\documentclass[reqno]{amsart}

\makeatletter

\renewcommand\subsection{\@startsection{subsection}{2}%
  \z@{-\baselineskip}{-0.3em}%
  {\normalfont\bfseries}}

\renewcommand\subsubsection{\@startsection{subsubsection}{3}%
  \z@{-\baselineskip}{-0.3em}%
  {\normalfont\bfseries}}

\makeatother

\usepackage{hyperref}

\usepackage{graphicx} 
\usepackage{amsmath}
\usepackage{amsthm}
\usepackage{hyperref}
\usepackage{amssymb}
\usepackage{cases}
\usepackage{multicol}
\usepackage{comment}
\usepackage{cases}
\usepackage{dsfont}
\usepackage{dsfont}
\usepackage{blindtext}
\usepackage{multirow}
\usepackage{geometry}
\usepackage{tikz-cd}
 
\usepackage[colorinlistoftodos]{todonotes}
\newcommand{\dd}{\mathrm{d}}
\newcommand{\DD}{\mathrm{D}}
\newcommand{\curl}{\mathrm{curl}}
\newcommand{\ddiv}{\mathrm{div}}

\newcommand{\IId}{\mathrm{Id}}

\newcommand{\supp}{\text{supp}}
\newcommand{\R}{\mathbb{R}}

\newcommand{\cn}{\cdot\nabla}

\theoremstyle{plain}
\newtheorem{theorem}{Theorem}[section]

\newtheorem{prop}[theorem]{Proposition}
\newtheorem{lemma}[theorem]{Lemma}

\theoremstyle{definition}
\newtheorem{definition}{Definition}[section]
\newtheorem{assumptions}{Standing Assumptions}
\newtheorem{remark}{Remark}

\counterwithin{equation}{section}

\title{$C^{1/5^{-}}$ Convex Integration Solutions of Ideal MHD}
\author{Matteo Giardi}
\thanks{MG was supported by the International Max Planck Research School Mathematics in the Sciences (IMPRS MiS)}
\author{L\'aszl\'o Sz\'ekelyhidi Jr.}
\thanks{LSz gratefully acknowledges the support of the Deutsche Forschungsgemeinschaft (DFG, German Research Foundation) through GZ SZ 325/2-1.}
\date{}

\begin{document}

\begin{abstract} 
\noindent For any $0\leq \gamma < 1/5$, we construct weak solutions $(v, B, p )$ of the Ideal MHD Equations which do not conserve the total kinetic energy, the cross-helicity and lie in $C^\gamma(\mathbb{T}^3\times\mathbb{R})$. In the spirit of Arnold's formulation of ideal hydrodynamics, a solution is thought of as a path of volume-preserving diffeomorphisms; the proof is then based on the interplay between classical convex integration techniques and geometric constructions at the level of the Lie algebra of this Lie group. Our work substantially extends the recent work of and building on the recent work of  Enciso, Peñafiel-Tomás and Peralta-Salas.
\end{abstract}

\maketitle

{
\makeatletter
\renewcommand{\l@section}{\@dottedtocline{1}{0em}{2.3em}}
\renewcommand{\l@subsection}{\@dottedtocline{2}{1.5em}{3em}}
\renewcommand{\l@subsubsection}{\@dottedtocline{3}{3em}{4em}}
\makeatother
\small
\tableofcontents
}

\section{Introduction}
We study the ideal incompressible magneto-hydrodynamic
equations (MHD) in three dimensions. This is a classical macroscopic model
describing the evolution of an electrically conducting fluids such as plasmas
and liquid metals. The MHD system describes the simultaneous evolution of a
\textit{velocity field} $v$ and a \textit{magnetic field} $B$ which are
divergence free. The evolution of $v$ is described by the Cauchy momentum
equation with an external force given by the Lorentz force induced by $B$.
Incompressibility is ensured by a pressure gradient $\nabla p$. In turn, the
evolution of $B$ is described by the induction equation which couples
Maxwell-Faraday law, with the \textit{electric field} $E$ given by Ohm's law:
\begin{equation}\label{FH}
    \begin{cases} 
        \partial_tB+\curl\ E=0,\\
        \ddiv \ B=0,\\
        E+v\times B=0.\\
    \end{cases}
\end{equation}
The resulting system of equations is then:
\begin{equation}\label{MHD}
    \begin{cases}
        \partial_tv+v\cdot\nabla v+\nabla p-B\cdot\nabla B=0,\\
        \partial_tB+\curl[B\times v]=0,\\
        \ddiv \ v=\ddiv \ B=0.\\
    \end{cases} \
\end{equation}
We consider this system in the periodic setting, in other words in the periodic spatial domain $\mathbb{T}^3=[-\pi,\pi]^3$, so that the unknowns are the vector fields $v, B:\mathbb{T}^3\times[0, T]\to\R^3$ and the scalar field $p:\mathbb{T}^3\times[0, T]\to\R$. We call \textit{weak solution} a triple $(v, B, p)$ satisfying \eqref{MHD} in the sense of distributions. 

\noindent In this paper we are concerned with \textit{weak solutions} of the MHD system, which are triples $(v,B,p)$ satisfying \eqref{MHD} in the sense of distributions. Our motivation comes from magneto-hydrodynamic turbulence and the search for an analogue of the famous K41 theory describing hydrodynamic turbulence. In the latter case a theory of weak solutions for the incompressible Euler equations, describing ideal hydrodynamics, arose from the work of Onsager in 1949 \cite{Onsager}. This work has received considerable attention in the past 20 years under the name of `Onsager's Conjecture', after its initial revival in the early 1990s \cite{Eyink1994,CET}. The conjecture was eventually resolved in the work of P. Isett \cite{Isett2018}, see also \cite{BuckmasterDeLellisSzekelyhidiVicol2019}, using the technique of convex integration. For a comprehensive review of this subject we refer to the review \cite{EyinkReview}. Subsequently, the technique has been adapted to various other fluid-dynamic models, including the ideal MHD setting in \cite{FLSz2021,BBV,FLSz2024,MiaoNieYe,EnPePe}. For a more in-depth review of these works, see Subsection \ref{ss:literature} below. In particular, our work is a substantial refinement of \cite{EnPePe} and can be stated as follows.

\begin{theorem}[Main Theorem]\label{main} For any $0\leq \gamma<1/5$, there exist, weak solutions $(v, B, p)$ of \eqref{MHD} with $u, B\in C^\gamma(\mathbb{T}^3_x\times \mathbb{R}_t)$ with nontrivial magnetic helicity, which do not conserve the total kinetic energy and the cross-helicity. 
\end{theorem}

\subsection{Context}\label{ss:context}
In order to explain the physical context and motivation of our work, let us first recall the closely related hydrodynamic case and weak solutions of the incompressible Euler equations. In a broad-brush picture (ignoring corrections due to intermittency), there are two classical questions defining the contours of the theory: the first is about the scale-invariant power law characterizing the cascade of kinetic energy - in the hydrodynamic context this is the famous $k^{-5/3}$ law of Kolmogorov (the K41 theory). The second is about the critical regularity required for validity of the ideal conservation of energy - in the hydrodynamic context this is precisely the question raised by L. Onsager, with the by now verified $C^{1/3}$ threshold on the H\"older scale (let us call this the O49 theory). Sufficiency of $1/3^+$ H\"older regularity was proved in \cite{CET}, whereas the necessity was shown in the work of \cite{Isett2018}. It appears to be a stroke of coincidence that the K41 exponent $-5/3$ and the O49 exponent $1/3$ are both uniquely determined and agree in terms of scaling. In particular, the K41 theory does not require any dynamical insight into what is going on in the inertial range, but the need for such an insight would arise if dimensional analysis fails to give a unique answer - as it happens for MHD turbulence below. We also remark that the O49 exponent also conveniently agrees with the natural bound reachable by convex integration as introduced by Nash \cite{Nash1954}, see \cite{DLSz2016} for a discussion of this fact. 

Now let us turn to MHD turbulence. Since any solution of the incompressible Euler equations is trivially a solution of the ideal MHD system (with zero magnetic field $B$), a fortiori both the $k^{-5/3}$ spectrum and the $C^{1/3}$ threshold appear as candidates. However, the presence of a large-scale non-zero magnetic field $B_0$ substantially changes the situation: firstly, it introduces an inherent anisotropy, in stark contrast with one of the basic assumptions underlying the K41 theory, leading to a distinction between $k_{\parallel}$ and $k_{\perp}$, representing the components of the wave-vector $k$ of fluctuations of the turbulent fields parallel and perpendicular to $B_0$. Secondly, it leads to the presence of Alfv\'en waves, propagating in the direction $k$ with a characteristic speed proportional to $B_0\cdot k=|B_0|k_\parallel$ - see Subsection \ref{overview} below. In turn, this leads to the introduction of a second characteristic time scale, the Alfv\'en time $\tau^A\sim (v_Ak_\parallel)^{-1}$, where $v_A=\frac{B_0}{\sqrt{4\pi\rho_0}}$ denotes the Alfv\'en speed and $\rho_0$ is the plasma density (assumed to be constant). Both the anisotropy and the additional characteristic scale mean that additional hypothesis are needed to deduce scaling laws on the energy spectrum. Kraichnan \cite{Kraichnan1965} argued that in homogeneous isotropic MHD turbulence one should have equipartition of energy in the inertial range in both kinetic and magnetic energy spectra and the scaling law should depend on $v_A$ - leading to a $k^{-3/2}$ spectrum (known as the Iroshnikov-Kraichnan spectrum). In contrast, in \cite{Goldreich1995} Goldreich and Sridhar argued that the strong anisotropy should be reflected by different scaling laws for the $k_{\parallel}$ and $k_{\perp}$ spectra, and in particular, assuming that the Alfv\'en time equals the eddy turnover time (known as the critical balance hypothesis), obtained $E(k_\perp)\sim \epsilon^{2/3}k_\perp^{-5/3}$. The nature of the inertial range in plasma turbulence remains a subject of intensive debate, we refer to the excellent survey \cite{Schekochihin}.

In terms of analytic properties, it is known that a non-vanishing large-scale magnetic field and the resulting hyperbolic structure can have a regularizing property on solutions, we refer to \cite{BSS1988,ZHuali2026}. Further, ideal MHD with non-vanishing magnetic field possesses two more ideal conserved quantities besides the total energy, namely cross-helicity and magnetic helicity (and in fact these are all the ideal invariants, see \cite{KhesinPeraltaYang}). We recall that the \textit{total energy} and the \textit{cross-helicity} are given by
\begin{align*}
	\mathcal{E}_{(v,B)}(t)&:=\int_{\mathbb{T}^3}\left(|v(x,t)|^2+|B(x,t)|^2\right)\dd x,\\ 
	\mathcal{H}^\times_{(v,B)}(t)&:=\int_{\mathbb{T}^3}B(x,t)\cdot v(x,t)\dd x,
\end{align*}
whereas \textit{magnetic helicity} can be defined via a vector potential $A$ such that $\curl A=B$ as
\begin{equation*}
	\mathcal{H}_{B}(t):=\int_{\mathbb{T}^3}A(x,t)\cdot B(x,t)\dd x.
\end{equation*}
Magnetic helicity was first studied by Woltjer in~\cite{Woltjer} and interpreted topologically as the linkage and twist of magnetic field lines in the highly influential work of Moffatt~\cite{Moffatt}, see also~\cite{ArnoldKhesinBook}. In analogy with the K41 hypothesis of anomalous dissipation in the infinite Reynolds number limit, in MHD turbulence the rate of total energy dissipation in viscous, resistive MHD seems not to tend to zero when the Reynolds number and magnetic Reynolds number tend to infinity, see~\cite{DallasAlexakis,Linkmann2015nonuniversality,Mininni2009finite}. On the other hand magnetic helicity is a rather robust conserved quantity even in turbulent regimes, and J.B.Taylor conjectured that  magnetic helicity is approximately conserved for small resistivities~\cite{Taylor1974} - this conjecture was recently verified in \cite{FaracoLindberg-Taylor}. 

The analogue of the Onsager conjecture, i.e.~critical threshold for conservation of the above quantities was first addressed by Caflisch, Klapper and Steele in \cite{CKS}, building upon the work of Constantin, E and Titi \cite{CET}, and endpoint cases were subsequently studied in \cite{KangLee2007}. It turns out that total energy and cross-helicity are conserved, provided $v\in C([0, T], C^{\gamma_1}(\mathbb T^3))$ and $B\in C([0, T], C^{\gamma_2}(\mathbb T^3))$ with $\gamma_1+2\gamma_2>0$, whereas conservation of magnetic helicity merely requires $(v, B)\in L^3(\mathbb{T}^3\times(0, T))$. Thus, a natural conjecture, put forward for instance in \cite{rev2}, is the following direct analogue of Onsager's conjecture: weak solutions $v,B\in C([0, T], C^{\gamma}(\mathbb T^3))$ conserve total energy and cross-helicity if $\gamma>1/3$, but may fail to do so if $\gamma<1/3$. 

\subsection{Overview of the Literature}\label{ss:literature}
 
As already pointed out above, any weak solution of the 3D incompressible Euler equations is trivially a weak solution of \eqref{MHD} with vanishing magnetic field $B\equiv 0$, so that Onsager's conjecture with exponent $1/3$ is valid in this special case. More generally, weak ``$2\frac{1}{2}$-dimensional'' solutions have been constructed in~\cite{Bronzi2015} Bronzi, Lopes Filho and Nussenzveig Lopes in the class $L^\infty$ and, recently by Miao, Nie, Ye in \cite{MiaoNieYe} in the class $C^{1/3-}(\mathbb T^3\times\R)$. These are solutions of the form  $v(x_1,x_2,x_3,t) = (v_1(x_1,x_2,t),v_2(x_1,x_2,t),0)$ and $B(x,t) = (0,0,b(x_1,x_2,t))$. Under this structural restriction the 3D MHD system \eqref{MHD} reduces to the 2D Euler equations for $(v_1,v_2)$ coupled with a passive scalar equation for $b$. Evidently, for such solutions magnetic helicity as well as cross-helicity are identically zero. 

The first results for genuinely 3D weak solutions appeared in the work of Beekie, Buckmaster and Vicol \cite{BBV} as well as in the work of Faraco, Lindberg and the second author \cite{FLSz2021}. In \cite{BBV} the authors constructed unbounded weak solutions in the energy class $v,B\in L^\infty_t L^2_x$ which do not preserve magnetic helicity nor energy or cross-helicity. Compared to the thresholds for conservation of energy and magnetic helicity obtained in \cite{CKS,KangLee2007} mentioned above, we see that the space $L^\infty_t L^2_x$ is super-critical with respect to magnetic helicity conservation. In view of Taylor's conjecture proved in \cite{FaracoLindberg-Taylor}, such solutions cannot arise in the infinite Reynolds number limit and thus are not valid representatives of fully developed turbulence. In \cite{FLSz2021}
the authors constructed bounded (but discontinuous) weak solutions which do not preserve energy and cross-helicity - but of course, a fortiori, magnetic helicity is conserved, and indeed vanishes identically. In subsequent work \cite{FLSz2024} the authors extended their construction to bounded weak solutions with non-vanishing (constant) magnetic helicity. In both works \cite{FLSz2021,FLSz2024} a key point was to introduce a relaxation of the MHD system \eqref{MHD} which decouples the effects of hydrodynamic turbulence (in the form of the appearance of a Reynolds stress term in the momentum equation) from possible small scale effects in the Faraday-Ohm system \eqref{FH} in a manner consistent with magnetic helicity robustness and conservation. In a nutshell, the relaxation involves replacing Ohm's law (the third equation in \eqref{FH}) by the nonlinear constraint $E\cdot B=0$. 

Very recently, Enciso, Peñafiel-Tomás and Peralta-Salas \cite{EnPePe} succeeded in constructing H\"older continuous weak solutions that do not conserve energy, cross-helicity, and have non-trivial magnetic helicity. The weak solution obtained is in the class $v,B\in C^\gamma(\mathbb T^3\times [0, T])$ with $\gamma\leq 1/200$. Their construction is based on convex integration applied to a more restrictive, partial relaxation of \eqref{MHD}, where the Faraday-Ohm system \eqref{FH} is solved exactly, leaving only the appearance of a Reynolds stress term in the momentum equation as the only small-scale effect, see \eqref{relaxedMHD}. In this way the magnetic field $B$ is ``slaved" to the velocity via Lie transport. Our paper is bulding on this same ansatz and can be considered as a substantial extension of the work \cite{EnPePe}, improving the exponent $1/200$ to $1/5$. In subsection \ref{overview} we recall the main ideas of this construction and point out the main new ideas and techniques introduced in our work. 

In the remaining part of the introduction, we will give an overview of the convex integration scheme implemented in this work (Subsection \ref{overview}), then state the main Iterative Proposition and the Iterative Assumptions, and provide a flowchart of its proof (Subsection \ref{iteratives}). Finally, in Subsection \ref{choiceofparameters} we discuss the parameters and set their values.


\subsection{Overview of the Construction} \label{overview} In this section we briefly sketch the main ideas in the approach used both in our work and the work \cite{EnPePe} and point out the key difficulties and points of departure. For the notation and differential-geometric concepts we refer to Appendix \ref{diffgeom}. 

\noindent \textbf{Alfv\'en Waves.} Let us start by considering small fluctuations around a constant state $(\bar{v}+w,\bar{B}+b)$, and linearizing. We obtain
\begin{equation*}
    \begin{split}
        \partial_tw+\bar v\cn w-\bar B\cn b+\nabla p&=0,\\
        \partial_tb+\bar v\cn b-\bar B\cn w&=0.
    \end{split}
\end{equation*}
Equivalently written for the Els\"asser variables $z^{\pm}=w\pm b$:
\begin{align*}
\partial_tz^++\bar{v}\cdot\nabla z^+-\bar{B}\cdot z^++\nabla p&=0,\\
\partial_tz^-+\bar{v}\cdot\nabla z^-+\bar{B}\cdot z^-+\nabla p&=0,
\end{align*}
which admit plane-wave solutions of the form
\begin{equation*}
    z^{\pm}=\delta^{1/2}e^{i\lambda(x-t\bar{v}\pm t\bar{B})\cdot k}\zeta^{\pm},
\end{equation*}
with amplitude $\delta^{1/2}$, wave vector $k$ and some state vectors $\zeta^\pm$ with $\zeta^{\pm}\cdot k=0$. Whilst one can assume $\bar{v}=0$ by Galilean invariance, the differing wave velocities $\pm\bar{B}\cdot k$ reflect fundamental properties of Alfv\'en wave packets, most pronounced in the case $k\parallel \bar{B}$. On the other hand the nonlinear term in Els\"asser variables takes the form $\textrm{div}(z^+\otimes z^-)$, and recall that a key aspect of convex integration constructions is that spatial averages of the nonlinear term are able to balance the Reynolds term. That is, we would like to have $\langle z^+\otimes z^-\rangle=R$, uniformly in time. However, this is only possible if the Alfv\'enic wave structure is suppressed, i.e.~provided $k_\parallel=0$. This leads to the ansatz 
\begin{equation}\label{e:ansatz1}
w(x,t)=\delta^{1/2}h(\lambda(x-t\bar{v})\cdot k)\zeta,\quad\textrm{ with }k\cdot \bar{B}=0,\,k\cdot \zeta=0.
\end{equation}
In the case where $\bar{v}$ and $\bar{B}$ are no longer constant but slowly varying, the plane-wave ansatz above needs to be replaced by $e^{i\lambda\phi(x,t)}\zeta^{\pm}$ with phase function $\phi(x,t)$ which needs to simultaneously solve
$$
\partial_t\phi+\bar{v}\cdot\nabla\phi=0,\textrm{ and }\bar{B}\cdot\nabla \phi=0.
$$
The differential-geometric structure that makes this construction possible at least locally in space-time, when $(\bar B, \ \bar v)$ are not constant anymore, is the commutation relation
$$[\partial_t+\bar v\cn,\bar B \cn]=\partial_t \bar B+\bar v\cn \bar B-\bar B\cn \bar v=0,$$
which is a consequence of \eqref{FH}. In terms of Els{\"a}sser variables and the Alfv\'en transport operators 
$$\mathcal{A}^\pm=\partial_t+\bar z^\pm\cn, \ \ \ \ \bar z^\pm=\bar v\pm\bar B,$$
the above is equivalent to
\begin{equation}\label{commutationintro}
    \mathcal{A}^+\mathcal{A}^-=\mathcal{A}^-\mathcal{A}^+.
\end{equation}
We will use this commutation relation to reduce locally any slowly varying solution of \eqref{FH} $\bar{v},\bar{B}$ to the constant case above, introducing a local chart construction, this is the content of Section \ref{chart}. The price to pay is that we will need to keep track not only of Lagrangian transport derivatives (as in the Euler case), but also of the Alfv\'en transport derivatives; see the iterative assumptions \eqref{inductiveassumptionsgeneral}. The parameter $\tau^c\sim 1/(\lambda_q\delta_q^{1/2})$ will quantify the size of the support of the space-time cutoffs used.

\noindent \textbf{Partial Relaxation and Lie-Taylor Expansion.} 
The first step in any convex integration construction is to identify a suitable relaxed system of equations. This is usually achieved by local averaging or filtering, and is also a key process in identifying the threshold regularity for energy conservation \cite{CET,Eyink1994,EyinkReview,AluieEyink}. Such filtering process identifies additional `defect' terms arising from the effect of small-scale fluctuations on the large scale dynamics via nonlinear interactions. In the case of \eqref{MHD} such relaxation takes the form
\begin{equation}
    \begin{cases}
        \partial_t\bar v+\ddiv\left[\bar v\otimes \bar v-\bar B\otimes \bar B\right]+\nabla \bar p=\ddiv R,\\
        \partial_t \bar B+\curl [\bar B\times \bar v]=\curl \ M,\\
        \ddiv \ \bar v= \ddiv \bar B =0
    \end{cases}
\end{equation}
where $R$ is a 2-tensor (usually called the Reynolds stress) and $M$ is a vector field (called the subscale electromotive force in \cite{AluieEyink}). An additional feature not seen by pure filtering, but resulting from the differential (div-curl) structure (first observed by L. Tartar \cite{Tartar}) is that, at least in the presence of large scale-separation, $M$ necessarily satisfies a geometric constraint of the type:
$$
M\cdot \bar B=0.
$$
This constraint is a local expression of magnetic helicity conservation \cite{FLSz2024}. In \cite{EnPePe} the authors consider a more restrictive, partial relaxation, where $M\equiv 0$, so that $\bar{B}$ is given exactly as the solution of the Faraday-Ohm system \eqref{FH}, i.e. the Lie transport of the initial field via the vector field $\bar{v}$.  As is well known, the solution can be then written as
$$
\bar B(\cdot,t)=(\bar X_t)_*\bar B|_{t=0},
$$
see \eqref{duhamel}, where $\bar X_t$ denotes the Lagrangian flow of $\bar v$. Thus a natural way to produce a perturbation is by precomposition of $\bar X$ with a given volume preserving map $X$ close to the identity - this is done in \cite{EnPePe}. In our work, instead of constructing $X$ directly, we construct it as the time one map of a vector field $\xi$, that is $X=\exp_{\IId}^\xi$. Motivated by the work of Vladimirov, Moffat Ilin \cite{VMI}, we call $\xi$ Lagrangian Displacement Field (LDF). Then, similarly to \cite{VMI}, we define the perturbed Lagrangian map as $X'=X\circ \bar X=\exp_{\IId}^\xi\circ \bar X$.
A key point is then to ensure good control of the flow map over long time scales, and this is done by introducing a fast time-scale $\tau^a<\tau^c$ on which $\xi$ (and consequently $w$) will oscillate. To see how this works, we first observe that the perturbation $w$ of $\bar v$ associated with a Lagrangian Deformation Field $\xi$ is given by: 
\begin{equation}\label{wintro}
\begin{split}
    w&=\partial_tX\circ X^{-1}+(X_*-\IId_*)\bar v\\
    &=(\partial_t+\mathcal{L}_{\bar{v}})\xi+\sum_{k=1}^{k_0}\frac{(-1)^k}{(k+1)!}\mathcal{L}_{\xi}^k(\partial_t+\mathcal{L}_{\bar{v}})\xi+r^{k_0}_w,
\end{split}
\end{equation}
which we call a Lie-Taylor expansion for $w$. In particular, a first approximation for $\xi$ corresponding to $w$ in \eqref{e:ansatz1} is given by 
\begin{equation}\label{e:ansatz2}
\xi(x,t)=\delta^{1/2}\tau^ag^{[1]}(t/\tau^a)h(\lambda(x-t\bar v)\cdot k)\zeta,
\end{equation}
leading to 
\begin{equation}\label{e:ansatz3}
    w(x,t)=\delta^{1/2}g(t/\tau^a)h(\lambda(x-t\bar v)\cdot k)\zeta+O(\tau^a)
\end{equation}
with $g=(g^{[1]})'$.
We also remark at this stage that this type of velocity perturbation is also compatible with the Nash-Newton iteration introduced by Giri-Radu in  
\cite{GR}, which is useful when a Mikado-type building block is unavaliable and one has to resort to simple plane-waves. The case in which $(\bar v, \bar B)$ are not constant is dealt with a chart construction, see Lemma \ref{chartprop}, and the geometric decomposition in Lemma \ref{geomlemma}.

\noindent Our construction pairs the LDF and chart construction toolbox mentioned above with this observation. Compared to \cite{EnPePe}, the use of a LDF moves the problem from the Lie Group of volume-preserving diffeomorphisms to its Lie-algebra of divergence-free vector fields and effectively transforms compositions into addition, allowing for better estimates in terms of the Lie derivatives of the various objects. The price to pay is to work with Lie-Taylor expansions \eqref{wintro}, which are explicit in terms of $\xi$ up to an arbitrarily small remainder. The study of such objects is the content of Section \ref{lt}, and the key tool will be the Lagrangian Perturbation Lemma \ref{lpl}. 

\noindent The parameters $k_0^g, \ k_0^p$ tracking the depth of the expansions, in the two stages of our construction, the Galbrun stage and Nash stage, and $\mathcal{T}_g, \ \mathcal{T}_p$ measuring the gain corresponding to each additional order of expansion in the Lie-Taylor series, \eqref{wintro}, will play a fundamental role. Typical estimates will look like:
\begin{equation*}
    \begin{split}
        &||(\partial_t+\mathcal{L}_{\bar{v}})\xi||_0\lesssim \delta^{1/2},\\
        &||\mathcal{L}_{\xi}^k(\partial_t+\mathcal{L}_{\bar{v}})\xi||_0\lesssim \mathcal{T}^k\delta^{1/2},\\
        &||r_w^{k_0}||_0\lesssim \mathcal{T}^{k_0+1}\delta^{1/2}
    \end{split}
\end{equation*}
where the gain $\mathcal{T}$ is heuristically given by:
\begin{equation}\label{Tintro}
    \mathcal{T}\simeq (\text{good derivative})\cdot ( C^0 \text{ size of } \xi)\ \text{ e.g.} \ \mathcal{T}_p=\lambda_q\tau^a\delta_{q+1}^{1/2}.
\end{equation}
For the exact estimates we are referring to, see, for example, Lemma \ref{estimatesperturbationnash}. The order $k_0$ will then be chosen so that, despite not being able to control the geometry of the remainder, the factor $\mathcal{T}^{k_0+1}$ is so small that we can still ensure transport properties matching the ones of the explicit terms. The technical tools needed to prove such bounds are developed in Section \ref{inductives} and are collected in what we call the Inductive Lemma \ref{inductive}. 

\noindent 

\noindent \textbf{Galbrun's Equation and Stages.} Since the geometry of the problem mentioned above effectively reduces the available types of oscillations to 2D constructions in the plane orthogonal to $\bar B$, which cannot be made disjoint in space, we pair the tools just described with a preliminary step in which we `well prepare' the Reynolds stress to look like a sum of disjoint simple components: 
\begin{equation}\label{wellprepared}
    \sum_I g^2_IA_I,
\end{equation}
as in Giri-Radu \cite{GR}, where $g_I$ are fast oscillating periodic disjoint $L^2$ normalized functions of time. Note that this goes hand in hand with the observations from the previous paragraph about fast-time oscillations. Pairing the LDF construction with having to solve a linearised ideal MHD equation at $(\bar v, \bar B, \bar p)$ with right-hand side: 
$$F=\ddiv \sum_I(1-g_I^2)A_I,$$ 
as in \cite{GR}, results in a second-order equation for the LDF $\xi$ and pressure $\pi$, which is similar to Galbrun's Equation in acoustics, see also the work of Lindblad \cite{Lindblad2005FreeBoundary} on the free boundary problem for the Euler equation. Namely,
\begin{equation}\label{Galbrunintro}
\begin{cases}
    \mathcal{A}^+\mathcal{A}^-\xi+\xi \cn (\nabla\bar p)+\nabla \pi=F,\\
    \ddiv \ \xi =0.
\end{cases}
\end{equation}
The commutation of the Alfv\'{e}n transport operators \eqref{commutationintro} will be crucial for propagating $C^{r+\alpha}$ H{\"o}lder estimates. The study of the solutions of this equation, in fact, for the equation satisfied by the 1-form potential $\Theta$ of $\xi$ is the content of Section \ref{gn}. We implement only a single step in the so-called Newton iteration from \cite{GR}, which we call the Galbrun Stage. This will be followed by the Nash Stage, in which we correct the Reynolds stress using the quadratic term, as it is standard in convex integration schemes. 

\noindent \textbf{Mollification Issues.} A key gadget in every convex integration scheme is mollification. In our setting, this is problematic as it is an operation that does not commute with the Faraday-Ohm system \eqref{FH}. To address this, we keep the non-mollified objects as they are, but construct, in each Stage, the perturbations out of mollified objects only. To reinstate the desired commutation properties, we locally correct the mollified vector fields by solving the (partially relaxed) MHD system exactly, locally in time, with initial data the mollified fields, see \eqref{localcorreq}. Stability estimates will ensure that this correction is small enough. Schematically, this reads:
$$\begin{cases}
    v_q\\
    B_q\\
\end{cases}\overset{\text{mollification}}{\longrightarrow} \begin{cases}
    v_\ell\\
    B_\ell\\
\end{cases} \overset{\text{local in time recorrection}}{\longrightarrow}\begin{cases}
    v_{\ell,j}=v_\ell+\text{small}\\
    B_{\ell,j}=B_\ell +\text{small}\\
\end{cases}  $$
$$\overset{\text{LDF construction}}{\longrightarrow}\begin{cases}
    \xi\\
    X=\exp_{\IId}^\xi
\end{cases} \overset{\text{perturbation in Lagrangian coordinates}}{\longrightarrow} \begin{cases}
    v_q+w=\partial_tX\circ X^{-1}+X_*v_q\\
    B_q+b=X_*B_q\\
\end{cases}$$
In the terminology of \cite{BuckmasterDeLellisSzekelyhidiVicol2019}, this can be viewed as a sort of effective glueing procedure that is performed once per Stage.

\noindent This adds in the Lie-Taylor expansion \eqref{wintro} above with $\bar v=v_q$ error terms of the form: 
\begin{equation}\label{nnmoll}
    \mathcal{L}^k_{\xi}(v_q-v_\ell),
\end{equation}
these are still not mollified, and we thus control only a finite number of derivatives.

\noindent This is ultimately why we can propagate only $M\ll N$ good derivatives of the Reynolds stress and the fields' transport properties, see \eqref{inductiveassumptionsgeneral}. To tackle this issue, we use the fact that pure derivatives do not require such Lie-Taylor expansions to be read from the construction, and we can propagate as many as needed.

\noindent Assume we want to bound $0\leq r \leq M$ transport derivatives of $w$ as in \eqref{wintro} with the additional errors in \eqref{nnmoll} and we have $M$ space derivativatives on the transport derivatives of $v_q-v_\ell$ and $N$ space derivatives on the pure time-derivative of $v_q-v_\ell$. The key idea is that, for $r$ small, we can ensure the properties we want by expanding and choosing $M\gg k_0$, while for large $r$, we first ensure that we still have enough good derivatives, by choosing $N-k_0\gg M$, we then use the weaker pure time derivative bound instead of the transport derivative bound, and still win using the ratio: $$(\text{good derivative}/\text{ bad derivative})^{r}=(\lambda_q/\lambda_{q+1})^r,$$
as $r$ will be large enough. This kind of derivative and expansion level bookkeeping is the most technical part of this work. As an example, the need to keep track of second-order pure time derivatives on $(v_q, \ B_q)$ corresponds to the fact that for high space derivatives ($\geq M-1$), we have no information on $\mathcal{A}_q^\pm R_q$, see \eqref{inductiveassumptionsgeneral}.

\noindent \textbf{Regularity Cap.} We now give a heuristic explanation for the claimed $1/5^-$ regularity. We first note that, as in \cite{GR}, the need to cut off in time the solution of \eqref{Galbrunintro} introduces an error of size:
$$\frac{1}{\lambda_q}\frac{(\text{fast time-scale})}{(\text{cut-off size})}(\text{size of r.h.s.})=\frac{\tau^a}{\tau^c}\delta_{q+1},$$
here $\tau^a$ is the fast time-scale introduced above and $\tau^c$, for the time being can be thought of as the eddy turnover time $1/\lambda_q\delta_q^{1/2}$, the size of the right-hand side is $\lambda_q\delta_{q+1}\sim \ddiv A_I$ and the $1/\lambda_q$ gain happens because as in \cite{GR} we work at the level of the potential of $\xi$ and not directly on \eqref{Galbrunintro}. In \cite{GR}, the procedure is iterated to remove this error; here, we keep it as it is.

\noindent Secondly, going back to the simplified constant case in \eqref{e:ansatz2} above, together with the relation between the leading term of $w$ and $\xi$ in \eqref{wintro}, we end up with an ansatz for $\xi$ of the form:
\begin{equation}\label{xiintro}
    \xi(x,t)=\delta_{q+1}^{1/2}\chi(t) \tau^ag^{[1]}(t/\tau^a)h(\lambda_{q+1} (x-t\bar v)\cdot k)\zeta,
\end{equation}
where $(g^{[1]})'=g$, and $\chi$ is a cut of in time (in the actual construction will be in time and space) whose support has size $\tau^c$. From this ansatz, we get:
\begin{equation}\label{wfromxi}
    w\simeq(\partial_t+\mathcal{L}_{\bar v})\xi=\delta_{q+1}^{1/2}\chi g h \zeta+ \tau^ag^{[1]}\partial_t\chi \delta_{q+1}^{1/2} h \zeta=w^o+w^c,
\end{equation}
see \eqref{leadingtermnashexpansion} for the actual computation. On the one hand, assuming the normalization $\int h^2=1$ we have the required balance in the principal quadratic term:
$$
w^o\otimes w^o=\delta_{q+1}\chi^2 g^2 h^2 \zeta\otimes \zeta=\delta_{q+1}^2\chi^2 g^2 (h^2-1) \zeta\otimes \zeta+\delta_{q+1}^2\chi^2 g^2\zeta\otimes \zeta,
$$
which has space-average matching exactly the rewritten Reynold stress \eqref{wellprepared} as required. However, the next order terms lead to quadratic error terms of the form
$$
||w^o\otimes w^c+w^c\otimes w^o||_0 \lesssim \frac{\tau^a}{\tau^c}\delta_{q+1}.
$$
Happily, these are of the same size as the cut-off error above.

\noindent Lastly, we see that when the transport operator in the ideal MHD hits the fast time oscillating $g$ we lose a $1/\tau^a$, and upon inverting the divergence to rewrite it as a new Reynold stress, we expect this error to have size:
$$||\ddiv^{-1}(\partial_t+\bar v\cn)w||_0\lesssim \frac{\delta_{q+1}^{1/2}}{\lambda_{q+1}\tau^a}=(\tau^c/\tau^a)\frac{\lambda_q\delta_q^{1/2}\delta_{q+1}^{1/2}}{\lambda_{q+1}}.$$

\noindent Assuming no other larger term appears between the various errors in the construction, this leaves us with a new Reynolds stress of size:
$$(\tau^a/\tau^c)\delta_{q+1}+(\tau^c/\tau^a)\frac{\lambda_q\delta_q^{1/2}\delta_{q+1}^{1/2}}{\lambda_{q+1}},$$
optimising this expression in the ratio $(\tau^a/\tau^c)$ imposing the required smallness and using the usual relations $\lambda_q=a^{b^q}, \delta_q=\lambda_q^{-2\beta}$, we end up with:
$$(\tau^a/\tau^c)\delta_{q+1}+(\tau^c/\tau^a)\frac{\lambda_q\delta_q^{1/2}\delta_{q+1}^{1/2}}{\lambda_{q+1}}= 2\sqrt{\frac{\lambda_q\delta_q^{1/2}\delta_{q+1}^{3/2}}{\lambda_{q+1}}}\leq \delta_{q+2} \Longrightarrow \beta < \frac{1}{(4b+1)}.$$


\subsection{Subsolutions and Iterative Proposition}\label{iteratives} In this Subsection, we specify the precise convex integration framework we adopt and state the estimates we require on the sequence of subsolutions we are about to construct. According to the observations in Subsection \ref{overview}, we consider the following relaxation of the MHD system:
\begin{equation}\label{relaxedMHD}
    \begin{cases}
        \partial_tv+\ddiv\left[v\otimes v-B\otimes B\right]+\nabla p=\ddiv R,\\
        \partial_t B+\curl [B\times v]=0,\\
        \ddiv \ v= \ddiv B =0.
    \end{cases}
\end{equation}
We call a quadrupole $(v, \ B, \ p, \ R)\in C^\infty(\mathbb{T}^3\times\mathbb{R})$ solving \eqref{relaxedMHD} a \textit{subsolution}. Moreover, we say that a sequence of subsolutions $\{(v_q, \ B_q, \ p_q, \ R_q)\}_q$  is \textit{admissible}  if the following iterative assumptions are satisfied:
\begin{subequations}\label{inductiveassumptionsgeneral}
    \begin{align}
        &||v_q||_0, \ ||B_q||_0\leq C_0(1-\delta_q^{1/2}),\\
        &||\partial_t^j v_q||_r, \ ||\partial_t^j B_q||_r\leq C_0 \lambda_q^{r+j}\delta_q^{1/2} \ \text{ for } j=0,1,2 \ \ \& \ \ 0\leq r\leq N-j \ \ \& \ \ (r,j) \ \neq (0,0),\\
        &||\mathcal{A}^\pm_q v_q||_r, \ ||\mathcal{A}^\pm_q B_q||_r\leq C_0 \lambda_q^{r+1}\delta_q \ \text{ for } 0\leq r\leq M-1,\\
        &||p_q||_r\leq C_0^2 \lambda_q^{r}\delta_q \ \text{ for } 1\leq r\leq M,\\
        &||\mathcal{A}^\pm_q p_q||_r\leq C_0^2 \lambda_q^{r+1}\delta_q^{3/2} \ \text{ for } 1\leq r\leq M-1,\\ 
        &||R_q||_r\leq  \lambda_q^{r-\alpha}\delta_{q+1} \ \text{ for } \ 0\leq r\leq M,\label{iterativeRq}\\
        &||\mathcal{A}^\pm_q R_q||_r\leq  \lambda_q^{r+1-\alpha}\delta_q^{1/2}\delta_{q+1} \ \text{ for } \ 0\leq r\leq M-1,\\ 
        &\supp_t R_q\subset ((1+\delta_q^{1/2})1/2,(1-\delta_q^{1/2})5/2),\\
        &|B_q|\geq c_0(1+\delta_q^{1/2})
    \end{align}
\end{subequations}
where the various parameters will be defined in Subsection \ref{choiceofparameters}, $c_0, C_0$ are constants independent of $q$ and 
$$\mathcal{A}^\pm_q=\partial_t+z_q^\pm\cn \ \text{ with } \ z_q^\pm=v_q\pm B_q.$$
The following proposition, which is at the heart of every convex integration scheme, states that given a subsolution satisfying the inductive assumptions at step $q$, it is possible to construct a new subsolution satisfying the inductive assumptions at step $q+1$, effectively constructing an admissible sequence.

\begin{prop}[Iterative Proposition]\label{iterative} Given $0<\beta< 1/5$ and $\ b,\ \gamma_a, \ \gamma_\ell, \ \gamma_{CZ}, \ m_0, \ k_0^g, \ k_0^p, \ N, \ M$ from Subsection \ref{choiceofparameters} satisfying the constraints in Lemma \ref{parameterconstraints}, there exist: 
$$\bar a>1, \ \bar C_0>1$$ 
such that, for any $a> \bar a, \ C_0> \bar C_0$ the following holds:
for any subsolution  $(v_q, \ B_q, \ p_q,\ R_q)$ satisfying the iterative assumptions \eqref{inductiveassumptionsgeneral}, there exists another subsolution $(v_{q+1}, \ B_{q+1}, \ p_{q+1},\ R_{q+1})$ satisfying \eqref{inductiveassumptionsgeneral} with $q$ replaced by $q+1$ together with the estimates:
    \begin{equation}\label{iterativeeq}
        \begin{split}
            &||\partial_t^j(v_{q+1}-v_q)||_r, \ ||\partial_t^j(B_{q+1}-B_q)||_r\leq C_0\lambda_{q+1}^{r+j}\delta_{q+1}^{1/2} \ \text{ for }\ j=0,1,2 \ \text{ and }\ 0\leq r \leq N-j,\\
            &\supp_t(v_{q+1}-v_q), \supp_t(B_{q+1}-v_q) \subset (1/2, \ 5/2).
        \end{split}
    \end{equation}
\end{prop}
\noindent The proof of the Main Theorem \ref{main} follows easily from this one once the sequence is properly initialised, see Subsection \ref{proofmain}. 

\noindent \textbf{Flowchart of the Proof of the Iterative Proposition.} We split the proof of the Iterative Proposition \ref{iterative} as follows. We invite the reader to visit Subsection \ref{overview} where the ideas and part of the notation are explained first.

\noindent \textit{Galbrun Stage - Section \ref{Galbrun}.} The goal of this first stage is to rewrite the Reynolds stress $R_q$ in space-time disjoint components.
\begin{itemize}
    \item \textit{Subsection} \ref{moll1}: we mollify and local-in-time recorrect the vector fields to ensure the commutation of the Alfv\'{e}n transport operators. 
    $$
    (v_q, \ B_q, \ p_q)
\overset{\text{mollification}}{\longrightarrow} \begin{cases}
    v_{\ell,j}=v_\ell+\text{small}\\
    B_{\ell,j}=B_\ell +\text{small}\\
    p_{\ell,j}=p_\ell+\text{small}\\
\end{cases} \leadsto \mathcal{A}^+_{\ell,j}\mathcal{A}^-_{\ell,j}=\mathcal{A}^-_{\ell,j}\mathcal{A}^+_{\ell,j}$$
    \item \textit{Subsection} \ref{geomconstr}: we construct space-time local charts $\Psi_I$ adapted to $(v_{\ell,j}, B_{\ell,j})$ together with a space-time partition of unity $\{\chi_I^2\}_I$ and use these to decompose the mollified Reynold stress in simple tensors. $$ \delta_{q+1}\left(\IId-\frac{R_\ell}{\delta_{q+1}}\right)=\sum_Ia_I^2\chi_I^2\Psi^{2*}_I\zeta\otimes\Psi^{2*}_I\zeta=\sum_IA_I$$ 
    \item \textit{Galbrun LDF} \ref{gnldf}: we solve, locally in time, the linearised MHD systems with a fast-time oscillating right-hand side in terms of a first LDF $\xi^{g}=\sum_j\xi^g_j$ and pressure $\pi^g=\sum_j\pi^g_j$, namely 
    \begin{equation}
        \begin{cases}
    \mathcal{A}^+_{\ell,j}\mathcal{A}^-_{\ell_j}\xi_j^g+\xi_j^g \cn [\nabla p_{\ell,j}-\ddiv \ R_\ell]+\nabla \pi_j^g=\ddiv\left[\sum_{I:I_t=j}(1-g_I^2)A_I\right],\\
    \ddiv \ \xi_j^g =0
        \end{cases}
    \end{equation}
    and rewrite the decomposition $\sum_IA_I$ as a sum of disjoint parts $\sum g_I^2A_I$. As in \cite{GR}, this will be done at the level of the 1-form potential $\Theta^g_j$ of $\xi^g_j$, see \eqref{galbrunapplied}. We finally perturb $(v_q, B_q, p_q)$ to $(\tilde v_q, \tilde B_q, \tilde p_q)$ in Lagrangian coordinates.
    $$
        ( v_{\ell,j}, \ B_{\ell,j}, \ p_{\ell,j})+\begin{cases}
    \xi^g\\
    X^g=\exp_{\IId}^{\xi^g}\\
\end{cases} \overset{\text{Galbrun stage}}{\longrightarrow}\begin{cases}
    \tilde{v}_q=\partial_tX^g\circ (X^g)^{-1}+X^g_*v_q&=v_q +w^{g}\\
    \tilde{B}_q=X^g_* B_q&=B_q +b^{g}\\
    \tilde p_q=p_q+\pi^g
\end{cases}$$
    \item \textit{Checkpoint - Subsection \ref{check}:} we collect all the estimates from the Galbrun Stage in Proposition \ref{recap}. 
\end{itemize}

\noindent \textit{Nash Stage - Section \ref{Nash}.} We use the high-to-low cascade in the quadratic term to cancel the rewritten Reynolds stress.
\begin{itemize}
    \item \textit{Subsections} \ref{secmollnash}, \ref{geomconstrn}: to improve the transport properties we do a second mollification and local-in-time recorrection, with this, we update the charts to be adapted to $(\tilde v_{\ell,j}, \ \tilde B_{\ell,j})$ and finally update the geometric decomposition. This produces a chart update error $R^{crt}$ (flow error in \cite{GR}).
        $$( \tilde{v}_q, \ \tilde{B}_q, \ \tilde{p}_q)
\overset{\text{mollification}}{\longrightarrow}\begin{cases}
    \tilde v_{\ell,j}=\tilde v_\ell + \text{small}\\
    \tilde B_{\ell,j}= \tilde B_\ell +\text{small}\\
    \tilde p_{\ell,j}= \tilde p_\ell +\text{small}\\
\end{cases} \leadsto \ \Psi_I, \ 
        A_I\overset{\text{update}}{\longrightarrow} 
        \tilde \Psi_I, \ \tilde A_I \leadsto R^{crt}=\sum_Ig_I^2(\tilde A_I-A_I)$$
    \item \textit{Nash LDF} \ref{nashldf}: with the updated charts and decomposition, we construct a second Lagrangian displacement field $\xi^p$ as a sum of space-time-disjoint pieces that cancel the rewritten Reynolds stress, the usual step in convex integration. We perturb $(\tilde v_q, \tilde B_q, \tilde p_q)$ to $(v_{q+1}, B_{q+1}, p_{q+1})$ as before.
    $$( \tilde v_{\ell,j}, \ \tilde B_{\ell,j}, \ \tilde p_{\ell,j})+\begin{cases}
    \xi^p\\
    X^p=\exp_{\IId}^{\xi^p}\\
\end{cases}
\overset{\text{Nash stage}}{\longrightarrow} \begin{cases}
    v_{q+1}=\partial_tX^p\circ (X^p)^{-1}+X^p_*\tilde v_q&=v_q+w^{g}+w^p\\
    B_{q+1}=X^p_*\tilde v_q&=B_q+b^{g}+b^p\\
    p_{q+1}=p_q+\pi^g+\pi^p
\end{cases}$$
\end{itemize}

\noindent \textit{Conclusion - Section \ref{conclusion}.} We collect all the estimates from the two Stages and conclude the proof of the Iterative Proposition \ref{iterative}, this is Subsection \ref{proofiterative}. At the very end, in Subsection \ref{proofmain}, we prove the Main Theorem \ref{main} by constructing three distinct initial sub-solutions satisfying the Iterative Assumptions \eqref{inductiveassumptionsgeneral}, thereby demonstrating the variety of data the scheme can handle.

\noindent\textbf{Global Viewpoint and Error Decomposition.} Both stages will produce errors which prevent the perturbed fields from solving the Ideal MHD system exactly. These constitute the new Reynolds stress $R_{q+1}$. We now give a global view of its structure. According to the above decomposition of $(v_{q+1}, \ B_{q+1}, \ p_{q+1})$, the fact that by construction: $$\delta_{q+1}\left(\IId-\frac{R_\ell}{\delta_{q+1}}\right)=\sum_IA_I\Longrightarrow \ddiv \left[R_\ell\right]=-\ddiv \left[\sum_IA_I\right]$$
and that $(v_{q}, \ B_{q}, \ p_{q}, \ R_q)$ solves \eqref{relaxedMHD}, we can write:
\begin{equation}\label{Rq+1}
    \begin{split}
        \ddiv \ R_{q+1}&=\partial_tv_{q+1}+\ddiv(v_{q+1}\otimes v_{q+1}-B_{q+1}\otimes B_{q+1})+\nabla p_{q+1}\\
        &=\partial_tw^p+(v_q+w^{g})\cn w^p-(B_q+b^{g})\cn b^p+w^p\cn(v_q+w^{g})-b^p\cn(B_q+b^{g})\\
        &+\partial_tw^{g}+v_q\cn w^{g}-B_q\cn b^{g}+w^{g}\cn v_q-b^{g}\cn B_q+\nabla \pi^{gr}+\ddiv \left[R_\ell+\sum_I g_I^2 A_I\right]\\
        &+\ddiv [w^p\otimes w^p -b^p\otimes b^p-\sum_I g_I^2 A_I +\pi^p]\\
        &+\ddiv[w^{g}\otimes w^{g}-b^{g}\otimes b^{g}+R_q-R_\ell]\\
        &=\underbrace{\partial_tw^{g}+v_q\cn w^{g}-B_q\cn b^{g}+w^{g}\cn v_q-b^{g}\cn B_q+\nabla \pi^{gr}+\ddiv \left[\sum_I (g_I^2-1) A_I\right]}_{\text{forced linearized MHD + linear errors: Galbrun stage}}\\
        &+\underbrace{\ddiv\left[w^{g}\otimes w^{g}-b^{g}\otimes b^{g}\right]}_{\text{residual quadratic interactions}}+\ddiv[R_q-R_\ell]\\
        &+\underbrace{\partial_tw^p+(v_q+w^{g})\cn w^p-(B_q+b^{g})\cn b^p+w^p\cn(v_q+w^{g})-b^p\cn(B_q+b^{g})}_{\text{linear errors: Nash stage}}\\
        &+\underbrace{\ddiv \left[w^p\otimes w^p-b^p\otimes b^p -\sum_I g_I^2 \tilde{A}_I +\pi^p\IId\right]}_{\text{quadratic interaction + cancellation}}+\underbrace{\ddiv\left[\sum_I g_I^2 (\tilde{A}_I-A_I)\right]}_{\text{chart update error}}\\
        &=\ddiv [R^g+R^p]
    \end{split}
\end{equation}
where we defined $R_{q+1}=R^g+R^p$ and enforced the following identities
\begin{equation}\label{Rgp}
    \begin{split}
        \ddiv R^g&=\partial_tw^{g}+v_q\cn w^{g}-B_q\cn b^{g}+w^{g}\cn v_q-b^{g}\cn B_q+\nabla \pi^{gr}+\ddiv \left[\sum_I (g_I^2-1) A_I\right]\\
        &+\ddiv\left [w^{g}\otimes w^{g}-b^{g}\otimes b^{g}+R_q-R_\ell\right],\\
        \ddiv R^p&=\partial_tw^p+(v_q+w^{g})\cn w^p-(B_q+b^{g})\cn b^p+w^p\cn(v_q+w^{g})-b^p\cn(B_q+b^{g})\\
        &+\ddiv\left[w^p\otimes w^p-b^p\otimes b^p -\sum_I g_I^2 \tilde{A}_I +\pi^p\IId\right]+\ddiv\left[\sum_I g_I^2 (\tilde{A}_I-A_I)\right].
    \end{split}
\end{equation}
The precise definitions of $R^g$ and $R^p$ will be given in Subsections \ref{Rgs} and \ref{Rps} after appropriate rewriting. Note that what is called the flow error in \cite{GR} here is the chart update error.


\subsection{Parameters}\label{choiceofparameters} We now define all the parameters that will appear throughout this work. 

\noindent \textit{Frequency-Amplitude}: for any $q\geq0$ integer and $a,b> 1$ to be determined later, we let 
    $$\lambda_q=a^{b^q} \ \ \text{ and } \ \ \delta_q=\lambda_q^{-2\beta}$$
    note that, since the constructions will be local in space, we don't need to take the integer part of the frequencies.
    
\noindent \textit{Mollification}: $\ell^{-1}=\lambda_q^{1+(b-1)\gamma_\ell} $ where $\gamma_\ell$ quantifies the additional smallness of the mollification scale compared to the inverse of a `good derivative' $1/\lambda_q$.

\noindent \textit{Deep Mollification}: $m_0\geq 1$ corresponds to the number $\bar m_0=m_0-1$ of zero moments we require on the convolution kernel, sometimes referred to as the depth of the smoothing operator, see \cite{gromov1986partial} and \eqref{momentumconditions} for the notation used here. This allows via Proposition \eqref{deepmollification} to choose $\gamma_\ell>0$ as small as needed, that is $\ell\lambda_q$ arbitrarily close to $1$ and trade $m_0$  `good derivative' for the desired mollification-error estimate, see \eqref{constraintm0} below. 

\noindent \textit{Time scales:} we use the following `slow' and `fast' time scales   
    \begin{equation}\label{tauc}
        \text{slow time-scale } \leadsto\tau^c=\frac{1}{\lambda_q\delta_q^{1/2}\lambda_{q+1}^{\alpha}}
    \end{equation}
    this corresponds to the size of the support of the space-time cut-offs and partitions of unity and 
    \begin{equation}\label{taua}
         \text{fast time-scale } \leadsto \tau^a=\lambda_q^{-(b-1)\gamma_a}\tau^c
    \end{equation}
    where $\gamma_a$ quantifies the gain of the fast time frequency over the Alfv\'{e}n transport derivative of the slow coefficients, namely $1/\tau^c$. The letter $a$ stands for Alfv\'{e}n. We described the role of this parameter already in the introductory Subsection \ref{overview}.

\noindent \textit{Mollification along the Alfv\'en Directions}: $\ell_t=\lambda_q^{-(b-1)\gamma_t}\tau^c$ it corresponds to the time-scale used to mollify the slow-coefficients in the Nash stage to avoid a loss of transport derivatives. This procedure is classically called mollification along the flow, see \cite{Isettflow}. Here, we need to adapt the usual procedure to the case of multiple space-time-commuting vector fields.

\noindent \textit{CZ exponent}: $\alpha=\frac{(b-1)}{b}\gamma_{CZ}$, note that this choice gives $\lambda_{q+1}^\alpha=\lambda_q^{(b-1)\gamma_{CZ}}$.
    
\noindent \textit{Number of derivatives}: $M, \ N$ they quantify the number of derivatives propagated in the scheme, see \eqref{inductiveassumptionsgeneral}. 

\noindent \textit{Lie-Taylor Expansion Parameters:} $k_0^g, \ k_0^p$ they fix the order of the Lie-Taylor expansions used in the Galbrun and Nash ($p$ for principal) Stages. The following quantities are associated with the gain from each additional order of expansion in the Lie-Taylor series, see \eqref{gainTp} below and \eqref{Tintro} in the overview Subsection. 
\begin{equation*}
    \begin{split}
        \mathcal{T}_g&=\lambda_q^2\ell^{-\alpha}\tau^a\tau^c\delta_{q+1}=\ell^{-\alpha}\lambda_{q+1}^{-2\alpha}(\tau^a/\tau^c)(\delta_{q+1}/\delta_q)=\ell^{-\alpha}\left(\frac{\lambda_q}{\lambda_{q+1}}\right)^{\gamma_a+2\beta+2\gamma_{CZ}},\\
        \mathcal{T}_p&=\lambda_q\tau^a\delta_{q+1}^{1/2}=\lambda_{q+1}^{-\alpha}(\tau^a/\tau^c)(\delta_{q+1}/\delta_q)^{1/2}=\left(\frac{\lambda_q}{\lambda_{q+1}}\right)^{\gamma_a+\beta+\gamma_{CZ}}
    \end{split}
\end{equation*}
where we use the $g, \ p$ notation again to differentiate between the stages.

We now list all the constraints, provide brief explanations of their origins, and show that they can all be satisfied simultaneously. For convenience, we formulate this as a Lemma. The reader may wish to skip this part on a first reading.
\begin{lemma}[Choice of Parameters]\label{parameterconstraints} Given any $0<\beta <1/5$, there exists a choice of parameters above, such that the following constraints can be met simultaneously. 
\begin{itemize}
    \item \textbf{Gain in the Lie-Taylor expansion.} 
    \begin{equation}\label{gainTp}
        \gamma_a+\beta-\gamma_\ell>0 \Longrightarrow \lambda_q^{(b-1)\gamma_\ell}\mathcal{T}_p<1 \ \text{ and } \ \lambda_q^{(b-1)\gamma_\ell}\mathcal{T}_g<1.
    \end{equation}
    \item \textbf{Admissibility of the loss function.} Let $\bar{\bar L}=\delta_q^{-1}$ we require:
    \begin{equation}\label{constraintadmissibility}
        \left(\lambda_{q}^{(b-1)\gamma_\ell}\mathcal{T}_p\right)^{k_0^g}\bar{\bar L}\leq 1/2 \ \text{ and } \ 1-(\gamma_a+2\beta+\gamma_\ell+2\gamma_{CZ})\geq 0  .
    \end{equation}
    \item  \textbf{Moment of the mollification kernel.} We will need to choose $\gamma_\ell$ sufficiently small, and to retain the necessary smallness on the mollification error, we need to require:
    \begin{equation}\label{constraintm0}
        m_0\gamma_\ell \geq 1 \Longrightarrow (\ell\lambda_q)^{m_0}\leq \frac{\lambda_q}{\lambda_{q+1}}.
    \end{equation}
    \item \textbf{Mollification along the flow.} $\tau^a=\ell_t$ but this is a choice out of convenience.
    \item \textbf{Reynold Stress - Optimisation of $\gamma_a$.} We choose $\gamma_a$ to minimise the size of the new Reynold stress, which corresponds to:
    \begin{equation}\label{constraintR1}
        (\tau^a/\tau^c)\delta_{q+1}+(\tau^c/\tau^a)\frac{\lambda_q\delta_q^{1/2}\delta_{q+1}^{1/2}}{\lambda_{q+1}}=2\sqrt{\frac{\lambda_q\delta_q^{1/2}\delta_{q+1}^{3/2}}{\lambda_{q+1}}}.
    \end{equation}
    
    \item \textbf{Reynolds stress - $\alpha$ smallness.} This will appear as the final constraint to close the scheme, ensuring that we can rehabsorb the implicit constants and still guarantee an additional $\lambda_{q+1}^{-\alpha}$ smallness on $R_{q+1}$. 
    \begin{equation}\label{constraintR2}
        \lambda_{q+1}^{2\alpha}\sqrt{\frac{\lambda_q\delta_q^{1/2}\delta_{q+1}^{3/2}}{\lambda_{q+1}}}\leq \lambda_{q+1}^{-2\alpha}\delta_{q+2}
    \end{equation}
    \item \textbf{Constraint on $M$.} 
    \begin{equation}\label{M}
        M-2k_0^g-2m_0-11\geq 0
    \end{equation}
    \item \textbf{Constraint on $N$.} This will appear in the Nash stage to ensure that the remainders of the Lie-Taylor expansions are small enough to close the estimates.
    \begin{equation}\label{constraintremainder}
        \mathcal{M}_p^{M}\left(\lambda_q^{(b-1)\gamma_\ell}\mathcal{T}_p\right)^{k_0^p}\leq \left(\frac{\lambda_q}{\lambda_{q+1}}\right)^5\delta_q \ \text{ and } \ k_0^p=(N-m_0-5)-M
    \end{equation}
    Given $M$ and all the other parameters (for any $a$), this fixes a value for $N$. This also entails that for $N-m_0-2\leq r \leq N$ we have:
    \begin{equation}\label{highderivativesconstraint}
        \left(\lambda_{q}^{(b-1)\gamma_\ell}\mathcal{T}_p\right)^r\leq \frac{\lambda_q}{\lambda_{q+1}}.
    \end{equation}
    We also require:
    \begin{equation}\label{misc}
        (N-M+2)\gamma_\ell-(1+2\gamma_{CZ}+\gamma_a+2\beta)\geq 0, \ \ \ \ N-k_0^g-m_0-4\geq M.
    \end{equation}
    \item \textbf{Miscellaneous.} Those are constraints that allow us to simplify some inequalities but are not crucial to the scheme itself:
    \begin{equation}\label{misc1}
        \begin{split}
            &\ell^{-\alpha}(\delta_{q+1}/\delta_q)<1 \leadsto 2\beta\geq \gamma_{CZ},\\
            &\tau^a/\tau^c\lambda_{q+1}^\alpha<1\leadsto\gamma_a-\gamma_{CZ}\geq 0,\\
            &\ell^{-\alpha}\leq \lambda_q^{(b-1)\gamma_\ell}\leadsto \gamma_\ell-\gamma_{CZ}\geq 0,\\
            &\tau^a/\tau^c\lambda_q^{(b-1)\gamma_\ell}<1\leadsto\gamma_a-\gamma_\ell\geq 0,\\
            &\lambda_q/\lambda_{q+1}\leq 1/2,\\
            &\beta(2b+1)<1-2\gamma_{CZ},\\
            &m_0\geq2.\\
        \end{split}
    \end{equation} 
\end{itemize}
\end{lemma}

\begin{proof}[Proof of Lemma \ref{parameterconstraints}] Given $0<\beta <1/5$ we can find $b$ such that:
\begin{equation}\label{brange}
    1<b<\frac{1-\beta}{4\beta} \Longrightarrow \beta <\frac{1}{4b+1}.
\end{equation}

\noindent \textit{Choice of $\gamma_a$.} We now choose, $\gamma_a$ so that $\tau^a/\tau^c=\lambda_q^{-(b-1)\gamma_a}$ minimises the expression: 
$$(\tau^a/\tau^c)\delta_{q+1}+(\tau^c/\tau^a)\frac{\lambda_q\delta_q^{1/2}\delta_{q+1}^{1/2}}{\lambda_{q+1}}.$$
Recall that the minimisation problem: 
$f(s)=s A+\frac{1}{s}B$ for a positive parameter $0\leq s\leq 1$ and fixed values $A,B>0$  has a solution given by:
$$s_{min}=\sqrt{B/A} \ \text{ and } \ f(s_{min})=2\sqrt{AB},$$
we conclude that:
\begin{equation}\label{gammaa}
    \lambda_q^{(b-1)\gamma_a}=(\tau^c/\tau^a)=\sqrt{\frac{\lambda_{q+1}\delta_{q+1}^{1/2}}{\lambda_q\delta_q^{1/2}}}=\lambda_{q}^{\frac{1}{2}(b-1)(1-\beta)}\iff \gamma_a=\frac{1}{2}(1-\beta)
\end{equation}
with optimal value given by the geometric mean, namely:
\begin{equation}\label{optimisation}
    (\tau^a/\tau^c)\delta_{q+1}+(\tau^c/\tau^a)\frac{\lambda_q\delta_q^{1/2}\delta_{q+1}^{1/2}}{\lambda_{q+1}}=2\sqrt{\frac{\lambda_q\delta_q^{1/2}\delta_{q+1}^{3/2}}{\lambda_{q+1}}}.
\end{equation}
This also takes care of \eqref{constraintR1}.

\noindent \textit{Choice of $\gamma_\ell, \ \gamma_{CZ}$.} We now find $\gamma_{CZ}$ and $\gamma_\ell$ so that the following two constraints are satisfied. 

\noindent The first one is given by:
$$1-(\gamma_a+2\beta+\gamma_\ell+2\gamma_{CZ})\geq 0 \iff \gamma_\ell+2\gamma_{CZ}\leq 1-(\gamma_a+2\beta)=\frac{1-3\beta}{2},$$
where we used \eqref{gammaa}. Note that positive values of $\gamma_\ell, \ \gamma_{CZ}$ exist as soon as $\beta<1/3$ and, in particular, in our case. 

\noindent The second one reads:
$$\lambda_{q+1}^{2\alpha}\sqrt{\frac{\lambda_q\delta_q^{1/2}\delta_{q+1}^{3/2}}{\lambda_{q+1}}}\leq \lambda_{q+1}^{-2\alpha}\delta_{q+2} \iff \beta\leq\frac{1}{4b+1} -\frac{8\gamma_{CZ}}{4b+1},$$
 note that our choice of $b$ in \eqref{brange} is there to guarantee that this bound holds strictly for $\gamma_{CZ}=0$, and in particular, we have some room to fit an additional parameter. Namely, there exists $\gamma_{CZ}$ positive such that: 
 $$0<\gamma_{CZ}\leq \frac{1-(1+4b)\beta}{8}.$$
 
\noindent These two constraints take care of the second requirement in \eqref{constraintadmissibility} and \eqref{constraintR2}, adding the ones in \eqref{misc1}, we get bounds:
\begin{equation*}
        0<\gamma_{\ell}\leq  \frac{1-3\beta}{4},\ \ \ \
        0<\gamma_{CZ}\leq \min\left \{\frac{1-3\beta}{8}, \ \frac{1-(1+4b)\beta}{8}, \ 2\beta, \ \gamma_\ell\right\}
\end{equation*}
we can then fix $\gamma_\ell$ given $\beta$ and then $\gamma_{CZ}$ once fixed $\gamma_\ell$. Note that this choice also ensures that $\gamma_a+\beta-\gamma_\ell>0$ as required in \eqref{gainTp}. 

\noindent \textit{Choice of $m_0$.} With this set, we can choose $m_0>1$, so that
$m_0\gamma_\ell \geq 1$, this takes care of \eqref{constraintm0} and of the remaining constraint in \eqref{misc1}. 

\noindent \textit{Choice of $k_0^g$ and $M$.} Recall that $\bar{\bar L}= 1/\delta_q$ and the definition of $\mathcal{T}_p$ in \eqref{Tp}. We choose $k_0^g$ to satisfy:
\begin{equation}\label{lowerk0g}
    \left(\lambda_{q}^{(b-1)\gamma_\ell}\mathcal{T}_p\right)^{k_0^g}\bar{\bar L}\leq 1/2\Longleftarrow \lambda_q^{2\beta-k_0^g(b-1)(\beta +\gamma_a-\gamma_\ell)}\leq 1/2\Longleftarrow k_0^g> \frac{2\beta}{(b-1)(\beta+\gamma_a-\gamma_\ell)} \text{ and } \ a\gg 0,
\end{equation}
this takes care of the remaining half of \eqref{constraintadmissibility}, and fixes a lower bound for $M$ via \eqref{M}, namely
$$M-2k_0^g-2m_0-11\geq 0.$$

\noindent \textit{Choice of $N$.} We are left with the constraints in \eqref{constraintremainder} and \eqref{misc}. Since all the parameters but not $N$ (and $a$) are set, we postulate:
\begin{equation}\label{Nansatz}
    N=n_0(M+5)+k_0^g+m_0+5
\end{equation}
and find an integer $n_0\geq 1$ so large that:
$$\mathcal{M}_p^{M}\left(\lambda_q^{(b-1)\gamma_\ell}\mathcal{T}_p\right)^{k_0^p}=\mathcal{M}_p^{M}\left(\lambda_q^{(b-1)\gamma_\ell}\mathcal{T}_p\right)^{N-(M+m_0+5)}\leq \left(\frac{\lambda_q}{\lambda_{q+1}}\right)^5\delta_q$$
and 
$$(N-M+2)\gamma_\ell-(1+2\gamma_{CZ}+\gamma_a+2\beta)\geq 0,$$
note that we have a chance to do so because our choice of $\gamma_a$ and $\gamma_\ell$ above guarantees \eqref{gainTp}. 

\noindent For the first constraint, we have: 
$$\mathcal{M}_p=\lambda_{q+1}/\lambda_q\mathcal{T}_p \ \ \text{ and } \ \ \mathcal{T}_p\leq (\lambda_q/\lambda_{q+1})^{\gamma_a+\beta},$$
see \eqref{Tp}, and we can rewrite it as:
$$\lambda_q^{2\beta+ (b-1)(M+5-(N-m_0-5)(\gamma_a+\beta-\gamma_\ell))}\leq 1 \iff (b-1)(M+5-(N-m_0-5)(\gamma_a+\beta-\gamma_\ell))+2\beta\leq 0,$$
we now plug in \eqref{Nansatz} and end up with a lower bound for $n_0$ given by:
$$n_0\geq \frac{1}{\gamma_a+\beta-\gamma_\ell}-\frac{k_0^g}{(M+5)}+\underbrace{\frac{2\beta}{(M+5)(b-1)(\gamma_a+\beta-\gamma_\ell)}}_{\simeq1}.$$
We remark that, given our lower bound on $M$, the third term is of order one, which is why we chose this Ansatz for $N$. 

\noindent From the second constraint, we get: 
$$n_0\geq \frac{(1+\gamma_a+2\beta+2\gamma_{CZ})}{\gamma_\ell(M+5)}+\frac{M-2-(k_0^g+m_0+5)}{M+5}$$
and upon taking $n_0$ large enough to satisfy both lower bounds, the proof is complete. 
\end{proof}

\begin{remark}[Use of Lemma \ref{parameterconstraints}] We collect all the constraints here, but they may appear in different proofs throughout this work. Moreover, different versions of the same constraint may arise; we list only the tightest, and one can derive the needed version from that. We will refer to this section when needed.
\end{remark}


\section{Galbrun Stage}\label{Galbrun}
\subsection{Space Mollification \texorpdfstring{$(v_q, \ B_q, \ R_q)\leadsto (v_{\ell,j},\ B_{\ell,j},\ R_{\ell})$}{(vq, Bq, Rq) to (v l,j, B l,j, R l)}}\label{moll1}
To control an arbitrary large number of derivatives, we need to mollify $(v_q, \ B_q, \ R_q)$. We will use the definitions and notation as in Proposition \ref{deepmollification}. For the reasons explained in the introductory Subsection \ref{overview}, we first construct a partition of unity of $\mathbb{R}$. We use the parameter $\tau^c$ to control the size of the partition elements. We set $t_j=j\tau^c$ for $j\in\mathbb{Z}$, standard techniques allow us to construct a sequence of cut-offs $\eta_j$ and auxiliary cut-offs $\tilde\eta_j$, such that:
\begin{itemize}
    \item The squares give us a partition of unity of $\mathbb{R}$ that is $\sum_{j\in\mathbb{Z}}\eta_j^2\equiv1$.
    \item For all $j\in \mathbb{Z}$ we have support restrictions:
    \begin{equation}\label{supporteta}
            \supp \ \eta_j\subset (t_j-2/3\tau,t_j+2/3\tau),\ \
            \tilde \eta_j \equiv 1 \text{ on } \ (t_j-2/3\tau,t_j+2/3\tau),\ \
            \supp \ \tilde \eta_j\subset (t_j-\tau,t_j+\tau),
    \end{equation}
    in particular, there are at most two active cut-offs at the same time, namely
    $$\supp \ \eta_j\cap \sup\ \eta_{j\pm 2}=\emptyset \ \text{ and } \ \supp\  \tilde \eta_j\cap \sup\ \eta_{j\pm 2}=\emptyset.$$
    \item We have estimates:
    \begin{equation}\label{timecutoffproperties}
        \sup_t|\partial_t^r\eta_j|,\ \sup_t|\partial_t^r\tilde\eta_j|\lesssim (\tau^c)^{-r}
    \end{equation}
    where the implicit constant depends on $r$.
\end{itemize}

\noindent We think of the elements as belonging to two families $\{e,o\}$ for even and odd intervals. We are ready to begin with the mollification. We adopt the notation from Proposition \ref{deepmollification} and set:
$$v_\ell=(v_q)_\ell, \ \ \ p_\ell=(p_q)_\ell, \ \ \ B_\ell=(B_q)_\ell, \ \ \ R_{\ell}=(R_q)_{\ell}$$
together with Els{\"a}sser fields and associated Alfv\'en transport operators:
$$z_q^\pm=v_q\pm B_q, \ \ \ z_{\ell}^\pm=v_\ell\pm B_\ell, \ \ \ \mathcal{A}^\pm=\partial_t+z^\pm_q\cn, \ \ \ \mathcal{A}^\pm_{\ell}=\partial_t+z^\pm_{\ell}\cn.$$
Commuting mollification with linear differential operators, we see that $B_\ell$ is not a solution of \eqref{FH}, indeed
\begin{equation*}
    \begin{cases}
        \partial_tB_\ell+\curl \left[B_\ell\times v_\ell\right]=\curl [M_\ell],\\
        \ddiv \ B_\ell=0\\
    \end{cases}
    \ \text{ where }\ M_\ell=B_\ell\times v_\ell-(B_q\times v_q)_\ell.
\end{equation*}
Similarly, we see that $v_\ell$ solves:
\begin{equation*}
    \begin{cases}
        \partial_t v_\ell+\ddiv\left[v_\ell\otimes v_\ell-B_{\ell}\otimes B_{\ell}\right]+\nabla p_\ell= \ddiv [R_\ell+R_{\ell}^c],\\
        \ddiv   \ v_\ell=0 
    \end{cases}
    \ \text{ where } \ R_\ell^c =\left(v_\ell\otimes v_\ell -B_{\ell}\otimes B_{\ell}\right)-\left(v_q\otimes v_q-B_q\otimes B_q\right)_\ell.
\end{equation*}
We gather the mollification and mollification-error estimates in the following lemma.
\begin{lemma}[Mollification Estimates]\label{standardmollgalbrun}Let $r\geq 0$ be an integer. We have:
    \begin{subequations}
        \begin{align}
            &||\partial_t^jv_\ell||_r,||\partial_t^jB_\ell||_r\lesssim \lambda_q^{r+j}\lambda_q^{[r+j-N]^+(b-1)\gamma_\ell}\delta_q^{1/2} \ \text{ for } j=0,1,2, \ \text{ and }\ (r,j)\neq (0,0),\\
            &||\mathcal{A}^\pm v_\ell||_r,\ ||\mathcal{A}^\pm B_\ell||_r\lesssim \lambda_q^{r+1}\lambda_q^{[r-(M-1)]^+(b-1)\gamma_\ell}\delta_q \ \text{ for } \ 0\leq r\leq N-1,\\
            &||\mathcal{A}^\pm_\ell v_\ell||_r,\ ||\mathcal{A}^\pm_\ell B_\ell||_r\lesssim \lambda_q^{r+1}\lambda_q^{[r-(M-1)]^+(b-1)\gamma_\ell}\delta_q, \\
            &||p_\ell||_r\lesssim \lambda_q^{r}\lambda_q^{[r-M]^+(b-1)\gamma_\ell}\delta_{q},\\
            &||\mathcal{A}^\pm p_\ell||_r\lesssim \lambda_q^{r+1}\lambda_q^{[r-(M-1)]^+(b-1)\gamma_\ell}\delta_{q}^{3/2} \ \text{ for } \ 1\leq r\leq N-1, \\
            &||R_\ell||\lesssim \lambda_q^{r-\alpha}\lambda_q^{[r-M]^+(b-1)\gamma_\ell}\delta_{q+1},\\
            &||\mathcal{A}^\pm R_\ell||_r\lesssim \lambda_q^{r+1-\alpha}\lambda_q^{[r-(M-1)]^+(b-1)\gamma_\ell}\delta_q^{1/2}\delta_{q+1} \ \text{ for } \ 0\leq r\leq N-1.
        \end{align}
    \end{subequations}
    Together with the following error estimates:
    \begin{subequations}
        \begin{align}
            &||\partial_t^j(v_q-v_\ell)||_r, ||\partial_t^j(B_q-B_\ell)||_r\lesssim \lambda_q^{r+j}(\ell\lambda_q)^{m_0}\delta_q^{1/2} \ \text{ for } j=0,1,2, \ \text{ and }\ 0\leq r\leq N-j-m_0, \\
            &||\mathcal{A}^\pm(v_q-v_\ell)||_r, ||\mathcal{A}^\pm(B_q-B_\ell)||_r\lesssim \lambda_q^{r+1}(\ell\lambda_q)^{m_0}\delta_q \ \text{ for } \ 0\leq r\leq M-1-m_0, \\
            &||p_q-p_\ell||_r\lesssim \lambda_q^{r}(\ell\lambda_q)^{m_0}\delta_{q} \ \text{ for } \ 1\leq r \leq M-m_0, \\
            &||\mathcal{A}^\pm (p_q-p_\ell)||_r\lesssim \lambda_q^{r+1}(\ell\lambda_q)^{m_0}\delta_{q}^{3/2} \ \text{ for } \ 1\leq r \leq M-1-m_0, \\
            &||R_q-R_\ell||_r\lesssim \lambda_q^{r}(\ell\lambda_q)^{m_0}\delta_{q+1}\ \text{ for } \ 0\leq r\leq M-m_0,\\
            &||\mathcal{A}^\pm(R_q-R_\ell)||_r\lesssim \lambda_q^{r+1}(\ell\lambda_q)^{m_0}\delta_q^{1/2}\delta_{q+1}\ \text{ for } \ 0\leq r\leq M-m_0-1
        \end{align}
    \end{subequations}
    and bounds for the forcings:
    \begin{subequations}
        \begin{align}
            &||R_\ell^c||_r,\ ||M_\ell||_r\lesssim \lambda_q^r\lambda_q^{[r-(N-m_0)]^+(b-1)\gamma_\ell}(\ell\lambda_q)^{m_0}\delta_q,\\
            &||\mathcal{A}^\pm_\ell R_\ell^c||_r, \ ||\mathcal{A}^\pm_\ell M_\ell||_r \lesssim \lambda_q^{r+1}\lambda_q^{[r-(M-1-m_0)]^+(b-1)\gamma_\ell}(\ell\lambda_q)^{m_0}\delta_q^{3/2}.
        \end{align}
    \end{subequations}
    The implicit constants depend on $r, \ C_0$ and the specific choice of the mollification kernel.
\end{lemma}
The lemma is a direct consequence of Proposition \ref{deepmollification}, we omit the details of the proof. 

\noindent We are now ready to correct $(v_\ell, \ B_\ell, \ p_\ell)$ so that, at least locally in time, we still have a solution to \eqref{relaxedMHD}. Let $j$ be any index of the time partition $\{\eta_j\}_j$ we constructed. We define the corrected vector fields and pressure $(v_{\ell,j}, \ B_{\ell,j}, \ p_{\ell,j})$ as solutions on $\mathbb{T}^3\times (t_j-\tau^c, \ t_j+\tau^c)$ of the full relaxed MHD system:
\begin{equation}\label{localcorreq}
    \begin{cases}
        \partial_t v_{\ell,j}+\ddiv\left[v_{\ell,j}\otimes v_{\ell,j}-B_{\ell,j}\otimes B_{\ell,j}\right]+\nabla p_{\ell,j}=\ddiv\ R_\ell,\\
        \partial_tB_{\ell,j}+\curl \left[B_{\ell,j}\times v_{\ell,j}\right]=0,\\
        \ddiv \ B_{\ell,j}=\ddiv \ v_{\ell,j}=0,\\
        (v_{\ell,j},B_{\ell,j})|_{t=t_j}=(v_\ell,B_\ell)|_{t=t_j}.
    \end{cases}
\end{equation}
Existence and uniqueness follow from the restriction of the time lifespan to $|t-t_j|<\tau^c$. Indeed, Lemma \ref{classicalgalbrun} below gives in \eqref{T*} a minimal time $T^*$ for well-posedness, the definition of $\tau^c$ in \eqref{tauc} and the bounds in Lemma \ref{standardmollgalbrun} then guarantee that:
$$\frac{\tau^c}{T^*}\lesssim \left(\frac{\ell^{-1}}{\lambda_{q+1}}\right)^\alpha$$
and in particular, given $\alpha, \ b$, we can choose $a$ large enough so that $\tau^c /T^*<1$. 

\begin{lemma}[Classical Solutions with Forcing]\label{classicalgalbrun} For any $\alpha\in (0,1)$, there exist a constant $c$ such that for 
\begin{equation}\label{T*}
    T^*\leq c\frac{1}{||v_\ell|_{t_j}||_{1+\alpha}+||B_\ell|_{t_j}||_{1+\alpha}+\int_{t_j-T^*}^{t_j+T^*}||\ddiv R_\ell(\cdot,s)||_{1+\alpha}\dd s}
\end{equation}
the problem \eqref{localcorreq} admits a unique solution $(v_{\ell,j}, \ B_{\ell,j})$ on $\mathbb{T}^3\times (t_j-T^*, \ t_j+T^*)$ and
\begin{equation}\label{aprioriclassical}
        \begin{split}
            &||v_{\ell,j}||_{r+\alpha}+ ||B_{\ell,j}||_{r+\alpha}\lesssim ||v_\ell|_{t_j}||_{r+\alpha}+||B_\ell|_{t_j}||_{r+\alpha} +\int_{t_j-T^*}^{t_j+T^*}||\ddiv R_\ell(\cdot,s)||_{r+\alpha}\dd s\ \text{ for }\ r\geq 1\\
        \end{split}
    \end{equation}
where the implicit constant depends on $r, \ \alpha$. In particular,
    \begin{equation}\label{classicalgalbrunbounds}
        \begin{split}
            &||v_{\ell,j}||_{r+\alpha}, \ ||B_{\ell,j}||_{r+\alpha}\lesssim \lambda_q^r\lambda_q^{[r-(M-1)]^+(b-1)\gamma_\ell}\ell^{-\alpha}\delta_q^{1/2} \ \text{ for }\ r\geq 1.\\
        \end{split}
    \end{equation}
\end{lemma}
The result and its proof are simple adaptations of the corresponding results for the incompressible Euler equations. This will be apparent after passing to the Els{\"a}sser fields formulation of ideal MHD, namely: 
\begin{equation}\label{alfvenoperators}
    z_{\ell,j}^\pm=v_{\ell,j}\pm B_{\ell,j}, \ \ \ \mathcal{A}^\pm_{\ell,j}=\partial_t+z^\pm_{\ell,j}\cn
\end{equation}
in this variables \eqref{localcorreq} becomes:
\begin{equation}\label{localcorreqz}
    \begin{cases}
        \partial_tz_{\ell,j}^\pm+z_{\ell,j}^\mp\cn z_{\ell,j}^\pm+\nabla p_{\ell,j}=\ddiv \ R_\ell,\\
        \ddiv \ z_{\ell,j}^\pm=0,\\
        z_{\ell,j}^\pm|_{t=t_j}=z_{\ell}^\pm|_{t=t_j}
    \end{cases}
\end{equation}
which indeed has a structure similar to that of the Incompressible Euler equations.

\begin{proof}[Proof of Lemma \ref{classicalgalbrun}] Existence and uniqueness follow from the classical contraction mapping argument used for the incompressible Euler equations, see for example \cite[Chapter 4]{bertozzimajda}. 

\noindent We now move to the a priori estimates. In what follows, for notational convenience, we assume $t_j\leq t<t_j+T^*$ for $T^*$ to be determined; the bounds for $t_j-T^*<t\leq t_j$ follow similarly. We first take $\ddiv$ of \eqref{localcorreqz} and deduce the identity: 
$$\nabla p_{\ell,j}=\nabla\Delta^{-1}\ddiv\left[\ddiv \ R_\ell- z_{\ell,j}^\mp\cn z_{\ell,j}^\pm\right],$$
with this at hand, we rewrite \eqref{localcorreqz} as: 
$$(\partial_t+z_{\ell,j}^\mp\cn) z_{\ell,j}^\pm=\mathbb{P}\ddiv R_\ell+\nabla\Delta^{-1}\ddiv[z_{\ell,j}^\mp\cn z_{\ell,j}^\pm].$$
Proposition \ref{standardtransportestimate} gives bounds for forced transport equations under the assumption $|t-t_j||z^\pm_{\ell,j}||_1\leq 1$ while Proposition \ref{czstuff} the boundedness of $CZ$ operators. We deduce that: 
\begin{equation*}
    \begin{split}
        ||z^\pm_{\ell,j}(\cdot,t)||_\alpha&\lesssim ||z^\pm_{\ell}(\cdot,t_j)||_\alpha+\int_{t_j}^t||(\mathbb{P}\ddiv \ R_\ell)(\cdot,s)||_\alpha+||(\nabla\Delta^{-1}\ddiv[z_{\ell,j}^\mp\cn z_{\ell,j}^\pm])(\cdot,s)||_\alpha\dd s\\
        &\lesssim ||z^\pm_{\ell}(\cdot,t_j)||_\alpha +\int_{t_j}^t||(\ddiv \ R_\ell)(\cdot,s)||_\alpha+||z_{\ell,j}^\mp(\cdot,s)||_\alpha[ z_{\ell,j}^\pm(\cdot,s)]_{1+\alpha}\dd s.\\
    \end{split}
\end{equation*}
Summing the two estimates, we obtain:
$$\sum_{\sigma\in \{\pm\}}||z^\pm_{\ell,j}(\cdot,t)||_\alpha\lesssim \sum_{\sigma\in \{\pm\}}||z^\sigma_{\ell}(\cdot,t_j)||_\alpha +\int_{t_j}^t||(\ddiv \ R_\ell)(\cdot,s)||_\alpha+\sum_{\sigma\in \{\pm\}}||z_{\ell,j}^\sigma(\cdot,s)||_\alpha\sum_{\sigma\in \{\pm\}}[ z_{\ell,j}^\sigma(\cdot,s)]_{1+\alpha}\dd s.$$
Now, let $i=1,2,3$ be any coordinate index, differentiating \eqref{localcorreqz} along $\partial_i$ and using that CZ operators commute with pure derivatives, we get:
\begin{equation}\label{1+alpha}
    \begin{split}
        (\partial_t+z_{\ell,j}^\mp\cn)(\partial_i z_{\ell,j}^\pm)&=\mathbb{P}\ddiv \partial_i R_\ell-\mathbb{P}\partial_i z_{\ell,j}^\mp\cn z_{\ell,j}^\pm+\nabla \Delta^{-1}\ddiv[z_{\ell,j}^\mp\cn \partial_i z_{\ell,j}^\pm]\\
        &=\mathbb{P}\ddiv \partial_i R_\ell-\mathbb{P}\partial_i z_{\ell,j}^\mp\cn z_{\ell,j}^\pm+\nabla \Delta^{-1}\ddiv[\partial_i z_{\ell,j}^\pm\cn  z_{\ell,j}^\mp],\\
    \end{split}
\end{equation}
arguing as above, we deduce:
\begin{equation*}
    \begin{split}
        \sum_{\sigma\in \{\pm\}}[z^\pm_{\ell,j}(\cdot,t)]_{1+\alpha}&\lesssim \sum_{\sigma\in \{\pm\}}[z^\sigma_{\ell}(\cdot,t_j)]_{1+\alpha}+\int_{t_j}^t[(\ddiv R_\ell)(\cdot,s)]_{1+\alpha}+[ z_{\ell,j}^-(\cdot,s)]_{1+\alpha}[ z_{\ell,j}^+(\cdot,s)]_{1+\alpha}\dd s\\
        &\lesssim \sum_{\sigma\in \{\pm\}}[z^\sigma_{\ell}(\cdot,t_j)]_{1+\alpha}+\int_{t_j}^t[(\ddiv R_\ell)(\cdot,s)]_{1+\alpha}+\sum_{\sigma\in \{\pm\}}[ z_{\ell,j}^\sigma(\cdot,s)]_{1+\alpha}\sum_{\sigma\in \{\pm\}}[ z_{\ell,j}^\sigma(\cdot,s)]_{1+\alpha}\dd s\\
    \end{split}
\end{equation*}
Summing the $\alpha$ and $1+\alpha$ semi-norm estimates and taking $\sup_{t_j\leq t'\leq t}$, we obtain:
\begin{equation*}
    \begin{split}
        \sum_{\sigma\in \{\pm\}}\sup_{t_j\leq t'\leq t}||z^\sigma_{\ell,j}(\cdot,t')||_{1+\alpha}&\lesssim \sum_{\sigma\in \{\pm\}}||z^\sigma_{\ell}(\cdot,t_j)||_{1+\alpha}+\int_{t_j}^t||(\mathbb{P}\ddiv \ R_\ell)(\cdot,s)||_{1+\alpha}\dd s\\
        &+(t-t_j)\sup_{t_j\leq t'\leq t}\sum_{\sigma\in \{\pm\}}[ z_{\ell,j}^\sigma(\cdot,t')]_{1+\alpha}\sum_{\sigma\in \{\pm\}}||z^\sigma_{\ell,j}(\cdot,t')||_{1+\alpha}\\
        &\leq \sum_{\sigma\in \{\pm\}}||z^\sigma_{\ell}(\cdot,t_j)||_{1+\alpha}+\int_{t_j}^t||(\mathbb{P}\ddiv \ R_\ell)(\cdot,s)||_{1+\alpha}\dd s\\
        &+(t-t_j)\left[\sum_{\sigma\in \{\pm\}}\sup_{t_j\leq t'\leq t}||z^\sigma_{\ell,j}(\cdot,t')||_{1+\alpha}\right]^2.\\
    \end{split}
\end{equation*}
Now, restricting the time interval to ensure $$(t-t_j)C\sum_{\sigma\in \{\pm\}}\sup_{t_j\leq t'\leq t}||z^\sigma_{\ell,j}(\cdot,t')||_{1+\alpha}<1/2$$ where $C$ is the implicit constant in the bound above and depends on $\alpha$, we obtain:
$$\sum_{\sigma\in \{\pm\}}\sup_{t_j\leq t'\leq t}||z^\sigma_{\ell,j}(\cdot,t')||_{1+\alpha}\leq 2C\left[ \sum_{\sigma\in \{\pm\}}||z^\sigma_{\ell}(\cdot,t_j)||_{1+\alpha}+\int_{t_j}^t||(\mathbb{P}\ddiv \ R_\ell)(\cdot,s)||_{1+\alpha}\dd s\right]$$
and the estimate holds, in particular, under the constraint $|t-t_j|< T^*$ with: 
$$T^*\leq \frac{1}{(2C)^2\left[ \sum_{\sigma\in \{\pm\}}||z^\sigma_{\ell}(\cdot,t_j)||_{1+\alpha}+\int_{t_j}^{t_j+T^*}||(\mathbb{P}\ddiv \ R_\ell)(\cdot,s)||_{1+\alpha}\dd s\right]}.$$

\noindent Higher-order derivative estimates follow by arguing as above, see in particular \eqref{1+alpha}, but now applying $\partial_\theta$ to \eqref{localcorreqz} for an arbitrary multi-index $\theta$ with $|\theta|>1$, and using the bounds for forced transport equations in Proposition \ref{standardtransportestimate} again. 

\noindent The bound \eqref{classicalgalbrunbounds} then comes from our choice of $\tau^c$ in \eqref{tauc}, the estimates in Lemma \ref{standardmollgalbrun} and the a priori estimate \eqref{aprioriclassical} we just obtained. 
\end{proof}

We now prove stability estimates for this local correction procedure. This is also a standard result, apart from the Alfv\'en transport and the second-order time derivative bounds.
\begin{lemma}[Stability Estimates - Galbrun Stage]\label{stabilitygalbrun} Let $j\in \mathbb{Z}$, $\sigma=0,1,2$ and $r\geq 0$ an integer. The following bounds hold:
\begin{subequations}
        \begin{align}
        &||\partial_t^\sigma (v_\ell-v_{\ell,j})||_{r+\alpha}, \ ||\partial_t^{\sigma}(B_\ell-B_{\ell,j})||_{r+\alpha}\lesssim \lambda_q^{r+\sigma+1}\lambda_q^{[r+\sigma -(M-m_0-1)]^+(b-1)\gamma_\ell}\ell^{-\alpha}(\ell\lambda_q)^{m_0}\tau^c\delta_{q} ,\\
        &||\mathcal{A}^\pm _{\ell,j}(v_\ell-v_{\ell,j})||_{r+\alpha}, \ ||\mathcal{A}^\pm _{\ell,j}(B_\ell-B_{\ell,j})||_{r+\alpha}\lesssim \lambda_q^{r+1}\lambda_q^{[r-(M-m_0-2)]^+(b-1)\gamma_\ell}\ell^{-\alpha}(\ell\lambda_q)^{m_0}\delta_q ,\\
        &||p_\ell-p_{\ell,j}||_{r+\alpha}\lesssim \lambda_q^{r}\lambda_q^{[r-(M-m_0-1)]^+(b-1)\gamma_\ell}\ell^{-\alpha}(\ell\lambda_q)^{m_0}\delta_{q} \ \text{ for } \ r\geq 1,\\
        &||\mathcal{A}^\pm _{\ell,j}\left(p_\ell-p_{\ell,j}\right)||_{r+\alpha}\lesssim \lambda_q^{r+1}\lambda_q^{[r-(M-m_0-2)]^+(b-1)\gamma_\ell}\ell^{-2\alpha}(\ell\lambda_q)^{m_0}\delta_{q}^{3/2} \ \text{ for } \ r\geq 1.
    \end{align}
\end{subequations}
The implicit constants depend on $r,  \ \alpha,  \ C_0$.
\end{lemma}
\begin{remark}
    In the applications, we will often use these estimates after substituting the bound $\lambda_q\ell^{-\alpha}\delta_q^{1/2}\leq 1/\tau^c$.
    We also note that the additional $\alpha$ loss in the Alfv\'en transport estimate for the pressure arises from commuting the transport operator with a non-local CZ-type operator.
\end{remark}
\begin{proof}[Proof of Lemma \ref{stabilitygalbrun}] For notational convenience, we always assume $t_j\leq t<t_j+\tau^c$ the estimates for $t_j-\tau^c<t\leq t_j$ can be done identically. As in the proof of Lemma \ref{classicalgalbrun} it is convenient to pass to the Els{\"a}sser variables formulation. We let:
\begin{equation*} 
        \Delta^\pm=z_{\ell}^\pm-z_{\ell,j}^\pm, \ \Delta^p= p_\ell-p_{\ell,j}, \
        R_\ell^{c,\pm}=R_\ell^c\pm[M_\ell]_\times,
\end{equation*}
where $[M_\ell]_\times$ is the unique antisymmetric matrix associated with $M_\ell$. Note that $z_\ell^\pm$ solves:
\begin{equation}\label{mollifiedeq}
    \begin{cases}
        \partial_tz_{\ell}^\pm+z_{\ell}^\mp\cn z_{\ell}^\pm+\nabla p_{\ell}=\ddiv[R_\ell+R_\ell^{c,\pm }],\\
        \ddiv \ z_{\ell}^\pm=0.\\
    \end{cases}
\end{equation}
With the notation set, we are ready to begin with the estimates.

\noindent \textbf{Pure space derivatives and standard transport estimates.} From \eqref{localcorreqz} and \eqref{mollifiedeq} we deduce that the differences $\Delta^\pm$ solve:
\begin{equation}\label{Deltaequation}
    \begin{split}
        \ddiv\ [R_\ell^{c,\pm}]&=\partial_t\Delta^\pm+z_{\ell,j}^\mp\cn \Delta^\pm+\Delta^\mp\cn z_{\ell}^\pm +\nabla\Delta^p\\
        &\overset{\Delta \text{ estimates}}{=}\partial_t\Delta^\pm+z_{\ell,j}^\mp\cn \Delta^\pm+\Delta^\mp\cn z_{\ell}^\pm\\
        &+\nabla\Delta^{-1}\ddiv\left[\ddiv[\ R_\ell^{c,\pm}]-(\Delta^\pm\cn z^\mp_{\ell,j}+\Delta^\mp\cn z^\pm_{\ell})\right]\\
        &\overset{\mathcal{A}\Delta \text{ estimates }}{=}\mathcal{A}^\mp_{\ell,j}\Delta^\pm+\Delta^\mp\cn z_{\ell}^\pm\\
        &+\nabla\Delta^{-1}\ddiv\left[\ddiv[R_\ell^{c,\pm}]-(-\Delta^\pm\cn \Delta^\mp +\Delta^\pm\cn z^\mp_{\ell}+\Delta^\mp\cn z^\pm_{\ell})\right],
    \end{split}
\end{equation}
with zero initial conditions. To remove the pressure term, we used the fact that taking $\ddiv$ of the equation we have:
\begin{equation}\label{stabilitypressure}
    \begin{split}
    \Delta\Delta^p&=\ddiv\left[\ddiv[\ R_\ell^{c,\pm}]-(z_{\ell,j}^\mp\cn \Delta^\pm  +\Delta^\mp\cn z_{\ell}^\pm)\right]\\
    &=\ddiv\left[\ddiv[\ R_\ell^{c,\pm}]-(\Delta^\pm\cn z_{\ell,j}^\mp +\Delta^\mp\cn z_{\ell}^\pm)\right].\\
    \end{split}
\end{equation}
Since we work, on time-scales of order $\tau^c$, with $$\tau_c ||z^\pm_{\ell,j}||_\alpha\lesssim \left(\frac{\ell^{-1}}{\lambda_{q+1}}\right)^\alpha<1$$ we can provide estimates by means of standard transport $\&$ Gr{\"o}nwall arguments. We write $\mathbb{P}=\IId-\nabla\Delta^{-1}\ddiv$ for the Leray projector and read \eqref{Deltaequation} as:
$$\mathcal{A}^\mp_{\ell,j} \Delta^\pm=\mathbb{P}\ddiv\ [R_\ell^{c,\pm}]-\mathbb{P}[\Delta^\mp\cn z^\pm_{\ell}]+\nabla\Delta^{-1}\ddiv\left[\Delta^\pm\cn z^\mp_{\ell,j}\right]=:F^\pm.$$
Here we see that $\Delta^\pm$ loses one good derivative on $z^\pm_{\ell,j}, \ z^\pm_\ell$ and $R_\ell^{c,\pm}$ and thus according to Lemmas \ref{classicalgalbrun} and \ref{standardmollgalbrun} we will have only $M-m_0-1$ left. Using Lemmas \ref{standardmollgalbrun}, \ref{classicalgalbrun}, and Proposition \ref{czstuff} we can estimate the right-hand side as:
\begin{equation}\label{Fpm}
    \begin{split}
        ||F^\pm||_{r+\alpha}&\lesssim \lambda_q^{r+1}\lambda_q^{[r-(M-m_0-1)]^+(b-1)\gamma_\ell}\ell^{-\alpha}(\ell\lambda_q)^{m_0}\delta_q\\
        &+(||\Delta^+||_0+||\Delta^-||_0)\lambda_q^{r+1}\lambda_q^{[r-(M-2)]^+(b-1)\gamma_\ell}\ell^{-\alpha}\delta_q^{1/2}+(||\Delta^+||_{r+\alpha}+||\Delta^-||_{r+\alpha})\lambda_q\ell^{-\alpha}\delta_q^{1/2}\\
        &\leq \lambda_q^{r+1}\lambda_q^{[r-(M-m_0-1)]^+(b-1)\gamma_\ell}\ell^{-\alpha}(\ell\lambda_q)^{m_0}\delta_q\\
        &+\sum_{\sigma\in\{+,-\}}\left[||\Delta^\sigma||_\alpha\lambda_q^{r+1}\lambda_q^{[r-(M-2)]^+(b-1)\gamma_\ell}\ell^{-\alpha}\delta_q^{1/2}+||\Delta^\sigma||_{r+\alpha}\lambda_q\ell^{-\alpha}\delta_q^{1/2}\right].
    \end{split}
\end{equation}
Now using the estimate on the forced transport equation contained in Proposition \ref{standardtransportestimate}, with 
$$\Delta^\pm|_{t=t_j}=(z^\pm_\ell-z^\pm_{\ell,j})|_{t=t_j}=0,$$
we deduce:
\begin{equation*}
\begin{split}
    \sum_{\sigma\in\{+,-\}}||\Delta^\sigma(\cdot,t)||_{\alpha}&\lesssim \sum_{\sigma\in\{+,-\}}\int_{t_j}^t[F^\sigma(\cdot,s)]_{\alpha}\dd s\\
    &\lesssim \tau^c\lambda_q\ell^{-\alpha}(\ell\lambda_q)^{m_0}\delta_q+\lambda_q\ell^{-\alpha}\delta_q^{1/2}\int_{t_j}^t\sum_{\sigma\in\{+,-\}}||\Delta^\sigma(\cdot,s)||_\alpha\dd s,
    \\
    \sum_{\sigma\in\{+,-\}}||\Delta^\sigma||_{r+\alpha}&\lesssim \sum_{\sigma\in\{+,-\}}\int_{t_j}^t[F^\sigma(\cdot,s)]_{r+\alpha}+(s-t_j)[z^{-\sigma}_{\ell,j}]_{r+\alpha}[F^\sigma(\cdot,s)]_1\dd s\\
    &\lesssim \tau^c\lambda_q^{r+1}\lambda_q^{[r+1-(M-m_0-1)]^+(b-1)\gamma_\ell}\ell^{-\alpha}(\ell\lambda_q)^{m_0}\delta_q\\
    &+\lambda_q^{r+1}\lambda_q^{[r+1-(M-2)]^+(b-1)\gamma_\ell}\delta_q^{1/2}\ell^{-\alpha}\int_{t_j}^t\sum_{\sigma\in\{+,-\}}||\Delta^\sigma||_{\alpha}\dd s+\lambda_q\delta_q^{1/2}\ell^{-\alpha}\int_{t_j}^t\sum_{\sigma\in\{+,-\}}||\Delta^\sigma||_{r+\alpha}\dd s\\
    &+\lambda_q^{r}\lambda_q^{[r-(M-1)]^+(b-1)\gamma_\ell}\ell^{-\alpha}\delta_q^{1/2}\tau^c\int_{t_j}^t\sum_{\sigma\in\{+,-\}}\lambda_q^2\ell^{-\alpha}\delta_q^{1/2}||\Delta^\sigma(\cdot,s)||_{\alpha}+\lambda_q\ell^{-\alpha}\delta_q^{1/2}||\Delta^\sigma(\cdot,s)||_{1+\alpha})\dd s.\\
\end{split}
\end{equation*}
Now we use Gr\"onwall Lemma on the first estimate, plug the estimate in the second one for $r=1$, plug the result in the same estimate but for $r>1$ and conclude:
\begin{equation}\label{Deltastep1}
    \begin{split}
        ||\Delta^\pm||_{r+\alpha}&\lesssim \tau^c\lambda_q^{r+1}\lambda_q^{[r-(M-m_0-1)]^+(b-1)\gamma_\ell}\ell^{-\alpha}(\ell\lambda_q)^{m_0}\delta_q\\
        &\leq \lambda_q^{r}\lambda_q^{[r-(M-m_0-1)]^+(b-1)\gamma_\ell}(\ell\lambda_q)^{m_0}\delta_q^{1/2},\\
        ||\mathcal{A}^{\mp}_{\ell,j}\Delta^\pm||_{r+\alpha}&\lesssim \lambda_q^{r+1}\lambda_q^{[r-(M-m_0-1)]^+(b-1)\gamma_\ell}\ell^{-\alpha}(\ell\lambda_q)^{m_0}\delta_q
    \end{split}
\end{equation}
where the second estimate can be deduced by plugging the first one into \eqref{Fpm}. The implicit constant depends on $r, \ \alpha, \  C_0$. 

\noindent \textbf{Special transport estimates.} Using the commutation $\mathcal{A}^{+}_{\ell,j}\mathcal{A}^{-}_{\ell,j}=\mathcal{A}^{-}_{\ell,j}\mathcal{A}^{+}_{\ell,j}$
and the estimates on $\mathcal{A}^{\pm}_\ell z^\pm_\ell $ from Lemma \ref{standardmollgalbrun} we will deduce estimates on $\mathcal{A}^{\pm}_{\ell,j}\Delta^\pm$.\\
Applying $\mathcal{A}^{\pm}_{\ell,j}$ to \eqref{Deltaequation}, we get:
\begin{equation*}
    \begin{split}
        \mathcal{A}^{\mp}_{\ell,j}\mathcal{A}^{\pm}_{\ell,j}\Delta^\pm&=\mathcal{A}^{\pm}_{\ell,j}\mathbb{P}\ddiv R_\ell^{c,\pm}-\mathcal{A}^{\pm}_{\ell,j}\mathbb{P}[\Delta^\mp\cn z^\pm_{\ell}]+\mathcal{A}^{\pm}_{\ell,j}\nabla\Delta^{-1}\ddiv\left[ -\Delta^\pm\cn \Delta^\mp+\Delta^\pm\cn z^\mp_{\ell}\right]\\
        &=\underbrace{[\mathcal{A}^{\pm}_{\ell,j},\mathbb{P}\ddiv] R_\ell^{c,\pm}-[\mathcal{A}^{\pm}_{\ell,j},\mathbb{P}][\Delta^\mp\cn z^\pm_{\ell}]+[\mathcal{A}^{\pm}_{\ell,j},\nabla\Delta^{-1}\ddiv]\left[ -\Delta^\pm\cn \Delta^\mp+\Delta^\pm\cn z^\mp_{\ell}\right]}_{T_1}\\
        &+\mathbb{P}\underbrace{\left[\ddiv \mathcal{A}^{\pm}_{\ell,j}R_\ell^{c,\pm}-\mathcal{A}^{\pm}_{\ell,j}[\Delta^\mp\cn z^\pm_{\ell}]\right]}_{T_3}+\nabla\Delta^{-1}\ddiv\underbrace{\mathcal{A}^{\pm}_{\ell,j}\left[ -\Delta^\pm\cn \Delta^\mp+\Delta^\pm\cn z^\mp_{\ell}\right]}_{T_2}\\
    \end{split}
\end{equation*}
Here we see that we lose one derivative on $\mathcal{A}^\pm_{\ell,j}\Delta^\mp, \ \mathcal{A}^\pm_\ell z^\pm_\ell, \ \mathcal{A}^\pm_\ell R^{c,\pm}_\ell$ and thus, according to Lemma \ref{standardmollgalbrun} and \eqref{Deltastep1}, we have only $M-m_0-2$ left. The commutator terms can be estimated using Proposition \ref{czstuff}, the bounds in \eqref{Deltastep1} and Lemma \ref{standardmollgalbrun}:
\begin{equation*}
    \begin{split}
        ||T_1||_{r+\alpha} &\lesssim ||[\mathcal{A}^{\pm}_{\ell,j},\ddiv\mathbb{P}] R_\ell^{c,\pm}||_{r+\alpha}+||[\mathcal{A}^{\pm}_{\ell,j},\mathbb{P}][\Delta^\mp\cn z^\pm_{\ell}]||_{r+\alpha}+||[\mathcal{A}^{\pm}_{\ell,j},\nabla\Delta^{-1}\ddiv]\left[ -\Delta^\pm\cn \Delta^\mp+\Delta^\pm\cn z^\mp_{\ell}\right]||_{r+\alpha}\\
        &\lesssim\underbrace{||z^{\pm}_{\ell,j}||_{1+\alpha}||R_\ell^{c,\pm}||_{r+1+\alpha}+||z^{\pm}_{\ell,j}||_{r+1+\alpha}||R_\ell^{c,\pm}||_{1+\alpha}}_{\ell^{-2\alpha}\lambda_q^{r+2}(\ell\lambda_q)^{m_0}\delta_q^{3/2}}+||z^{\pm}_{\ell,j}||_{1+\alpha}||\Delta^\mp\cn z^\pm_{\ell}||_{r+\alpha}+||z^{\pm}_{\ell,j}||_{r+1+\alpha}||\Delta^\mp\cn z^\pm_{\ell}||_{\alpha}\\
        &+\underbrace{||[z^{\pm}_{\ell,j}||_{1+\alpha}(||\Delta^\pm\cn \Delta^\mp||_{r+\alpha}+||\Delta^\pm\cn z^\mp_{\ell}||_{r+\alpha})+||[z^{\pm}_{\ell,j}||_{r+1+\alpha}(||\Delta^\pm\cn \Delta^\mp||_{\alpha}+||\Delta^\pm\cn z^\mp_{\ell}||_{\alpha})}_{\lambda_q\delta_q^{1/2}\ell^{-\alpha}\lambda_q^{r+1}(\ell\lambda_q)^{2m_0}\delta_q+(\lambda_q\delta_q^{1/2}\ell^{-\alpha})^2\lambda_q^{r}(\ell\lambda_q)^{m_0}\delta_q^{1/2}}\\
        &\lesssim\lambda_q^{r+2} \lambda_q^{[r+1-(M-m_0-1)]^+(b-1)\gamma_\ell}\ell^{-2\alpha}(\ell\lambda_q)^{m_0}\delta_q^{3/2}\\
        &\leq \lambda_q^{r+1} \lambda_q^{[r-(M-m_0-1)]^+(b-1)\gamma_\ell}\ell^{-\alpha}(\ell\lambda_q)^{m_0}(1/\tau^c)\delta_q\\
    \end{split}
\end{equation*}
where we used $\tau^c\lambda_q\ell^{-\alpha}\delta_q^{1/2}\leq 1$ and $(\ell\lambda_q)^{m_0}\lambda_q^{(b-1)\gamma_\ell}\leq (\lambda_q/\lambda_{q+1})^{1-\gamma_  \ell}\leq 1$. In what follows, we will do so without mentioning it.

\noindent We now rewrite:
\begin{equation*}
    \begin{split}
        T_2&=-\underbrace{(\mathcal{A}^{\pm}_{\ell,j}\Delta^\pm)\cn \Delta^\mp}_{\text{Gr{\"o}nwall}}-\underbrace{\Delta^\pm\cn (\mathcal{A}^{\pm}_{\ell,j}\Delta^\mp)}_{\ell^{-\alpha}\lambda_q^{r+2}(\ell\lambda_q)^{2m_0}\delta_q^{3/2}}+\underbrace{(\Delta^\pm\cn z^\pm_{\ell,j})\cn \Delta^\mp}_{\ell^{-2\alpha}\lambda_q^{r+2}(\ell\lambda_q)^{2m_0}\delta_q^{3/2}}\\
        &+\underbrace{(\mathcal{A}^{\pm}_{\ell,j}\Delta^\pm)\cn z^\mp_{\ell}}_{\text{Gr{\"o}nwall}}+\underbrace{\Delta^\pm\cn(\mathcal{A}^{\pm}_\ell z^\mp_{\ell}+(z^\pm_{\ell,j}-z_\ell^\pm)\cn z^\mp_\ell)}_{\lambda_q^r\ell^{-\alpha}(\ell\lambda_q)^{m_0}\lambda_q^2\delta_q^{3/2}+\lambda_q^r\ell^{-\alpha}((\ell\lambda_q)^{m_0}\delta_q^{1/2})^2\lambda_q^2\delta_q^{1/2}}-\underbrace{(\Delta^\pm\cn z_{\ell,j}^\pm)\cn z^\mp_{\ell}}_{\ell^{-2\alpha}\lambda_q^r(\ell\lambda_q)^{m_0}\delta_q^{1/2}(\lambda_q\delta_q^{1/2})^2}\\
    \end{split}
\end{equation*}
and 
\begin{equation*}
    \begin{split}
        T_3&=\underbrace{\ddiv \left[\mathcal{A}^{\pm}_{\ell}R_\ell^{c,\pm}+(z^\pm_{\ell,j}-z^\pm_\ell)\cn R_\ell^{c,\pm} \right]}_{\lambda_q^{r+2}\ell^{-\alpha}(\ell\lambda_q)^{m_0}\delta_q^{3/2}+\lambda_q^{r+2}\ell^{-\alpha}(\ell\lambda_q)^{2m_0}\delta_q^{3/2}}\\
        &-\underbrace{(\mathcal{A}^{\pm}_{\ell,j}\Delta^\mp)\cn z^\pm_{\ell}}_{\lambda_q^{r+1}\ell^{-\alpha}(\ell\lambda_q)^{m_0}\delta_q}-\underbrace{\Delta^\mp\cn(\mathcal{A}^{\pm}_\ell z^\pm_{\ell}+(z^\pm_{\ell,j}-z_\ell^\pm)\cn z^\pm_\ell)}_{\lambda_q^r\ell^{-\alpha}(\ell\lambda_q)^{m_0}\lambda_q^2\delta_q^{3/2}+\lambda_q^r\ell^{-\alpha}((\ell\lambda_q)^{m_0}\delta_q^{1/2})^2\lambda_q^2\delta_q^{1/2}}+\underbrace{(\Delta^\mp\cn z_{\ell,j}^\pm)\cn z^\pm_{\ell}}_{\ell^{-2\alpha}\lambda_q^r(\ell\lambda_q)^{m_0}\delta_q^{1/2}(\lambda_q\delta_q^{1/2})^2}
    \end{split}
\end{equation*}
and from the bounds in \eqref{Deltastep1}, Lemma \eqref{standardmollgalbrun} and the definition of $\tau^c$ \eqref{tauc}, we conclude:
\begin{equation*}
    \begin{split}
        ||\nabla\Delta^{-1}\ddiv \ T_2||_{r+\alpha}
        &\lesssim ||\mathcal{A}^{\pm}_{\ell,j}\Delta^\pm||_\alpha\lambda_q^{r+1}\lambda_q^{[r-(M-m_0-1)]^+(b-1)\gamma_\ell}\ell^{-\alpha}\delta_q^{1/2}+||\mathcal{A}^{\pm}_{\ell,j}\Delta^\pm||_{r+\alpha}\lambda_q\ell^{-\alpha}\delta_q^{1/2}\\
        &+\ell^{-\alpha}(\ell\lambda_q)^{m_0}\lambda_q^{r+1}\lambda_q^{[r-(M-m_0-1)]^+(b-1)\gamma_\ell}(1/\tau^c)\delta_q,\\
        ||\mathbb{P} T_3||_{r+\alpha}&\lesssim \lambda_q^{r+2}\lambda_q^{[r-(M-m_0-2)]^+(b-1)\gamma_\ell}\ell^{-\alpha}(\ell\lambda_q)^{m_0}\delta_q^{3/2}
    \end{split}
\end{equation*}
where, in addition to the observations above, we used that $m_0>1$. We are now ready to do a second transport $\&$ Gr\"onwall argument. Gathering the bounds, we deduce:
\begin{equation}\label{puretimecorrection}
\begin{split}
    ||\mathcal{A}^{\mp}_{\ell,j}\mathcal{A}^{\pm}_{\ell,j}\Delta^\pm||_{r+\alpha}&\lesssim ||T_1||_{r+\alpha}+||\nabla\Delta^{-1}\ddiv \ T_2||_{r+\alpha}+||\mathbb{P} T_3||_{r+\alpha}\\
    &\lesssim ||\mathcal{A}^{\pm}_{\ell,j}\Delta^\pm||_\alpha\lambda_q^{r+1}\lambda_q^{[r-(M-m_0-1)]^+(b-1)\gamma_\ell}\ell^{-\alpha}\delta_q^{1/2}+||\mathcal{A}^{\pm}_{\ell,j}\Delta^\pm||_{r+\alpha}\lambda_q\ell^{-\alpha}\delta_q^{1/2}\\
    &+\ell^{-\alpha}(\ell\lambda_q)^{m_0}\lambda_q^{r+1}\lambda_q^{[r-(M-m_0-2)]^+(b-1)\gamma_\ell}(1/\tau^c)\delta_q,\\
\end{split}
\end{equation}
we now apply Proposition \ref{standardtransportestimate}, with
\begin{equation*}
    \begin{split}
        \mathcal{A}^{\pm}_{\ell,j}\Delta^\pm|_{t=t_j}&=\mathcal{A}^{\mp}_{\ell,j}\Delta^\pm|_{t=t_j}+2B_{\ell,j}\cn \Delta^\pm|_{t=t_j}\\
        &=\mathcal{A}^{\mp}_{\ell,j}\Delta^\pm|_{t=t_j}+2B_{\ell,j}(\cdot,t_j)\cn [\underbrace{(z^\pm_\ell(\cdot,t_j)-z^\pm_{\ell,j}(\cdot, t_j)}_{\equiv0}]\\
        &=\mathcal{A}^{\mp}_{\ell,j}\Delta^\pm|_{t=t_j}
    \end{split}
\end{equation*}
and deduce from \eqref{puretimecorrection} the bounds:
\begin{equation*}
    \begin{split}
        ||\mathcal{A}^{\pm}_{\ell,j}\Delta^\pm||_{\alpha}&\lesssim ||\mathcal{A}^{\mp}_{\ell,j}\Delta^\pm||_\alpha+\int_{t_j}^t[\mathcal{A}^{\mp}_{\ell,j}\mathcal{A}^{\pm}_{\ell,j}\Delta^\pm(\cdot,s)]_\alpha\dd s\\
        &\lesssim \lambda_q\ell^{-\alpha}(\ell\lambda_q)^{m_0}\delta_q+\lambda_q\ell^{-\alpha}\delta_q^{1/2}\int_{t_j}^t||\mathcal{A}^{\pm}_{\ell,j}\Delta^\pm||_{\alpha}(\cdot,s)\dd s,
        \\
        ||\mathcal{A}^{\pm}_{\ell,j}\Delta^\pm||_{r+\alpha}&\lesssim ||\mathcal{A}^{\mp}_{\ell,j}\Delta^\pm||_{r+\alpha}+|t-t_j|[z^\mp_{\ell,j}]_{r+\alpha}[\mathcal{A}^{\mp}_{\ell,j}\Delta^\pm]_{1}\\
        &+\int_{t_j}^t[\mathcal{A}^{\mp}_{\ell,j}\mathcal{A}^{\pm}_{\ell,j}\Delta^\pm(\cdot,s)]_{r+\alpha}+(s-t_j)[z^{\mp}_{\ell,j}]_{r+\alpha}[\mathcal{A}^{\mp}_{\ell,j}\mathcal{A}^{\pm}_{\ell,j}\Delta^\pm(\cdot,s)]_1\dd s\\
        &\lesssim \lambda_q^{r+1}\lambda_q^{[r-(M-m_0-2)]^+(b-1)\gamma_\ell}\ell^{-\alpha}(\ell\lambda_q)^{m_0}\delta_q\\
        &+\int_{t_j}^t(||\mathcal{A}^{\pm}_{\ell,j}\Delta^\pm(\cdot,s)||_\alpha\lambda_q^{r+1}\lambda_q^{[r-(M-m_0-1)]^+(b-1)\gamma_\ell}\ell^{-\alpha}\delta_q^{1/2}+||\mathcal{A}^{\pm}_{\ell,j}\Delta^\pm(\cdot,s)||_{r+\alpha}\lambda_q\ell^{-\alpha}\delta_q^{1/2})\dd s\\
        &+\tau^c\lambda_q^r\lambda_q^{[r-(M-1)]^+(b-1)\gamma_\ell}\ell^{-\alpha}\delta_q^{1/2}\int_{t_j}^t(||\mathcal{A}^{\pm}_{\ell,j}\Delta^\pm(\cdot,s)||_\alpha\lambda_q^{2}\ell^{-\alpha}\delta_q^{1/2}+||\mathcal{A}^{\pm}_{\ell,j}\Delta^\pm(\cdot,s)||_{1+\alpha}\lambda_q\ell^{-\alpha}\delta_q^{1/2})\dd s
    \end{split}
\end{equation*}
where we used $M-m_0-2\geq 0$. Now we conclude as before, we use Gr{\"o}nwall Lemma on the first estimate, plug the estimate in the second one for $r=1$, plug the result in the same estimate but for $r>1$ and we conclude that:
\begin{equation}\label{Deltastep2}
    \begin{split}
        ||\mathcal{A}^{\pm}_{\ell,j}\Delta^\pm||_{r+\alpha}&\lesssim \lambda_q^{r+1}\lambda_q^{[r-(M-m_0-2)]^+(b-1)\gamma_\ell}\ell^{-\alpha}(\ell\lambda_q)^{m_0}\delta_q
    \end{split}
\end{equation}
where the implicit constant depends on $r, \ \alpha, \  C_0$. The stated estimates then follow from unwinding the definitions, namely
$$v_\ell-v_{\ell,j}=\frac{\Delta^++\Delta^-}{2} \ \text{ and } \ B_\ell-B_{\ell,j}=\frac{\Delta^+-\Delta^-}{2}.$$

\noindent \textbf{Pressure estimates.} To prove bounds on $\nabla \Delta^p$ we use that $\nabla\Delta^{-1}\ddiv$ is a CZ-type operator together with the result in Proposition \ref{czstuff} and the bounds we just proved. Indeed from \eqref{stabilitypressure} and the bounds in Lemma \ref{standardmollgalbrun} and \eqref{Deltastep1} we have:
\begin{equation*}
    \begin{split}
        ||\nabla \Delta^p||_{r+\alpha}&\lesssim||R_\ell^{c,\pm}||_{r+1+\alpha}+||z_{\ell,j}^\mp||_{r+1+\alpha}||\Delta^\pm||_0+||z_{\ell,j}^\mp||_{1}||\Delta^\pm||_{r+\alpha}+||z_{\ell}^\pm||_{r+1+\alpha}||\Delta^\mp||_0+||z_{\ell}^\pm||_{1}||\Delta^\mp||_{r+\alpha}\\
        &\lesssim\lambda_q^{r+1}\lambda_q^{[r+1-(M-m_0-1)]^+(b-1)\gamma_\ell}\ell^{-\alpha}(\ell\lambda_q)^{m_0}\delta_q
    \end{split}
\end{equation*}
and with the help of the commutator estimates in Proposition \ref{czstuff}, the bounds in Lemma \ref{standardmollgalbrun} and \eqref{Deltastep1}, \eqref{Deltastep2}, we deduce:
\begin{equation*}
    \begin{split}
        ||\mathcal{A}^\pm_{\ell,j} \nabla\Delta^p||_{r+\alpha}&\lesssim ||[z_{\ell,j},\nabla\Delta^{-1}\ddiv][\ddiv[\ R_\ell^{c,\pm}]-(\Delta^\pm\cn z_{\ell,j}^\mp+\Delta^\mp\cn z_{\ell}^\pm)]||_{r+\alpha}\\
        &+||\mathcal{A}^\pm_{\ell,j}[\ddiv[\ R_\ell^{c,\pm}]-(\Delta^\pm\cn z_{\ell,j}^\mp+\Delta^\mp\cn z_{\ell}^\pm)]||_{r+\alpha}\\
        &\lesssim \underbrace{||z_{\ell,j}^\pm||_{r+1+\alpha}(||\Delta^\pm\cn z_{\ell,j}^\mp||_\alpha+||\Delta^\pm\cn z_\ell^\mp||_\alpha+||\ddiv \ R_\ell^{c,\pm}||_\alpha)}_{\lambda_q^{r+2}\ell^{-2\alpha}(\ell\lambda_q)^{m_0}\delta_q^{3/2}}\\
        &+||z_{\ell,j}^\pm||_{1+\alpha}(||\Delta^\pm\cn z_{\ell,j}^\mp||_{r+\alpha}+||\Delta^\pm\cn z_\ell^\mp||_{r+\alpha}+||\ddiv \ R_\ell^{c,\pm}||_{r+\alpha})\\
        &+\underbrace{||(\mathcal{A}^\pm_{\ell,j}\Delta^\pm)\cn z_{\ell,j}^\mp||_{r+\alpha}}_{\lambda_q^{r+2}\ell^{-2\alpha}(\ell\lambda_q)^{m_0}\delta_{q}^{3/2}}+\underbrace{||\Delta^\pm\cn \mathcal{A}^\pm_{\ell,j}z_{\ell,j}^\mp||_{r+\alpha}}_{\lambda_q^{r+2}\ell^{-\alpha}(\ell\lambda_q)^{m_0}\delta_{q}^{3/2}}+\underbrace{||(\Delta^\pm\cn z^\pm_{\ell,j})\cn z_{\ell,j}^\mp||_{r+\alpha}}_{\lambda_q^{r+2}\ell^{-2\alpha}(\ell\lambda_q)^{m_0}\delta_{q}^{3/2}}\\
        &+||(\mathcal{A}^\pm_{\ell,j}\Delta^\mp)\cn z_\ell^\pm||_{r+\alpha}+||\Delta^\mp\cn [(z^\pm_{\ell,j}-z^\pm_\ell)\cn z_\ell^\pm+\mathcal{A}_\ell^\pm z_\ell^\pm]||_{r+\alpha}+||(\Delta^\mp\cn z^\pm_{\ell,j})\cn z_\ell^\mp||_{r+\alpha}\\
        &+\underbrace{||\ddiv\mathcal{A}^\pm_\ell R_\ell^{c,\pm}||_{r+\alpha}+||[z^\pm_\ell\cn ,\ddiv] R_\ell^{c,\pm}||_{r+\alpha}}_{\lambda_q^{r+2}\ell^{-\alpha}(\ell\lambda_q)^{m_0}\delta_q^{3/2}}+\underbrace{||(z^\pm_{\ell,j}-z^\pm_\ell)\cn\ddiv\ R_\ell^{c,\pm}||_{r+\alpha}}_{\lambda_q^{r+2}\ell^{-\alpha}(\ell\lambda_q)^{2m_0}\delta_q^{3/2}}\\
        &\lesssim \lambda_q^{r+2}\lambda_q^{[r+1-(M-m_0-2)]^+(b-1)\gamma_\ell}\ell^{-2\alpha}(\ell\lambda_q)^{m_0}\delta_q^{3/2}
    \end{split}
\end{equation*}
and we conclude:
$$||\nabla\mathcal{A}^\pm_{\ell,j} \Delta^p||_{r+\alpha}\lesssim ||\mathcal{A}^\pm_{\ell,j}\nabla \Delta^p||_{r+\alpha}+||(\DD z^\pm_{\ell,j})^\top\nabla \Delta^p||_{r+\alpha}\lesssim \lambda_q^{r+2}\lambda_q^{[r+1-(M-m_0-2)]^+(b-1)\gamma_\ell}\ell^{-2\alpha}(\ell\lambda_q)^{m_0}\delta_q^{3/2}$$
where the implicit constant depends on $r, \ \alpha, \  C_0$.

\noindent \textbf{Pure time derivatives.} We write the first time derivative as:
\begin{equation*}
    \begin{split}
        \partial_t\Delta^\pm=\mathcal{A}^{\mp}_{\ell,j}\Delta^\pm-z^{\mp}_{\ell,j}\cn \Delta^\pm
    \end{split}
\end{equation*}
and we deduce from the bounds shown above that:
\begin{equation*}
    \begin{split}
        ||\partial_t\Delta^\pm||_{r+\alpha}&\lesssim \lambda_q^{r+2}\lambda_q^{[r+1-(M-m_0-1)]^+(b-1)\gamma_\ell}\ell^{-\alpha}(\ell\lambda_q)^{m_0}\tau^c\delta_{q}.\\
    \end{split}
\end{equation*}
We read the second-order pure time derivative bound simply by applying $\partial_t$ to equation \eqref{Deltaequation}:
\begin{equation*}
    \begin{split}
         \partial_t^2\Delta^\pm=-\mathbb{P}\left[(\partial_tz^\mp_{\ell,j})\cn \Delta^\pm+z^\mp_{\ell,j}\cn(\partial_t \Delta^\pm)+(\partial_t\Delta^\mp)\cn z^\pm_{\ell}+\Delta^\mp\cn (\partial_t z^\pm_{\ell})\right]+\mathbb{P}\ddiv\ [\partial_t R_\ell^{c,\pm}]
    \end{split}
\end{equation*}
and the first-order bound above, the ones in \eqref{Deltastep1} and Lemma \ref{standardmollgalbrun}, then give: 
\begin{equation*}
    \begin{split}
        ||\partial_t^2\Delta^\pm||_{r+\alpha}
        &\lesssim \lambda_q^{r+3}\lambda_q^{[r+2 -(M-m_0-1)]^+(b-1)\gamma_\ell}\ell^{-\alpha}(\ell\lambda_q)^{m_0}\tau^c\delta_{q}. 
    \end{split}
\end{equation*}
The implicit constants in the above statements depend on $r, \ \alpha, \  C_0$.
\end{proof}


\begin{remark}[A Useful Identity] Let 
$$\Delta^v=v_\ell-v_{\ell,j}, \ \ \ \ \Delta^B=B_\ell-B_{\ell,j}, \ \ \ \ \Delta^p=p_\ell-p_{\ell,j},$$
we now rewrite
\begin{equation*}
    \begin{split}
        \ddiv \ R_\ell^c&=\partial_t\Delta^v+\ddiv\left[\left[v_\ell\otimes v_\ell-B_\ell\otimes B_\ell\right]-\left[v_{\ell,j}\otimes v_{\ell,j}-B_{\ell,j}\otimes B_{\ell,j}\right]\right]+\nabla\Delta^p\\
        &=\partial_t\Delta^v+\ddiv\left[\left[v_\ell\otimes \Delta^v+\Delta^v\otimes v_{\ell,j}\right]-\left[B_\ell\otimes \Delta^B+\Delta^B\otimes B_{\ell,j}\right]\right]+\nabla\Delta^p\\
        &=\partial_t\Delta^v+\ddiv\left[\Delta^v\otimes\Delta^v-\Delta^B\otimes\Delta^B+\left[\Delta^v\otimes v_{\ell,j}\right]^{sym}-\left[\Delta^B\otimes B_{\ell,j}\right]^{sym}\right]\\
        &+\nabla\Delta^p\\
    \end{split}
\end{equation*}
and we conclude that:
\begin{equation}\label{specialdoublegalbrun}
    \begin{split}
        (\partial_t+\mathcal{L}_{v_{\ell,j}})\Delta^v-\mathcal{L}_{B_{\ell,j}}\Delta^B&=\ddiv \ R_\ell^c-2\left(\Delta^v\cn v_{\ell,j}-\Delta^B\cn B_{\ell,j}\right)\\
        &-\ddiv\left[\Delta^v\otimes\Delta^v-\Delta^B\otimes\Delta^B\right]-\nabla\Delta^p\\
    \end{split}
\end{equation}
This will be useful later when computing the transport part of the new Reynolds stress.
\end{remark}

We now collect the work done in this Subsection. 
\begin{lemma}[Local Correction]\label{localcorrg} Let $j\in \mathbb{Z}$, $r\geq 0$ an integer and $\underline{r}=M-m_0-1$. The following bounds hold:
\begin{subequations}
        \begin{align}
        &||v_{\ell,j}||_{0}, \ ||B_{\ell,j}||_{0}\leq C_0 \ \text{ and } |B_{\ell,j}|\geq c_0 \ \text{ everywhere on} \ \mathbb{T}^3\times \mathbb{R}, \\
        &||\partial_t^\sigma v_{\ell,j}||_{r}, \ ||\partial_t^{\sigma}B_{\ell,j}||_{r}\lesssim \lambda_q^{r+\sigma}\lambda_q^{[r+\sigma-\underline{r}]^+(b-1)\gamma_\ell}\delta_{q}^{1/2} \ \text{ for } \ \sigma =0,1,2  \ \& \ (\sigma ,r)\neq (0,0),\\
        &||\mathcal{A}^\pm _{\ell,j}v_{\ell,j}||_{r}, \ ||\mathcal{A}^\pm _{\ell,j}B_{\ell,j}||_{r}\lesssim \lambda_q^{r+1}\lambda_q^{[r-(\underline{r}-1)]^+(b-1)\gamma_\ell}\delta_q,\\
        &||p_{\ell,j}||_{r}\lesssim \lambda_q^{r}\lambda_q^{[r-\underline{r}]^+(b-1)\gamma_\ell}\delta_{q} \ \text{ for }  \ r\geq 1,\\
        &||\mathcal{A}^\pm _{\ell,j}p_{\ell,j}||_{r}\lesssim \lambda_q^{r+1}\lambda_q^{[r-(\underline{r}-1)]^+(b-1)\gamma_\ell}\delta_{q}^{3/2} \ \text{ for }  \ r\geq 1.
    \end{align}
\end{subequations}
The implicit constants depend on $r, \ C_0, \ \alpha$.
\end{lemma}
The bounds follow immediately from writing:
\begin{equation*}
    \begin{split}
        &\partial_t^\sigma v_{\ell,j}=\begin{cases}
            (v_{\ell,j}-v_\ell)+(v_\ell-v_q)+v_q &\text{ for } r=\sigma=0,\\
            \partial_t^\sigma(v_{\ell,j}-v_\ell)+\partial_t^\sigma v_\ell &\text{ else } \\
        \end{cases}\\
        &\mathcal{A}^\pm _{\ell,j}v_{\ell,j}=\mathcal{A}^\pm _{\ell,j}(v_{\ell,j}-v_\ell)+\mathcal{A}^\pm _{\ell}v_\ell+(v_{\ell,j}-v_\ell)\cn v_\ell,
    \end{split}
\end{equation*}
together with the estimates in Lemmas \ref{standardmollgalbrun}, \ref{stabilitygalbrun}, and the fact that $(\ell\lambda_q)^{m_0}\ell^{-2\alpha}\leq 1$, see \eqref{constraintadmissibility}. Similar decompositions hold for $B_{\ell,j}$ and $p_{\ell,j}$.


\subsection{Geometric Constructions}\label{geomconstr}
With the help of Proposition \ref{minimal} and the explicit construction afterwards, see the definitions and notation there, we build a minimal covering $\mathcal{C}=\{U_{j'}\}_{j'}$ of $\mathbb{T}^3$ with disjoint families $\{\mathcal{F}_l\}$ for $l=1,\dots,4$ and an associated partition of unity $\{\theta_{j'}\}_{j'}$ with diameter $4/3\tau^c$ and granularity $1/3\tau^c$. Namely, for each $U_{j'}\in \mathcal{C}$ we have $\theta_{j'}\in C_c^{\infty}(U_{j'},[0,1])$ satisfying:
\begin{itemize}
    \item $\sum_{j'}\theta_{j'}^2\equiv1$ and for each point in $\mathbb{T}^3$ at most four $\theta_{j'}$ are non-zero.
    \item For each $j'$ there exists $x_{j'}$ such that $U_{j'}\subset B_{2/3\tau^c}(x_{j'})$ and  $$U_{j'_1},U_{j'_2}\in \mathcal{F}_l \Longrightarrow \text{dist}(U_{j'_1}, \ U_{j'_2})> 1/3\tau^c.$$
    \item $||\theta_{j'}||_r\lesssim (\tau^c)^{-r}$ where the implicit constant depends on $r$ and the specific choice of the partition of unity.
\end{itemize}
Now set $J=(j,j')$ and define:
\begin{equation}\label{J}
    \chi_J=\chi_{j,j'}=\eta_j(t)\theta_{j'}(x) \ \text{ and } \ Q_J=B_{\tau^c}(x_{j'})\times B_{\tau^c}(t_{j}),
\end{equation}
where $\{\eta_j\}_j$ was constructed in Subsection \ref{moll1}, note that $\{\chi_J^2\}_{J}$ is now a space-time partition of unity, in the sense that:
$$\sum_J\chi_J^2\equiv1 \ \ \ \text{ everywhere on } \mathbb{T}^3\times \mathbb{R}$$
and additionally satisfies:
\begin{equation}\label{cutoffbounds}
    \supp_{x,t}\chi_J \subset B_{2/3\tau^c}(x_j)\times B_{2/3\tau^c}(t_j), \ \ \ ||\partial_t^{r'}\chi_J||_r\lesssim (\tau^c)^{-(r+r')}.
\end{equation}
We proceed by decomposing matrices close to the identity into simple tensors satisfying an orthogonality constraint, which will play a key role in the Nash stage when the fast coefficients are added (see the introductory Subsection \ref{overview}).
\begin{lemma}[Geometric Decomposition with Constraint]\label{geomlemma}
  Any matrix $M\in B_{1/2}(\IId)$ can be written as:
    \begin{equation}\label{matrixdecomp}
        M=\sum_{\zeta\in \Lambda} \gamma_\zeta^2(M) \zeta\otimes\zeta
    \end{equation}
where $\Lambda$ is a set of 6 unit vectors and
$\gamma_\zeta:B_{1/2}(\IId)\to \mathbb{R}$
are smooth functions. Moreover, each $\zeta\in \Lambda$ is associated with an oriented orthonormal basis $(k,\ \nu,\ \zeta)$, where $k$ satisfies $k\cdot e_3=0$.
\end{lemma}

This decomposition is a standard tool of convex integration; see, for example \cite[Appendix B]{BuckmasterDeLellisIsettSzekelyhidi2015}, for a proof and an account of its development. This gives the claim for a collection $\Lambda$ of six unit vectors $\zeta$. The fact that we can extend these to ordered orthonormal bases $(k, \ \nu, \ \zeta)$ with $k\cdot e_3=0$ follows from the fact that the intersection of planes $e_3^\perp\cap\zeta^\perp$ contains at least a line spanned by a unit vector $k$. We can finally set $\nu=\zeta\times k$.


\noindent 
Note that the commutation:
$$\partial_tB_{\ell,j}+\curl \left[B_{\ell,j}\times v_{\ell,j}\right]=0\leadsto \mathcal{A}^+_{\ell,j}\mathcal{A}^-_{\ell,j}=\mathcal{A}^-_{\ell,j}\mathcal{A}^+_{\ell,j}$$
coming from \eqref{localcorreq} (see the explicit computation in \eqref{pmcommutation}), together with the bounds in Lemma \ref{localcorrg}, ensures that the assumptions of Lemma \ref{chartconstr} are satisfied. This Lemma constructs for us charts $\Psi_J:Q_J\to \mathbb{R}^3$ adapted to $(v_{\ell,j},\ B_{\ell,j})$. We summarise their properties in the following Lemma, which is just a rephrasing of \ref{chartconstr} under the above notation. We refer the reader to Appendix \ref{diffgeom}, where we introduce and motivate the notational conventions for the differential-geometric objects.
\begin{lemma}[Adapted Charts]\label{chartpropg} Let $J=(j,j')$ a chart label as in \eqref{J} and $\underline{r}=M-m_0-2$. For any fixed $t\in B_{\tau^c}(t_j)$ the map: 
$$\Psi_J(\cdot,t):B_{\tau^c}(x_{j'})\to \mathbb{R}^3$$
is a diffeomorphism onto its image, with estimates:
\begin{equation*}
    \begin{split}
        &||\DD \Psi_J-\IId||_0,||(\DD\Psi_J)^{-1}-\IId||_0,||\det[\DD\Psi_J]-1||_0 \lesssim \lambda_{q+1}^{-\alpha},\\
        &||\DD \Psi_J||_r,||(\DD\Psi_J)^{-1}||_{r}\lesssim \lambda_q^{r}\lambda_q^{[r-\underline{r}]^+(b-1)\gamma_\ell} \ \text{ for } r \geq 1
    \end{split}
\end{equation*}
where the H{\"o}lder norms are understood to be taken on $Q_J$ and the implicit constants depend on $r,\ C_0,\ c_0, \ \alpha$.

\noindent Moreover, let $\zeta$ be a fixed unit vector, which we think of as a 2-form.  We have:
    \begin{equation*}
    (\partial_t+\mathcal{L}_{v_{\ell,j}})\Psi_J^{2*}\zeta=0 \ \text{ and }\ \mathcal{L}_{B_{\ell,j}}\Psi_J^{2*}\zeta=0.
\end{equation*}
\end{lemma}

\noindent \textbf{Decomposition of the Reynolds stress.} With the space-time charts and the geometric decomposition of matrices close to the identity, we are ready to rewrite the mollified Reynolds stress $R_\ell$. Using the matrix decomposition \eqref{matrixdecomp} we compute:
\begin{equation*}
    \begin{split}
        (\DD \Psi_J)^{-1}M(\DD \Psi_J)^{-T}=\sum_{\zeta\in \Lambda} \gamma_\zeta^2(M) (\DD \Psi_J)^{-1}\zeta\otimes (\DD \Psi_J)^{-1}\zeta=\sum_{\zeta\in \Lambda} \gamma_\zeta^2(M) \det[\DD\Psi_J]^{-2}\Psi_J^{2*}\zeta\otimes \Psi_J^{2*}\zeta\\
    \end{split}
\end{equation*}
and we deduce:
$$M=\sum_{\zeta\in \Lambda} \gamma_\zeta^2\left(\det[\DD\Psi_J]^{-2}\DD \Psi_J M\DD \Psi^{T}_J\right) \Psi_J^{2*}\zeta\otimes \Psi_J^{2*}\zeta$$
as soon as the composition with $\gamma_\zeta$ is well defined. This is indeed the case. Computations as in the proof of Lemma \ref{slowcoeffestimates} show that for $M=\IId-\frac{R_{\ell}}{\delta_{q+1}}$ we have:
$$||\det[\DD\Psi_J]^{-2}\DD \Psi_J M\DD \Psi_J^{T}-\IId||_0\lesssim \lambda_{q}^{-\alpha}$$
and in particular, for any given $\alpha, \ b$ we can choose $a$ large enough so that: 
$$\det[\DD\Psi_J]^{-2}\DD \Psi_J M\DD \Psi_J^{T}\in B_{1/2}(\IId)$$
and Lemma \ref{geomlemma} can be applied. With this and the partition of unity in \eqref{J} at hand, we can write:
\begin{equation}\label{reynoldsdecomposition}
    \begin{split}
        \delta_{q+1}\left(\IId-\frac{R_\ell}{\delta_{q+1}}\right)&=\sum_J\delta_{q+1}\chi_J^2\left(\IId-\frac{R_\ell}{\delta_{q+1}}\right)\\
        &=\sum_J\sum_{\zeta\in \Lambda} \delta_{q+1}\chi_J^2\gamma_\zeta^2\left(\det[\DD\Psi]^2\DD \Psi_J \left(\IId-R_\ell/\delta_{q+1}\right)\DD \Psi_J^T\right) \Psi_J^{2*}\zeta\otimes \Psi_J^{2*}\zeta\\
        &=\sum_J\sum_{\zeta\in \Lambda} a_{I}^2\Psi_{J}^{2*}\zeta\otimes \Psi_{J}^{2*}\zeta\\
        &=\sum_{I\in \mathcal{I}}A_I
    \end{split}
\end{equation}
where we defined $$I=(\zeta,J)=(I_\zeta,I_t,I_x)=(\zeta,j,j')$$ and $\mathcal{I}$ to be the finite set containing all the possible labels and $A_I$ the tensor associated with the index $I$ having coefficients $a_I$, namely
\begin{equation}\label{slowcoeffgalb}
\begin{split}
    a_I&=\delta_{q+1}^{1/2} \chi_J\gamma_\zeta\left(\det[\DD\Psi_J]^{-2}\DD \Psi_J \left(\IId-R_{\ell}/\delta_{q+1}\right)\DD \Psi_J^T\right),\\
    A_I&=A_{\zeta,J}=a_I^2\Psi_J^{2*}\zeta\otimes \Psi_J^{2*}\zeta.\\
\end{split}
\end{equation}
Note that $\Psi_j^{2*}\zeta$ is exactly the type of divergence-free vector field we constructed in Lemma \ref{chartpropg}. We now provide estimates on the decomposition.
\begin{lemma}[Estimates on the Decomposition]\label{slowcoeffgalbrun} Let $r\geq 0$ an integer and $\underline{r}=M-m_0-2$. For any $I=(\zeta,j,j')\in \mathcal{I}$ we have:
\begin{equation*}
    \begin{split}
        &||\partial_t^j a_I||_r\lesssim \lambda_q^{r+j}\lambda_q^{[r+j-\underline{r}]^+(b-1)\gamma_\ell}\delta_{q+1}^{1/2} \  \text{ for } \ j=0,1,\\
        &||\mathcal{A}^\pm_{\ell,I} a_I||_r\lesssim \lambda_q^{r}\lambda_q^{[r-\underline{r}]^+(b-1)\gamma_\ell}1/\tau^c\delta_{q+1}^{1/2}.
    \end{split}
\end{equation*}
Analogously, the tensor decomposition satisfies: 
\begin{equation*}
    \begin{split}
        &||\partial_t^j A_I||_r\lesssim \lambda_q^r\lambda_q^{[r-\underline{r}]^+(b-1)\gamma_\ell}\delta_{q+1} \  \text{ for } \ j=0,1,\\
        &||\mathcal{A}^\pm A_I||_r\lesssim \lambda_q^{r}\lambda_q^{[r-\underline{r}]^+(b-1)\gamma_\ell}1/\tau^c\delta_{q+1}.
    \end{split}
\end{equation*}
The implicit constants in both statements depend on $r,\ C_0,\ c_0, \ \alpha$.

\noindent Finally, we have the following support localisation property:
\begin{equation*} 
    \supp_{x,t} a_I,\ \supp_{x,t} A_I\Subset B_{2/3\tau^c}(x_{j'})\times B_{2/3\tau^c}(t_j)\cap \supp_{x,t}R_\ell.
\end{equation*}
\end{lemma} 

The proof of this Lemma is similar to the proof of Lemma \ref{slowcoeffestimates}, and we omit it. Here, in contrast with \ref{slowcoeffestimates}, we do not claim bounds for higher-order transport and pure time derivatives; this is why a preliminary mollification along the Alfv\'en directions, see Subsection \ref{timemollification}, is not necessary. This suffices, as we will use these objects only as the right-hand side of \eqref{galbrunapplied}, and the transport properties of the solution will be guaranteed by the equation itself.


\subsection{Galbrun LDF}\label{gnldf}
With the geometric constructions at hand, we proceed to the construction of the LDF. To better identify objects, we define the type map:
\begin{equation}\label{typemap}
    p: I \mapsto p(I)\in \Lambda \times \{e,o\}\times \{0,1,2,3\}
\end{equation}
distinguishing if the index is associated with an even or odd time interval, the family in the space partition, and the direction in the matrix decomposition. This is important because we will assign disjoint time-oscillating profiles to different types. Note that if two indices have the same type, then the associated cut-offs are disjoint in space by construction of the minimal covering, even if they share the same time interval and direction.  
\begin{lemma}[Time Profiles]\label{timeprofiles} There exists a collection of functions $\{g_p\}_p\subset C_c^\infty(0,1)$ indexed by
$$p\in \Lambda\times\{e,o\}\times \{0,1,2,3\}$$
such that: 
\begin{equation*}
        \int_0^1g_p=0, \ \int_0^1g^2_p=1, \ 
        \supp \ g_p \cap \supp \  g_{p'}=\emptyset \text{ for } p\neq p'.
\end{equation*}
We think of each $g_p$ as a periodic function $\mathbb{R}\to \mathbb{R}$ and define: 
\begin{equation*}
    \begin{cases}
        \alpha_{p}(t)= g^{[1]}_{p}(t)=\int_0^t g_{p}(s)\dd s,\\
        f_{p}=(1-g_{p}^2).
    \end{cases}
\end{equation*}
Additionally, the collection  $\{g_p\}$ can be chosen so that the resulting $\{\alpha_p\}_p$ satisfy for $p\neq p'$: 
$$\supp \ \alpha_p \cap \supp \ \alpha_{p'}=\emptyset, \ \supp \ \alpha_p \cap \supp \ g_{p'}=\emptyset.$$
\end{lemma}
The proof of the Lemma is a quick adaptation of \cite[Lemma 3.3]{GR}, and we omit it.  Without renaming, we rescale each profile by the parameter $\tau^a$ in \eqref{taua}, namely
$$g_p=g_p(t/\tau^a),\ \ \ f_p=f_p(t/\tau^a),\ \ \ \alpha_p(t)=\alpha_p(t/\tau^a).$$
As planned, we now engineer the LDF $\xi^{g}$ so that the leading order of the associated perturbation is a solution of Galbrun's equation with right-hand side
$$\ddiv\sum_{I}f_{p(I)}\mathbb{T} A_I,$$
we construct the potential of $\xi^{g}$ and work locally in time, that is, on $(t_j-\tau^c,t_j+\tau^c)$. We invite the reader to visit Section \ref{gn} before proceeding. 

\noindent According to the results in Subsection \ref{rew}, this corresponds to solving on $\mathbb{T}^3\times(t_j-\tau^c,t_j+\tau^c)$ the problem:
\begin{equation}\label{galbrunapplied}
\begin{cases}
    \mathcal{A}^-_{\ell,j}\mathcal{A}^+_{\ell,j}\Theta_j^{g}+\mathbb{H}_1[\Theta_j^{g}]+\mathbb{H}_2[\mathcal{A}^+_{\ell,j}\Theta_j^{g},\mathcal{A}^-_{\ell,j}\Theta_j^{g}]=\sum_{I:I_t=j}f_{p(I)}\mathbb{T} A_I,\\
    \Theta_j^{g}|_{t=t_j}=0,\\
    \mathcal{A}^+_{\ell,j}\Theta_j^{g}|_{t=t_j}=\tau^a\sum_{I:I_t=j}f_{p(I)}^{[1]}\mathbb{T} A_I|_{t=t_j}\\
\end{cases}
\end{equation}
where the initial data is chosen to simplify the a priori estimates from Subsection \ref{apriorie}, but plays no role here. Since $\supp \ \tilde \eta_j \subset (t_j-\tau^c,t_j+\tau^c)$, we can cut off the construction in time and define:
$$\xi_j^{g}=\curl [\tilde\eta_{j}\Theta_j^{g}]=\tilde\eta_{j}\curl [\Theta_j^{g}].$$
We now compute the cut-off error, since 
$$\tilde \eta_j\equiv 1 \ \text{ on } \  (t_j-2/3\tau^c,t_j+2/3\tau^c) \ \ \text{ and } \ \ \supp_t A_I\subset \supp \ \eta_j\subset (t_j-2/3\tau^c,t_j+2/3\tau^c),$$ 
see \eqref{slowcoeffgalb}, \eqref{supporteta},  we get from \eqref{galbrunapplied}: 
\begin{equation*}
\begin{split}
   \sum_{I:I_t=j}f_{p(I)}\mathbb{T} A_I&=\tilde{\eta}_j\sum_{I:I_t=j}f_{p(I)}\mathbb{T} A_I\\
    &=\tilde{\eta}_j\mathcal{A}^-_{\ell,j}\mathcal{A}^+_{\ell,j}\Theta_j^g+\tilde{\eta}_j\mathbb{H}_1[\Theta_j^g]+\tilde{\eta}_j\mathbb{H}_2[\mathcal{A}^+_{\ell,j}\Theta_j^g,\mathcal{A}^-_{\ell,j}\Theta_j^g]\\
    &=\mathcal{A}^-_{\ell,j}\mathcal{A}^+_{\ell,j}(\tilde{\eta}_j\Theta_j^{g})+\mathbb{H}_1[\tilde{\eta}_j\Theta_j^{g}]+\mathbb{H}_2[\mathcal{A}^+_{\ell,j}(\tilde{\eta}_j\Theta_j^{g}),\mathcal{A}^-_{\ell,j}(\tilde{\eta}_j\Theta_j^{g})]\\
    &-(\partial_t\tilde{\eta}_j)(\mathcal{A}^+_{\ell,j}\Theta_j^{g}+\mathcal{A}^-_{\ell,j}\Theta_j^{g})-(\partial_t^2\tilde{\eta}_j)\Theta_j^{g}-\partial_t\tilde{\eta}_j\mathbb{H}_2[\Theta_j^{g},\Theta_j^{g}].\\
\end{split}
\end{equation*}
We now set for $\mathcal{R}$ as in \eqref{invdiv}: 
\begin{equation}\label{Rcutj}
    R^{cut}_j=\mathcal{R}\curl\left[(\partial_t\tilde{\eta}_j)(\mathcal{A}^+_{\ell,j}\Theta_j^{g}+\mathcal{A}^-_{\ell,j}\Theta_j^{g})+(\partial_t^2\tilde{\eta}_j)\Theta_j^{g}+\partial_t\tilde{\eta}_j\mathbb{H}_2[\Theta_j^{g},\Theta_j^{g}]\right]\\
\end{equation}
and define, according to \eqref{firstapprox}, the principal part of the Galbrun perturbation to be:
\begin{equation*}
        w_j^{g,p}=(\partial_t+\mathcal{L}_{v_{\ell,j}})\xi_j^{g}, \ \ \ \ 
        b_j^{g,p}=\mathcal{L}_{B_{\ell,j}}\xi_j^{g},
\end{equation*}
where $p$ stands for principal. 

\noindent We now unwind the rewriting done in Subsection \ref{rew}. We take $\curl$ of \eqref{galbrunapplied} and solve for the associated Poisson equation to recover $\pi^{g}_j$ and remove the Leray projectors. We obtain that $(w_j^{g,p}, \ b_j^{g,p}, \ \pi^{g}_j)$ is a solution of: 
\begin{equation*}
\partial_tw^{g,p}_j+v_{\ell,j}\cn w^{g,p}-B_{\ell,j}\cn b^{g,p}_j+w^{g,p}\cn v_{\ell,j}-b^{g,p}_j\cn B_{\ell,j}+\nabla \pi^{gr}_j=\ddiv\left[\sum_{I:I_t=j}f_{p(I)} A_I+R_j^{cut}\right],\\
\end{equation*}
with associated Poisson problem:
\begin{equation}\label{poissong}
    \begin{cases}
        \Delta \pi^{gr}_j=-2\ddiv\left[w^{g,p}_j\cn v_{\ell,j}-b^{g,p}_j\cn B_{\ell,j}\right]+\ddiv \ddiv \left[\sum_{I:I_t=j}f_{p(I)}\ddiv A_I\right],\\
     \int_{\mathbb{T}^3} \pi^{g}_j=0.
    \end{cases}
\end{equation}
and we finally collect all the constructed LDFs, pressure terms and cut-off errors:  
\begin{equation*}
        \xi^{g}=\sum_j\xi^{g}_j, \ (w^{g,p},\ b^{g,p}, \ \pi^{g})=\sum_j (w^{g,p}_j,\ b^{g,p}_j,\ \pi_j^{g}), \ R^{cut}=\sum_j R_j^{cut},
\end{equation*}
note that the pressure is also localised in time, as the rhs of \eqref{poissong} is and abusing notation we also redefine:
$$\Theta^g_j:=\tilde\eta_j\Theta^g_j \leadsto \xi^g_j=\curl [\Theta^g_j].$$

The work done in Subsection \ref{apriorie} allows us to conclude that the following estimates hold.
\begin{lemma}[Galbrun LDF]\label{estgalbrun}Let $j\in \mathbb{Z}$, $r\geq0$ be an integer and $\underline{r}=M-m_0-4$. We have:
    \begin{subequations}\label{estgalbrun1}
        \begin{align}
            &||\partial_t^\sigma\Theta_j^{g}||_{r+\alpha}\lesssim \lambda_q^{r+\sigma}\lambda_{q}^{[r+\sigma-\underline{r}]^+(b-1)\gamma_\ell}\ell^{-\alpha}\tau^c\tau^a\delta_{q+1} \ \text{ for } \ \sigma=0,1,2,\\
            &||\partial_t^\sigma\mathcal{A}^\pm_{\ell,j}\Theta_j^{g}||_{r+\alpha}\lesssim \lambda_q^{r+\sigma}\lambda_{q}^{[r+\sigma-\underline{r}]^+(b-1)\gamma_\ell}\ell^{-\alpha}\tau^a\delta_{q+1} \ \text{ for } \ \sigma=0,1,2,\\
            &||\mathcal{A}^\pm_{\ell,j}\mathcal{A}^\mp_{\ell,j}\Theta_j^{g}||_{r+\alpha}, \ ||(\mathcal{A}^{\pm}_{\ell,j})^2\Theta_j^{g}||_{r+\alpha}\lesssim \lambda_q^r\lambda_{q}^{[r-\underline{r}]^+(b-1)\gamma_\ell}\ell^{-\alpha}\delta_{q+1}.
        \end{align}
    \end{subequations}
    In particular, 
    \begin{subequations}
    \begin{align}
        &||\partial_t^\sigma w_j^{g,p}||_{r+\alpha}, \ ||\partial_t^\sigma b_j^{g,p}||_{r+\alpha}\lesssim \lambda_q^{r+1+\sigma}\lambda_{q}^{[r+\sigma-(\underline{r}-1)]^+(b-1)\gamma_\ell}\ell^{-\alpha}\tau^a\delta_{q+1} \ \text{ for } \ \sigma=0,1,2,\\
        &||\mathcal{A}^\pm_{\ell,j}w_j^{g,p}||_{r+\alpha}, \ ||\mathcal{A}^\pm_{\ell,j}b_j^{g,p}||_{r+\alpha}\lesssim \lambda_q^{r+1}\lambda_{q}^{[r-(\underline{r}-1)]^+(b-1)\gamma_\ell}\ell^{-\alpha}\delta_{q+1},\\
        &||\pi^{gr}_j||_{r+\alpha}\lesssim \lambda_q^{r}\lambda_{q}^{[r-\underline{r}]^+(b-1)\gamma_\ell}\ell^{-\alpha}\delta_{q+1}\ \text{ for } \ r\geq 1,\\
        &||\mathcal{A}^\pm_{\ell,j} \pi^{gr}_j||_{r+\alpha}\lesssim \lambda_q^{r}\lambda_{q}^{[r-\underline{r}]^+(b-1)\gamma_\ell}\ell^{-\alpha}(1/\tau^a)\delta_{q+1}\ \text{ for } \ r\geq 1.
    \end{align}
    \end{subequations} 

    \noindent The implicit constants depend on $r, \ C_0, \ c_0, \ \alpha$.
\end{lemma}

\begin{remark}[Good Derivatives Loss]
    The equation loses two derivatives on the ones we have on the data, namely $M-m_0-2$, see Lemma \ref{slowcoeffgalbrun}, and there is an additional loss in passing to the vector fields because we have to move first from $\Theta_j^g$ to $\xi_j^g$, and this loses one derivative. 
\end{remark}

\begin{proof}[Proof of Lemma \ref{estgalbrun}] We want to apply Lemma \ref{fullestimates} to get estimates on the solution $\Theta_j^g$ of \eqref{galbrunapplied}. To do so, we need to verify that the Standing Assumptions \ref{standingass} are met. We notice immediately that the right-hand side of \eqref{galbrunapplied} is: 
$$\sum_{I:I_t=j}f_{p(I)}\mathbb{T} A_I,$$
which, strictly speaking, is not of the form $f\mathbb{T}F$; however, at each space-time point, by construction, at most one of the $f_{p(I)}\mathbb{T} A_I$ is non-zero, and the proof of Lemma \ref{fullestimates} can be adapted to this case. Indeed, the integration by parts in time, type of arguments in the proof of \ref{fullestimates} would, in this case, involve finitely many time profiles, the number of which is fixed independently of $q$, namely $48$. The bounds on $(v_{\ell,j}, \ B_{\ell,j}, \ p_{\ell,j}, \ R_\ell)$ can be found in Lemmas \ref{standardmollgalbrun} and \ref{localcorrg}. They solve the relaxed MHD equation by construction, see \eqref{localcorreq}. The bounds on the right-hand side then follow from Lemma \ref{slowcoeffgalbrun}. We set the parameter $\underline{r}$ in the Standing Assumptions \ref{standingass} to be $\underline{r}=M-m_0-2$. The fact that we actually redefined $\Theta_j^g$ to be $\tilde\eta_j\Theta_j^g$ does not change the estimates, because of the bound in \eqref{timecutoffproperties}. We conclude that the bounds in \eqref{estgalbrun1} hold.

\noindent \textbf{From $\Theta_j^g$ to $(w_j^{g,p}, \ b_j^{g,p})$.} We now use the definition of $\underline{r}$ in the statement of the Lemma that is $\underline{r}=M-m_0-4$. We can transfer the estimates from the potential to the vector field with the help of the following computation, which relies on \eqref{DG1}:
\begin{equation*}
    \begin{split}
        w_j^{g,p}&=(\partial_t+\mathcal{L}_{v_{\ell,j}})\xi_j^{g}\\
                &=(\partial_t+\mathcal{L}_{v_{\ell,j}})\curl[\tilde\eta_j\Theta_j^g]\\
                &=\curl\left[(\partial_t\tilde\eta_j)\Theta_j^g+1/2\tilde\eta_j(\mathcal{A}^+_{\ell,j}+\mathcal{A}^-_{\ell,j})\Theta_j^g+(\DD v^\pm_{\ell,j})^\top[\tilde\eta_j\Theta_j^g]\right].
    \end{split}
\end{equation*}
Similarly, for $b_j^{g,p}$. From this commuting the time derivative with $\curl$, and the bounds in \eqref{estgalbrun1} we deduce:
$$||\partial_t^\sigma w_j^{g,p}||_{r+\alpha}, \ ||\partial_t^\sigma b_j^{g,p}||_{r+\alpha}\lesssim \lambda_q^{r+1+\sigma}\lambda_{q}^{[r+\sigma-(\underline{r}-1)]^+(b-1)\gamma_\ell}\ell^{-\alpha}\tau^a\delta_{q+1} \ \text{ for } \ \sigma=0,1,2 \ \text{ and } \ r\geq 0.$$
Note the loss of one good derivative. We use the same idea to deduce the Alfv\'en transport bound.
\begin{equation*}
    \begin{split}
        \mathcal{A}^\pm_{\ell,j}w_j^{g,p}&=(\partial_t+\mathcal{L}_{z_{\ell,j}^\pm})w_j^{g,p}+w_j^{g,p}\cn z^\pm_{\ell,j}\\
        &=(\partial_t+\mathcal{L}_{z_{\ell,j}^\pm})\curl\left[(\partial_t\tilde\eta_j)\Theta_j^g+1/2\tilde\eta_j(\mathcal{A}^+_{\ell,j}+\mathcal{A}^-_{\ell,j})\Theta_j^g+(\DD v^\pm_{\ell,j})^\top[\tilde\eta_j\Theta_j^g]\right]\\
        &+w_j^{g,p}\cn z^\pm_{\ell,j}\\
        &=\curl\left[(\partial_t^2\tilde\eta_j)\Theta_j^g+(\partial_t\tilde\eta_j)\mathcal{A}^\pm_{\ell,j}\Theta_j^g+1/2\tilde\eta_j\mathcal{A}^\pm_{\ell,j}(\mathcal{A}^+_{\ell,j}+\mathcal{A}^-_{\ell,j})\Theta_j^g+1/2(\partial_t\tilde\eta_j)(\mathcal{A}^+_{\ell,j}+\mathcal{A}^-_{\ell,j})\Theta_j^g\right]\\
        &+\curl \left[(\DD \mathcal{A}^\pm_{\ell,j}v^\pm_{\ell,j})^\top[\tilde\eta_j\Theta_j^g]+(\DD v^\pm_{\ell,j})^\top[(\partial_t\tilde\eta_j)\Theta_j^g+\tilde\eta_j\mathcal{A}^\pm_{\ell,j}\Theta_j^g]\right]\\
        &+\curl\left[(\DD z^\pm_{\ell,j})^\top\left[(\partial_t\tilde\eta_j)\Theta_j^g+1/2\tilde\eta_j(\mathcal{A}^+_{\ell,j}+\mathcal{A}^-_{\ell,j})\Theta_j^g+(\DD v^\pm_{\ell,j})^\top[\tilde\eta_j\Theta_j^g]\right]\right]\\
        &+w_j^{g,p}\cn z^\pm_{\ell,j}\\
    \end{split}
\end{equation*}
Similarly, for $b_j^{g,p}$. We deduce from \eqref{estgalbrun1} the bound above and the estimates in Lemma \ref{localcorrg} that:
\begin{equation}\label{leadinggalbrun}
    ||\mathcal{A}^\pm_{\ell,j}w_j^{g,p}||_{r+\alpha}, \ ||\mathcal{A}^\pm_{\ell,j}b_j^{g,p}||_{r+\alpha}\lesssim \lambda_q^{r+1}\lambda_{q}^{[r-(\underline{r}-1)]^+(b-1)\gamma_\ell}\ell^{-\alpha}\delta_{q+1} \ \text{ for } \ r\geq 0.
\end{equation}

\noindent \textbf{Estimates on $\nabla \pi^g$.} We will use the bounds in \eqref{leadinggalbrun} without mention. By definition, we have:
\begin{equation*}
    \begin{split}
        \nabla\pi^{gr}_j=-2\nabla\Delta^{-1}\ddiv\left[w^{g,p}_j\cn v_{\ell,j}-b^{g,p}_j\cn B_{\ell,j}\right]+\nabla\Delta^{-1}\ddiv \ddiv \underbrace{\sum_{I:I_t=j}f_{p(I)}A_I}_{=:F}\\
    \end{split}
\end{equation*}
We will now use that $\nabla\Delta^{-1}\ddiv$ is a CZ-type operator together with Proposition \ref{czstuff} to conclude estimates on $\nabla \pi^{gr}_j$. We have:
\begin{equation*}
    \begin{split}
        ||\nabla \pi^{gr}_j||_{r+\alpha} &\lesssim ||v_{\ell,j}||_{r+1+\alpha}||w^{g,p}_j||_0+ ||v_{\ell,j}||_{1}||w^{g,p}_j||_{r+\alpha}+||B_{\ell,j}||_{r+1+\alpha}||b^{g,p}_j||_0+||B_{\ell,j}||_1||b^{g,p}_j||_{r+\alpha}\\
        &+||F||_{1+\alpha}\\
        &\lesssim\lambda_q^{r+2}\lambda_{q}^{[r-(\underline{r}-1)]^+(b-1)\gamma_\ell}\ell^{-2\alpha}\tau^a\delta_q^{1/2}\delta_{q+1}+\lambda_q^{r+1}\lambda_{q}^{[r+1-(M-m_0-2)]^+(b-1)\gamma_\ell}\ell^{-\alpha}\delta_{q+1}\\
        &\leq\lambda_q^{r+1}\lambda_{q}^{[r-(\underline{r}-1)]^+(b-1)\gamma_\ell}\ell^{-\alpha}\delta_{q+1}.\\
    \end{split}
\end{equation*}
By means of Lemmas \ref{slowcoeffgalbrun}, \ref{stabilitygalbrun} and Proposition \ref{czstuff} to deal with the commutators, we now bound:
\begin{equation*}
    \begin{split}
        ||\mathcal{A}^\pm_{\ell,j} \nabla \pi^{gr}_j||_{r+\alpha}&\lesssim ||[z_{\ell,j}^\pm,\nabla\Delta^{-1}\ddiv][-2w^{g,p}_j\cn v_{\ell,j}+2b^{g,p}_j\cn B_{\ell,j}+\ddiv F]||_{r+\alpha}\\
        &+||\mathcal{A}^\pm_{\ell,j}[-2w^{g,p}_j\cn v_{\ell,j}+2b^{g,p}_j\cn B_{\ell,j}+\ddiv F]||_{r+\alpha}\\
        &\lesssim \underbrace{||z_{\ell,j}^\pm||_{r+1+\alpha}(||w^{g,p}_j\cn v_{\ell,j}||_\alpha+||\ddiv F||_\alpha)+||z_{\ell,j}^\pm||_{1+\alpha}(||w^{g,p}_j\cn v_{\ell,j}||_{r+\alpha}+||\ddiv F||_{r+\alpha})}_{\lambda_q^{r+3}\ell^{-3\alpha}\tau^a\delta_q\delta_{q+1}+\lambda_q^{r+2}\ell^{-3\alpha}\tau^a\delta_q^{1/2}\delta_{q+1}}\\
        &+\underbrace{||(\mathcal{A}^\pm_{\ell,j}w^{g,p}_j)\cn v_{\ell,j}||_{r+\alpha}}_{\lambda_q^{r+3}\ell^{-2\alpha}\delta_q^{1/2}\delta_{q+1}}+\underbrace{||w^{g,p}_j\cn \mathcal{A}^\pm_{\ell,j}v_{\ell,j}||_{r+\alpha}}_{\lambda_q^{r+3}\ell^{-2\alpha}\tau^a\delta_q\delta_{q+1}}+\underbrace{||(w^{g,p}_j\cn z^\pm_{\ell,j})\cn v_{\ell,j}||_{r+\alpha}}_{\lambda_q^{r+3}\ell^{-2\alpha}\tau^a\delta_q\delta_{q+1}}\\
        &+\underbrace{||\ddiv\mathcal{A}^\pm_{\ell,j} F||_{r+\alpha}}_{\lambda_q^{r+1}\ell^{-\alpha}(1/\tau^a)\delta_{q+1}}+\underbrace{||[z^\pm_{\ell,j}\cn ,\ddiv] F||_{r+\alpha}}_{\lambda_q^{r+2}\ell^{-\alpha}\delta_q^{1/2}\delta_{q+1}}\\
        &+\text{magnetic field terms}\\
        &\lesssim \lambda_q^{r+1}\lambda_{q}^{[r-(\underline{r}-1)]^+(b-1)\gamma_\ell}\ell^{-\alpha}(1/\tau^a)\delta_{q+1},
    \end{split}
\end{equation*}
and from the identity
$$\nabla \mathcal{A}^\pm_{\ell,j} \pi^{gr}_j=\mathcal{A}^\pm_{\ell,j} \nabla \pi^{gr}_j+(\DD z^\pm_{\ell,j})^\top \nabla \pi^{gr}_j,$$ 
we conclude that:
$$||\nabla \mathcal{A}^\pm_{\ell,j} \pi^{gr}_j||_r\lesssim \lambda_q^{r+1}\lambda_{q}^{[r-(r-1)]^+(b-1)\gamma_\ell}\ell^{-\alpha}(1/\tau^a)\delta_{q+1} \ \text{ for } \ r\geq 0.$$

\noindent The implicit constants in the above bounds depend on $r, \ C_0, \ c_0, \ \alpha$.
\end{proof}

\begin{lemma}[Estimates on $R^{cut}$]\label{Rcut} Let $r\geq 0$ be an integer and $\underline{r}=M-m_0-4$. We have:
    \begin{equation*}
        \begin{split}
            &||R^{cut}||_{r+\alpha}\lesssim \lambda_q^{r}\lambda_q^{[r-\underline{r}]^+(b-1)\gamma_\ell}\ell^{-\alpha}(\tau^a/\tau^c)\delta_{q+1},\\
            &||\mathcal{A}^\pm R^{cut}||_{r+\alpha}\lesssim \lambda_q^{r}\lambda_q^{[r-\underline{r}]^+(b-1)\gamma_\ell}\ell^{-\alpha}(1/\tau^c)\delta_{q+1} \ \text{ for } \ 0\leq r\leq N-m_0-1.
        \end{split}
    \end{equation*}
The implicit constants depend on $r, \ C_0, \ c_0, \ \alpha$.
\end{lemma}
\begin{remark}
    Note that with the current Iterative Assumptions \eqref{inductiveassumptionsgeneral}, we are not able to guarantee the improved transport estimates as in \cite{GR}. 
\end{remark}

\begin{proof}[Proof of Lemma \ref{Rcut}] We first use the definition of $\mathbb{H}_2$ and $\mathbb{T}$ given in \eqref{Hdefinition} to rewrite $R^{cut}_j$ as:
\begin{equation*}
    \begin{split}
        R^{cut}_j&=\mathcal{R}\curl\left[(\partial_t\tilde{\eta}_j)(\mathcal{A}^+_{\ell,j}\Theta_j^{g}+\mathcal{A}^-_{\ell,j}\Theta_j^{g})+(\partial_t^2\tilde{\eta}_j)\Theta_j^{g}+\partial_t\tilde{\eta}_j\mathbb{H}_2[\Theta_j^{g},\Theta_j^{g}]\right]\\
        &=\mathcal{R}\curl\left[(\partial_t\tilde{\eta}_j)(\mathcal{A}^+_{\ell,j}\Theta_j^{g}+\mathcal{A}^-_{\ell,j}\Theta_j^{g})+(\partial_t^2\tilde{\eta}_j)\Theta_j^{g}\right]\\
        &+\mathcal{R}\curl\left[\partial_t\tilde{\eta}_j(\DD z_{\ell,j}^-)^\top [\Theta_j^{g}] +\partial_t\tilde{\eta}_j(\DD z_{\ell,j}^+)^\top [\Theta_j^{g}]+\partial_t\tilde{\eta}_j\mathbb{T}\left[\Theta_j^{g} \times \nabla z_{\ell,j}^-+\Theta_j^{g} \times \nabla z_{\ell,j}^+\right]^\top\right]\\
        &=\underbrace{\mathcal{R}\curl\left[(\partial_t\tilde{\eta}_j)(\mathcal{A}^+_{\ell,j}\Theta_j^{g}+\mathcal{A}^-_{\ell,j}\Theta_j^{g})+(\partial_t^2\tilde{\eta}_j)\Theta_j^{g}\right]}_{T_1}\\
        &+\underbrace{\mathcal{R}\curl\left[\partial_t\tilde{\eta}_j(\DD z_{\ell,j}^-)^\top [\Theta_j^{g}] +\partial_t\tilde{\eta}_j(\DD z_{\ell,j}^+)^\top [\Theta_j^{g}] \right]+\mathcal{R}\mathbb{P}\ddiv\left[\partial_t\tilde{\eta}_j\left[\Theta_j^{g} \times \nabla z_{\ell,j}^-+\Theta_j^{g} \times \nabla z_{\ell,j}^+\right]^\top\right]}_{T_2}.\\
    \end{split}
\end{equation*}

\noindent \textit{Estimates on $T_1$.} Since $\mathcal{R}\curl$ is a CZ-type operator, we can use Proposition \ref{czstuff} to deal with it. The bounds in Lemmas \ref{localcorrg}, \ref{estgalbrun} and \eqref{timecutoffproperties}, then give:
\begin{equation*}
    \begin{split}
        ||T_1||_{r+\alpha}&\lesssim ||\partial_t\tilde{\eta}_j\mathcal{A}^+_{\ell,j}\Theta_j^{g}||_{r+\alpha}+||\partial_t\tilde{\eta}_j\mathcal{A}^-_{\ell,j}\Theta_j^{g}||_{r+\alpha}+||(\partial_t^2\tilde{\eta}_j)\Theta_j^{g}||_{r+\alpha}\\
        &\lesssim \lambda_q^{r}\lambda_q^{[r-\underline{r}]^+(b-1)\gamma_\ell}\ell^{-\alpha}(\tau^a/\tau^c)\delta_{q+1}.
    \end{split}
\end{equation*}
Using in adddition the commutator estimates in Proposition \ref{czstuff} and the bounds in Lemmas \ref{standardmollgalbrun}, \ref{stabilitygalbrun} to correct the transport operator, we obtain:
\begin{equation*}
    \begin{split}
        ||\mathcal{A}^\pm T_1||_{r+\alpha}&\lesssim ||[z_{\ell,j}^\pm,\mathcal{R}\curl][(\partial_t\tilde{\eta}_j)(\mathcal{A}^+_{\ell,j}\Theta_j^{g}+\mathcal{A}^-_{\ell,j}\Theta_j^{g})+(\partial_t^2\tilde{\eta}_j)\Theta_j^{g}]||_{r+\alpha}\\
        &+||(z_q^\pm-z_{\ell,j}^\pm)\cn\mathcal{R}\curl[(\partial_t\tilde{\eta}_j)(\mathcal{A}^+_{\ell,j}\Theta_j^{g}+\mathcal{A}^-_{\ell,j}\Theta_j^{g})+(\partial_t^2\tilde{\eta}_j)\Theta_j^{g}]||_{r+\alpha}\\
        &+ \underbrace{||\partial_t^2\tilde{\eta}_j(\mathcal{A}^+_{\ell,j}\Theta_j^{g}+\mathcal{A}^-_{\ell,j}\Theta_j^{g})||_{r+\alpha}}_{\lambda_q^{r}\ell^{-\alpha}(\tau^a/\tau^c)(1/\tau^c)\delta_{q+1}}+\underbrace{||\partial_t\tilde{\eta}_j\mathcal{A}^\pm_{\ell,j}\mathcal{A}^+_{\ell,j}\Theta_j^{g}||_{r+\alpha}+||\partial_t\tilde{\eta}_j\mathcal{A}^\pm_{\ell,j}\mathcal{A}^-_{\ell,j}\Theta_j^{g}||_{r+\alpha}}_{\lambda_q^{r}\ell^{-\alpha}(1/\tau^c)\delta_{q+1}}\\
        &+\underbrace{||(\partial_t^3\tilde{\eta}_j)\Theta_j^{g}||_{r+\alpha}+||(\partial_t^2\tilde{\eta}_j)\mathcal{A}^\pm_{\ell,j}\Theta_j^{g}||_{r+\alpha}}_{\lambda_q^{r}\ell^{-\alpha}(\tau^a/\tau^c)(1/\tau^c)\delta_{q+1}}\\
        &\lesssim\lambda_q^{r}\lambda_q^{[r-\underline{r}]^+(b-1)\gamma_\ell}\ell^{-\alpha}(1/\tau^c)\delta_{q+1}\\
    \end{split}
\end{equation*}
for $0\leq r \leq N-m_0-1$.

\noindent \textit{Estimates on $T_2$.} We show explicit bounds only for the first term, namely $\mathcal{R}\curl\partial_t\tilde{\eta}_j(\DD z_{\ell,j}^-)^\top [\Theta_j^{g}]$, the others can be treated similarly. As above, we deduce:
\begin{equation*}
        ||\mathcal{R}\curl\partial_t\tilde{\eta}_j(\DD z_{\ell,j}^-)^\top [\Theta_j^{g}]||_{r+\alpha}\lesssim ||\partial_t\tilde{\eta}_j(\DD z_{\ell,j}^-)^\top [\Theta_j^{g}]||_{r+\alpha} \lesssim \lambda_q^{r}\lambda_q^{[r-\underline{r}]^+(b-1)\gamma_\ell}\ell^{-\alpha}(\tau^a/\tau^c)\delta_{q+1}.
\end{equation*}
Now note that:
\begin{equation*}
    \begin{split}
        \mathcal{A}^\pm_{\ell,j}\left[\partial_t\tilde{\eta}_j(\DD z_{\ell,j}^-)^\top [\Theta_j^{g}]\right]&=\partial_t^2\tilde{\eta}_j(\DD z_{\ell,j}^-)^\top [\Theta_j^{g}]+\partial_t\tilde{\eta}_j(\DD \mathcal{A}^\pm_{\ell,j}z_{\ell,j}^-)^\top [\Theta_j^{g}]-\partial_t\tilde{\eta}_j(\DD z_{\ell,j}^\pm)^\top(\DD z_{\ell,j}^-)^\top [\Theta_j^{g}]\\
        &+\partial_t\tilde{\eta}_j(\DD z_{\ell,j}^-)^\top [\mathcal{A}^\pm_{\ell,j}\Theta_j^{g}]
    \end{split}
\end{equation*}
and as above, we conclude that:
\begin{equation*}
    \begin{split}
        &||\mathcal{A}^\pm \mathcal{R}\curl\partial_t\tilde{\eta}_j(\DD z_{\ell,j}^-)^\top [\Theta_j^{g}]||_{r+\alpha}\\
        &\lesssim ||[z_{\ell,j}^\pm,\mathcal{R}\curl][\partial_t\tilde{\eta}_j(\DD z_{\ell,j}^-)^\top [\Theta_j^{g}]]||_{r+\alpha}\\
        &+||(z_q^\pm-z_{\ell,j}^\pm)\cn\mathcal{R}\curl[\partial_t\tilde{\eta}_j(\DD z_{\ell,j}^-)^\top [\Theta_j^{g}]]||_{r+\alpha}\\
        &+ \underbrace{||\partial_t^2\tilde{\eta}_j(\DD z_{\ell,j}^-)^\top [\Theta_j^{g}]||_{r+\alpha}}_{\lambda_q^{r}\ell^{-\alpha}(\tau^a/\tau^c)(1/\tau^c)\delta_{q+1}}+\underbrace{||\partial_t\tilde{\eta}_j(\DD \mathcal{A}^\pm_{\ell,j}z_{\ell,j}^-)^\top [\Theta_j^{g}]||_{r+\alpha}}_{\lambda_q^{r+2}\ell^{-2\alpha}\tau^a\delta_q\delta_{q+1}}+\underbrace{||\partial_t\tilde{\eta}_j(\DD z_{\ell,j}^\pm)^\top(\DD z_{\ell,j}^-)^\top [\Theta_j^{g}]||_{r+\alpha}}_{\lambda_q^{r+2}\ell^{-3\alpha}\tau^a\delta_q\delta_{q+1}}\\
        &\lesssim ||z_{\ell,j}^\pm||_{r+\alpha}||\partial_t\tilde{\eta}_j(\DD z_{\ell,j}^-)^\top [\Theta_j^{g}]||_\alpha+||z_{\ell,j}^\pm||_{\alpha}||\partial_t\tilde{\eta}_j(\DD z_{\ell,j}^-)^\top [\Theta_j^{g}]||_{r+\alpha}\\
        &+||z_q^\pm-z_{\ell,j}^\pm||_{r+\alpha}||\partial_t\tilde{\eta}_j(\DD z_{\ell,j}^-)^\top [\Theta_j^{g}]||_{1+\alpha}+||z_q^\pm-z_{\ell,j}^\pm||_{\alpha}||\partial_t\tilde{\eta}_j(\DD z_{\ell,j}^-)^\top [\Theta_j^{g}]||_{r+1+\alpha}\\
        &+\lambda_q^{r}\lambda_q^{[r-\underline{r}]^+(b-1)\gamma_\ell}\ell^{-\alpha}(\tau^a/\tau^c)(1/\tau^c)\delta_{q+1}\\
        &\lesssim\lambda_q^{r}\lambda_q^{[r-\underline{r}]^+(b-1)\gamma_\ell}\ell^{-\alpha}(\tau^a/\tau^c)(1/\tau^c)\delta_{q+1}\\
    \end{split}
\end{equation*}
for $0\leq r \leq N-m_0-1$.

\noindent \textit{Conclusion.} From the fact that at most two $\tilde{\eta}_j$ are non-zero at the same time we conclude that:
\begin{equation*}
    \begin{split}
        &||R^{cut}||_{r+\alpha}\lesssim \sup_j||R_j^{cut}||_{r+\alpha}\lesssim \lambda_q^{r}\lambda_q^{[r-\underline{r}]^+(b-1)\gamma_\ell}\ell^{-\alpha}(\tau^a/\tau^c)\delta_{q+1} \ \text{ for } \ r\geq 0,\\
        &||\mathcal{A}^\pm R^{cut}||_{r+\alpha}\lesssim \sup_j||\mathcal{A}^\pm R_j^{cut}||_{r+\alpha}\lesssim \lambda_q^{r}\lambda_q^{[r-\underline{r}]^+(b-1)\gamma_\ell}\ell^{-\alpha}(1/\tau^c)\delta_{q+1} \ \text{ for } \ 0\leq r \leq N-m_0-1.\\
    \end{split}
\end{equation*}
The implicit constants depend on $r, \ C_0, \ c_0, \ \alpha$.
\end{proof}

\subsubsection{Splitting of the Perturbation.}
Following Lemma \ref{lpl} we construct $(w^{g}, b^{g})$ as a perturbation of $(v_q, B_q)$ along the LDF $\xi^{g}$. In the following, we split $(w^{g}, b^{g})$ in simpler components.

\noindent\textit{Velocity Field.} Since $\xi^g$ is designed to behave correctly with $(v_{\ell,j}, \ B_{\ell,j})$, which are not globally defined in time, we need to be extra careful. Let us first separate the non-mollified part:
\begin{equation}\label{split1g}
    \begin{split}
        w^{g}&=\partial_tX^{g}\circ (X^{g})^{-1}+(X^{g}_*-\IId_*)v_q\\
        &=\underbrace{\partial_tX^{g}\circ (X^{g})^{-1}+(X^{g}_*-\IId_*)v_\ell}_{\bar w^g}+\underbrace{(X^{g}_*-\IId_*)(v_q-v_\ell)}_{\mathring w^g}
    \end{split}
\end{equation}
and according to Lemma \ref{lpl} we have:
\begin{equation}\label{decompgalbrunw}
    \begin{split}
        \bar w^g 
            &=\curl\left[\sum_j\sum_{k=0}^{k_0^g}\frac{(-1)^k}{(k+1)!}\mathcal{L}_{\xi^g}^k(\partial_t+\mathcal{L}_{v_{\ell,j}})\Theta^g_j+\sum_j\sum_{k=0}^{k_0^g}\frac{(-1)^k}{(k+1)!}\mathcal{L}_{\xi^g}^k\mathcal{L}_{v_{\ell}- v_{\ell,j}}\Theta_j^g+\theta^g_w\right],\\
        \mathring w^g
                    &=\curl\left[\sum_j\sum_{k=0}^{k_0^g}\frac{(-1)^k}{(k+1)!}\mathcal{L}_{\xi^g}^k\mathcal{L}_{v_q- v_\ell}\Theta_j^g+\mathring\theta^g_w\right]\\
    \end{split}
\end{equation}
where
\begin{equation*}
    \begin{split}
        \theta^g_w
        &=\sum_j\frac{(-1)^{k_0^g+1}}{(k_0^g+1)!}\int_0^{1}(X_{s}^g)_*\left[\mathcal{L}_{\xi^{g}}^{k_0^g+1}(\partial_t+\mathcal{L}_{v_{\ell,j}})\Theta^{g}_j\right](1-s)^{k_0^g+1}\dd s\\
        &+\sum_j\frac{(-1)^{k_0^g+1}}{(k_0^g+1)!}\int_0^{1}(X_{s}^g)_*\left[\mathcal{L}_{\xi^{g}}^{k_0^g+1}\mathcal{L}_{v_\ell-v_{\ell,j}}\Theta^{g}_j\right](1-s)^{k_0^g+1}\dd s,\\
        \mathring\theta^g_w&=\sum_j\frac{(-1)^{k_0^g+1}}{(k_0^g+1)!}\int_0^{1}(X_{s}^g)_*\left[\mathcal{L}_{\xi^{g}}^{k_0^g+1}\mathcal{L}_{v_q-v_\ell}\Theta^{g}_j\right](1-s)^{k_0^g+1}\dd s.\\
    \end{split}
\end{equation*}

\noindent \textit{Magnetic field.} We proceed similarly, we first split:
\begin{equation*}
    \begin{split}
        &b^g=\underbrace{(X^{g}_*-\IId_*)B_\ell}_{\bar{b}^p}+\underbrace{(X^{g}_*-\IId_*)(B_q-B_\ell)}_{\mathring b^p},
    \end{split}
\end{equation*}
then Lemma \ref{lpl} gives:
\begin{equation}\label{decompgalbrunb}
    \begin{split}
        \bar b^g
            &=\curl\left[\sum_j\sum_{k=0}^{k_0^g}\frac{(-1)^k}{(k+1)!}\mathcal{L}_{\xi^g}^k\mathcal{L}_{B_{\ell,j}}\Theta^g_j+\sum_j\sum_{k=0}^{k_0^g}\frac{(-1)^k}{(k+1)!}\mathcal{L}_{\xi^g}^k\mathcal{L}_{B_{\ell}- B_{\ell,j}}\Theta_j^g+\theta^g_b\right],\\
        \mathring b^g
                    &=\curl\left[\sum_j\sum_{k=0}^{k_0^g}\frac{(-1)^k}{(k+1)!}\mathcal{L}_{\xi^g}^k\mathcal{L}_{B_q- B_\ell}\Theta_j^g+\mathring\theta^g_b\right]\\
    \end{split}
\end{equation}
where 
\begin{equation*}
    \begin{split}
        \theta^g_b&=\sum_j\frac{(-1)^{k_0^g+1}}{(k_0^g+1)!}\int_0^{1}(X_{s}^g)_*\left[\mathcal{L}_{\xi^{g}}^{k_0^g+1}\mathcal{L}_{B_{\ell,j}}\Theta^{g}_j\right](1-s)^{k_0^g+1}\dd s\\
        &+\sum_j\frac{(-1)^{k_0^g+1}}{(k_0^g+1)!}\int_0^{1}(X_{s}^g)_*\left[\mathcal{L}_{\xi^{g}}^{k_0^g+1}\mathcal{L}_{B_\ell-B_{\ell,j}}\Theta^{g}_j\right](1-s)^{k_0^g+1}\dd s,\\
        \mathring\theta^g_b&=\sum_j\frac{(-1)^{k_0^g+1}}{(k_0^g+1)!}\int_0^{1}(X_{s}^g)_*\left[\mathcal{L}_{\xi^{g}}^{k_0^g+1}\mathcal{L}_{B_q-B_\ell}\Theta^{g}_j\right](1-s)^{k_0^g+1}\dd s.\\
    \end{split}
\end{equation*}

\subsection{Estimates on the Perturbation}\label{estperturbationgalbrun}
In this subsection, we study the higher-order explicit terms and remainders of the Lie-Taylor series associated with the perturbation $(w^g, \ b^g)$ of $(v_q, \ B_q)$ as in \eqref{decompgalbrunw}, \eqref{decompgalbrunb}.

\noindent Key quantities will be: 
\begin{equation}\label{Tg}
    \mathcal{T}_g=\lambda_q^2\ell^{-\alpha}\tau^a\tau^c\delta_{q+1}=\ell^{-\alpha}\left(\frac{\lambda_q}{\lambda_{q+1}}\right)^{\gamma_a+2\beta+2\gamma_{CZ}}
\end{equation}
which measures the smallness gain corresponding to one Lie-derivation, and objects which we call loss functions, $L, \ L_{\mathcal{A}}:\mathbb{N}_{\geq 0}\to \mathbb{R}_{\geq 1}$, taking the form:
\begin{equation*}
    \begin{split}
        &L(r)=\lambda_q^{[r-\underline{r}](b-1)\gamma_\ell},\\
        &L_{\mathcal{A}}(r)=1_{r\leq \underline{r}}+1_{r\geq \underline{r}+1} \lambda_q^{[r-\underline{r}](b-1)\gamma_\ell}\bar L,
    \end{split}
\end{equation*}
for parameters $\underline{r}$ and $\bar L$ e.g. $\bar L=\delta_q^{-1/2}$. These keep track of the range of derivatives $0\leq r\leq \underline{r}$ over which we have the desired estimates and the loss that occurs afterwards. The subscript $\mathcal{A}$ will indicate that the loss function is related to an Alfv\'en transport estimate.

\subsubsection{Estimates on the Lagrangian Perturbation}
Let $X^g:[-1,1]_s\times\mathbb{T}^3_x\times \mathbb{R}_t\to \mathbb{T}^3$ solve:
\begin{equation}\label{lagrangiang}
    \begin{cases}
        \partial_sX_s^g(x,t)=\xi^g(X_s^g(x,t),t),\\
        X_0^g(x,t)=x.
    \end{cases}
\end{equation}
In Lemma \ref{estgalbrun} we gave estimates for $\xi^g=\curl \sum_j\Theta^g_j$. After choosing $a$ sufficiently large to reabsorb the implicit constant, we can assume 
\begin{equation}\label{derivativelagrangiang}
    ||\DD \xi||_0\leq C\mathcal{T}_g\leq C\left(\frac{\lambda_q}{\lambda_{q+1}}\right)^{\gamma_a+2\beta}\leq 1
\end{equation}
and from this we deduce the following Lemma.

\begin{lemma}[Estimates on $X^{g}$]\label{estgalbrunflow} Let $r\geq 0$ integer and $\underline{r}=M-m_0-6$. The Lagrangian perturbation $X^g$ associated with $\xi^g$ satisfies:
\begin{equation*}
    \begin{split}
        &||\DD X^g_s-\IId||_0\lesssim \mathcal{T}_g \ \ \text{ for } \ \ |s|\leq 1,\\
        &||\partial_t^j X^g_s||_0, \lesssim\lambda_q^{j-1}\mathcal{T}_g \ \ \text{ for } \ \ j= 1,2 \ \ \&  \ \  |s|\leq 1, \\
        &||\partial_t^j \DD X^{g}_s||_r\lesssim\lambda_q^{r+j}\lambda_q^{[r+j-\underline{r}]^+(b-1)\gamma_\ell}\mathcal{T}_g \ \ \text{ for } \ j= 0,1,2 \ \ \&  \ \ r\geq 1 \ \ \& \ \ |s|\leq 1 \\
    \end{split}
\end{equation*}
where the implicit constants depend on $r,\ C_0,\ c_0,\ \alpha$. 
\end{lemma}

\begin{remark}[Loss of Good Derivatives]
    The additional loss of two derivatives in \ref{estgalbrunflow} when compared to Lemma \ref{estgalbrun} is there because we need to pass from the potential to the actual LDF and then control the derivative of the Lagrangian flow with the derivative of the LDF as in \eqref{derivativelagrangiang}.
\end{remark}

\begin{proof}[Proof of Lemma \ref{estgalbrunflow}] The estimates on the pure space derivatives are standard, see Proposition \ref{standardlagrangianestimate}, treating time as a parameter or adapt the argument we are about to give. 

\noindent We will deduce the pure time derivative bounds from those for $\xi^g$, which can be read immediately from the ones on $\Theta_j^g$ in Lemma \ref{estgalbrun}. For notational convenience, we assume $0\leq s \leq 1$; the bounds for $-1\leq s \leq 0$ follow identically. From the Lagrangian flow equation for $X^g_s,$ \eqref{lagrangiang}, we deduce:
\begin{equation*}
    \begin{split}
        \partial_t X^g_s(x,t)&=\partial_t\int_0^s\xi^g(X^g_s(x,t),t)\dd s\\
        &=\int_0^s\partial_t(\xi^g)(X^g_s(x,t),t)+\DD \xi^g(X^g_s(x,t),t)[\partial_t X^g_s]\dd s\\
    \end{split}
\end{equation*}
Using this formula and the composition estimates in Proposition \ref{compestimates}, we deduce:
\begin{equation*}
    \begin{split}
        ||\partial_t X^g_s||_r&\leq \int_0^s||\partial_t(\xi^g)\circ X^g_s)||_r+||\DD \xi^g\circ X^g_s||_r||\partial_t X^g_s||_0+||\DD \xi^g\circ X^g_s||_0||\partial_t X^g_s||_r\dd s\\
        &\overset{r=0}{\leq}\int_0^s||\partial_t\xi^g||_0+||\DD \xi^g||_0||\partial_t X^g_s||_0\dd s\\
        &\overset{r\geq 1}{\lesssim} \int_0^s||\partial_t\xi^g||_1||\DD X^g_s||_{r-1}+||\partial_t\xi^g||_r||\DD X^g_s||_{0}^r\dd s\\
        &+\int_0^s\left(||\DD \xi^g||_1||\DD X^g_s||_{r-1}+||\DD \xi^g||_r||\DD X^g_s||_0^r\right])||\partial_t X^g_s||_0+||\DD \xi^g||_0||\partial_t X^g_s||_r\dd s\\
    \end{split}
\end{equation*}
Given that $||\DD \xi^g||_0\lesssim \mathcal{T}^g\leq 1$ we can run a standard Gr{\"o}nwall argument. We first deduce the estimate for $r=0$, plug this one into the bound for $r\geq 1$ and conclude: 
$$||\partial_t X^g_s||_r\lesssim \lambda_q^{r}\lambda_q^{[r-\underline{r}]^+(b-1)\gamma_\ell}\mathcal{T}_g \ \text{ for }\ |s|\leq 1 \ \ \& \ \ r\geq 0.$$

\noindent We now run the same argument to deduce the second-order pure-time-derivative bound. We first compute:
\begin{equation*}
    \begin{split}
        \partial_t^2 X^g_s(x,t)&=\int_0^s\partial_t^2(\xi^g)(X^g_s(x,t),t)+2\DD (\partial_t\xi^g)(X^g_s(x,t),t)[\partial_t X^g_s]\dd s\\
        &+\int_0^s\DD^2\xi^g(X^g_s(x,t),t)[\partial_t X^g_s,\partial_t X^g_s]+\DD\xi^g(X^g_s(x,t),t)[\partial_t^2 X^g_s]\dd s,
    \end{split}
\end{equation*}
from this obtain
\begin{equation*}
    \begin{split}
        ||\partial_t^2 X^g_s||_r&\lesssim\int_0^s||\partial_t^2\xi^g\circ X^g_s||_r+||\DD (\partial_t\xi^g)\circ X^g_s||_r||\partial_t X^g_s||_0+||\DD (\partial_t\xi^g)\circ X^g_s||_0||\partial_t X^g_s||_r\dd s\\
        &+\int_0^s||\DD^2\xi^g\circ X^g_s||_r||\partial_t X^g_s||_0^2+||\DD^2\xi^g\circ X^g_s||_0||\partial_t X^g_s||_0||\partial_t X^g_s||_r\dd s\\
        &+\int_0^s||\DD\xi^g\circ X^g_s||_r||\partial_t^2 X^g_s||_0+||\DD\xi^g\circ X^g_s||_0||\partial_t^2 X^g_s||_r\dd s\\
        &\lesssim \lambda_q^{r+1}\lambda_q^{[r+2-\underline{r}]^+(b-1)\gamma_\ell}\mathcal{T}_g+\int_0^s||\DD\xi^g\circ X^g_s||_r||\partial_t^2 X^g_s||_0+||\DD\xi^g||_0||\partial_t^2 X^g_s||_r\dd s\\
    \end{split}
\end{equation*}
and conclude that:
$$||\partial_t^2 X^g_s||_r\lesssim \lambda_q^{r+1}\lambda_q^{[r+1-\underline{r}]^+(b-1)\gamma_\ell}\mathcal{T}_g \ \text{ for }\ |s|\leq 1 \ \ \& \ \ r\geq 0.$$
The implicit constants in the above bounds depend on $r, \ C_0, \ c_0, \ \alpha$.
\end{proof}

\subsubsection{Estimates on the Lie-Taylor Expansion}
We split the estimates between the mollified terms and the non-mollified terms. Out of convenience, we state bounds only for the potentials $\Theta_j^g$. The bounds for $\xi_j^g$ follow immediately by commuting the Lie-derivatives with the $\curl$ operator by means of the identity \eqref{DG1}. The key technical tool is the Inductive Lemma \ref{inductive} as in Remark \ref{slowinductive}, which allows us to deduce the estimates for an arbitrary term of the Lie-Taylor expansion from those of the leading term.

\begin{lemma}[Lie-Derivative and Transport Estimates - Mollified Terms]\label{estimatesperturbationgalbrun} Fix $\bar N$ a non-negative integer. Let $r,\ \sigma, \ k\geq 0$ integers, $j\in \mathbb{Z}$ and $\underline{r}=M-m_0-5$. There exist constants $C, C'$ which depend on $\bar N, \ C_0, \ c_0, \ \alpha $ but are uniform in the parameters $j, \ k, \ \sigma, \ r$, such that the following estimates hold.

\noindent\textbf{Lie-derivatives bounds.}  For $\sigma=0,1,2$ and $0\leq r+k+\sigma \leq \bar N$, we have:
\begin{equation*}
    \begin{split}
        &||\partial_t^\sigma\mathcal{L}_{\xi^{g}}^k(\partial_t+\mathcal{L}_{v_{\ell,j}})\Theta^{g}_j||_{r+\alpha}, \ ||\partial_t^\sigma\mathcal{L}_{\xi^{g}}^k\mathcal{L}_{B_{\ell,j}}\Theta^{g}_j||_{r+\alpha}\leq C'(C)^k \lambda_q^{r+\sigma}\lambda_q^{[r+k+\sigma-\underline{r}]^+(b-1)\gamma_\ell}\mathcal{T}_g^k\ell^{-\alpha}\tau^a \delta_{q+1},\\
        &||\partial_t^\sigma\mathcal{L}_{\xi^{g}}^k\mathcal{L}_{v_\ell-v_{\ell,j}}\Theta^{g}_j||_{r+\alpha}, \ ||\partial_t^\sigma\mathcal{L}_{\xi^{g}}^k\mathcal{L}_{B_\ell-B_{\ell,j}}\Theta^{g}_j||_{r+\alpha}\leq C'(C)^k\lambda_q^{r+\sigma}\lambda_q^{[r+k+\sigma-\underline{r}]^+(b-1)\gamma_\ell}\mathcal{T}_g^k\ell^{-\alpha}(\ell\lambda_q)^{m_0}\tau^a\delta_{q+1}.\\
    \end{split}
\end{equation*}

\noindent \textbf{Alfv\'en-transport bounds.} For $0\leq r+k \leq \bar N-1$, we have:
\begin{equation*}
    \begin{split}
        &||\mathcal{A}^\pm_{\ell,j}\mathcal{L}_{\xi^{g}}^k(\partial_t+\mathcal{L}_{v_{\ell,j}})\Theta^{g}_j||_{r}, \ ||\mathcal{A}^\pm_{\ell,j}\mathcal{L}_{\xi^{g}}^k\mathcal{L}_{B_{\ell,j}}\Theta^{g}_j||_{r}\leq C'(C)^k \lambda_q^{r}\lambda_q^{[r+k-\underline{r}]^+(b-1)\gamma_\ell}\mathcal{T}_g^k\ell^{-\alpha}\delta_{q+1},\\
        &||\mathcal{A}^\pm_{\ell,j}\mathcal{L}_{\xi^{g}}^k\mathcal{L}_{v_\ell-v_{\ell,j}}\Theta^{g}_j||_{r}, \ ||\mathcal{A}^\pm_{\ell,j}\mathcal{L}_{\xi^{g}}^k\mathcal{L}_{B_\ell-B_{\ell,j}}\Theta^{g}_j||_{r}\leq C'(C)^k\lambda_q^{r}\lambda_q^{[r+k-\underline{r}]^+(b-1)\gamma_\ell}\mathcal{T}^k_g(\ell\lambda_q)^{m_0}\ell^{-\alpha}\lambda_q\delta_q^{1/2}\tau^a\delta_{q+1}.\\
    \end{split}
\end{equation*}
\end{lemma}

\begin{proof}[Proof of Lemma \ref{estimatesperturbationgalbrun}] We want to apply the Inductive Lemma \ref{inductive} as in Remark \ref{slowinductive}, with the following $F_0^j$:
$$(\partial_t+\mathcal{L}_{v_{\ell,j}})\Theta^{g}_j, \ \mathcal{L}_{B_{\ell,j}}\Theta^{g}_j, \ \mathcal{L}_{v_\ell-v_{\ell,j}}\Theta^{g}_j, \ \mathcal{L}_{B_\ell-B_{\ell,j}}\Theta^{g}_j.$$
We let: 
$$\mathcal{J}=\left\{J=(j_1,\dots,j_k):j_i\sim j\right\}$$ 
be the set of all ordered $k$-tuples of indices $j_i\in \{j,\ j+1,\ j-1\}$ neighbouring $j$, note that the cardinality of this set is $2^k$, which is in line with the implicit constant scaling exponentially in $k$. Since the support of the other $\xi^g_{j_i}$ do not intersect that of $\xi_j^g$ by construction, we can write:
$$\mathcal{L}_{\xi^g}^kF_0^j=\sum_{J\in \mathcal{J}}\mathcal{L}_{\xi^{g}_{j_1}}\dots \mathcal{L}_{\xi^{g}_{j_k}}F_0^j$$
In the notation of the Inductive Lemma, we are thus interested in sequences of vector fields $(\sigma_i)_i:=(\xi_{j_i}^p)_i$ with $j_i\sim j$. We now need to provide estimates on $F_0^j$ and $\xi_{j_i}^g$, we remark that we are not trying to be precise with the exact loss $\lambda_q\leadsto \ell^{-1}$ quantified by $\lambda_q^{[r-\underline{r}]^+(b-1)\gamma_\ell}$, we will pick the worst lower bound, namely $\underline{r}=M-m_0-5$ coming from the mollification correction terms, between all the estimates even if a specific term admits a better bound.

\noindent \textit{$F_0^j$ of transport type.} We first compute:
\begin{equation*}
    \begin{split}
        (\partial_t+\mathcal{L}_{v_{\ell,j}})\Theta^{g}_j&= \frac{1}{2}(\mathcal{A}_{\ell,j}^++\mathcal{A}_{\ell,j}^-)\Theta^{g}_j+(\DD v_{\ell,j})^\top\Theta^{g}_j,\\
        \partial_t[(\partial_t+\mathcal{L}_{v_{\ell,j}})\Theta^{g}_j]&= \frac{1}{2}\partial_t(\mathcal{A}_{\ell,j}^++\mathcal{A}_{\ell,j}^-)\Theta^{g}_j+(\DD\partial_tv_{\ell,j})^\top\Theta^{g}_j+(\DD v_{\ell,j})^\top\partial_t\Theta^{g}_j,\\
        \partial_t^2[(\partial_t+\mathcal{L}_{v_{\ell,j}})\Theta^{g}_j]&= \frac{1}{2}\partial_t^2(\mathcal{A}_{\ell,j}^++\mathcal{A}_{\ell,j}^-)\Theta^{g}_j+(\DD\partial_t^2v_{\ell,j})^\top\Theta^{g}_j+(\DD v_{\ell,j})^\top\partial_t^2\Theta^{g}_j+2(\DD\partial_tv_{\ell,j})^\top\partial_t\Theta^{g}_j\\
    \end{split}
\end{equation*}
and by means of the bounds in Lemmas \ref{localcorrg} and \ref{estgalbrun} we deduce:
$$||\partial_t^\sigma(\partial_t+\mathcal{L}_{v_{\ell,j}})\Theta^{g}_j||_r\lesssim \lambda_q^{r+\sigma}\lambda_q^{[r+\sigma-(\underline{r}+1)](b-1)\gamma_\ell}\ell^{-\alpha}\tau^a\delta_{q+1} \ \text{ for }\ \sigma=0,1,2 \ \ \& \ \ r\geq 0.$$

\noindent We move to the transport estimates. We first compute:
\begin{equation*}
    \begin{split}
        (\partial_t+\mathcal{L}_{z^\pm_{\ell,j}})(\partial_t+\mathcal{L}_{v_{\ell,j}})\Theta^{g}_j&=\frac{1}{2}(\partial_t+\mathcal{L}_{z^\pm_{\ell,j}})(\mathcal{A}_{\ell,j}^++\mathcal{A}_{\ell,j}^-)\Theta^{g}_j+(\partial_t+\mathcal{L}_{z^\pm_{\ell,j}})\DD v_{\ell,j}^\top\Theta^{g}_j\\
        &=\frac{1}{2}\mathcal{A}_{\ell,j}^\pm(\mathcal{A}_{\ell,j}^++\mathcal{A}_{\ell,j}^-)\Theta^{g}_j+(\DD z^{\pm}_{\ell,j})^\top (\mathcal{A}_{\ell,j}^++\mathcal{A}_{\ell,j}^-)\Theta^{g}_j\\
        &+\DD (\mathcal{A}^\pm_{\ell,j}v_{\ell,j})^\top\Theta^{g}_j+\DD (v_{\ell,j})^\top\mathcal{A}^\pm_{\ell,j}\Theta^{g}_j
    \end{split}
\end{equation*}
and by means of Lemmas \ref{classicalgalbrun}, \ref{stabilitygalbrun} and \ref{estgalbrun} we deduce:
\begin{equation*}
    \begin{split}
        &||(\partial_t+\mathcal{L}_{z^\pm_{\ell,j}})(\partial_t+\mathcal{L}_{v_{\ell,j}})\Theta^{g}_j||_{r+\alpha}\\
        &\leq||\mathcal{A}_{\ell,j}^\pm(\mathcal{A}_{\ell,j}^++\mathcal{A}_{\ell,j}^-)\Theta^{g}_j||_{r+\alpha}+||(\DD z^{\pm}_{\ell,j})^\top (\mathcal{A}_{\ell,j}^++\mathcal{A}_{\ell,j}^-)\Theta^{g}_j||_{r+\alpha}\\
        &+||\DD (\mathcal{A}^\pm_{\ell,j}v_{\ell,j})^\top\Theta^{g}_j||_{r+\alpha}+||\DD (v_{\ell,j})^\top\mathcal{A}^\pm_{\ell,j}\Theta^{g}_j||_{r+\alpha}\\
        &\leq\underbrace{||\mathcal{A}_{\ell,j}^\pm(\mathcal{A}_{\ell,j}^++\mathcal{A}_{\ell,j}^-)\Theta^{g}_j||_{r+\alpha}}_{\lambda_q^{r}\ell^{-\alpha}\delta_{q+1}}+\underbrace{||z^{\pm}_{\ell,j}||_1||(\mathcal{A}_{\ell,j}^++\mathcal{A}_{\ell,j}^-)\Theta^{g}_j||_{r+\alpha}+||z^{\pm}_{\ell,j}||_{r+1+\alpha}||(\mathcal{A}_{\ell,j}^++\mathcal{A}_{\ell,j}^-)\Theta^{g}_j||_0}_{\lambda_q^{r+1}\delta_q^{1/2}\ell^{-2\alpha}\tau^a\delta_{q+1}}\\
        &+\underbrace{||\mathcal{A}^\pm_{\ell,j}v_{\ell,j}||_{1}||\Theta^{g}_j||_{r+\alpha}+||\mathcal{A}^\pm_{\ell,j}v_{\ell,j}||_{r+1+\alpha}||\Theta^{g}_j||_{0}}_{\lambda_q^{r+2}\ell^{-2\alpha}\delta_q\tau^c\tau^a\delta_{q+1}}+\underbrace{||v_{\ell,j}||_1||\mathcal{A}^\pm_{\ell,j}\Theta^{g}_j||_{r+\alpha}+||v_{\ell,j}||_{r+1+\alpha}||\mathcal{A}^\pm_{\ell,j}\Theta^{g}_j||_{0}}_{\lambda_q^{r+1}\delta_q^{1/2}\ell^{-2\alpha}\tau^a\delta_{q+1}}\\
        &\lesssim \lambda_q^{r}\lambda_q^{[r-(\underline{r}+1)](b-1)\gamma_\ell}\ell^{-\alpha}\delta_{q+1}\\
    \end{split}
\end{equation*}
Here we used that: $$\tau^c\ell^{-\alpha}\lambda_q\delta_q^{1/2}\leq 1$$
and will do so in what follows, without further mention. Similar computation for $\mathcal{L}_{B_{\ell,j}}\Theta^g_j$ lead to:
\begin{equation*}
    \begin{split}
        ||\partial_t^{\sigma}\mathcal{L}_{B_{\ell,j}}\Theta^g_j||_{r+\alpha}&\lesssim \lambda_q^{r+\sigma}\lambda_q^{[r+\sigma-(\underline{r}+1)](b-1)\gamma_\ell}\ell^{-\alpha}\tau^a\delta_{q+1} \ \text{ for } \ \sigma=0,1,2 \ \ \& \ \ r\geq 0,\\
        ||(\partial_t+\mathcal{L}_{z^\pm_{\ell,j}})\mathcal{L}_{B_{\ell,j}}\Theta^g_j||_{r+\alpha}&\lesssim \lambda_q^{r}\lambda_q^{[r-(\underline{r}+1)](b-1)\gamma_\ell}\ell^{-\alpha}\delta_{q+1}\  \text{ for }\ r\geq 0.
    \end{split}
\end{equation*}

\noindent \textit{$F_0^j$ coming from the mollification correction.} We first compute:
\begin{equation*}
    \begin{split}
        \mathcal{L}_{v_\ell-v_{\ell,j}}\Theta^{g}_j&= (v_\ell-v_{\ell,j})\cn \Theta^{g}_j+\DD(v_\ell-v_{\ell,j})^\top\Theta^{g}_j,\\
        \partial_t[\mathcal{L}_{v_\ell-v_{\ell,j}}\Theta^{g}_j]&= \partial_t(v_\ell-v_{\ell,j})\cn \Theta^{g}_j+(v_\ell-v_{\ell,j})\cn \partial_t\Theta^{g}_j+\DD\partial_t(v_\ell-v_{\ell,j})^\top\Theta^{g}_j+\DD(v_\ell-v_{\ell,j})^\top\partial_t\Theta^{g}_j,\\
        \partial_t^2[\mathcal{L}_{v_\ell-v_{\ell,j}}\Theta^{g}_j]&= \partial_t^2(v_\ell-v_{\ell,j})\cn \Theta^{g}_j+(v_\ell-v_{\ell,j})\cn \partial_t^2\Theta^{g}_j+2\partial_t(v_\ell-v_{\ell,j})\cn \partial_t\Theta^{g}_j\\
        &+\DD\partial_t^2(v_\ell-v_{\ell,j})^\top\Theta^{g}_j+\DD(v_\ell-v_{\ell,j})^\top\partial_t^2\Theta^{g}_j+2\DD\partial_t(v_\ell-v_{\ell,j})^\top\partial_t\Theta^{g}_j,\\
    \end{split}
\end{equation*}
and by means of the bounds in Lemmas \ref{stabilitygalbrun}, \ref{estgalbrun} we deduce:
$$||\partial_t^\sigma\mathcal{L}_{v_\ell-v_{\ell,j}}\Theta^{g}_j||_r\lesssim \lambda_q^{r+\sigma}\lambda_q^{[r+\sigma-\underline{r}](b-1)\gamma_\ell}\ell^{-\alpha}\tau^a\delta_{q+1} \ \text{ for }\ \sigma=0,1,2 \ \text{ and }\ r\geq 0.$$

\noindent We now use the identity in \eqref{trick} to write:
$$(\partial_t+\mathcal{L}_{z^\pm_{\ell,j}})\mathcal{L}_{v_\ell-v_{\ell,j}}\Theta^{g}_j=\mathcal{L}_{(\partial_t+\mathcal{L}_{z^\pm_{\ell,j}})(v_\ell-v_{\ell,j})}\Theta^{g}_j+\mathcal{L}_{v_\ell-v_{\ell,j}}(\partial_t+\mathcal{L}_{z^\pm_{\ell,j}})\Theta^{g}_j.$$
Lemmas \ref{classicalgalbrun} and \ref{stabilitygalbrun} now give: 
\begin{equation*}
    \begin{split}
        ||(\partial_t+\mathcal{L}_{z^\pm_{\ell,j}})(v_\ell-v_{\ell,j})||_{r+\alpha}&\lesssim ||\mathcal{A}^\pm_{\ell,j}(v_\ell-v_{\ell,j})||_{r+\alpha}+||(v_\ell-v_{\ell,j})\cn z^\pm_{\ell,j}||_{r+\alpha}\\
        &\lesssim||\mathcal{A}^\pm_{\ell,j}(v_\ell-v_{\ell,j})||_{r+\alpha}+||v_\ell-v_{\ell,j}||_0||z^\pm_{\ell,j}||_{r+1+\alpha}+||v_\ell-v_{\ell,j}||_{r+\alpha}||z^\pm_{\ell,j}||_1\\
        &\lesssim \lambda_q^{r+1}\lambda_q^{[r-\underline{r}]^+(b-1)\gamma_\ell}\ell^{-\alpha}(\ell\lambda_q)^{m_0}\delta_q,
        \\
    \end{split}
\end{equation*}
using in addition the bounds in Lemma \ref{estgalbrun}, we deduce:
\begin{equation*}
    \begin{split}
        &||(\partial_t+\mathcal{L}_{z^\pm_{\ell,j}})\mathcal{L}_{v_\ell-v_{\ell,j}}\Theta^{g}_j||_{r+\alpha}\\
        &\leq||\mathcal{L}_{(\partial_t+\mathcal{L}_{z^\pm_{\ell,j}})(v_\ell-v_{\ell,j})}\Theta^{g}_j||_{r+\alpha}+||\mathcal{L}_{v_\ell-v_{\ell,j}}(\partial_t+\mathcal{L}_{z^\pm_{\ell,j}})\Theta^{g}_j||_{r+\alpha}\\
        &\leq||[(\partial_t+\mathcal{L}_{z^\pm_{\ell,j}})(v_\ell-v_{\ell,j})]\cn\Theta^{g}_j||_{r+\alpha}+||\DD[(\partial_t+\mathcal{L}_{z^\pm_{\ell,j}})(v_\ell-v_{\ell,j})]^\top\Theta^{g}_j||_{r+\alpha}\\
        &+||(v_\ell-v_{\ell,j})\cn(\partial_t+\mathcal{L}_{z^\pm_{\ell,j}})\Theta^{g}_j||_{r+\alpha}+||\DD(v_\ell-v_{\ell,j})^\top(\partial_t+\mathcal{L}_{z^\pm_{\ell,j}})\Theta^{g}_j||_{r+\alpha}\\
        &\lesssim \lambda_q^{r+2}\lambda_q^{[r-\underline{r}]^+(b-1)\gamma_\ell}\ell^{-2\alpha}(\ell\lambda_q)^{m_0}\tau^a\tau^c\delta_q\delta_{q+1}+\lambda_q^{r+1}\lambda_q^{[r-\underline{r}]^+(b-1)\gamma_\ell}(\ell\lambda_q)^{m_0}\ell^{-2\alpha}\tau^a\delta_q^{1/2}\delta_{q+1}\\
        &\lesssim\lambda_q^{r+1}\lambda_q^{[r-\underline{r}]^+(b-1)\gamma_\ell}(\ell\lambda_q)^{m_0}\ell^{-\alpha}\tau^a\delta_q^{1/2}\delta_{q+1}.
    \end{split}
\end{equation*}
Similarly, one can show:
$$||(\partial_t+\mathcal{L}_{z^\pm_{\ell,j}})\mathcal{L}_{B_\ell-B_{\ell,j}}\Theta^{g}_j||_{r+\alpha}\lesssim \lambda_q^{r+1}\lambda_q^{[r-\underline{r}]^+(b-1)\gamma_\ell}(\ell\lambda_q)^{m_0}\ell^{-\alpha}\tau^a\delta_q^{1/2}\delta_{q+1}.$$

\noindent \textit{Estimates on the $\sigma_i$.} As the estimates on $\xi_{j_i}^p$ can be deduced from \ref{estgalbrun} upon using the identity \eqref{DG1}. We need only deal with the error in the transport, namely 
$$(\partial_t+\mathcal{L}_{z^\pm_{\ell,j}})\xi^{g}_{j_i}=(\partial_t+\mathcal{L}_{z^\pm_{\ell,j_i}})\xi^{g}_{j_i}+\mathcal{L}_{z^\pm_{\ell,j}-z^{\pm}_{\ell,j_i}}\xi^g_{j_i},$$
by means of Lemmas \ref{stabilitygalbrun}, \ref{estgalbrun} we first deduce:
\begin{equation}\label{transporthighliegalbrun3}
    \begin{split}
        ||(\partial_t+\mathcal{L}_{z^\pm_{\ell,j_i}})\xi^g_{j_i}||_{r+\alpha}&\lesssim\lambda_q^{r+1}\lambda_q^{[r-\underline{r}]^+(b-1)\gamma_\ell}\ell^{-\alpha}\tau^a\delta_{q+1},
        \\
        ||\mathcal{L}_{z^\pm_{\ell,j}-z^{\pm}_{\ell,j_i}}\xi^g_{j_i}||_{r+\alpha}&\leq||(z^\pm_{\ell,j}-z^\pm_{\ell,j_i})\cn\xi^{g}_j||_{r+\alpha}+||\xi^{g}_j\cn (z^\pm_{\ell,j}-z^\pm_{\ell,j_i})||_{r+\alpha}\\
        &\lesssim(||z^\pm_{\ell,j}-z_\ell^\pm||_0+||z^\pm_{\ell,j_i}-z_\ell^\pm||_0)||\xi^{g}_j||_{r+1+\alpha}\\
        &+(||z^\pm_{\ell,j}-z_\ell^\pm||_{r+\alpha}+||z^\pm_{\ell,j_i}-z_\ell^\pm||_{r+\alpha})||\xi^{g}_j||_{1}\\
        &+||\xi^{g}_j||_0(||z^\pm_{\ell,j}-z_\ell^\pm||_{r+1+\alpha}+||z^\pm_{\ell,j_i}-z_\ell^\pm||_{r+1})\\
        &+||\xi^{g}_j||_{r+\alpha}(||z^\pm_{\ell,j}-z_\ell^\pm||_1+||z^\pm_{\ell,j_i}-z_\ell^\pm||_1)\\
        &\lesssim \lambda_q^{r+2}\lambda_q^{[r-\underline{r}]^+(b-1)\gamma_\ell}\ell^{-2\alpha}(\ell\lambda_q)^{m_0}\tau^c\tau^a\delta_q^{1/2}\delta_{q+1}\\
        &\lesssim\lambda_q^{r+1}\lambda_q^{[r-\underline{r}]^+(b-1)\gamma_\ell}(\ell\lambda_q)^{m_0}\ell^{-\alpha}\tau^a\delta_{q+1}\\
    \end{split}
\end{equation}
and we conclude from $(\ell\lambda_q)^{m_0}<1$ that:
$$||(\partial_t+\mathcal{L}_{z^\pm_{\ell,j}})\xi^{g}_{j_i}||_{r+\alpha}\lesssim \lambda_q^{r+1}\lambda_q^{[r-\underline{r}]^+(b-1)\gamma_\ell}\ell^{-\alpha}\tau^a\delta_{q+1} \ \text{ for } \ r\geq 0.$$

\noindent \textit{Conclusion.} With the estimates on the various $F_0^j$ and $\sigma_i=\xi_{j_i}$ at hand we can apply the Inductive Lemma \ref{inductive} as in Remark \ref{slowinductive}, with:
\begin{equation*}
        \begin{split}
            &A=\ell^{-\alpha}\tau^a\delta_{q+1} \ \text{ or } \ell^{-\alpha}(\ell\lambda_q)^{m_0}\tau^a\delta_{q+1},\\ 
            &A_{\mathcal{A}}=\ell^{-\alpha}\delta_{q+1} \ \text{ or } \ell^{-\alpha}(\ell\lambda_q)^{m_0}\ell^{-\alpha}\lambda_q\delta_q^{1/2}\tau^a\delta_{q+1},\\
            &\bar\varsigma_i= \lambda_q\ell^{-\alpha}\tau^a\tau^c\delta_{q+1}, \ \bar\varsigma_{i,\mathcal{A}}=\lambda_q\ell^{-\alpha}\tau^a\delta_{q+1},\\
            &L(r)=L_{\mathcal{A}}(r)=\lambda_q^{[r-\underline{r}]^+(b-1)\gamma_\ell}.\\
        \end{split}
    \end{equation*}
Note that $L=L_\mathcal{A}$, we do not lose a derivative because the objects were constructed by solving \eqref{galbrunapplied} exactly. We conclude that:
\begin{equation*}
    \begin{split}
        ||\partial_t^{\sigma}\mathcal{L}_{\xi^{g}_{j_1}}\dots \mathcal{L}_{\xi^{g}_{j_k}}(\partial_t+\mathcal{L}_{v_{\ell,j}})\Theta^{g}_j||_{r+\alpha}&\leq C' (C)^k \lambda_q^{r}\lambda_q^{[r+\sigma+k-\underline{r}]^+(b-1)\gamma_\ell}\ell^{-\alpha}\mathcal{T}^k_g\tau^a\delta_{q+1} \\ 
        &\text{ for } \sigma=0,1,2 \ \text{ and } 0\leq r+k +\sigma\leq \bar N\\
        ||\mathcal{A}^\pm_{\ell,j}\mathcal{L}_{\xi^{g}_{j_1}}\dots \mathcal{L}_{\xi^{g}_{j_k}}(\partial_t+\mathcal{L}_{v_{\ell,j}})\Theta^{g}_j||_{r+\alpha}&\leq C' (C)^k\lambda_q^{r}\lambda_q^{[r-\underline{r}]^+(b-1)\gamma_\ell}\ell^{-\alpha}\mathcal{T}^k_g\delta_{q+1} \\
        &\text{ for } \ 0\leq r+k \leq \bar N-1\\ 
    \end{split}
\end{equation*}
and
\begin{equation*}
    \begin{split}
        ||\partial_t^\sigma\mathcal{L}_{\xi^{g}_{j_1}}\dots \mathcal{L}_{\xi^{g}_{j_k}}\mathcal{L}_{v_\ell-v_{\ell,j}}\Theta^{g}_j||_{r+\alpha}&\leq C' (C)^k \lambda_q^{r}\lambda_q^{[r+k-\underline{r}]^+(b-1)\gamma_\ell}\ell^{-\alpha}\mathcal{T}^k_g(\ell\lambda_q)^{m_0}\tau^a \delta_{q+1} \\ 
        &\text{ for } \sigma=0,1,2 \ \text{ and } 0\leq r +k+\sigma\leq \bar N,\\
        ||\mathcal{A}^\pm_{\ell,j}\mathcal{L}_{\xi^{g}_{j_1}}\dots \mathcal{L}_{\xi^{g}_{j_k}}\mathcal{L}_{v_\ell-v_{\ell,j}}\Theta^{g}_j||_{r+\alpha}&\leq C' (C)^k\lambda_q^{r}\lambda_q^{[r+k-\underline{r}]^+(b-1)\gamma_\ell}\ell^{-\alpha}\mathcal{T}^k_g(\ell\lambda_q)^{m_0} \lambda_q\delta_q^{1/2}\tau^a\delta_{q+1}\\ 
        &\text{ for } 0\leq r+k \leq \bar N-1.\\ 
    \end{split}
\end{equation*}
The implicit constants and $C$ depend on $\bar N,\ C_0,\ c_0, \ \alpha$ but not on $r,\ k$. Similar bounds hold for the terms involving the magnetic field.

\noindent The claimed estimates then follow by summing over $J=(j_i,\dots,j_k)\in \mathcal{J}$ and redefining $C, \ C'$ accordingly.
\end{proof}


\begin{lemma}[Remainders - Galbrun Stage]\label{remaindersg} For $k_0\geq 0$ let: 
\begin{equation*}
    \begin{split}
        \theta^{k_0}&=
        \begin{cases}
            \sum_j\frac{(-1)^{k_0+1}}{(k_0+1)!}\int_0^{1}(X_{s}^g)_*\left[\mathcal{L}_{\xi^g}^{k_0+1}(\partial_t+\mathcal{L}_{ v_{\ell,j}})\Theta^g_j\right](1-s)^{k_0+1}\dd s, \\
            \sum_j\frac{(-1)^{k_0+1}}{(k_0+1)!}\int_0^{1}(X_{s}^g)_*\left[\mathcal{L}_{\xi^g}^{k_0+1}\mathcal{L}_{ B_{\ell,j}}\Theta^g_j\right](1-s)^{k_0+1}\dd s
        \end{cases}
        \\
        \theta^{k_0}_c&=
        \begin{cases}
            \sum_j\frac{(-1)^{k_0+1}}{(k_0+1)!}\int_0^{1}(X_{s}^g)_*\left[\mathcal{L}_{\xi^g}^{k_0+1}\mathcal{L}_{B_\ell- B_{\ell,j}}\Theta_j^g\right](1-s)^{k_0+1}\dd s,\\
            \sum_j\frac{(-1)^{k_0+1}}{(k_0+1)!}\int_0^{1}(X_{s}^g)_*\left[\mathcal{L}_{\xi^g}^{k_0+1}\mathcal{L}_{v_\ell- v_{\ell,j}}\Theta_j^g\right](1-s)^{k_0+1}\dd s.
        \end{cases}
    \end{split}
\end{equation*}
Fix $\bar N$ a non-negative integer and $\underline{r}=M-m_0-5$. There exist constants $C, \ C'$, which depend on $\bar N, \ C_0, \ c_0, \ \alpha $ such that for $r\geq 0$ integer and $\sigma=0,1,2$ satisfying $0\leq r +(k_0+1)+\sigma\leq \bar N$ we have:
\begin{equation*}
    \begin{split}
        &||\partial_t^\sigma \theta^{k_0}||_{r+\alpha} \leq C' \frac{(C)^{k_0+1}}{(k_0+2)!}\lambda_q^{r+\sigma}\lambda_q^{[r+\sigma+(k_0+1)-\underline{r}]^+(b-1)\gamma_\ell}\mathcal{T}_g^{k_0+1}\ell^{-\alpha}\tau^a\delta_{q+1},\\
        &||\partial_t^\sigma\theta^{k_0}_c||_{r+\alpha}\leq C' \frac{(C)^{k_0+1}}{(k_0+2)!}\lambda_q^{r+\sigma}\lambda_q^{[r+\sigma+(k_0+1)-\underline{r}]^+(b-1)\gamma_\ell}\mathcal{T}_g^{k_0+1}(\ell\lambda_q)^{m_0}\ell^{-\alpha}\tau^a\delta_{q+1}.\\
    \end{split}
\end{equation*}
\end{lemma}
\noindent The proof of this Lemma is a variation of the proof of Lemma \ref{remaindersmoll} using the estimates in Lemmas \ref{localcorrg} and \ref{estgalbrun} this time, and we omit it.
\begin{remark}From Lemma \ref{remaindersg} with $k_0=k_0^g, \ \bar N=N$, we deduce that:
\begin{equation*}
    ||\partial_t^\sigma \theta^g_w||_{r+\alpha}, \ ||\partial_t^\sigma \theta^g_b||_{r+\alpha} \leq C' \frac{(C)^{k_0^g+1}}{(k_0^g+2)!}\lambda_q^{r+\sigma}\lambda_q^{[r+\sigma+(k_0^g+1)-\underline{r}]^+(b-1)\gamma_\ell}\mathcal{T}_g^{k_0^g+1}\ell^{-\alpha}\tau^a\delta_{q+1}
\end{equation*}
for $0\leq r +(k_0^g+1)+\sigma\leq N$, where $\theta_w^g, \ \theta_b^g$ are as in \eqref{decompgalbrunw} and \eqref{decompgalbrunb}.
\end{remark}


\begin{lemma}[Lie-Derivative and Transport Estimates - Non Mollified Terms]\label{estimatesperturbationgalbrun2} Let $j\in \mathbb{Z}$ and $r, \ k, \ \sigma\geq 0$ integers. Moreover, let $L, \ L_{\mathcal{A}}:\mathbb{N}_{\geq 0}\to \mathbb{R}_{\geq 1}$ be the loss functions: 
\begin{equation}\label{LAgalbrun}
    \begin{split}
        &L(r)=\lambda_q^{[r-\underline{r}](b-1)\gamma_\ell},\\
        &L_{\mathcal{A}}(r)=1_{r\leq \underline{r}-1}+1_{r\geq \underline{r}} \lambda_q^{[r-(\underline{r}-1)](b-1)\gamma_\ell}\bar L
    \end{split}
\end{equation}
with $\underline{r}=M-m_0-5, \  \bar{L}=\delta_q^{-1/2}$. There exist constants $C, \ C'$ which depend on all the parameters but not on $a$ and are uniform in $k, \ r, \ j, \ \sigma$, such that the following bounds hold.

\noindent\textbf{Lie-derivatives bounds.}  For $\sigma=0,1,2$ and $0\leq r+k+\sigma\leq N-m_0-2$, we have:
\begin{equation*}
    ||\partial_t^\sigma\mathcal{L}_{\xi^{g}}^k\mathcal{L}_{v_q-v_\ell}\Theta^{g}_j||_{r+\alpha}, \ ||\partial_t^\sigma\mathcal{L}_{\xi^{g}}^k\mathcal{L}_{B_q-B_\ell}\Theta^{g}_j||_{r+\alpha}\leq C'(C)^k\lambda_q^{r+\sigma}L(r+k+\sigma)\mathcal{T}_g^k\ell^{-\alpha}(\ell\lambda_q)^{m_0}\tau^a\delta_{q+1}.
\end{equation*}

\noindent \textbf{Alfv\'en-transport bounds.} For $0\leq r+k\leq N-m_0-3$, we have:
\begin{equation*}
    ||\mathcal{A}^\pm_{\ell,j}\mathcal{L}_{\xi^{g}}^k\mathcal{L}_{v_q-v_\ell}\Theta^{g}_j||_{r+\alpha}, \ ||\mathcal{A}^\pm_{\ell,j}\mathcal{L}_{\xi^{g}}^k\mathcal{L}_{B_q-B_\ell}\Theta^{g}_j||_{r+\alpha}\leq C'(C)^k\lambda_q^{r}L_{\mathcal{A}}(r+k)\mathcal{T}^k_g(\ell\lambda_q)^{m_0}\ell^{-\alpha}\lambda_q\delta_q^{1/2}\tau^a\delta_{q+1}.
\end{equation*}

\noindent \textbf{Remainder bounds.} For $\sigma=0,1,2$ and $0\leq r+\sigma\leq N-(k_0^g+1)-m_0-2$, we have:
\begin{equation*}
\begin{split}
    ||\partial_t^\sigma\mathring \theta^g_w||_{r+\alpha}, \ ||\partial_t^\sigma\mathring \theta^g_b||_{r+\alpha}\leq\frac{C'(C)^{k_0^g+1}}{(k_0^g+2)!} \lambda_{q}^r\mathcal{T}_g^{k_0^p+1}L(k_0^g+1+\sigma)(\ell\lambda_q)^{m_0}\ell^{-\alpha}\tau^a\delta_{q+1}.\\
\end{split}
\end{equation*}
\end{lemma}

The proof of this Lemma is mutatis mutandis as that of Lemma \ref{estimatesperturbationgalbrun}. We remark, however, the following key differences:
\begin{itemize}
    \item We cannot pick $\bar N$ arbitrarily large, the presence of $v_q-v_\ell,\ B_q-B_\ell$ forces us to set $\bar N=N-m_0-2$, the loss of $m_0$ coming from the mollification error estimate, and the additional two are one to interpolate the $C^{\alpha}$ norm and the other because even for $k=0$ we have one Lie-derivative acting. This is also why $\underline{r}$ is smaller compared to the one in Lemma \ref{estgalbrun}.
    \item We state the bound for the transport derivative this way out of convenience, the loss $L_{\mathcal{A}}$ just saying that if $r+k> \underline{r}-1$ we don't have information on the transport derivative of $v_q-v_\ell, \ B_q-B_\ell$ and we switch to the pure time derivative estimate losing a $\delta_q^{1/2}$. To match the loss function of the pure derivative bound, we also shift $\underline{r}$ to $\underline{r}-1$, compare this to the Alfv\'en transport bound in Lemma \ref{estimatesperturbationgalbrun}.
\end{itemize}


\subsubsection{Estimates on the Vector Fields.} With all the estimates for the intermediate objects making the Lie-Taylor expansions of the Galbrun perturbation, we are ready to prove estimates on the actual vector fields.

\begin{lemma}[Estimates on $(\mathring w^g,\mathring b^g)$]\label{estimatesfieldsringg} Let $r\geq 0$ integer and $\underline{r}=M-m_0-6, \ \bar L=1/\delta_q^{1/2}$. Let $L_{g,\mathcal{A}}:\mathbb{N}_{\geq 0}\to \mathbb{R}_{\geq 1}$ be one of the following admissible loss functions:
\begin{subnumcases}{L_{g,\mathcal{A}}(r)=}
    1_{r\leq \underline{r}-k_0^g-2}+1_{\underline{r}-k_0^g-1\leq r\leq \underline{r}-1} \left[1+\left(\lambda_q^{(b-1)\gamma_\ell}\mathcal{T}_g\right)^{\underline{r}-1-r}\bar L\right]+1_{ r\geq \underline{r}}\lambda_q^{[r-(\underline{r}-1)](b-1)\gamma_\ell}\bar L \label{precise}
    \\
    1_{r\leq \underline{r}-k_0^g-1}+1_{r\geq \underline{r}-k_0^g}\left(\frac{\lambda_{q+1}}{\lambda_{q}}\right)^{r-(\underline{r}-k_0^g-1)}=\lambda_q^{[r-(\underline{r}-k_0^g-1)]^+(b-1)}\label{lossy}
\end{subnumcases}
Under the choice of parameters in \ref{choiceofparameters}, we have: 
\begin{equation}\label{boundlossfunction}
    \eqref{precise} \leq 2\eqref{lossy}
\end{equation}
and 
\begin{equation*}
    \begin{split}
        &||\partial_t^j\mathring w^g||_{r},\ ||\partial_t^j \mathring b^g||_{r}\lesssim \lambda_q^{r+j}\lambda_q^{[r+j-\underline{r}]^+(b-1)\gamma_\ell}(\ell\lambda_q)^{m_0}\mathcal{T}_g\delta_q^{1/2} \ \text{ for } \ j=0,1,2 \ \text{ and } \ 0\leq r\leq N-j,\\
        &||\mathcal{A}^\pm \mathring w^g||_{r}, \ ||\mathcal{A}^\pm \mathring b^g||_{r}\lesssim \lambda_q^{r}L_{g,\mathcal{A}}(r)(\ell\lambda_q)^{m_0}\mathcal{T}_g1/\tau^c\delta_q^{1/2} \ \text{ for }  0\leq r\leq N-1.
    \end{split}
\end{equation*}
The implicit constants depend on $r$ and all the other parameters, but not on $a$.
\end{lemma}
\begin{remark}(Loss Function Heuristics) The loss function $L_{g,\mathcal{A}}$ keeps track of three intervals of derivatives. In the first one, we have the expected estimate. In the second one, we lose the good transport bound, $\bar L$ appears, but we retain some smallness from the Lie-Taylor expansion, quantified by $\mathcal{T}_g$. In the third one, we have no smallness left, and the mollification parameter also starts to appear. The claim in \eqref{boundlossfunction} tells us that, upon paying the ratio $\lambda_{q+1}/\lambda_q$, we can remove the loss, and still ensure an estimate compatible with our Iterative Assumption \eqref{inductiveassumptionsgeneral}.   
\end{remark}

\begin{proof}[Proof of Lemma \ref{estimatesfieldsringg}] We will only provide bounds for $\mathring w^g$, those for $\mathring b^g$ follow similarly.

\noindent \textbf{Pure derivative bounds.} Using the definition of $X^g$ in \eqref{lagrangiang} and $\mathring w^g$ in \eqref{decompgalbrunb}, and the property of pullbacks in \eqref{lietransport}, we can write:
\begin{equation}\label{directexpansion}
    \begin{split}
        \mathring w^g&=X^g_*(v_q-v_\ell)-(v_q-v_\ell)\\
        &=(X^g_1)_*(v_q-v_\ell)-(X^g_0)_*(v_q-v_\ell)\\
        &=\int_0^1\partial_s[(X_s^g)_*(v_q-v_\ell)]\dd s\\
        &=\int_0^1\partial_s[(X_{-s}^g)^*(v_q-v_\ell)]\dd s\\
        &=-\int_0^1(X_{-s}^g)^*\mathcal{L}_{\xi^g}(v_q-v_\ell)\dd s\\
        &=\int_0^1(X_{-s}^g)^*\mathcal{L}_{v_q-v_\ell}\xi^g\dd s\\
        &=\int_0^1(X_{-s}^g)^*\mathcal{L}_{v_q-v_\ell}\bigg[\curl\sum_j\Theta^g_j\bigg]\dd s\\
        &=\curl\sum_j\int_0^1(X_s^g)_*\left[\mathcal{L}_{v_q-v_\ell}\Theta^g_j\right]\dd s
    \end{split}
\end{equation}
here we also used that $(X^g_s)^{-1}=X^g_{-s}$ as $\xi^g$ does not depend on $s$ and the identity in \eqref{DG0}. From the bounds in Lemmas \ref{estgalbrun}, \ref{standardmollgalbrun} and \eqref{Tg} it follows that:
\begin{equation*}
    \begin{split}
        ||\partial_t^j\mathcal{L}_{v_q-v_\ell}\Theta^g_j||_r&\lesssim \lambda_q^{r+j}\lambda_q^{[r+j-(\underline{r}+1)]^+(b-1)\gamma_\ell}(\ell\lambda_q)^{m_0}\lambda_q\ell^{-\alpha}\delta_q^{1/2}\tau^a\tau^c\delta_{q+1}\\
        &=\lambda_q^{r+j}\lambda_q^{[r+j-(\underline{r}+1)]^+(b-1)\gamma_\ell}(\ell\lambda_q)^{m_0}\frac{1}{\lambda_q}\mathcal{T}_g\delta_q^{1/2}
    \end{split}
\end{equation*}
for $j=0,1,2$ and $0\leq r\leq N-m_0-j-1$. Estimating the integral as in the proof of Lemma \ref{remaindersg}, we deduce:
\begin{equation}\label{wcirclowderivative}
        ||\partial_t^j\mathring w^g||_{r}\lesssim \lambda_q^{r+j}\lambda_q^{[r+j-\underline{r}]^+(b-1)\gamma_\ell}(\ell\lambda_q)^{m_0}\mathcal{T}_g\delta_q^{1/2} \ \text{ for }\ j=0,1,2 \ \text{ and }\ 0\leq r \leq N-m_0-j-2
\end{equation}
where the implicit constants depend on $r$ and all the other parameters, but not on $a$. Note that the need to still expand up to the first order gives the loss of two additional good derivatives compared to the ones we have on the potentials $\Theta_j^g$, from $M-m_0-4\leadsto M-m_0-6$, which matches the loss we have on $\DD X^g-\IId$ from Lemma \ref{estgalbrunflow}.

\noindent For $N-m_0-j-2\leq r \leq N$, we cannot use the expansion, and we estimate the terms separately. We will need that $(X^g)^{-1}=X^g_{-1}$ the estimates in Lemmas \ref{estgalbrunflow} and \ref{standardmollgalbrun} and the composition estimates in Proposition \ref{compestimates}.
\begin{equation}\label{largepure}
    \begin{split}
        ||\mathring w^g||_r&=||X^g_*(v_q-v_\ell)||_r+||v_q-v_\ell||_r\\
        &\lesssim ||\DD X^g \circ (X^g)^{-1}||_r||(v_q-v_\ell)\circ (X^g)^{-1}||_0+||\DD X^g \circ (X^g)^{-1}||_0||(v_q-v_\ell)\circ (X^g)^{-1}||_r\\
        &+||v_q-v_\ell||_r\\
        &\lesssim ||v_q-v_\ell||_0\left[||\DD X^g ||_{1}||\DD X^g_{-1}||_{r-1}+||\DD X^g ||_{r}||\DD X_{-1}^g||_0^{r}\right]\\
        &+||\DD X^g||_0\left[||v_q-v_\ell||_1||\DD X^g_{-1}||_{r-1}+||v_q-v_\ell||_{r}||\DD X^g_{-1}||_0^r\right]\\
        &+||v_q||+||v_\ell||_r\\
        &\lesssim\lambda_q^r\lambda_q^{[r-\underline{r}]^+(b-1)\gamma_\ell}\mathcal{T}_g(\ell\lambda_q)^{m_0}\delta_q^{1/2}+\lambda_q^r\delta_q^{1/2}.\\
    \end{split}
\end{equation}

\noindent Similarly, for high space derivatives of the time derivative, we cannot use the Lie-Taylor expansion. Set $M=\DD X^g\circ(X^g)^{-1}$. We begin with the following calculations:
\begin{equation}\label{firsttimederivate}
    \begin{split}
       \partial_t \left[X^g_*(v_q-v_\ell)\right]&=\partial_tM(v_q-v_\ell) \circ(X^g)^{-1}\\
       &+M\left[\partial_t(v_q-v_\ell) \circ(X^g)^{-1}       +\DD(v_q-v_\ell) \circ(X^g)^{-1}\partial_t (X^g)^{-1}\right],
       \\
       \partial_t M&=\DD \partial_t X^g\circ(X^g)^{-1}+\DD^2X^g \circ(X^g)^{-1}[\partial_t(X^g)^{-1}, \cdot ]
    \end{split}
\end{equation}
and 
\begin{equation}\label{secondtimederivative}
    \begin{split}
        \partial_t^2 \left[X^g_*(v_q-v_\ell)\right]&=\partial_t^2M(v_q-v_\ell) \circ(X^g)^{-1}\\
        &+\partial_tM\left[\partial_t(v_q-v_\ell) \circ(X^g)^{-1}+\DD(v_q-v_\ell) \circ(X^g)^{-1}\partial_t (X^g)^{-1}\right]\\
       &+M\left[\partial_t^2(v_q-v_\ell) \circ(X^g)^{-1}+2\DD[\partial_t(v_q-v_\ell)] \circ(X^g)^{-1}\partial_t (X^g)^{-1}\right]\\
       &+M\left[\DD^2[v_q-v_\ell] \circ(X^g)^{-1}[\partial_t (X^g)^{-1},\partial_t (X^g)^{-1}]+\DD(v_q-v_\ell) \circ(X^g)^{-1}\partial_t^2 (X^g)^{-1}\right],
       \\
       \partial_t^2 M&=\DD \partial_t^2 X^g\circ(X^g)^{-1}+2\DD^2\partial_t X^g \circ(X^g)^{-1}[\partial_t(X^g)^{-1}, \cdot ]\\
       &+\DD^2X^g \circ(X^g)^{-1}[\partial_t^2(X^g)^{-1}, \cdot ]+\DD^3X^g \circ(X^g)^{-1}[\partial_t(X^g)^{-1}, \ \partial_t(X^g)^{-1}, \cdot ].
    \end{split}
\end{equation}
Arguing as in \eqref{largepure}, we deduce:
\begin{equation}\label{puretimering}
\begin{split}
    ||\partial_t^j\mathring w^g||_{r}&\leq ||\partial_t^j \left[X^g_*(v_q-v_\ell)\right]||_{r}+||\partial_t^j(v_q-v_\ell)||_{r}\\
    &\lesssim \lambda_q^{r+j}\lambda_q^{[r+j-\underline{r}]^+(b-1)\gamma_\ell}\mathcal{T}_g(\ell\lambda_q)^{m_0}\delta_q^{1/2} +\lambda_q^{r+j}\delta_q^{1/2}
\end{split}
\end{equation}
for $j=0,1,2$ and $N-m_0-j-2\leq r \leq N-j$. 

\noindent Now, given our choice of parameters in \ref{choiceofparameters}, see \eqref{misc} and \eqref{Tg}, we have:
$$(N-m_0-2-j-\underline{r})\gamma_\ell-(1+2\gamma_{CZ}+\gamma_a+2\beta)\geq 0 \Longrightarrow \lambda_q^{r+j}\delta_q^{1/2}\leq \lambda_q^{r+j}\lambda_q^{[r+j-\underline{r}]^+(b-1)\gamma_\ell}(\ell\lambda_q)^{m_0}\mathcal{T}_g\delta_q^{1/2}$$
and from \eqref{wcirclowderivative}, \eqref{puretimering} we conclude:
\begin{equation*}
        ||\partial_t^j\mathring w^g||_{r}\lesssim \lambda_q^{r+j}\lambda_q^{[r+j-\underline{r}]^+(b-1)\gamma_\ell}(\ell\lambda_q)^{m_0}\ell^{-\alpha}\lambda_q\tau^a\delta_{q+1} \ \text{ for }\ j=0,1,2 \ \text{ and }\ 0\leq r \leq N-j
\end{equation*}
where the implicit constant depends on $r$ and all the parameters, but not on $a$. 

\noindent \textbf{Alfv\'en transport bounds.} Recall the full Lie-Taylor expansion given in \eqref{decompgalbrunw}, namely 
$$ \mathring w^g=\curl\left[\sum_j\sum_{k=0}^{k_0^g}\frac{(-1)^k}{(k+1)!}\mathcal{L}_{\xi^g}^k\mathcal{L}_{v_q- v_\ell}\Theta_j^g+\mathring\theta^g_w\right],$$
we will use without explicit mention that thanks to the commutation \eqref{DG1}, given any 1-form $F$, the following identity holds:
\begin{equation}\label{identityfieldsg}
    \begin{split}
        \mathcal{A}^\pm_{\ell,j}\curl F&=(\partial_t+\mathcal{L}_{z^\pm_{\ell,j}})\curl F+\curl F\cn z^\pm_{\ell,j}\\
        &=\curl (\partial_t+\mathcal{L}^1_{z^\pm_{\ell,j}})F+\curl F\cn z^\pm_{\ell,j}\\
        &=\curl[\mathcal{A}^\pm_{\ell,j}F]+\curl[(\DD z^\pm_{\ell,j})^\top F]+\curl [F]\cn z^\pm_{\ell,j}
    \end{split}
\end{equation}
and we can thus deduce the required transport bounds from those for the potentials $\Theta^g_j$ in Lemma \ref{estimatesperturbationgalbrun2} and the ones in Lemmas \ref{standardmollgalbrun}, \ref{stabilitygalbrun} to correct the transport operator. We omit the calculations and state the final bound:
$$||\mathcal{A}^\pm\curl\mathcal{L}_{\xi^g}^k\mathcal{L}_{v_q- v_\ell}\Theta_j^g||_{r}\lesssim (C)^k\lambda_q^{r}L_{\mathcal{A}}(r+k+1)\mathcal{T}_g^k(\ell\lambda_q)^{m_0}\mathcal{T}_g(1/\tau^c)\delta_q^{1/2}$$
for $0\leq r +k+m_0+3\leq N$.

\noindent We now deal with the remainder. We are forced to trade its smallness for a poor pure time-derivative estimate, since we don't control its geometric properties. Given that our choice of parameters in \ref{choiceofparameters}, see \eqref{constraintadmissibility}, in particular ensures 
$$(\lambda_q^{(b-1)\gamma_\ell}\mathcal{T}_g)^{k_0^g}\leq \delta_q^{1/2},$$ 
using the bounds in Lemma \ref{estimatesperturbationgalbrun2} and in the Iterative Assumptions \ref{inductiveassumptionsgeneral}, we deduce:
\begin{equation*}
    \begin{split}
        ||\mathcal{A}^\pm\curl \ \mathring \theta^{g}_w||_r&\lesssim ||\partial_t\mathring \theta^{g}_w||_{r+1}+||z^\pm_q||_0||\mathring \theta^{g}_w||_{r+2}+||z^\pm_q||_r||\mathring \theta^{g}_w||_2\\
        &\lesssim \frac{(C)^{k_0^g+1}}{(k_0^g+2)!}\lambda_q^{r+2}\lambda_q^{[r+2+k_0^g+1-(\underline{r}+1)](b-1)\gamma_\ell}\mathcal{T}^{k_0^g+1}_g(\ell\lambda_q)^{m_0}\ell^{-\alpha}\tau^a\delta_{q+1}\\
        &\lesssim \lambda_q^{r}\lambda_q^{[r+1-\underline{r}](b-1)\gamma_\ell}(\ell\lambda_q)^{m_0}\ell^{-\alpha}\tau^a\lambda_q^2\delta_q^{1/2}\delta_{q+1}\\
        &=\lambda_q^{r}\lambda_q^{[r-(\underline{r}-1)]^+(b-1)\gamma_\ell}(\ell\lambda_q)^{m_0}\mathcal{T}_g(1/\tau^c)\delta_q^{1/2}
    \end{split}
\end{equation*}
for $0\leq r +(k_0^g+1)+m_0+3\leq N$. Note the difference in the definitions of $\underline{r}$.

\noindent Lemma \ref{estimatesperturbationgalbrun2} also gives estimates for the explicit part of the expansion. We will split the sum in the parameter $k$ according to the number of derivatives $r$ and the definition of the loss function $L_\mathcal{A}$ in \eqref{LAgalbrun} to isolate the bad regime where $\bar L=\delta_q^{-1/2}$ appears.  Using in addition the fact that at most two $\Theta_j^g$ are non-zero, we deduce:
\begin{equation*}
\begin{split}
    ||\mathcal{A}^\pm \mathring w^g||_{r} &\leq \sup_j\sum_{k=0}^{k_0^g}\frac{1}{(k+1)!}||\mathcal{A}^\pm\curl\mathcal{L}_{\xi^g}^k\mathcal{L}_{v_q- v_\ell}\Theta_j^g||_{r}+||\mathcal{A}^\pm\curl \ \mathring \theta ^g_w||_{r}\\
    &\lesssim\sum_{k=0}^{k_0^g}\frac{(C)^k}{(k+1)!}\lambda_q^{r}L_{\mathcal{A}}(r+k+1)\mathcal{T}_g^k(\ell\lambda_q)^{m_0}\mathcal{T}_g(1/\tau^c)\delta_q^{1/2}+\lambda_q^{r}\lambda_q^{[r-(\underline{r}-1)]^+(b-1)\gamma_\ell}(\ell\lambda_q)^{m_0}\mathcal{T}_g(1/\tau^c)\delta_q^{1/2}\\
    &= \lambda_q^{r}(\ell\lambda_q)^{m_0}\mathcal{T}_g(1/\tau^c)\delta_q^{1/2}\left[\sum_{k+r\leq \underline{r}-2,0\leq k\leq k_0^g}\frac{\left(C\mathcal{T}_g\right)^k}{(k+1)!}\right]\\
    &+ \lambda_q^{r}(\ell\lambda_q)^{m_0}\mathcal{T}_g(1/\tau^c)\delta_q^{1/2}\left[\bar L\sum_{k+r\geq \underline{r}-1,0\leq k\leq k_0^g}\frac{\left(C\lambda_q^{(b-1)\gamma_\ell}\mathcal{T}_g\right)^k}{(k+1)!}\lambda_q^{[r-(\underline{r}-1)]^+(b-1)\gamma_\ell}\right]\\
    &+\lambda_q^{r}\lambda_q^{[r-(\underline{r}-1)]^+(b-1)\gamma_\ell}(\ell\lambda_q)^{m_0}\mathcal{T}_g(1/\tau^c)\delta_q^{1/2}\\
    &= \lambda_q^{r}(\ell\lambda_q)^{m_0}\mathcal{T}_g(1/\tau^c)\delta_q^{1/2}\Biggr[\underbrace{\sum_{k=0}^{[\underline{r}-1-r]^+-1}\frac{\left(C\mathcal{T}_g\right)^k}{(k+1)!}}_{\lesssim 1}+ \bar L\underbrace{\sum_{k=[\underline{r}-1-r]^+}^{ k_0^g}\frac{\left(C\lambda_q^{(b-1)\gamma_\ell}\mathcal{T}_g\right)^k}{(k+1)!}}_{\lesssim \left(\lambda_q^{(b-1)\gamma_\ell}\mathcal{T}_g\right)^{[\underline{r}-1-r]^+}}\lambda_q^{[r-(\underline{r}-1)]^+(b-1)\gamma_\ell}\Biggr]\\
    &+\lambda_q^{r}\lambda_q^{[r-(\underline{r}-1)]^+(b-1)\gamma_\ell}(\ell\lambda_q)^{m_0}\mathcal{T}_g(1/\tau^c)\delta_q^{1/2}\\
    &\lesssim \lambda_q^r(\ell\lambda_q)^{m_0}\mathcal{T}_g(1/\tau^c)\delta_q^{1/2}\underbrace{\left[1_{[\underline{r}-1-r]^+\geq 1 }+1_{[\underline{r}-1-r]^+\leq k_0^g} \lambda_q^{[r-(\underline{r}-1)]^+(b-1)\gamma_\ell}\left(\lambda_q^{(b-1)\gamma_\ell}\mathcal{T}_g\right)^{[\underline{r}-1-r]^+}\bar L\right]}_{=f(r)}\\
    &+\lambda_q^{r}\lambda_q^{[r-(\underline{r}-1)]^+(b-1)\gamma_\ell}(\ell\lambda_q)^{m_0}\mathcal{T}_g(1/\tau^c)\delta_q^{1/2}\\
\end{split}
\end{equation*}
for $0\leq r \leq N-(k_0^g+1)-m_0-3$ and note the different definition of $\underline{r}$ compared to Lemma \ref{estimatesperturbationgalbrun2}. The implicit constant depends on $r$ and all the parameters, but not on $a$. Note, however, that for $N-(k_0^g+1)-m_0-3\leq r\leq N-1$ we have: 
$$f(r)=\lambda_q^{[r-(\underline{r}-1)](b-1)\gamma_\ell}\bar L$$
which matches up to a $\lambda_{q+1}^\alpha$ loss the pure time derivative estimates \eqref{puretimering} above. 

\noindent We now bound $f$. For $r$ sufficiently large, each new derivative produces a loss of $\lambda_q^{(b-1)\gamma_\ell}\mathcal{T}_g$, we will compensate this by trading a good derivative $\lambda_q$, for a bad derivative $\lambda_{q+1}$. Having $\underline{r}-1$ good derivatives available, with $\underline{r}$ sufficiently large, will save the day. First note that:
$$1_{[\underline{r}-1-r]^+\geq 1}=1_{r\leq \underline{r}-2} \ \text{ and } \ 1_{[\underline{r}-1-r]^+\leq k_0^g}=1_{r\geq \underline{r}-1-k_0^g},$$
we can then rewrite:
\begin{equation}\label{f(r)}
    \begin{split}
        f(r)&=1_{[\underline{r}-1-r]^+\geq 1}+1_{[\underline{r}-1-r]^++\leq k_0^g} \lambda_q^{[r-(\underline{r}-1)]^+(b-1)\gamma_\ell}\left(\lambda_q^{(b-1)\gamma_\ell}\mathcal{T}_g\right)^{[\underline{r}-1-r]^+}\bar L\\
        &=1_{r\leq \underline{r}-2}+1_{r\geq \underline{r}-k_0^g-1} \lambda_q^{[r-(\underline{r}-1)]^+(b-1)\gamma_\ell}\left(\lambda_q^{(b-1)\gamma_\ell}\mathcal{T}_g\right)^{[\underline{r}-1-r]^+}\bar L\\
        &\overset{(a)}{=}1_{r\leq \underline{r}-2}+\underbrace{1_{\underline{r}-k_0^g-1\leq r\leq \underline{r}-1} \left(\lambda_q^{(b-1)\gamma_\ell}\mathcal{T}_g\right)^{\underline{r}-1-r}\bar L}_{f_1(r)}+\underbrace{1_{ r\geq \underline{r}}\lambda_q^{[r-(\underline{r}-1)](b-1)\gamma_\ell}\bar L}_{f_2(r)},\\
    \end{split}
\end{equation}
from the choice of parameters \ref{choiceofparameters}, see \eqref{constraintadmissibility} and \eqref{M}, we have: 
$$\left(\lambda_q^{(b-1)\gamma_\ell}\mathcal{T}_g\right)^{k_0^g}\bar L\leq 1, \ 1-(2\beta +\gamma_a+\gamma_\ell+2\gamma_{CZ})\geq 0, \ \underline{r}-k_0^g\geq 0$$
using in addition the definition of $\mathcal{T}_g$ in \eqref{Tg} we can bound:
\begin{equation}\label{f1}
    \begin{split}
        f_1(r)&=1_{\underline{r}-k_0^g-1\leq r\leq \underline{r}-1}\left(\lambda_q^{(b-1)\gamma_\ell}\mathcal{T}_g\right)^{k_0^g}\bar L\left(\lambda_q^{(b-1)\gamma_\ell}\mathcal{T}_g\right)^{\underline{r}-1-r-k_0^g}\\
        &=1_{\underline{r}-k_0^g-1\leq r\leq \underline{r}-1}\underbrace{\left(\lambda_q^{(b-1)\gamma_\ell}\mathcal{T}_g\right)^{k_0^g}\bar L}_{\leq 1}\left(\lambda_q^{(b-1)\gamma_\ell}\mathcal{T}_g\right)^{-(r-(\underline{r}-1-k_0^g))}\\
        &\leq 1_{\underline{r}-1-k_0^g\leq r\leq \underline{r}-1}\left(\frac{\lambda_q}{\lambda_{q+1}}\right)^{\underbrace{[1-(2\beta+\gamma_a+\gamma_\ell+2\gamma_{CZ})]}_{\geq 0}(r-(\underline{r}-1-k_0^g))}\left(\frac{\lambda_{q+1}}{\lambda_{q}}\right)^{r-(\underline{r}-1-k_0^g)}\\
        &\leq 1_{\underline{r}-1-k_0^g\leq r\leq \underline{r}-1}\left(\frac{\lambda_{q+1}}{\lambda_{q}}\right)^{r-(\underline{r}-1-k_0^g)},
        \\
        f_2(r)&=1_{ r\geq \underline{r}}\lambda_q^{[r-(\underline{r}-1)](b-1)\gamma_\ell}\underbrace{\bar L\left(\frac{\lambda_{q}}{\lambda_{q+1}}\right)^{k_0^g}}_{\leq 1} \left(\frac{\lambda_{q+1}}{\lambda_{q}}\right)^{k_0^g}\\
        &\leq 1_{ r\geq \underline{r}}\left(\frac{\lambda_{q+1}}{\lambda_{q}}\right)^{r-(\underline{r}-k_0^g-1)}
    \end{split}
\end{equation}
and collecting also the more precise estimate $(a)$ in \eqref{f(r)} the inequality in \eqref{boundlossfunction} follows, indeed
\begin{equation*}
\begin{split}
    f(r)&\leq1_{r\leq \underline{r}-k_0^g-2}+1_{\underline{r}-k_0^g-1\leq r\leq \underline{r}-1} \left[1+\left(\lambda_q^{(b-1)\gamma_\ell}\mathcal{T}_g\right)^{\underline{r}-1-r}\bar L\right]+1_{ r\geq \underline{r}}\lambda_q^{[r-(\underline{r}-1)](b-1)\gamma_\ell}\bar L\\
    &\leq 2\left[1_{r\leq \underline{r}-k_0^g-1}+1_{r\geq \underline{r}-k_0^g}\left(\frac{\lambda_{q+1}}{\lambda_{q}}\right)^{r-(\underline{r}-k_0^g-1)}\right].\\
\end{split}
\end{equation*}
In particular, given the definition of $L_{g,\mathcal{A}}$ in the statement of the Lemma, and the remark about extending the estimate up to $r\leq N-1$ made above, we conclude that:
$$||\mathcal{A}^\pm \mathring w^g||_{r}\lesssim \lambda_q^{r}L_{g,\mathcal{A}}(r)(\ell\lambda_q)^{m_0}\mathcal{T}_g(1/\tau^c)\delta_q^{1/2} \ \text{ for } \ 0\leq r\leq N-1$$
where the implicit constant depends on $r$ and all the parameters, but not on $a$.
\end{proof}


\begin{lemma}[Estimates on $(\bar w^g, \ \bar b^g)$]\label{estimatesfieldsbarg} Let $r\geq0$ integer and set $\underline{r}=M-m_0-6$. We have:
\begin{equation*}
    \begin{split}
        &||\partial_t^j\bar w^g||_{r},\ ||\partial_t^j \bar b^g||_{r}\lesssim \lambda_q^{r+j}\lambda_q^{[r-(\underline{r}-j)]^+(b-1)\gamma_\ell}\ell^{-\alpha}\tau^a\lambda_q\delta_{q+1} \ \text{ for } \ j=0,1,2 \ \text{ and } \ 0\leq r\leq N-j,\\
        &||\mathcal{A}^\pm \bar w^g||_{r}, \ ||\mathcal{A}^\pm \bar b^g||_{r}\lesssim \lambda_q^{r}\lambda_q^{[r-(\underline{r}-1)]^+(b-1)\gamma_\ell}\ell^{-\alpha}\lambda_q\delta_{q+1} \ \text{ for }  0\leq r\leq N-m_0.
    \end{split}
\end{equation*}
    The implicit constants depend on $r$ and all the other parameters, but not on $a$.
\end{lemma}
\begin{remark}The maximum number of derivatives on the Alfv\'en transport bounds comes from the fact that: 
$$\mathcal{A}^\pm=\partial_t+z_q^\pm\cn$$ 
where $z_q^\pm$ is not mollified, and we have to locally correct it to $\mathcal{A}^\pm_{\ell,j}=\partial_t+z_{\ell,j}^\pm\cn$ to match the transport on each component of the perturbation. This uses the mollification-error and stability estimates in Lemmas \ref{standardmollgalbrun}, \ref{stabilitygalbrun} which lose $m_0$ derivatives.
\end{remark}

\begin{proof}[Proof of Lemma \ref{estimatesfieldsbarg}.] We will only prove bounds for $\bar w^g$, those for $\bar b^g$ follow mutatis mutandis.

\noindent \textbf{Preliminaries.} Recall the decomposition \eqref{decompgalbrunw}, namely
\begin{equation}\label{expansionrecallg}
    \bar w^g=\curl\left[\sum_{k=0}^{k_0^g}\frac{(-1)^k}{(k+1)!}\mathcal{L}_{\xi^g}^k(\partial_t+\mathcal{L}_{v_{\ell,j}})\Theta^g_j+\sum_{k=0}^{k_0^g}\frac{(-1)^k}{(k+1)!}\mathcal{L}_{\xi^g}^k\mathcal{L}_{v_{\ell}- v_{\ell,j}}\Theta_j^g+\theta^g_w\right],
\end{equation}
we now let $\bar N = N+k_0^g+2$ and apply Lemmas \ref{estimatesperturbationgalbrun}, \ref{remaindersg}, to get the estimates for the various terms in the expansion. This suffices as we can use the commutation \eqref{identityfieldsg} as in the proof of Lemma \ref{estimatesfieldsringg}. In what follows, we use the notation from \ref{estimatesperturbationgalbrun}, \ref{remaindersg}. 

\noindent \textbf{Pure derivative bounds.} From the fact that we have at most two non-zero $\Theta_j^g$ at the same time, we deduce:
\begin{equation*}
    \begin{split}
        ||\partial_t^j\bar w^g||_{r}&\lesssim \sum_{k=0}^{k_0^g}\frac{1}{(k+1)!}\sup_j\left[||\partial_t^j\mathcal{L}_{\xi^g}^k(\partial_t+\mathcal{L}_{v_{\ell,j}})\Theta^g_j||_{r+1}+||\partial_t^j\mathcal{L}_{\xi^g}^k\mathcal{L}_{v_{\ell}- v_{\ell,j}}\Theta_j^g||_{r+1}\right]\\
        &+||\partial_t^j \theta^{g}_w||_{r+1}\\
        &\lesssim \lambda_q^{r+1+j}\lambda_q^{[r+j-\underline{r}]^+(b-1)\gamma_\ell}\ell^{-\alpha}\tau^a\delta_{q+1}\left[1+(\ell\lambda_q)^{m_0}\right]\underbrace{\sum_{k=0}^{k_0^g}\frac{\left(C\lambda_q^{(b-1)\gamma_\ell}\mathcal{T}_g\right)^k}{(k+1)!}}_{\lesssim 1}\\
        &+\lambda_q^{r+1+j}\lambda_q^{[r+j+k_0^g+1-\underline{r}]^+(b-1)\gamma_\ell}\mathcal{T}_g^{k_0^g+1}\ell^{-\alpha}\tau^a \delta_{q+1}\\
        &\lesssim \lambda_q^{r+1+j}\lambda_q^{[r-(\underline{r}-j)]^+(b-1)\gamma_\ell}\ell^{-\alpha}\tau^a\delta_{q+1}
    \end{split}
\end{equation*}
for $0\leq r \leq \bar N-k_0^g-1-j=N-j$ and the implicit constant depends on $r, \ \bar N$ and all the other parameters, but not on $a$. Here we used $(\lambda_q^{(b-1)\gamma_\ell}\mathcal{T}_p)^{k_0^g+1}<1$.

\noindent  \textbf{Alfv\'en transport bounds.} Lemma \ref{estimatesperturbationgalbrun} also proves Alfv\'en transport bounds for the explicit part of the Lie-Taylor expansion. Note, however, that we are forced to use pure first-order time derivatives for the remainder. Our choice of parameters in \ref{choiceofparameters}, see \eqref{constraintadmissibility}, guarantees: $$\left(\lambda_q^{(b-1)\gamma_\ell}\mathcal{T}_p\right)^{k_0^g}\leq \delta_q^{1/2}$$ 
and from the bounds in Lemma \ref{remaindersg} and in the Iterative Assumptions \ref{inductiveassumptionsgeneral}, we deduce:
\begin{equation*}
    \begin{split}
        ||\mathcal{A}^\pm \curl \ \theta^{g}_w||_r&\lesssim ||\partial_t \theta ^{g}_w||_{r+1}+||z^\pm_q||_r||\theta^{g}_w||_2+||z^\pm_q||_0||\theta^{g}_w||_{r+2}\\
        &\lesssim \frac{(C)^{k_0^g+1}}{(k_0^g+2)!}\lambda_q^{r+2}\lambda_q^{[r-(\underline{r}-k_0^g-2)]^+(b-1)\gamma_\ell}\mathcal{T}_g^{k_0^g+1}\ell^{-\alpha}\tau^a \delta_{q+1}\\
        &\leq \lambda_q^{r+1}\lambda_q^{[r-(\underline{r}-1)]^+(b-1)\gamma_\ell}\ell^{-\alpha}\tau^a\lambda_q\delta_q^{1/2} \delta_{q+1}\\
        &\leq \lambda_q^{r+1}\lambda_q^{[r-(\underline{r}-1)]^+(b-1)\gamma_\ell}\ell^{-\alpha}\delta_{q+1}
    \end{split}
\end{equation*}
for $0\leq r \leq \bar N-(k_0^g+1)-2=N-1$, note the different definitions of $\underline{r}$. We will use this in the following computations. We now argue as in the corresponding estimate in the proof of Lemma \ref{estimatesfieldsringg}. We take the transport derivative of \eqref{expansionrecallg}, use the commutation \eqref{identityfieldsg}, and correct the transport operator, we obtain:
\begin{equation*}
    \begin{split}
        &||\mathcal{A}^\pm\bar w^g||_{r}\\
        &\lesssim\sup_j\sum_{k=0}^{k_0^g}\frac{1}{(k+1)!}\Biggr[\underbrace{||\mathcal{A}^\pm_{\ell,j}\curl \ \mathcal{L}_{\xi^g}^k(\partial_t+\mathcal{L}_{v_{\ell,j}})\Theta^g_j||_{r}}_{\lambda_q^{r+1}\ell^{-\alpha}\delta_{q+1}}+\underbrace{||\mathcal{A}^\pm_{\ell,j}\curl \ \mathcal{L}_{\xi^g}^k\mathcal{L}_{v_{\ell}- v_{\ell,j}}\Theta_j^g||_{r}}_{\lambda_q^{r+1}(\ell\lambda_q)^{m_0}(\tau^a/\tau^c)\delta_{q+1}}\Biggr]\\
        &+\sup_j\sum_{k=0}^{k_0^g}\frac{1}{(k+1)!}\Biggr[\underbrace{||(z^\pm_q-z^\pm_{\ell,j})\cn\curl \ \mathcal{L}_{\xi^g}^k(\partial_t+\mathcal{L}_{v_{\ell,j}})\Theta^g_j||_{r}}_{\lambda_q^{r+1}(\ell\lambda_q)^{m_0}(\tau^a/\tau^c)\delta_{q+1}}+\underbrace{||(z^\pm_q-z^\pm_{\ell,j})\cn\curl \ \mathcal{L}_{\xi^g}^k\mathcal{L}_{v_{\ell}- v_{\ell,j}}\Theta_j^g||_{r}}_{\lambda_q^{r+1}(\ell\lambda_q)^{2m_0}(\tau^a/\tau^c)\delta_{q+1}}\Biggr]\\
        &+||\mathcal{A}^\pm \curl \ \theta^{g}_w||_{r}\\
        &\lesssim \lambda_q^{r+1}\lambda_q^{[r-(\underline{r}-1)]^+(b-1)\gamma_\ell}\ell^{-\alpha}\delta_{q+1}\left[1+(\tau^a/\tau^c)(\ell\lambda_q)^{m_0}+(\tau^a/\tau^c)(\ell\lambda_q)^{2m_0}\right]\underbrace{\sum_{k=0}^{k_0^g}\frac{\left(C\lambda_q^{(b-1)\gamma_\ell}\mathcal{T}_g\right)^k}{(k+1)!}}_{\lesssim 1}\\
        &+\lambda_q^{r+1}\lambda_q^{[r-(\underline{r}-1)]^+(b-1)\gamma_\ell}\ell^{-\alpha}\delta_{q+1}\\
        &\lesssim \lambda_q^{r+1}\lambda_q^{[r-(\underline{r}-1)]^+(b-1)\gamma_\ell}\ell^{-\alpha}\delta_{q+1}
    \end{split}
\end{equation*}
for $0\leq r \leq \min\{N-1, \ N-m_0\}=N-m_0$, the maximum number of derivatives coming from the mollification error estimates, given that $m_0>1$, see \eqref{misc1}. The implicit constant depends on $r, \ \bar N$ and all the parameters, but not on $a$. 
\end{proof}


\subsection{Reynolds Stress: \texorpdfstring{$R^g$}{Rg}}\label{Rgs}
We first split
\begin{equation}\label{Rg}
    R^g=R^{g,lin}+R^{g,qua}
\end{equation}
where $R^{g,lin}$ is required to satisfy
\begin{equation}\label{Rglin1}
    \begin{split}
        \ddiv \ R^{g,lin}&=\partial_tw^{g}+v_q\cn w^{g}-B_q\cn b^{g}+w^{g}\cn v_q-b^{g}\cn B_q+\nabla \pi^{gr}+\ddiv \left[\sum_I (g_I^2-1) A_I\right]\\
        &+\ddiv [R_q-R_\ell]\\
    \end{split}
\end{equation}
and 
$$R^{g, qua}=w^{g}\otimes w^{g}-b^{g}\otimes b^{g}$$
this ensure that \eqref{Rgp} is satisfied. The precise definition of $R^{g,lin}$ will be given after further rewriting. This denotes the errors which are linear with respect to the perturbation $(w^g,b^g)$ of $(v_q, \ B_q)$, although part of it is due to the non-linear way we adopt to construct the perturbation, that is, Lemma \ref{lpl}, while $R^{g,qua}$ comes from the quadratic errors. 

\noindent \textbf{Decomposition of Linear Errors.} We now split the linear error terms into simpler components. We recall that in Section \ref{gnldf} we constructed for each element $j$ of the time partition, the principal part of the Galbrun perturbation $(w_j^{g,p}, \ b_j^{g,p},\ \pi_j^{g})$ together with their cut-off error $R_j^{cut}$ so that they solve:
\begin{equation*}
    \begin{split}
        \partial_tw^{g,p}_j+v_{\ell,j}\cn w^{g,p}_j+w^{g,p}_j\cn v_{\ell,j}-B_{\ell,j}\cn b^{g,p}_j-b^{g,p}_j\cn B_{\ell,j}+\nabla \pi^{gr}_j&=\ddiv \left[\sum_{I:I_t=j} (1-g_I^2) A_I+R^{cut}_j\right].\\
    \end{split}
\end{equation*}
Similarly as in \eqref{decompgalbrunw}, \eqref{decompgalbrunb}, we decompose $(w^{g},b^{g})$ in:
\begin{equation*}
    \begin{split}
         w^g&=\sum_j w_j^{g,p}+w_j^{h}+r_w\\
         &=\sum_j w_j^{g,p}+\curl[\Theta_{w,j}^h]+\curl[\vartheta_w],\\
    \end{split}
    \ \ 
    \begin{split}
        b^g&=\sum_j b_j^{g,p}+b_j^h+r_b\\
        &=\sum_j b_j^{g,p}+\curl[\Theta_{b,j}^h]+\curl[\vartheta_b].\\
    \end{split}
\end{equation*}
Lemma \ref{lpl} gives explicit expressions:
\begin{equation*}
    \begin{cases}
        \Theta_{w,j}^h=\sum_{k=1}^{k_0^g}\frac{(-1)^k}{(k+1)!}\mathcal{L}_{\xi^{g}}^k(\partial_t+\mathcal{L}_{v_{\ell,j}})\Theta_j^{g}+\sum_{k=0}^{k_0^g}\frac{(-1)^k}{(k+1)!}\mathcal{L}_{\xi^{g}}^k\mathcal{L}_{v_q-v_{\ell,j}}\Theta_j^{g},\\
        \vartheta_w=\sum_j\frac{(-1)^{k_0^g+1}}{(k+1)!}\int_0^{1}(X_{s}^g)_*\left[\mathcal{L}_{\xi^{g}}^{k_0^g+1}(\partial_t+\mathcal{L}_{ v_q})\Theta^{p}_j\right](1-s)^{k_0^g+1}\dd s\\
    \end{cases}
\end{equation*}
and
\begin{equation*}
    \begin{split}
        \begin{cases}
        \Theta_{b,j}^h=\sum_{k=1}^{k_0^g}\frac{(-1)^k}{(k+1)!}\mathcal{L}_{\xi^{p}_j}^k\mathcal{L}_{B_{\ell,j}}\Theta_j^{p}+\sum_{k=0}^{k_0^g}\frac{(-1)^k}{(k+1)!}\mathcal{L}_{\xi^{p}_j}^k\mathcal{L}_{B_q-B_{\ell,j}}\Theta_j^{p},\\
        \vartheta_b=\sum_j\frac{(-1)^{k_0^g+1}}{(k+1)!}\int_0^{1}(X_{s}^g)_*\left[\mathcal{L}_{\xi^{g}}^{k_0^g+1}\mathcal{L}_{B_q}\Theta^{p}_j\right](1-s)^{k_0^g+1}\dd s.
    \end{cases}
    \end{split}
\end{equation*}

\noindent We now manipulate \eqref{Rglin1} further:
\begin{equation}\label{decompling}
    \begin{split}
        &\ddiv \ R^{g,lin}\\
        &=\partial_tw^{g}+\ddiv\left[v_q\otimes w^{g}+w^{g}\otimes v_q-B_q\otimes b^{g}-b^{g}\otimes B_q\right]+\nabla \pi^{gr}+\ddiv \left[\sum_I (g_I^2-1) A_I \right]\\
        &+\ddiv [R_q-R_\ell]\\
        &=\sum_j\left[\partial_tw^{g,p}_j+v_{\ell,j}\cn w^{g,p}_j+w^{g,p}_j\cn v_{\ell,j}-B_{\ell,j}\cn b^{g,p}_j-b^{g,p}_j\cn B_{\ell,j}+\nabla \pi^{gr}_j+\sum_{I:j_I=j}(g_I^2-1)\ddiv A_I\right]\\
        &+\sum_j\left[(\partial_t+\mathcal{L}_{v_{\ell,j}})w^{h}_j-\mathcal{L}_{B_{\ell,j}}b_j^{h}\right]+2\sum_j\left[w^{h}_j\cn v_{\ell,j}-b^{h}_j\cn B_{\ell,j}\right]\\
        &+\ddiv\left[\sum_j\left[(v_q-v_{\ell,j})\otimes (w_j^{g,p}+w_j^{h})\right]^{sym}-\left[(B_q-B_{\ell,j})\otimes (b_j^{g,p}+b_j^{h})\right]^{sym}+R_q-R_\ell\right]\\
        &+\partial_tr_{w}+\ddiv\left[v_q\otimes r_{w}+r_{w}\otimes v_q-B_q\otimes r_{b}-r_{b}\otimes B_q\right]\\
        &= \ddiv\left[R^{cut}+R^{g,tr}+ R^{g,na}+ R^{g,mo}+ R^{g,rem}\right]
    \end{split}
\end{equation}
where, explicitly, we set 
$$R^{g,lin}=R^{cut}+R^{g,tr}+ R^{g,na}+ R^{g,mo}+ R^{g,rem}$$
and defined $R^{g,mo}$ to be the error due to mollification and local recorrection
\begin{equation*}
    \begin{split}
        R^{g,mo}&=\sum_j\left[(v_q-v_{\ell,j})\otimes (w_j^{g,p}+w_j^{h})\right]^{sym}-\left[(B_q-B_{\ell,j})\otimes (b_j^{g,p}+b_j^{h})\right]^{sym}+R_q-R_\ell,\\
    \end{split}
\end{equation*}
$R^{g,tr}$ to be the error coming from the transport of the high-order terms in the Lie-Taylor expansion
\begin{equation*}
    \begin{split}
        R^{g,tr}&=\sum_j\mathcal{R}\left[(\partial_t+\mathcal{L}_{v_{\ell,j}})w^{h}_j-\mathcal{L}_{B_{\ell,j}}b_j^{h}\right]=\sum_j\mathcal{R}\curl\left[(\partial_t +\mathcal{L}_{ v_{\ell,j}})\Theta_{w,j}^h-\mathcal{L}_{ B_{\ell,j}}\Theta_{b,j}^h\right],\\
    \end{split}
\end{equation*}
$R^{g,na}$ to be the error coming from the Nash terms arising from higher-order Lie derivatives
\begin{equation*}
    \begin{split}
        R^{g,na}&=\sum_j2\mathcal{R}\left[w^{h}_j\cn v_{\ell,j}-b^{h}_j\cn B_{\ell,j}\right]=\sum_j2\mathcal{R}\ddiv\left[\Theta_{w,j}^h\times\nabla v_{\ell,j}-\Theta_{b,j}^h\times\nabla  B_{\ell,j}\right]^\top\\
    \end{split}
\end{equation*}
and $R^{g,rem}$ to be the error coming from the remainders
\begin{equation*}
    \begin{split}
        R^{g,rem}&=\mathcal{R}\left[\partial_tr_{w}+\ddiv\left[v_q\otimes r_{w}+r_{w}\otimes v_q-B_q\otimes r_{b}-r_{b}\otimes B_q\right]\right]\\
        &=\mathcal{R}\curl[\partial_t\theta_w]+\mathcal{R}\ddiv\left[v_q\otimes \curl \ \theta_{w}+\curl \ \theta_{w}\otimes v_q-B_q\otimes \curl \ \theta_{b}-\curl \ \theta_{b}\otimes B_q\right],
    \end{split}
\end{equation*}
the error coming from the time cut-off of the Galbrun's equation $R^{cut}=\sum_jR_j^{cut}$ was defined in \eqref{Rcutj}. In the above, we used the identities in \eqref{DG2} and \eqref{DG1}, and $\mathcal{R}$ is the classical inverse divergence operator, see \eqref{invdiv}.

\subsubsection{Estimates on the Transport Error Term}\label{linearerrorsg} 

\begin{lemma}[Estimates on $R^{g,tr}$]\label{Rgtr} Let $\underline{r}=M-k_0^g-m_0-6$. Under the choice of parameters \ref{choiceofparameters}, the following bounds hold:
\begin{equation*}
\begin{split}
    &||R^{g,tr}||_{r+\alpha}\lesssim \lambda_{q}^r\lambda_q^{[r-(\underline{r}-k_0^g)]^+(b-1)}\mathcal{T}_g\ell^{-\alpha}\delta_{q+1} \ \text{ for } \ 0\leq r \leq N-m_0-k_0^g-3,\\
    &||\mathcal{A}^\pm R^{g,tr}||_{r+\alpha}\lesssim \lambda_{q}^r\lambda_q^{[r-(\underline{r}-k_0^g-1)]^+(b-1)}(1/\tau^a)\mathcal{T}_g\ell^{-\alpha}\delta_{q+1} \ \text{ for } \ 0\leq r \leq N-m_0-k_0^g-4.
\end{split}
\end{equation*}
The implicit constants depend on $r$ and all the other parameters, but not on $a$. In particular, the bounds hold for $0\leq r \leq M$ and $0\leq r \leq M-1$, respectively.
\end{lemma}
\begin{proof}[Proof of Lemma \ref{Rgtr}] We first highlight the mollification and local re-correction differences $v_q-v_{\ell,j}, \ B_q-B_{\ell,j}$ in the definition of $R^{g,tr}$. All the Lie derivatives here are in the sense of 1-forms.
\begin{equation}\label{reynoldstransportgalbrun}
    \begin{split}
        R^{g,tr}&=\sum_j\mathcal{R}\curl\left[(\partial_t +\mathcal{L}_{ v_{\ell,j}})\Theta_{w,j}^h-\mathcal{L}_{ B_{\ell,j}}\Theta_{b,j}^h\right]\\
        &=\underbrace{\sum_j\mathcal{R}\curl\sum_{k=1}^{k_0}\frac{(-1)^k}{(k+1)!}\left[(\partial_t +\mathcal{L}_{ v_{\ell,j}})\mathcal{L}_{\xi^{g}}^k(\partial_t+\mathcal{L}_{v_{\ell,j}})\Theta_j^{g}-\mathcal{L}_{ B_{\ell,j}}\mathcal{L}_{\xi^{p}_j}^k\mathcal{L}_{ B_{\ell,j}}\Theta_j^{p}\right]}_{R_1}\\
        &+\underbrace{\sum_j\mathcal{R}\curl\sum_{k=0}^{k_0}\frac{(-1)^k}{(k+1)!}\left[(\partial_t +\mathcal{L}_{ v_{\ell,j}})\mathcal{L}_{\xi^{g}}^k\mathcal{L}_{v_q-v_{\ell,j}}\Theta_j^{g}-\mathcal{L}_{ B_{\ell,j}}\mathcal{L}_{\xi^{p}_j}^k\mathcal{L}_{B_q-B_{\ell,j}}\Theta_j^{p}\right]}_{R_2}\\
    \end{split}
\end{equation}

\noindent We now use the trick in \eqref{trick} to rewrite $R_1, \ R_2,$ by commuting the outside transport operators with the successive Lie derivations. We obtain:
\begin{equation*}
    \begin{split}
        R_1&=\mathcal{R}\curl\sum_j\sum_{k=1}^{k_0^g}\frac{(-1)^k}{(k+1)!} \mathcal{L}_{\xi^{g}}^k\left[\left(\partial_t+\mathcal{L}_{v_{\ell,j}}\right)^2\Theta_j^{g}-\left(\mathcal{L}_{B_{\ell,j}}\right)^2\Theta_j^{g}\right]=:T_1\\
        &+\mathcal{R}\curl\sum_j\sum_{k=1}^{k_0^g}\frac{(-1)^k}{(k+1)!}\sum_{i=0}^{k-1}\left[\mathcal{L}_{\xi^{g}}^i\mathcal{L}_{\partial_t\xi^{g}+[v_{\ell,j},\xi^{g}]} \mathcal{L}_{\xi^{g}}^{(k-1)-i}\left(\partial_t+\mathcal{L}_{v_{\ell,j}}\right)\Theta_j^{g} \right]=:T_{2,w}\\
        &-\mathcal{R}\curl\sum_j\sum_{k=1}^{k_0^g}\frac{(-1)^k}{(k+1)!}\sum_{i=0}^{k-1}\left[\mathcal{L}_{\xi^{g}}^i\mathcal{L}_{[B_{\ell,j},\xi^{g}]} \mathcal{L}_{\xi^{g}}^{(k-1)-i}\mathcal{L}_{B_{\ell,j}}\Theta_j^{g}\right]=:T_{2,b}\\
    \end{split}
\end{equation*}
and
\begin{equation*}
    \begin{split}
        R_2=&\mathcal{R}\curl\sum_j\sum_{k=0}^{k_0^g}\frac{(-1)^k}{(k+1)!}\mathcal{L}_{\xi^{g}}^k\left[\mathcal{L}_{v_q-v_{\ell,j}}\left(\partial_t+\mathcal{L}_{v_{\ell,j}}\right)\Theta_j^{g}-\mathcal{L}_{B_q-B_{\ell,j}}\mathcal{L}_{B_{\ell,j}}\Theta_j^{g}\right]=:T_3\\
        &+\mathcal{R}\curl\sum_j\sum_{k=0}^{k_0^g}\frac{(-1)^k}{(k+1)!} \mathcal{L}_{\xi^{g}}^k\left[\mathcal{L}_{\partial_t(v_q-v_{\ell,j})+[v_{\ell,j},v_q-v_{\ell,j}]}\Theta_j^{g}-\mathcal{L}_{[B_{\ell,j},B_q-B_{\ell,j}]}\Theta_j^{g}\right]=:T_4\\
        &+\mathcal{R}\curl\sum_j\sum_{k=1}^{k_0^g}\frac{(-1)^k}{(k+1)!}\sum_{i=0}^{k-1}\left[\mathcal{L}_{\xi^{g}}^i\mathcal{L}_{\partial_t\xi^{g}+[v_{\ell,j},\xi^{g}]} \mathcal{L}_{\xi^{g}}^{(k-1)-i}\mathcal{L}_{v_q-v_{\ell,j}}\Theta_j^{g}\right]=:T_{5,w}\\
        &-\mathcal{R}\curl\sum_j\sum_{k=1}^{k_0^g}\frac{(-1)^k}{(k+1)!}\sum_{i=0}^{k-1}\left[\mathcal{L}_{\xi^{g}}^i\mathcal{L}_{[B_{\ell,j},\xi^{g}]} \mathcal{L}_{\xi_j}^{(k-1)-i}\mathcal{L}_{B_q-B_{\ell,j}}\Theta_j^{g}\right]=:T_{5,b}\\
    \end{split}
\end{equation*}

\noindent We work on each $T_i$ separately, provide estimates for the top term, and then deal with the multiple Lie-differentiations by means of the Inductive Lemma \ref{inductive} as in Remark \ref{slowinductive}. In what follows, we use most of the time the worst value $\underline{r}$ for the loss function $\lambda_q^{[r-\underline{r}]^+(b-1)\gamma_\ell}$ rather than changing it according to the precise estimate at hand. The key values are: $$\underline{r}=M-m_0-5,\ \ \ \  \underline{r}'=M-m_0-1,$$ 
coming from Lemmas \ref{estimatesperturbationgalbrun}, \ref{estimatesperturbationgalbrun2} and Lemma \ref{localcorrg} respectively (note that in the statement of the Lemma, for notational convenience, we set $\underline{r}=M-m_0-6$).

\noindent\textbf{Estimates on $T_1$.} The top term reads:
$$S=\left(\partial_t+\mathcal{L}_{v_{\ell,j}}\right)^2\Theta_j^{g}-\left(\mathcal{L}_{B_{\ell,j}}\right)^2\Theta_j^{g}.$$
Although we have explicit second-order estimates for $\Theta_j^g$ in Lemma \ref{estgalbrun}, to prove transport estimates for $R^{g,tr}$ we will need a `third-order' one. The explicit formula for this second-order Lie-transport will save the day. From \eqref{secondlietransportgalbrun} we have:
\begin{equation*}
    \begin{split}
        \left(\partial_t+\mathcal{L}_{v_{\ell,j}}\right)^2\Theta_j^{g}-\left(\mathcal{L}_{B_{\ell,j}}\right)^2\Theta_j^{g}&=-\mathbb{T}\tilde H_1(\Theta_j^g)-\mathbb{T}\tilde H_2(\mathcal{A}^+_{\ell,j}\Theta_j^g,\mathcal{A}^-_{\ell,j}\Theta_j^g)+\mathbb{T}F_j
    \end{split}
\end{equation*}
where $$F_j=\sum_{I:I_t=j}f_{p(I)}A_I.$$
To save notation, we omit the $j$ dependence in $\tilde H_1, \ \tilde H_2$. The bounds in Lemmas \ref{slowcoeffgalbrun},  \ref{estgalbrun} together with Lemma \ref{operatornorms} to deal with the operators $\mathbb{T}\tilde H_i$, allow us to compute:
\begin{equation}\label{T1galbrun}
    \begin{split}
        ||\mathbb{T}\tilde H_1(\Theta_j^g)||_{r+\alpha}&\lesssim h_1^r||\Theta^g_j||_0+h_1||\Theta^g_j||_{r+\alpha}\\
        &\lesssim \lambda_q^{r+2}\lambda_q^{[r-\underline{r}]^+(b-1)\gamma_\ell}\lambda_q^{}\ell^{-2\alpha}\tau^c\tau^a\delta_q\delta_{q+1}\\
        &\leq \lambda_q^{r}\lambda_q^{[r-\underline{r}]^+(b-1)\gamma_\ell}\delta_{q+1},
        \\
        ||\mathbb{T}\tilde H_2(\mathcal{A}^+_{\ell,j}\Theta_j^g,\mathcal{A}^-_{\ell,j}\Theta_j^g)||_{r+\alpha}&\lesssim h_2^r(||\mathcal{A}^+_{\ell,j}\Theta^g_j||_0+||\mathcal{A}^-_{\ell,j}\Theta^g_j||_0)+h_2(||\mathcal{A}^+_{\ell,j}\Theta^g_j||_{r+\alpha}+||\mathcal{A}^-_{\ell,j}\Theta^g_j||_{r+\alpha})\\
        &\lesssim \lambda_q^{r+1}\lambda_q^{[r-\underline{r}]^+(b-1)\gamma_\ell}\ell^{-2\alpha}\tau^a\delta_q^{1/2}\delta_{q+1}\\
        &\leq \lambda_q^{r}\lambda_q^{[r-\underline{r}]^+(b-1)\gamma_\ell}\ell^{-\alpha}\delta_{q+1},
        \\
        ||\mathbb{T}F_j||_{r+\alpha}&\lesssim \lambda_q^{r}\lambda_q^{[r-\underline{r}]^+(b-1)\gamma_\ell}\ell^{-\alpha}\delta_{q+1}.\\
    \end{split}
\end{equation}
Here, we used the fact that $\tau^a\ell^{-\alpha}\lambda_q\delta_q^{1/2}\leq 1$, and we will use this several times in the following without further mention. Collecting the bounds, we get: 
$$||S||_{r+\alpha}\lesssim \lambda_q^{r}\lambda_q^{[r-\underline{r}]^+(b-1)\gamma_\ell}\ell^{-\alpha}\delta_{q+1} \ \text{ for } r\geq 0.$$
Concerning Alfv\'en transport, we first compute:
\begin{equation*}
    \begin{split}
        \mathcal{A}^\pm_{\ell,j}S&=\mathcal{A}^\pm_{\ell,j}\left[-\mathbb{T}\tilde H_1(\Theta_j^g)-\mathbb{T}\tilde H_2(\mathcal{A}^+_{\ell,j}\Theta_j^g,\mathcal{A}^-_{\ell,j}\Theta_j^g)+\mathbb{T}F_j\right]\\
        &=-[\mathcal{A}^\pm_{\ell,j},\mathbb{T}\tilde H_1](\Theta_j^g)-[\mathcal{A}^\pm_{\ell,j}, \mathbb{T}\tilde H_2](\mathcal{A}^+_{\ell,j}\Theta_j^g,\mathcal{A}^-_{\ell,j}\Theta_j^g)+[\mathcal{A}^\pm_{\ell,j},\mathbb{T}]F_j\\
        &-\mathbb{T}\tilde H_1(\mathcal{A}^\pm_{\ell,j}\Theta_j^g)-\mathbb{T}\tilde H_2(\mathcal{A}^\pm_{\ell,j}\mathcal{A}^+_{\ell,j}\Theta_j^g,\mathcal{A}^\pm_{\ell,j}\mathcal{A}^-_{\ell,j}\Theta_j^g)+\mathbb{T}\mathcal{A}^\pm_{\ell,j}F_j.
    \end{split}
\end{equation*}
The estimates for the terms obtained after commuting the transport can be carried out in the same way as in \eqref{T1galbrun}, together with the bounds in Lemma \ref{estgalbrun}. We get:
\begin{equation*}
    \begin{split}
        &||\mathbb{T}\tilde H_1(\mathcal{A}^\pm_{\ell,j}\Theta_j^g)-\mathbb{T}\tilde H_2(\mathcal{A}^\pm_{\ell,j}\mathcal{A}^+_{\ell,j}\Theta_j^g,\mathcal{A}^\pm_{\ell,j}\mathcal{A}^-_{\ell,j}\Theta_j^g)+\mathbb{T}\mathcal{A}^\pm_{\ell,j}F_j||_{r+\alpha}\\\
        &\lesssim\lambda_q^{r+2}\lambda_q^{[r-\underline{r}]^+(b-1)\gamma_\ell}\ell^{-2\alpha}\tau^a\delta_q\delta_{q+1}+\lambda_q^{r+1}\lambda_q^{[r-\underline{r}]^+(b-1)\gamma_\ell}\ell^{-2\alpha}\delta_q^{1/2}\delta_{q+1}+\lambda_q^{r}\lambda_q^{[r-\underline{r}]^+(b-1)\gamma_\ell}(1/\tau^a)\ell^{-\alpha}\delta_{q+1}\\
        &\lesssim \lambda_q^{r}\lambda_q^{[r-\underline{r}]^+(b-1)\gamma_\ell}(1/\tau^a)\ell^{-\alpha}\delta_{q+1}.\\
    \end{split}
\end{equation*}
The bounds in Lemmas \ref{slowcoeff},  \ref{estgalbrun}, together with Lemma \ref{operatornorms} to deal with the commutators, give:
\begin{equation*}
    \begin{split}
        ||[\mathcal{A}^\pm_{\ell,j},\mathbb{T}\tilde H_1](\Theta_j^g)||_{r+\alpha}&\lesssim \lambda_q\delta_q^{1/2}\ell^{-\alpha}h_1^r||\Theta^g_j||_\alpha+\lambda_q\delta_q^{1/2}\ell^{-\alpha}h_1||\Theta^g_j||_{r+\alpha}\\
        &\lesssim \lambda_q^{r+3}\lambda_q^{[r-\underline{r}]^+(b-1)\gamma_\ell}\ell^{-3\alpha}\tau^c\tau^a\delta_q^{3/2}\delta_{q+1}\\
        &\leq \lambda_q^{r+1}\lambda_q^{[r-\underline{r}]^+(b-1)\gamma_\ell}\ell^{-\alpha}\delta_q^{1/2}\delta_{q+1},
        \\
        ||[\mathcal{A}^\pm_{\ell,j}, \mathbb{T}\tilde H_2](\mathcal{A}^+_{\ell,j}\Theta_j^g,\mathcal{A}^-_{\ell,j}\Theta_j^g)||_{r+\alpha}&\lesssim \lambda_q\delta_q^{1/2}\ell^{-\alpha}h_2^r(||\mathcal{A}^+_{\ell,j}\Theta^g_j||_\alpha+||\mathcal{A}^-_{\ell,j}\Theta^g_j||_\alpha)\\
        &+\lambda_q\delta_q^{1/2}\ell^{-\alpha}h_2(||\mathcal{A}^+_{\ell,j}\Theta^g_j||_{r+\alpha}+||\mathcal{A}^-_{\ell,j}\Theta^g_j||_{r+\alpha})\\
        &\lesssim \lambda_q^{r+2}\lambda_q^{[r-\underline{r}]^+(b-1)\gamma_\ell}\ell^{-3\alpha}\tau^a\delta_q\delta_{q+1}\\
        &\leq \lambda_q^{r+1}\lambda_q^{[r-\underline{r}]^+(b-1)\gamma_\ell}\ell^{-\alpha}\delta_q^{1/2}\delta_{q+1},
        \\
        ||[\mathcal{A}^\pm_{\ell,j},\mathbb{T}]F_j||_{r+\alpha}&\lesssim\lambda_q^{r+1}\lambda_q^{[r-\underline{r}]^+(b-1)\gamma_\ell}\ell^{-2\alpha}\delta_q^{1/2}\delta_{q+1}.
    \end{split}
\end{equation*}
where we used that $(\tau^a/\tau^c)\ell^{-\alpha}<1$ see \eqref{misc1}. Collecting the bounds, we conclude:
$$||\mathcal{A}^\pm_{\ell,j}S||_{r+\alpha}\lesssim \lambda_q^{r}\lambda_q^{[r-\underline{r}]^+(b-1)\gamma_\ell}(1/\tau^a)\ell^{-\alpha}\delta_{q+1} \ \text{ for } \ r\geq 0.$$

\noindent \textit{Conclusion.} We can now argue by means of Remark \ref{slowinductive}, with: 
\begin{equation*}
        \begin{split}
            &A=\ell^{-\alpha}\delta_{q+1}, \ A_{\mathcal{A}}=(1/\tau^a)\ell^{-\alpha}\delta_{q+1},\\
            &\bar\varsigma_i= \lambda_q\ell^{-\alpha}\tau^a\tau^c\delta_{q+1}, \ \bar\varsigma_{i,\mathcal{A}}=\lambda_q\ell^{-\alpha}\tau^a\delta_{q+1},\\
            &L(r)=L_{\mathcal{A}}(r)=\lambda_q^{[r-\underline{r}]^+(b-1)\gamma_\ell},\\
            &\bar N=N+k_0^g.
        \end{split}
    \end{equation*}
Note that $L=L_{\mathcal{A}}$ here because $\Theta_j^g$ is constructed as a solution to \eqref{galbrunapplied}, and thus we do not lose any good derivative in the transport estimate. For $k\geq 1$ and $0\leq r+k\leq \bar N$, we conclude that: 
\begin{equation}\label{T1kg}
\begin{split}
    &||\mathcal{L}^k_{\xi^g}S||_{r+\alpha}\leq C' (C)^k\lambda_q^{r}\lambda_q^{[r+k-\underline{r}]^+(b-1)\gamma_\ell}\mathcal{T}^k_g\ell^{-\alpha}\delta_{q+1},\\
    &||\mathcal{A}^\pm_{\ell,j}\mathcal{L}^k_{\xi^g}S||_{r+\alpha}\leq C' (C)^k\lambda_q^{r}\lambda_q^{[r+k-\underline{r}]^+(b-1)\gamma_\ell}\mathcal{T}^k_g\ell^{-\alpha}1/\tau^a\delta_{q+1}.
\end{split}
\end{equation}
The constants $C,\ C'$ depend on $\bar N$ and all the other parameters, but not on $a$ and are uniform in $r,\ k, \ j$.

\noindent We will argue as in the proof of the transport bounds in Lemma \ref{estimatesfieldsringg}.  By means of Proposition \ref{czstuff} to deal with $\mathcal{R}\curl$ and the fact that there are at most two non-zero $\Theta^g_j$ at each time, we conclude that:
\begin{equation}\label{T1g}
    \begin{split}
        ||T_1||_{r+\alpha}&\lesssim \sup_j\sum_{k=1}^{k_0^g} \frac{1}{(k+1)!}||\mathcal{R}\curl \mathcal{L}_{\xi^{g}}^kS||_{r+\alpha}\lesssim\sup_j\sum_{k=1}^{k_0^g} \frac{1}{(k+1)!}|| \mathcal{L}_{\xi^{g}}^kS||_{r+\alpha}\\
        &\lesssim\sum_{k=1}^{k_0^g} \frac{(C)^k}{(k+1)!}\lambda_q^{r}\lambda_q^{[r+k-\underline{r}]^+(b-1)\gamma_\ell}\mathcal{T}^k_g\ell^{-\alpha}\delta_{q+1}\\
        &=\lambda_q^{r}\ell^{-\alpha}\delta_{q+1}\Biggr[\underbrace{\sum_{k=1}^{[\underline{r}-r]^+} \frac{\left(C\mathcal{T}_g\right)^k}{(k+1)!}}_{\lesssim 1_{[\underline{r}-r]^+\geq 1}\mathcal{T}_g}+\underbrace{\sum_{k=[\underline{r}-r]^++1}^{k_0^g} \frac{\left(C\lambda_q^{(b-1)\gamma_\ell}\mathcal{T}_g\right)^k}{(k+1)!}\lambda_q^{[r-\underline{r}]^+(b-1)\gamma_\ell}}_{\lesssim 1_{[\underline{r}-r]^++1\leq k_0^g}\left(\lambda_q^{(b-1)\gamma_\ell}\mathcal{T}_g\right)^{[\underline{r}-r]^++1}}\Biggr]\\
        &\lesssim \lambda_q^{r}\mathcal{T}_g\ell^{-\alpha}\delta_{q+1}\underbrace{\left[1_{[\underline{r}-r]^+\geq 1}+1_{[\underline{r}-r]^++1\leq k_0^g}\lambda_q^{[r-\underline{r}]^+(b-1)\gamma_\ell}\lambda_q^{(b-1)\gamma_\ell}\left(\lambda_q^{(b-1)\gamma_\ell}\mathcal{T}_g\right)^{[\underline{r}-r]^+}\right]}_{f(r)}\\
    \end{split}
\end{equation}
for $0\leq r\leq \bar N-k_0^g=N$. 

\noindent Now note that:
$$1_{[\underline{r}-r]^+\geq 1}=1_{r\leq \underline{r}-1} \ \text{ and } \ 1_{[\underline{r}-r]^++1\leq k_0^g}=1_{r\geq \underline{r}-k_0^g+1},$$
we can then rewrite
\begin{equation*}
    \begin{split}
        f(r)&=1_{[\underline{r}-r]^+\geq 1}+1_{[\underline{r}-r]^++1\leq k_0^g} \lambda_q^{[r-\underline{r}]^+(b-1)\gamma_\ell}\lambda_q^{(b-1)\gamma_\ell}\left(\lambda_q^{(b-1)\gamma_\ell}\mathcal{T}_g\right)^{[\underline{r}-r]^+}\\
        &=1_{r\leq \underline{r}-1}+1_{r\geq \underline{r}-k_0^g+1} \lambda_q^{[r-\underline{r}]^+(b-1)\gamma_\ell}\lambda_q^{(b-1)\gamma_\ell}\left(\lambda_q^{(b-1)\gamma_\ell}\mathcal{T}_g\right)^{[\underline{r}-r]^+}\\
        &=1_{r\leq \underline{r}-1}+1_{\underline{r}-k_0^g+1\leq r\leq \underline{r}} \lambda_q^{(b-1)\gamma_\ell}\left(\lambda_q^{(b-1)\gamma_\ell}\mathcal{T}_g\right)^{\underline{r}-r}+1_{ r\geq \underline{r}+1}\lambda_q^{[r-\underline{r}](b-1)\gamma_\ell}\lambda_q^{(b-1)\gamma_\ell}\\
        &\leq 2\lambda_q^{[r-(\underline{r}-k_0^g)]^+(b-1)\gamma_\ell}
    \end{split}
\end{equation*}
and conclude:
$$||T_1||_{r+\alpha}\lesssim \lambda_q^{r}\lambda_q^{[r-(\underline{r}-k_0^g)]^+(b-1)\gamma_\ell}\mathcal{T}_g\ell^{-\alpha}\delta_{q+1} \ \text{ for } \ r\geq 0.$$
where the implicit constant depends on $r$ and all the other parameters, but not on $a$. This bound is not optimal with respect to the loss function, but other terms in $R^{g,tr}$ will yield worse estimates, see $T_4$ below.

\noindent From the bounds in Lemmas \ref{localcorrg}, \ref{stabilitygalbrun} \ref{standardmollgalbrun}, to correct the transport, \eqref{T1kg} above and Proposition \ref{czstuff} to deal with the commutator, we deduce:
\begin{equation*}
    \begin{split}
        &||\mathcal{A}^\pm T_1||_{r+\alpha}\\
        &\lesssim \sup_j\sum_{k=1}^{k_0^g} \frac{1}{(k+1)!} \left[||\mathcal{R}\curl \mathcal{A}^\pm_{\ell,j}\mathcal{L}_{\xi^{g}}^kS||_{r+\alpha}+||[\mathcal{A}^\pm_{\ell,j},\mathcal{R}\curl] \mathcal{L}_{\xi^{g}}^kS||_{r+\alpha}+||(z_q^\pm-z^\pm_{\ell,j})\cn\mathcal{R}\curl \mathcal{L}_{\xi^{g}}^kS||_{r+\alpha}\right]\\
        &\lesssim \sup_j\sum_{k=1}^{k_0^g} \frac{1}{(k+1)!} \underbrace{\left[||z^\pm_{\ell,j}||_{r+1+\alpha}||\mathcal{L}_{\xi^{g}}^kS||_{\alpha}+||z^\pm_{\ell,j}||_{1+\alpha}||\mathcal{L}_{\xi^{g}}^kS||_{r+\alpha}\right]}_{\lambda_q^{r+1}\mathcal{T}^k_g\ell^{-2\alpha}\delta_q^{1/2}\delta_{q+1}}\\
        &+ \sup_j\sum_{k=1}^{k_0^g} \frac{1}{(k+1)!}\underbrace{\left[||z_q^\pm-z^\pm_{\ell,j}||_{r+\alpha}||\mathcal{L}_{\xi^{g}}^kS||_{1+\alpha}+||z_q^\pm-z^\pm_{\ell,j}||_{0}|| \mathcal{L}_{\xi^{g}}^kS||_{r+1+\alpha}\right]}_{\lambda_q^{r+1}\mathcal{T}^k_g(\ell\lambda_q)^{m_0}\ell^{-2\alpha}\delta_q^{1/2}\delta_{q+1}}\\
        &+ \sup_j\sum_{k=1}^{k_0^g} \frac{1}{(k+1)!}\underbrace{|| \mathcal{A}^\pm_{\ell,j}\mathcal{L}_{\xi^{g}}^kS||_{r+\alpha}}_{\lambda_q^{r}\mathcal{T}_g^k(1/\tau^a)\ell^{-\alpha}\delta_{q+1}}\\
        &\lesssim\sum_{k=1}^{k_0^g} \frac{(C)^k}{(k+1)!}\lambda_q^{r}\lambda_q^{[r+k-\underline{r}]^+(b-1)\gamma_\ell}\mathcal{T}^k_g\ell^{-\alpha}1/\tau^a\delta_{q+1}\\
    \end{split}
\end{equation*}
for $0\leq r\leq \min\{N-m_0-1,\ \bar N-k_0^p\}=N-m_0-1$, where we used $1-\gamma_\ell>0$ see \eqref{constraintadmissibility}. 

\noindent To deal with the sum, we now argue as before and split it according to the number of derivatives $r$. We omit the details. 
\begin{equation*}
    \begin{split}
        ||\mathcal{A}^\pm T_1||_{r+\alpha}&\lesssim\lambda_q^{r}\ell^{-\alpha}1/\tau^a\delta_{q+1}\left[\sum_{k=1}^{[\underline{r}-r]^+} \frac{\left(C\mathcal{T}_g\right)^k}{(k+1)!}+\sum_{k=[\underline{r}-r]^++1}^{k_0^g} \frac{\left(C\lambda_q^{(b-1)\gamma_\ell}\mathcal{T}_g\right)^k}{(k+1)!}\lambda_q^{[r-\underline{r}]^+(b-1)\gamma_\ell}\right]\\
        &\lesssim \lambda_q^{r}\mathcal{T}_g\ell^{-\alpha}(1/\tau^a)\delta_{q+1}\left[1_{[\underline{r}-r]^+\geq 1}+1_{[\underline{r}-r]^++1\leq k_0^g}\lambda_q^{[r-\underline{r}]^+(b-1)\gamma_\ell}\left(\lambda_q^{(b-1)\gamma_\ell}\mathcal{T}_g\right)^{[\underline{r}-r]^+}\right]\\
        &\lesssim \lambda_q^{r}\lambda_q^{[r-(\underline{r}-k_0^g)]^+(b-1)\gamma_\ell}\mathcal{T}_g\ell^{-\alpha}(1/\tau^a)\delta_{q+1}
    \end{split}
\end{equation*}
for $0\leq r\leq N-m_0-1$. The implicit constant depends on $r$ and all the other parameters, but not on $a$.

\noindent \textbf{Estimates on $T_{4}$.} The top term can be written as:
$$\mathcal{L}_{(\partial_t+\mathcal{L}_{v_{\ell,j}})(v_q-v_{\ell,j})-\mathcal{L}_{B_{\ell,j}}(B_q-B_{\ell,j})}\Theta_j^{g}$$
and adopting a similar notation to Section \ref{secmollnash}, we define:
\begin{equation*}
    \begin{split}
        &\Delta^v= v_\ell - v_{\ell,j}, \ 
        \Delta^B= B_\ell - B_{\ell,j}, \ \Delta^p=p_\ell-p_{\ell,j}, \\
        &S=(\partial_t+\mathcal{L}_{v_{\ell,j}})(v_q-v_{\ell,j})-\mathcal{L}_{B_{\ell,j}}(B_q-B_{\ell,j})
    \end{split}
\end{equation*}
where the Lie derivative here is just the commutator of vector fields. We begin with the following computation:
\begin{equation}\label{nastytrasnporttermg}
    \begin{split}
    S&=(\partial_t+\mathcal{L}_{ v_{\ell,j}})( v_q- v_\ell)+(\partial_t+\mathcal{L}_{ v_{\ell,j}})( v_\ell- v_{\ell,j})-\mathcal{L}_{ B_{\ell,j}}( B_q- B_\ell)-\mathcal{L}_{ B_{\ell,j}}( B_\ell- B_{\ell,j})\\
    &=\underbrace{(\partial_t+\mathcal{L}_{ v_{\ell,j}})\Delta^v-\mathcal{L}_{ B_{\ell,j}}\Delta^B}_{S_1}+\underbrace{(\partial_t+ v_q\cn)( v_q- v_\ell)- B_q\cn( B_q- B_\ell)}_{S_2}\\
    &\underbrace{-( v_q- v_\ell)\cn  v_{\ell,j}+( B_q- B_\ell)\cn  B_{\ell,j}+( v_{\ell,j}- v_q)\cn ( v_q- v_\ell)-( B_{\ell,j}- B_q)\cn ( B_q- B_\ell)}_{S_3}\\
    &=S_1+S_2+S_3
    \end{split}
\end{equation}
and we now treat each $S_i$ separately.

\noindent \textit{Estimates on $S_1$.} We could estimate $S_1$ using Lemmas \ref{localcorrg} and \ref{stabilitygalbrun}. Note, however, that we must consider the transport of the new Reynold stress itself; the explicit representation \eqref{specialdoublegalbrun} will save us the day. This reads:
\begin{equation*}
    \begin{split}
        (\partial_t+\mathcal{L}_{v_{\ell,j}})\Delta^v+\mathcal{L}_{B_{\ell,j}}\Delta^B&=\ddiv \ R_\ell^c-2\left(\Delta^v\cn  v_{\ell,j}-\Delta^B\cn  B_{\ell,j}\right)-\ddiv\left[\Delta^v\otimes\Delta^v-\Delta^B\otimes\Delta^B\right]-\nabla\Delta^p\\
        &=\underbrace{\mathbb{P}\ddiv \ R_\ell^c}_{\lambda_q^{r+1}\ell^{-\alpha}(\ell\lambda_q)^{m_0}\delta_q}-\underbrace{2\mathbb{P}\left(\Delta^v\cn  v_{\ell,j}-\Delta^B\cn  B_{\ell,j}\right)}_{\lambda_q^{r+1}\ell^{-\alpha}(\ell\lambda_q)^{m_0}\delta_q}-\underbrace{\mathbb{P}\ddiv\left[\Delta^v\otimes\Delta^v-\Delta^B\otimes\Delta^B\right]}_{\lambda_q^{r+1}\ell^{-\alpha}(\ell\lambda_q)^{2m_0}\delta_q}
    \end{split}
\end{equation*}
This is the main reason for constructing the velocity and magnetic-field re-correction by solving the full nonlinear equation. By means of Lemmas \ref{standardmollgalbrun}, \ref{stabilitygalbrun} and \ref{stabilitygalbrun} we bound:
\begin{equation*}
    \begin{split}
        ||S_1||_{r+\alpha}
        &\lesssim\lambda_q^{r+1}\lambda_q^{[r-\underline{r}']^+(b-1)\gamma_\ell}\ell^{-\alpha}(\ell\lambda_q)^{m_0}\delta_q \ \text{ for } \ r\geq 0\\
    \end{split}
\end{equation*}
where we used $1-\gamma_\ell>0$, see \eqref{constraintadmissibility}. The transport estimate now follows immediately from the same Lemmas and the commutator estimate in Proposition \ref{czstuff}, indeed: 
\begin{equation*}
    \begin{split}
        \mathcal{A}^\pm_{\ell,j}S_1&=\mathcal{A}^\pm_{\ell,j}\left[\mathbb{P}\ddiv \ R_\ell^c-2\mathbb{P}\left(\Delta^v\cn  v_{\ell,j}-\Delta^B\cn  B_{\ell,j}\right)-\mathbb{P}\ddiv\left[\Delta^v\otimes\Delta^v-\Delta^B\otimes\Delta^B\right]\right]\\
        &=\underbrace{[\mathcal{A}^\pm_{\ell,j},\mathbb{P}\ddiv](R^c_\ell-\Delta^v\otimes\Delta^v+\Delta^B\otimes\Delta^B)}_{\lambda_q^{r+2}\ell^{-2\alpha}(\ell\lambda_q)^{m_0}\delta_q^{3/2}+\lambda_q^{r+2}\ell^{-2\alpha}(\ell\lambda_q)^{2m_0}\delta_q^{3/2}}-\underbrace{2[\mathcal{A}^\pm_{\ell,j},\mathbb{P}](\Delta^v\cn  v_{\ell,j}-\Delta^B\cn  B_{\ell,j})}_{\lambda_q^{r+2}\ell^{-2\alpha}(\ell\lambda_q)^{m_0}\delta_q^{3/2}}\\
        &+\underbrace{\mathbb{P}\ddiv\left[\mathcal{A}^\pm_\ell R_\ell^c+( z_{\ell,j}^\pm- z_\ell^\pm)\cn R^c_\ell-[\mathcal{A}^\pm_{\ell,j}\Delta^v\otimes\Delta^v]^{sym}+[\mathcal{A}^\pm_{\ell,j}\Delta^B\otimes\Delta^B]^{sym}\right]}_{\lambda_q^{r+2}\ell^{-\alpha}(\ell\lambda_q)^{m_0}\delta_q^{3/2}+\lambda_q^{r+2}\ell^{-\alpha}(\ell\lambda_q)^{2m_0}\delta_q^{3/2}}\\
        &-\underbrace{2\mathbb{P}\left[(\mathcal{A}^\pm_{\ell,j}\Delta^v)\cn  v_{\ell,j}+\Delta^v\cn \mathcal{A}^\pm_{\ell,j} v_{\ell,j}-(\Delta^v\cn  z^\pm_{\ell,j})\cn  v_{\ell,j}\right]}_{\lambda_q^{r+2}\ell^{-\alpha}(\ell\lambda_q)^{m_0}\delta_q^{3/2}+\lambda_q^{r+2}\ell^{-\alpha}(\ell\lambda_q)^{2m_0}\delta_q^{3/2}}\\
        &+2\mathbb{P}\left[(\mathcal{A}^\pm_{\ell,j}\Delta^B)\cn  B_{\ell,j}+\Delta^B\cn \mathcal{A}^\pm_{\ell,j} B_{\ell,j}-(\Delta^B\cn  z^\pm_{\ell,j})\cn  B_{\ell,j}\right],\\
    \end{split}
\end{equation*}
we conclude that:
\begin{equation*}
    \begin{split}
        ||\mathcal{ A}^\pm_{\ell,j}S_1||_{r+\alpha}
        &\lesssim\lambda_q^{r+2}\lambda_q^{[r-(\underline{r}'-1)]^+(b-1)\gamma_\ell}\ell^{-2\alpha}(\ell\lambda_q)^{m_0}\delta_q^{3/2} \ \text{ for } \ r\geq 0.
    \end{split}
\end{equation*}

\noindent \textit{Estimates on $S_2$.} From Lemma \ref{standardmollgalbrun} we immediately deduce:
\begin{equation}
    \begin{split}
    ||S_2||_r&\lesssim ||(\partial_t+ v_q\cn)( v_q- v_\ell)||_r+ ||B_q\cn( B_q- B_\ell)||_r\\
       &=\frac{1}{2}||(\mathcal{A}^++\mathcal{A}^-)( v_q- v_\ell)||_r+ \frac{1}{2}||(\mathcal{A}^+-\mathcal{A}^-)( B_q- B_\ell)||_r\\
       &\lesssim\lambda_q^{r+1}(\ell\lambda_q)^{m_0} \delta_{q}
    \end{split}
\end{equation}
for $0\leq r\leq M-m_0-1$. For any derivative in the range $M-m_0-1\leq r \leq N- m_0-1$, we are forced to bound:
\begin{equation*}
    \begin{split}
        ||S_2||_{r}&\lesssim ||\partial_tv_q||_{r}+||v_q\cn (v_q- v_\ell)||_{r}+||B_q\cn (B_q- B_\ell)||_{r}\\
        &\lesssim ||\partial_tv_q||_{r}+||v_q||_{r}||v_q- v_\ell||_{1}+||v_q||_{0}||v_q- v_\ell||_{r+1}\\
        &+||B_q||_{r}||B_q- B_\ell||_{1}+||B_q||_{0}||B_q- B_\ell||_{r+1}\\
        &\lesssim \lambda_q^{r+1}(\ell\lambda_q)^{m_0}\delta_q^{1/2}.
    \end{split}
\end{equation*}
To estimate $\mathcal{A}^\pm_{\ell,j}S_2$, we first rewrite $S_2$ as a commutator between mollification and the MHD equation. Namely:
\begin{equation}\label{nasty2g}
    \begin{split}
    S_2&=(\partial_t+ v_q\cn)( v_q- v_\ell)- B_q\cn( B_q- B_\ell)\\
        &=(\partial_t+ v_q\cn)( v_q-( v_q)_\ell)- B_q\cn ( B_q-( B_q)_\ell)\\
        &=\left[(\partial_t+ v_q\cn) v_q- B_q\cn  B_q\right]-\left[(\partial_t+ v_q\cn) v_q- B_q\cn  B_q\right]_\ell\\
        &-([ v_q\cn,*\rho_\ell] v_q-[ B_q\cn,*\rho_\ell] B_q)\\
        &=\ddiv\left[R_q-R_\ell\right]-\nabla (p_q- p_\ell)-([ v_q\cn,*\rho_\ell] v_q-[ B_q\cn,*\rho_\ell] B_q)\\
    \end{split}
\end{equation}
where we used that $(v_q, B_q, p_q, R_q)$ solves \eqref{relaxedMHD}. We now apply $\mathcal{A}^\pm_{\ell,j}$ to \eqref{nasty2g} an write it as:
\begin{equation*}
\begin{split}
    \mathcal{A}^\pm_{\ell,j}S_2&=\mathcal{A}^\pm_{\ell,j}\ddiv[R_q-R_\ell]+\mathcal{A}^\pm_{\ell,j}\nabla(p_q-p_\ell)-(\mathcal{A}^\pm_{\ell,j}[ v_q\cn,*\rho_\ell] v_q-\mathcal{A}^\pm_{\ell,j}[ B_q\cn,*\rho_\ell] B_q)\\
    &=\underbrace{\ddiv\mathcal{A}^\pm\left[R_q-R_\ell+(p_q-p_\ell)\IId\right]}_{\lambda_q^{r+2}(\ell\lambda_q)^{m_0} \delta_q^{1/2}\delta_{q+1}+\lambda_q^{r+2}(\ell\lambda_q)^{m_0} \delta_{q}^{3/2}}-\underbrace{\left([ v_q\cn,*\rho_\ell] \mathcal{A}^\pm v_q-[ B_q\cn,*\rho_\ell] \mathcal{A}^\pm B_q\right)}_{\lambda_q^{r+2}(\ell\lambda_q)^{m_0}\delta_q^{3/2}}\\
    &+\underbrace{(z^\pm_{\ell,j}-z_q^\pm)\cn\left[\ddiv[R_q-R_\ell]+\nabla(p_q-p_\ell)-([ v_q\cn,*\rho_\ell] v_q-[ B_q\cn,*\rho_\ell] B_q)\right]}_{\lambda_q^{r+2}(\ell\lambda_q)^{2m_0}\delta_q^{1/2}\delta_{q+1}+\lambda_q^{r+2}(\ell\lambda_q)^{2m_0}\delta_q^{3/2}}\\
    &+\underbrace{[\mathcal{A}^\pm,\ddiv][R_q-R_\ell +(p_q-p_\ell)\IId]}_{\lambda_q^{r+2}(\ell\lambda_q)^{m_0}\delta_q^{1/2}\delta_{q+1}+\lambda_q^{r+2}(\ell\lambda_q)^{m_0}\delta_q^{3/2}}-\underbrace{([\mathcal{A}^\pm,[ v_q\cn,*\rho_\ell]] v_q-[\mathcal{A}^\pm,[ B_q\cn,*\rho_\ell]] B_q)}_{\lambda_q^{r+2}(\ell\lambda_q)^{m_0}\delta_q^{3/2}}\\
\end{split}        
\end{equation*} 
Note the double commutator of transport with mollification, according to Proposition \ref{deepmollification} and the Iterative Assumptions \eqref{inductiveassumptionsgeneral}, we have: 
\begin{equation*}
    \begin{split}
        &||[\mathcal{A}^\pm,[ v_q\cn,*\rho_\ell] v_q]||_{r}, \ ||[\mathcal{A}^\pm,[ B_q\cn,*\rho_\ell] B_q]||_{r}\lesssim\lambda_q^{r} (\ell\lambda_q)^{m_0}\lambda_q^2\delta_q^{3/2} \ \text{ for } \ 0\leq r\leq \underline{r}'-1,\\
        &||[ v_q\cn,*\rho_\ell] \mathcal{A}^\pm v_q]||_{r}, \ ||[ B_q\cn,*\rho_\ell] \mathcal{A}^\pm B_q]||_{r}\lesssim\lambda_q^{r} (\ell\lambda_q)^{m_0}\lambda_q^2\delta_q^{3/2} \ \text{ for } \ 0\leq r\leq \underline{r}'-1\\
    \end{split}
\end{equation*}
and using in addition the bounds in Lemmas \ref{standardmollgalbrun}, \ref{stabilitygalbrun}, \ref{localcorrg}, we obtain:
\begin{equation*}
    \begin{split}
        ||\mathcal{A}^\pm_{\ell,j}S_2||_{r}
        &\lesssim\lambda_q^{r+2}(\ell\lambda_q)^{m_0}\delta_q^{3/2} \ \text{ for } \ 0\leq r\leq \underline{r}'-1.
    \end{split}
\end{equation*}
For any derivative in the range $\underline{r}'-1\leq r \leq N- m_0-2$, we are forced to bound:
\begin{equation*}
        ||\mathcal{A}^\pm_{\ell,j}S_2||_{r}\lesssim ||\partial_tS_2||_{r}+||z^\pm_{\ell,j}||_0||S_2||_{r+1}+||S_2||_{1}||z^\pm_{\ell,j}||_{r+1}\\
\end{equation*}
where
\begin{equation}\label{puresecondT4g}
    \begin{split}
        \partial_tS_2&=\partial_t[(\partial_t+v_q\cn)(v_q-v_\ell)]\\
        &=\partial_t^2(v_q-v_\ell)+\partial_t v_q\cn (v_q-v_\ell)+v_q\cn \partial_t(v_q-v_\ell),
    \end{split}
\end{equation}
similarly for the magnetic field, the bounds in Lemma \ref{standardmollgalbrun} then give:
\begin{equation*}
    ||\partial_t S_2||_r\lesssim \lambda_q^{r+2}(\ell\lambda_q)^{m_0}\delta_q^{1/2}
\end{equation*}
and we conclude:
\begin{equation*}
    \begin{split}
        ||\mathcal{A}^\pm_{\ell,j}S_2||_{r}\lesssim  \lambda_q^{r+2}\lambda_q^{[r-(\underline{r}'-1)]^+(b-1)\gamma_\ell}(\ell\lambda_q)^{m_0}\delta_q^{1/2} \ \text{ for } \ \underline{r}'-1\leq r \leq N- m_0-2
    \end{split}
\end{equation*}
where, the loss function appears because of the estimates on $z^\pm_{\ell,j}$ in Lemma \ref{localcorrg}.

\noindent \textit{Estimates on $S_3$.} We have:
\begin{equation*}
    \begin{split}
        S_3&=-[\underbrace{( v_q- v_\ell)\cn  v_{\ell,j}-( B_q- B_\ell)\cn  B_{\ell,j}}_{\lambda_q^{r+1}\ell^{-\alpha}(\ell\lambda_q)^{m_0}\delta_q}]+\underbrace{(v_{\ell,j}- v_q)\cn ( v_q- v_\ell)-( B_{\ell,j}- B_q)\cn ( B_q- B_\ell)}_{\lambda_q^{r+1}(\ell\lambda_q)^{2m_0}\delta_q}
    \end{split}
\end{equation*}
and Lemmas \ref{standardmollgalbrun}, \ref{stabilitygalbrun} and \ref{localcorrg}  give the necessary bounds, we conclude:
\begin{equation*}
    \begin{split}
        ||S_3||_{r+\alpha}
        &\lesssim \lambda_q^{r+1}\lambda_q^{[r-(\underline{r}'-1)]^+(b-1)\gamma_\ell}\ell^{-\alpha}(\ell\lambda_q)^{m_0}\delta_q \ \text{ for } \ 0\leq r \leq N-m_0-2
    \end{split}
\end{equation*}
where the additional loss of one derivative comes from interpolating the $C^{r+\alpha}$ norm in the non-mollified terms. 

\noindent Commuting the transport operator with $\cn$, we rewrite:
\begin{equation*}
    \begin{split}
        \mathcal{A}^\pm_{\ell,j}S_3&=-(v_q-v_\ell)\cn \mathcal{A}^\pm_{\ell,j} v_{\ell,j}+( B_q- B_\ell)\cn  \mathcal{A}^\pm_{\ell,j}B_{\ell,j}\\
        &-[\mathcal{A}^\pm(v_q-v_\ell)-(v_q-v_\ell)\cn z^{\pm}_{\ell,j}+(z^\pm_{\ell,j}-z^\pm_q)\cn(v_q-v_\ell)]\cn  v_{\ell,j}\\
        &+[\mathcal{A}^\pm(B_q-B_\ell)-(B_q-B_\ell)\cn z^{\pm}_{\ell,j}+(z^\pm_{\ell,j}-z^\pm_q)\cn(B_q-B_\ell)]\cn  B_{\ell,j}\\
        &+[\mathcal{A}^\pm_{\ell,j}(v_{\ell,j}- v_\ell)-(v_{\ell,j}- v_\ell)\cn z^\pm_{\ell,j}]\cn ( v_q- v_\ell)\\
        &+(v_\ell-v_q)\cn [\mathcal{A}^\pm( v_q- v_\ell)+(z^\pm_{\ell,j}-z^\pm_q)\cn(v_q-v_\ell)]\\
        &-[\mathcal{A}^\pm_{\ell,j}(B_{\ell,j}- B_\ell)-(B_{\ell,j}- B_\ell)\cn z^\pm_{\ell,j}]\cn (B_q- B_\ell)\\
        &+(B_\ell-B_q)\cn [\mathcal{A}^\pm( B_q- B_\ell)+(z^\pm_{\ell,j}-z^\pm_q)\cn(B_q-B_\ell)]\\
    \end{split}
\end{equation*}
and we conclude from the Lemmas above that:
\begin{equation*}
    \begin{split}
        ||\mathcal{A}^\pm_{\ell,j}S_3||_{r+\alpha}
        &\lesssim\lambda_q^{r+2}\ell^{-\alpha}(\ell\lambda_q)^{m_0}\delta_q^{3/2} \ \text{ for } \ 0\leq r \leq \underline{r}'-2.
    \end{split}
\end{equation*}
For any derivative in the range $\underline{r}'-2\leq r \leq N-m_0-3$, as in $S_2$, we are forced to bound the transport derivatives with pure derivatives and get:
\begin{equation*}
    \begin{split}
        ||\mathcal{A}^\pm_{\ell,j}S_3||_{r+\alpha}&\lesssim\lambda_q^{r+2}\lambda_q^{[r-(\underline{r}'-2)](b-1)\gamma_\ell}\ell^{-\alpha}(\ell\lambda_q)^{m_0}\delta_q.
    \end{split}
\end{equation*}

\noindent \textit{Conclusion.} Let us now set loss functions: 
\begin{equation}\label{T4g}
\begin{split}
    &\bar L=1/\delta_q^{1/2} \ \leadsto \ L=1_{r\leq \underline{r}-1}+1_{r\geq \underline{r}}\lambda_q^{([r-(\underline{r}-1)](b-1)\gamma_\ell}\bar L,\\
    &\bar{\bar L}=1/\delta_q,\ \leadsto \ L_{\mathcal{A}}=1_{r\leq \underline{r}-2}+1_{r\geq \underline{r}-1}\lambda_q^{([r-(\underline{r}-2)](b-1)\gamma_\ell}\bar{ \bar L}.\\
\end{split}
\end{equation}
The additional shift of $-1,-2$ is to get a cleaner formula afterwards and guarantee that
$$L(r+1)\leq L_{\mathcal{A}}(r),$$
note that in the definition we use $\underline{r}$, which is associated with derivatives of $\Theta_j^g$, not $\underline{r}'$. 

\noindent We now collect all the bounds for $S$ and use interpolation to deduce the $C^{r+\alpha}$ norm of $S_2$. Since $\underline{r}'-3\geq  \underline{r}$, we can write:
\begin{equation*}
    \begin{split}
        ||S||_{r+\alpha}&\lesssim||S_1||_{r+\alpha}+||S_2||_{r+\alpha}+||S_3||_{r+\alpha}\\
        &\lesssim \lambda_q^{r+1}L(r-1)\ell^{-\alpha}(\ell\lambda_q)^{m_0}\delta_q \ \text{ for } \ 0\leq N-m_0-2,
        \\
        ||\mathcal{A}^\pm_{\ell,j}S||_{r+\alpha}&\lesssim||\mathcal{A}^\pm_{\ell,j}S_1||_{r+\alpha}+||\mathcal{A}^\pm_{\ell,j}S_2||_{r+\alpha}+||\mathcal{A}^\pm_{\ell,j}S_3||_{r+\alpha}\\
        &\lesssim \lambda_q^{r+2} L_{\mathcal{A}}(r-1)\ell^{-2\alpha}(\ell\lambda_q)^{m_0}\delta_q^{3/2} \ \text{ for } \ 0\leq r\leq  N- m_0-3.\\
    \end{split}
\end{equation*}

\noindent We conclude that the top term in $T_4$ satisfies:
\begin{equation*}
    \begin{split}
        ||\mathcal{L}_S\Theta_j^g||_{r+\alpha}&\lesssim ||S||_{r+\alpha}||\Theta_j^g||_{1}+||S||_{0}||\Theta_j^g||_{r+1+\alpha}+||S||_{r+1+\alpha}||\Theta_j^g||_{0}+||S||_{1}||\Theta_j^g||_{r+\alpha}\\
        &\lesssim\lambda_q^{r+2}L(r)\ell^{-2\alpha}(\ell\lambda_q)^{m_0}\tau^c\tau^a\delta_q\delta_{q+1}\\
        &\leq\lambda_q^{r}L(r)(\ell\lambda_q)^{m_0}(\tau^a/\tau^c)\delta_{q+1}\\
    \end{split}
    \end{equation*}
for $0\leq r\leq N-m_0-3$. Here we used $\tau^c\ell^{-\alpha}\lambda_q\delta_q^{1/2}\leq 1$.  Moreover, 
\begin{equation*}
    \begin{split}
        ||\mathcal{A}^\pm_{\ell,j}\mathcal{L}_S\Theta_j^g||_{r+\alpha}&\lesssim ||\mathcal{L}_{(\partial_t+\mathcal{L}_{z^\pm_{\ell,j}})S}\Theta_j^g||_{r+\alpha}+||\mathcal{L}_{S}(\partial_t+\mathcal{L}_{z^\pm_{\ell,j}})\Theta_j^g||_{r+\alpha}+||\mathcal{L}_S\Theta_j^g\cn z^\pm_{\ell,j}||_{r+\alpha}\\
        &\lesssim\lambda_q^{r}L_{\mathcal{A}}(r)(\ell\lambda_q)^{m_0}(\tau^a/\tau^c)(1/\tau^c)\delta_{q+1}
    \end{split}
\end{equation*}
for $0\leq r \leq N- m_0-4$. 

\noindent We can now apply Lemma \ref{inductive} as in Remark \ref{slowinductive} with $L,L_{\mathcal{A}}$ as in \eqref{T4g} and
\begin{equation*}
        \begin{split}
            &A=(\ell\lambda_q)^{m_0}(\tau^a/\tau^c)\delta_{q+1}, \ A_{\mathcal{A}}=(\ell\lambda_q)^{m_0}(\tau^a/\tau^c)(1/\tau^c)\delta_{q+1,}\\
            &\bar\varsigma_i= \lambda_q\ell^{-\alpha}\tau^a\tau^c\delta_{q+1}, \ \bar\varsigma_{i,\mathcal{A}}=\lambda_q\ell^{-\alpha}\tau^a\delta_{q+1},\\
            &\bar N=N-m_0-3,
        \end{split}
\end{equation*}
we deduce:
\begin{equation*}
    \begin{split}
        &||\mathcal{L}_{\xi^g}^k[\mathcal{L}_{S}\Theta_j^g]||_{r+\alpha}
        \leq C'(C)^k\lambda_q^{r}L(r+k)\mathcal{T}^k_g(\ell\lambda_q)^{m_0}(\tau^a/\tau^c)\delta_{q+1} \ \text{ for } \ 0\leq r+k \leq N- m_0-3, \\
        &||\mathcal{A}^\pm_{\ell,j}\mathcal{L}_{\xi^g}^k[\mathcal{L}_{S}\Theta_j^g]||_{r+\alpha}\leq C' (C)^k\lambda_q^{r}L_{\mathcal{A}}(r+k)\mathcal{T}^k_g(\ell\lambda_q)^{m_0}(\tau^a/\tau^c)(1/\tau^c)\delta_{q+1}\ \text{ for } \ 0\leq r+k \leq N- m_0-4 \\
    \end{split}
\end{equation*}
where $C, \ C'$ depend on all the parameters but not on $a$ and are uniform in $k, \ r, \ j$. 

\noindent To bound $T_4$ we need to sum up these estimates over $k$, we will split the sum according to the definitions of $L,\ L_{\mathcal{A}}$  in \eqref{T4g} and the number of derivatives $r$ to isolate the bad range in which the losses $\bar L, \ \bar{\bar L}$ appear. We proceed as in the proof of Lemma \ref{estimatesfieldsringg} and skip the details, see \eqref{f(r)} in particular.

\noindent We begin with the pure derivatives. By means of Proposition \ref{czstuff} to deal with $\mathcal{R}\curl$ and the fact that there are at most two non-zero $\Theta^g_j$ at each time, we deduce:
\begin{equation}\label{transportg1}
    \begin{split}
        &||T_4||_{r+\alpha}\\
        &\lesssim \sup_j\sum_{k=0}^{k_0^g}\frac{1}{(k+1)!}||\mathcal{R}\curl \mathcal{L}_{\xi^{g}}^k[\mathcal{L}_{S}\Theta_j^{g}||_{r+\alpha}\lesssim\sup_j\sum_{k=0}^{k_0^g}\frac{1}{(k+1)!}||\mathcal{L}_{\xi^{g}}^k[\mathcal{L}_{S}\Theta_j^{g}||_{r+\alpha}\\
        &\lesssim \sum_{k=0}^{k_0^g}\frac{1}{(k+1)!}(C)^k\lambda_q^{r}L(r+k)\mathcal{T}^k_g(\ell\lambda_q)^{m_0}(\tau^a/\tau^c)\delta_{q+1}\\
        &\lesssim \lambda_q^{r}(\ell\lambda_q)^{m_0}(\tau^a/\tau^c)\delta_{q+1}\left[\sum_{k+r\leq \underline{r}-2, \ 0\leq k\leq k_0^g}\frac{\left(C\mathcal{T}_g\right)^k}{(k+1)!}\right]\\
        &\lesssim \lambda_q^{r}(\ell\lambda_q)^{m_0}(\tau^a/\tau^c)\delta_{q+1}\left[\bar L\sum_{k+r\geq \underline{r}-1,\ 0\leq k\leq k_0^g}\frac{\left(C\lambda_q^{(b-1)\gamma_\ell}\mathcal{T}_g\right)^k}{(k+1)!}\lambda_q^{[r-(\underline{r}-1)]^+(b-1)\gamma_\ell}\right]\\
        &= \lambda_q^{r}(\ell\lambda_q)^{m_0}(\tau^a/\tau^c)\delta_{q+1}\Biggr[\underbrace{\sum_{k=0}^{[\underline{r}-1-r]^+-1}\frac{\left(C\mathcal{T}_g\right)^k}{(k+1)!}}_{\lesssim1_{[\underline{r}-1-r]^+\geq 1}}+ \bar L\underbrace{\sum_{k=[\underline{r}-1-r]^+}^{ k_0^g}\frac{\left(C\lambda_q^{(b-1)\gamma_\ell}\mathcal{T}_g\right)^k}{(k+1)!}}_{\lesssim 1_{[\underline{r}-1-r]^+\leq k_0^g}\left(\lambda_q^{(b-1)\gamma_\ell}\mathcal{T}_g\right)^{[\underline{r}-1-r]^+}}\lambda_q^{[r-(\underline{r}-1)]^+(b-1)\gamma_\ell}\Biggr]\\
        &\lesssim  \lambda_q^{r}(\ell\lambda_q)^{m_0}(\tau^a/\tau^c)\delta_{q+1}\left[1_{[\underline{r}-1-r]^+\geq 1}+1_{[\underline{r}-1-r]^+\leq k_0^g} \lambda_q^{[r-(\underline{r}-1)]^+(b-1)\gamma_\ell}\left(\lambda_q^{(b-1)\gamma_\ell}\mathcal{T}_g\right)^{[\underline{r}-1-r]^+}\bar L\right]\\
        &\lesssim \lambda_q^{r}\lambda_q^{[r-(\underline{r}-k_0^g-1)]^+(b-1)}(\ell\lambda_q)^{m_0}(\tau^a/\tau^c)\delta_{q+1}
    \end{split}
\end{equation}
for $0\leq r\leq N-m_0-k_0^g-3$. The implicit constant depends on $r$ and all the parameters, but not on $a$.

\noindent We repeat the argument for the Alfv\'en transport bound. The additional use of Lemmas \ref{standardmollgalbrun}, \ref{stabilitygalbrun} to correct the transport operator and Proposition \ref{czstuff} to deal with the commutator, gives:
\begin{equation*}
    \begin{split}
        &||\mathcal{A}^\pm T_4||_{r+\alpha}\\
        &\lesssim \sup_j\sum_{k=0}^{k_0^g}\frac{1}{(k+1)!} \left[||\mathcal{R}\curl \mathcal{A}^\pm_{\ell,j}\mathcal{L}_{\xi^{g}}^k\mathcal{L}_{S}\Theta_j^{g}||_{r+\alpha}+||[\mathcal{A}^\pm_{\ell,j},\mathcal{R}\curl] \mathcal{L}_{\xi^{g}}^k\mathcal{L}_{S}\Theta_j^{g}||_{r+\alpha}+||(z_q^\pm-z^\pm_{\ell,j})\cn\mathcal{R}\curl \mathcal{L}_{\xi^{g}}^k\mathcal{L}_{S}\Theta_j^{g}||_{r+\alpha}\right]\\
        &\lesssim \sup_j\sum_{k=0}^{k_0^g}\frac{1}{(k+1)!} \underbrace{\left[||z^\pm_{\ell,j}||_{r+1+\alpha}||\mathcal{L}_{\xi^{g}}^k\mathcal{L}_{S}\Theta_j^{g}||_{\alpha}+||z^\pm_{\ell,j}||_{1+\alpha}||\mathcal{L}_{\xi^{g}}^k\mathcal{L}_{S}\Theta_j^{g}||_{r+\alpha}\right]}_{\lambda_q^{r+1}\ell^{-\alpha}L(r+k+1)\mathcal{T}^k_g(\ell\lambda_q)^{m_0}(\tau^a/\tau^c)\delta_q^{1/2}\delta_{q+1}}\\
        &+ \sup_j\sum_{k=0}^{k_0^g}\frac{1}{(k+1)!}\underbrace{\left[||z_q^\pm-z^\pm_{\ell,j}||_{r+\alpha}||\mathcal{L}_{\xi^{g}}^k\mathcal{L}_{S}\Theta_j^{g}||_{1+\alpha}+||z_q^\pm-z^\pm_{\ell,j}||_{0}|| \mathcal{L}_{\xi^{g}}^k\mathcal{L}_{S}\Theta_j^{g}||_{r+1+\alpha}\right]}_{\lambda_q^{r+1}\ell^{-\alpha}L(r+k+1)\mathcal{T}^k_g(\ell\lambda_q)^{2m_0}(\tau^a/\tau^c)\delta_q^{1/2}\delta_{q+1}}\\
        &+ \sup_j\sum_{k=0}^{k_0^g}\frac{1}{(k+1)!}\underbrace{|| \mathcal{A}^\pm_{\ell,j}\mathcal{L}_{\xi^{g}}^k\mathcal{L}_{S}\Theta_j^{g}||_{r+\alpha}}_{\lambda_q^{r}L_{\mathcal{A}}(r+k)\mathcal{T}^k_g(\ell\lambda_q)^{m_0}(\tau^a/\tau^c)(1/\tau^c)\delta_{q+1}}\\
        &\lesssim\sum_{k=0}^{k_0^g}\frac{(C)^k}{(k+1)!}\lambda_q^{r}L_{\mathcal{A}}(r+k)\mathcal{T}^k_g(\ell\lambda_q)^{m_0}(\tau^a/\tau^c)(1/\tau^c)\delta_{q+1}\\
    \end{split}
\end{equation*}
for $0\leq r\leq N-m_0-k_0^g-4$, where we used that $L(r+1)\leq L_{\mathcal{A}}(r)$. 

\noindent We now argue as before and split the sum according to $L_{\mathcal{A}}$ and $r$. Given that our assumptions, see \eqref{constraintadmissibility}, in fact guarantee that $(\lambda_q^{(b-1)\gamma_\ell}\mathcal{T}_g)^{k_0^g}\bar{\bar L}\leq 1$, we conclude:
\begin{equation*}
    \begin{split}
        &||\mathcal{A}^\pm T_4||_{r+\alpha}\\
        &\lesssim \sum_{k=0}^{k_0^g}\frac{(C)^k}{(k+1)!}\lambda_q^{r}L_{\mathcal{A}}(r+k)\mathcal{T}^k_g(\ell\lambda_q)^{m_0}(\tau^a/\tau^c)(1/\tau^c)\delta_{q+1}\\
        &=\lambda_q^r(\ell\lambda_q)^{m_0}(\tau^a/\tau^c)(1/\tau^c)\delta_{q+1}\left[\sum_{k=0}^{[\underline{r}-2-r]^+-1}\frac{\left(C_r\mathcal{T}_g\right)^k}{(k+1)!}+ \bar{\bar L}\sum_{k=[\underline{r}-2-r]^+}^{ k_0^g}\frac{\left(C\lambda_q^{(b-1)\gamma_\ell}\mathcal{T}_g\right)^k}{(k+1)!}\lambda_q^{[r-(\underline{r}-2)]^+(b-1)\gamma_\ell}\right]\\
        &\lesssim \lambda_q^r(\ell\lambda_q)^{m_0}(\tau^a/\tau^c)(1/\tau^c)\delta_{q+1}\left[1_{[\underline{r}-2-r]^+\geq 1}+1_{[\underline{r}-2-r]^+\leq k_0^g} \lambda_q^{[r-(\underline{r}-2)]^+(b-1)\gamma_\ell}\left(\lambda_q^{(b-1)\gamma_\ell}\mathcal{T}_g\right)^{[\underline{r}-2-r]^+}\bar{\bar L
        }\right]\\
        &\lesssim \lambda_q^r\lambda_q^{[r-(\underline{r}-k_0^g-2)]^+(b-1)}(\ell\lambda_q)^{m_0}(\tau^a/\tau^c)(1/\tau^c)\delta_{q+1}
    \end{split}
\end{equation*}
for $0\leq r\leq N-m_0-k_0^g-4$. The implicit constants in the above bounds depend on $r$ and all the parameters, but not on $a$.


\noindent \textbf{Estimates on $T_{2,w}, \ T_{2,b},\ T_3, \ T_{5,w}, \ T_{5,b}$.} The proof of the following bounds is just a slight modification of the techniques and ideas contained in the proof of Lemma \ref{estimatesperturbationgalbrun} and above. We refer to that and state the final bounds directly. 
\begin{equation*}
    \begin{split}
        &||T_{2,w}||_{r+\alpha}, \ ||T_{2,b}||_{r+\alpha}\lesssim \lambda_{q}^r\lambda_q^{[r-(\underline{r}-1)]^+(b-1)\gamma_\ell}\mathcal{T}_g(\tau^a/\tau^c)\ell^{-\alpha}\delta_{q+1} \ \text{ for } r\geq 0,\\
        &||\mathcal{A}^\pm T_{2,w}||_{r+\alpha}, \ ||\mathcal{A}^\pm T_{2,b}||_{r+\alpha}\lesssim  \lambda_{q}^r\lambda_q^{[r-(\underline{r}-1)]^+(b-1)\gamma_\ell}\mathcal{T}_g(1/\tau^c)\ell^{-\alpha}\delta_{q+1} \ \text{ for } 0\leq r \leq N-m_0-1,\\
        &||T_{3}||_{r+\alpha}\lesssim \lambda_{q}^r\lambda_q^{[r-\underline{r}]^+(b-1)\gamma_\ell}(\tau^a/\tau^c)(\ell\lambda_q)^{m_0}\delta_{q+1} \ \text{ for } \ 0\leq r\leq N-k_0^g-m_0-2,\\
        &||\mathcal{A}^\pm T_{3}||_{r+\alpha}\lesssim\lambda_{q}^r\lambda_q^{[r-(\underline{r}-k_0^g-1)]^+(b-1)}(1/\tau^c)(\ell\lambda_q)^{m_0}\delta_{q+1} \ \text{ for } \ 0\leq r\leq N-k_0^g-m_0-3,\\
        &||T_{5,w}||_{r+\alpha}, \ ||T_{5,b}||_{r+\alpha}\lesssim \lambda_{q}^r\lambda_q^{[r-(\underline{r}-1)]^+(b-1)\gamma_\ell}\mathcal{T}_g(\tau^a/\tau^c)(\ell\lambda_q)^{m_0}\ell^{-\alpha}\delta_{q+1} \ \text{ for } \ 0\leq r\leq N-k_0^g-m_0-2,\\
        &||\mathcal{A}^\pm T_{5,w}||_{r+\alpha}, \ ||\mathcal{A}^\pm T_{5,b}||_{r+\alpha}\lesssim \lambda_{q}^r\lambda_q^{[r-(\underline{r}-k_0^g-1)]^+(b-1)}\mathcal{T}_g(1/\tau^c)(\ell\lambda_q)^{m_0}\ell^{-\alpha}\delta_{q+1} \ \text{ for } \ 0\leq r\leq N-k_0^g-m_0-3.\\
    \end{split}
\end{equation*}
where the implicit constants depend on $r$ and all the other parameters, but not on $a$. 

\begin{remark}The estimate for $T_5$ is somewhat different, indeed, the leading term in the sum ($k=1, i=0$), already has two Lie derivations acting and this is the reason for the $\underline{r}-1$ in the pure derivative loss function, moreover since the summation in $i$ goes up to $k_0^g-1$ in the loss function for the transport bound we get $\underline{r}-k_0^g-1$ and not $\underline{r}-k_0^g-2$, which would look more natural given the pure derivative bound. Similar considerations hold for $T_2$, but there all the terms are mollified, and we have no loss in the transport bound.
\end{remark}

\noindent \textbf{Conclusion.} Gathering the estimates for each of the $T_i$ we get:
\begin{equation*}
    \begin{split}
        ||R^{g,tr}||_{r+\alpha}&\lesssim ||T_1||_{r+\alpha}+||T_{2,a}||_{r+\alpha}+||T_{2,b}||_{r+\alpha}+||T_3||_{r+\alpha}+||T_4||_{r+\alpha}+||T_{5,a}||_{r+\alpha}+||T_{5,b}||_{r+\alpha}\\
        &\lesssim \lambda_q^{r}\lambda_q^{[r-(\underline{r}-k_0^g)]^+(b-1)\gamma_\ell}\mathcal{T}_g\ell^{-\alpha}\delta_{q+1}+\lambda_q^{r}\lambda_q^{[r-(\underline{r}-k_0^g-1)]^+(b-1)\gamma_\ell}\mathcal{T}_g(\tau^a/\tau^c)\ell^{-\alpha}\delta_{q+1}\\
        &+\lambda_q^{r}\lambda_q^{[r-(\underline{r}-k_0^g-1)]^+(b-1)\gamma_\ell}(\ell\lambda_q)^{m_0}(\tau^a/\tau^c)\delta_{q+1}+\lambda_q^{r}\lambda_q^{[r-(\underline{r}-k_0^g-1)]^+(b-1)}(\ell\lambda_q)^{m_0}(\tau^a/\tau^c)\delta_{q+1}\\
        &+\lambda_{q}^r\lambda_q^{[r-(\underline{r}-k_0^g-1)]^+(b-1)\gamma_\ell}\mathcal{T}_g(\tau^a/\tau^c)(\ell\lambda_q)^{m_0}\ell^{-\alpha}\delta_{q+1}\\
        &\lesssim \lambda_{q}^r\lambda_q^{[r-(\underline{r}-k_0^g-1)]^+(b-1)}\mathcal{T}_g\ell^{-\alpha}\delta_{q+1}
    \end{split}
\end{equation*}
for $0\leq r \leq N-m_0-k_0^g-3$. Here we used $(\ell\lambda_q)^{m_0}\leq \lambda_q/\lambda_{q+1}\leq \mathcal{T}_g$, see \eqref{admissibilityofloss}. 

\noindent Moreover, we have:
\begin{equation*}
    \begin{split}
    ||\mathcal{A}^\pm R^{g,tr}||_{r+\alpha}&\leq ||\mathcal{A}^\pm T_1||_{r+\alpha}+||\mathcal{A}^\pm T_{2,a}||_{r+\alpha}+||\mathcal{A}^\pm T_{2,b}||_{r+\alpha}\\
    &+||\mathcal{A}^\pm T_3||_{r+\alpha}+||\mathcal{A}^\pm T_4||_{r+\alpha}+||\mathcal{A}^\pm T_{5,a}||_{r+\alpha}+||\mathcal{A}^\pm T_{5,b}||_{r+\alpha}\\
    &\lesssim\lambda_q^{r}\lambda_q^{[r-(\underline{r}-k_0^g)]^+(b-1)\gamma_\ell}\mathcal{T}_g(1/\tau^a)\ell^{-\alpha}\delta_{q+1}+\lambda_{q}^r\lambda_q^{[r-(\underline{r}-k_0^g-1)]^+(b-1)\gamma_\ell}\mathcal{T}_g(1/\tau^c)\ell^{-\alpha}\delta_{q+1}\\
    &+\lambda_{q}^r\lambda_q^{[r-(\underline{r}-k_0^g-1)]^+(b-1)}(\ell\lambda_q)^{m_0}(1/\tau^c)\delta_{q+1}+\lambda_{q}^r\lambda_q^{[r-(\underline{r}-k_0^g-2)]^+(b-1)}(\ell\lambda_q)^{m_0}(\tau^a/\tau^c)(1/\tau^c)\delta_{q+1}\\
    &+\lambda_{q}^r\lambda_q^{[r-(\underline{r}-k_0^g-1)]^+(b-1)}\mathcal{T}_g(1/\tau^c)(\ell\lambda_q)^{m_0}\ell^{-\alpha}\delta_{q+1}\\
    &\lesssim \lambda_{q}^r\lambda_q^{[r-(\underline{r}-k_0^g-2)]^+(b-1)}\mathcal{T}_g(1/\tau^a)\ell^{-\alpha}\delta_{q+1}
    \end{split}
\end{equation*}
for $0\leq r \leq N-m_0-k_0^g-4$. 

\noindent The proof of the Lemma is complete upon noticing that our choice of parameters in \ref{choiceofparameters}, see \eqref{misc}, guarantees that
$$N-k_0^g-m_0-3\geq M$$ 
and thus the bounds hold for $0\leq r \leq M$ and $0\leq r \leq M-1$ respectively.
\end{proof}


\subsubsection{Estimates on the Nash, Mollification, Remainder and Quadratic Error Terms}

\begin{lemma}[Estimates on the Nash Error]\label{Rgna} Let $\underline{r}=M-m_0-5$. Under the choice of parameters \ref{choiceofparameters}, the following bounds hold:
\begin{equation*}
    \begin{split}
        &||R^{g,na}||_{r+\alpha}\lesssim \lambda_q^{r}\lambda_q^{[r-(\underline{r}-k_0^g)]^+(b-1)\gamma_\ell}\mathcal{T}_g(\tau^a/\tau^c)\ell^{-\alpha}\delta_{q+1} \ \text{ for } \ 0\leq r \leq N-k_0^g-m_0-2,\\
        &||\mathcal{A}^\pm R^{g,na}||_{r+\alpha}\lesssim \lambda_q^{r}\lambda_q^{[r-(\underline{r}-k_0^g-1)]^+(b-1)}\mathcal{T}_g(1/\tau^c)\ell^{-\alpha}\delta_{q+1} \ \text{ for } \ 0\leq r \leq N-k_0^g-m_0-3.
    \end{split} 
\end{equation*}
The implicit constants depend on $r$ and all the other parameters, but not on $a$. In particular, the bounds hold for $0\leq r\leq M$ and $0\leq r\leq M-1$, respectively.
\end{lemma}
\begin{proof}[Proof of Lemma \ref{Rgna}] The Reynolds stress coming from the Nash-type terms $R^{g,na}$ is given by:
\begin{equation*}
    \begin{split}
        R^{g,na}&=\sum_j2\mathcal{R}\ddiv\left[\Theta_{w,j}^h\times\nabla v_{\ell,j}-\Theta_{b,j}^h\times\nabla  B_{\ell,j}\right]^\top\\
        &=2\sum_j\sum_{k=1}^{k_0^g}\frac{(-1)^k}{(k+1)!}\mathcal{R}\ddiv\Biggr[\underbrace{\left[[\mathcal{L}_{\xi^{g}}^k(\partial_t+\mathcal{L}_{v_{\ell,j}})\Theta_j^{g}]\times \nabla v_{\ell,j}\right]^\top}_{T_{1,w}^k}-\underbrace{\left[[\mathcal{L}_{\xi^{g}}^k\mathcal{L}_{B_{\ell,j}}\Theta_j^{g}]\times \nabla B_{\ell,j}\right]^\top}_{T_{1,b}^k}\Biggr]\\
        &+2\sum_j\sum_{k=0}^{k_0^g}\frac{(-1)^k}{(k+1)!}\mathcal{R}\ddiv\Biggr[\underbrace{\left[[\mathcal{L}_{\xi^{g}}^k\mathcal{L}_{v_q-v_{\ell,j}}\Theta_j^{g}]\times \nabla v_{\ell,j}\right]^\top}_{T_{2,w}^k}-\underbrace{\left[[\mathcal{L}_{\xi^{g}}^k\mathcal{L}_{B_q-B_{\ell,j}}\Theta_j^{g}]\times \nabla B_{\ell,j}\right]^\top}_{T_{2,b}^k}\Biggr]\\
    \end{split}
\end{equation*}
We omit the dependence on $j$ in $T_{1,w}^k,\ T_{2,w}^k, \ T_{1,b}^k, \ T_{2,b}^k$, as the bounds will be uniform in that parameter.

\noindent We apply Lemma \ref{estimatesperturbationgalbrun} with $\bar N =N+k_0^g$ and Lemma \ref{estimatesperturbationgalbrun2}. We let $C, \ C'$ be the implicit constants in the statements, which, without loss of generality, we might assume are the same. The bounds for $R^{g,na}$ follow from those in \ref{estimatesperturbationgalbrun}, \ref{estimatesperturbationgalbrun2}, \ref{stabilitygalbrun} and Proposition \ref{czstuff} to deal with $\mathcal{R}\ddiv$. We won't mention this again. We will argue, as in the proof of Lemma \ref{Rgtr}, to handle the lossy derivatives in the transport of the non-mollified terms; we omit some details and refer to that proof. 

\noindent \textbf{Pure derivative bounds.} We have:
\begin{equation*}
    \begin{split}
        ||\mathcal{R}\ddiv T_{1,w}^k||_{r+\alpha}&\lesssim||\mathcal{L}_{\xi^{g}}^k(\partial_t+\mathcal{L}_{v_{\ell,j}})\Theta_j^{g}\times \nabla v_{\ell,j}||_{r+\alpha}\\
        &\lesssim ||\mathcal{L}_{\xi^{g}}^k(\partial_t+\mathcal{L}_{v_{\ell,j}})\Theta_j^{g}||_{r+\alpha}||v_{\ell,j}||_1+||\mathcal{L}_{\xi^{g}}^k(\partial_t+\mathcal{L}_{v_{\ell,j}})\Theta_j^{g}||_{0}||v_{\ell,j}||_{r+1+\alpha}\\
        &\lesssim (C)^k \lambda_q^{r}\lambda_q^{[r+k-\underline{r}]^+(b-1)\gamma_\ell}\mathcal{T}^k_g(\tau^a/\tau^c)\ell^{-\alpha}\delta_{q+1} \ \text{ for } \ 0\leq r+k\leq \bar N,
        \\
        ||\mathcal{R}\ddiv T_{2,w}^k||_{r+\alpha}&\lesssim ||\mathcal{L}_{\xi^{g}}^k\mathcal{L}_{v_q-v_{\ell,j}}\Theta_j^{g}\times \nabla v_{\ell,j}||_{r+\alpha}\\
        &\lesssim||\mathcal{L}_{\xi^{g}}^k\mathcal{L}_{v_q-v_{\ell,j}}\Theta_j^{g}||_{r+\alpha}||v_{\ell,j}||_{1}+||\mathcal{L}_{\xi^{g}}^k\mathcal{L}_{v_q-v_{\ell,j}}\Theta_j^{g}||_{0}||v_{\ell,j}||_{r+\alpha+1}\\
        &\lesssim (C)^k\lambda_q^r\lambda_q^{[r+k-\underline{r}]^+(b-1)\gamma_\ell}\mathcal{T}^k_g(\ell\lambda_q)^{m_0}(\tau^a/\tau^c)\ell^{-\alpha}\delta_{q+1} \ \text{ for } 0\leq r+k\leq N-m_0-2.\ 
    \end{split}
\end{equation*}
The bounds for the terms involving the magnetic field follow similarly.
\begin{equation*}
    \begin{split}
    ||\mathcal{R}\ddiv T_{1,b}^k||_{r+\alpha}&\lesssim (C)^k\lambda_q^{r}\lambda_q^{[r+k-\underline{r}]^+(b-1)\gamma_\ell}\mathcal{T}^k_g(\tau^a/\tau^c)\ell^{-\alpha}\delta_{q+1} \ \text{ for } \ 0\leq r+k\leq \bar N,
    \\
    ||\mathcal{R}\ddiv T_{2,b}^k||_{r+\alpha}&\lesssim (C)^k\lambda_q^r\lambda_q^{[r+k-\underline{r}]^+(b-1)\gamma_\ell}\mathcal{T}^k_g(\ell\lambda_q)^{m_0}(\tau^a/\tau^c)\ell^{-\alpha}\delta_{q+1} \ \text{ for } \ 0\leq r+k\leq N-m_0-2.
    \end{split}
\end{equation*}
The implicit constant and $C$ in the bounds above depend on $r$ and all the parameters, but not on $a, \ k$.

\noindent \textbf{Alfv\'en transport bounds.} Using Proposition \ref{czstuff} to deal with the commutator and Lemmas \ref{standardmollgalbrun}, \ref{stabilitygalbrun} to correct the transport operators, we deduce:
\begin{equation*}
    \begin{split}
        ||\mathcal{A}^\pm \mathcal{R}\ddiv T^k_{1,w}||_{r+\alpha}&\lesssim||[\mathcal{A}^\pm_{\ell,j},\mathcal{R}\ddiv]T^k_{1,w}||_{r+\alpha}+||(z^\pm_q-z^\pm_{\ell,j})\cn \mathcal{R}\ddiv T^k_{1,w}||_{r+\alpha}+||\mathcal{R}\ddiv\mathcal{A}^\pm_{\ell,j}T^k_{1,w}||_{r+\alpha}\\
        &\lesssim \underbrace{||z^\pm_{\ell,j}||_{r+1+\alpha}||T^k_{1,w}||_{\alpha}+||z^\pm_{\ell,j}||_{1+\alpha}||T^k_{1,w}||_{r+\alpha}}_{\lambda_q^{r+1} \ell^{-2\alpha}\mathcal{T}^k_g(\tau^a/\tau^c)\delta_q^{1/2}\delta_{q+1}}\\
        &+\underbrace{||z^\pm_q-z^\pm_{\ell,j}||_{r+\alpha}||T^k_{1,w}||_{1+\alpha}+||z^\pm_q-z^\pm_{\ell,j}||_{0}||T^k_{1,w}||_{r+1+\alpha}}_{\lambda_q^{r+1} \mathcal{T}^k_g(\ell\lambda_q)^{m_0}(\tau^a/\tau^c)\ell^{-2\alpha}\delta_q^{1/2}\delta_{q+1}}\\
        &+\underbrace{||[\mathcal{A}^\pm_{\ell,j}\mathcal{L}_{\xi^{g}}^k(\partial_t+\mathcal{L}_{v_{\ell,j}})\Theta_j^{g}]\times \nabla v_{\ell,j}||_{r+\alpha}}_{\lambda_q^{r}\mathcal{T}^k_g(1/\tau^c)\ell^{-\alpha}\delta_{q+1}}+\underbrace{||\mathcal{L}_{\xi^{g}}^k(\partial_t+\mathcal{L}_{v_{\ell,j}})\Theta_j^{g}\times \nabla \mathcal{A}^\pm_{\ell,j}v_{\ell,j}||_{r+\alpha}}_{\lambda_q^{r+1}\delta_q^{1/2}\mathcal{T}^k_g(\tau^a/\tau^c)\ell^{-\alpha}\delta_{q+1}}\\
        &+\underbrace{||\mathcal{L}_{\xi^{g}}^k(\partial_t+\mathcal{L}_{v_{\ell,j}})\Theta_j^{g}\times (\DD z^\pm_{\ell,j})^\top(\DD v_{\ell,j})^\top||_{r+\alpha}}_{\lambda_q^{r+1} \mathcal{T}^k(\tau^a/\tau^c)\delta_q^{1/2}\ell^{-\alpha}\delta_{q+1}}\\
        &\lesssim (C)^k\lambda_q^{r}\lambda_q^{[r+k-\underline{r}]^+(b-1)\gamma_\ell}\mathcal{T}^k_g(1/\tau^c)\ell^{-\alpha}\delta_{q+1}\\
    \end{split}
\end{equation*}
for $0\leq r\leq N-m_0-1$ and $0\leq r+k\leq \bar N-1$. Moreover, 
\begin{equation*}
    \begin{split}
        ||\mathcal{A}^\pm \mathcal{R}\ddiv T^k_{2,w}||_{r+\alpha}&\lesssim||[\mathcal{A}^\pm_{\ell,j},\mathcal{R}\ddiv]T^k_{2,w}||_{r+\alpha}+||(z^\pm_q-z^\pm_{\ell,j})\cn \mathcal{R}\ddiv T^k_{2,w}||_{r+\alpha}+||\mathcal{R}\ddiv\mathcal{A}^\pm_{\ell,j}T^k_{2,w}||_{r+\alpha}\\
        &\lesssim \underbrace{||z^\pm_{\ell,j}||_{r+1+\alpha}||T^k_{2,w}||_{\alpha}+||z^\pm_{\ell,j}||_{1+\alpha}||T^k_{2,w}||_{r+\alpha}}_{\lambda_q^{r+1}\delta_q^{1/2}\mathcal{T}^k_g\ell^{-2\alpha}(\ell\lambda_q)^{m_0}(\tau^a/\tau^c)\delta_{q+1}}\\
        &+\underbrace{||z^\pm_q-z^\pm_{\ell,j}||_{r+\alpha}||T^k_{2,w}||_{1+\alpha}+||z^\pm_q-z^\pm_{\ell,j}||_{0}||T^k_{2,w}||_{r+1+\alpha}}_{\lambda_q^{r+1}\delta_q^{1/2}\mathcal{T}^k_g(\ell\lambda_q)^{2m_0}(\tau^a/\tau^c)\ell^{-2\alpha}\delta_{q+1}}\\
        &+\underbrace{||[\mathcal{A}^\pm_{\ell,j}\mathcal{L}_{\xi^{g}}^k\mathcal{L}_{v_q-v_{\ell,j}}\Theta_j^{g}]\times \nabla v_{\ell,j}||_{r+\alpha}}_{\lambda_q^{r+1}\delta_q^{1/2}\mathcal{T}^k_g(\ell\lambda_q)^{m_0}(\tau^a/\tau^c)\ell^{-\alpha}\delta_{q+1}}+\underbrace{||\mathcal{L}_{\xi^{g}}^k\mathcal{L}_{v_q-v_{\ell,j}}\Theta_j^{g}\times \nabla \mathcal{A}^\pm_{\ell,j}v_{\ell,j}||_{r+\alpha}}_{\lambda_q^{r+1}\delta_q^{1/2}\mathcal{T}^k_g(\ell\lambda_q)^{m_0}(\tau^a/\tau^c)\ell^{-\alpha}\delta_{q+1}}\\
        &+\underbrace{||\mathcal{L}_{\xi^{g}}^k\mathcal{L}_{v_q-v_{\ell,j}}\Theta_j^{g}\times (\DD z^\pm_{\ell,j})^\top(\DD v_{\ell,j})^\top||_{r+\alpha}}_{\lambda_q^{r+1}\delta_q^{1/2}\mathcal{T}^k_g(\ell\lambda_q)^{m_0}(\tau^a/\tau^c)\ell^{-\alpha}\delta_{q+1}}\\
        &\lesssim (C)^k\lambda_q^{r}L_{\mathcal{A}}(r+k)\mathcal{T}^k_g(\ell\lambda_q)^{m_0}(\tau^a/\tau^c)(1/\tau^c)\ell^{-\alpha}\delta_{q+1}
    \end{split}
\end{equation*}
for $0\leq r+k\leq N-m_0-3$. 

\noindent Similarly, one can show:
\begin{equation*}
    \begin{split}
    ||\mathcal{A}^\pm \mathcal{R}\ddiv T_{1,b}^k||_{r+\alpha}&\lesssim (C)^k\lambda_q^{r}\lambda_q^{[r+k-\underline{r}]^+(b-1)\gamma_\ell}\mathcal{T}^k_g(1/\tau^c)\ell^{-\alpha}\delta_{q+1}\ \text{ for } \ 0\leq r\leq N-m_0-1, \ 0\leq r+k\leq \bar N-1,\\
    ||\mathcal{A}^\pm \mathcal{R}\ddiv T_{2,b}^k||_{r+\alpha}&\lesssim(C)^k\lambda_q^{r}L_{\mathcal{A}}(r+k)\mathcal{T}^k_g(\ell\lambda_q)^{m_0}(\tau^a/\tau^c)(1/\tau^c)\ell^{-\alpha}\delta_{q+1} \ \text{ for } \ 0\leq r+k\leq N-m_0-3.\\
    \end{split}
\end{equation*}
The implicit constants and $C$ in the above depend on $r$ and the other parameters, but not on $a, \ k$.

\noindent \textbf{Conclusion.} To complete the proof we need to sum up these estimates over $k$, we will split the sum according to the definitions of the loss functions $\lambda_q^{[r+k-\underline{r}]^+(b-1)\gamma_\ell},\ L_{\mathcal{A}}(r+k)$, see \eqref{LAgalbrun}, and the number of derivatives $r$ to isolate the bad range in which the losses $\lambda_q^{(b-1)\gamma_\ell}, \ \bar L$ appears. We proceed as in the proof of the Alfv\'en transport bound in Lemma \ref{estimatesfieldsringg}, see \eqref{f(r)}.

\noindent Collecting the bounds above together with the fact that at most two $\Theta_j^g$ are non-zero at the same time, we get:
\begin{equation}\label{nashg1}
    \begin{split}
        ||R^{g,na}||_{r+\alpha}&\leq \sum_{k=1}^{k_0^g}\frac{1}{(k+1)!}\left[||\mathcal{R}\ddiv T_{1,2}||_{r+\alpha}+||\mathcal{R}\ddiv T_{2,w}||_{r+\alpha}+||\mathcal{R}\ddiv T_{1,b}||_{r+\alpha}+||\mathcal{R}\ddiv T_{2,b}||_{r+\alpha}\right]\\
        &\lesssim\sum_{k=1}^{k_0^g}\frac{(C)^k}{(k+1)!}\lambda_q^{r}\lambda_q^{[r+k-\underline{r}]^+(b-1)\gamma_\ell}\mathcal{T}^k_g(\tau^a/\tau^c)\ell^{-\alpha}\delta_{q+1}\\
        &+\sum_{k=0}^{k_0^g}\frac{(C)^k}{(k+1)!}\lambda_q^{r}\lambda_q^{[r+k-\underline{r}]^+(b-1)\gamma_\ell}\mathcal{T}^k_g(\ell\lambda_q)^{m_0}(\tau^a/\tau^c)\ell^{-\alpha}\delta_{q+1}\\
        &\lesssim\lambda_q^{r}\lambda_q^{[r-\underline{r}]^+(b-1)\gamma_\ell}(\tau^a/\tau^c)\ell^{-\alpha}\delta_{q+1}\left[\sum_{k=1}^{[\underline{r}-r]^+}\frac{\left(C\mathcal{T}_g\right)^k}{(k+1)!}+\sum_{k=[\underline{r}-r]^++1}^{k_0^g}\frac{\left(C\lambda_q^{(b-1)\gamma_\ell}\mathcal{T}_g\right)^k}{(k+1)!}\right]\\
        &+\lambda_q^{r}\lambda_q^{[r-\underline{r}]^+(b-1)\gamma_\ell}(\ell\lambda_q)^{m_0}(\tau^a/\tau^c)\ell^{-\alpha}\delta_{q+1}\left[\sum_{k=0}^{[\underline{r}-r]^+-1}\frac{\left(C\mathcal{T}_g\right)^k}{(k+1)!}+\sum_{k=[\underline{r}-r]^+}^{k_0^g}\frac{\left(C\lambda_q^{(b-1)\gamma_\ell}\mathcal{T}_g\right)^k}{(k+1)!}\right]\\
        &\lesssim \lambda_q^{r}\mathcal{T}_g(\tau^a/\tau^c)\ell^{-\alpha}\delta_{q+1}\left[1_{[\underline{r}-r]^+\geq 1}+1_{[\underline{r}-r]^++1\leq k_0^g}\lambda_q^{[r-\underline{r}]^+(b-1)\gamma_\ell}\lambda_q^{(b-1)\gamma_\ell}\left(\lambda_q^{(b-1)\gamma_\ell}\mathcal{T}_g\right)^{[\underline{r}-r]^+}\right]\\
        &+ \lambda_q^{r}(\ell\lambda_q)^{m_0}(\tau^a/\tau^c)\ell^{-\alpha}\delta_{q+1}\left[1_{[\underline{r}-r]^+\geq 1}+1_{[\underline{r}-r]^+\leq k_0^g}\lambda_q^{[r-\underline{r}]^+(b-1)\gamma_\ell}\left(\lambda_q^{(b-1)\gamma_\ell}\mathcal{T}_g\right)^{[\underline{r}-r]^+}\right]\\
        &\lesssim \lambda_q^{r}\lambda_q^{[r-(\underline{r}-k_0^g)]^+(b-1)\gamma_\ell}\mathcal{T}_g(\tau^a/\tau^c)\ell^{-\alpha}\delta_{q+1}
    \end{split} 
\end{equation}
for $0\leq r\leq N-k_0^g-m_0-2$. The implicit constant depends on $r$ and all the other parameters, but not on $a$. Here we used that$(\ell\lambda_q)^{m_0}\leq \lambda_q/\lambda_{q+1}\leq \mathcal{T}_g$, see \eqref{constraintadmissibility}. 

\noindent We proceed to the Alfv\'en transport estimate. 
\begin{equation*}
    \begin{split}
        &||\mathcal{A}^\pm R^{g,na}||_{r+\alpha}\\
        &\leq \sum_{k=1}^{k_0^g}\frac{1}{(k+1)!}\left[||\mathcal{A}^\pm \mathcal{R}\ddiv T_{1,2}||_{r+\alpha}+||\mathcal{A}^\pm \mathcal{R}\ddiv T_{2,w}||_{r+\alpha}+||\mathcal{A}^\pm \mathcal{R}\ddiv T_{1,b}||_{r+\alpha}+||\mathcal{A}^\pm \mathcal{R}\ddiv T_{2,b}||_{r+\alpha}\right]\\\\
        &\lesssim \underbrace{\sum_{k=1}^{k_0^g}\frac{(C)^k}{(k+1)!}\left[ \lambda_q^{r}\lambda_q^{[r-\underline{r}]^+(b-1)\gamma_\ell}\mathcal{T}^k_g(1/\tau^c)\ell^{-\alpha}\delta_{q+1}\right]}_{S_1}+ \underbrace{\sum_{k=0}^{k_0^g}\frac{(C)^k}{(k+1)!}\left[\lambda_q^{r}L_{\mathcal{A}}(r+k)\mathcal{T}^k_g(\ell\lambda_q)^{m_0}(\tau^a/\tau^c)(1/\tau^c)\ell^{-\alpha}\delta_{q+1}\right]}_{S_2}\\
    \end{split} 
\end{equation*}
for $0\leq r\leq N-k_0^g-m_0-3$. 

\noindent We can bound $S_1$ as in \eqref{nashg1}, and deduce:
$$S_1\lesssim \lambda_q^{r}\lambda_q^{[r-(\underline{r}-k_0^g)]^+(b-1)\gamma_\ell}\mathcal{T}_g(1/\tau^c)\ell^{-\alpha}\delta_{q+1} \ \text{ for } \ 0\leq r\leq \bar N-k_0^g=N$$ 

\noindent To deal with the loss function $L_{\mathcal{A}}$ in $S_2$, we proceed as in \eqref{transportg1}; we omit the details.
\begin{equation*}
    \begin{split}
        S_2&= \lambda_q^{r}(\ell\lambda_q)^{m_0}(\tau^a/\tau^c)(1/\tau^c)\ell^{-\alpha}\delta_{q+1}\sum_{k=0}^{k_0^g}\frac{(C)^k}{(k+1)!}\left[L_{\mathcal{A}}(r+k)\mathcal{T}^k_g\right]\\
        &\lesssim \lambda_q^{r}(\ell\lambda_q)^{m_0}(\tau^a/\tau^c)(1/\tau^c)\ell^{-\alpha}\delta_{q+1}\Biggr[\underbrace{\sum_{k=0}^{[\underline{r}-1-r]^+-1}\frac{\left(C\mathcal{T}_g\right)^k}{(k+1)!}}_{\lesssim1_{[\underline{r}-1-r]^+\geq 1}\mathcal{T}_g}+ \bar L\underbrace{\sum_{k=[\underline{r}-1-r]^+}^{ k_0^g}\frac{\left(C\lambda_q^{(b-1)\gamma_\ell}\mathcal{T}_g\right)^k}{(k+1)!}}_{\lesssim 1_{[\underline{r}-1-r]^+\leq k_0^g}\left(\lambda_q^{(b-1)\gamma_\ell}\mathcal{T}_g\right)^{[\underline{r}-1-r]^+}}\lambda_q^{[r-(\underline{r}-1)]^+(b-1)\gamma_\ell}\Biggr]\\
        &\lesssim  \lambda_q^{r}\mathcal{T}_g(\ell\lambda_q)^{m_0}(\tau^a/\tau^c)(1/\tau^c)\ell^{-\alpha}\delta_{q+1}\left[1_{[\underline{r}-1-r]^+\geq 1}+1_{[\underline{r}-1-r]^+\leq k_0^g} \lambda_q^{[r-(\underline{r}-1)]^+(b-1)\gamma_\ell}\left(\lambda_q^{(b-1)\gamma_\ell}\mathcal{T}_g\right)^{[\underline{r}-1-r]^+}\bar L\right]\\
        &\lesssim \lambda_q^{r}\lambda_q^{[r-(\underline{r}-k_0^g-1)]^+(b-1)}(\ell\lambda_q)^{m_0}(\tau^a/\tau^c)(1/\tau^c)\ell^{-\alpha}\delta_{q+1}
    \end{split}
\end{equation*}
for $0\leq r\leq N-k_0^g-m_0-3$. 

\noindent We conclude that:
\begin{equation*}
    \begin{split}
        ||\mathcal{A}^\pm R^{g,na}||_{r+\alpha}&\lesssim \lambda_q^{r}\lambda_q^{[r-(\underline{r}-k_0^g)]^+(b-1)\gamma_\ell}\mathcal{T}_g(1/\tau^c)\ell^{-\alpha}\delta_{q+1}\\
        &+\lambda_q^{r}\lambda_q^{[r-(\underline{r}-k_0^g-1)]^+(b-1)}(\ell\lambda_q)^{m_0}(\tau^a/\tau^c)(1/\tau^c)\ell^{-\alpha}\delta_{q+1}\\
        &\lesssim \lambda_q^{r}\lambda_q^{[r-(\underline{r}-k_0^g-1)]^+(b-1)}\mathcal{T}_g(1/\tau^c)\ell^{-\alpha}\delta_{q+1}
    \end{split}
\end{equation*}
for $0\leq r\leq N-k_0^g-m_0-3$, where we used that$(\ell\lambda_q)^{m_0}\leq \lambda_q/\lambda_{q+1}\leq \mathcal{T}_g$, see \eqref{constraintadmissibility}.  The implicit constant depends on $r$ and all the other parameters, but not on $a$. 

\noindent The proof of the Lemma is complete upon noticing that according to our choice of parameters in \ref{choiceofparameters}, see \eqref{misc}, we have: 
$$N-k_0^g-m_0-2\geq M$$ 
and thus the bounds hold for $0\leq r \leq M$ and $0\leq r \leq M-1$ respectively.
\end{proof}

\begin{lemma}[Estimates on the Remainder Error]\label{Rgrem} Let $\underline{r}=M-m_0-6$. Under the choice of parameters \ref{choiceofparameters}, the following bounds hold:
\begin{equation*}
    \begin{split}
        ||R^{g,rem}||_{r+\alpha}&\lesssim \lambda_q^{r}\lambda_q^{[r-\underline{r}]^+(b-1)\gamma_\ell}(\tau^a/\tau^c) \delta_{q+1}\ \text{ for } 0\leq r \leq N-(k_0^g+1)-m_0-3,\\
        ||\mathcal{A}^\pm R^{g,rem}||_{r+\alpha}&\lesssim\lambda_q^{r}\lambda_q^{[r-(\underline{r}-1)]^+(b-1)\gamma_\ell}\lambda_q\delta_q^{1/2}(\tau^a/\tau^c)\delta_{q+1} \ \text{ for } \ 0\leq r \leq  N-(k_0^g+1)-m_0-4.
    \end{split}
\end{equation*}
The implicit constants depend on $r$ and all the other parameters, but not on $a$. In particular, the bounds hold for $0\leq r\leq M$ and $0\leq r\leq M-1$, respectively.
\end{lemma}
\begin{proof}[Proof of Lemma \ref{Rgrem}] The key ingredient in the proof is the fact that from the choice of parameters in \ref{choiceofparameters}, see \eqref{constraintadmissibility}, we have: 
$$\left(\lambda_q^{(b-1)\gamma_\ell}\mathcal{T}_g\right)^{k_0^g}\leq \delta_q.$$

\noindent According to the estimates in Lemmas \ref{remaindersg} with $k_0=k_0^g$ and \ref{estimatesperturbationgalbrun2} and Proposition \ref{czstuff} to deal with $\mathcal{R}\curl$ and $\mathcal{R}\ddiv$, we can bound:
\begin{equation*}
    \begin{split}
        ||R^{g,rem}||_{r+\alpha}&\leq||\mathcal{R}\curl[\partial_t\theta_w]||_{r+\alpha}+||\mathcal{R}\ddiv\left[v_q\otimes r_{w}+r_{w}\otimes v_q-B_q\otimes r_{b}-r_{b}\otimes B_q\right]||_{r+\alpha}\\
        &\lesssim||\partial_t\theta_w||_{r+\alpha}+||v_q\otimes r_{w}+r_{w}\otimes v_q-B_q\otimes r_{b}-r_{b}\otimes B_q||_{r+\alpha}\\
        &\lesssim||\partial_t\theta_w||_{r+\alpha}+||v_q||_{r+\alpha}||\theta_{w}||_{1}+||v_q||_{0}||\theta_{w}||_{r+1+\alpha}+||B_q||_{r+\alpha}||\theta_{b}||_{1}+||B_q||_{0}||\theta_{b}||_{r+1+\alpha}\\
        &\lesssim \lambda_q^{r+1}\lambda_q^{[r-\underline{r}]^+(b-1)\gamma_\ell}\left(\lambda_q^{(b-1)\gamma_\ell}\mathcal{T}_g\right)^{k_0^g+1}\ell^{-\alpha}\tau^a \delta_{q+1}\\
        &\leq \lambda_q^{r}\lambda_q^{[r-\underline{r}]^+(b-1)\gamma_\ell}\tau^a/\tau^c\delta_{q+1}
    \end{split}
\end{equation*}
for $0\leq r\leq N-(k_0^g+1)-m_0-3$. 

\noindent We estimate the transport derivatives using only first- and second-order pure time derivatives. We obtain:
\begin{equation*}
    \begin{split}
        ||\mathcal{A}^\pm R^{g,rem}||_{r+\alpha}&\leq||\partial_t R^{g,rem}||_{r+\alpha}+||z^\pm_q\cn R^{g,rem}||_{r+\alpha}\\
        &\lesssim ||\partial_t^2\theta_w||_{r+\alpha}+||[\partial_t v_q\otimes \curl \ \theta_{w}]^{sym}||_{r+\alpha}+||[v_q\otimes \curl\ \partial_t\theta_{w}]^{sym}||_{r+\alpha}\\
        &+||[\partial_t B_q\otimes \curl \ \theta_{b}]^{sym}||_{r+\alpha}+||[B_q\otimes \curl \ \partial_t \theta_{b}]^{sym}||_{r+\alpha}\\
        &+||z^\pm_q||_{r+\alpha}||R^{g,rem}||_{1}+||z^\pm_q||_{0}||R^{g,rem}||_{r+\alpha+1}\\
        &\lesssim\lambda_q^{r+2}\lambda_q^{[r-(\underline{r}-1)]^+(b-1)\gamma_\ell}\left(\lambda_q^{(b-1)\gamma_\ell}\mathcal{T}_g\right)^{k_0^g+1}\ell^{-\alpha}\tau^a \delta_{q+1}+\lambda_q^{r+1}\lambda_q^{[r-(\underline{r}-1)]^+(b-1)\gamma_\ell}(\tau^a/\tau^c)\ell^{-\alpha}\delta_q^{1/2}\delta_{q+1}\\
        &\leq \lambda_q^{r}\lambda_q^{[r-(\underline{r}-1)]^+(b-1)\gamma_\ell}\lambda_q\delta_q^{1/2}(\tau^a/\tau^c) \delta_{q+1}
    \end{split}
\end{equation*}
for $0\leq r\leq N-(k_0^g+1)-m_0-4$. 

\noindent Given our choice of parameters in \ref{choiceofparameters}, see \eqref{misc}, the bound above hold in particular for $0\leq r\leq M$ and $0\leq r\leq M-1$ respectively. Moreover, the implicit constants depend on $r$ and all the other parameters, but not on $a$.
\end{proof}


\begin{lemma}[Estimates on the Mollification Error]\label{Rgmoll} Let $\underline{r}=M-m_0-6$. Under the choice of parameters in \ref{choiceofparameters}, the following bounds hold:
\begin{equation*}
    \begin{split}
        ||R^{g,mo}||_r&\lesssim\lambda_q^{r}\lambda_q^{[r-\underline{r}]^+(b-1)}\frac{\lambda_q\delta_{q+1}}{\lambda_{q+1}} \ \text{ for } \ 0\leq r \leq M,\\
        ||\mathcal{A}^\pm R^{g,mo}||_r&\lesssim\lambda_q^{r}\lambda_q^{[r-(\underline{r}-k_0^g-1)]^+(b-1)}1/\tau^c\frac{\lambda_q\delta_{q+1}}{\lambda_{q+1}} \ \text{ for } \ 0\leq r \leq M-1.
    \end{split}
\end{equation*}
The implicit constants depend on $r$ and all the other parameters, but not on $a$.
\end{lemma}
\begin{proof}[Proof of Lemma \ref{Rgmoll}] During the proof, we will need the estimates in Lemma \ref{estimatesperturbationgalbrun2}. These come with the loss function:  
$$L_{g,\mathcal{A}}(r)=1_{r\leq \underline{r}-k_0^g-1}+1_{r\geq \underline{r}-k_0^g}\left(\frac{\lambda_{q+1}}{\lambda_{q}}\right)^{r-(\underline{r}-k_0^g-1)}=\lambda_q^{[r-(\underline{r}-k_0^g-1)]^+(b-1)}.$$
We will structure the estimates on $R_q-R_\ell, \ v_q-v_\ell, \ B_q-B_\ell$ in the same way. We now recall the definition of $R^{g,mo}$ and split it in the following way: 
\begin{equation*}
    \begin{split}
        R^{g,mo}&=\sum_j[\underbrace{(v_q-v_{\ell,j})\otimes (w_j^{g,p}+w_j^{h})}_{T_{1,w}}]^{sym}-[\underbrace{(B_q-B_{\ell,j})\otimes (b_j^{g,p}+b_j^{h})}_{T_{1,b}}]^{sym}+\underbrace{R_q-R_\ell}_{T_2},\\
    \end{split}
\end{equation*}
we omit the dependence on $j$ of $T_{1,w}, \ T_{1,b}$.

\noindent \textbf{Estimates on $T_2$.} For $0\leq r \leq M-m_0$, the required bounds are contained in Lemmas \ref{standardmollgalbrun}, \ref{stabilitygalbrun}. For $M-m_0\leq r \leq M$, we bound each term separately and trade a bad derivative for a good one to retain the smallness. We write this in a compact form as:
\begin{equation*}
    \begin{split}
        ||R_q-R_\ell||_r&\lesssim \lambda_q^{r}(\ell\lambda_q)^{m_0}\delta_{q+1}\left[1_{r\leq M-m_0}+1_{r\geq M-m_0+1}\frac{\lambda_{q+1}}{\lambda_{q}}\frac{\lambda_{q}}{\lambda_{q+1}(\ell\lambda_q)^{m_0}}\right]\\ 
        &\leq  \lambda_q^{r}\lambda_q^{[r-(M-m_0)]^+(b-1)}\frac{\lambda_q\delta_{q+1}}{\lambda_{q+1}} \ \text{ for } \ 0\leq r \leq M,\\
        ||\mathcal{A}^\pm(R_q-R_\ell)||_r&\lesssim \lambda_q^{r}(\ell\lambda_q)^{m_0}\lambda_q\delta_q^{1/2}\delta_{q+1}\left[1_{r\leq M-m_0-1}+1_{r\geq M-m_0}\frac{\lambda_{q+1}}{\lambda_{q}}\frac{\lambda_{q}}{\lambda_{q+1}(\ell\lambda_q)^{m_0}}\right]\\
        &\leq  \lambda_q^{r}\lambda_q^{[r-(M-m_0-1)]^+(b-1)}\lambda_q\delta_q^{1/2}\frac{\lambda_q\delta_{q+1}}{\lambda_{q+1}} \ \text{ for } \ 0\leq r \leq M-1.
    \end{split}
\end{equation*}

\noindent \textbf{Estimates on $T_{1,w}$.} From the bounds in Lemmas  \ref{standardmollgalbrun}, \ref{stabilitygalbrun}, \ref{estimatesfieldsringg}, \ref{estimatesfieldsbarg}, we deduce:
\begin{equation*}
    \begin{split}
        ||T_{1,w}||_r&\lesssim(||v_q-v_\ell||_r+||v_\ell-v_{\ell,j}||_{r})||w_j^{g,p}+w_j^{h}||_0+(||v_q-v_\ell||_0+||v_\ell-v_{\ell,j}||_{0})||w_j^{g,p}+w_j^{h}||_{r}\\
        &\lesssim \lambda_q^{r}\lambda_q^{[r-\underline{r}]^+(b-1)\gamma_\ell}(\ell\lambda_q)^{m_0}(\tau^a/\tau^c)\delta_{q+1}
    \end{split}
\end{equation*}
for $0\leq r\leq N-m_0$.

\noindent We move to the Alfv\'en transport estimates.  We first rewrite:
\begin{equation*}
    \begin{split}
        \mathcal{A}^\pm T_{1,w}&=\left[\mathcal{A}^\pm (v_q- v_\ell)+(z^\pm_q-z^\pm_{\ell,j})\cn(v_\ell-v_{\ell,j})+\mathcal{A}^\pm_{\ell,j}(v_\ell-v_{\ell,j})\right]\otimes (w_j^{g,p}+w_j^{h})\\
        &+(v_q-v_{\ell,j})\otimes\left[ \mathcal{A}^\pm_{\ell,j}(w_j^{g,p}+w_j^{h})+(z^\pm_q-z^\pm_{\ell,j})\cn (w_j^{g,p}+w_j^{h})\right]\\
    \end{split}
\end{equation*}
Reasoning as in $T_2$, using the bounds in Lemma \ref{standardmollgalbrun} and in the Iterative Assumptions \eqref{inductiveassumptionsgeneral}, we first write for $0\leq r\leq M-1$:
\begin{equation*}
    \begin{split}
        ||[\mathcal{A}^\pm (v_q- v_\ell)||_r&\lesssim \lambda_q^{r+1}(\ell\lambda_q)^{m_0}\delta_q\left[1_{r\leq M-m_0-1}+1_{r\geq M-m_0}\frac{\lambda_{q+1}}{\lambda_{q}}\frac{\lambda_{q}}{\lambda_{q+1}(\ell\lambda_q)^{m_0}}\right]\\
        &\leq\lambda_q^{r}\lambda_q^{[r-(M-m_0-1)]^+(b-1)}\frac{\lambda_q^2\delta_q}{\lambda_{q+1}}
    \end{split}
\end{equation*}
and with the additional use of the bounds in Lemmas \ref{stabilitygalbrun}, \ref{estimatesperturbationgalbrun}, \ref{estimatesperturbationgalbrun2}, we deduce:
\begin{equation*}
    \begin{split}
        ||\mathcal{A}^\pm T_{1,w}||_r&\lesssim \left\|\left[\mathcal{A}^\pm (v_q- v_\ell)+(z^\pm_q-z^\pm_{\ell,j})\cn(v_\ell-v_{\ell,j})+\mathcal{A}^\pm_{\ell,j}(v_\ell-v_{\ell,j})\right]\otimes (w_j^{g,p}+w_j^{h}) \right\|_r \\
        &+\left\|(v_q-v_{\ell,j})\otimes\left[ \mathcal{A}^\pm_{\ell,j}(w_j^{g,p}+w_j^{h})+(z^\pm_q-z^\pm_{\ell,j})\cn (w_j^{g,p}+w_j^{h})\right]\right\|_r\\
        &\lesssim\lambda_q^{r}\lambda_q^{[r-(M-m_0-1)]^+(b-1)}\frac{\lambda_q^2\ell^{-\alpha}\tau^a\delta_q\delta_{q+1}}{\lambda_{q+1}}\\
        &+\lambda_q^{r}\lambda_q^{[r-(\underline{r}-1)]^+(b-1)\gamma_\ell}\ell^{-\alpha}(\ell\lambda_q)^{m_0}\lambda_q\delta_q^{1/2}\delta_{q+1}+\lambda_q^{r}\lambda_q^{[r-(\underline{r}-k_0^g-1)]^+(b-1)}\mathcal{T}_g(\ell\lambda_q)^{m_0}(1/\tau^c)\delta_q\\\\
        &\lesssim\lambda_q^{r}\lambda_q^{[r-(\underline{r}-k_0^g-1)]^+(b-1)}1/\tau^c\frac{\lambda_q\delta_{q+1}}{\lambda_{q+1}}
    \end{split}
\end{equation*}
for $0\leq r \leq M-1$. Here we used the definition of $\mathcal{T}_g=\lambda_q^2\tau^a\tau^c\ell^{-\alpha}\delta_{q+1}$ and the fact that $M-1\leq N-k_0^g-m_0-3$, see \eqref{misc}.

\noindent Using the fact that at most two labels $j$ are active at each time and the fact that $T_{1,b}$ enjoys the same bounds as $T_{1,w}$, we conclude that:
\begin{equation*}
    \begin{split}
        ||R^{g,mo}||_r&\lesssim \lambda_q^{r}\lambda_q^{[r-\underline{r}]^+(b-1)\gamma_\ell}(\ell\lambda_q)^{m_0}(\tau^a/\tau^c)\delta_{q+1}\\
        &+\lambda_q^{r}\lambda_q^{[r-(M-m_0)]^+(b-1)}\frac{\lambda_q\delta_{q+1}}{\lambda_{q+1}}\\
        &\lesssim \lambda_q^{r}\lambda_q^{[r-\underline{r}]^+(b-1)}\frac{\lambda_q\delta_{q+1}}{\lambda_{q+1}} \ \text{ for } \ 0\leq r \leq M,
        \\
        ||\mathcal{A}^\pm R^{g,mo}||_r&\lesssim\lambda_q^{r}\lambda_q^{[r-(\underline{r}-k_0^g-1)]^+(b-1)}1/\tau^c\frac{\lambda_q\delta_{q+1}}{\lambda_{q+1}}\\
        &+\lambda_q^{r}\lambda_q^{[r-(M-m_0-1)]^+(b-1)}\lambda_q\delta_q^{1/2}\frac{\lambda_q\delta_{q+1}}{\lambda_{q+1}}\\
        &\lesssim \lambda_q^{r}\lambda_q^{[r-(\underline{r}-k_0^g-1)]^+(b-1)}1/\tau^c\frac{\lambda_q\delta_{q+1}}{\lambda_{q+1}} \ \text{ for } \ 0\leq r \leq M-1.\\
    \end{split}
\end{equation*}
The implicit constants depend on $r$ and all the other parameters, but not on $a$. 
\end{proof}

\begin{lemma}[Estimates on the Quadratic Error]\label{Rgqua} Let $\underline{r}=M-m_0-6$. Under the choice of parameters in \ref{choiceofparameters}, the following bounds hold:
\begin{equation*}
    \begin{split}
        ||R^{g,qua}||_{r}&\lesssim \lambda_q^{r}\lambda_q^{[r-(\underline{r}-k_0^g)]^+(b-1)\gamma_\ell}\left(\lambda_q\ell^{-\alpha}\tau^a\delta_{q+1}\right)^2 \ \text{ for } \ 0\leq r \leq M,\\
        ||\mathcal{A}^\pm R^{g,qua}||_{r}&\lesssim \lambda_q^{r}\lambda_q^{[r-(\underline{r}-k_0^g-1)]^+(b-1)}(1/\tau^a)\left(\lambda_q\ell^{-\alpha}\tau^a\delta_{q+1}\right)^2 \ \text{ for } \ 0\leq r \leq M-1.\\
    \end{split}
\end{equation*}
The implicit constants depend on $r$ and all the other parameters, but not on $a$.
\end{lemma}
\noindent The estimates follow immediately from its definition in \eqref{Rg}, the decomposition in \eqref{split1g} and the bounds Lemmas \ref{estimatesfieldsringg}, \ref{estimatesfieldsbarg}. We also remark that:
$$\left(\lambda_q\ell^{-\alpha}\tau^a\delta_{q+1}\right)^2\leq (\tau^a/\tau^c)^2(\delta_{q+1}/\delta_q)\delta_{q+1}.$$


\subsection{Checkpoint \texorpdfstring{$(v_q, \ B_q, \ R_q)\leadsto (\tilde v_q,\ \tilde B_q,\ \tilde R_q)$}{(vq, Bq, Rq) to (tilde v, tilde B, tilde Rq)}}\label{check}
Before proceeding with the Nash stage, we summarise the work and estimates shown so far. We denote the objects resulting from the Galbrun stage either with a superscript $g$ or with a tilde. The vector fields and pressure are perturbed in the following way:
$$(\tilde v_q,\ \tilde B_q, \ \tilde p_q)= (v_q+w^g, \ B_q+b^g, \ p_q+\pi^g)$$
while the right-hand side of \eqref{relaxedMHD} is now:
$$\tilde R_q:= R^g+\sum_I g_I^2 A_I.$$
According to the decomposition \eqref{Rq+1}, with $(w^p,b^p,\pi^p)$ and $\tilde A_I$ equal to zero, and the fact that our way of constructing the perturbation ensures that the Faraday-Ohm system \eqref{FH} is solved exactly, we have:
\begin{equation*}
    \begin{cases}
        \partial_t\tilde v_q+\tilde v_q\cn \tilde v_q-\tilde B_q\cn \tilde B_q+\nabla \tilde p_q= \ddiv \ \tilde R_q,\\
        \partial_t \tilde B_q +\curl [\tilde B_q\times \tilde v_q]=0,\\
        \ddiv \ \tilde v_q=\ddiv \ \tilde B_q=0.\\
    \end{cases}
\end{equation*}
We now gather the estimates on these intermediate objects.
\begin{prop}[Galbrun Stage]\label{recap} Set $\underline{r}=M-m_0-6$ and let $L_{g,\mathcal{A}}:\mathbb{N}_{\geq 0}\to \mathbb{R}_{\geq 1}$ be the admissible loss function from Lemma \ref{estimatesfieldsringg}. The following bounds hold:
\begin{subequations}
    \begin{align}
        & |\tilde v_q |,\ |\tilde B_q|\leq C_0(1-3\delta_{q+1}^{1/2}) \ \text{ and } \ |\tilde B_q|\geq c_0(1+2\delta_{q+1}^{1/2}) \ \text{ everywhere on } \ \mathbb{T}^3\times \mathbb{R},\\
        &||\partial_t^j\tilde v_q ||_{r},\ ||\partial_t^j\tilde B_q||_{r}\lesssim \lambda_q^{r+j}\lambda_{q}^{[r+j-\underline{r}]^+(b-1)\gamma_\ell}\delta_q^{1/2}  \ \text{ for } \ j=0,1,2, \ \& \  0\leq r \leq N-j \ \& \ (j,r)\neq (0,0),\\
        &||\tilde p ||_{r}\lesssim \lambda_q^{r}\lambda_{q}^{[r-\underline{r}]^+(b-1)\gamma_\ell}\delta_q \ \text{ for } \  1\leq r \leq M,\\
        &||\mathcal{\tilde A}^\pm\tilde v_q ||_{r},\ ||\mathcal{\tilde A}^\pm\tilde B_q||_{r}\lesssim \lambda_q^{r+1}L_{g,\mathcal{A}}(r)\delta_q   \ \text{ for }\  0\leq r \leq N-1,\\
        & ||\mathcal{\tilde A}^\pm\tilde p ||_{r}\lesssim \lambda_q^{r}\lambda_{q}^{[r-\underline{r}]^+(b-1)\gamma_\ell}\ell^{-\alpha}(1/\tau^a)\delta_{q+1}  \ \text{ for } \  1\leq r \leq M-1,\\
        &||R^g||_{r}\lesssim \lambda_q^{r}\lambda_q^{[r-(\underline{r}-k_0^g)]^+(b-1)}\ell^{-\alpha}(\tau^a/\tau^c)\delta_{q+1}  \ \text{ for }\  0\leq r \leq M,\\
        &||\mathcal{\tilde A}^\pm R^g||_{r}\lesssim \lambda_q^{r}\lambda_q^{[r-(\underline{r}-k_0^g-1)]^+(b-1)}\ell^{-\alpha}(1/\tau^c)\delta_{q+1}   \ \text{ for }\  0\leq r \leq M-1.
    \end{align}
\end{subequations}
\end{prop}
\begin{remark}[New Pressure Bound]\label{newpressurebound}Here, we opt for a simplified notation and very suboptimal estimates. Note, however, the loss in the transport bound for the pressure, when compared to the one in \eqref{inductiveassumptionsgeneral}, this is because $\tilde R_q$ contains $\sum_I g_I^2 A_I$, which has no additional smallness, and $\ddiv \ddiv \sum_I g_I^2 A_I\neq 0$, so the pressure sees the fast oscillations we introduced in the Galbrun Stage. In the vector fields' bounds, the additional smallness of the Galbrun perturbation compensates for this fact.
\end{remark}
\begin{remark}[Total Good Derivative Loss $\underline{r}$] We start with $M$ and then lose: 
\begin{itemize}
    \item $m_0$ in the deep mollification.
    \item 1 in the local recorrection as we solve a first-order (non-linear) PDE.
    \item 1 for the construction of the charts in the Gr{\"o}nwall estimate for the flow.
    \item 2 in Galbrun's equation as it is of second-order.
    \item 1 in the correction terms of the Lie series $\mathcal{L}_{v_q-v_\ell}\Theta_j^g$.
    \item 1 passing from $\Theta_j^g\leadsto \xi_j^g$.
\end{itemize}
In the Alfv\'en transport and Reynold stress estimates, we lose an additional $k_0^g$ due to the Lie-Taylor expansion needed to read the geometric properties.
\end{remark}


\begin{proof}[Proof of Proposition \ref{recap}] We deal with the vector fields and the pressure first, and gather the bounds for $R^g$ at the end.

\noindent \textbf{Recap vector fields.} We will show the computations only for the velocity field; the bounds for the magnetic field follow similarly. The estimates in Lemmas \ref{estimatesfieldsringg} and \ref{estimatesfieldsbarg} give:
\begin{equation*}
    \begin{split}
        ||\partial_t^jw^g||_{r}&\leq ||\partial_t^j\bar w^g||_{r}+||\partial_t^j\mathring w^g||_r\\
        &\lesssim\lambda_q^{r+j}\lambda_q^{[r+j-\underline{r}]^+(b-1)\gamma_\ell}\ell^{-\alpha}\tau^a\lambda_q\delta_{q+1}+\lambda_q^{r+j}\lambda_q^{[r+j-\underline{r}]^+(b-1)\gamma_\ell}\mathcal{T}_g(\ell\lambda_q)^{m_0}\delta_q^{1/2}\\
        &\lesssim \lambda_q^{r+j}\lambda_q^{[r+j-\underline{r}]^+(b-1)\gamma_\ell}\ell^{-\alpha}\tau^a\lambda_q\delta_{q+1}
    \end{split}
\end{equation*}
for $0\leq r \leq N-j$, where we used that:
$$\ell^{-\alpha}\tau^a\lambda_q\delta_{q+1}=\lambda_{q+1}^{\alpha}\mathcal{T}_g\delta_q^{1/2}>  \mathcal{T}_g(\ell\lambda_q)^{m_0}\delta_q^{1/2}.$$
Moreover,
\begin{equation*}
    \begin{split}
        ||\mathcal{A}^\pm w^g||_{r}&\leq ||\mathcal{A}^\pm\bar w^g||_{r}+||\mathcal{A}^\pm \mathring w^g||_r\\
        &\lesssim \lambda_q^{r}\lambda_q^{[r-(\underline{r}-1)]^+(b-1)\gamma_\ell}\ell^{-\alpha}\lambda_q\delta_{q+1}+\lambda_q^{r}L_{g,\mathcal{A}}(r)\mathcal{T}_g(\ell\lambda_q)^{m_0}(1/\tau^c)\delta_q^{1/2}\\
        &\lesssim\lambda_q^{r}L_{g,\mathcal{A}}(r)\ell^{-\alpha}\lambda_q\delta_{q+1}.
    \end{split}
\end{equation*}
From this, we also deduce:
\begin{equation*}
    \begin{split}
        ||\mathcal{\tilde A}^\pm w^g||_{r}&\leq ||\mathcal{A}^\pm w^g||_r+||(w^g\pm b^g)\cn w^g||_r\\
        &\lesssim \lambda_q^{r}L_{g,\mathcal{A}}(r)\ell^{-\alpha}\lambda_q\delta_{q+1}+\lambda_q^{r+1}\lambda_q^{[r-(\underline{r}-1)]^+(b-1)\gamma_\ell}(\lambda_q\ell^{-\alpha}\tau^a\delta_{q+1})^2\\
        &\lesssim \lambda_q^{r}L_{g,\mathcal{A}}(r)\ell^{-\alpha}\lambda_q\delta_{q+1}
    \end{split}
\end{equation*}
for $0\leq r \leq N-1$. The implicit constants depend on $r$ and all the parameters, but not on $a$.

\noindent Given our Iterative Assumptions \eqref{inductiveassumptionsgeneral} we conclude, for $j=0,1,2$, $0\leq r \leq N-j$ and $(j,r)\neq (0,0)$, that:
\begin{equation*}
    \begin{split}
        ||\partial_t^j\tilde v_q||_{r}&\leq ||\partial_t^jv_q||_{r}+||\partial_t^jw^g||_{r}\\
        &\lesssim \lambda_q^{r+j}\delta_q^{1/2}+\lambda_q^{[r+j-\underline{r}]^+(b-1)\gamma_\ell}\ell^{-\alpha}\tau^a\lambda_q\delta_{q+1}\\
        &\lesssim \lambda_q^{r+j}\lambda_q^{[r+j-\underline{r}]^+(b-1)\gamma_\ell}\delta_q^{1/2}
    \end{split}
\end{equation*}
For $r=j=0$, we have for some constant $C$ depending on the parameters but not on $a$:
$$||\tilde v_q||_0\leq ||v_q||_{0}+||w^g||_{0}\leq C_0(1-\delta_q^{1/2})+C\lambda_q\ell^{-\alpha}\tau^a\delta_{q+1}\leq C_0(1-3\delta_{q+1}^{1/2}),$$
where the last inequality follows by choosing $a$ sufficiently large, given all the other parameters, so that:
$$\delta_q^{1/2}/\delta_{q+1}^{1/2}\geq3+C/C_0(\tau^a/\tau^c)(\delta_{q+1}^{1/2}/\delta_q^{1/2}).$$
We now check the non-vanishing condition. For some constant $C$ depending on the parameters but not on $a$, we have:
\begin{equation*}
    \begin{split}
        |\tilde B_q|\geq |B_q|-|b^g|\geq c_0(1+\delta_q^{1/2})- C\lambda_q\ell^{-\alpha}\tau^a\delta_{q+1}\geq c_0(1+2\delta_{q+1}^{1/2}),
    \end{split}
\end{equation*}
where the last inequality follows again by choosing $a$ large enough. 

\noindent Note that by bridging the Alfv\'en transport estimate with the pure time derivative bound, we can write:
$$||\mathcal{A}^\pm v_q||_r\lesssim \lambda_q^{r+1}\delta_q\left[1_{r\leq M-1}+\bar L1_{r\geq M}\right] \ \text{ for } \ 0\leq r \leq N-1$$
where $\bar L=\delta_q^{1/2}$. By definition, see \eqref{precise}, we have:
$$1_{r\leq M-1}+\bar L1_{r\geq M}\leq L_{g,\mathcal{A}}(r)$$
and thus
\begin{equation*}
    \begin{split}
        ||\mathcal{\tilde A}^\pm \tilde v_q||_{r}&\leq ||\mathcal{A}^\pm v_q||_r+||\mathcal{\tilde A}^\pm w^g||_r+||(w^g\pm b^g)\cn v_q||_r\\
        &\lesssim \lambda_q^{r+1}\delta_q\left[1_{r\leq M-1}+\bar L1_{r\geq M}\right]+\lambda_q^{r}L_{g,\mathcal{A}}(r)\ell^{-\alpha}\lambda_q\delta_{q+1}+\lambda_q^{r}\lambda_q^{[r-\underline{r}]^+(b-1)\gamma_\ell}\ell^{-\alpha}\lambda_q\delta_q^{1/2}\tau^a\delta_{q+1}\\
        &\leq \lambda_q^{r+1}L_{g,\mathcal{A}}(r)\delta_q
    \end{split}
\end{equation*}
for $0\leq r\leq N-1$. The implicit constants depend on $r$ and all the other parameters, but not on $a$. Here we used that $2\beta\geq \gamma_{CZ}$, see \eqref{misc1}.

\noindent \textbf{Recap pressure.} The claimed estimates follow directly from the ones in the Iterative Assumptions \eqref{inductiveassumptionsgeneral} and Lemma \ref{estgalbrun}.
\begin{equation}
    \begin{split}
        ||\tilde p_q||_r&\leq ||p_q||_r+||\pi^g||_r\\
        &\lesssim \lambda_q^{r}\delta_q+\lambda_q^{r}\lambda_{q}^{[r-\underline{r}]^+(b-1)\gamma_\ell}\ell^{-\alpha}\delta_{q+1}\\
        &\lesssim\lambda_q^{r}\lambda_{q}^{[r-\underline{r}]^+(b-1)\gamma_\ell}\delta_q 
    \end{split}
\end{equation}
for $0\leq r \leq M$. Moreover,
\begin{equation*}
    \begin{split}
        ||\mathcal{\tilde A}^\pm\tilde p_q||_r&\lesssim ||\mathcal{A}^\pm p_q||_r+||\mathcal{A}^\pm \pi^g||_r+||w^g\pm b^g||_{r}||p_q+\pi^g||_{1}+||w^g\pm b^g||_{0}||p_q+\pi^g||_{r+1}\\
        &\lesssim \lambda_q^{r+1}\delta_{q}^{3/2}+\lambda_q^{r}\lambda_{q}^{[r-\underline{r}]^+(b-1)\gamma_\ell}\ell^{-\alpha}(1/\tau^a)\delta_{q+1}+\lambda_q^{r+2}\lambda_{q}^{[r-\underline{r}]^+(b-1)\gamma_\ell}\ell^{-\alpha}\tau^a\delta_q\delta_{q+1}\\
        &\lesssim\lambda_q^{r}\lambda_{q}^{[r-\underline{r}]^+(b-1)\gamma_\ell}\ell^{-\alpha}(1/\tau^a)\delta_{q+1}
    \end{split}
\end{equation*}
for $0\leq r \leq M-1$, where we used that thanks to our choice of parameters in \eqref{choiceofparameters}, see \eqref{gammaa}, we have: 
$$ \gamma_a=\frac{1}{2}(1-\beta) \text{ and } \beta < 1/5 \Longrightarrow 2\beta\leq \gamma_a \Longrightarrow\tau^c/\tau^a\delta_{q+1}/\delta_q\geq 1 \Longrightarrow\ell^{-\alpha}(1/\tau^a)\delta_{q+1}\geq \lambda_q\delta_q^{3/2}.$$


\noindent \textbf{Recap Reynold stress.} Collecting the bounds from Lemmas \ref{Rcut}, \ref{Rgtr}, \ref{Rgna}, \ref{Rgrem}, \ref{Rgmoll}, \ref{Rgqua}, we deduce:
\begin{equation*}
    \begin{split}
        ||R^g||_{r}&\leq ||R^{g,tr}||_r+||R^{cut}||_r+ ||R^{g,moll}||_r+ ||R^{g,na}||_r+ ||R^{g,rem}||_r+||R^{g,qua}||_r\\
        &\lesssim \lambda_q^{r}\lambda_q^{[r-(\underline{r}-k_0^g)]^+(b-1)}\mathcal{T}_g\ell^{-\alpha}\delta_{q+1}+\lambda_q^{r}\lambda_q^{[r-\underline{r}]^+(b-1)\gamma_\ell}(\tau^a/\tau^c)\ell^{-\alpha}\delta_{q+1}+\lambda_q^r\lambda_q^{[r-\underline{r}]^+(b-1)}\frac{\lambda_q\delta_{q+1}}{\lambda_{q+1}}\\
        &+\lambda_q^{r}\lambda_q^{[r-(\underline{r}-k_0^g)]^+(b-1)\gamma_\ell}\mathcal{T}_g(\tau^a/\tau^c)\ell^{-\alpha}\delta_{q+1}+\lambda_q^{r}\lambda_q^{[r-\underline{r}]^+(b-1)\gamma_\ell}(\tau^a/\tau^c)\delta_{q+1}\\
        &+\lambda_q^{r}\lambda_q^{[r-\underline{r}]^+(b-1)\gamma_\ell}(\tau^a/\tau^c)^2(\delta_{q+1}/\delta_q)\delta_{q+1}\\
        &\lesssim\lambda_q^{r}\lambda_q^{[r-(\underline{r}-k_0^g)]^+(b-1)}(\tau^a/\tau^c)\ell^{-\alpha}\delta_{q+1}\\
    \end{split}
\end{equation*}
for $0\leq r \leq M$. The implicit constant depends on $r$ and all the other parameters, but not on $a$. Here we used that:
$$\frac{\lambda_q}{\lambda_{q+1}}\leq \mathcal{T}^g=\lambda_q^2\tau^c\tau^a\ell^{-\alpha}\delta_{q+1}\leq (\tau^a/\tau^c)\lambda_{q+1}^{-\alpha}(\delta_{q+1}/\delta_q)\leq(\tau^a/\tau^c)$$
see \eqref{constraintadmissibility}, \eqref{Tg}.

\noindent From the same Lemmas, we also get:
\begin{equation*}
    \begin{split}
        ||\mathcal{A}^\pm R^g||_{r}&\leq ||\mathcal{A}^\pm R^{g,tr}||_r+||\mathcal{A}^\pm R^{cut}||_r+ ||\mathcal{A}^\pm R^{g,na}||_r+ ||\mathcal{A}^\pm R^{g,moll}||_r+ ||\mathcal{A}^\pm R^{g,rem}||_r+||\mathcal{A}^\pm R^{g,qua}||_r\\
        &\lesssim \lambda_q^{r}\lambda_q^{[r-(\underline{r}-k_0^g-1)]^+(b-1)}\mathcal{T}_g(1/\tau^a)\ell^{-\alpha}\delta_{q+1}+\lambda_q^{r}\lambda_q^{[r-\underline{r}]^+(b-1)\gamma_\ell}\ell^{-\alpha}(1/\tau^c)\delta_{q+1}\\
        &+\lambda_q^{r}\lambda_q^{[r-(\underline{r}-k_0^g-1)]^+(b-1)}\mathcal{T}_g(1/\tau^c)\ell^{-\alpha}\delta_{q+1}+ \lambda_q^r\lambda_q^{[r-(\underline{r}-k_0^g-1)]^+(b-1)}1/\tau^c\frac{\lambda_q\delta_{q+1}}{\lambda_{q+1}}\\
        &+\lambda_q^{r}\lambda_q^{[r-(\underline{r}-1)]^+(b-1)\gamma_\ell}\lambda_q\delta_q^{1/2}(\tau^a/\tau^c)\delta_{q+1}+\lambda_q^{r}\lambda_q^{[r-(\underline{r}-k_0^g-1)]^+(b-1)}(1/\tau^c)(\tau^a/\tau^c)(\delta_{q+1}/\delta_q)\delta_{q+1}\\
        &\lesssim \lambda_q^{r}\lambda_q^{[r-(\underline{r}-k_0^g-1)]^+(b-1)}\ell^{-\alpha}(1/\tau^c)\delta_{q+1}
    \end{split}
\end{equation*}
and we conclude that:
\begin{equation*}
        ||\mathcal{\tilde A}^\pm R^g||_{r}\leq ||\mathcal{A}^\pm R^g||_{r}+||(w^g\pm b^g)\cn R^g||_r\lesssim \lambda_q^{r}\lambda_q^{[r-(\underline{r}-k_0^g-1)]^+(b-1)}1/\tau^c\ell^{-\alpha}\delta_{q+1}
\end{equation*}
for $0\leq r \leq M-1$. The implicit constant depends on $r$ and all the other parameters, but not on $a$.
\end{proof}

\section{Nash Stage}\label{Nash}
We now construct the main part of the perturbation. This will erase the `well-prepared' Reynold stress from the Galbrun stage. Many of the techniques and ideas have been discussed in previous sections; we will provide only brief summaries.

\subsection{Space Mollification \texorpdfstring{$(\tilde v_q, \ \tilde B_q)\leadsto (\tilde v_{\ell,j},\ \tilde B_{\ell,j})$}{(tilde vq, tilde Bq) to (tilde v l,j, tilde B l,j)}}\label{secmollnash}
The vector fields $(\tilde B_q, \ \tilde v_q)$ resulting from the Galbrun stage contain terms which are not mollified, and thus we need to perform a second mollification to avoid a possible loss of derivatives. As before, since this procedure breaks the commutativity of the Alfv\'en transport operators, we will construct a secondary local recorrection. The ideas are the same as in Section \ref{moll1}, and we omit the proofs.

\noindent Recall the notation from Proposition \ref{recap}, we now set: 
\begin{equation*}
    \begin{cases}
        \tilde v_\ell=(\tilde v_q)_\ell,\\
        \tilde B_\ell=(\tilde B_q)_\ell,\\
        \tilde z^\pm_\ell = \tilde v_\ell\pm \tilde B_\ell
    \end{cases} \text{ and }
    \begin{cases}
        \tilde R_\ell=(\tilde R_q)_\ell,\\
        \tilde M_\ell=(\tilde B_\ell\times \tilde v_\ell)-(\tilde B_q\times \tilde v_q)_\ell,\\
        \tilde R_\ell^c=\left(\tilde v_\ell\otimes \tilde v_\ell-\tilde B_\ell\otimes \tilde B_\ell\right)-\left(\tilde v_q\otimes \tilde v_q-\tilde B_q\otimes \tilde B_q\right)_\ell\\
    \end{cases}
\end{equation*}
together with Alfv\'en transport operators:
\begin{equation}\label{operatorsn}
    \mathcal{\tilde A}^\pm=\partial_t+\tilde z^\pm_q \cn,\
        \mathcal{\tilde A}^\pm_\ell=\partial_t+\tilde z^\pm_\ell \cn, \
        \mathcal{\tilde A}^\pm_{\ell,j}=\partial_t+\tilde z^\pm_{\ell,j} \cn.
\end{equation}
We gather the estimates concerning the mollification in the following lemma.
\begin{lemma}[Mollification Estimates]\label{standardmollnash}Let $r\geq 0$ be an integer and $\underline{r}=M-m_0-k_0^g-6$. We have:
    \begin{subequations}
        \begin{align}
            &||\partial_t^j\tilde v_\ell||_r,||\partial_t^j\tilde B_\ell||_r\lesssim \lambda_q^{r+j}\lambda_q^{[r+j-\underline{r}]^+(b-1)\gamma_\ell}\delta_q^{1/2} \ \text{ for } j=0,1,2 \ \& \ (r,j)\neq (0,0),\\
            &||\mathcal{\tilde A}^\pm \tilde v_\ell||_r,\ ||\mathcal{\tilde A}^\pm \tilde B_\ell||_r\lesssim \lambda_q^{r+1}\lambda_q^{[r-(\underline{r}-1)]^+(b-1)\gamma_\ell}\delta_q \ \text{ for } \ 0\leq r\leq N-1,\\
            &||\tilde p_\ell||_r\lesssim \lambda_q^{r}\lambda_q^{[r-\underline{r}]^+(b-1)\gamma_\ell}\delta_{q} \ \text{ for } \ r\geq 1, \\
            &||\mathcal{\tilde A}^\pm \tilde p_\ell||_r\lesssim \lambda_q^{r}\lambda_q^{[r-(\underline{r}-1)]^+(b-1)\gamma_\ell}\ell^{-\alpha}(1/\tau^a)\delta_{q+1} \ \text{ for } \ 1\leq r\leq N-1,\\
            &||\tilde R_\ell||\lesssim \lambda_q^{r}\lambda_q^{[r-\underline{r}]^+(b-1)\gamma_\ell}\delta_{q+1}, \\
            &||\mathcal{\tilde A}^\pm \tilde R_\ell||_r\lesssim \lambda_q^{r}\lambda_q^{[r-(\underline{r}-1)]^+(b-1)\gamma_\ell}1/\tau^a\delta_{q+1} \ \text{ for } \ 0\leq r\leq N-1.
        \end{align}
    \end{subequations}
    Moreover, let $L_{p,\mathcal{A}}:\mathbb{N}_{\geq 0}\to \mathbb{R}_{\geq 1}$ be the admissible loss function in \eqref{lossparameters}, the following error bounds hold:
    \begin{subequations}
        \begin{align}
            &||\partial_t^j(\tilde v_q-\tilde v_\ell)||_r\lesssim \lambda_q^{r+j}\lambda_q^{[r+j-(\underline{r}-m_0)]^+(b-1)\gamma_\ell}(\ell\lambda_q)^{m_0}\delta_q^{1/2} \ \text{ for } j=0,1,2, \ \&\ 0\leq r\leq N-m_0-j,\\
            &||\partial_t^j(\tilde B_q-\tilde B_\ell)||_r\lesssim \lambda_q^{r+j}\lambda_q^{[r+j-(\underline{r}-m_0)]^+(b-1)\gamma_\ell}(\ell\lambda_q)^{m_0}\delta_q^{1/2} \ \text{ for } j=0,1,2, \ \&\ 0\leq r\leq N-m_0-j,\\
            &||\mathcal{\tilde A}^\pm(\tilde v_q-\tilde v_\ell)||_r, ||\mathcal{\tilde A}^\pm(\tilde B_q-\tilde B_\ell)||_r\lesssim \lambda_q^{r+1}L_{p,\mathcal{A}}(r)(\ell\lambda_q)^{m_0}\delta_q \ \text{ for } \ 0\leq r\leq N-m_0-1,\\
            &||\tilde p_q-\tilde p_\ell||_r\lesssim \lambda_q^{r}\lambda_q^{[r-(\underline{r}-m_0)]^+(b-1)\gamma_\ell}(\ell\lambda_q)^{m_0}\delta_q \ \text{ for } \ 1\leq r\leq M-m_0,\\
            &||\mathcal{\tilde A}^\pm(\tilde p_q-\tilde p_\ell)||_r\lesssim \lambda_q^{r}\lambda_{q}^{[r-(\underline{r}-m_0-1)]^+(b-1)\gamma_\ell}\ell^{-\alpha}(1/\tau^a)\delta_{q+1} \ \text{ for } \ 0\leq r\leq M-m_0-1,\\
            &||\tilde R_q-\tilde R_\ell||_r\lesssim \lambda_q^{r}(\ell\lambda_q)^{m_0}\delta_{q+1}\ \text{ for } \ 0\leq r\leq \underline{r}-m_0,\\
            &||\mathcal{\tilde A}^\pm(\tilde R_q-\tilde R_\ell)||_r\lesssim \lambda_q^{r}(\ell\lambda_q)^{m_0}(1/\tau^a)\delta_{q+1}\ \text{ for } \ 0\leq r\leq \underline{r}-m_0-1.
        \end{align}
    \end{subequations}
    Finally, we have the following bounds for the forcings: 
    \begin{equation}
        \begin{split}
            &||\tilde R_\ell^c||_r,\ ||\tilde M_\ell||_r\lesssim \lambda_q^r\lambda_q^{[r-(\underline{r}-m_0)]^+(b-1)\gamma_\ell}(\ell\lambda_q)^{m_0}\delta_q,\\
            &||\mathcal{\tilde A}^\pm_\ell \tilde R_\ell^c||_r, \ ||\mathcal{\tilde A}^\pm_\ell M_\ell||_r \lesssim \lambda_q^{r+1}\lambda_q^{[r-(\underline{r}-m_0-1)]^+(b-1)\gamma_\ell}(\ell\lambda_q)^{m_0}\delta_q^{3/2}.\\
        \end{split}
    \end{equation}
    The implicit constants depend on $r$, all the parameters, and the specific choice of the mollification kernel, but not on $a$.
\end{lemma}

The bounds follow immediately from the ones in Propositions \ref{recap} and \ref{deepmollification}. It is important to remark here that since we proved `good properties' up to $r=\underline{r}$ derivatives, after mollification, we retain them for all $r$ upon paying $\ell^{-1}$ for each additional derivative. 

\noindent We are now ready to construct the corrections. Let $j$ be a fixed index of the time partition $\{\eta_j\}_j$ from Subsection \ref{moll1}. We now solve on $\mathbb{T}^3\times (t_j-\tau^c, t_j+\tau^c)$ the full relaxed MHD system: 
\begin{equation*}
    \begin{cases}
        \partial_t \tilde v_{\ell,j}+\ddiv[\tilde v_{\ell,j}\otimes \tilde v_{\ell,j}-\tilde B_{\ell,j}\otimes \tilde B_{\ell,j}]+\nabla \tilde p_{\ell,j}=\ddiv\ \tilde R_\ell,\\
        \partial_t\tilde B_{\ell,j}+\curl [\tilde B_{\ell,j}\times \tilde v_{\ell,j}]=0,\\
        \ddiv \ \tilde B_{\ell,j}=\ddiv \ \tilde v_{\ell,j}=0,\\
        (\tilde v_{\ell,j},\tilde B_{\ell,j})|_{t=t_j}=(\tilde v_\ell,\tilde B_\ell)|_{t=t_j}.
    \end{cases}
\end{equation*}

Existence and uniqueness follow from $||\tilde v_q||_1\tau^c,\ ||\tilde B_q||_1\tau^c\leq 1$. 

\begin{lemma}[Classical Solutions with Forcing]\label{classicaln} Let $\underline{r}=M-2m_0-k_0^g-7$, we have:
\begin{equation*}
        \begin{split}
            &||\tilde v_{\ell,j}||_{r+\alpha}+ ||\tilde B_{\ell,j}||_{r+\alpha}\lesssim \lambda_q^{[r-\underline{r}]^+(b-1)\gamma_\ell}\delta_q^{1/2} \ \text{ for }\ r\geq 1.\\
        \end{split}
    \end{equation*}
    The implicit constants depend on $r$ and all the parameters, but not on $a$.
\end{lemma}

We gather the  estimates for: $$(\tilde \Delta^v, \ \tilde\Delta^B, \ \tilde \Delta^p)=(\tilde v_\ell-\tilde v_{\ell,j}, \ \tilde B_\ell-\tilde B_{\ell,j},\ \tilde p_\ell- \tilde p_{\ell,j}),$$ in the following Lemma.
\begin{lemma}[Stability Estimates]\label{stabilitynash}Let $j\in\mathbb{Z}$, $\sigma=0,1,2$, $r\geq0$ integer and $\underline{r}=M-2m_0-k_0^g-7$. We have the following stability estimates:
    \begin{subequations}
    \begin{align}
        &||\partial_t^\sigma (\tilde v_\ell-\tilde v_{\ell,j})||_{r+\alpha}, \ ||\partial_t^{\sigma}(\tilde B_\ell-\tilde B_{\ell,j})||_{r+\alpha}\lesssim \lambda_q^{r+\sigma+1}\lambda_q^{[r-(\underline{r}-\sigma)]^+(b-1)\gamma_\ell}\ell^{-\alpha}(\ell\lambda_q)^{m_0}\tau^c\delta_{q},\\
        &||\mathcal{\tilde A}^\pm _{\ell,j}(\tilde v_\ell-\tilde v_{\ell,j})||_{r+\alpha}, \ ||\mathcal{\tilde A}^\pm _{\ell,j}(\tilde B_\ell-\tilde B_{\ell,j})||_{r+\alpha}\lesssim \lambda_q^{r+1}\lambda_q^{[r-(\underline{r}-1)]^+(b-1)\gamma_\ell}\ell^{-\alpha}(\ell\lambda_q)^{m_0}\delta_q,\\
        &||\tilde p_\ell-\tilde p_{\ell,j}||_{r+\alpha}\lesssim \lambda_q^{r}\lambda_q^{[r-\underline{r}]^+(b-1)\gamma_\ell}\ell^{-\alpha}(\ell\lambda_q)^{m_0}\delta_{q} \ \text{ for } \ r\geq 1,\\
        &||\mathcal{\tilde A}^\pm _{\ell,j}\left(\tilde p_\ell-\tilde p_{\ell,j}\right)||_{r+\alpha}\lesssim \lambda_q^{r+1}\lambda_q^{[r-(\underline{r}-1)]^+(b-1)\gamma_\ell}(\ell\lambda_q)^{m_0}\ell^{-2\alpha}\delta_{q}^{3/2} \ \text{ for } \ r\geq 1.
    \end{align}
\end{subequations}
The implicit constants depend only on $r$ and all the other parameters, but not on $a$.
\end{lemma}

Note that there is a discrepancy in the Alfv\'en transport estimate for $\tilde p_q-\tilde p_\ell$ and $\tilde \Delta^p$. This can be explained by the fact that the right-hand side in the stability equation has transport estimates inherited from those of the vector fields on which we have the expected one, namely
$$\partial_t\tilde \Delta^\pm+\tilde z_{\ell,j}^\mp\cn \tilde \Delta^\pm+\tilde \Delta^\mp\cn \tilde z_{\ell}^\pm +\nabla\tilde \Delta^p=\ddiv\ [\tilde R_\ell^{c,\pm}]$$
where $\tilde R_\ell^{c,\pm}=\tilde R_\ell^c\pm[\tilde M_\ell]_\times$, see \eqref{Deltaequation} for the details, while $\tilde p_q$ sees the new fast-in-time oscillating right-hand side we introduced in the Galbrun stage, see Remark \ref{newpressurebound}.

 We now gather the resulting estimates on the corrected vector fields and pressure.
\begin{lemma}[Local Correction]\label{localcorrn} Let $\underline{r}=M-2m_0-k_0^g-7$ and $j\in \mathbb{Z}$. The following bounds hold:
    \begin{subequations}
    \begin{align}
        &||\tilde v_{\ell,j}||_{0}, \ ||\tilde B_{\ell,j}||_{0}\leq C_0, \ \ \ \ |\tilde B_{\ell,j}|\geq c_0 \ \text{ everywhere on } \ \mathbb{T}^3\times \mathbb{R},\\
        &||\partial_t^\sigma \tilde v_{\ell,j}||_{r}, \ ||\partial_t^{\sigma}\tilde B_{\ell,j}||_{r}\lesssim \lambda_q^{r+\sigma}\lambda_q^{[r-(\underline{r}-\sigma)]^+(b-1)\gamma_\ell}\delta_{q}^{1/2} \ \text{ for } \ \sigma =0,1,2 \ \text{ and } r\geq 0,\\
        &||\mathcal{\tilde A}^\pm _{\ell,j}\tilde v_{\ell,j}||_{r+\alpha}, \ ||\mathcal{\tilde A}^\pm _{\ell,j}\tilde B_{\ell,j}||_{r}\lesssim \lambda_q^{r+1}\lambda_q^{[r-(\underline{r}-1)]^+(b-1)\gamma_\ell}\delta_q \ \text{ for } \ r\geq 0,\\
        &||\tilde p_{\ell,j}||_{r}\lesssim \lambda_q^r\lambda_{q}^{[r-\underline{r}]^+(b-1)\gamma_\ell}\delta_q \ \text{ for } \ r\geq 1,\\
        &||\mathcal{\tilde A}^\pm _{\ell,j}\tilde p_{\ell,j}||_{r}\lesssim \lambda_q^{r}\lambda_q^{[r-(\underline{r}-1)]^+(b-1)\gamma_\ell}\ell^{-\alpha}(1/\tau^a)\delta_{q+1} \ \text{ for } \ r\geq 1.
    \end{align}
\end{subequations}
The implicit constants depend on $r$ and all the other parameters, but not on $a$.
\end{lemma}


\subsection{Mollification along the Alfv\'en Directions}\label{timemollification} 
Using the result of Section \ref{secmollnash}, we now construct an operator regularising functions along the Alfv\'en directions. This kind of convolution along the flow is a generalisation of the construction first introduced by Isett in \cite{Isettflow} to the case of multiple, commuting, space-time flows. This will be useful later to avoid a loss of transport derivatives on the slow coefficients, see \eqref{slowcoeff}.

\noindent \textbf{Preliminaries.} Let $\mathcal{\tilde A}^\pm_{\ell,j}$ as in \eqref{operatorsn} but seen as vector fields $(\tilde z_{\ell,j}^{\pm},1)^\top$ in $\mathbb{T}^3\times B_{\tau^c}(t_j)\to \mathbb{T}^3\times \mathbb{R}$. We associate to these the Lagrangian flow maps $\mathbb{X}^{j,\pm}$, that is, we solve:
\begin{equation*}
    \begin{cases}
        \partial_s \mathbb{X}^{j,\pm}_s(x,t)=\mathcal{\tilde A}^\pm_{\ell,j}(\mathbb{X}^{j,\pm}_s(x,t)),\\
        \mathbb{X}^{j,\pm}_0(x,t)=(x,t).
    \end{cases}
\end{equation*}
The relation to the standard flow map is as follows. Let $X^{j,\pm}_{t,s}$ solve:
\begin{equation*}
    \begin{cases}
        \partial_s X^{j,\pm}_{s,t}(x)=\tilde z^{\pm}_{\ell,j}(X^{j,\pm}_{s,t}(x),s),\\
        X^{j,\pm}_{t,t}=\IId
    \end{cases}
\end{equation*}
then, by the uniqueness of the Cauchy problem, we have:
\begin{equation}\label{totvsnormal}
    \mathbb{X}^{j,\pm}_s(x,t)=(X^{j,\pm}_{t+s,t}(x),t+s).
\end{equation}
This identity will be handy in a moment. Now, recall that given any two vector fields $v, w$ with Lagrangian flow maps $X^v, X^w$, we have: $$[v,w]=0\iff X^{v}_s\circ X^{w}_{s'}=X^{w}_{s'}\circ X^{v}_s$$
see for example \cite[Chapter 10]{lee2013smooth} for a proof, taking $\partial_s,\partial_{s'}$ and then letting $s=0,s'=0$ respectively one finds the identities:
\begin{equation}\label{totalflowcommutation}
    \DD X^{v}_s [w]=w(X^{v}_s) \ \text{ and } \ \DD X^{w}_{s'} [v]=v(X^{w}_{s'}),
\end{equation}
we will apply these identities with $\mathcal{\tilde A}^\pm_{\ell,j},\ \mathbb{X}^{j,\pm}$ as: 
$$\frac{1}{2}[\mathcal{\tilde A}^+_{\ell,j},\mathcal{\tilde A}^-_{\ell,j}]=\partial_t \tilde B_{\ell,j}+[\tilde v_{\ell,j},\tilde B_{\ell,j}]=0, \ [\mathcal{\tilde A}^\pm_{\ell,j},\mathcal{\tilde A}^\pm_{\ell,j}]=0.$$

\noindent \textbf{Mollification along the Alfv\'en directions.} Let $F$ be a possibly tensor valued function: 
$$F:\mathbb{T}^3\times \mathbb{R}\to \mathbb{R}^N$$
with $\supp_{t} F\subset B_{2/3\tau^c}(t_j)$, and fix $\ell_t$ with $\ell_t< \frac{1}{6}\tau^c$. We define the mollification along the Alfv\'en directions of $F$ at scale $\ell_t$ to be: 
$$(F)_{\ell_t}^j:\mathbb{T}^3\times \mathbb{R}\to \mathbb{R}^N$$
given by
\begin{equation}\label{definitiontimemoll}
    \begin{split}
        (F)_{\ell_t}^j(x,t)=\int_{-\ell_t}^{\ell_t}\int_{-\ell_t}^{\ell_t}F(\mathbb{X}^{j,+}_s\circ \mathbb{X}^{j,-}_{s'}(x,t))\rho_{\ell_t}(s)\rho_{\ell_t}(s')\dd s'\dd s,
    \end{split}
\end{equation}
where $\rho\in C_c^\infty(\mathbb{R})$ is a convolution kernel which we require to satisfy: 
\begin{equation}
    \int_{-1}^1 \rho(s)\dd s=1, \ \ \ \ \rho\geq 0, \ \ \ \ \supp \ \rho \Subset (-1,1).
\end{equation}
The fact that $\supp_t F\subset B_{2/3\tau^c}(t_j)$ makes this object well-defined. To see this, we can argue as follows. The identity \eqref{totvsnormal} shows that any point $t\in B_{2/3\tau^c}(t_j)$ gets translated by $\mathbb{X}^{j,+}_s \circ\mathbb{X}^{j,-}_{s'}$ at most of $2\ell_t$, but since $\ell_t<  1/6\tau^c$ we get that $B_{2/3\tau^c}(t_j)\pm  2\ell_t\subset B_{\tau^c}(t_j)$, here both flows are well defined and we have good control on their derivatives as $||\tilde z^\pm_{\ell,j}||_1\tau^c\leq 1$ upon choosing $a$ sufficiently large given $b, \ \beta$.

We gather the properties of this space-time regularisation in the following Lemma.
\begin{lemma}[Mollification along the Alfv\'en directions] \label{mollalongtheflow} Let $\underline{r}=M-2m_0-k_0^g-8$, fix $0<\ell_t< 1/6\tau^c$ and let $f$ be a smooth function 
$$f:\mathbb{T}^3\times \mathbb{R}\to \mathbb{R} \ \text{ with } \ \supp_{t} f\subset B_{2/3\tau^c}(t_j).$$
The mollification along the Alfv\'en directions satisfies:
\begin{subequations}
    \begin{align}
        & ||(f)_{\ell_t}^j||_0\leq ||f||_0,\\
        &||(f)_{\ell_t}^j||_r\lesssim \lambda_q^{r-1}\lambda_q^{[r-1-\underline{r}](b-1)\gamma_\ell}||f||_1+||f||_r \ \text{ for } \ r\geq 1,\\
        &||(f)_{\ell_t}^j-f||_0\leq \ell_t\sum_{\sigma\in\{+, \ -\}}||\mathcal{\tilde A}_{\ell,j}^\sigma f||_0,\\
        &||(f)_{\ell_t}^j-f||_r\lesssim \ell_t\sum_{\sigma\in\{+, \ -\}}(||\mathcal{A}_{\ell,j}^\sigma f||_1\lambda_q^{r-1}\lambda_q^{[r-1-\underline{r}](b-1)\gamma_\ell}+||\mathcal{A}_{\ell,j}^\sigma f||_r)  \ \text{ for } \ r\geq 1.
    \end{align}
\end{subequations}
The implicit constants depend on $r$ and all the parameters, but not on $a$. Moreover, the following properties hold:
    \begin{enumerate}
    \item \textbf{Commutation:} $$\mathcal{\tilde A}^\pm_{\ell,j} (f)^j_{\ell_t}= (\mathcal{\tilde A}^\pm_{\ell,j}f)^j_{\ell_t}.$$
    \item \textbf{Integration by parts:}
    $$\int_{-\ell_t}^{\ell_t}\int_{-\ell_t}^{\ell_t}(\mathcal{\tilde A}^+_{\ell,j}f)(\mathbb{X}^{j,+}_s\circ \mathbb{X}^{j,-}_{s'})\rho_{\ell_t}(s)\rho_{\ell_t}(s')\dd s'\dd s=-\frac{1}{\ell_t}\int_{-\ell_t}^{\ell_t}\int_{-\ell_t}^{\ell_t}f(\mathbb{X}^{j,+}_s\circ \mathbb{X}^{j,-}_{s'})\rho_{\ell_t}'(s)\rho_{\ell_t}(s')\dd s'\dd s$$
    and an analogous identity holds for $\mathcal{\tilde A}^-_{\ell,j}$.
    \item \textbf{Support preservation:} $$\supp_{x,t} (f)_{\ell_t}^j\subset B_{2(C_0+1)\ell_t}(\supp_{x,t}\ f)$$
    where $C_0$ is as in \eqref{inductiveassumptionsgeneral}.
\end{enumerate}
In particular, the Lemma holds for $\ell_t=\tau^a$ upon choosing $a$ large enough.
\end{lemma}
\begin{remark}[Commutation + IBP] Combining the commutation and integration by parts properties allows one to control as many Alfv\'en transport derivatives as one desires upon paying $\ell_t^{-1}$ for each additional one, see \eqref{secondtransport3} for the explicit computation.
\end{remark}
\begin{remark}
    Although in the application we will use this operator only for scalar functions, the Lemma holds for the general vector-valued case.
\end{remark}
\begin{proof}[Proof of Lemma \ref{mollalongtheflow}.] The proof is an adaptation of the corresponding one in \cite{Isettflow}. The claims hold for $\ell_t=\tau^a$ since by definition we have $\tau^a=\lambda_q^{-(b-1)\gamma_a}\tau^c$ and given any $\gamma_a>0, \ b$ we can choose $a$ so large that $\tau^a<\frac{1}{6}\tau^c$.

\noindent \textbf{Commutation property.} We first use the properties in \eqref{totalflowcommutation}, to compute:
\begin{equation*}
    \begin{split}
        \mathcal{\tilde A}^+_{\ell,j}[f(\mathbb{X}^{j,+}_s\circ \mathbb{X}^{j,-}_{s'})]&=(\DD_{x,t}f)(\mathbb{X}^{j,+}_s\circ \mathbb{X}^{j,-}_{s'})[(\DD _{x,t}\mathbb{X}^{j,+}_s)(\mathbb{X}^{j,-}_{s'})[\DD_{x,t}\mathbb{X}^{j,-}_{s'}[\mathcal{\tilde A}^+_{\ell,j}]]]\\
        &=(\DD_{x,t}f)(\mathbb{X}^{j,+}_s\circ \mathbb{X}^{j,-}_{s'})[(\DD _{x,t}\mathbb{X}^{j,+}_s)(\mathbb{X}^{j,-}_{s'})[\mathcal{\tilde A}^+_{\ell,j}(\mathbb{X}^{j,-}_{s'})]\\
        &=(\DD_{x,t}f)(\mathbb{X}^{j,+}_s\circ \mathbb{X}^{j,-}_{s'})[(\DD _{x,t}\mathbb{X}^{j,+}_s[\mathcal{\tilde A}^+_{\ell,j}])(\mathbb{X}^{j,-}_{s'})]\\
        &=(\DD_{x,t}f)(\mathbb{X}^{j,+}_s\circ \mathbb{X}^{j,-}_{s'})[\mathcal{A}^+_{\ell,j}(\mathbb{X}^{j,+}_s\circ \mathbb{X}^{j,-}_{s'})]\\
        &=(\DD_{x,t}f) [\mathcal{\tilde A}^+_{\ell,j}])(\mathbb{X}^{j,+}_s\circ \mathbb{X}^{j,-}_{s'})\\
        &=(\mathcal{\tilde A}^+_{\ell,j}f)(\mathbb{X}^{j,+}_s\circ \mathbb{X}^{j,-}_{s'})
    \end{split}
\end{equation*}
similar calculations hold for $\mathcal{\tilde A}^-_{\ell,j}$ and we conclude that: 
\begin{equation}\label{firsttransport}
    \begin{split}
        &\mathcal{\tilde A}^\pm_{\ell,j} (f)^j_{\ell_t}=\int_{-\ell_t}^{\ell_t}\int_{-\ell_t}^{\ell_t}(\mathcal{\tilde A}^\pm_{\ell,j}f)(\mathbb{X}^{j,+}_s\circ \mathbb{X}^{j,-}_{s'})\rho_{\ell_t}(s)\rho_{\ell_t}(s')\dd s'\dd s= (\mathcal{\tilde A}^\pm_{\ell,j}f)^j_{\ell_t}.
    \end{split}
\end{equation}

\noindent \textbf{Integration by parts property.} We first compute:
\begin{equation}\label{secondtransport}
    \begin{split}
        &\int_{-\ell_t}^{\ell_t}\int_{-\ell_t}^{\ell_t}(\mathcal{\tilde A}^+_{\ell,j}f)(\mathbb{X}^{j,+}_s\circ \mathbb{X}^{j,-}_{s'})\rho_{\ell_t}(s)\rho_{\ell_t}(s')\dd s'\dd s\\
        &=\int_{-\ell_t}^{\ell_t}\int_{-\ell_t}^{\ell_t}\partial_s[(f)(\mathbb{X}^{j,+}_s)](\mathbb{X}^{j,-}_{s'})\rho_{\ell_t}(s)\rho_{\ell_t}(s')\dd s'\dd s\\
        &=\int_{-\ell_t}^{\ell_t}\int_{-\ell_t}^{\ell_t}\partial_s[(f)(\mathbb{X}^{j,+}_s\circ \mathbb{X}^{j,-}_{s'})]\rho_{\ell_t}(s)\rho_{\ell_t}(s')\dd s'\dd s\\
        &=\int_{-\ell_t}^{\ell_t}f(\mathbb{X}^{j,+}_{s}\circ \mathbb{X}^{j,-}_{s'})\underbrace{\rho_{\ell_t}(s)}_{=0}|_{s=-\ell_t}^{s=+\ell_t}\rho_{\ell_t}(s')\dd s'\\
        &-\int_{-\ell_t}^{\ell_t}\int_{-\ell_t}^{\ell_t} f(\mathbb{X}^{j,-}_{s'}\circ \mathbb{X}^{j,+}_s)\partial_{s}(\rho_{\ell_t})(s)\rho_{\ell_t}(s')\dd s'\dd s\\
        &=-\frac{1}{\ell_t}\int_{-\ell_t}^{\ell_t}\int_{-\ell_t}^{\ell_t}f(\mathbb{X}^{j,-}_{s'}\circ \mathbb{X}^{j,+}_s)\rho_{\ell_t}'(s)\rho_{\ell_t}(s')\dd s'\dd s.\\
    \end{split}
\end{equation}
Commuting the flows, the same calculations can be done for $\mathcal{\tilde A}^-_{\ell,j}$. In particular, we are free to take an additional Alfv\'en derivatives at a cost of $\ell_t^{-1}$. Indeed, we can use the identity \eqref{secondtransport} we just found and \eqref{totalflowcommutation} again to compute:
\begin{equation}\label{secondtransport3}
    \begin{split}
        \mathcal{\tilde A}^\sigma_{\ell,j}\mathcal{\tilde A}^\pm_{\ell,j}(f)^j_{\ell_t}&=\mathcal{\tilde A}^\pm_{\ell,j}\int_{-\ell_t}^{\ell_t}\int_{-\ell_t}^{\ell_t}(\mathcal{\tilde A}^\pm_{\ell,j}f)(\mathbb{X}^{j,+}_s\circ \mathbb{X}^{j,-}_{s'})\rho_{\ell_t}(s)\rho_{\ell_t}(s')\dd s'\dd s\\
        &=-\frac{1}{\ell_t}\mathcal{\tilde A}^\sigma_{\ell,j}\int_{-\ell_t}^{\ell_t}\int_{-\ell_t}^{\ell_t}[f(\mathbb{X}^{j,-}_{s'}\circ \mathbb{X}^{j,+}_s)]\rho_{\ell_t}'(s)\rho_{\ell_t}(s')\dd s'\dd s\\
        &=-\frac{1}{\ell_t}\int_{-\ell_t}^{\ell_t}\int_{-\ell_t}^{\ell_t}(\mathcal{\tilde A}^\sigma_{\ell,j}f)(\mathbb{X}^{j,-}_{s'}\circ \mathbb{X}^{j,+}_s)]\rho_{\ell_t}'(s)\rho_{\ell_t}(s')\dd s'\dd s
    \end{split}
\end{equation}
where $\sigma \in\{+,\ -\}$. Iterating, one can take any number of Alfv\'en transport derivatives. 

\noindent \textbf{Support property.} We first note that thanks to the identity \eqref{totvsnormal} for $|s'|\leq \ell_t$ we have:
$$\mathbb{X}^{j,-}_{s'}(\mathbb{T}^3\times B_{2/3\tau^c}(t_j))\subseteq \mathbb{T}^3\times B_{2/3\tau^c+\ell_t}(t_j)$$
and thus for $|s|,|s'|\leq \ell_t, \ |t-t_j| \leq 2/3\tau^c$ we deduce:
\begin{equation}\label{supportpm}
    \begin{split}
        |(\mathbb{X}^{j,+}_s\circ \mathbb{X}^{j,-}_{s'})-\IId_{x,t}|&\leq|(\mathbb{X}^{j,+}_s-\IId_{x,t})(\mathbb{X}^{j,-}_{s'})|+|\mathbb{X}^{j,-}_{s'}-\IId_{x,t}|\\
        &\leq\left[\sup_{(x,t)\in\mathbb{T}^3\times B_{2/3\tau^c+\ell_t}(t_j)}|\mathbb{X}^{j,+}_s-\IId_{x,t}|\right]+|\mathbb{X}^{j,-}_{s'}-\IId_{x,t}|\\
        &\leq |s|(||z^{+}_{\ell,j}||_0+1)+|s'|(||z^{+}_{\ell,j}||_0+1)\\
        &\leq 2(C_0+1)\ell_t\\
    \end{split}
\end{equation}
where we write $\IId_{x,t}$ to mean the identity map $(x,t)\mapsto (x,t)$. We can now compute:
\begin{equation*}
    \begin{split}
        \supp_{x,t}(f)^j_{\ell_t}&\subset\bigcup_{|s|,|s'|\leq \ell_t}\supp_{x,t}[(f)(\mathbb{X}^{j,+}_s\circ \mathbb{X}^{j,-}_{s'})]\\
        &=\bigcup_{|s|,|s'|\leq \ell_t}(\mathbb{X}^{j,+}_s\circ \mathbb{X}^{j,-}_{s'})^{-1}(\supp_{x,t}\ f)\\
        &=\bigcup_{|s|,|s'|\leq \ell_t}(\mathbb{X}^{j,-}_{s'})^{-1}\circ (\mathbb{X}^{j,+}_s)^{-1}(\supp_{x,t}\ f)\\
        &=\bigcup_{|s|,|s'|\leq \ell_t}\mathbb{X}^{j,-}_{-s'}\circ \mathbb{X}^{j,+}_{-s}(\supp_{x,t}\ f)\\
        &=\bigcup_{|s|,|s'|\leq \ell_t}\mathbb{X}^{j,-}_{s'}\circ \mathbb{X}^{j,+}_{s}(\supp_{x,t}\ f)\\
        &\subset B_{2(C_0+1)\ell_t}(\supp_{x,t}\ f)
    \end{split}
\end{equation*}
where we used that $\supp_t f\subset B_{2/3\tau^c}(t_j)$, and thus the bound in \eqref{supportpm} holds. The support property follows.

\noindent \textbf{Error computation.} Using the fact that the convolution kernel $\rho$ has average one, we rewrite:
\begin{equation}\label{molltime1}
    \begin{split}
        f_{\ell_t}-f&=\int_{-\ell_t}^{\ell_t}\int_{-\ell_t}^{\ell_t}\left[f(\mathbb{X}^{j,+}_s\circ \mathbb{X}^{j,-}_{s'}(x,t))-f(x,t)\right]\rho_{\ell_t}(s)\rho_{\ell_t}(s')\dd s'\dd s\\
        &=\int_{-\ell_t}^{\ell_t}\int_{-\ell_t}^{\ell_t}\left[f(\mathbb{X}^{j,+}_s\circ \mathbb{X}^{j,-}_{s'}(x,t))-f(\mathbb{X}^{j,+}_{s}(x,t))\right]\rho_{\ell_t}(s)\rho_{\ell_t}(s')\dd s'\dd s\\
        &+\int_{-\ell_t}^{\ell_t}\int_{-\ell_t}^{\ell_t}\left[f(\mathbb{X}^{j,+}_{s}(x,t))-f(x,t)\right]\rho_{\ell_t}(s)\rho_{\ell_t}(s')\dd s'\dd s\\
    \end{split}
\end{equation}
and we then express these differences as integrals of transport derivatives, the key point being that we want to compare the scale $\ell_t$ with the transport estimates of $f$, namely
\begin{equation*}
\begin{split}
    f(\mathbb{X}^{j,+}_s\circ \mathbb{X}^{j,-}_{s'}(x,t))-f(\mathbb{X}^{j,+}_{s}(x,t))&=\int_0^{s'}\partial_{s''}[f(\mathbb{X}^{j,+}_{s}\circ \mathbb{X}^{j,-}_{s''}(x,t))]\dd s''\\
    &=\int_0^{s'}\partial_{s''}[f(\mathbb{X}^{j,-}_{s''}\circ \mathbb{X}^{j,+}_{s}(x,t))]\dd s''\\
    &=\int_0^{s}\partial_{s''}[f(\mathbb{X}^{j,-}_{s''})]\circ \mathbb{X}^{j,+}_{s}(x,t))]\dd s''\\
    &=\int_0^{s'}(\mathcal{A}^-_{\ell,j}f)(\mathbb{X}^{j,-}_{s''}\circ\mathbb{X}^{j,+}_{s})\dd s''.\\
\end{split}
\end{equation*}
We conclude that for any $|s|, \ |s'|\leq \ell_t$, we have:
\begin{equation}\label{molltime2}
    \begin{split}
        ||f(\mathbb{X}^{j,+}_s\circ \mathbb{X}^{j,-}_{s'})-f(\mathbb{X}^{j,+}_{s})||_r\leq &\ell_t ||(\mathcal{\tilde A}_{\ell,j}^-f)(\mathbb{X}^{j,+}_s\circ \mathbb{X}^{j,-}_{s''})||_r\\
        &\overset{r=0}{\leq}\ell_t||\mathcal{\tilde A}_{\ell,j}^-f||_0\\
        &\overset{r\geq1}{\lesssim} \ell_t(||\mathcal{\tilde A}_{\ell,j}^-f||_1||\mathbb{X}^{j,+}_s\circ \mathbb{X}^{j,-}_{s''}||_r+||\mathcal{\tilde A}_{\ell,j}^-f||_r||\mathbb{X}^{j,+}_s\circ \mathbb{X}^{j,-}_{s''}||_1^r)\\
        &\lesssim \ell_t(||\mathcal{\tilde A}_{\ell,j}^-f||_1\lambda_q^{r-1}\lambda_q^{[r-1-\underline{r}](b-1)\gamma_\ell}+||\mathcal{\tilde A}_{\ell,j}^-f||_r).
    \end{split}
\end{equation}
Similarly, one can rewrite:
$$f(\mathbb{X}^{j,+}_{s}(x,t))-f(x,t)=\int_0^{s}(\mathcal{\tilde A}^+_{\ell,j}f)(\mathbb{X}^{j,+}_{s''}(x,t))\dd s''$$
and deduce the bound:
\begin{equation}\label{molltime3}
\begin{split}
    ||f(\mathbb{X}^{j,+}_{s})-f||_r
        &\overset{r=0}{\leq}\ell_t ||\mathcal{\tilde A}_{\ell,j}^+f||_0\\
        &\overset{r\geq 1}{\lesssim} \ell_t(||\mathcal{\tilde A}_{\ell,j}^+f||_1\lambda_q^{r-1}\lambda_q^{[r-1-\underline{r}](b-1)\gamma_\ell}+||\mathcal{\tilde A}_{\ell,j}^+f||_r).
\end{split}
\end{equation}
We now use \eqref{molltime2} and \eqref{molltime3} to bound \eqref{molltime1} and conclude that:
\begin{equation}\label{errormollflow1}
\begin{split}
    ||(f)_{\ell_t}^j-f||_r&\overset{r=0}{\leq}\ell_t \sum_{\sigma\in\{+, \ -\}}||\mathcal{A}_{\ell,j}^\sigma f||_0\\
    &\overset{r\geq 1}{\lesssim} \ell_t\sum_{\sigma\in\{+, \ -\}}(||\mathcal{A}_{\ell,j}^\sigma f||_1\lambda_q^{r-1}\lambda_q^{[r-1-\underline{r}](b-1)\gamma_\ell}+||\mathcal{A}_{\ell,j}^\sigma f||_r).\\
\end{split}
\end{equation}

\noindent \textbf{Derivative estimates.} We have:
\begin{equation}\label{derivativemollflow}
    \begin{split}
    ||(f)_{\ell_t}^j||_r&\leq \int_{-\ell_t}^{\ell_t}\int_{-\ell_t}^{\ell_t}||f(\mathbb{X}^{j,+}_s\circ \mathbb{X}^{j,-}_{s'})||_r\rho_{\ell_t}(s)\rho_{\ell_t}(s')\dd s'\dd s\\
    &\leq \sup_{|s|,|s'|\leq \ell_t} ||f(\mathbb{X}^{j,+}_s\circ \mathbb{X}^{j,-}_{s'})||_r\\
    &\overset{r=0}{\leq} ||f||_0\\
    &\overset{r\geq 1}{\lesssim} \sup_{|s|,|s'|\leq \ell_t}\left[||f||_1||\DD(\mathbb{X}^{j,+}_s\circ \mathbb{X}^{j,-}_{s'})||_{r-1}+||f||_r||\DD(\mathbb{X}^{j,+}_s\circ \mathbb{X}^{j,-}_{s'})||_{0}^r\right]\\
    &\lesssim \lambda_q^{r-1}\lambda_q^{[r-1-\underline{r}](b-1)\gamma_\ell}||f||_1+||f||_r.
    \end{split}
\end{equation}

\noindent In the above, the bounds on the flows' composition come from the following computation:
\begin{equation*}
    \begin{split}
        \sup_{|s|,|s'|\leq \ell_t}||\DD(\mathbb{X}^{j,+}_s\circ \mathbb{X}^{j,-}_{s'})||_{0}&\leq ||\DD\mathbb{X}^{j,+}_s-\IId_{x,t}||_0||\DD\mathbb{X}^{j,+}_{s'}-\IId_{x,t}||_0+||\DD\mathbb{X}^{j,+}_s-\IId_{x,t}||_0+||\DD\mathbb{X}^{j,+}_{s'}-\IId_{x,t}||_0+|\IId_{x,t}|_{op}\\
        &\leq 1+ C\left(\ell_t^2\lambda_q^2\delta_q+ \ell_t\lambda_q\delta_q^{1/2}\right)\\
        &\leq 2
    \end{split}
\end{equation*}
where we used the estimates in Lemma \ref{localcorrn} on $\tilde z^\pm_{\ell,j}$ and Proposition \ref{standardlagrangianestimate} to deduce the ones on $\mathbb{X}^\pm$ and the last inequality follows by choosing $a$ sufficiently large given the other parameters and our definition of $\tau^c$ in \eqref{tauc}. For higher derivatives, by means of the composition estimates in Proposition \ref{compestimates}, we get:
\begin{equation*}
    \begin{split}
        \sup_{|s|,|s'|\leq \ell_t}||\DD(\mathbb{X}^{j,+}_s\circ \mathbb{X}^{j,-}_{s'})||_{r}&\lesssim\sup_{|s|,|s'|\leq \ell_t}\left[||\DD(\mathbb{X}^{j,+}_s)\circ \mathbb{X}^{j,-}_{s'})||_r||\DD \mathbb{X}^{j,-}_{s'}||_0+||\DD\mathbb{X}^{j,+}_s||_0||\DD \mathbb{X}^{j,-}_{s'}||_r\right]\\
        &\lesssim \sup_{|s|,|s'|\leq \ell_t}\left[||\DD\mathbb{X}^{j,+}_s||_{1}||\mathbb{X}^{j,-}_{s'}||_r+||\DD\mathbb{X}^{j,+}_s||_{r}||\DD \mathbb{X}^{j,-}_{s'}||_0^r\right]+\lambda_q^{r}\lambda_q^{[r-\underline{r}](b-1)\gamma_\ell}\\
        &\lesssim \lambda_q^{r}\lambda_q^{[r-\underline{r}](b-1)\gamma_\ell}.
    \end{split}
\end{equation*}

\noindent The implicit constants in the above statements depend on $r$ and the other parameters, but not on $a$.
\end{proof}

\subsection{Update of the Decomposition}\label{geomconstrn}
The goal of this section is to `update' the geometric decomposition of Section \ref{geomconstr} to take care of the additional Galbrun perturbation. We construct new charts $\tilde \Psi_I$ adapted to $(\tilde v_{\ell,j},\tilde B_{\ell,j})$ and obtain a new decomposition $\tilde A_I$. We then show that the difference $\tilde A_I-A_I$ between the updated one and the old one is small, which corresponds to the flow error in \cite{GR}. This procedure is necessary to enhance the transport properties of the Nash LDF $\xi^p$, which will be constructed next. 

\noindent We use the same space-time cut-off functions and cubes constructed in Section \ref{geomconstr}, namely $\eta_j,\theta_{j'}$ and $Q_J$. Proceeding as in Section \ref{geomconstr} the commutation, $\mathcal{\tilde A}^+_{\ell,j}\mathcal{\tilde A}^-_{\ell,j}=\mathcal{\tilde A}^-_{\ell,j}\mathcal{\tilde A}^+_{\ell,j}$ the bounds in Lemma \ref{localcorrn}
and the results in Section \ref{chart}, namely Lemma \ref{chartconstr}, allow us to construct charts $\tilde \Psi_J:Q_J\to \mathbb{R}^3$ adapted to $(\tilde v_{\ell,j},\ \tilde B_{\ell,j})$. The conventions and notation used for the differential geometric objects can be found in Appendix \ref{diffgeom}.

\begin{lemma}[Chart Properties]\label{chartprop} Let $J=(j,j')$ a chart label as in \eqref{J} and $\underline{r}=M-2m_0-k_0^g-8$. For any fixed $t\in B_{\tau^c}(t_j)$ the map: 
$$\tilde \Psi_J(\cdot,t): B_{\tau^c}(x_{j'})\to \mathbb{R}^3$$
is a diffeomorphism onto its image, with estimates:
\begin{equation*}
    \begin{split}
        &||\DD \tilde \Psi_J-\IId||_0, \ ||(\DD\tilde \Psi_J)^{-1}-\IId||_0, \ ||\det[\DD\tilde \Psi_J]-1||_0\lesssim \lambda_{q+1}^{-\alpha},\\
        &||\partial_t^\sigma\DD \tilde \Psi_J||_r,||\partial_t^\sigma(\DD\tilde \Psi_J)^{-1}||_r\lesssim  \lambda_q^{r+\sigma}\lambda_q^{[r+\sigma-\underline{r}]^+(b-1)\gamma_\ell} \ \text{ for } \  \sigma=0,1,2,  \ r\geq 0\\
    \end{split}
\end{equation*}
where the H{\"o}lder norms are understood to be taken on $Q_J$ and the implicit constants depend on $r$ and all the other parameters, but not on $a$.

\noindent Moreover, for any scalar function $\varphi:\mathbb{R}\to \mathbb{R}$, $k\in e_3^\perp$ define 
$\varphi_{k,\lambda_{q+1}}(x)=\varphi(\lambda_{q+1}x\cdot k)$, which we regard as a function $\mathbb{R}^3\to\mathbb{R}$, and for any fixed unit vector $\nu$ which we think of as a 1-form, we have:
    \begin{equation*}
\begin{cases}
    (\partial_t+\mathcal{L}_{\tilde v_{\ell,j}}^1)\tilde \Psi_J^{1*}\nu=0,\\
    \mathcal{L}_{\tilde B_{\ell,j}}^1\tilde \Psi_J^{1*}\nu=0\\
\end{cases} \ \text{ and } \
    \begin{cases}
    (\partial_t+\tilde v_{\ell,j}\cn)\tilde \Psi_J^*\varphi_{k,\lambda_{q+1}}=0,\\
    \tilde B_{\ell,j}\cn \tilde \Psi_J^{*}\varphi_{k,\lambda_{q+1}}=0.\\
\end{cases}
\end{equation*}
In particular, if $(k,\nu,\zeta)\in \Lambda$ from Lemma \ref{geomlemma}, then 
$$\xi=\frac{1}{\lambda_{q+1}}\curl[\tilde \Psi_J^{1*}(\varphi_k\nu)]=\tilde \Psi_J^{2*}(\varphi_{k,\lambda_{q+1}}'\zeta)$$
is a well-defined vector field on $Q_J$ satisfying:
\begin{equation*}
    (\partial_t+\mathcal{L}_{\tilde v_{\ell,j}})\xi=0, \
    \mathcal{L}_{\tilde B_{\ell,j}}\xi=0, \
    \ddiv \ \xi=0.
\end{equation*}

\noindent Finally, given any positive integer $\bar N$, there exists $a$ sufficiently large so that:
$$||\tilde \Psi_J^{1*}(\varphi_{\lambda,k}\nu)||_r\lesssim \lambda_{q+1}^r \ \text{ for } \ 0\leq r \leq \bar N,$$
where the implicit constant depends only on $r$ but not on $C_0, \ a$ and the other parameters.
\end{lemma}

The Lemma is just a reformulation of \ref{chartconstr} in the current notation. For the final additional estimate, see \eqref{principaloscillation} in the proof of Lemma \ref{leadingtermsnash}.

\noindent The key takeaway here is that the geometric constraint in Lemma \ref{geomlemma} pairs with the fact that the charts map $B_{\ell,j}$ to a vector parallel to $e_3$. We will use this in the construction of the Nash LDF, see \eqref{ansatz}.


\noindent We can now update $A_I$ to $\tilde A_I$ and $a_I$ to $\tilde a_I$. A key difference here compared to the construction in the Galbrun stage is that we also mollify along the Alfv\'en direction the slow coefficients by means of the operator $(\dots)_{\ell_t}^{j}$ defined in \eqref{definitiontimemoll}. Namely, for $I=(\zeta, J)$, we set:
\begin{equation}\label{slowcoeff}
\begin{split}
\tilde a_I&=\delta_{q+1}^{1/2} \left(\chi_I\gamma_\zeta(\det[\DD \tilde\Psi_J]^{-2}\DD \tilde \Psi_J \left(\IId-R_\ell/\delta_{q+1}\right)\DD \tilde \Psi_J^T)\right)^{I_t}_{\ell_t},\\
    \tilde A_I&=\delta_{q+1}\tilde a_I^2\tilde\Psi_J^{2*}\zeta\otimes \tilde\Psi_J^{2*}\zeta.\\
\end{split}
\end{equation}
In the classical constructions, e.g. \cite{Isettflow} one mollifies only the Reynolds stress; this is not sufficient because we can provide transport estimates for the charts up to second-order, but in the transport part of the new Reynolds stress, we will need a third-order one (see $T_1$ in the proof of Lemma \ref{Rptr}). We now collect all the bounds relative to the slow coefficients. We remark that, despite the mollification along the Alfv\'en directions, the coefficients associated with charts of the same type still have disjoint supports.
\begin{lemma}[Estimates on the Decomposition]\label{slowcoeffestimates}Let $p$ be the type map in \eqref{typemap} and $\underline{r}=M-2m_0-k_0^g-8$. For any $I=(\zeta,j,j')\in \mathcal{I}$ the slow coefficients $\tilde a_I$ satisfy:
$$||\tilde a_I||_0\lesssim \delta_{q+1}^{1/2}$$
where the implicit constant is independent of $C_0$. Moreover,
\begin{equation*}
    \begin{split}
        &||\partial_t^j\tilde a_I||_r\lesssim \lambda_q^{r+j}\lambda_q^{[r+j-\underline{r}]^+(b-1)\gamma_\ell}\delta_{q+1}^{1/2} \ \text{ for } \  j=0,1,2 \ \text{ and } \ r\geq 0,\\
        &||(\mathcal{A}^\pm_{\ell,I})^{j+1}\tilde a_I||_r\lesssim 1/\tau^c\ell_t^{-j}\lambda_q^r \lambda_q^{[r-\underline{r}]^+(b-1)\gamma_\ell}\delta_{q+1}^{1/2} \ \text{ for } \ r,j\geq 0\\
    \end{split}
\end{equation*}
and we have stability estimates: 
\begin{equation}\label{stabilitycoeff}
    \begin{split}
        &||\tilde a_I-a_I||_r\lesssim \lambda_q^r \lambda_q^{[r-(\underline{r}-1)]^+(b-1)\gamma_\ell}[\ell_t/\tau^c+\mathcal{T}_g]\delta_{q+1}^{1/2} \ \text{ for } \  r\geq 0,\\
        &||\mathcal{\tilde A}^{\pm}(\tilde a_I-a_I)||_r\lesssim 1/\tau^c\lambda_q^r \lambda_q^{[r-\underline{r}]^+(b-1)\gamma_\ell}\delta_{q+1}^{1/2} \ \text{ for } \ 0\leq r \leq N-m_0\\
    \end{split}
\end{equation}
where the implicit constants depend on $r, \ j$ and all the other parameters, but not on $a$.

\noindent Finally, we have the following support localisation properties:
\begin{equation} \label{supportslow}
    \begin{split}
        &\supp_{x,t} \tilde a_I\Subset B_{\tau^c}(x_{j'})\times B_{\tau^c}(t_j),\\
        &\supp_{x,t} \tilde a_I\subset B_{C\ell_t}(\supp_{x,t}R_\ell),\\
        &p(I_1)=p(I_2) \Longrightarrow \supp_{x,t} \tilde a_{I_1} \cap \supp_{x,t} \tilde a_{I_2}=\emptyset
    \end{split}
\end{equation}
for any $I, \ I_1, \ I_2\in\mathcal{I}$ and some constant $C$ depending on the parameters but not on $a$.
\end{lemma} 
\begin{remark} The fact that we don't lose one good derivative in the transport bound comes from the fact that we lost one already in the construction of the chart, compared to the ones we have on the vector fields, but note that we lose one in the stability estimate \eqref{stabilitycoeff} i.e. two on the vector fields, as per Lemma \ref{chartstability}.
    
\end{remark}
\begin{remark}
    From the $j$ th-order Alfv\'en transport bound arguing as in \eqref{puretimeslowcoeff}, follows that:
    $$||\partial_t^j(\mathcal{A}^\pm_{\ell,I})\tilde a_I||_r\lesssim 1/\tau^c\lambda_q^{r+j} \lambda_q^{[r+j-\underline{r}]^+(b-1)\gamma_\ell}\delta_{q+1}^{1/2} \ \text{ for } \ j=0,1,2 \ \text{ and } \ r\geq 0.$$
    This will be needed later. 
\end{remark}

\begin{proof}[Proof of Lemma \ref{slowcoeffestimates}] We will apply the properties of the mollification along the Alfv\'en directions studied in Section \ref{timemollification}, see Lemma \ref{mollalongtheflow}, and the estimates in Lemmas \ref{standardmollnash} and \ref{chartprop}. In what follows, abusing notation, we always index $I$, even if the object in question does not depend on all its components. Moreover, the $C^0_tC^r_x$ norms are thought on the space-time domain $Q_I$ of the chart.

\noindent \textbf{Pure derivative estimates.} We first note that for $a$  sufficiently large $\DD \tilde\Psi_I\in B_{1/2}(\IId)$, moreover, for an implicit constant which depends on $r$ we have:
\begin{equation}\label{determinant}
    ||\det[\cdots]^{-2}|_{B_{1/2}(\IId)}||_r\lesssim 1
\end{equation}
the composition estimates in Proposition \ref{compestimates} for $r\geq 1$ give:
\begin{equation*}
    \begin{split}
        ||\det[\DD \tilde\Psi_I]^{-2}||_r\lesssim ||\det[\cdots]^{-2}|_{B_{1/2}(\IId)}||_1||\DD \tilde\Psi_I||_{r}+||\det[\cdots]^{-2}|_{B_{1/2}(\IId)}||_r||\DD \tilde\Psi_I||_{1}^r\lesssim \lambda_q^r \lambda_q^{[r-\underline{r}]^+(b-1)\gamma_\ell}
    \end{split}
\end{equation*}
and we deduce: 
\begin{equation*}
    \begin{split}
        ||\det[\DD \tilde\Psi_I]^{-2}\DD \tilde \Psi_I \left(\IId-R_\ell/\delta_{q+1}\right)\DD \tilde \Psi_I^T||_r&\lesssim ||\det[\DD \tilde\Psi_I]^{-2}||_r||\DD \tilde \Psi_I||^2_0||\IId-R_\ell/\delta_{q+1}||_0\\
        &+||\det[\DD \tilde\Psi_I]^{-2}||_0||\DD \tilde \Psi_I||^2_0||\IId-R_\ell/\delta_{q+1}||_r\\
        &+||\det[\DD \tilde\Psi_I]^{-2}||_0||\DD \tilde \Psi_I||_0||\DD \tilde \Psi_I||_r||\IId-R_\ell/\delta_{q+1}||_0\\
        &\lesssim \lambda_q^r \lambda_q^{[r-\underline{r}]^+(b-1)\gamma_\ell},
    \end{split}
\end{equation*}
from \eqref{derivativemollflow} we conclude:
\begin{equation}\label{pureslow}
    \begin{split}
        ||\tilde a_I||_r&=||\left(\chi_I\gamma_\zeta(\det[\DD \tilde\Psi_I]^{-2}\DD \tilde \Psi_I \left(\IId-R_\ell/\delta_{q+1}\right)\DD \tilde \Psi_I^T)\right)^{I_t}_{\ell_t}||_r\delta_{q+1}^{1/2} \\
        &\overset{r=0}{\leq}||\chi_I\gamma_\zeta(\det[\DD \tilde\Psi_I]^{-2}\DD \tilde \Psi_I \left(\IId-R_\ell/\delta_{q+1}\right)\DD \tilde \Psi_I^T)||_0\delta_{q+1}^{1/2}\\
        &\leq ||\gamma_\zeta|_{B_{1/2}(\IId)}||_0\delta_{q+1}^{1/2}\lesssim \delta_{q+1}^{1/2}\\
        &\overset{r\geq 1}{\lesssim}||\chi_I\gamma_\zeta(\det[\DD \tilde\Psi_I]^{-2}\DD \tilde \Psi_I \left(\IId-R_\ell/\delta_{q+1}\right)\DD \tilde \Psi_I^T)||_1\lambda_q^{r-1} \lambda_q^{[r-1-\underline{r}]^+(b-1)\gamma_\ell}\delta_{q+1}^{1/2}\\
        &+||\chi_I\gamma_\zeta(\det[\DD \tilde\Psi_I]^{-2}\DD \tilde \Psi_I \left(\IId-R_\ell/\delta_{q+1}\right)\DD \tilde \Psi_I^T)||_r\delta_{q+1}^{1/2}\\
        &\lesssim \lambda_q^{r} \lambda_q^{[r-\underline{r}]^+(b-1)\gamma_\ell}\delta_{q+1}^{1/2}
    \end{split}
\end{equation}
where the implicit constant for $r\geq 0$ depends on $r$ and all the other parameters. Note, however, that the bound for $r=0$ depends only on $\beta, \ b$ but not on $C_0$ and the other parameters.

\noindent \textbf{Alfv\'en transport estimates.} From \eqref{transportdeterminant} in Remark \ref{interestin} we have that: 
$$\mathcal{\tilde A}^\pm_{\ell,I}\det[\DD \tilde\Psi_I]=0.$$
The transport estimates for the charts follow from the Lie-transport properties in Lemma \ref{chartconstr}, namely
\begin{equation*}
    \begin{split}
        &\DD \tilde\Psi_I[B_{\ell,I}]=c e_3\Longrightarrow 0=\DD(\DD \tilde\Psi_I[B_{\ell,I}])=\DD \tilde\Psi_I\DD B_{\ell,I}+B_{\ell,I}\cn\DD \tilde\Psi_I,\\
        &\partial_t\tilde\Psi_I+\DD \tilde\Psi_I[v_{\ell,I}]=0\Longrightarrow 0=\DD(\partial_t\tilde\Psi_I+\DD \tilde\Psi_I[v_{\ell,I}])=(\partial_t+v_{\ell,I}\cn) \DD\tilde\Psi_I+\DD\tilde\Psi_I\DD v_{\ell,I}
    \end{split}
\end{equation*}
and thus since $z^\pm_{\ell,I}=v_{\ell,I}\pm B_{\ell,I}$ combination of the two gives:
$$\mathcal{\tilde A}^\pm_{\ell,I}\DD \tilde\Psi_I=(\partial_t+z^\pm_{\ell,I})\DD \tilde\Psi_I= \DD \tilde\Psi_I\DD z^{\pm}_{\ell,I},$$
by taking the transpose, we also see that: 
$$\mathcal{\tilde A}^\pm_{\ell,I}(\DD \tilde\Psi_I)^\top=(\mathcal{\tilde A}^\pm_{\ell,I}\DD \tilde\Psi_I)^\top=(\DD \tilde\Psi_I\DD z^{\pm}_{\ell,I})^\top=(\DD z^{\pm}_{\ell,I})^\top(\DD \tilde\Psi_I)^\top$$
and from the bounds in \ref{chartprop} and \ref{localcorrn} we deduce:
\begin{equation*}
    \begin{split}
        ||\mathcal{\tilde A}^\pm_{\ell,I}\DD \tilde\Psi_I||_r, \ ||\mathcal{\tilde A}^\pm_{\ell,I}(\DD \tilde\Psi_I)^\top||_r&\lesssim ||z^\pm_{\ell,I}||_{r+1}||\DD \tilde\Psi_I||_0+||z^\pm_{\ell,I}||_1||\DD \tilde\Psi_I||_r\\
        &\lesssim \lambda_q^{r+1}\lambda_q^{[r-\underline{r}]^+(b-1)\gamma_\ell}\delta_q^{1/2}.\\
    \end{split}
\end{equation*}
With this at hand and the bounds in Lemmas \ref{chartprop}, \ref{standardmollgalbrun} and \eqref{cutoffbounds}, we first obtain:
\begin{equation}\label{transportinside}
    \begin{split}
        &||\mathcal{\tilde A}^\pm_{\ell,I}\left[\chi_I\gamma_\zeta(\det[\DD \tilde\Psi_I]^{-2}\DD \tilde \Psi_I \left(\IId-R_\ell/\delta_{q+1}\right)\DD \tilde \Psi_I^T)\right]||_r\\
        &\lesssim ||(\mathcal{\tilde A}^\pm_{\ell,I}\chi_I)\det[\DD \tilde\Psi_I]^{-2}\gamma_\zeta\DD \tilde \Psi_I\left(\IId-R_\ell/\delta_{q+1}\right)\DD \tilde \Psi_I^T||_r\\
        &+||\chi_I\det[\DD \tilde\Psi_I]^{-2}\DD \gamma_\zeta(\mathcal{\tilde A}^\pm_{\ell,I}\DD \tilde \Psi_I) \left(\IId-R_\ell/\delta_{q+1}\right)\DD \tilde \Psi_I^T||_r\\
        &+||\chi_I\det[\DD \tilde\Psi_I]^{-2}\DD \gamma_\zeta\DD \tilde \Psi_I \left(\mathcal{\tilde A}^\pm_{\ell,I}R_\ell/\delta_{q+1}\right)\DD \tilde \Psi_I^T||_r\\
        &+||\chi_I\det[\DD \tilde\Psi_I]^{-2}\DD \gamma_\zeta\DD \tilde \Psi_I\left(\IId-R_\ell/\delta_{q+1}\right)(\mathcal{\tilde A}^\pm_{\ell,I}\DD \tilde \Psi_I^T)||_r\\
        &\lesssim \lambda_q^{r}\lambda_q^{[r-\underline{r}]^+(b-1)\gamma_\ell}1/\tau^c\\
    \end{split}
\end{equation}
and conclude, using \eqref{firsttransport} to commute the mollification along the flow and the transport operator and \eqref{derivativemollflow} to compute the $C^r$ norm, that:
\begin{equation*}
    \begin{split}
        ||\mathcal{\tilde A}^\pm_{\ell,I}\tilde a_I||_r&=||\mathcal{\tilde A}^\pm_{\ell,I}\left(\chi_I\gamma_\zeta(\det[\DD \tilde\Psi_I]^{-2}\DD \tilde \Psi_I \left(\IId-R_\ell/\delta_{q+1}\right)\DD \tilde \Psi_I^T)\right)^{I_t}_{\ell_t}||_r\delta_{q+1}^{1/2}\\
        &=||\left(\mathcal{\tilde A}^\pm_{\ell,I}\left[\chi_I\gamma_\zeta(\det[\DD \tilde\Psi_I]^{-2}\DD \tilde \Psi_I \left(\IId-R_\ell/\delta_{q+1}\right)\DD \tilde \Psi_I^T)\right]\right)^{I_t}_{\ell_t}||_r\delta_{q+1}^{1/2}\\
        &\lesssim ||\mathcal{\tilde A}^\pm_{\ell,I}\left[\chi_I\gamma_\zeta(\det[\DD \tilde\Psi_I]^{-2}\DD \tilde \Psi_I \left(\IId-R_\ell/\delta_{q+1}\right)\DD \tilde \Psi_I^T)\right]||_1\lambda_q^{r-1} \lambda_q^{[r-1-\underline{r}]^+(b-1)\gamma_\ell}\delta_{q+1}^{1/2}\\
        &+||\mathcal{\tilde A}^\pm_{\ell,I}\left[\chi_I\gamma_\zeta(\det[\DD \tilde\Psi_I]^{-2}\DD \tilde \Psi_I \left(\IId-R_\ell/\delta_{q+1}\right)\DD \tilde \Psi_I^T)\right]||_r\delta_{q+1}^{1/2}\\
        &\lesssim\lambda_q^{r}\lambda_q^{[r-\underline{r}]^+(b-1)\gamma_\ell}1/\tau^c\delta_{q+1}^{1/2}. 
    \end{split}
\end{equation*}
We can move any additional derivative to the kernel in the mollification along the flow, see \eqref{secondtransport}, and we conclude:
\begin{equation}\label{transportslownash}
    ||(\mathcal{\tilde A}^\pm_{\ell,I})^{j+1}\tilde a_I||_r\lesssim\lambda_q^{r}\lambda_q^{[r-\underline{r}]^+(b-1)\gamma_\ell}1/\tau^c\ell_t^{-j}\delta_{q+1}^{1/2} \ \text{ for } \ r,j\geq 0.
\end{equation}
The implicit constant depends on $r,j$ and all the other parameters, but not on $a$.

\noindent \textbf{Stability estimates.} According to Lemma \ref{chartstability} we have:
$$||\DD\tilde\Psi_I-\DD\Psi_I||_r\lesssim \lambda_q^{r}\lambda_q^{[r-(\underline{r}-1)]^+(b-1)\gamma_\ell}\mathcal{T}_g,$$
we now split:
\begin{equation*}
    \begin{split}
        \tilde a_I-a_I&=\underbrace{\delta_{q+1}^{1/2}\left(\chi_I\gamma_\zeta(\det[\DD \tilde\Psi_I]^{-2}\DD \tilde \Psi_I \left(\IId-R_\ell/\delta_{q+1}\right)\DD \tilde \Psi_I^T)\right)^{I_t}_{\ell_t}-\delta_{q+1}^{1/2}\chi_I\gamma_\zeta(\det[\DD \tilde\Psi_I]^{-2}\DD \tilde \Psi_I \left(\IId-R_\ell/\delta_{q+1}\right)\DD \tilde \Psi_I^T)}_{T_1}\\
        &+\underbrace{\delta_{q+1}^{1/2}\chi_I\gamma_\zeta(\det[\DD \tilde\Psi_I]^{-2}\DD \tilde \Psi_I \left(\IId-R_\ell/\delta_{q+1}\right)\DD \tilde \Psi_I^T)-\delta_{q+1}^{1/2}\chi_I\gamma_\zeta(\det[\DD \Psi_I]^2\DD \Psi_I \left(\IId-R_\ell/\delta_{q+1}\right)\DD \Psi_I^T)}_{T_2}
    \end{split}
\end{equation*}
The estimate on $T_1$ follows immediately from \eqref{errormollflow1} and \eqref{transportinside}:
$$||T_1||_r\lesssim  \lambda_q^{r}\lambda_q^{[r-\underline{r}]^+(b-1)\gamma_\ell}\ell_t/\tau^c\delta_{q+1}^{1/2}.$$
We now estimate $T_2$. Define for any matrices $A, M\in B_{1/2}(\IId)$:
$$f(A,M)=\gamma_\zeta(\det[A]^{-2}A MA^T),$$
from \eqref{determinant} and Lemma \ref{matrixdecomp} it follows that for some implicit constant which depends on $r$ we have
\begin{equation}\label{updatecoeff0}
    \begin{split}
        ||f|_{B_{1/2}(\IId)\times B_{1/2}(\IId)}||_r \lesssim 1.
    \end{split}
\end{equation}
We now compute:
\begin{equation}\label{updatecoeff}
    \begin{split}
        f(A_1,M)-f(A_0,M)&=\int_0^1\partial_s\left[f(sA_1+(1-s)A_0,M)\right] \dd s\\
        &=\int_0^1(\partial_Af)(sA_1+(1-s)A_0,M)[A_1-A_0]\dd s\\
    \end{split}
\end{equation}
and we use this with $A_1=\DD\tilde\Psi_I, \ A_2=\DD\Psi_I$ and $M=\IId-R_\ell/\delta_{q+1}$ to rewrite:
\begin{equation}\label{updatecoeff1}
    T_2=f(A_1,M)-f(A_0,M).
\end{equation}
From the bound in \eqref{updatecoeff0}, the composition estimates in Proposition \ref{compestimates} and the expression for \eqref{updatecoeff1} in \eqref{updatecoeff}, we obtain:
\begin{equation*}
    \begin{split}
        ||T_2||_r&\leq\int_0^1||(\partial_Af)(s\DD\tilde\Psi_I+(1-s)\DD\Psi_I,\IId-R_\ell/\delta_{q+1})[\DD\tilde\Psi_I-\DD\Psi_I]||_r\dd s\\
        &\lesssim ||\DD\tilde\Psi_I-\DD\Psi_I]||_0\int_0^1||(\partial_Af)(s\DD\tilde\Psi_I+(1-s)\DD\Psi_I,\IId-R_\ell/\delta_{q+1})||_r\dd s\\
        &+||\DD\tilde\Psi_I-\DD\Psi_I]||_r\int_0^1||(\partial_Af)(s\DD\tilde\Psi_I+(1-s)\DD\Psi_I,\IId-R_\ell/\delta_{q+1})||_0\dd s\\
        &\lesssim||\DD\tilde\Psi_I-\DD\Psi_I]||_0||f|_{B_{1/2}(\IId)\times B_{1/2}(\IId)}||_{2}(||\DD\tilde\Psi_I||_r+||\DD\Psi_I||_r+||R_\ell/\delta_{q+1}||_r)\\
        &+||\DD\tilde\Psi_I-\DD\Psi_I]||_0||f|_{B_{1/2}(\IId)\times B_{1/2}(\IId)}||_{r+1}(||\DD\tilde\Psi_I||_1+||\DD\Psi_I||_1+||R_\ell/\delta_{q+1}||_1)^r\\
        &+||\DD\tilde\Psi_I-\DD\Psi_I]||_r||f|_{B_{1/2}(\IId)\times B_{1/2}(\IId)}||_1\\
        &\lesssim \lambda_q^{r}\lambda_q^{[r-(\underline{r}-1)]^+(b-1)\gamma_\ell}\mathcal{T}_g\delta_{q+1}^{1/2}
    \end{split}
\end{equation*}
and we conclude:
\begin{equation*}
    \begin{split}
        ||\tilde a_I-a_I||_r\leq ||T_1||_r+||T_2||_r\lesssim \lambda_q^{r}\lambda_q^{[r-(\underline{r}-1)]^+(b-1)\gamma_\ell}\delta_{q+1}^{1/2}\left[(\ell_t/\tau^c)+\mathcal{T}_g\right].
    \end{split}
\end{equation*}

\noindent For the Alfv\'en transport estimate, we consider each term separately. We first rewrite:
\begin{equation*}
    \begin{split}
        \tilde z_q^\pm&=\tilde z^\pm_{\ell,I}+(\tilde z^\pm_{\ell}-\tilde z^\pm_{\ell,I})+(\tilde z_q^\pm-\tilde z_\ell^\pm),\\
        \tilde z_q^\pm&=z^\pm_{\ell,I}+(z^\pm_{\ell}-z^\pm_{\ell,I})+(z^\pm_q-z^\pm_\ell)+(\tilde z_q^\pm-z^\pm_q).
    \end{split}
\end{equation*}
Now, by means of Lemmas \ref{standardmollgalbrun}, \ref{stabilitygalbrun}, \ref{slowcoeffgalb}, \ref{standardmollnash} and \ref{stabilitynash}, and the estimate in \eqref{transportslownash} we proved above, we conclude:
\begin{equation*}
    \begin{split}
        ||\mathcal{\tilde A}^\pm(\tilde a_I-a_I)||_r&\leq||\mathcal{\tilde A}^\pm\tilde a_I||_r+||\mathcal{\tilde A}^\pm a_I||_r\\
        &\leq ||\mathcal{\tilde A}^\pm_{\ell,I}\tilde a_I||_r+||(\tilde z^\pm_{\ell}-\tilde z^\pm_{\ell,I})\cn\tilde a_I||_r+||(\tilde z_q^\pm-\tilde z_\ell^\pm)\cn\tilde a_I||_r\\
        &+||\mathcal{A}^\pm_{\ell,I}a_I||_r+||(z^\pm_{\ell}-z^\pm_{\ell,I})\cn a_I||_r+||(z^\pm_q-z^\pm_\ell)\cn a_I||_r+||(\tilde z_q^\pm-z^\pm_q)\cn a_I||_r\\
        &\lesssim \lambda_q^{r}\lambda_q^{[r-\underline{r}]^+(b-1)\gamma_\ell}1/\tau^c\delta_{q+1}^{1/2}+\lambda_q^{r+1}\lambda_q^{[r-(\underline{r}-1)]^+(b-1)\gamma_\ell}(\ell\lambda_q)^{m_0}\delta_q^{1/2}\delta_{q+1}^{1/2}\\
        &+\lambda_q^{r}\lambda_q^{[r-(\underline{r}-1)]^+(b-1)\gamma_\ell}\underbrace{\lambda_q^2\ell^{-\alpha}\tau^a\delta_{q+1}}_{1/\tau^c\mathcal{T}_g}\delta_{q+1}^{1/2}\\
        &\lesssim \lambda_q^{r}\lambda_q^{[r-\underline{r}]^+(b-1)\gamma_\ell}1/\tau^c\delta_{q+1}^{1/2}
    \end{split}
\end{equation*}
for $0\leq r\leq N-m_0$. Here we used $\lambda_q^{(b-1)\gamma_\ell}(\ell\lambda_q)^{m_0}, \ \lambda_q^{(b-1)\gamma_\ell}\mathcal{T}_g<1$, see \eqref{constraintadmissibility}, \eqref{constraintm0}. The implicit constant depends on $r$ and all the other parameters, but not on $a$.

\noindent \textbf{Pure time derivatives.} We write the pure time derivatives as:
\begin{equation}\label{puretimeslowcoeff}
    \begin{split}
        \partial_t \tilde a_I&=\mathcal{\tilde A}^{\pm}_{\ell,I}\tilde a_I-z^\pm_{\ell,I}\cn\tilde a_I,\\
        \partial_t^2\tilde a_I&=\mathcal{\tilde A}^{\mp}_{\ell,I}\mathcal{\tilde A}^{\pm}_{\ell,I}\tilde a_I-(\partial_t\tilde z^\pm_{\ell,I})\cn\tilde a_I-2\tilde z^\pm_{\ell,I}\cn\partial_t\tilde a_I-(\tilde z^\mp_{\ell,I}\cn)^2\tilde a_I\\
    \end{split}
\end{equation}
and from the bounds in Lemma \ref{localcorrn} and \eqref{pureslow}, \eqref{transportslownash} above we conclude:
$$||\partial_t^j \tilde a_I||_r\lesssim \lambda_q^{r+j}\lambda_q^{[r+j-\underline{r}]^+(b-1)\gamma_\ell}\delta_{q+1}^{1/2} \ \text{ for } j=0,1,2 \ \text{ and } \ r\geq 0.$$

\noindent \textbf{Support Property.} Note that by construction of the space partition of unity in Section \ref{geomconstr} with the obvious notation, we have: 
\begin{equation}\label{minimaldistance}
    p(I_1)=p(I_2)\Longrightarrow \text{dist}(U_{j'_1}, \ U_{j'_2})> 1/3\tau^c
\end{equation}
and
$$\supp_x \chi_{I_1}\subset U_{j'_1}, \ \supp_x \chi_{I_2}\subset U_{j'_2}.$$
The issue here is that we mollify the slow coefficients along the Alfv\'en direction through the operator $(\dots)^{I_t}_{\ell_t}$, namely
$$\tilde a_I=\delta_{q+1}^{1/2} \left(\chi_I\gamma_\zeta(\det[\DD \tilde\Psi_J]^{-2}\DD \tilde \Psi_J \left(\IId-R_\ell/\delta_{q+1}\right)\DD \tilde \Psi_J^T)\right)^{I_t}_{\ell_t}$$
and this translates the support of the function. The last property in Lemma \ref{timemollification} ensures that:
\begin{equation*}
    \begin{split}
        \supp_{t,x}\tilde a_{I_i}&\subset B_{2(C_0+1)\ell_t}\left(\supp_{x,t}\left(\chi_I\gamma_\zeta(\det[\DD \tilde\Psi_J]^{-2}\DD \tilde \Psi_J \left(\IId-R_\ell/\delta_{q+1}\right)\DD \tilde \Psi_J^T\right)\right)\\
        &\subset B_{2(C_0+1)\ell_t}\left(U_{j'_i}\times B_{2/3\tau^c}(t_{j_i})\cap \supp_{x,t}R_\ell\right)\\
    \end{split}
\end{equation*}
and given \eqref{minimaldistance}, upon choosing $a$ sufficiently large so that: $$4(C_0+1)\ell_t=4(C_0+1)\tau^a< 1/3\tau^c \ \text{ and } \ 2/3\tau^c+2(C_0+1)\tau^a<\tau^c$$ 
we can ensure:
$$\supp_{x}\ \tilde a_{I_1}\cap \supp_{x}\ \tilde a_{I_2}=\emptyset \ \text{ and } \ \supp_{x,t}\ \tilde a_{I}\Subset B_{\tau^c}(x_{j'})\times B_{\tau^c}(t_j)  \ \text{ and } \ \supp_{x,t}\ \tilde a_{I}\subset B_{2(C_0+1)\ell_t}(\supp_{x,t}R_\ell)$$
for $I_i$ as above and $I$ any index. The support properties are still guaranteed, and the claim in the Lemma follows with $C=2(C_0+1)$.
\end{proof}


\begin{lemma}[Chart Update Error]\label{Rcrt} Let $\underline{r}=M-2m_0-k_0^g-9$. Define the chart update error to be
    $$R^{p,crt}=\sum_Ig_I^2\left[\tilde A_I-A_I\right]$$
    then, the following bounds hold:
    \begin{equation*}
    \begin{split}
       & ||R^{p,crt}||_r\lesssim \lambda_q^r\lambda_q^{[r-\underline{r}]^+(b-1)\gamma_\ell}(\tau^a/\tau^c)\delta_{q+1} \ \text{ for }\  r \geq 0,\\
        &||\mathcal{\tilde A}^{\pm} R^{p,crt}||_r\lesssim \lambda_q^r\lambda_q^{[r-\underline{r}]^+(b-1)\gamma_\ell}(1/\tau^c)\delta_{q+1} \  \text{ for }\ 0\leq r \leq N-m_0.\\
    \end{split}        
    \end{equation*}
    The implicit constants depend on $r$ and all the other parameters, but not on $a$.
\end{lemma}
\begin{remark} The bounds given here differ from the ones for the `flow error' in \cite{GR} because of our choice $\gamma_t=\gamma_a$.
\end{remark}

\begin{proof}[Proof of Lemma \ref{Rcrt}] In what follows the $C^0_tC^r_x$ norms are thought on the space-time domain $Q_I$ of the chart. Recall that indexing $I$ all the objects we have:
\begin{equation*}
    \begin{split}
        A_I&=\delta_{q+1}a_I^2\Psi_I^{2*}\zeta\otimes \Psi_I^{2*}\zeta=\delta_{q+1} a_I^2\det[\DD \tilde \Psi_I]^2\DD \Psi_I^{-1}\zeta\otimes \DD \Psi_I^{-1}\zeta,\\
        \tilde A_I&=\delta_{q+1}\tilde a_I^2\tilde\Psi_I^{2*}\zeta\otimes \tilde\Psi_I^{2*}\zeta=\delta_{q+1}\tilde a_I^2\det[\DD \tilde \Psi_I]^2\DD \tilde \Psi_I^{-1}\zeta\otimes \DD \tilde \Psi_I^{-1}\zeta.
    \end{split}
\end{equation*}
The proof of this Lemma is an application of the work done in Section \ref{chart}, in particular of the chart stability Lemma \ref{chartstability} from which we deduce:
$$||(\DD\tilde\Psi_I)^{-1}-(\DD\Psi_I)^{-1}||_r, \ ||\det[\DD \tilde \Psi_I]-\det[\DD  \Psi_I]||_r\lesssim \lambda_q^{r}\lambda_q^{[r-\underline{r}]^+(b-1)\gamma_\ell}\mathcal{T}_g.$$
Given that for each space-time point at most one $g_I^2(\tilde A_I-A_I)$ is non-zero, see Lemma \ref{disjoint}, and the bounds in Lemma \ref{slowcoeffestimates}, we conclude:
\begin{equation*}
    \begin{split}
        ||R^{p,crt}||_r&\lesssim \lambda_q^{r}\lambda_q^{[r-\underline{r}]^+(b-1)\gamma_\ell}\delta_{q+1}\left[(\ell_t/\tau^c)+\mathcal{T}_g\right]\\
        & \lesssim \lambda_q^{r}\lambda_q^{[r-\underline{r}]^+(b-1)\gamma_\ell}(\tau^a/\tau^c)\delta_{q+1},
    \end{split}
\end{equation*}
for $r\geq 0$, where we used $\ell_t=\tau^a$ and $(\tau^a/\tau^c)\geq \mathcal{T}_g$.

\noindent Using the fact that $\mathcal{\tilde A}^\pm_{\ell,I}\det[\DD\tilde \Psi_I]=\mathcal{A}^\pm_{\ell,I}\det[\DD\Psi_I]=0$ and that $(\partial_t+\mathcal{L}_{\tilde z^\pm_{\ell,I}})\tilde \Psi^{*}_I\zeta=(\partial_t+\mathcal{L}_{ z^\pm_{\ell,I}})\Psi^{*}_I\zeta=0$ from Remark \ref{interestin} and Lemmas \ref{chartpropg}, \ref{chartprop}, and arguing as in the proof of Lemma  \ref{slowcoeffestimates}, we first deduce:
\begin{equation*}
    \begin{split}
        ||\mathcal{\tilde A}^\pm [g^2(A_I-\tilde A_I)]||_r & \leq||(\mathcal{\tilde A}^\pm g^2)(A_I-\tilde A_I)]||_r+||\mathcal{\tilde A}^\pm A_I||_r+ ||\mathcal{\tilde A}^\pm \tilde A_I||_r\\
        &\lesssim \lambda_q^{r}\lambda_q^{[r-\underline{r}]^+(b-1)\gamma_\ell}(1/\tau^a)\delta_{q+1}\left[(\ell_t/\tau^c)+\mathcal{T}_g\right]\\
        &+\lambda_q^{r}\lambda_q^{[r-\underline{r}]^+(b-1)\gamma_\ell}1/\tau^c\delta_{q+1}\\
        &\lesssim \lambda_q^{r}\lambda_q^{[r-\underline{r}]^+(b-1)\gamma_\ell}\delta_{q+1}[\underbrace{ \ell_t/(\tau^a\tau^c)}_{\leq1/\tau^c}+\underbrace{1/\tau^a\mathcal{T}_g}_{\leq 1/\tau^c\delta_{q+1}/\delta_q}]\\
        &+\lambda_q^{r}\lambda_q^{[r-\underline{r}]^+(b-1)\gamma_\ell}1/\tau^c\delta_{q+1}\\
        &\lesssim \lambda_q^{r}\lambda_q^{[r-\underline{r}]^+(b-1)\gamma_\ell}1/\tau^c\delta_{q+1},
    \end{split}
\end{equation*}
for $0\leq r\leq N-m_0$ and from the fact that at most one $g_I^2(\tilde A_I-A_I)$ is non-zero et each space-time point we conclude that:
\begin{equation*}
        ||\mathcal{\tilde A}^\pm R^{p,crt}||_r\leq \sup_I||\mathcal{\tilde A}^\pm [g^2(A_I-\tilde A_I)]||_r\lesssim \lambda_q^{r}\lambda_q^{[r-\underline{r}]^+(b-1)\gamma_\ell}1/\tau^c\delta_{q+1}.
\end{equation*}
The implicit constants in the above bounds depend on $r$ and all the parameters, but not on $a$.
\end{proof}


\subsection{Nash LDF}\label{nashldf}
With the time-oscillating profiles from Lemma \ref{timeprofiles}, the updated coefficients from \ref{slowcoeffestimates}, charts from \ref{chartprop} and cut-offs from \eqref{cutoffbounds}, we are ready to define the Nash LDF $\xi^p$. Using a similar notation to Lemma \ref{chartprop} given any $\varphi:\mathbb{R}\to\mathbb{R}$, any $\lambda>0$ and unit vector $k$ we define:
$$\varphi_{\lambda,k}(x)=\varphi(\lambda x\cdot k)$$ and think of $\varphi_{\lambda,k}$ as a function $\mathbb{R}^3\to \mathbb{R}$. We make the explicit choice $\varphi=\sqrt{2}\sin$, but this plays no role. It follows that: 
$$||\varphi||_r\leq \sqrt{2}\lesssim 1.$$
Given an index $I=(\zeta,J)\in \mathcal{I}$ and $(k, \nu, \zeta)$ the orthonormal basis associated to $\zeta$ from the Geometric Lemma \ref{geomlemma} we define label maps: 
$$I\to \nu_I,\ I\to k_I,\ I\to \zeta_I$$
associating the index $I$ to its orthonormal basis. Abusing of notation we will use  the label $I$ even if the labeled object depends only on some components of the index e.g. we write $\ \tilde \Psi_I,\ \tilde B_{\ell,I}, \ \alpha_I$ instead of $\tilde \Psi_{J}, \ \tilde B_{\ell,j}, \ \alpha_{p(I)}$. We now define $\Theta^p, \ \Theta^p_I$ to be the 1-forms:
\begin{equation}\label{ansatz}
        \Theta^p_I=\frac{\tau^a}{\lambda_{q+1}}\tilde a_I  \alpha_{I}\tilde\Psi_I^{1*}(\varphi_{\lambda_{q+1},k_I}\nu_I), \ \ \ \ \Theta^p=\sum_{I\in\mathcal{I}}\Theta_I^p\\
\end{equation}
and $\xi^p$ to be the divergence-free vector field having $\Theta^p$ as potential, namely,
\begin{equation*}
        \xi_I^p=\curl [\Theta_I^p], \ \ \ \ \xi^p=\curl [\Theta^p]=\sum_{I\in\mathcal{I}} \xi_I^p.
\end{equation*}
Computing the $\curl$ explicitly, we obtain the following representation:
\begin{equation}\label{xinash}
    \begin{split}
        \xi^p&=\curl[\Theta^p]\\
        &=\sum_{I\in\mathcal{I}}\frac{\tau^a}{\lambda_{q+1}}\tilde a_I\alpha_I\tilde \Psi^{2*}_I\curl\left [\varphi_{\lambda_{q+1},k_I}\nu_I\right]+\frac{\tau^a}{\lambda_{q+1}}\alpha_I\nabla\tilde a_I\times\tilde\Psi_I^{1*}(\varphi_{\lambda_{q+1},k_I}\nu_I)\\
        &=\sum_{I\in\mathcal{I}}\tau^a\tilde a_I\alpha_I\tilde \Psi^{2*}_I\left[\varphi_{\lambda_{q+1},k_I}'\zeta_I\right]+\frac{\tau^a}{\lambda_{q+1}}\alpha_I\nabla\tilde a_I\times\tilde\Psi_I^{1*}(\varphi_{\lambda_{q+1},k_I}\nu_I)\\
    \end{split}
\end{equation}
The definition is motivated by the Lie-transport properties of Lemma \ref{chartprop} and the geometric decomposition \ref{slowcoeff}. 
\begin{lemma}[Disjointedness]\label{disjoint} Let $I,I'\in \mathcal{I}$ with $I\neq I'$, then
\begin{equation*}
     \supp_{x,t}\Theta_I\cap \supp_{x,t} \Theta_{I'}=\emptyset,
\end{equation*}
since the property is retained for all space and time derivatives, in particular, we have $$\supp_{x,t}\xi_I\cap \supp_{x,t}\Theta_{I'}=0.$$ 
\end{lemma}

This is a simple consequence of our Ansatz \eqref{ansatz} and Lemmas \ref{timeprofiles}, \ref{slowcoeffestimates}, and ultimately due to the presence of the fast-oscillating time profiles $\{\alpha_I\}_I$ and a minimal partition in space $\{\theta_I\}_I$. Let $p$ be the type map in \eqref{typemap}, we distinguish two cases:
\begin{itemize}
    \item $p(I)=p(I')$ in this case the indices are associated to the same family in the covering $\mathcal{C}$ of $\mathbb{T}^3$ and are thus disjoint by construction. See Lemma \ref{slowcoeffestimates}.\\
    \item $p(I)\neq p(I')$ the time profiles $\alpha_I=\alpha_{p(I)}$ and $\alpha_{I'}=\alpha_{p(I')}$ are disjoint in time. See lemma \ref{timeprofiles}. No interactions also in this case.
\end{itemize} 

\subsubsection{Splitting of the Perturbation} \label{splittingnash}
In the language of Lemma \ref{lpl}, we now construct a perturbation $(w^p,b^p)$ of $(\tilde v_q,\tilde B_q)$ along $\xi^{p}$; the disjointness properties from Lemma \ref{disjoint} will play a crucial role. We follow the same strategy as for the Galbrun perturbations, see Section \ref{gnldf} for more details.

\noindent\textit{Velocity Field.} According to Lemma \ref{lpl} we now expand the perturbation $w^{g}$ along $\xi^{p}$ of $\tilde v_q$. Since $\xi^p_I$ is designed to behave correctly with $(\tilde v_{\ell,I}, \ \tilde B_{\ell,I})$, which are not globally defined in time, we need to be extra careful. Let us first separate the non-mollified part:
\begin{equation}\label{split1p}
    \begin{split}
        w^p&=\partial_tX^p\circ (X^p)^{-1}+(X^p_*-\IId_*)\tilde v_q\\
        &=\underbrace{\partial_tX^p\circ (X^p)^{-1}+(X^p_*-\IId_*)\tilde v_\ell}_{\bar w^p}+\underbrace{(X^p_*-\IId_*)(\tilde v_q-\tilde v_\ell)}_{\mathring w^p}
    \end{split}
\end{equation}
and according to Lemma \ref{lpl} we have:
\begin{equation}\label{decompnashw}
    \begin{split}
        \bar w^p 
            &=\curl\left[\sum_I\sum_{k=0}^{k_0^p}\frac{(-1)^k}{(k+1)!}\mathcal{L}_{\xi^p_I}^k(\partial_t+\mathcal{L}_{\tilde v_{\ell,I}})\Theta^p_I+\sum_I\sum_{k=0}^{k_0^p}\frac{(-1)^k}{(k+1)!}\mathcal{L}_{\xi^p_I}^k\mathcal{L}_{\tilde v_\ell-\tilde v_{\ell,I}}\Theta_I^p+\theta_w^p\right],\\
        \mathring w^p
                    &=\curl\left[\sum_I\sum_{k=0}^{k_0^p}\frac{(-1)^k}{(k+1)!}\mathcal{L}_{\xi^p_I}^k\mathcal{L}_{\tilde v_q-\tilde v_\ell}\Theta_I^p+\mathring \theta_w^p\right]\\
    \end{split}
\end{equation}
where 
\begin{equation*}
    \begin{split}
        \theta_w^p&= \sum_I\frac{(-1)^{k_0^p+1}}{(k_0^p+1)!}\int_0^{1}(X_{s}^p)_*\left[\mathcal{L}_{\xi^p_I}^{k_0^p+1}(\partial_t+\mathcal{L}_{\tilde v_{\ell,I}})\Theta^p_I\right](1-s)^{k_0^p+1}\dd s\\
        &+\sum_I\frac{(-1)^{k_0^p+1}}{(k_0^p+1)!}\int_0^{1}(X_{s}^p)_*\left[\mathcal{L}_{\xi^p_I}^{k_0^p+1}\mathcal{L}_{\tilde v_\ell-\tilde v_{\ell,I}}\Theta_I^p\right](1-s)^{k_0^p+1}\dd s,\\
        \mathring \theta_w^p&=\sum_I\frac{(-1)^{k_0^p+1}}{(k_0^p+1)!}\int_0^{1}(X_{s}^p)_*\left[\mathcal{L}_{\xi^p_I}^{k_0^p+1}\mathcal{L}_{\tilde v_q-\tilde v_\ell}\Theta_I^p\right](1-s)^{k_0^p+1}\dd s.
    \end{split}
\end{equation*}

\noindent \textit{Magnetic field.} We proceed similarly, we first split:
\begin{equation*}
    \begin{split}
        &b^g=\underbrace{(X^p_*-\IId_*)\tilde B_\ell}_{\bar{b}^p}+\underbrace{(X^p_*-\IId_*)(\tilde B_q-\tilde B_\ell)}_{\mathring b^p}
    \end{split}
\end{equation*}
and Lemma \ref{lpl} gives:
\begin{equation}\label{decompnashb}
    \begin{split}
        \bar b^p 
            &=\curl\left[\sum_I\sum_{k=0}^{k_0^p}\frac{(-1)^k}{(k+1)!}\mathcal{L}_{\xi^p_I}^k\mathcal{L}_{\tilde B_{\ell,I}}\Theta^p_I+\sum_I\sum_{k=0}^{k_0^p}\frac{(-1)^k}{(k+1)!}\mathcal{L}_{\xi^p_I}^k\mathcal{L}_{\tilde B_\ell-\tilde B_{\ell,I}}\Theta_I^p+\theta^p_b\right],\\
        \mathring b^p
                    &=\curl\left[\sum_I\sum_{k=0}^{k_0^p}\frac{(-1)^k}{(k+1)!}\mathcal{L}_{\xi^p_I}^k\mathcal{L}_{\tilde B_q-\tilde B_\ell}\Theta_I^p+\mathring\theta^p_b\right]\\
    \end{split}
\end{equation}
where 
\begin{equation*}
    \begin{split}
        \theta^p_b&=\sum_I\frac{(-1)^{k_0^p+1}}{(k_0^p+1)!}\int_0^{1}(X_{s}^p)_*\left[\mathcal{L}_{\xi^p_I}^{k_0^p+1}\mathcal{L}_{\tilde B_{\ell,I}}\Theta^p_I\right](1-s)^{k_0^p+1}\dd s\\
        &+\sum_I\frac{(-1)^{k_0^p+1}}{(k_0^p+1)!}\int_0^{1}(X_{s}^p)_*\left[\mathcal{L}_{\xi^p_I}^{k_0^p+1}\mathcal{L}_{\tilde B_\ell-\tilde B_{\ell,I}}\Theta_I^p\right](1-s)^{k_0^p+1}\dd s,\\
        \mathring\theta^p_b&=\sum_I\frac{(-1)^{k_0^p+1}}{(k_0^p+1)!}\int_0^{1}(X_{s}^p)_*\left[\mathcal{L}_{\xi^p_I}^{k_0^p+1}\mathcal{L}_{\tilde B_q-\tilde B_\ell}\Theta_I^p\right](1-s)^{k_0^p+1}\dd s.\\
    \end{split}
\end{equation*}

\noindent Note that in the higher-order terms of the Lie-Taylor expansion, there are no cross terms involving different indices; this is a consequence of Lemma \ref{disjoint}. 

\noindent We now study the leading terms:
\begin{equation}\label{leadingdefinition}
    w^{p,p}_I=\curl\left[(\partial_t+\mathcal{L}_{\tilde v_{\ell,I}})\Theta_I^{p}\right], \ b_I^{p,p}=\curl\left[\mathcal{L}_{\tilde B_{\ell,I}}\Theta_I^{p}\right], \ w^{p,p}=\sum_I w^{p,p}_I, \ b^{p,p}=\sum_I b^{p,p}_I.
\end{equation}
Using the definition of $\Theta_I^p$ in \eqref{ansatz} and the commutation \eqref{DG1} we compute:
\begin{equation}\label{leadingtermnashexpansion}
        \begin{split}
            w^{p,p}_I&=\frac{\tau^a}{\lambda_{q+1}}\curl\left[(\partial_t+\mathcal{L}_{\tilde v_{\ell,I}})\left[ \tilde a_I \alpha_{I}\Psi_I^{1*}(\varphi_{\lambda_{q+1},k_I}\nu_I)\right]\right]\\
            &=\frac{1}{\lambda_{q+1}}\curl\Biggr[ \tilde a_I \underbrace{\alpha_I'}_{=g_{I}}\tilde \Psi_I^{1*}(\varphi_{\lambda_{q+1},k_I}\nu_I)+\tau^a \alpha_I(\partial_t+\tilde v_{\ell,I}\cn)( \tilde a_I)\tilde \Psi_I^{1*}(\varphi_{\lambda_{q+1},k_I}\nu_I)\Biggr]\\
            &=\underbrace{ \tilde a_I g_{I}\tilde \Psi_I^{2*}(\varphi_{\lambda_{q+1},k_I}'\zeta_I)}_{w_I^{p,o}}+\frac{1}{\lambda_{q+1}}g_I\nabla\tilde a_I\times\Psi_I^{1*}(\varphi_{\lambda_{q+1},k_I}\nu_I)\\
            &+\tau^a \alpha_I\frac{1}{2}(\mathcal{A}^+_{\ell,I}+\mathcal{A}^-_{\ell,I})( \tilde a_I)\tilde \Psi_I^{2*}(\varphi_{\lambda_{q+1},k_I}'\zeta_I)+\frac{\tau^a}{\lambda_{q+1}} \alpha_I\frac{1}{2}\nabla(\mathcal{A}^+_{\ell,I}+\mathcal{A}^-_{\ell,I})( \tilde a_I)\times \tilde \Psi_I^{1*}(\varphi_{\lambda_{q+1},k_I}\nu_I)\\
            &=w_I^{p,o}+w_I^{p,c}\\
            b_I^{p,p}&=\frac{\tau^a}{\lambda_{q+1}}\curl\left[\mathcal{L}_{\tilde B_{\ell,I}}\left[ \tilde a_I \alpha_{I}\Psi_I^{1*}(\varphi_{\lambda_{q+1},k_I}\nu_I)\right]\right]=\frac{\tau^a}{\lambda_{q+1}}\curl\left[\tilde B_{\ell,I}\cn( \tilde a_I) \alpha_{I}\Psi_I^{1*}(\varphi_{\lambda_{q+1},k_I}\nu_I)\right]\\
            &=\tau^a\alpha_{I}\frac{1}{2}(\mathcal{\tilde A}^+_{\ell,I}-\mathcal{\tilde A}^-_{\ell,I})(\tilde a_I) \Psi_I^{2*}(\varphi_{\lambda_{q+1},k_I}'\zeta_I)+\frac{\tau^a}{\lambda_{q+1}}\alpha_{I}\frac{1}{2}\nabla(\mathcal{\tilde A}^+_{\ell,I}-\mathcal{\tilde A}^-_{\ell,I})\tilde a_I \Psi_I^{1*}(\varphi_{\lambda_{q+1},k_I}\nu_I)
        \end{split}
\end{equation}
where we used the transport properties of the charts in Lemma \ref{chartprop} and the fact that by definition $\alpha_I'=g_I$, see Lemma \ref{timeprofiles}. This makes the computation \eqref{wfromxi} in the overview Subsection \ref{overview} rigorous. In light of the geometric decomposition in \eqref{reynoldsdecomposition} it is clear how we will use the term 
$$w^{p,o}_I=  \tilde a_I g_{I}\tilde \Psi_I^{2*}(\varphi_{\lambda_{q+1},k_I}'\zeta_I)$$ 
to correct the `well prepared' Reynolds stress coming from the Galbrun stage. As it is often done in convex integration, we refer to $w_I^{p,o}$ as the oscillation part and $w_I^{p,c}$ as the corrector part of the leading term of the velocity field perturbation.
\begin{lemma}[Leading Terms]\label{leadingtermsnash} For any $I\in \mathcal{I}$,  $r\geq 0$ integer and $j=0,1,2$ we have:
    \begin{equation*}
        \begin{split}
            &||\partial_t^jw^{p,p}_I||_r\lesssim \lambda_{q+1}^{r+j}\delta_{q+1}^{1/2},\\
            &||\mathcal{\tilde A}^\pm_{\ell,I} w^{p,p}_I||_r\lesssim\lambda_{q+1}^r(1/\tau^a)\delta_{q+1}^{1/2}.\\
        \end{split}
    \end{equation*}
    The implicit constants depend on $r$ and all the other parameters, but not on $I$ and $a$. Moreover, given any $\bar N$ and $\beta, \ b$, there exists $a$ sufficiently large such that the implicit constant is fixed for  $0\leq r \leq \bar N$ and independent of $C_0$ and the other parameters.

    \noindent We also have the following bounds for the corrector and magnetic terms:
    \begin{equation*}
        \begin{split}
            &||\partial_t^jw^{p,c}_I||_r, ||\partial_t^j b^{p,p}_I||_r\lesssim \lambda_{q+1}^{r+j}(\tau^a/\tau^c)\delta_{q+1}^{1/2},\\
            &||\mathcal{\tilde A}^\pm_{\ell,I} w^{p,c}_I||_r, ||\mathcal{\tilde A}^\pm_{\ell,I} b^{p,p}_I||_r\lesssim \lambda_{q+1}^r(\tau^a/\tau^c)\delta_{q+1}^{1/2}\left[\frac{1}{\tau^a}+\frac{1}{\ell_t}\right].\\
        \end{split}
    \end{equation*}
    The implicit constants depend on $r$ and all the other parameters, but not on $a$.
\end{lemma}
\begin{proof}[Proof of Lemma \ref{leadingtermsnash}] We will prove the bounds for the velocity field only; the ones for the magnetic field follow as those for $w^{p,c}$. We set $\underline{r}=M-2m_0-k_0^g-8$.

\noindent \textbf{Oscillation term $w^{p,o}_I$.} From the definition in  \eqref{leadingtermnashexpansion}, the bounds in Lemmas \ref{chartprop}, \ref{slowcoeffestimates} and the composition estimates in Proposition \ref{compestimates}, we deduce that:
\begin{equation}\label{principaloscillation}
    \begin{split}
        ||w^{p,o}_I||_r&\leq \sup_t|g_I| \ ||\tilde a_I\det[\DD\tilde \Psi_I]\varphi'(\lambda_{q+1}\tilde \Psi_I\cdot k_I)\tilde \Psi_I^{-1}\zeta_I||_r\\
        &\lesssim\sup_t|g_I|\left[||\tilde a_I||_r||\det[\DD\tilde \Psi_I]||_0||\DD\tilde \Psi_I^{-1}\zeta||_0||\varphi'||_0+||\tilde a_I||_0||\det[\DD\tilde \Psi]||_r||\DD\tilde \Psi_I^{-1}\zeta||_0||\varphi'||_0\right]\\
        &+\sup_t|g_I|||\tilde a_I||_0||\det[\DD\tilde \Psi]||_0||\DD\tilde \Psi_I^{-1}\zeta||_r||\varphi'||_0\\
        &+\sup_t|g_I|||\tilde a_I||_0||\det[\DD\tilde\Psi]||_0||\DD\tilde \Psi_I^{-1}\zeta||_0||\varphi'(\lambda_{q+1}\tilde \Psi_I\cdot k_I)||_r\\
        &\overset{r=0}{\lesssim}\delta_{q+1}^{1/2}\\
        &\overset{r\geq 1}{\leq} C_r \lambda_q^r\lambda_q^{[r-\underline{r}]^+(b-1)\gamma_\ell}\delta_{q+1}^{1/2}+C'_r\delta_{q+1}^{1/2}\left[\lambda_{q+1}||\DD \tilde \Psi_I||_{r-1}+\lambda_{q+1}^r||\DD\tilde\Psi_I||_0^r\right]\\
        &\leq C_r \lambda_q^r\lambda_q^{[r-\underline{r}]^+(b-1)\gamma_\ell}\delta_{q+1}^{1/2}+C'_r\delta_{q+1}^{1/2}\left[C_r''\lambda_{q+1}\lambda_q^{r-1}\lambda_q^{[r-1-\underline{r}]^+(b-1)\gamma_\ell}+\lambda_{q+1}^r2^r\right]\\
        &\leq(C'_r2^r)\lambda_{q+1}^r\delta_{q+1}^{1/2}\left[1+\frac{C_r}{(C'_r2^r)}\left(\frac{\ell^{-1}}{\lambda_{q+1}}\right)^r+\frac{C''_r}{2^r}\left(\frac{\ell^{-1}}{\lambda_{q+1}}\right)^{r-1}\right]
    \end{split}
\end{equation}
Here, the implicit constant in the bound for $r=0$ is independent of $C_0$. Moreover, $C_r, \ C''_r$ depend on $r, \ C_0$ and the other parameters, while $C_r'$ depends on $r$ but not on $C_0$ and the other parameters. Upon choosing $a$ sufficiently large, we can ensure that: 
$$1+\frac{C_r}{(C'_r2^r)}\left(\frac{\ell^{-1}}{\lambda_{q+1}}\right)^r+\frac{C''_r}{2^r}\left(\frac{\ell^{-1}}{\lambda_{q+1}}\right)^{r-1}\leq 2 \ \text{ for } \ 0\leq r \leq \bar N$$
and conclude that: 
$$||w^{p,o}_I||_r\leq C_{\bar N}\lambda_{q+1}^r\delta_{q+1}^{1/2} \ \text{ for } \ 0\leq r \leq \bar N$$
for some fixed $C_{\bar N}$ which depends on $\bar N$ but is independent of $C_0$ and the other parameters. Arguing similarly for the $\partial_t, \ \partial_t^2$ case, the claimed bound follows.

\noindent To be precise, one should distinguish the case $r=1$ as well to ensure that $C''_1$ does not depend on $C_0$, this is the case because according to Lemma \ref{chartprop} we have: 
$$||\DD \tilde \Psi_I^{-1}-\IId||_0\lesssim \lambda_{q+1}^{-\alpha},$$ 
this is also what is behind the claim for $r=0$.

\noindent We now recall Remark \ref{interestin} which gives $\mathcal{\tilde A}^\pm_{\ell,I}\det[\DD\tilde \Psi_I]=0$. The estimates and transport properties in Lemmas \ref{chartprop}, \ref{slowcoeffestimates}, then give:
\begin{equation*}
    \begin{split}
        ||\mathcal{\tilde A}^\pm_{\ell,I} w^{p,p}_I||_r&\lesssim 1/\tau^a\sup_t|g_I'| \ ||\tilde a_I\tilde \Psi_I^{2*}(\varphi_{\lambda_{q+1},k_I}'\zeta_I)||_r+\sup_t|g_I|\ ||\mathcal{\tilde A}^\pm_{\ell,I}(\tilde a_I)\tilde \Psi_I^{2*}(\varphi_{\lambda_{q+1},k_I}'\zeta_I)||_r\\
        &+\sup_t|g_I|\ ||\tilde a_I\tilde \Psi_I^{2*}(\varphi_{\lambda_{q+1},k_I}'\zeta_I)\cn \tilde z^\pm_{\ell,I}||_r.
    \end{split}
\end{equation*}
We can now argue as above and deduce that upon taking $a$ sufficiently large, we have:
$$||\mathcal{\tilde A}^\pm_{\ell,I} w^{p,p}_I||_r\leq C_{\bar N} \lambda_{q+1}^r1/\tau^a\delta_{q+1}^{1/2} \ \text{ for } \ 0\leq r \leq \bar N$$
where $C_{\bar N}$ is a possibly different constant compared to the one above and is independent of $C_0$ and the other parameters.

\noindent \textbf{Corrector term $w_I^{p,c}$.} According to \eqref{leadingtermnashexpansion} we have:
\begin{equation*}
    \begin{split}
        w_I^{p,c}&=\underbrace{\tau^a \alpha_I\frac{1}{2}(\mathcal{A}^+_{\ell,I}+\mathcal{A}^-_{\ell,I})( \tilde a_I)\tilde \Psi_I^{2*}(\varphi_{\lambda_{q+1},k_I}'\zeta_I)}_{T_1}\\
        &+\underbrace{\frac{1}{\lambda_{q+1}}g_I\nabla \tilde a_I\times\Psi_I^{1*}(\varphi_{\lambda_{q+1},k_I}\nu_I)}_{T_{2,a}}+\underbrace{\frac{\tau^a}{\lambda_{q+1}} \alpha_I\frac{1}{2}\nabla(\mathcal{A}^+_{\ell,I}+\mathcal{A}^-_{\ell,I})( \tilde a_I)\times \tilde \Psi_I^{1*}(\varphi_{\lambda_{q+1},k_I}\nu_I)}_{T_{2,b}}.
    \end{split}
\end{equation*}
The estimates on $T_1$ follow exactly as the ones for $w_I^{p,o}$ and we omit the details:
\begin{equation*}
    \begin{split}
        &||T_1||_r\lesssim\lambda_{q+1}^r(\tau^a/\tau^c)\delta_{q+1}^{1/2} \ \text{ for } \ r\geq 0,\\
        &||\mathcal{A}^\pm_{\ell,I}T_1||_r\lesssim \lambda_{q+1}^r\delta_{q+1}^{1/2}\left[\frac{1}{\tau^c}+\frac{\tau^a}{\tau^c\ell_t}\right] \ \text{ for } \ r\geq 0,
    \end{split}
\end{equation*}
We now prove bounds for $T_{2,a}$. $T_{2,b}$ has the same structure, and we omit the details. From the estimates in Lemmas \ref{chartprop} and \ref{slowcoeffestimates}, we deduce:
\begin{equation*}
    \begin{split}
        ||T_{2,a}||_r&\lesssim\frac{1}{\lambda_{q+1}}\sup_I|g_I|\left[||\tilde a_I||_{r+1}||\Psi_I^{1*}\nu_I||_0||\Psi_I^{*}\varphi_{\lambda_{q+1},k_I}||_0+||\tilde a_I||_{1}||\Psi_I^{1*}\nu_I||_r||\Psi_I^{*}\varphi_{\lambda_{q+1},k_I}||_0\right]\\
        &+\frac{1}{\lambda_{q+1}}\sup_I|g_I|\ ||\tilde a_I||_1||\Psi_I^{1*}\nu_I||_0||\Psi_I^{*}\varphi_{\lambda_{q+1},k_I}||_r\\
        &\lesssim \lambda_q^r\lambda_q^{[r-\underline{r}]^+(b-1)\gamma_\ell}\frac{\lambda_q}{\lambda_{q+1}}\delta_{q+1}^{1/2}+\lambda_{q+1}^r\frac{\lambda_q}{\lambda_{q+1}}\delta_{q+1}^{1/2}\\
        &\lesssim\lambda_{q+1}^r\frac{\lambda_q}{\lambda_{q+1}}\delta_{q+1}^{1/2}.\\
    \end{split}
\end{equation*}

\noindent Given the properties of the chart from Lemma \ref{chartprop}, we deduce:
\begin{equation*}
    \begin{split}
        ||\mathcal{\tilde A}^\pm_{\ell,I}T_{2,a}||_r&\lesssim \frac{1}{\lambda_{q+1}}\sup_I|g_I|\left[||[\nabla \mathcal{\tilde A}^\pm_{\ell,I}\tilde a_I]\times\Psi_I^{1*}(\varphi_{\lambda_{q+1},k_I}\nu_I)||_r+||(\DD z^\pm_{\ell,I})^\top\nabla \tilde a_I\times\Psi_I^{1*}(\varphi_{\lambda_{q+1},k_I}\nu_I)||_r\right]\\
        &+\frac{1}{\lambda_{q+1}}\sup_I|g_I| \ ||\nabla \tilde a_I\times(\DD z^{\pm}_{\ell,I})^\top[\Psi_I^{1*}(\varphi_{\lambda_{q+1},k_I}\nu_I)]||_r+\frac{1}{\lambda_{q+1}\tau^a}\sup_I|\partial_tg_I| \ ||\nabla \tilde a_I\times\Psi_I^{1*}(\varphi_{\lambda_{q+1},k_I}\nu_I)||_r\\
        &\lesssim \lambda_{q+1}^r\frac{\lambda_q}{\lambda_{q+1}}(1/\tau^a)\delta_{q+1}^{1/2}+\lambda_{q+1}^r\frac{\lambda_q}{\lambda_{q+1}}(1/\tau^c)\delta_{q+1}^{1/2}\\
        &\lesssim \lambda_{q+1}^r\frac{\lambda_q}{\lambda_{q+1}}(1/\tau^a)\delta_{q+1}^{1/2}.
    \end{split}
\end{equation*}
Similarly, one can show:
\begin{equation*}
    \begin{split}
        &||T_{2,b}||_r\lesssim \lambda_{q+1}^r\frac{\lambda_q}{\lambda_{q+1}}(\tau^a/\tau^c)\delta_{q+1}^{1/2},\\
        &||\mathcal{\tilde A}^\pm_{\ell,I}T_{2,b}||_r\lesssim\lambda_{q+1}^r\frac{\lambda_q}{\lambda_{q+1}}\delta_{q+1}^{1/2}\left[\frac{1}{\tau^c}+\frac{\tau^a}{\tau^c\ell_t}\right].
    \end{split}
\end{equation*}

\noindent Gathering the estimates above, we conclude: 
\begin{equation*}
    \begin{split}
        ||w^{p,c}_I||_r&\lesssim ||T_1||_r+||T_{2,a}||_r+||T_{2,b}||_r\\
        &\lesssim   \lambda_{q+1}^r(\tau^a/\tau^c)\delta_{q+1}^{1/2}+ \lambda_{q+1}^r\frac{\lambda_q}{\lambda_{q+1}}\delta_{q+1}^{1/2}+\lambda_{q+1}^r\frac{\lambda_q}{\lambda_{q+1}}(\tau^a/\tau^c)\delta_{q+1}^{1/2}\\
        &\lesssim\lambda_{q+1}^r(\tau^a/\tau^c)\delta_{q+1}^{1/2},\\
        ||\mathcal{\tilde A}^\pm_{\ell,I}w^{p,c}_I||_r&\lesssim ||\mathcal{\tilde A}^\pm_{\ell,I}T_1||_r+||\mathcal{\tilde A}^\pm_{\ell,I}T_{2,a}||_r+||\mathcal{\tilde A}^\pm_{\ell,I}T_{2,b}||_r\\
        &\lesssim \lambda_{q+1}^r\delta_{q+1}^{1/2}\left[\frac{1}{\tau^c}+\frac{\tau^a}{\tau^c\ell_t}\right]+\lambda_{q+1}^r\frac{\lambda_q}{\lambda_{q+1}}(1/\tau^a)\delta_{q+1}^{1/2}+\lambda_{q+1}^r\frac{\lambda_q}{\lambda_{q+1}}\delta_{q+1}^{1/2}\left[\frac{1}{\tau^c}+\frac{\tau^a}{\tau^c\ell_t}\right]\\
        &\lesssim \lambda_{q+1}^r(\tau^a/\tau^c)\delta_{q+1}^{1/2}\left[\frac{1}{\tau^a}+\frac{1}{\ell_t}\right]
    \end{split}
\end{equation*}
for $r\geq 0$, where used that $\tau^a=\ell_t$. The implicit constants depend on $r$ and all the other parameters, but not on $a$.
\end{proof}

\begin{remark}[Inductive Lemma]\label{inductiveapplication} We remark that this is exactly the type of vector fields we will use in combination with the Inductive Lemma \ref{inductive}. Indeed, in the computation of the new Reynolds stress, see \eqref{Rptrequation}, we will have a Lie-Taylor expansion in which the Lie-derivation happens repeatedly along:
\begin{equation*}
        \begin{split}
            \sigma_1&=\xi_I^p=\frac{1}{\lambda_{q+1}}\curl[\underbrace{\tau^a\alpha_I \tilde a_I }_{=:\varsigma_1}\tilde\Psi_I^{1*}(\varphi_{\lambda_{q+1},k_I}\nu_I)]\\
            \sigma_2&=(\partial_t+\mathcal{L}_{\tilde v_{\ell,I}})\xi_I^p=\frac{1}{\lambda_{q+1}}\curl[\underbrace{(\alpha'_I\tilde a_I+\tau^a\alpha_I (\partial_t+\tilde v_{\ell,I}\cn)(\tilde a_I))}_{=:\varsigma_2} \Psi_I^{1*}(\varphi_{\lambda_{q+1},k_I}\nu_I)]\\
        \sigma_3&=\mathcal{L}_{\tilde B_{\ell,I}}\xi_I^p=\frac{1}{\lambda_{q+1}}\curl[\underbrace{\tau^a\alpha_I \tilde B_{\ell,I}\cn(\tilde a_I)}_{=:\varsigma_3}\tilde\Psi_I^{1*}(\varphi_{\lambda_{q+1},k_I}\nu_I)]
        \end{split}
\end{equation*}
In the notation of Lemma \ref{inductive}, we then have:
\begin{equation*}
    \begin{split}
        &\bar \varsigma_1=\tau^a\delta_{q+1}^{1/2}, \ \bar \varsigma_{1,\mathcal{A}}=\delta_{q+1}^{1/2},\\
        &\bar \varsigma_2=\delta_{q+1}^{1/2}, \ \bar \varsigma_{2,\mathcal{A}}=[1/\tau^a+\tau^a/(\tau^c\ell_t)]\delta_{q+1}^{1/2},\\
        &\bar \varsigma_3=(\tau^a/\tau^c)\delta_{q+1}^{1/2}, \ \bar \varsigma_{2,\mathcal{A}}=[1/\tau^c+\tau^a/(\tau^c\ell_t)]\delta_{q+1}^{1/2},\\
        &L^\varsigma(r)=\lambda_q^{[r-\underline{r}]^+(b-1)\gamma_\ell}\\
    \end{split}
\end{equation*}
where $\underline{r}=M-2m_0-k_0^g-8$.
\end{remark}


\subsection{Estimates on the Perturbation}\label{leadings}

To prove estimates on the perturbation $(w^p,\ b^p)$, we first need bounds on the various terms appearing in their Lie-Taylor expansion coming from Lemma \ref{lpl}. As for the Galbrun stage, a key quantity is the smallness gain $\mathcal{T}_p$ corresponding to each additional Lie-derivative. This gain occurs only when the Lie derivative acts on objects with a specific, well-designed geometry. We define a corresponding quantity $\mathcal{M}_p$ telling us how much we lose when we take a Lie derivative in a direction on which we have no geometric control, e.g. $\tilde v_q-\tilde v_\ell$. Note that this is an issue we did not have in the Galbrun stage, where all the quantities were `slow' in space and thus did not produce gradients of order $\lambda_{q+1}$. Not surprisingly, $\mathcal{M}_p$ corresponds to the $C^0$ size of $\DD\xi^p$. More precisely, we set:
\begin{equation}\label{Tp}
    \begin{split}
        &\mathcal{T}_p=\lambda_q\tau^a\delta_{q+1}^{1/2}=\lambda_{q+1}^{-\alpha}(\tau^a/\tau^c)(\delta_{q+1}/\delta_q)^{1/2}=\left(\frac{\lambda_q}{\lambda_{q+1}}\right)^{\gamma_a+\beta+\gamma_{CZ}},\\
        &\mathcal{M}_p=\frac{\lambda_{q+1}}{\lambda_q}\mathcal{T}_p=\left(\frac{\lambda_{q+1}}{\lambda_{q}}\right)^{1-(\gamma_a+\beta+\gamma_{CZ})}.
    \end{split}
\end{equation}
Moreover, to deal with the degrading bounds on $\tilde v_q, \ \tilde B_q$, we define for $\underline{r}=M-2m_0-k_0^g-8$ loss functions $L_p,L_{p,\mathcal{A}}:\mathbb{N}_{\geq 0}\to \mathbb{R}_{\geq 1}$ with:
\begin{equation} \label{lossparameters}
    \begin{split}
        L_p(r)&=\lambda_q^{[r-\underline{r}]]^+(b-1)\gamma_\ell},
        \\
        L_{p,\mathcal{A}}(r)
        &=1_{r\leq \underline{r}-k_0^g-2}+1_{\underline{r}-k_0^g-1\leq r\leq \underline{r}-1} \left[1+\left(\lambda_q^{(b-1)\gamma_\ell}\mathcal{T}_g\right)^{\underline{r}-1-r}\bar L\right]+1_{ r\geq \underline{r}}\lambda_q^{[r-(\underline{r}-1)](b-1)\gamma_\ell}\bar L.\\
    \end{split}
\end{equation}
These are nothing but a shift of $m_0+k_0^g+2$ of the loss functions appearing in Lemma \ref{recap}. This is not strictly necessary, but it simplifies the estimates by ensuring that:
\begin{equation}\label{lossyvsmoll}
    L_p(r)\leq L_{p,\mathcal{A}}(r)
\end{equation}
where the left-hand side is the loss function associated with the Alfv\'en transport of the slow-coefficients $\mathcal{\tilde A}_{\ell,I}^\pm \tilde a_I$ from Lemma \ref{slowcoeffestimates}.

\noindent We now show that $L_{p,\mathcal{A}}$, is an admissible loss function in the sense of Definition \ref{admissible}, we need only verify that:
$$\left(\frac{\lambda_q}{\lambda_{q+1}}\right)^{r'}L_{p,\mathcal{A}}(r+r')\leq L_{p,\mathcal{A}}(r) \ \text{ for } \ r,r'\geq 0.$$
Let:
$$I_1=[0,\ \underline{r}-k_0^g-2]\cap \mathbb{N}, \ \ \ \ I_2=[\underline{r}-k_0^g-1,\  \underline{r}-1]\cap \mathbb{N}, \ \ \ \ I_3=\{r\geq \underline{r}\}\cap \mathbb{N}$$
and we recall that because of our choice of parameters in \ref{choiceofparameters}, see \eqref{constraintadmissibility}, \eqref{misc1}, we have:
$$\left(\lambda_q^{(b-1)\gamma_\ell}\mathcal{T}_g\right)^{k_0^g}\bar L\leq 1/2, \ \ \ \ 1-(\gamma_a+2\beta+\gamma_\ell+2\gamma_{CZ})\geq 0,\ \ \ \ \frac{\lambda_q}{\lambda_{q+1}}\leq \frac{1}{2}$$
We now distinguish the following cases:

\noindent \textbf{Case $r\in I_1$.} If $r+r'\in I_1$, the claim is trivial, if $r+r'\in I_2$  we have $r'\geq 1, \ \underline{r}-1-r\geq k_0^g+1\geq k_0^g$ and we can bound:
\begin{equation*}
    \begin{split}
        \left(\frac{\lambda_q}{\lambda_{q+1}}\right)^{r'}L_{p,\mathcal{A}}(r+r')&=\left(\frac{\lambda_q}{\lambda_{q+1}}\right)^{r'}\left[1+\left(\lambda_q^{(b-1)\gamma_\ell}\mathcal{T}_g\right)^{\underline{r}-1-(r+r')}\bar L\right]\\
        &\leq\left[\left(\frac{\lambda_q}{\lambda_{q+1}}\right)^{r'}+\left(\frac{\lambda_q}{\lambda_{q+1}}\right)^{r'\underbrace{[1-(\gamma_a+2\beta+\gamma_\ell+2\gamma_{CZ})]}_{\geq 0}}\left(\lambda_q^{(b-1)\gamma_\ell}\mathcal{T}_g\right)^{\underline{r}-1-r}\bar L\right]\\
        &\leq\left[\frac{\lambda_q}{\lambda_{q+1}}+\left(\lambda_q^{(b-1)\gamma_\ell}\mathcal{T}_g\right)^{k_0^g}\bar L\right]\\
        &\leq 1=L_{p,\mathcal{A}}(r).
    \end{split}
\end{equation*}
If $r+r'\in I_3$ since $r\in I_1$ we have $r'\geq k_0^g, \ r\leq \underline{r}-k_0^g-2\leq \underline{r}-1$ and thus:
\begin{equation*}
    \begin{split}
        \left(\frac{\lambda_q}{\lambda_{q+1}}\right)^{r'}L_{p,\mathcal{A}}(r+r')&=\left(\frac{\lambda_q}{\lambda_{q+1}}\right)^{r'}\lambda_q^{[r+r'-(\underline{r}-1)](b-1)\gamma_\ell}\bar L\\
        &=\left(\frac{\lambda_q}{\lambda_{q+1}}\right)^{r'(1-\gamma_\ell)}\lambda_q^{\underbrace{[r-(\underline{r}-1)](b-1)\gamma_\ell}_{\leq 0}}\bar L\\
        &\leq \left(\frac{\lambda_q}{\lambda_{q+1}}\right)^{k_0^g(\gamma_a+2\beta+2\gamma_{CZ})}\bar L\\
        &\leq \mathcal{T}_g^{k_0^g}\bar L\\
        &\leq 1=L_{p,\mathcal{A}}(r)
    \end{split}
\end{equation*}
where we used that by definition, see \eqref{Tg}, we have: 
$$\mathcal{T}_g=\lambda_q^2\tau^a\tau^c\ell^{-\alpha}\delta_{q+1}\geq \left(\frac{\lambda_q}{\lambda_{q+1}}\right)^{\gamma_a+2\beta+2\gamma_{CZ}}.$$

\noindent \noindent \textbf{Case $r\in I_2$.} If $r+r'\in I_2$, similarly as above we compute:
\begin{equation*}
    \begin{split}
        \left(\frac{\lambda_q}{\lambda_{q+1}}\right)^{r'}L_{p,\mathcal{A}}(r+r')&=\left(\frac{\lambda_q}{\lambda_{q+1}}\right)^{r'}\left[1+\left(\lambda_q^{(b-1)\gamma_\ell}\mathcal{T}_g\right)^{\underline{r}-1-(r+r')}\bar L\right]\\
        &\leq\left[\left(\frac{\lambda_q}{\lambda_{q+1}}\right)^{r'}+\left(\frac{\lambda_q}{\lambda_{q+1}}\right)^{r'\underbrace{[1-(\gamma_a+2\beta+\gamma_\ell+2\gamma_{CZ})]}_{\geq 0}}\left(\lambda_q^{(b-1)\gamma_\ell}\mathcal{T}_g\right)^{\underline{r}-1-r}\bar L\right]\\
        &\leq\left[1+\left(\lambda_q^{(b-1)\gamma_\ell}\mathcal{T}_g\right)^{\underline{r}-1-r}\bar L\right]\\
        &=L_{p,\mathcal{A}}(r).
    \end{split}
\end{equation*}
If $r+r'\in I_3$ since $r\in I_2$ we have $r'\geq \underline{r}-r\geq \underline{r}-1-r, \ r\leq \underline{r}-1$ and we compute:
\begin{equation*}
    \begin{split}
        \left(\frac{\lambda_q}{\lambda_{q+1}}\right)^{r'}L_{p,\mathcal{A}}(r+r')&=\left(\frac{\lambda_q}{\lambda_{q+1}}\right)^{r'}\lambda_q^{[r+r'-(\underline{r}-1)](b-1)\gamma_\ell}\bar L\\
        &=\left(\frac{\lambda_q}{\lambda_{q+1}}\right)^{r'(1-\gamma_\ell)}\lambda_q^{\underbrace{[r-(\underline{r}-1)](b-1)\gamma_\ell}_{\leq 0}}\bar L\\
        &\leq \left(\frac{\lambda_q}{\lambda_{q+1}}\right)^{(\underline{r}-1-r)(\gamma_a+2\beta+2\gamma_{CZ})}\bar L\\
        &\leq \mathcal{T}_g^{\underline{r}-1-r}\bar L\\
        &\leq 1+\left(\lambda_q^{(b-1)\gamma_\ell}\mathcal{T}_g\right)^{\underline{r}-1-r}\bar L=L_{p,\mathcal{A}}(r).
    \end{split}
\end{equation*}

\subsubsection{Estimates on the Lagrangian Perturbation.} 
Let $X^p:[-1,1]_s\times\mathbb{T}^3_x\times\mathbb{R}_t\to \mathbb{T}^3$ solve:
\begin{equation}\label{lagrangianp}
    \begin{cases}
        \partial_sX_s^p(x,t)=\xi^p(X_s^p(x,t),t),\\
        X_0^p(x,t)=x.
    \end{cases}
\end{equation} 
In this Subsection, we want to improve on the usual Gr{\"o}nwall bound:
$$||\DD X^p_s-\IId||_0\leq e^{\int_0^s||\DD \xi||_{0}\dd s}-1$$
this is needed because, $|s|\leq 1$ and $||\DD \xi^p||_0\lesssim \lambda_{q+1}\tau^a\delta_{q+1}^{1/2}$ blows-up along the scheme 
$$\tau^a\delta_{q+1}^{1/2}\gg \frac{1}{\lambda_{q+1}},$$
the key idea is that we have better control over the $k$-th momentum of $\xi^p$ for any $k\geq 1$, relative to its full differential, namely
$$||(\xi^p\cn)^k\xi^p||_0\ll ||\DD^k\xi^p||_0$$
this observation is useful because a Taylor expansion of $X^p_s$ in the parameter $s$ at $s=0$, gives:
\begin{equation}\label{taylorflow}
        X_s^p=\IId+\sum_{k=1}^{k_0}\frac{s^k}{k!}\partial_s^kX_s^p|_{s=0}+r^{k_0,s}_X=\IId+s\xi^p+\sum_{k=2}^{k_0}\frac{s^k}{k!}(\xi^p\cn)^{k-1}\xi^p+r^{k_0,s}_X
\end{equation}
where the remainder can be explicitly written as:
\begin{equation*}
        r^{k_0,s}_X=\frac{1}{k_0!}\int_0^s\partial_s^{k_0+1}X_{s'}^p(s-s')^{k_0}\dd s'=\frac{1}{k_0!}\int_0^s[(\xi^p\cn)^{k_0}\xi^p](X_{s'}^p)(s-s')^{k_0}\dd s'.
\end{equation*}
We now show that we indeed have better control over the $k$-th momentum. 
\begin{lemma}[Momentum Estimates] \label{momentumestimates} 
Let $\zeta$ be any constant unit vector and $\bar N\geq 2$ a non-negative integer. Set $\underline{r}=M-2m_0-k_0^g-8$. Then there exist constants $C, \ C'$ such that for 
$$j=0,1,2, \ \ \ \ 0\leq r+j+k \leq \bar N$$ 
the following bounds hold:
    \begin{equation*}
    \begin{split}
        &||\partial_t^j(\xi^p\cn)^k\xi^p||_r\leq C'(C)^k \lambda_{q+1}^{r+j}\lambda_q^{[k-\underline{r}]^+(b-1)\gamma_\ell}\mathcal{T}_p^k\tau^a\delta_{q+1}^{1/2}\\
        &||\partial_t^j\mathcal{L}^k_\xi\mathcal{L}_\xi\zeta||_r\leq C'(C)^k \lambda_{q+1}^{r+j}\lambda_q^{[k-\underline{r}]^+(b-1)\gamma_\ell}\mathcal{T}_p^k\mathcal{M}_p\\
    \end{split}
    \end{equation*}
    The constants, $C, \ C'$ depend on $\bar N$ and all the parameters but not on $a$ and are uniform in $k, \ r, \ j, \ \zeta $.
\end{lemma}
\begin{proof}[Proof of Lemma \ref{momentumestimates}] The second statement is a simple application of the Inductive Lemma \ref{inductive}. We omit the details, as the same type of argument will appear several times in the proof of Lemma \ref{estimatesperturbationnash}. 

\noindent The first statement follows from a slight variation of the proof of Lemma \ref{inductive}. By means of the disjointedness properties from Lemma \ref{disjoint} and the definition of $\xi_I^p$ in \eqref{xinash}, we first write:
\begin{equation*}
    \begin{split}
        &(\xi^p\cn)^k \xi^p=\sum_I (\xi_I^p\cn)^k \xi_I^p,\\
        &\xi_I^p=\tau^a\tilde a_I\alpha_{I}\tilde \Psi^{2*}_I\left[\varphi_{\lambda_{q+1},k_I}'\zeta_I\right]+\frac{\tau^a}{\lambda_{q+1}}\alpha_{I}\nabla\tilde a_I\times\tilde\Psi_I^{1*}(\varphi_{\lambda_{q+1},k_I}\nu_I).
    \end{split}
\end{equation*}
We prove the claim by induction; the $k=0$ case follows from Lemmas \ref{slowcoeffestimates}, \ref{chartprop} and the representation above, now assume that: 
$$F_{I,k}=(\xi^p_I\cn)^k \xi^p_I\simeq \tilde\Psi_I^{*}[(\phi_k)_{\lambda_{q+1},k_I}]Y_{I,k}=\phi_k(\lambda_{q+1}\tilde\Psi_I\cdot k_I)Y_{I,k}$$
with
\begin{equation*}||\partial_t^jY_{I,k}||_r\lesssim \lambda_q^{r+j}\lambda_q^{[r+k+j-\underline{r}]^+(b-1)\gamma_\ell}\mathcal{T}_p^k\tau^a\delta_{q+1}^{1/2} \ \text{ and } \ ||\phi_k||_r\lesssim 1
\end{equation*}
where we use the same convention for $\simeq$ used in the proof of Lemma \ref{inductive} to indicate a finite sum of terms of that structure, and obeying the same estimates. It will be apparent that the number of these grows at most exponentially in $k$, and this is the reason for the shape of the implicit constant, claimed in the statement of the Lemma.

\noindent We are interested in the bounds for $F_{I,k+1}$, we compute:
\begin{equation*}
    \begin{split}
        F_{I,k+1}&=\xi_I\cn F_{I,k}\\
        &=\left[\tau^a\tilde a_I\alpha_{I}\tilde \Psi^{2*}_I\left[\varphi_{\lambda_{q+1},k_I}'\zeta_I\right]+\frac{\tau^a}{\lambda_{q+1}}\alpha_{I}\nabla\tilde a_I\times\tilde\Psi_I^{1*}(\varphi_{\lambda_{q+1},k_I}\nu_I)\right]\cn\left[ \phi_k(\lambda_{q+1}\tilde\Psi_I\cdot k_I)Y_{I,k}\right]\\
        &=\left[\lambda_{q+1}\tilde \Psi_I^*[(\phi_k')_{\lambda_{q+1},k_I}]Y_{I,k}\otimes \tilde \Psi_I^{1*} k_I+\Psi_I^*[(\phi_k)_{\lambda_{q+1},k_I}]\DD Y_{I,k}\right]\left[\tau^a\tilde a_I\alpha_{I}\tilde \Psi^{2*}_I\left[\varphi_{\lambda_{q+1},k_I}'\zeta_I\right]\right]\\
        &+\left[\lambda_{q+1}\tilde \Psi_I^*[(\phi_k')_{\lambda_{q+1},k_I}]Y_{I,k}\otimes \tilde \Psi_I^{1*} k_I+\Psi_I^*[(\phi_k)_{\lambda_{q+1},k_I}]\DD Y_{I,k}\right]\left[\frac{\tau^a}{\lambda_{q+1}}\alpha_{I}\nabla\tilde a_I\times\tilde\Psi_I^{1*}(\varphi_{\lambda_{q+1},k_I}\nu_I)\right]\\
        &=\Psi_I^*[(\varphi'\phi_k)_{\lambda_{q+1},k_I}]\tau^a\tilde a_I\alpha_{I}(\tilde \Psi^{2*}_I\zeta_I\cn) Y_{I,k}\\
        &+\Psi_I^*[(\varphi\phi_k')_{\lambda_{q+1},k_I}]\tau^a\alpha_{I}(\nabla\tilde a_I\times\tilde\Psi_I^{1*}\nu_I)\cdot \Psi_I^{1*}k_IY_{I,k}+\Psi_I^*[(\varphi\phi_k)_{\lambda_{q+1},k_I}]\frac{\tau^a}{\lambda_{q+1}}\alpha_{I}(\nabla\tilde a_I\times\tilde\Psi_I^{1*}\nu_I)\cn Y_{I,k}
    \end{split}
\end{equation*}
The key point here is that the main term in $\xi_I^p$, because of its geometry, kills the $\lambda_{q+1}$ term in $\DD F_k$, while the corrector term has a $\lambda_q/\lambda_{q+1}$ embedded smallness compensating for the loss. From the computation above, it is clear that $F_{I,k+1}$ can be written as:
$$F_{I,k+1}\simeq \tilde\Psi_I^{*}[(\phi_{k+1})_{\lambda_{q+1},k_I}]Y_{I,k+1}$$
and from Lemmas \ref{chartprop}, \ref{slowcoeffestimates} it follows that:
\begin{equation*}
    ||Y_{I,k+1}||_r\lesssim \lambda_q^r\lambda_q^{[k+1-\underline{r}]^+(b-1)\gamma_\ell}\mathcal{T}_p^{k+1}\tau^a\delta_{q+1}^{1/2} \ \text{ and } \ ||\phi_{k+1}||_r\lesssim 1.
\end{equation*}
Commuting time and space derivatives from the calculations above and the estimates in Lemmas \ref{chartprop}, \ref{slowcoeffestimates}, again, we conclude:
$$||\partial_t^j Y_{I,k+1}||_r\lesssim \lambda_q^{r+j}\lambda_q^{[k+1+j-\underline{r}]^+(b-1)\gamma_\ell}\mathcal{T}_p^{k+1}\tau^a\delta_{q+1}^{1/2} \ \text{ for }\ j=0,1,2 \ \ \& \  \ 0\leq r+(k+1)+j\leq \bar N.$$
In particular, by means of the composition estimates in Proposition \ref{compestimates} and Lemma \ref{chartprop}, for $r\geq 1$ we deduce:
\begin{equation*}
    \begin{split}
        ||F_{k,I}||_r&\lesssim ||\tilde\Psi_I^{*}[(\phi_k)_{\lambda_{q+1},k_I}]Y_{I,k}||_r\\
        &\lesssim ||\tilde\Psi_I^{*}[(\phi_k)_{\lambda_{q+1},k_I}]||_r||Y_{I,k}||_0+||\tilde\Psi_I^{*}[(\phi_k)_{\lambda_{q+1},k_I}]||_0||Y_{I,k}||_r\\
        &\lesssim \left[\lambda_{q+1}||\phi_k||_1||\DD\tilde \Psi_I||_{r-1}+\lambda_{q+1}^r||\phi_k||_r||\DD\tilde \Psi_I||_0^r\right]||Y_{I,k}||_0\\
        &\lesssim \lambda_{q+1}^r\lambda_q^{[k-\underline{r}]^+(b-1)\gamma_\ell}\mathcal{T}^k_p\tau^a\delta_{q+1}^{1/2},
    \end{split}
\end{equation*}
similarly, one can show:
$$||\partial_t^j F_{k,I}||_r\lesssim \lambda_{q+1}^{r+j}\lambda_q^{[k-\underline{r}]^+(b-1)\gamma_\ell}\mathcal{T}^k_p\tau^a\delta_{q+1}^{1/2}.$$

\noindent From the fact that for each spacetime point there is at most one $\xi_I^p$ which is not zero, and the bound we just showed, we conclude that:
\begin{equation*}
    \begin{split}
        ||\partial_t^j(\xi^p\cn)^k \xi^p||_r&\leq \sup_I ||\partial_t^j(\xi_I^p\cn)^k \xi_I^p||_r\\
        &=\sup_I ||\partial_t^jF_{k,I}||_r\\
        &\lesssim \lambda_{q+1}^{r+j}\lambda_q^{[k-\underline{r}]^+(b-1)\gamma_\ell}\mathcal{T}^k_p\tau^a\delta_{q+1}^{1/2}
    \end{split}
\end{equation*}
where the implicit constant is of the form $C'(C)^{k+1}$ for some $C, \ C'$ which depend on $\bar N$ and all the other parameters but not on $a$ and are uniform in $k, \ r, \ j$.
\end{proof}

We are ready to prove bounds for the Lagrangian perturbation.
\begin{lemma}[Estimates on $X^{p}$] \label{estimatesflowp} Let $\bar N$ be a non-negative integer and $\zeta$ any constant unit vector. For $j=0,1,2$ and any $|s|\leq 1$, we have:
    \begin{equation*}
        \begin{split}
            &||\partial_t^j X^{p}_s||_r\lesssim\lambda_{q+1}^{r+j}\tau^a\delta_{q+1}^{1/2} \ \ \text{ for } \ \ 0\leq r\leq \bar N+1-j \ \ \& \ \ (r,j)\neq (0,0)\\
            &||\partial_t^j\left[\DD X^p\circ (X^p)^{-1}[\zeta]\right]||_r\lesssim \lambda_{q+1}^{r+j+1}\tau^a\delta_{q+1}^{1/2}=\lambda_{q+1}^{r+j}\mathcal{M}_p \ \text{ for } \ 0\leq r\leq \bar N-j\\
        \end{split}
    \end{equation*}
\noindent The implicit constants depend on $\bar N, \ r$ and all the other parameters but not on $a$.
\end{lemma}

What happens here is that, thanks to Lemma \ref{momentumestimates} and the assumptions, we have:
\begin{equation*}
    \begin{split}
       ||(\xi^p\cn)^{k_0}\xi^p||_{1}\lesssim\left(\lambda_q^{(b-1)\gamma_\ell}\mathcal{T}_p\right)^{k_0}\mathcal{M}_p
    \end{split}
\end{equation*}
and in particular, choosing $k_0$ sufficiently large, we can ensure
$$\frac{1}{k_0!}||(\xi^p\cn)^{k_0}\xi^p||_1\leq 1,$$ 
this allows us to run the usual Gr{\"o}nwall argument to control derivatives of the Lagrangian flow of a vector field, but at the level of the first derivative of the $k_0$ momentum of the vector field instead of the derivative of the vector field itself. 

\begin{proof}[Proof of Lemma \ref{estimatesflowp}] We first recall the expressions of $\mathcal{T}_p$ and $\mathcal{M}_p$ from \eqref{Tp}. Now, given any $b, \ \gamma_a, \ \gamma_\ell, \ \bar N$ we can find $k_0$ so large that for $0\leq r \leq \bar N$ we have:
\begin{equation}\label{proliminarieslagrangiann}
    \left(\lambda_q^{(b-1)\gamma_\ell}\mathcal{T}_p\right)^{k_0}\mathcal{M}_p^{r+2}\leq\left(\lambda_q^{(b-1)\gamma_\ell}\mathcal{T}_p\right)^{k_0}\mathcal{M}_p^{\bar N+2}\leq \left(\frac{\lambda_q}{\lambda_{q+1}}\right)^{k_0(\beta+\gamma_a-\gamma_\ell)+(\bar N+2)(\beta+\gamma_a-\gamma_\ell-1)}\overset{!}{\leq}1
\end{equation}
We now let $C, C'$ be the implicit constants from Lemma \ref{momentumestimates}, applied with $k_0+\bar N+3$ being the total number of derivatives on which we want the estimates (this number is called $\bar N$ in the statement of \ref{momentumestimates}, not to be confused with the one we are given as data now). In this sense, we say that $C, \ C'$ depend on $\bar N$ and all the other parameters.

\noindent For notational convenience we assume $0\leq s\leq 1$ the case $-1\leq s\leq 0$ can be handled similarly. The composition estimates in Proposition \ref{compestimates} will be used repeatedly without mention.

\noindent \textit{Estimates on the Lagrangian Perturbation.} By taking one derivative of the expression of $X_s^p-\IId$ given in \eqref{taylorflow} we can estimate:
\begin{equation}\label{lagrangianbasecase}
    \begin{split}
        ||\DD X_s^p-\IId||_0
        &\leq \sum_{k=1}^{k_0}\frac{s^{k}}{k!}\left[(\xi^p\cn)^{k-1}\xi^p\right]_1+\frac{1}{k_0!}\int_0^s\left[(\xi^p\cn)^{k_0}\xi^p\right]_1||\DD X_{s'}^p-\IId||_0(s-s')^{k_0}\dd s'\\
        &+\frac{1}{k_0!}\int_0^s\left[(\xi^p\cn)^{k_0}\xi^p\right]_1(s-s')^{k_0}\dd s'\\
        &\lesssim \mathcal{M}_p\sum_{k=1}^{k_0+1}\frac{s^k}{k!}(C)^{k-1}\lambda_q^{[k-1-\underline{r}]^+(b-1)\gamma_\ell}\mathcal{T}_p^{k-1}\\
        &+\lambda_q^{[k_0-\underline{r}]^+(b-1)\gamma_\ell}\underbrace{\frac{1}{k_0!}(C)^{k_0}\mathcal{T}_p^{k_0}\mathcal{M}_p}_{\lesssim 1}\int_0^s||\DD X_{s'}^p-\IId||_0(s-s')^{k_0}\dd s'\\
        &\lesssim \mathcal{M}_p\underbrace{\sum_{k=+1}^{k_0}\frac{\left(C\lambda_q^{(b-1)\gamma_\ell}\mathcal{T}_p\right)^{k-1}}{k!}}_{\lesssim 1} +\int_0^s\left[(\xi^p\cn)^{k_0}\xi^p\right]_1||\DD X_{s'}^p-\IId||_0(s-s')^{k_0}\dd s'\\
        &\lesssim \mathcal{M}_p+\int_0^s||\DD X_{s'}^p-\IId||_0(s-s')^{k_0}\dd s'\\
    \end{split}
\end{equation}
An application of Gr{\"o}nwall Lemma gives:
\begin{equation*}
    \begin{split}
        ||\DD X_s^p-\IId||_0&\lesssim \mathcal{M}_p
    \end{split}
\end{equation*}
and the case $r=0$ for pure space derivatives is proved.

\noindent For the $r\geq 1$ case, we can bound the remainder as follows:
\begin{equation*}
    \begin{split}
        ||r_X^{k_0,s}||_{r}&\leq\frac{1}{k_0!}\int_0^s||(\xi^p\cn)^{k_0}\xi^p](X_{s'}^p)||_{r}\dd s'\\
        &\leq\frac{1}{k_0!}\int_0^s\left[[(\xi^p\cn)^{k_0}\xi^p]_1||\DD X_{s'}^p||_{r-1}+||(\xi^p\cn)^{k_0}\xi^p||_{r+1}||\DD X_{s'}^p||_0^r\right]\dd s'\\
        &\lesssim \frac{1}{k_0!}(C)^{k_0}\mathcal{T}_p^{k_0}\mathcal{M}_p\int_0^s||X_{s'}^p||_{r}\dd s'+\frac{1}{k_0!}(C)^{k_0}\lambda_{q+1}^{r}\left(\lambda_q^{(b-1)\gamma_\ell}\mathcal{T}_p\right)^{k_0}\mathcal{M}_p^{r+1}\\
        &\leq \frac{1}{k_0!}(C)^{k_0}\int_0^s||X_{s'}^p||_{r}\dd s'+\frac{1}{k_0!}(C)^{k_0}\lambda_{q+1}^{r}\mathcal{M}_p
    \end{split}
\end{equation*}
where we used Lemma \ref{momentumestimates} together with the case $r=0$ and \eqref{proliminarieslagrangiann}, we conclude that:
\begin{equation*}
    \begin{split}
        ||\DD X_s^p||_r&\leq \sum_{k=1}^{k_0}\frac{1}{k!}||(\xi^p\cn)^{k-1}\xi^p||_{r+1}+||r^{k_0,s}_X||_{r+1}\\
        &\lesssim \lambda_{q+1}^{r}\mathcal{M}_p\underbrace{\sum_{k=1}^{k_0}\frac{1}{k!}(C)^{k-1}\lambda_q^{[k-1-\underline{r}]^+(b-1)\gamma_\ell}\mathcal{T}^{k-1}_p}_{\lesssim 1}+\int_0^s||X_{s'}^p||_{r}\dd s'\\
        &\lesssim \lambda_{q+1}^{r}\mathcal{M}_p+\int_0^s||X_{s'}^p||_{r}\dd s'\\
    \end{split}
\end{equation*}
for $0\leq r \leq \bar N$, where we argued as above, to address the sum. We now apply Gr{\"o}nwall Lemma and conclude: 
\begin{equation*}
        ||\DD X_s^p||_r\lesssim \lambda_{q+1}^{r}\mathcal{M}_p \ \text{ for } \ 0\leq s\leq 1 \ \ \& \ \ 1\leq r\leq \bar N.
\end{equation*}

\noindent For the pure time derivatives case, we first compute:
\begin{equation}\label{timederivativeremindernash}
    \begin{split}
        \partial_tr^{k_0,s}_X&=\frac{1}{k_0!}\int_0^s\left[\partial_t[(\xi^p\cn)^{k_0}\xi^p](X_{s'}^p)+\DD[(\xi^p\cn)^{k_0}\xi^p](X_{s'}^p)\left[\partial_t X_{s'}^p\right]\right] (s-s')^{k_0}\dd s',
        \\
        \partial_t^2r^{k_0,s}_X&=\frac{1}{k_0!}\int_0^s\left[\partial_t^2[(\xi^p\cn)^{k_0}\xi^p](X_{s'}^p)+\DD^2[(\xi^p\cn)^{k_0}\xi^p](X_{s'}^p)\left[\partial_tX_{s'}^p,\partial_tX_{s'}^p\right]\right](s-s')^{k_0}\dd s'\\
        &+\frac{1}{k_0!}\int_0^s\left[2\DD[\partial_t[(\xi^p\cn)^{k_0}\xi^p]](X_{s'}^p)\left[\partial_t X_{s'}^p\right]+\DD[(\xi^p\cn)^{k_0}\xi^p](X_{s'}^p)\left[\partial_t^2 X_{s'}^p\right]\right](s-s')^{k_0}\dd s'.
    \end{split}
\end{equation}
From the bounds in Lemma \ref{momentumestimates}, we deduce:
\begin{equation*}
    \begin{split}
        ||\partial_tr^{k_0,s}_X||_{r}&\lesssim \frac{1}{k_0!}\int_0^s||\partial_t[(\xi^p\cn)^{k_0}\xi^p](X_{s'}^p)||_r(s-s')^{k_0}\dd s'\\
        &+\frac{1}{k_0!}\int_0^s\left[||\DD[(\xi^p\cn)^{k_0}\xi^p](X_{s'}^p)||_{r}||\partial_t X_{s'}^p||_0+||\DD[(\xi^p\cn)^{k_0}\xi^p](X_{s'}^p)||_{0}||\partial_t X_{s'}^p||_{r}\right] (s-s')^{k_0}\dd s'\\
        &\overset{r=0}{\leq}\frac{s^{k_0+1}}{(k_0+1)!}C^{k_0}\left(\lambda_q^{(b-1)\gamma_\ell}\mathcal{T}_p\right)^{k_0}\mathcal{M}_p+\frac{1}{k_0!}C^{k_0}\left(\lambda_q^{(b-1)\gamma_\ell}\mathcal{T}_p\right)^{k_0}\mathcal{M}_p\int_0^s||\partial_t X_{s'}^p||_0(s-s')^{k_0}\dd s\\
        &\overset{r\geq 1}{\lesssim}\frac{1}{k_0!}\int_0^s\left[||\partial_t[(\xi^p\cn)^{k_0}\xi^p||_{1}||\DD X_{s'}^p||_{r}+||\partial_t[(\xi^p\cn)^{k_0}\xi^p||_{r}||\DD X_{s'}^p||_{0}^{r}\right](s-s')^{k_0}\dd s'\\
        &+\frac{1}{k_0!}\int_0^s\left[||\DD[(\xi^p\cn)^{k_0}\xi^p]||_1||\DD X_{s'}^p||_{r-1}+||\DD[(\xi^p\cn)^{k_0}\xi^p]||_r||\DD X_{s'}^p||_{0}^r\right]||\partial_t X_{s'}^p||_0 (s-s')^{k_0}\dd s'\\
        &+\frac{1}{k_0!}\int_0^s||\DD[(\xi^p\cn)^{k_0}\xi^p]||_0||\partial_t X_{s'}^p||_r (s-s')^{k_0}\dd s'\\
        &\lesssim \frac{s^{k_0+1}}{(k_0+1)!}C^{k_0}\lambda_{q+1}^{r}\left(\lambda_q^{(b-1)\gamma_\ell}\mathcal{T}_p\right)^{k_0}\mathcal{M}_p^{r+1}+\frac{1}{k_0!}C^{k_0}\lambda_{q+1}^{r}\left(\lambda_q^{(b-1)\gamma_\ell}\mathcal{T}_p\right)^{k_0}\mathcal{M}_p^{r+1}\int_0^s||\partial_t X_{s'}^p||_0 (s-s')^{k_0}\dd s'\\
        &+\frac{1}{k_0!}C^{k_0}\left(\lambda_q^{(b-1)\gamma_\ell}\mathcal{T}_p\right)^{k_0}\mathcal{M}_p\int_0^s||\partial_t X_{s'}^p||_r (s-s')^{k_0}\dd s'\\
    \end{split}
\end{equation*}
and given our choice of $k_0$ in \eqref{proliminarieslagrangiann} we conclude:
\begin{equation*}
    \begin{split}
        ||\partial_tr^{k_0,s}_X||_{r}&\lesssim \begin{cases}
        \mathcal{M}_p+\int_0^s||\partial_t X_{s'}^p||_0 (s-s')^{k_0}\dd s' & \ \text{ for } \ r=0,\\
            \lambda_{q+1}^r\mathcal{M}_p+\lambda_{q+1}^r\int_0^s||\partial_t X_{s'}^p||_0 (s-s')^{k_0}\dd s'+\int_0^s||\partial_t X_{s'}^p||_r (s-s')^{k_0}\dd s' & \ \text{ for } \ r\geq 1.
        \end{cases}
    \end{split}
\end{equation*}
We now take the time derivative of the full expression in \eqref{taylorflow}, we plug the bound just found and the ones in Lemma \ref{momentumestimates} and deduce:
\begin{equation*}
    \begin{split}
        ||\partial_tX_s^p||_r&\leq \sum_{k=1}^{k_0}\frac{s^k}{k!}||\partial_t\left[(\xi^p\cn)^{k-1}\xi^p\right]||_r+||\partial_tr^{k_0,s}_X||_r\\
        &\lesssim \lambda_{q+1}^r\mathcal{M}_p+\lambda_{q+1}^{r}\int_0^s||\partial_t X_{s'}^p||_0 (s-s')^k\dd s'+\int_0^s||\partial_t X_{s'}^p||_r (s-s')^k\dd s'\\
    \end{split}
\end{equation*}
where we argued, as in \eqref{lagrangianbasecase}, to deal with the sum in the parameter $k$. An application of Gr{\"o}nwall Lemma, first for the case $r=0, \ 0\leq s\leq 1$ and then for the case $0\leq r \leq \bar N-1, \ 0\leq s\leq 1$,  allows us to conclude that:
$$||\partial_tX_s^p||_r\lesssim \lambda_{q+1}^{r}\mathcal{M}_p=\lambda_{q+1}^{r+1}\tau^a\delta_{q+1}^{1/2} \ \text{ for } \ 0\leq s\leq 1 \ \ \& \ \ 1\leq r\leq \bar N-1,$$
now arguing as above, but making use of the expression for $\partial_t^{2}$ in \eqref{timederivativeremindernash} instead of $\partial_t$ one can show:
$$||\partial_t^2X_s^p||_r\lesssim \lambda_{q+1}^{r+1}\mathcal{M}_p=\lambda_{q+1}^{r+2}\tau^a\delta_{q+1}^{1/2} \ \text{ for } \ 0\leq s\leq 1 \ \ \& \ \ 1\leq r\leq \bar N-2$$
where the implicit constants in the bounds above depend on $\bar N, \ r$ and all the other parameters, but not on $a$ and $s$ for $0\leq s\leq 1$.

\noindent \textit{Estimates on the push-forward of a constant vector.} To prove the second assertion of the Lemma, we argue similarly. For any constant vector $\zeta$ we have:
\begin{equation*}
\begin{split}
    \DD X^p((X^p)^{-1})[\zeta]&=(X^p)_*\zeta\\
    &=\zeta+\sum_{k=1}^{k_0}\frac{(-1)^k}{k!}\mathcal{L}_{\xi^p}^k\zeta+\frac{(-1)^{k_0+1}}{k_0!}\int_{0}^1(X_{s}^p)_*\mathcal{L}_{\xi^p}^{k_0+1}\zeta (1-s)^{k_0}\dd s\\
    &=\zeta+\sum_{k=1}^{k_0}\frac{(-1)^k}{k!}\mathcal{L}_{\xi^p}^k\zeta+\frac{(-1)^{k_0+1}}{k_0!}\int_{0}^1\DD X_s^p(X_{-s}^p)[\mathcal{L}_{\xi^p}^{k_0+1}\zeta (X_{-s}^p)] (1-s)^{k_0}\dd s\\
\end{split}
\end{equation*}
where the expansion can be deduced from the ideas behind Lemma \ref{lpl} and we used that since $\xi^p$ is independent of $s$, $X_s$ has the semi-group property $X_s\circ X_{s'}=X_{s+s'}$ and in particular we have $X_s^{-1}=X_{-s}$.  

\noindent From the results above for $1\leq r \leq \bar N$ we deduce:
\begin{equation*}
    \begin{split}
        ||[\DD X_s^p(X_{-s}^p)]||_r&\lesssim[\DD X^p]_1||\DD X_{-s}^p||_{r-1}+||\DD X^p||_r||\DD X_{-s}^p||_0^{r}\\
        &\lesssim \lambda_{q+1}^{r}\mathcal{M}^{r+1}_p,
        \\
        ||(\mathcal{L}_{\xi^p}^{k_0+1}\zeta)(X^p_{-s})||_r&\lesssim [\mathcal{L}_{\xi}^{k_0+1}\zeta]_1||\DD X^p_{-s}||_{r-1}+||\mathcal{L}_{\xi}^{k_0+1}\zeta||_r||\DD X^p_{-s}||_0^{r}\\
        &\lesssim \lambda_{q+1}^{r}\left(\lambda_q^{(b-1)\gamma_\ell}\mathcal{T}_p\right)^{k_0}\mathcal{M}^{r+1}_p,\\
    \end{split}
\end{equation*}
moreover,
\begin{equation*}
    \begin{split}
        ||\partial_t\DD X_s^p(X_{-s}^p)||_r&\leq ||\DD (\partial_tX_s^p)(X_{-s}^p)||_r+||\DD^2 X_s^p(X_{-s}^p)[\partial_t X_{-s}^p,\cdot]||_r\\
        &\lesssim[\DD (\partial_tX_s^p)]_1||\DD X_{-1}^p||_{r-1}+||\DD (\partial_tX_s^p)||_r||\DD X_{-1}^p||_0^{r}\\
        &+||\DD^2 X_s^p(X_{-s}^p)||_r||\partial_t X_{-s}^p||_0+||\DD^2 X_s^p(X_{-s}^p)||_0||\partial_t X_{-s}^p||_r\\
        &\lesssim \lambda_{q+1}^{r+1}\mathcal{M}_p^{r+1}+\left[||\DD^2 X_s^p||_1||\DD X_{-s}^p||_{r-1}+||\DD^2 X_s^p||_r||\DD X_{-s}^p||_0^r\right]\mathcal{M}_p\\
        &\lesssim \lambda_{q+1}^{r+1}\mathcal{M}_p^{r+2},
        \\
        ||\partial_t[(\mathcal{L}_{\xi^p}^{k_0+1}\zeta)(X^p_{-s})]||_r&\lesssim [\partial_t\mathcal{L}_{\xi}^{k_0+1}\zeta]_1||\DD X^p_{-s}||_{r-1}+||\partial_t\mathcal{L}_{\xi}^{k_0+1}\zeta||_r||\DD X^p_{-s}||_0^{r}+||\DD[\mathcal{L}_{\xi}^{k_0+1}\zeta](X_{-s}^p)[\partial_t X_{-s}^p]||_r\\
        &\lesssim \lambda_{q+1}^{r+1}(C)^{k_0}\mathcal{M}^{r+1}_p\left(\lambda_q^{(b-1)\gamma_\ell}\mathcal{T}_p\right)^{k_0}+||\DD[\mathcal{L}_{\xi}^{k_0+1}\zeta||_0||\partial_t X_{-s}^p||_r\\
        &+\left[[\DD[\mathcal{L}_{\xi}^{k_0+1}\zeta]_1||\DD X_{-s}^p||_{r-1}+||\DD[\mathcal{L}_{\xi}^{k_0+1}\zeta||_r||\DD X_{-s}^p||_0^r\right]||\partial_t X_{-s}^p||_0\\
        &\lesssim (C)^{k_0}\lambda_{q+1}^{r+1}\mathcal{M}^{r+2}_p\left(\lambda_q^{(b-1)\gamma_\ell}\mathcal{T}_p\right)^{k_0},\\
    \end{split}
\end{equation*}
similarly,
\begin{equation*}
    \begin{split}
        ||\partial_t^2\DD X_s^p(X_{-s}^p)||_r&\lesssim \lambda_{q+1}^{r+2}\mathcal{M}_p^{r+3},\\
        ||\partial_t^2(\mathcal{L}_{\xi^p}^{k_0+1}\zeta)(X^p_{-s})||_r&\lesssim (C)^{k_0}\lambda_{q+1}^{r+2}\mathcal{M}^{r+3}_p\left(\lambda_q^{(b-1)\gamma_\ell}\mathcal{T}_p\right)^{k_0}.\\
    \end{split}
\end{equation*}
Given our choice of $k_0$ in \eqref{proliminarieslagrangiann}, we deduce that:
\begin{equation*}
    \begin{split}
        ||\partial_t^j[\DD X_s^p(X_{-s}^p)(\mathcal{L}_{\xi^p}^{k_0+1}\zeta)(X^p_{-s})]||_r\lesssim (C)^{k_0}\lambda_{q+1}^{r+j}\mathcal{M}_p \ \text{ for } \ 0\leq r\leq \bar N-j.
    \end{split}
\end{equation*}
We are ready to estimate:
\begin{equation*}
    \begin{split}
        ||\partial_t^j[\DD X^p \circ (X^p)^{-1}\zeta-\zeta]||_r&\lesssim\sum_{k=1}^{k_0}\frac{1}{k!}||\partial_t^j[\mathcal{L}_{\xi^p}^k\zeta]||_r+\frac{1}{k_0!}\int_{0}^1||\partial_t^j[\DD X_s^p(X_{-s}^p)[\mathcal{L}_{\xi^p}^{k_0+1}\zeta (X_{-s}^p)]]||_r (1-s)^{k_0}\dd s\\
        &\lesssim \lambda_{q+1}^{r+j}\mathcal{M}_p\sum_{k=0}^{k_0-1}\frac{1}{(k+1)!}(C)^{k}\lambda_q^{[k-(\underline{r}-1)]^+(b-1)\gamma_\ell}\mathcal{T}^k_p+\frac{1}{(k_0+1)!}C^{k_0}\lambda_{q+1}^{r+j}\mathcal{M}_p\\
        &\lesssim \lambda_{q+1}^{r+j}\mathcal{M}_p\\
    \end{split}
\end{equation*}
for $0\leq r \leq \bar N-j$. The implicit constant depends on $\bar N, \ r$ and all the other parameters, but not on $a$. 
\end{proof}

\begin{remark}
    The relevance of the second claim in Lemma \ref{estimatesflowp} can be seen from the following computation. We will write $|\cdot|_{op}$ for the operator norm on matrices, we have:
\begin{equation*}
    \begin{split}
        [\DD X \circ X^{-1}]_r&=\sum_{|\theta|=r}\sup_{t\in \mathbb{R}}\sup_{x\in\mathbb{T}^3}|\partial_\theta\DD X \circ X^{-1}|_{op}\\
        &=\sum_{|\theta|=r}\sup_{t\in \mathbb{R}}\sup_{x\in\mathbb{T}^3}\sup_{|\zeta|=1}|\partial_\theta\DD X \circ X^{-1}[\zeta]|\\
        &=\sum_{|\theta|=r}\sup_{t\in \mathbb{R}}\sup_{x\in\mathbb{T}^3}\sup_{|\zeta|=1}|\partial_\theta(\DD X \circ X^{-1}[\zeta])|.\\
    \end{split}
\end{equation*}
It thus suffices to estimate $|\partial_\theta(\DD X \circ X^{-1}\zeta)|$ independently of $\zeta$ with $|\zeta|=1$ as in the Lemma \ref{estimatesflowp} to deduce bounds for $||\DD X \circ X^{-1}||_r$.
\end{remark}

\subsubsection{Estimates on the Lie-Taylor Expansion.}
In this Subsection, we provide bounds on all the objects in the decomposition of the perturbation from Subsection \ref{splittingnash}. The key tool for studying these Lie-Taylor expansions is the Inductive Lemma \ref{inductive}, which guarantees estimates for all $k\geq 1$ from those for $k=0$ and careful bookkeeping of the geometry. We invite the reader to visit Section \ref{inductives}, where this machinery is developed, before proceeding.

\begin{lemma}[Lie-Derivative and Transport Estimates - Mollified Terms]\label{estimatesperturbationnash} Fix $\bar N\geq 0$ an integer, let $j, \ k, \ r\geq 0$ also integers, $I\in \mathcal{I}$ and $\underline{r}=M-2m_0-k_0^g-8$. There exist constants $C,\ C'$ which depend on $\bar N$ and all the parameters but not on $a$ and are uniform in $I, \ j, \ k, \ r$, such that the following estimates hold.

\noindent\textbf{Lie-derivatives bounds.}  For $j=0,1,2$ and $0\leq r+k+j \leq \bar N$, we have:
\begin{equation*}
    \begin{split}
        &||\partial_t^j\mathcal{L}_{\xi^{p}_I}^k(\partial_t+\mathcal{L}_{\tilde v_{\ell,I}})\Theta_I^{p}||_{r}\leq C'(C)^k \lambda_{q+1}^{r+j-1}\lambda_q^{[k-\underline{r}]^+(b-1)\gamma_\ell}\mathcal{T}_p^k\delta_{q+1}^{1/2}, \\
        &||\partial_t^j\mathcal{L}_{\xi^{p}_I}^k\mathcal{L}_{\tilde B_{\ell,I}}\Theta_I^{p}||_r\leq C'(C)^k \lambda_{q+1}^{r+j-1}\lambda_q^{[k-\underline{r}]^+(b-1)\gamma_\ell}\mathcal{T}_p^k(\tau^a/\tau^c)\delta_{q+1}^{1/2}, \\
        &||\partial_t^j\mathcal{L}_{\xi^{p}_I}^k\mathcal{L}_{\tilde v_\ell-\tilde v_{\ell,I}}\Theta_I^{p}||_{r}, \ ||\partial_t^j\mathcal{L}_{\xi^{p}_I}^k\mathcal{L}_{\tilde B_\ell-\tilde B_{\ell,I}}\Theta_I^{p}||_{r}\leq C'(C)^k\lambda_{q+1}^{r+j}\lambda_q^{[k-\underline{r}]^+(b-1)\gamma_\ell}\mathcal{T}_p^k(\ell\lambda_q)^{m_0}\tau^a\delta_q^{1/2}\delta_{q+1}^{1/2}.
    \end{split}
\end{equation*}

\noindent \textbf{Alfv\'en-transport bounds.} For  $0\leq r+k \leq \bar N-1$, we have:
\begin{equation*}
    \begin{split}
        &||\mathcal{\tilde A}_{\ell,I}^\pm\mathcal{L}_{\xi^{p}_I}^k(\partial_t+\mathcal{L}_{\tilde v_{\ell,I}})\Theta_I^{p}||_r\leq C'(C)^k \lambda_{q+1}^{r-1}\lambda_q^{[k-\underline{r}]^+(b-1)\gamma_\ell}\mathcal{T}_p^k\delta_{q+1}^{1/2}\left[\frac{1}{\tau^a}+\frac{1}{\ell_t}(\tau^a/\tau^c)\right], \\
        &||\mathcal{\tilde A}_{\ell,I}^\pm \mathcal{L}_{\xi^{p}_I}^k\mathcal{L}_{\tilde B_{\ell,I}}\Theta_I^{p}||_r\leq C'(C)^k \lambda_{q+1}^{r-1}\lambda_q^{[k-\underline{r}]^+(b-1)\gamma_\ell}\mathcal{T}_p^k(\tau^a/\tau^c)\delta_{q+1}^{1/2}\left[\frac{1}{\tau^a}+\frac{1}{\ell_t}\right], \\
        &||\mathcal{\tilde A}_{\ell,I}^\pm\mathcal{L}_{\xi^{p}_I}^k\mathcal{L}_{\tilde v_\ell-\tilde v_{\ell,I}}\Theta_I^{p}||_r, \ ||\mathcal{\tilde A}_{\ell,I}^\pm\mathcal{L}_{\xi^{p}_I}^k\mathcal{L}_{\tilde B_\ell-\tilde B_{\ell,I}}\Theta_I^{p}||_r\leq C'(C)^k\lambda_{q+1}^{r}\lambda_q^{[k-\underline{r}]^+(b-1)\gamma_\ell}\mathcal{T}_p^k(\ell\lambda_q)^{m_0}\delta_q^{1/2}\delta_{q+1}^{1/2}.
    \end{split}
\end{equation*}
\end{lemma}

\begin{proof}[Proof of Lemma \ref{estimatesperturbationnash}] We first study the transport terms and then move to the mollification correction terms. We will adopt a notation similar to the one appearing in the Inductive Lemma \ref{inductive}.

\noindent \textbf{Transport Terms.}  Given the Ansatz \eqref{ansatz}, and the properties in Lemma \ref{chartprop} we can explicitly compute:
\begin{equation*}
        \begin{split}
            (\partial_t+\mathcal{L}_{\tilde v_{\ell,I}})\Theta_I^{p}&=\frac{\tau^a}{\lambda_{q+1}}(\partial_t+\mathcal{L}_{\tilde v_{\ell,I}})\left[ \tilde a_I \alpha_{I}\Psi_I^{1*}(\varphi_{\lambda_{q+1},k_I}\nu_I)\right]\\
            &=\frac{1}{\lambda_{q+1}}\left[ \tilde a_I g_{I}\tilde \Psi_I^{1*}(\varphi_{\lambda_{q+1},k_I}\nu_I)+\tau^a \alpha_I(\partial_t+\tilde v_{\ell,I}\cn)( \tilde a_I)\tilde \Psi_I^{1*}(\varphi_{\lambda_{q+1},k_I}\nu_I)\right]\\
        \end{split}
\end{equation*}
where we used $\alpha_I'=g_I$, see Lemma \ref{timeprofiles}. In particular, in the notation of Lemma \ref{inductive}, we see that: 
$$(\partial_t+\mathcal{L}_{\tilde v_{\ell,I}})\Theta_I^{p}=F_0=a\Psi^{1*}(\phi\nu)$$
where 
$$a=\frac{1}{\lambda_{q+1}}\left[ \tilde a_I g_{I}+\tau^a\alpha_I(\partial_t+\tilde v_{\ell,I}\cn)( \tilde a_I)\right] \ \text{ and } \ \Psi^{1*}(\phi\nu)=\Psi_I^{1*}(\varphi_{\lambda_{q+1},k_I}\nu_I),$$ 
the transport properties and estimates on the chart are guaranteed by Lemma \ref{chartprop}, which in combination with Lemma \ref{slowcoeffestimates} also gives: 
\begin{equation*}
    \begin{split}
        \lambda_{q+1}||\partial_t^j a||_r&\lesssim \lambda_q^{r+j+[r+j-\underline{r}]^+(b-1)\gamma_\ell}\delta_{q+1}^{1/2},
        \\
        \lambda_{q+1}||\mathcal{\tilde A}^\pm_{\ell,I} a||_r&\lesssim\sup_t|g_{I}|\ ||\mathcal{\tilde A}^\pm_{\ell,I}\tilde a_I||_r+1/\tau^a\sup_t|g_{I}'| \ || \tilde a_I||_r+\frac{\tau^a}{2}\sup_t| \alpha_I| \ ||\mathcal{\tilde A}^\pm_{\ell,I}(\mathcal{\tilde A}^+_{\ell,I}+\mathcal{\tilde A}^-_{\ell,I})( \tilde a_I)||_r\\
        &+\frac{1}{2}|\alpha_I'| \ ||(\mathcal{\tilde A}^+_{\ell,I}+\mathcal{\tilde A}^-_{\ell,I})( \tilde a_I)||_r\\
        &\lesssim\lambda_q^{r+[r-\underline{r}]^+(b-1)\gamma_\ell}(1/\tau^c)\delta_{q+1}^{1/2}+(1/\tau^a)\lambda_q^{r+[r-\underline{r}]^+(b-1)\gamma_\ell}\delta_{q+1}^{1/2}+\tau^a/(\tau^c\ell_t)\lambda_q^{r+[r-\underline{r}]^+(b-1)\gamma_\ell}\delta_{q+1}^{1/2}\\
        &\lesssim \left[\frac{1}{\tau^a}+\frac{\tau^a}{\tau^c \ell_t}\right]\lambda_q^{r+[r-\underline{r}]^+(b-1)\gamma_\ell}\delta_{q+1}^{1/2},
    \end{split}
\end{equation*}
we now apply Lemma \ref{inductive} with: 
\begin{equation*}
    \begin{split}
        &z^\pm=\tilde z^\pm_{\ell,I},\\ 
        &\sigma_i=\sigma=\xi_{I},\\
        &\bar\varsigma_i=\bar\varsigma=\tau^a\delta_{q+1}^{1/2}, \ \bar\varsigma_{i,\mathcal{A}}=\bar\varsigma_{\mathcal{A}}=\delta_{q+1}^{1/2},\\
        &L^1(r)=L^\varsigma(r)=L^1_{\mathcal{A}}(r)=\lambda_q^{[r-\underline{r}]^+(b-1)\gamma_\ell},\\
        & A_1=\delta_{q+1}^{1/2},\ A_{1,\mathcal{A}}=\frac{1}{\lambda_{q+1}}\left[\frac{1}{\tau^a}+\frac{\tau^a}{\tau^c \ell_t}\right]\delta_{q+1}^{1/2},\\
        &A_3=A_{3,\mathcal{A}}=0.
    \end{split}
\end{equation*}
We are allowed to choose $\bar N$ in the statement of the Lemma as large as we like, as all the terms are mollified, and we can set $\bar C$ to be the largest of the implicit constants in the estimates above and in Lemma \ref{slowcoeffestimates}, for $0\leq r \leq \bar N$ multiplied by the number of terms in the definition of $a$, that is, two. The Lemma gives decompositions:
\begin{equation*}
    \begin{split}
        \mathcal{L}_{\xi^{p}_I}^k(\partial_t+\mathcal{L}_{\tilde v_{\ell,I}})\Theta_I^{p}&=F_k\simeq a_{k}\tilde\Psi^{1*}_I(\phi_{k}\nu_I),
        \\
        (\partial_t+\mathcal{L}_{\tilde z_{\ell,I}^\pm})\mathcal{L}_{\xi^{p}_I}^k(\partial_t+\mathcal{L}_{\tilde v_{\ell,I}})\Theta_I^{p}&=(\partial_t+\mathcal{L}_{\tilde z_{\ell,I}^\pm}) F_k\simeq a_{k,\mathcal{A}}\tilde\Psi^{1*}_I(\phi_{k}\nu_I)\\
    \end{split}
\end{equation*}
where abusing notation, we do not rename the fast coefficients, together with the estimates:
\begin{equation*}
    \begin{split}
        ||\partial_t^ja_k||_r&\leq C'(C)^k \lambda_{q}^{r+k+j}L(r+k+j)A_1\prod_{i=1}^k\bar \varsigma_i= C'(C)^k\frac{1}{\lambda_{q+1}}\lambda_q^{r+k+j+[r+k+j-\underline{r}]^+(b-1)\gamma_\ell}(\tau^a)^k\delta_{q+1}^{(k+1)/2}\\
        &=C'(C)^k\frac{1}{\lambda_{q+1}}\lambda_q^{r+j+[r+k+j-\underline{r}]^+(b-1)\gamma_\ell}\mathcal{T}_p^k\delta_{q+1}^{1/2},\\
    \end{split}
\end{equation*}
for $0\leq r+k+j\leq \bar N$ and $j=0,1,2$, moreover
\begin{equation*}
    \begin{split}
        ||a_{k,\mathcal{A}}||_r&\leq C'(C^k) \lambda_{q}^{r+k}\left[L_{\mathcal{A}}(k+r)A_{1,\mathcal{A}}\prod_{i=1}^k\bar\varsigma_i +L(k+r)A_1\max_{j\in\{1,\dots,k\}}\bar\varsigma_{j,\mathcal{A}}\prod_{i=1,\dots,k \ i\neq j}\bar\varsigma_i\right]\\
        &\lesssim (C)^k\frac{1}{\lambda_{q+1}}\left[\frac{1}{\tau^a}+\frac{\tau^a}{\tau^c \ell_t}\right]\lambda_q^{r+k+[r+k-\underline{r}]^+(b-1)\gamma_\ell}(\tau^a)^k\delta_{q+1}^{(k+1)/2}\\
        &+(C)^k\frac{1}{\lambda_{q+1}}\lambda_q^{r+k+[r+k-\underline{r}]^+(b-1)\gamma_\ell}(\tau^a)^{k-1}\delta_{q+1}^{(k+1)/2}\\
        &\lesssim (C)^k\frac{1}{\lambda_{q+1}}\left[\frac{1}{\tau^a}+\frac{\tau^a}{\tau^c \ell_t}\right]\lambda_q^{r+[r+k-\underline{r}]^+(b-1)\gamma_\ell}\mathcal{T}_p^k\delta_{q+1}^{1/2},\\
    \end{split}
\end{equation*}
for $0\leq r+k\leq \bar N-1$ and finally:
\begin{equation*}
    ||\partial_t^j\Psi^{1*}_I(\phi_{k}\nu_I)||_r \lesssim \lambda_{q+1}^{r+j} \ \text{ for } \ r\geq 0 \text{ and } j=0,1,2.
\end{equation*}
We conclude that:
\begin{equation*}
    \begin{split}
        ||\mathcal{L}_{\xi^{p}_I}^k(\partial_t+\mathcal{L}_{\tilde v_{\ell,I}})\Theta_I^{p}||_r&\lesssim ||a_{k}||_r||\tilde\Psi^{1*}_I(\phi_{k}\nu_I)||_0+||a_{k}||_0||\tilde\Psi^{1*}_I(\phi_{k}\nu_I)||_r\\
        &\lesssim (C)^k\lambda_q^{r+[r+k-\underline{r}]^+(b-1)\gamma_\ell}\frac{1}{\lambda_{q+1}}\mathcal{T}_p^k\delta_{q+1}^{1/2}+(C^k)\lambda_{q+1}^r\lambda_q^{[k-\underline{r}]^+(b-1)\gamma_\ell}\frac{1}{\lambda_{q+1}}\mathcal{T}_p^k\delta_{q+1}^{1/2}\\
        &\leq (C)^k\lambda_{q+1}^r\underbrace{\left(\frac{\lambda_q}{\lambda_{q+1}}\right)^{r(1-\gamma_\ell)}}_{\leq 1}\lambda_q^{[k-\underline{r}]^+(b-1)\gamma_\ell}\frac{1}{\lambda_{q+1}}\mathcal{T}_p^k\delta_{q+1}^{1/2}\\
        &+(C)^k\lambda_{q+1}^r\lambda_q^{[k-\underline{r}]^+(b-1)\gamma_\ell}\frac{1}{\lambda_{q+1}}\mathcal{T}_p^k\delta_{q+1}^{1/2}\\
        &\lesssim (C)^k\lambda_{q+1}^{r-1}\lambda_q^{[k-\underline{r}]^+(b-1)\gamma_\ell}\mathcal{T}_p^k\delta_{q+1}^{1/2},\\
    \end{split}
\end{equation*}
for $0\leq r+k\leq \bar N$. Similarly, applying $\partial_t^j$ for $j=1,2$ and using the product rule, together with the estimates above, one can show:
$$||\partial_t^j\mathcal{L}_{\xi^{p}_I}^k(\partial_t+\mathcal{L}_{\tilde v_{\ell,I}})\Theta_I^{p}||_r\lesssim (C^k)\lambda_{q+1}^{r+j-1}\lambda_q^{[k-\underline{r}]^+(b-1)\gamma_\ell}\mathcal{T}_p^k\delta_{q+1}^{1/2}$$
for $0\leq r+k+j\leq \bar N$, moreover
\begin{equation*}
    \begin{split}
        ||\mathcal{\tilde A}^\pm_{\ell,I}\mathcal{L}_{\xi^{p}_I}^k(\partial_t+\mathcal{L}_{\tilde v_{\ell,I}})\Theta_I^{p}||_r
        &\leq||(\partial_t+\mathcal{L}_{\tilde z_{\ell,I}^\pm})\mathcal{L}_{\xi^{p}_I}^k(\partial_t+\mathcal{L}_{\tilde v_{\ell,I}})\Theta_I^{p}||_r+||\mathcal{L}_{\xi^{p}_I}^k(\partial_t+\mathcal{L}_{\tilde v_{\ell,I}})\Theta_I^{p}\cn \tilde z^\pm_{\ell,I}||_r\\
        &\lesssim ||a_{k,\mathcal{A}}||_r||\tilde\Psi^{1*}_I(\phi_{k}\nu_I)||_0+||a_{k,\mathcal{A}}||_0||\tilde\Psi^{1*}_I(\phi_{k}\nu_I)||_r\\
        &+||\mathcal{L}_{\xi^{p}_I}^k(\partial_t+\mathcal{L}_{\tilde v_{\ell,I}})\Theta_I^{p}||_r||\tilde z^\pm_{\ell,I}||_1+||\mathcal{L}_{\xi^{p}_I}^k(\partial_t+\mathcal{L}_{\tilde v_{\ell,I}})\Theta_I^{p}||_0||\tilde z^\pm_{\ell,I}||_{r+1}\\
        &\lesssim(C)^k\frac{1}{\lambda_{q+1}}\left[\frac{1}{\tau^a}+\frac{\tau^a}{\tau^c \ell_t}\right]\lambda_q^{r+[r+k-\underline{r}]^+(b-1)\gamma_\ell}\mathcal{T}_p^k\delta_{q+1}^{1/2}\\
        &+(C)^k\left[\frac{1}{\tau^a}+\frac{\tau^a}{\tau^c \ell_t}\right]\lambda_{q+1}^{r-1}\lambda_q^{[k-\underline{r}]^+(b-1)\gamma_\ell}\mathcal{T}_p^k\delta_{q+1}^{1/2}\\
        &+(C)^k\lambda_{q+1}^{r-1}\lambda_q^{[k-\underline{r}]^+(b-1)\gamma_\ell}\mathcal{T}_p^k(1/\tau^c)\delta_{q+1}^{1/2}\\
        &\lesssim(C)^k\lambda_{q+1}^{r-1}\lambda_q^{[k-\underline{r}]^+(b-1)\gamma_\ell}\mathcal{T}_p^k\left[\frac{1}{\tau^a}+\frac{\tau^a}{\tau^c \ell_t}\right]\delta_{q+1}^{1/2},
    \end{split}
\end{equation*}
for $0\leq r+k\leq \bar N-1$. The implicit constants and $C$ in the above statements depend on $\bar N$ and all the other parameters, but not on $a$ and are uniform in $I, \ r, \ j, \ k,$ satisfying the above constraints.

\noindent Proceeding exactly as above and given the extra $\tau^a/\tau^c$ smallness in $\mathcal{L}_{\tilde B_{\ell,I}}\Theta_I^{p}$ one can show that:
\begin{equation*}
    \begin{split}
         ||\partial_t^j\mathcal{L}_{\xi^{p}_I}^k\mathcal{L}_{\tilde B_{\ell,I}}\Theta_I^{p}||_r&\lesssim (C)^k\lambda_{q+1}^{r+j-1}\lambda_q^{[k-\underline{r}]^+(b-1)\gamma_\ell}\mathcal{T}_p^k(\tau^a/\tau^c)\delta_{q+1}^{1/2},
         \\
         ||\mathcal{\tilde A}^\pm_{\ell,I}\mathcal{L}_{\xi^{p}_I}^k\mathcal{L}_{\tilde B_{\ell,I}}\Theta_I^{p}||_r&\lesssim (C)^k\lambda_{q+1}^{r-1}\lambda_q^{[k-\underline{r}]^+(b-1)\gamma_\ell}\mathcal{T}_p^k\left[\frac{1}{\tau^c}+\frac{\tau^a}{\tau^c \ell_t}\right]\delta_{q+1}^{1/2}
    \end{split}
\end{equation*}
where the same ranges of parameters and remarks about the constants apply.

\noindent \textbf{Mollification correction terms.} The proof of these bounds is contained in that of Lemma \ref{estimatesperturbationnash2}, with the following differences:
\begin{itemize}
    \item Here we can set $\bar N$ as large as we like, while in \ref{estimatesperturbationnash2} we pick the largest possible number of derivatives given our inductive assumptions, namely $\bar N=N-m_0-1$.
    \item In \ref{estimatesperturbationnash2} the loss function for the Alfv\'en transport is more complex due to the presence of the non-mollified terms, here we can set 
    $$L_{\mathcal{A}}(r)=\lambda_q^{[r-\underline{r}]^+(b-1)\gamma_\ell},$$ 
    which can be handled as above.
\end{itemize}
\end{proof}

\begin{lemma}[Remainders - Nash Stage]\label{remaindersmoll} For $k_0\geq 0$ let: 
\begin{equation*}
    \begin{split}
        \theta^{k_0}_w&=\sum_I\frac{(-1)^{k_0+1}}{(k_0+1)!}\int_0^{1}(X_{s}^p)_*\left[\mathcal{L}_{\xi^p_I}^{k_0+1}(\partial_t+\mathcal{L}_{\tilde v_{\ell,I}})\Theta^p_I\right](1-s)^{k_0+1}\dd s,\\
        \theta^{k_0}_b&=\sum_I\frac{(-1)^{k_0+1}}{(k_0+1)!}\int_0^{1}(X_{s}^p)_*\left[\mathcal{L}_{\xi^p_I}^{k_0+1}\mathcal{L}_{ \tilde B_{\ell,I}}\Theta^p_I\right](1-s)^{k_0+1}\dd s,\\
        \theta^{k_0}_c&=
        \begin{cases}
            \sum_I\frac{(-1)^{k_0+1}}{(k_0+1)!}\int_0^{1}(X_{s}^p)_*\left[\mathcal{L}_{\xi^p_I}^{k_0+1}\mathcal{L}_{\tilde B_\ell-\tilde B_{\ell,I}}\Theta_I^p\right](1-s)^{k_0+1}\dd s, \\
            \sum_I\frac{(-1)^{k_0+1}}{(k_0+1)!}\int_0^{1}(X_{s}^p)_*\left[\mathcal{L}_{\xi^p_I}^{k_0+1}\mathcal{L}_{\tilde v_\ell-\tilde v_{\ell,I}}\Theta_I^p\right](1-s)^{k_0+1}\dd s.\\
        \end{cases}
    \end{split}
\end{equation*}
Fix $\bar N$ a non-negative integer and $\underline{r}=M-2m_0-k_0^g-8$. There exist constants $C,C'$, such that, for $r\geq 0$ integer and $j=0,1,2$ satisfying $0\leq r +(k_0+1)+j\leq \bar N$, we have:
\begin{equation*}
    \begin{split}
        &||\partial_t^j\theta^{k_0}_w||_{r} \leq C'\frac{(C)^{k_0+1}}{(k_0+2)!}\lambda_{q+1}^{r+j}\mathcal{M}_p^{r+j+1}\left[\frac{1}{\lambda_{q+1}}\lambda_q^{[k_0+1-\underline{r}]^+(b-1)\gamma_\ell}\mathcal{T}_p^{k_0+1}\delta_{q+1}^{1/2}\right],\\
        &||\partial_t^j\theta^{k_0}_b||_{r}\leq C'\frac{(C)^{k_0+1}}{(k_0+2)!}\lambda_{q+1}^{r+j}\mathcal{M}_p^{r+j+1}\left[\frac{1}{\lambda_{q+1}}\lambda_q^{[k_0+1-\underline{r}]^+(b-1)\gamma_\ell}\mathcal{T}_p^{k_0+1}(\tau^a/\tau^c)\delta_{q+1}^{1/2}\right],\\
        &||\partial_t^j\theta^{k_0}_c||_{r}\leq C'\frac{(C)^{k_0+1}}{(k_0+2)!}\lambda_{q+1}^{r+j}\mathcal{M}_p^{r+j+1}\left[\lambda_q^{[k_0+1-\underline{r}]^+(b-1)\gamma_\ell}\mathcal{T}_p^{k_0+1}(\ell\lambda_q)^{m_0}\tau^a\delta_q^{1/2}\delta_{q+1}^{1/2}\right].\\
    \end{split}
\end{equation*}
The constants $C, \ C'$ depend on all the parameters but not on $a$ and are uniform in $r, \ j, \ k_0$ satisfying the constraint above.
\end{lemma}
\begin{remark}From Lemma \ref{remaindersmoll} with $k_0=k_0^p, \ \bar N=N$, we deduce that:
\begin{equation*}
    \begin{split}
        ||\partial_t^j\theta^p||_{r} &\leq C'\frac{(C)^{k_0^p+1}}{(k_0^p+2)!}\lambda_{q+1}^{r+j}\mathcal{M}_p^{r+j+1}\left[\frac{1}{\lambda_{q+1}}\lambda_q^{[k_0^p+1-\underline{r}]^+(b-1)\gamma_\ell}\mathcal{T}_p^{k_0^p+1}\delta_{q+1}^{1/2}\right],\\
        ||\partial_t^j\theta^p_b||_{r}&\leq C'\frac{(C)^{k_0^p+1}}{(k_0^p+2)!}\lambda_{q+1}^{r+j}\mathcal{M}_p^{r+j+1}\left[\frac{1}{\lambda_{q+1}}\lambda_q^{[k_0^p+1-\underline{r}]^+(b-1)\gamma_\ell}\mathcal{T}_p^{k_0^p+1}(\tau^a/\tau^c)\delta_{q+1}^{1/2}\right]
    \end{split}
\end{equation*}
for $0\leq r +(k_0^p+1)+j\leq N$, where $\theta_w^p, \ \theta_b^p$ are as in \eqref{decompnashw} and \eqref{decompnashb}.
\end{remark}

\begin{proof}[Proof of Lemma \ref{remaindersmoll}]According to the bounds in Lemmas \ref{estimatesflowp}, \ref{estimatesperturbationnash} and the composition estimates in Proposition \ref{compestimates}, we have:
\begin{equation*}
    \begin{split}
        &||(X_s^p)_*[\mathcal{L}_{\xi^{p}_I}^{k_0+1}(\partial_t+\mathcal{L}_{\tilde v_{\ell,I}})\Theta^{p}_I]||_r\\
        &=||\DD X_s^p((X_{s}^p)^{-1})[\mathcal{L}_{\xi^{p}_I}^{k_0+1}(\partial_t+\mathcal{L}_{\tilde v_{\ell,I}})\Theta^{p}_I]((X_{s}^p)^{-1})||_r\\
        &\overset{r=0}{\leq}||\DD X_s^p||_0||\mathcal{L}_{\xi^{p}_I}^{k_0+1}(\partial_t+\mathcal{L}_{\tilde v_{\ell,I}})\Theta^{p}_I]||_0\\
        &\lesssim\frac{1}{\lambda_{q+1}}\mathcal{M}_p\lambda_q^{[k_0+1-\underline{r}]^+(b-1)\gamma_\ell}\mathcal{T}^{k_0+1}_p\delta_{q+1}^{1/2}\\
        &\overset{r\geq 1}{\lesssim} ||\DD X_s^p((X_{s}^p)^{-1})||_r||\mathcal{L}_{\xi^{p}_I}^{k_0+1}(\partial_t+\mathcal{L}_{\tilde v_{\ell,I}})\Theta^{p}_I||_0+||\DD X_s^p((X_{s}^p)^{-1})||_0||[\mathcal{L}_{\xi^{p}_I}^{k_0+1}(\partial_t+\mathcal{L}_{\tilde v_{\ell,I}})\Theta^{p}_I]((X_{s}^p)^{-1})||_r\\
        &\lesssim ||\mathcal{L}_{\xi^{p}_I}^{k_0+1}(\partial_t+\mathcal{L}_{\tilde v_{\ell,I}})\Theta^{p}_I||_0\left[||\DD X_s^p||_1||X_{-s}^p||_{r-1}+||\DD X_s^p||_r||\DD X_{-s}^p||_{0}^r\right]\\
        &+||\DD X_s^p||_0\left[||\mathcal{L}_{\xi^{p}_I}^{k_0+1}(\partial_t+\mathcal{L}_{\tilde v_{\ell,I}})\Theta^{p}_I||_1||\DD X_{-s}^p||_r+||[\mathcal{L}_{\xi^{p}_I}^{k_0+1}(\partial_t+\mathcal{L}_{\tilde v_{\ell,I}})\Theta^{p}_I]||_r||\DD X_{-s}^p||_0^r\right]\\
        &\lesssim (C)^{k_0+1}\lambda_{q+1}^{r-1}\mathcal{M}_p^{r+1}\lambda_q^{[k_0+1-\underline{r}]^+(b-1)\gamma_\ell}\mathcal{T}^{k_0+1}_p\delta_{q+1}^{1/2},
    \end{split}
\end{equation*}
for $0\leq s\leq 1$, $0\leq r+k_0\leq \bar N$. Recall that even if everything here is mollified, Lemma \ref{estimatesflowp} guarantees good bounds only up to an arbitrarily large, but fixed, $\bar N$. 

\noindent By means of the identities \eqref{firsttimederivate} and \eqref{secondtimederivative} and estimating as above with the bounds from Lemma \ref{estimatesflowp} and \ref{estimatesperturbationnash} one can show:
$$||\partial_t^j[(X_s^p)_*[\mathcal{L}_{\xi^{p}_I}^{k_0+1}(\partial_t+\mathcal{L}_{\tilde v_{\ell,I}})\Theta^{p}_I]]||_r\lesssim (C)^{k_0+1}\lambda_{q+1}^{r+j-1}\mathcal{M}_p^{r+j+1}\lambda_q^{[k_0+1-\underline{r}]^+(b-1)\gamma_\ell}\mathcal{T}^{k_0+1}_p\delta_{q+1}^{1/2},$$
for $0\leq s\leq 1$, $0\leq r+k_0+j\leq \bar N$. The same remarks about the constants apply.

\noindent The same argument given the estimates for $\mathcal{L}_{\xi^{p}_I}^{k_0+1}\mathcal{L}_{\tilde B_{\ell,I}}\Theta^{p}_I, \ \mathcal{L}_{\xi^{p}_I}^{k_0+1}\mathcal{L}_{\tilde B_\ell-\tilde B_{\ell,I}}\Theta^{p}_I, \mathcal{L}_{\xi^{p}_I}^{k_0+1}\mathcal{L}_{\tilde v_\ell-\tilde v_{\ell,I}}\Theta^{p}_I\ $ from Lemma \ref{estimatesperturbationnash}, shows:
\begin{equation*}
    \begin{split}
        &||\partial_t^j[(X_s^p)_*[\mathcal{L}_{\xi^{p}_I}^{k_0+1}\mathcal{L}_{\tilde B_{\ell,I}}\Theta^{p}_I]]||_r\lesssim (C)^{k_0+1}\lambda_{q+1}^{r+j-1}\mathcal{M}_p^{r+j+1}\lambda_q^{[k_0+1-\underline{r}]^+(b-1)\gamma_\ell}\mathcal{T}^{k_0+1}_p(\tau^a/\tau^c)\delta_{q+1}^{1/2},\\
        &||\partial_t^j[(X_s^p)_*[\mathcal{L}_{\xi^{p}_I}^{k_0+1}\mathcal{L}_{\tilde v_\ell-\tilde v_{\ell,I}}\Theta^{p}_I]]||_r\lesssim (C)^{k_0+1} \lambda_{q+1}^{r+j}\mathcal{M}_p^{r+j+1}\lambda_q^{[k_0+1-\underline{r}]^+(b-1)\gamma_\ell}\mathcal{T}_p^{k_0+1}(\ell\lambda_q)^{m_0}\tau^a\delta_q^{1/2}\delta_{q+1}^{1/2},\\
        &||\partial_t^j[(X_s^p)_*[\mathcal{L}_{\xi^{p}_I}^{k_0+1}\mathcal{L}_{\tilde B_\ell-\tilde B_{\ell,I}}\Theta^{p}_I]]||_r\lesssim (C)^{k_0+1} \lambda_{q+1}^{r+j}\mathcal{M}_p^{r+j+1}\lambda_q^{[k_0+1-\underline{r}]^+(b-1)\gamma_\ell}\mathcal{T}_p^{k_0+1}(\ell\lambda_q)^{m_0}\tau^a\delta_q^{1/2}\delta_{q+1}^{1/2}.
    \end{split}
\end{equation*}
The same range of indexes and the remarks about the constants apply.

\noindent Since by construction at most one index $\Theta_I^p$ is non-zero at each space-time point, see Lemma \ref{disjoint}, we conclude:
\begin{equation*}
    \begin{split}
        ||\partial_t^j\theta^{k_0}_w||_r&\lesssim \frac{1}{(k_0+1)!}\left[\int_0^1(1-s)^{k_0+1}\dd s\right]\sup_{I\in\mathcal{I}}\max_{0\leq s\leq 1}||\partial_t^j[(X_s^p)_*[\mathcal{L}_{\xi^{p}_I}^{k_0+1}(\partial_t+\mathcal{L}_{\tilde v_{\ell,I}})\Theta^{p}_I]]||_r\\
        &\lesssim\frac{(C)^{k_0+1}}{(k_0+2)!}\lambda_{q+1}^{r+j}\mathcal{M}_p^{r+j+1}\left[\frac{1}{\lambda_{q+1}}\mathcal{T}_p^{k_0+1}\lambda_q^{[k_0+1-\underline{r}]^+(b-1)\gamma_\ell}\delta_{q+1}^{1/2}\right],\\
    \end{split}
\end{equation*}
for $0\leq r +k_0+j\leq \bar N$, the implicit constant and $C$ depends on $\bar N, \ r$ and all the other parameters, but not on $a$ and is uniform in $k_0, \ j, \ r$ satisfying the above constraint. 

\noindent Similarly, one can show:
$$||\partial_t^j\theta^{k_0}_b||_r\lesssim \frac{(C)^{k_0+1}}{(k_0+2)!}\lambda_{q+1}^r\mathcal{M}_p^{r+j+1}\left[\frac{1}{\lambda_{q+1}}\lambda_q^{[k_0+1-\underline{r}]^+(b-1)\gamma_\ell}\mathcal{T}_p^{k_0+1}(\tau^a/\tau^c)\delta_{q+1}^{1/2}\right]$$
and
$$||\partial_t^j\theta^{k_0}_c||_r \lesssim \frac{(C)^{k_0+1}}{(k_0+2)!}\lambda_{q+1}^r\mathcal{M}_p^{r+j+1}\left[\lambda_q^{[k_0+1-\underline{r}]^+(b-1)\gamma_\ell}\mathcal{T}_p^{k_0+1}(\ell\lambda_q)^{m_0}\tau^a\delta_q^{1/2}\delta_{q+1}^{1/2}\right].$$
The same range of indexes and the remarks about the constants apply.
\end{proof}


\begin{lemma}[Lie-Derivative and Transport Estimates - Non Mollified Terms]\label{estimatesperturbationnash2} Let $I\in \mathcal{I}$ and $r, \ j, \ k\geq 0$ integers and $L_p, \ L_{p,\mathcal{A}}$ the admissible loss functions in \eqref{lossparameters}. There exist constants $C, \ C'$ which depend on all the parameters but not on $a$ and are uniform in $r, \ j, \ k, \ I$, such that the following estimates hold.

\noindent\textbf{Lie-derivatives bounds.} For $j=0,1,2$ and $0\leq r+k+j\leq N-m_0-1$, we have:
\begin{equation*}
        ||\partial_t^j\mathcal{L}_{\xi^{p}_I}^k\mathcal{L}_{\tilde v_q-\tilde v_\ell}\Theta_I^{p}||_{r}, \ ||\partial_t^j\mathcal{L}_{\xi^{p}_I}^k\mathcal{L}_{\tilde B_q-\tilde B_\ell}\Theta_I^{p}||_{r}\leq C'(C)^k\lambda_{q+1}^{r+j}L_p(k)\mathcal{T}_p^k(\ell\lambda_q)^{m_0}\tau^a\delta_q^{1/2}\delta_{q+1}^{1/2}. 
\end{equation*}

\noindent \textbf{Alfv\'en-transport bounds.} For $0\leq r+k\leq N-m_0-2$, we have:
\begin{equation*}
    ||\mathcal{\tilde A}_{\ell,I}^\pm\mathcal{L}_{\xi^{p}_I}^k\mathcal{L}_{\tilde v_q-\tilde v_\ell}\Theta_I^{p}||_r, \ ||\mathcal{\tilde A}_{\ell,I}^\pm\mathcal{L}_{\xi^{p}_I}^k\mathcal{L}_{\tilde B_q-\tilde B_\ell}\Theta_I^{p}||_r\leq C'(C)^k\lambda_{q+1}^{r}\left[L_p(k)+\tau^a/\tau^cL_{p,\mathcal{A}}(k)\right]\mathcal{T}_p^k(\ell\lambda_q)^{m_0}\delta_q^{1/2}\delta_{q+1}^{1/2}.
\end{equation*}

\noindent \textbf{remainder bounds.} For $j=0,1,2$ and $0\leq r+j\leq N-(k_0^p+1)-m_0-1$, we have:
\begin{equation*}
    ||\partial_t^j\mathring \theta^p_w||_{r}, \ ||\partial_t^j\mathring \theta^p_b||_{r}\leq C'\frac{(C)^{k_0^p+1}}{(k_0^p+2)!} \lambda_{q+1}^{r+j}\mathcal{M}_p^{r+j+1}\mathcal{T}_p^{k_0^p+1}L_p(k_0^p+1)(\ell\lambda_q)^{m_0}\tau^a\delta_q^{1/2}\delta_{q+1}^{1/2}.\\
\end{equation*}

\end{lemma}

\begin{proof}[Proof of Lemma \ref{estimatesperturbationnash2}] We proceed as in the proof of Lemma \ref{estimatesperturbationnash}; here, extra care is needed because terms with non-controlled geometry appear, and for high derivatives, the transport estimates degrade. We will adopt a notation similar to the one appearing in the Inductive Lemma \ref{inductive}.

\noindent\textbf{Estimates on the top terms.} We only show the bounds for the terms involving the velocity field; the estimates for the magnetic field follow similarly. From the definition of $\Theta_I^p$ in \eqref{ansatz}, it follows that:
\begin{equation*}
    \begin{split}
        \mathcal{L}_{\tilde v_q-\tilde v_\ell}\Theta_I^{p}&=\frac{\tau^a}{\lambda_{q+1}}(\tilde v_q-\tilde v_\ell)\cn\left[ \tilde a_I \alpha_{I}\Psi_I^{1*}(\varphi_{\lambda_{q+1},k_I}\nu_I)\right]+\frac{\tau^a}{\lambda_{q+1}}\DD(\tilde v_q-\tilde v_\ell)^\top \left[ \tilde a_I \alpha_{I}\Psi_I^{1*}(\varphi_{\lambda_{q+1},k_I}\nu_I)\right]\\
        &=\underbrace{\left[\tau^a \tilde a_I \alpha_{I}(\tilde v_q-\tilde v_\ell)\cdot(\Psi_I^{1*}k_I)\right]}_{a_1}\Psi_I^{1*}(\underbrace{\varphi_{\lambda_{q+1},k_I}'}_{\phi_1}\nu_I)\\
        &+\Psi_I^{*}(\underbrace{\varphi_{\lambda_{q+1},k_I}}_{\phi_3})\underbrace{\frac{\tau^a}{\lambda_{q+1}}\left[(\tilde v_q-\tilde v_\ell)\cn\left[ \tilde a_I \alpha_{I}\Psi_I^{1*}\nu_I\right]+\DD(\tilde v_q-\tilde v_\ell)^\top \left[ \tilde a_I \alpha_{I}\Psi_I^{1*}(\nu_I)\right]\right]}_{Y}\\
            &= a_1\Psi_I^{1*}(\phi_1\nu_I)+\Psi_I^{*}(\phi_3)Y=F_0
    \end{split}
\end{equation*}
Lemma \ref{chartprop} guarantees the transport properties and estimates on the chart, which, in combination with the bounds in Lemmas \ref{slowcoeffestimates}, \ref{standardmollnash}, and \ref{stabilitynash}, give: 
\begin{equation*}
    \begin{split}
        ||\partial_t^j a_1||_r
        &\lesssim \tau^a\lambda_q^{r+j}L_p(r+j)(\ell\lambda_q)^{m_0}\delta_q^{1/2}\delta_{q+1}^{1/2}
        \\
        ||\mathcal{\tilde A}^\pm_{\ell,I} a_1||_r&\lesssim\underbrace{|| \tilde a_I \alpha'_{I}(\tilde v_q-\tilde v_\ell)\cdot\Psi_I^{1*}k_I||_r}_{\lambda_q^{r}(\ell\lambda_q)^{m_0}L_p(r)\delta_q^{1/2}\delta_{q+1}^{1/2}}+\underbrace{\tau^a||\mathcal{\tilde A}^\pm_{\ell,I}( \tilde a_I) \alpha_{I}(\tilde v_q-\tilde v_\ell)\cdot\Psi_I^{1*}k_I||_r}_{(\tau^a/\tau^c)\lambda_q^{r}(\ell\lambda_q)^{m_0}L_p(r)\delta_q^{1/2}\delta_{q+1}^{1/2}}\\
        &+\underbrace{\tau^a|| \tilde a_I \alpha_{I}[[\mathcal{\tilde A}^\pm+(z^\pm_{\ell,I}-\tilde z^\pm_q)\cn](\tilde v_q-\tilde v_\ell)]\cdot\Psi_I^{1*}k_I||_r}_{\tau^a\lambda_q^{r+1}(\ell\lambda_q)^{m_0}L_{p,\mathcal{A}}(r)\delta_q\delta_{q+1}^{1/2}+\tau^a\lambda_q^{r+1}(\ell\lambda_q)^{2m_0}L_{p}(r)\delta_q\delta_{q+1}^{1/2}}+\underbrace{\tau^a|| \tilde a_I \alpha_{I}(\tilde v_q-\tilde v_\ell)\cdot(\mathcal{\tilde A}^\pm_{\ell,I}\Psi_I^{1*}k_I)||_r}_{\tau^a\lambda_q^{r+1}(\ell\lambda_q)^{m_0}L_p(r)\delta_q\delta_{q+1}^{1/2}}\\
        &\lesssim\lambda_q^r\left[L_p(r)+(\tau^a/\tau^c)L_{p,\mathcal{A}}(r)\right](\ell\lambda_q)^{m_0}\delta_q^{1/2}\delta_{q+1}^{1/2}
    \end{split}
\end{equation*}
for  $0\leq r+j\leq N-m_0$ and $0\leq r+j\leq N-m_0-1$ respectively and we used that according to Lemma \ref{chartprop} we have $$\mathcal{\tilde A}^\pm_{\ell,I}\Psi_I^{1*}k=-(\DD z^\pm_{\ell,I})^\top[\Psi_I^{1*}k]$$ together with $$L_p\leq L_{p,\mathcal{A}} \ \text{ and } L_{p,\mathcal{A}}(r)= L_{g,\mathcal{A}}(r+m_0+k_0^g+2)$$ 
see \eqref{lossyvsmoll} and their definitions in \eqref{lossparameters}. We will use these observations in the following without further mention. To estimate $Y$, we first split it into:
$$Y=\underbrace{\frac{\tau^a}{\lambda_{q+1}}(\tilde v_q-\tilde v_\ell)\cn\left[ \tilde a_I \alpha_{I}\Psi_I^{1*}\nu_I\right]}_{T_1}+\underbrace{\frac{\tau^a}{\lambda_{q+1}}\DD(\tilde v_q-\tilde v_\ell)^\top \left[ \tilde a_I \alpha_{I}\Psi_I^{1*}(\nu_I)\right]}_{T_2}$$  
According to the same Lemmas as above, we have: 
\begin{equation*}
    \begin{split}
        ||\partial_t^jT_1||_r
        &\lesssim \frac{\tau^a}{\lambda_{q+1}}\lambda_q^{r+j+1}L_p(r+j+1)(\ell\lambda_q)^{m_0}\delta_q^{1/2}\delta_{q+1}^{1/2},
        \\
        ||(\partial_t+\mathcal{L}_{\tilde z_{\ell,I}^\pm})T_1||_r&\lesssim\underbrace{\frac{1}{\lambda_{q+1}}||\alpha'_{I}(\tilde v_q-\tilde v_\ell)\cn( \tilde a_I\Psi_I^{1*}\nu_I)||_r}_{\frac{1}{\lambda_{q+1}}\lambda_q^{r+1}(\ell\lambda_q)^{m_0}L_p(r+1)\delta_q^{1/2}\delta_{q+1}^{1/2}}+\underbrace{\frac{\tau^a}{\lambda_{q+1}}||\alpha_{I}[\mathcal{\tilde A}^\pm(\tilde v_q-\tilde v_\ell)]\cn( \tilde a_I \Psi_I^{1*}\nu_I)||_r}_{\frac{\tau^a}{\lambda_{q+1}}\lambda_q^{r+2}(\ell\lambda_q)^{m_0}L_{p,\mathcal{A}}(r)\delta_q\delta_{q+1}^{1/2}}\\
        &+\underbrace{\frac{\tau^a}{\lambda_{q+1}}||\alpha_{I}[(\tilde z^\pm_{\ell,I}-\tilde z^\pm_q)\cn(\tilde v_q-\tilde v_\ell)]\cn( \tilde a_I \Psi_I^{1*}\nu_I)||_r}_{\frac{\tau^a}{\lambda_{q+1}}\lambda_q^{r+2}(\ell\lambda_q)^{2m_0}L_{p}(r+1)\delta_q\delta_{q+1}^{1/2}}\\
        &+\underbrace{\frac{\tau^a}{\lambda_{q+1}}||\alpha_I[(\tilde v_q-\tilde v_\ell)\cn z^\pm_{\ell,I}]\cn( \tilde a_I \Psi_I^{1*}\nu_I)||_r}_{\frac{\tau^a}{\lambda_{q+1}}\lambda_q^{r+2}(\ell\lambda_q)^{m_0}L_p(r+1)\delta_q\delta_{q+1}^{1/2}}+\underbrace{\frac{\tau^a}{\lambda_{q+1}}||\alpha_{I}(\tilde v_q-\tilde v_\ell)\cn\mathcal{\tilde A}^\pm_{\ell,I}( \tilde a_I \Psi_I^{1*}k)||_r}_{\frac{1}{\lambda_{q+1}}(\tau^a/\tau^c)\lambda_q^{r+1}(\ell\lambda_q)^{m_0}L_p(r+1)\delta_q^{1/2}\delta_{q+1}^{1/2}}\\
        &+\underbrace{||T_1\cn \tilde z^\pm_{\ell,I}||_r}_{\frac{\tau^a}{\lambda_{q+1}}\lambda_q^{r+2}L_p(r)(\ell\lambda_q)^{m_0}\delta_q\delta_{q+1}^{1/2}}\\
        &\lesssim\frac{1}{\lambda_{q+1}}\lambda_q^{r+1}\left[L_p(r+1)+\tau^a/\tau^cL_{p,\mathcal{A}}(r)\right](\ell\lambda_q)^{m_0}\delta_q^{1/2}\delta_{q+1}^{1/2}.
    \end{split}
\end{equation*}
Similar estimates hold for $T_2$; we only state the resulting bounds:
\begin{equation*}
    \begin{split}
        &||T_2||_r\lesssim\frac{\tau^a}{\lambda_{q+1}}\lambda_q^{r+j+1}L_p(r+j+1)(\ell\lambda_q)^{m_0}\delta_q^{1/2}\delta_{q+1}^{1/2},\\
        &||(\partial_t+\mathcal{L}_{\tilde z_{\ell,I}^\pm})T_2||_r\lesssim \frac{1}{\lambda_{q+1}}\lambda_q^{r+1}\left[L_p(r+1)+\tau^a/\tau^cL_{p,\mathcal{A}}(r)\right](\ell\lambda_q)^{m_0}\delta_q^{1/2}\delta_{q+1}^{1/2}\\
    \end{split}
\end{equation*}
and we conclude that:
\begin{equation*}
    \begin{split}
        ||\partial_t^jY||_r&\lesssim ||\partial_t^jT_1||_r+||\partial_t^jT_2||_r\\
            &\lesssim\frac{1}{\lambda_{q+1}}\lambda_q^{r+j+1}L_p(r+j+1)(\ell\lambda_q)^{m_0}\tau^a\delta_q^{1/2}\delta_{q+1}^{1/2} \ \text{ for } \ 0\leq r+j\leq N-m_0-1,
            \\
        ||(\partial_t+\mathcal{L}_{\tilde z_{\ell,I}^\pm}) Y||_r&\lesssim||(\partial_t+\mathcal{L}_{\tilde z_{\ell,I}^\pm})T_1||_r+||(\partial_t+\mathcal{L}_{\tilde z_{\ell,I}^\pm})T_2||_r\\
        &\lesssim\frac{1}{\lambda_{q+1}}\lambda_q^{r+1}\left[L_p(r+1)+\tau^a/\tau^cL_{p,\mathcal{A}}(r)\right](\ell\lambda_q)^{m_0}\delta_q^{1/2}\delta_{q+1}^{1/2} \ \text{ for } \ 0\leq r+j\leq N-m_0-2.
    \end{split}
\end{equation*}

\noindent \textbf{Full estimates.} We now follow the same strategy as before and apply Lemma \ref{inductive} with: 
\begin{equation*}
    \begin{split}
        &z^\pm=\tilde z^\pm_{\ell,I},\\ 
        &\sigma_i=\sigma=\xi_{I},\\
        &\bar\varsigma_i=\bar\varsigma=\tau^a\delta_{q+1}^{1/2}, \ \bar\varsigma_{i,\mathcal{A}}=\bar\varsigma_{\mathcal{A}}=\delta_{q+1}^{1/2},\\
        &L(r)=L_p(r+1), \ L_{\mathcal{A}}(r)=L_p(r+1)+\tau^a/\tau^cL_{p,\mathcal{A}}(r),\\
        &L^1(r)=L_p(r), \ L^1_{\mathcal{A}}(r)=L_p(r)+\tau^a/\tau^cL_{p,\mathcal{A}}(r),\\
        &L^\varsigma(r)=L_p(r),\\
        &A_1=(\ell\lambda_q)^{m_0}\tau^a\delta_q^{1/2}\delta_{q+1}^{1/2}, \ A_{1,\mathcal{A}}=(\ell\lambda_q)^{m_0}\delta_q^{1/2}\delta_{q+1}^{1/2},\\
        &A_3=\frac{\lambda_q}{\lambda_{q+1}}A_1, \ A_{3,\mathcal{A}}=\frac{\lambda_q}{\lambda_{q+1}}A_{1,\mathcal{A}}
    \end{split}
\end{equation*}
and let $\bar N$ be as large as possible, namely
$$\bar N =N-m_0-1$$
and we can set $\bar C$ to be the largest of the implicit constants in the estimates above and in Lemma \ref{slowcoeffestimates}, for $0\leq r \leq \bar N$ multiplied by the largest number of terms between the definitions of $\ a_1, \ a_2,\ Y$.

\noindent We conclude that: 
\begin{equation*} 
    \begin{split}
        \mathcal{L}_{\xi^{p}_I}^k\mathcal{L}_{\tilde v_q-\tilde v_\ell}\Theta_I^{p}&=F_k\simeq a_{1,k}\tilde \Psi^{1*}_I(\phi_{1,k}\nu_I)+a_{2,k} \tilde \Psi^{1*}_I(\phi_{2,k}k_I)+\tilde \Psi^{*}_I(\phi_{3,k}) Y_k,
        \\
        (\partial_t+\mathcal{L}_{\tilde z_{\ell,I}^\pm})\mathcal{L}_{\xi^{p}_I}^k\mathcal{L}_{\tilde v_q-\tilde v_\ell}\Theta_I^{p}&= (\partial_t+\mathcal{L}_{\tilde z_{\ell,I}^\pm})F_k\simeq a_{1,k,\mathcal{A}}\tilde \Psi^{1*}_I(\phi_{1,k}\nu_I)+a_{2,k,\mathcal{A}} \tilde \Psi^{1*}_I(\phi_{2,k}k_I)+\tilde \Psi^{*}_I(\phi_{3,k}) Y_{k,\mathcal{A}}
    \end{split}
\end{equation*}
where, abusing notation, we do not rename the fast coefficients. Moreover, for $k\geq 1, \ j=0,1,2$ and $0\leq r+j+k\leq \bar N$, we have:
\begin{equation*}
    \begin{split}
    ||\partial_t^j a_{1,k}||_r&\leq C'(C)^k\lambda_{q}^{r+k+j}L(r+j+k)A_1\prod_{i=1}^k\bar \varsigma_i\\
    &\lesssim\lambda_q^{r+j}L_p(r+j+k)\mathcal{T}_p^k(\ell\lambda_q)^{m_0}\tau^a\delta_q^{1/2}\delta_{q+1}^{1/2},
    \\
        ||\partial_t^j a_{2,k}||_r&\leq C'(C)^k \lambda_{q+1}\lambda_{q}^{r+k+j-1}L(r+j+k-1)A_3\prod_{i=1}^k\bar \varsigma_i\\
        &\lesssim(C)^k\lambda_q^{r+j}L_p(r+j+k)\mathcal{T}_p^k(\ell\lambda_q)^{m_0}\tau^a\delta_q^{1/2}\delta_{q+1}^{1/2},
        \\
        ||\partial_t^jY_k||_r&\leq C'(C)^k  \lambda_{q}^{r+k+j}L(r+j+k)A_3\prod_{i=1}^k\bar \varsigma_i\\
        &\lesssim (C)^k\left(\frac{\lambda_q}{\lambda_{q+1}}\right)\lambda_q^{r+j}L_p(r+j+k+1)\mathcal{T}_p^k(\ell\lambda_q)^{m_0}\delta_q^{1/2}\tau^a\delta_{q+1}^{1/2}\\
        &\leq(C)^k\lambda_q^{r+j}L_p(r+j+k)\mathcal{T}_p^k(\ell\lambda_q)^{m_0}\tau^a\delta_q^{1/2}\delta_{q+1}^{1/2},
    \end{split}
\end{equation*}
together with Alfv\'en transport estimates 
\begin{equation*}
    \begin{split}
        ||a_{1,k,\mathcal{A}}||_r&\leq C'(C)^k\lambda_{q}^{r+k}\left[L^1_{\mathcal{A}}(k+r)A_{1,\mathcal{A}}\prod_{i=1}^k\bar\varsigma_i +L^1(k+r)A_1\max_{j\in\{1,\dots,k\}}\bar\varsigma_{j,\mathcal{A}}\prod_{i=1,\dots,k \ i\neq j}\bar\varsigma_i\right]\\
        &\lesssim (C)^k \lambda_q^{r}\left[L_p(r+k)+\frac{\tau^a}{\tau^c}L_{p,\mathcal{A}}(r+k)\right]\mathcal{T}_p^k(\ell\lambda_q)^{m_0}\delta_q^{1/2}\delta_{q+1}^{1/2},
        \\
        ||a_{2,k,\mathcal{A}}||_r&\leq C'(C)^k \lambda_{q+1}\lambda_{q}^{r+k-1}\left[L_{\mathcal{A}}(r+k-1)A_{3,\mathcal{A}}\prod_{i=1}^k\bar\varsigma_i +L(r+k-1) A_3\max_{j\in\{1,\dots,k\}}\bar\varsigma_{j,\mathcal{A}}\prod_{i=1,\dots,k \ i\neq j}\bar\varsigma_i\right]\\
        &\lesssim (C)^k\lambda_q^{r}\left[L_p(r+k)+\frac{\tau^a}{\tau^c}L_{p,\mathcal{A}}(r+k-1)\right]\mathcal{T}_p^k(\ell\lambda_q)^{m_0}\delta_q^{1/2}\delta_{q+1}^{1/2},
        \\
        ||Y_{k,\mathcal{A}}||_r&\leq C'(C)^k\lambda_{q}^{r+k}\left[L_{\mathcal{A}}(k+r)A_{3,\mathcal{A}}\prod_{i=1}^k\bar\varsigma_i +L(k+r)A_3\max_{j\in\{1,\dots,k\}}\bar\varsigma_{j,\mathcal{A}}\prod_{i=1,\dots,k \ i\neq j}\bar\varsigma_i\right]\\
        &\lesssim (C)^k\lambda_q^{r}\left[L_p(r+k)+\frac{\tau^a}{\tau^c}L_{p,\mathcal{A}}(r+k-1)\right]\mathcal{T}_p^k(\ell\lambda_q)^{m_0}\delta_q^{1/2}\delta_{q+1}^{1/2}
    \end{split}
\end{equation*}
for $k \geq 1, \ 0\leq r+k\leq \bar N-1$ and finally
\begin{equation*}
    ||\partial_t^j\Psi^{1*}_I(\phi_{1,k}\nu_I)||_r, \ ||\partial_t^j\Psi^{1*}_I(\phi_{2,k}k_I)||_r, \ ||\partial_t^j\Psi^{1*}_I(\phi_{3,k})||_r \lesssim \lambda_{q+1}^{r+j}.
\end{equation*}
In the above bounds, $C, \ C'$ and the implicit constants depend on $\bar N$ and thus on all the parameters, but not on $a$ and are uniform in $r,\ k, \ j, \ I$ under the constraints above. Note that to simplify the expressions, we used the fact that $L_p$ and $ L_{p,\mathcal{A}}$ are admissible loss functions, see Definition \ref{admissible}, as shown after \eqref{lossparameters}. We will do this without further mention in the following. 

\noindent The bounds above, and the fact that $a_{1,k}$ has the worst estimates, give:
\begin{equation*}
    \begin{split}
        ||\mathcal{L}_{\xi^{p}_I}^k\mathcal{L}_{\tilde v_q-\tilde v_{\ell}}\Theta_I^{p}||_{r}&\lesssim ||a_{1,k}||_r||\tilde \Psi^{1*}_I(\phi_{1,k}\nu)||_0+||a_{1,k}||_0||\tilde \Psi^{1*}_I(\phi_{1,k}\nu)||_r+||a_{2,k}||_r||\tilde \Psi^{1*}_I(\phi_{2,k}k)||_0+||a_{2,k}||_0||\tilde \Psi^{1*}_I(\phi_{2,k}k)||_r\\
        &+||\tilde \Psi^{*}_I(\phi_{3,k})||_r||Y_k||_0+||\tilde \Psi^{*}_I(\phi_{3,k})||_0||Y_k||_r\\
        &\lesssim (C)^k\lambda_q^{r}L_p(k+r)\mathcal{T}_p^k(\ell\lambda_q)^{m_0}\tau^a\delta_q^{1/2}\delta_{q+1}^{1/2}+(C)^k\lambda_{q+1}^{r}L_p(k)\mathcal{T}_p^k(\ell\lambda_q)^{m_0}\tau^a\delta_q^{1/2}\delta_{q+1}^{1/2}\\
        &=(C)^k \lambda_{q+1}^{r}\underbrace{\left[\left(\frac{\lambda_q}{\lambda_{q+1}}\right)^rL_p(k+r)\right]}_{\leq L_p(k)}\mathcal{T}_p^k(\ell\lambda_q)^{m_0}\tau^a\delta_q^{1/2}\delta_{q+1}^{1/2}+(C)^k\lambda_{q+1}^{r}L_p(k)\mathcal{T}_p^k(\ell\lambda_q)^{m_0}\tau^a\delta_q^{1/2}\delta_{q+1}^{1/2}\\
        &\lesssim (C)^k\lambda_{q+1}^{r}L_p(k)\mathcal{T}_p^k(\ell\lambda_q)^{m_0}\tau^a\delta_q^{1/2}\delta_{q+1}^{1/2}
    \end{split}
\end{equation*}
for $0\leq r+k\leq \bar N$, now distributing the pure time derivatives between fast and slow coefficients and estimating in the same way, we also deduce:
$$||\partial_t^j\mathcal{L}_{\xi^{p}_I}^k\mathcal{L}_{\tilde v_q-\tilde v_{\ell}}\Theta_I^{p}||_r\lesssim (C)^k\lambda_{q+1}^{r+j}L_p(k)\mathcal{T}_p^k(\ell\lambda_q)^{m_0}\tau^a\delta_q^{1/2}\delta_{q+1}^{1/2} \ \text{ for } 0\leq r+j+k\leq \bar N.$$
Similarly, for $0\leq r+k \leq \bar N-1$, we have:
\begin{equation*}
    \begin{split}
        &||\mathcal{\tilde A}_{\ell,I}^\pm\mathcal{L}_{\xi^{p}_I}^k\mathcal{L}_{\tilde v_q-\tilde v_{\ell}}\Theta_I^{p}||_{r}\\
        &=||(\partial_t+\mathcal{L}_{\tilde z_{\ell,I}^\pm})\mathcal{L}_{\xi^{p}_I}^k\mathcal{L}_{\tilde v_q-\tilde v_{\ell}}\Theta_I^{p}||_{r}+||\mathcal{L}_{\xi^{p}_I}^k\mathcal{L}_{\tilde v_q-\tilde v_{\ell}}\Theta_I^{p}\cn \tilde z^\pm_{\ell,I}||_r\\
        &\lesssim ||a_{1,k,\mathcal{A}}||_r||\tilde \Psi^{1*}_I(\phi_{1,k}\nu_I)||_0+||a_{1,k,\mathcal{A}}||_0||\tilde \Psi^{1*}_I(\phi_{1,k}\nu_I)||_r+||a_{2,k,\mathcal{A}}||_r||\tilde \Psi^{1*}_I(\phi_{2,k}k_I)||_0+||a_{2,k,\mathcal{A}}||_0||\tilde \Psi^{1*}_I(\phi_{2,k}k_I)||_r\\
        &+||\tilde \Psi^{*}_I(\phi_{3,k})||_r||Y_{k,\mathcal{A}}||_0+||\tilde \Psi^{*}_I(\phi_{3,k})||_0||Y_{k,\mathcal{A}}||_r+||\mathcal{L}_{\xi^{p}_I}^k\mathcal{L}_{\tilde v_q-\tilde v_{\ell}}\Theta_I^{p}||_r|| \tilde z^\pm_{\ell,I}||_1+||\mathcal{L}_{\xi^{p}_I}^k\mathcal{L}_{\tilde v_q-\tilde v_{\ell}}\Theta_I^{p}||_0|| \tilde z^\pm_{\ell,I}||_{r+1}\\
        &\lesssim(C)^k \lambda_{q+1}^{r}\left[L_p(k)+(\tau^a/\tau^c)L_{p,\mathcal{A}}(k)\right]\mathcal{T}_p^k(\ell\lambda_q)^{m_0}\delta_q^{1/2}\delta_{q+1}^{1/2}\\
        &+(C)^k\lambda_q^{r}\left[L_p(k+r)+(\tau^a/\tau^c)L_{p,\mathcal{A}}(k+r)\right]\mathcal{T}_p^k(\ell\lambda_q)^{m_0}\delta_q^{1/2}\delta_{q+1}^{1/2}\\
        &+(C)^k\lambda_{q+1}^{r}L_p(k)\mathcal{T}_p^k(\ell\lambda_q)^{m_0}\lambda_q\tau^a\delta_q\delta_{q+1}^{1/2}\\
        &\lesssim(C)^k \lambda_{q+1}^{r}\left[L_p(k)+(\tau^a/\tau^c)L_{p,\mathcal{A}}(k)\right]\mathcal{T}_p^k(\ell\lambda_q)^{m_0}\delta_q^{1/2}\delta_{q+1}^{1/2}\\
        &+(C)^k\lambda_{q+1}^{r}\underbrace{\left(\frac{\lambda_q}{\lambda_{q+1}}\right)^r\left[L_p(k+r)+(\tau^a/\tau^c)L_{p,\mathcal{A}}(k+r)\right]}_{\leq L_p(k)+(\tau^a/\tau^c)L_{p,\mathcal{A}}(k)}\mathcal{T}_p^k(\ell\lambda_q)^{m_0}\delta_q^{1/2}\delta_{q+1}^{1/2}\\
        &+(C)^k\lambda_{q+1}^{r}L_p(k)\mathcal{T}_p^k(\ell\lambda_q)^{m_0}\lambda_q\tau^a\delta_q\delta_{q+1}^{1/2}\\
        &\lesssim (C)^k\lambda_{q+1}^{r}\mathcal{T}_p^k(\ell\lambda_q)^{m_0}\delta_q^{1/2}\delta_{q+1}^{1/2}\left[L_p(k)+(\tau^a/\tau^c)L_{p,\mathcal{A}}(k)\right]
    \end{split}
\end{equation*}
where the same remarks about the constants hold.

\noindent \textbf{Remainders.} The estimate on the remainders can be shown as in the proof of Lemma \ref{remaindersmoll}. In the notation of Lemma \ref{remaindersmoll}, we set $\bar N=N-m_0-1$, which is the largest number of derivatives at our disposal, and fix $k_0=k_0^p$, which is an admissible choice because of \eqref{constraintremainder}.
\end{proof}


\subsubsection{Estimates on the Vector Fields} 
With all the estimates for the intermediate objects making the Lie-Taylor expansions of the Nash perturbation, we are ready to prove bounds on the actual vector fields.

\begin{lemma}[Estimates on $(\mathring w^p, \ \mathring b^p)$]\label{estimatesfieldscircp} Under the choice of parameters in \ref{choiceofparameters}, the following bounds hold.
    \begin{equation*}
        \begin{split}
            &||\partial_t^j\mathring w^p||, \ ||\partial_t^j\mathring b^p||\lesssim  \lambda_{q+1}^{r+j}\mathcal{T}_p\delta_q^{1/2}  \ \text{ for } \ j=0,1,2 \ \text{ and } \ 0\leq r\leq N-j,\\
            &||\mathcal{\tilde A}^\pm\mathring w^p||_r, \ ||\mathcal{\tilde A}^{\pm}\mathring b^p||\lesssim \lambda_{q+1}^r\mathcal{T}_p(1/\tau^a)\delta_q^{1/2} \ \text{ for } \ 0\leq r\leq M-1.\\
        \end{split}
    \end{equation*}
    The implicit constants depend on $r$ and all the parameters, but not on $a$.
\end{lemma}
\begin{proof}[Proof of Lemma \ref{estimatesfieldscircp}] From Proposition \ref{recap} and \ref{standardmollnash} with $\underline{r}=M-m_0-6$ we have:
\begin{equation*}
    \begin{split}
        &||\partial_t^j\tilde v_q||_r, \ ||\partial_t^j\tilde B_q||_r \lesssim \lambda_q^{r+j}\lambda_q^{[r+j-\underline{r}]^+(b-1)\gamma_\ell}\delta_q^{1/2} \ \text{ for }\ j=0,1,2 \ \text{ and } 0\leq r \leq N-j,\\
        &||\partial_t^j(\tilde v_q-\tilde v_\ell)||_r, \ ||\partial_t^j(\tilde B_q-\tilde B_\ell)||_r \lesssim \lambda_q^{r+j}\lambda_q^{[r+j-(\underline{r}-m_0)]^+(b-1)\gamma_\ell}(\ell\lambda_q)^{m_0}\delta_q^{1/2} \ \text{ for }\ j=0,1,2 \ \text{ and } 0\leq r \leq N-m_0-j.\\
    \end{split}
\end{equation*}
We now apply Lemma \ref{estimatesflowp} with $\bar N=N$, to get good estimates on $X^p_s, \ \DD X^p \circ (X^p)^{-1}$ for $|s|\leq 1$. Moreover, we recall the estimates from Lemma \ref{estimatesperturbationnash2}, and let $C$ be the constant from its statement. We finally recall the decomposition in \eqref{decompnashw}, namely
\begin{equation*}
\begin{split}
    \mathring w^p&=
    \begin{cases}
        X^p_*(\tilde v_\ell-\tilde v_q)-(\tilde v_\ell-\tilde v_q),\\
        \curl\left[\sum_I\sum_{k=0}^{k_0^p}\frac{(-1)^k}{(k+1)!}\mathcal{L}_{\xi^p_I}^k\mathcal{L}_{\tilde v_q-\tilde v_\ell}\Theta_I^p+\mathring \theta_w^p\right],\\
    \end{cases}
\end{split}
\end{equation*}
we will use the first expression to prove bounds on pure-space and pure-time derivatives, thereby avoiding a loss. The second expression, together with the bounds in Lemma \ref{estimatesperturbationnash2} and the commutation \eqref{DG1}, namely
\begin{equation}\label{thetatoxi}
    \begin{split}
        \mathcal{A}^\pm_{\ell,I}\curl F&=(\partial_t+\mathcal{L}_{z^\pm_{\ell,I}})\curl F+\curl F\cn z^\pm_{\ell,I}\\
        &=\curl (\partial_t+\mathcal{L}^1_{z^\pm_{\ell,I}})F+\curl F\cn z^\pm_{\ell,I}\\
        &=\curl[\mathcal{A}^\pm_{\ell,I}F]+\curl[(\DD z^\pm_{\ell,I})^\top F]+\curl [F]\cn z^\pm_{\ell,I},
    \end{split}
\end{equation}
will allow us to access the geometric properties of the perturbation and prove the desired Alfv\'en transport bounds. We will use this without mention.

\noindent \textbf{Pure space derivatives.} The observations above together with the composition estimate in Proposition \ref{compestimates}, give:
\begin{equation*}
    \begin{split}
        ||\mathring w^p||_r&=||X^p_*(\tilde v_q-\tilde v_\ell)||_r+||\tilde v_q-\tilde v_\ell||_r\\
        &\lesssim ||\DD X^p \circ (X^p)^{-1}||_r||(\tilde v_q-\tilde v_\ell)\circ (X^p)^{-1}||_0+||\DD X^p \circ (X^p)^{-1}||_0||(\tilde v_q-\tilde v_\ell)\circ (X^p)^{-1}||_r+||\tilde v_q-\tilde v_\ell||_r,\\
    \end{split}
\end{equation*}
now, for $r=0$ we deduce:
\begin{equation*}
    \begin{split}
        ||\mathring w^p||_r&\lesssim ||\DD X^p ||_0||\tilde v_q-\tilde v_\ell||_0+||\tilde v_q-\tilde v_\ell||_0\\
        &\lesssim\mathcal{M}_p(\ell\lambda_q)^{m_0}\delta_q^{1/2}+(\ell\lambda_q)^{m_0}\delta_q^{1/2}\\
        &\leq \mathcal{T}_p\delta_q^{1/2}
    \end{split}
\end{equation*}
where we used  $\frac{\lambda_q}{\lambda_{q+1}}\mathcal{M}_p=\mathcal{T}_p$ and $(\ell \lambda_q)^{m_0}\leq \lambda_q/\lambda_{q+1}\leq \mathcal{T}_p$, while for $r\geq 1$, we have instead:
\begin{equation*}
    \begin{split}
        ||\mathring w^p||_r&\lesssim||\DD X^p \circ (X^p)^{-1}||_r||\tilde v_q-\tilde v_\ell||_0+||\tilde v_q-\tilde v_\ell||_r\\
        &+||\DD X^p \circ (X^p)^{-1}||_0\left[||\tilde v_q-\tilde v_\ell||_1||\DD X^p_{-1}||_{r-1}+||\tilde v_q-\tilde v_\ell||_{r}||\DD X^p_{-1}||_0^r\right],\\
    \end{split}
\end{equation*}
where we used that $(X^p)^{-1}=X^p_{-1}$, for $1\leq r \leq N-m_0$, we have access to the mollification estimates recalled above, and we conclude: 
\begin{equation*}
    \begin{split}
        ||\mathring w^p||_r&\lesssim \lambda_{q+1}^r\mathcal{M}_p(\ell\lambda_q)^{m_0}\delta_q^{1/2}+\lambda_{q+1}^{r-1}\lambda_q\mathcal{M}_p^2(\ell\lambda_q)^{m_0}\delta_q^{1/2}+\lambda_{q}^{r+[r-(\underline{r}-m_0)]^+(b-1)\gamma_\ell}\mathcal{M}_p^{r+1}(\ell\lambda_q)^{m_0}\delta_q^{1/2}\\
        &+\lambda_q^r(\ell\lambda_q)^{m_0}\delta_q^{1/2}\\
        &\lesssim\lambda_{q+1}^r\mathcal{M}_p(\ell\lambda_q)^{m_0}\delta_q^{1/2}+\lambda_{q}^{r+[r-(\underline{r}-m_0)]^+(b-1)\gamma_\ell}\mathcal{M}_p^{r+1}(\ell\lambda_q)^{m_0}\delta_q^{1/2}\\
        &\leq\lambda_{q+1}^r\mathcal{T}_p\delta_q^{1/2}+\lambda_{q}^{r+[r-(\underline{r}-m_0)]^+(b-1)\gamma_\ell}\mathcal{M}_p^{r}\mathcal{T}_p\delta_q^{1/2}\\
        &=\lambda_{q+1}^r\mathcal{T}_p\delta_q^{1/2}+\lambda_{q+1}^r\lambda_{q}^{[r-(\underline{r}-m_0)]^+(b-1)\gamma_\ell}\mathcal{T}_p^{r+1}\delta_q^{1/2}\\
        &= \lambda_{q+1}^r\mathcal{T}_p\delta_q^{1/2}\left[1+\left(\lambda_{q}^{(b-1)\gamma_\ell}\mathcal{T}_p\right)^r\right]\\
        &\lesssim \lambda_{q+1}^r\mathcal{T}_p\delta_q^{1/2}
    \end{split}
\end{equation*}
where we used that $(\ell\lambda_q)^{m_0}\leq \lambda_q/\lambda_{q+1}$ and $\lambda_q/\lambda_{q+1}\mathcal{M}_p=\mathcal{T}_p$ and $\lambda_{q}^{(b-1)\gamma_\ell}\mathcal{T}_p<1$, while for $N-m_0\leq r \leq N$ we can only bound the difference $\tilde v_q-\tilde v_\ell$ as two separate terms and incur in a loss: 
\begin{equation*}
    \begin{split}
        ||\mathring w^p||_r&\lesssim \lambda_{q+1}^r\mathcal{M}_p(\ell\lambda_q)^{m_0}\delta_q^{1/2}+\lambda_{q+1}^{r-1}\lambda_q\mathcal{M}_p^2(\ell\lambda_q)^{m_0}\delta_q^{1/2}+\lambda_{q}^{r+(r-\underline{r})(b-1)\gamma_\ell}\mathcal{M}_p^{r+1}\delta_q^{1/2}\\
        &\lesssim\lambda_{q+1}^r\mathcal{M}_p(\ell\lambda_q)^{m_0}\delta_q^{1/2}+\lambda_{q}^{r+(r-\underline{r})(b-1)\gamma_\ell}\mathcal{M}_p^{r+1}\delta_q^{1/2}\\
        &\lesssim\lambda_{q+1}^r\mathcal{T}_p\delta_q^{1/2}+\lambda_{q}^{r+(r-\underline{r})(b-1)\gamma_\ell}\mathcal{M}_p^{r+1}\delta_q^{1/2}\\
        &\leq\lambda_{q+1}^r\mathcal{T}_p\delta_q^{1/2}\left[1+\left(\lambda_{q}^{(b-1)\gamma_\ell}\mathcal{T}_p\right)^r\left(\frac{\lambda_{q+1}}{\lambda_q}\right)\right],\\
    \end{split}
\end{equation*}
our choice of parameters in \ref{choiceofparameters}, see \eqref{highderivativesconstraint}, guarantees, however, that:
\begin{equation}\label{constraintpurederivativenash}
    \left(\lambda_{q}^{(b-1)\gamma_\ell}\mathcal{T}_p\right)^{N-m_0-2}\leq \frac{\lambda_q}{\lambda_{q+1}}
\end{equation}
and we conclude:
$$||\mathring w^p||_r\lesssim \lambda_{q+1}^r\mathcal{T}_p\delta_q^{1/2} \text{ for } \ 0\leq r \leq N$$
where the implicit constant depends on $r$ and all the parameters, but not on $a$. The additional $-2$ in \eqref{constraintpurederivativenash} will be needed when considering $\partial_t, \ \partial_t^2$. 

\noindent \textbf{Pure time derivatives.} Set $M=\DD X^p\circ(X^p)^{-1}$, we begin with the following calculations:
\begin{equation*}
    \begin{split}
       \partial_t \left[X^p_*(\tilde v_q-\tilde v_\ell)\right]&=\partial_tM(\tilde v_q-\tilde v_\ell) \circ(X^p)^{-1}+M\left[\partial_t(\tilde v_q-\tilde v_\ell) \circ(X^p)^{-1}+\DD(\tilde v_q-\tilde v_\ell) \circ(X^p)^{-1}[\partial_t (X^p)^{-1}]\right]\\
    \end{split}
\end{equation*}
and 
\begin{equation*}
    \begin{split}
        \partial_t^2 \left[X^p_*(\tilde v_q-\tilde v_\ell)\right]&=\partial_t^2M(\tilde v_q-\tilde v_\ell) \circ(X^p)^{-1}+\partial_tM\left[\partial_t(\tilde v_q-\tilde v_\ell) \circ(X^p)^{-1}+\DD(\tilde v_q-\tilde v_\ell) \circ(X^p)^{-1}[\partial_t (X^p)^{-1}]\right]\\
       &+M\left[\partial_t^2(\tilde v_q-\tilde v_\ell) \circ(X^p)^{-1}+2\DD[\partial_t(\tilde v_q-\tilde v_\ell)] \circ(X^p)^{-1}\partial_t (X^p)^{-1}\right]\\
       &+M\left[\DD^2[\tilde v_q-\tilde v_\ell] \circ(X^p)^{-1}[\partial_t (X^p)^{-1},\partial_t (X^p)^{-1}]+\DD(\tilde v_q-\tilde v_\ell) \circ(X^p)^{-1}\partial_t^2 (X^p)^{-1}\right].
    \end{split}
\end{equation*}
We can now argue as above and deduce: 
\begin{equation*}
    \begin{split}
        &||\partial_t \left[X^p_*(\tilde v_q-\tilde v_\ell)\right]||_r\\
        &\leq ||\partial_tM||_r||(\tilde v_q-\tilde v_\ell) \circ(X^p)^{-1}||_0+||\partial_tM||_0||(\tilde v_q-\tilde v_\ell) \circ(X^p)^{-1}||_r\\
        &+||M||_r||\partial_t(\tilde v_q-\tilde v_\ell) \circ(X^p)^{-1}||_0+||M||_0||\partial_t(\tilde v_q-\tilde v_\ell) \circ(X^p)^{-1}||_r\\
        &+||M||_r||\DD(\tilde v_q-\tilde v_\ell) \circ(X^p)^{-1}\partial_t (X^p)^{-1}||_0+||M||_0||\DD(\tilde v_q-\tilde v_\ell) \circ(X^p)^{-1}[\partial_t (X^p)^{-1}]||_r\\
        &\overset{r=0}{\lesssim}\lambda_{q+1}\mathcal{T}_p\delta_q^{1/2}\\
        &\overset{r\geq 1}{\lesssim} \lambda_{q+1}^{r+1}\mathcal{M}_p(\ell\lambda_q)^{m_0}\delta_q^{1/2}+\lambda_{q+1}\mathcal{M}_p\left[||\tilde v_q-\tilde v_\ell||_1||\DD X^p_{-1}||_{r-1}+||\tilde v_q-\tilde v_\ell||_r||\DD X^p_{-1}||_{0}^r\right]\\
        &+\lambda_{q+1}^{r}\mathcal{M}_p\lambda_q(\ell\lambda_q)^{m_0}\delta_q^{1/2}+\mathcal{M}_p\left[||\partial_t(\tilde v_q-\tilde v_\ell)||_1||\DD X^p_{-1}||_{r-1}+||\partial_t(\tilde v_q-\tilde v_\ell)||_r||\DD X^p_{-1}||_0^r\right]\\
        &+\lambda_{q+1}^{r}\mathcal{M}_p^2\lambda_q(\ell\lambda_q)^{m_0}\delta_q^{1/2}\\
        &+\mathcal{M}_p\left[[||\DD(\tilde v_q-\tilde v_\ell)||_1 ||\DD X^p_{-1}||_{r-1}+||\DD(\tilde v_q-\tilde v_\ell)||_r ||\DD X^p_{-1}||_0^r]||\partial_t (X^p)^{-1}||_0+||\DD(\tilde v_q-\tilde v_\ell)||_0||\partial_t (X^p)^{-1}||_r\right]\\
        &\lesssim \lambda_{q+1}^{r+1}\mathcal{M}_p(\ell\lambda_q)^{m_0}\delta_q^{1/2}+\lambda_{q+1}^{r}\mathcal{M}_p^2\lambda_q(\ell\lambda_q)^{m_0}\delta_q^{1/2}+\lambda_q^{r+1}\lambda_q^{[r-(\underline{r}-m_0-1)]^+(b-1)\gamma_\ell}\mathcal{M}_p^{r+2}(\ell\lambda_q)^{m_0}\delta_q^{1/2}\\
        &\lesssim \lambda_{q+1}^{r+1}\mathcal{T}_p\delta_q^{1/2}+\lambda_{q+1}^{r+1}\underbrace{\lambda_q^{[r-(\underline{r}-m_0-1)]^+(b-1)\gamma_\ell}\mathcal{T}_p^{r+1}}_{\leq 1}\mathcal{T}_p\delta_q^{1/2}\\
        &\lesssim \lambda_{q+1}^{r+1}\mathcal{T}_p\delta_q^{1/2}\\
    \end{split}
\end{equation*}
for $0\leq r \leq N-m_0-1$, while for $N-m_0-1\leq r\leq N$, we do not have access to the improved mollification bounds on the difference $\tilde v_q-\tilde v_\ell$, and we are forced to estimate them separately. Proceeding as above, by means of \eqref{constraintpurederivativenash} to reabsorb the loss, we obtain:
\begin{equation*}
    \begin{split}
        ||\partial_t \left[X^p_*(\tilde v_q-\tilde v_\ell)\right]||_r&\lesssim\lambda_{q+1}^{r+1}\mathcal{T}_p\delta_q^{1/2}+\lambda_q^{r+1}\lambda_q^{[r+1-(\underline{r}-m_0)](b-1)\gamma_\ell}\mathcal{M}_p^{r+2}\delta_q^{1/2}\\
        &\leq \lambda_{q+1}^{r+1}\mathcal{T}_p\delta_q^{1/2}+\lambda_{q+1}^{r+1}\lambda_q^{[r+1](b-1)\gamma_\ell} \mathcal{T}^{r+2}_p\left(\frac{\lambda_{q+1}}{\lambda_q}\right)\delta_q^{1/2}\\
        &=\lambda_{q+1}^{r+1}\mathcal{T}_p\delta_q^{1/2}\left[1+\underbrace{\left(\lambda_q^{(b-1)\gamma_\ell} \mathcal{T}_p\right)^{r+1}\left(\frac{\lambda_{q+1}}{\lambda_q}\right)}_{\leq 1}\right]\\
        &\lesssim \lambda_{q+1}^{r+1}\mathcal{T}_p\delta_q^{1/2}.
    \end{split}
\end{equation*}
Similarly, one can show:
$$||\partial_t^2 \left[X^p_*(\tilde v_q-\tilde v_\ell)\right]||_r\lesssim \lambda_{q+1}^{r+2}\mathcal{T}_p\delta_q^{1/2} \ \text{ for } \ 0\leq r \leq N-2$$
and we conclude that:
$$||\partial_t^j\mathring w^p||\leq ||\partial_t^j \left[X^p_*(\tilde v_q-\tilde v_\ell)\right]||_r+||\partial_t^j(\tilde v_q-\tilde v_\ell)||  \lesssim \lambda_{q+1}^{r+j}\mathcal{T}_p\delta_q^{1/2}  \ \text{ for } \ j=0,1,2 \ \text{ and } \ 0\leq r\leq N-j$$
where the implicit constant depends on $r$ and all the parameters, but not on $a$. The estimates for $\mathring b^p$ can be obtained similarly.

\noindent \textbf{Alfv\'en transport estimates.} We first bound the Alfv\'en transport of the remainder with pure derivatives. According to the estimates in Lemma \ref{estimatesperturbationnash2}, we get: 
\begin{equation*}
    \begin{split}
        ||\mathcal{\tilde A}^\pm\curl[\mathring \theta ^p_w]||_{r}&\lesssim||\partial_t\mathring \theta^p_w||_{r+1}+||\tilde z^\pm||_r||\mathring \theta^p_w||_2+||\tilde z^\pm||_0||\mathring \theta^p_w||_{r+2}\\
        &\lesssim \frac{(C)^{k_0^p+1}}{(k_0^p+2)!}\lambda_{q+1}^{r+2}\mathcal{M}_p^{r+2}\mathcal{T}_p^{k_0^p+1}L_p(k_0^p+1)(\ell\lambda_q)^{m_0}\tau^a\delta_q^{1/2}\delta_{q+1}^{1/2},\\
    \end{split}
\end{equation*}
for $0\leq r+(k_0^p+1)+3+m_0\leq N$. Our choice of parameters in \ref{choiceofparameters}, see \eqref{constraintremainder}, guarantees in particular that
$$\lambda_{q+1}\mathcal{M}_p^{r+2}\left(\lambda_q^{(b-1)\gamma_\ell}\mathcal{T}_p\right)^{k_0^p+1}\leq (1/\tau^a) \ \text{ for } \ 0\leq r \leq M-1$$
and we conclude 
$$||\mathcal{\tilde A}^\pm\curl \ \mathring \theta^p_w||_{r}\lesssim \lambda_{q+1}^{r}\mathcal{T}_p(1/\tau^a)\delta_q^{1/2} \ \text{ for } \ 0\leq r \leq M-1.$$

\noindent We now deal with the explicit part of the series, let
$$S_k=\sum_I\mathcal{L}_{\xi^p_I}^k\mathcal{L}_{\tilde v_q-\tilde v_\ell}\Theta_I^p \ \ \text{ for } \ \ k\geq 0,$$
and according to Lemmas \ref{standardmollnash}, \ref{stabilitynash}, \ref{estimatesperturbationnash} and the remarks in the preliminary part of this proof, we have:
\begin{equation*}
    \begin{split}
        ||\mathcal{\tilde A}^\pm \curl \ S_k||_r&\lesssim \sup_I\left[||\mathcal{\tilde A}^\pm_{\ell,I}\curl\left[\mathcal{L}_{\xi^p_I}^k\mathcal{L}_{\tilde v_q-\tilde v_\ell}\Theta_I^p\right]||_r+ ||(\tilde z^\pm_q-z^\pm_{\ell,I})\cn \curl\left[\mathcal{L}_{\xi^p_I}^k\mathcal{L}_{\tilde v_q-\tilde v_\ell}\Theta_I^p\right]||_r\right]\\
        &\lesssim (C)^k\lambda_{q+1}^{r+1}\mathcal{T}_p^{k}(\ell\lambda_q)^{m_0}\delta_q^{1/2}\delta_{q+1}^{1/2}\left[L_p(k)+(\tau^a/\tau^c)L_{p,\mathcal{A}}(k)\right]\\
        &+(C)^k\lambda_{q+1}^{r+2}\mathcal{T}_p^kL_p(k)(\ell\lambda_q)^{2m_0}\tau^a\delta_q\delta_{q+1}^{1/2}\\
        &\lesssim(C)^k\lambda_{q+1}^{r}\mathcal{T}_p^{k}\lambda_q\delta_q^{1/2}\delta_{q+1}^{1/2}\left[L_p(k)+(\tau^a/\tau^c)L_{p,\mathcal{A}}(k)\right]
    \end{split}
\end{equation*}
for $0\leq r+k\leq N-m_0-3$, where we used that $(\ell\lambda_q)^{m_0}\leq \lambda_q/\lambda_{q+1}$ and $L_p\leq L_{p,\mathcal{A}}$, together with the fact that at most one $\Theta_I^p$ is non-zero for each space-time point, see Lemma \ref{disjoint}. 

\noindent Gathering all the estimates we obtain:
\begin{equation*}
    \begin{split}
        ||\mathcal{\tilde A}^\pm\mathring w^p||_r&\leq\sup_I\sum_{k=0}^{k_0^p}\frac{1}{(k+1)!}||\mathcal{\tilde A}^\pm \curl \ S_k||_r+||\mathcal{\tilde A}^\pm\curl[\mathring \theta^p_w]||_r\\
        &\lesssim  \underbrace{\sum_{k=0}^{k_0^p}\frac{(C)^k}{(k+1)!}\lambda_{q+1}^{r}(\tau^a/\tau^c)L_{p,\mathcal{A}}(k)\mathcal{T}_p^k\lambda_q\delta_q^{1/2}\delta_{q+1}^{1/2}}_{T_1}+\underbrace{\sum_{k=0}^{k_0^p}\frac{(C)^k}{(k+1)!}\lambda_{q+1}^{r}L_p(k)\mathcal{T}_p^k\lambda_q\delta_q^{1/2}\delta_{q+1}^{1/2}}_{T_2}\\
        &+\lambda_{q+1}^{r}\mathcal{T}_p(1/\tau^a)\delta_q^{1/2}\\
    \end{split}
\end{equation*}
for $0\leq r+(k_0^p+1)+3+m_0\leq N$.

\noindent We now split the sum in $T_1$ according with the explicit formula for $L_{p,\mathcal{A}}$ in \eqref{lossparameters}, namely
$$L_{p,\mathcal{A}}(k)=1_{k\leq \underline{r}-k_0^g-2}+1_{\underline{r}-k_0^g-1\leq k \leq \underline{r}-1}\left[1+\left(\lambda_q^{(b-1)\gamma_\ell}\mathcal{T}_g\right)^{\underline{r}-1-k}\bar L\right]+1_{k\geq \underline{r}}\lambda_q^{[k-(\underline{r}-1)](b-1)\gamma_\ell}\bar L$$
where $\underline{r}=M-2m_0-k_0^g-8$ and obtain:
\begin{equation}\label{circtransportestimate}
\begin{split}
    T_1&=\lambda_{q+1}^r(\tau^a/\tau^c)\lambda_q\delta_q^{1/2}\delta_{q+1}^{1/2}\sum_{k=0}^{k_0^p}\frac{(C)^k}{(k+1)!}\mathcal{T}_p^kL_{p,\mathcal{A}}(k)\\
    &=\lambda_{q+1}^r(\tau^a/\tau^c)\lambda_q\delta_q^{1/2}\delta_{q+1}^{1/2}\left[\sum_{k=0}^{\underline{r}-1}\frac{\left(C\mathcal{T}_p\right)^k}{(k+1)!}+\bar L\sum_{k=\underline{r}-k_0^g-1}^{\underline{r}-1}\frac{\left(C\mathcal{T}_p\right)^k}{(k+1)!}\left(\lambda_q^{(b-1)\gamma_\ell}\mathcal{T}_g\right)^{\underline{r}-1-k}\right]\\
    &+\lambda_{q+1}^r(\tau^a/\tau^c)\lambda_q\delta_q^{1/2}\delta_{q+1}^{1/2}\left[\lambda_q^{-(\underline{r}-1)(b-1)\gamma_\ell}\bar L\sum_{k=\underline{r}}^{ k_0^p}\frac{\left(C\lambda_q^{(b-1)\gamma_\ell}\mathcal{T}_p\right)^k}{(k+1)!}\right].
\end{split}
\end{equation}
The first summand is bounded by:
$$\sum_{k=0}^{\underline{r} -1}\frac{\left(C\mathcal{T}_p\right)^k}{(k+1)!}=\frac{1}{C\mathcal{T}_p}\sum_{k=1}^{\underline{r} }\frac{\left(C\mathcal{T}_p\right)^k}{k!}\lesssim 1,$$
we now rewrite and bound the second summand as:
\begin{equation*}
    \begin{split}
        \bar L\sum_{k=\underline{r}-k_0^g-1}^{\underline{r}-1}\frac{\left(C\mathcal{T}_p\right)^k}{(k+1)!}\left(\lambda_q^{(b-1)\gamma_\ell}\mathcal{T}_g\right)^{\underline{r}-1-k}&=\bar L\sum_{k=\underline{r}-k_0^g-1}^{\underline{r}-1}\frac{\left(C\mathcal{T}_p\right)^k}{(k+1)!}\left(\lambda_q^{(b-1)\gamma_\ell}\mathcal{T}_p\right)^{\underline{r}-1-k}\underbrace{(\mathcal{T}_g/\mathcal{T}_p)^{\underline{r}-1-k}}_{\leq 1}\\
        &=\bar L\left(\lambda_q^{(b-1)\gamma_\ell}\mathcal{T}_p\right)^{\underline{r}-1}\sum_{k=\underline{r}-k_0^g-1}^{\underline{r}-1}\frac{(C)^k}{(k+1)!}\\
        &\lesssim 1
    \end{split}
\end{equation*}
where, according to the choice of parameters in \ref{choiceofparameters}, see \eqref{constraintadmissibility} and \eqref{M}, we have: $$(\lambda_q^{(b-1)\gamma_\ell}\mathcal{T}_p)^{\underline{r}-1}\bar L\leq (\lambda_q^{(b-1)\gamma_\ell}\mathcal{T}_p)^{k_0^g}\bar L\leq 1,$$
the last summand can be dealt with similarly:
\begin{equation*}
    \begin{split}
        \lambda_q^{-(\underline{r}-1)(b-1)\gamma_\ell}\bar L\sum_{k=\underline{r}}^{ k_0^p}\frac{\left(C\lambda_q^{(b-1)\gamma_\ell}\mathcal{T}_p\right)^k}{(k+1)!}\lesssim \lambda_q^{-(\underline{r}-1)(b-1)\gamma_\ell}\bar L \left(\lambda_q^{(b-1)\gamma_\ell}\mathcal{T}_p\right)^{\underline{r}}\leq 1
    \end{split}
\end{equation*}
and we conclude:
$$T_1\lesssim \lambda_{q+1}^r(\tau^a/\tau^c)\lambda_q\delta_q^{1/2}\delta_{q+1}^{1/2}.$$

\noindent The bound for $T_2$ is easier, recall from \eqref{lossparameters} that $L_p(k)=\lambda_q^{[k-\underline{r}]^+(b-1)\gamma_\ell}$, plugging this in, we deduce:
\begin{equation*}
    \begin{split}
        T_2&=\sum_{k=0}^{k_0^p}\frac{C_r^k}{(k+1)!}\lambda_{q+1}^{r}L_p(k)\mathcal{T}_p^k\lambda_q\delta_q^{1/2}\delta_{q+1}^{1/2}\\
        &=\lambda_{q+1}^{r}\lambda_q\delta_q^{1/2}\delta_{q+1}^{1/2}\Biggr[\underbrace{\sum_{k=0}^{\underline{r}}\frac{\left(C_r\mathcal{T}_p\right)^k}{(k+1)!}}_{\lesssim 1}+\underbrace{\lambda_q^{-\underline{r}(b-1)\gamma_\ell}\sum_{k=\underline{r}+1}^{k_0^p}\frac{\left(C_r\lambda_q^{(b-1)\gamma_\ell}\mathcal{T}_p\right)^k}{k!}}_{\lesssim \lambda_q^{(b-1)\gamma_\ell}\mathcal{T}_p^{\underline{r}+1}\leq \mathcal{T}_p^{\underline{r}}}\Biggr]\\
        &\lesssim \lambda_{q+1}^{r}\lambda_q\delta_q^{1/2}\delta_{q+1}^{1/2}
    \end{split}
\end{equation*}
and we conclude:
$$||\mathcal{\tilde A}^\pm\mathring w^p||_r\lesssim \lambda_{q+1}^{r}\lambda_q\delta_q^{1/2}\delta_{q+1}^{1/2}= \lambda_{q+1}^r\mathcal{T}_p(1/\tau^a)\delta_q^{1/2},$$
for $0\leq r \leq (k_0^p+1)+3+m_0$ and thus in particular for $0\leq r \leq M-1$, see \eqref{constraintremainder}. The implicit constant depends on $r$ and all the parameters, but not on $a$. The estimates for the magnetic field can be obtained similarly.
\end{proof}

\begin{lemma}[Estimates on $(\bar w^p,\bar b^p)$]\label{estimatesfieldsp} Under the choice of parameters in \ref{choiceofparameters}, the following bounds hold. 
    \begin{equation*}
        \begin{split}
            &||\partial_t^j(\bar w^{p}-w^{p,p})||\lesssim \lambda_{q+1}^{r+j}\mathcal{T}_p\delta_{q}^{1/2}\ \text{ for } \ j=0,1,2 \ \text{ and }\ 0\leq r \leq N-j,\\
            &||\partial_t^j\bar b^{p}||\lesssim \lambda_{q+1}^{r+j}(\tau^a/\tau^c)\delta_{q+1}^{1/2}\ \text{ for } \ j=0,1,2 \ \text{ and }\ 0\leq r \leq N-j,\\
            &||\mathcal{\tilde A}^\pm[\bar w^{p}-w^{p,p}]||\lesssim \lambda_{q+1}^{r}\mathcal{T}_p(1/\tau^a)\delta_{q}^{1/2} \ \text{ for } \ 0\leq r \leq M-1,\\
            &||\mathcal{\tilde A}^{\pm}\bar b^p||\lesssim \lambda_{q+1}^{r}(1/\tau^c)\delta_{q+1}^{1/2} \ \text{ for } \ 0\leq r \leq M-1.\\
        \end{split}
    \end{equation*}
The implicit constants depend on $r$ and all the parameters, but not on $a$.
\end{lemma}
 
\begin{proof}[Proof of Lemma \ref{estimatesfieldsp}] According to the decomposition \eqref{decompnashw} and \eqref{leadingdefinition} we can write:
\begin{equation*}
    \begin{split}
        \bar w^p &=\curl\left[\sum_I\sum_{k=0}^{k_0}\frac{(-1)^k}{(k+1)!}\mathcal{L}_{\xi^p_I}^k(\partial_t+\mathcal{L}_{\tilde v_{\ell,I}})\Theta^p_I+\sum_I\sum_{k=0}^{k_0}\frac{(-1)^k}{(k+1)!}\mathcal{L}_{\xi^p_I}^k\mathcal{L}_{\tilde v_\ell-\tilde v_{\ell,I}}\Theta_I^p+\theta^p_w\right],\\
        w^{p,p}&=\curl\left[\sum_I(\partial_t+\mathcal{L}_{\tilde v_{\ell,I}})\Theta_I^{p}\right]
    \end{split}
\end{equation*}
and in the setting of Lemmas \ref{estimatesperturbationnash} and \ref{remaindersmoll}, given any $b, \  \gamma_a, \ \gamma_\ell$ we fix $\bar N \geq k_0\geq 0$ so that:
\begin{equation}\label{konash}
    \begin{split}
        & \bar N-k_0\geq N+1,\\
        &\mathcal{M}_p^{r+1+j}\left(\lambda_q^{(b-1)\gamma_\ell}\mathcal{T}_p\right)^{k_0+1}\leq \mathcal{T}_p \ \text{ for } \ 0\leq r \leq N-j,\\
        & \lambda_{q+1}\mathcal{M}_p^{r+2}\left(\lambda_q^{(b-1)\gamma_\ell}\mathcal{T}_p\right)^{k_0+1}\leq (1/\tau^a)\mathcal{T}_p \ \text{ for } \ 0\leq r \leq M-1
    \end{split}
\end{equation}
then the Lemmas give bounds up to constants $C, \ C'$ which w.l.o.g. we assume are the same. Note that these depend on $\bar{N}, \ k_0$, which come from the constraints above and thus ultimately depend on $b, \ \beta, \ \gamma_a, \ \gamma_\ell$.

\noindent \textbf{Pure derivatives bounds.} From Lemmas \ref{estimatesperturbationnash} and \ref{remaindersmoll} with $\bar N$ and $k_0$ as above, and the fact that at most one $\Theta_I^p$ is non-zero at each space-time point, see Lemma \ref{disjoint}, we deduce:
\begin{equation*}
\begin{split}
    &||\partial_t^j[\bar w^p-w^{p,p}]||_r\\
    &\leq \sup_I\sum_{k=1}^{k_0}\frac{1}{(k+1)!}||\partial_t^j[\mathcal{L}_{\xi^p_I}^k(\partial_t+\mathcal{L}_{\tilde v_{\ell,I}})\Theta^p_I]||_{r+1}+\sup_I\sum_{k=0}^{k_0}\frac{1}{(k+1)!}||\partial_t^j[\mathcal{L}_{\xi^p_I}^k\mathcal{L}_{\tilde v_\ell-\tilde v_{\ell,I}}\Theta_I^p]||_{r+1}+||\partial_t^j\theta^p_{w}||_{r+1}\\
    &\lesssim \underbrace{\sum_{k=1}^{k_0}\frac{(C)^k}{(k+1)!}\lambda_{q+1}^{r+j}\lambda_q^{[k-\underline{r}]^+(b-1)\gamma_\ell}\mathcal{T}_p^k\delta_{q+1}^{1/2}}_{T_1}+\underbrace{\sum_{k=0}^{k_0}\frac{(C)^k}{(k+1)!}\lambda_{q+1}^{r+1+j}\lambda_q^{[k-\underline{r}]^+(b-1)\gamma_\ell}\mathcal{T}_p^k(\ell\lambda_q)^{m_0}\tau^a\delta_q^{1/2}\delta_{q+1}^{1/2}}_{T_2}\\
    &+\frac{(C)^{k_0+1}}{(k_0+2)!}\lambda_{q+1}^{r+j}\mathcal{M}_p^{r+j+1}\mathcal{T}_p^{k_0+1}\left[\lambda_q^{[k_0+1-\underline{r}]^+(b-1)\gamma_\ell}[1+\tau^a/\tau^c]\delta_{q+1}^{1/2}+\lambda_{q+1}\lambda_q^{[k_0+1-\underline{r}]^+(b-1)\gamma_\ell}\mathcal{T}_p^{k_0+1}(\ell\lambda_q)^{m_0}\tau^a\delta_q^{1/2}\delta_{q+1}^{1/2}\right]\\
    &\lesssim T_1+T_2+\lambda_{q+1}^{r+j}\mathcal{T}_p\delta_{q+1}^{1/2},
\end{split}
\end{equation*}
for $0\leq r\leq N-j$, where we used $\tau^a=\ell_t$ and $(\ell\lambda_q)^{m_0}\leq \lambda_q/\lambda_{q+1}$ and our choice of $\bar N, \ k_0$ in \eqref{konash} to deal with the bound on the remainder. 

\noindent We can estimate $T_1$ by splitting the sum according to the $\gamma_\ell$ loss:
\begin{equation*}
    \begin{split}
        T_1&=\sum_{k=1}^{k_0}\frac{(C)^k}{(k+1)!}\lambda_{q+1}^{r+j}\lambda_q^{[k-\underline{r}]^+(b-1)\gamma_\ell}\mathcal{T}_p^k\delta_{q+1}^{1/2}\\
        &=\lambda_{q+1}^{r+j}\delta_{q+1}^{1/2}\left[\sum_{k=1}^{\underline{r}}\frac{\left(C\mathcal{T}_p\right)^k}{(k+1)!}\lambda_q^{[k-\underline{r}]^+(b-1)\gamma_\ell}\mathcal{T}_p^k+\lambda_q^{-\underline{r}(b-1)\gamma_\ell}\sum_{k=\underline{r}+1}^{k_0}\frac{\left(C\lambda_q^{(b-1)\gamma_\ell}\mathcal{T}_p\right)^k}{(k+1)!}\right]\\
        &\lesssim \lambda_{q+1}^{r+j}\delta_{q+1}^{1/2}\left[\mathcal{T}_p+\lambda_q^{(b-1)\gamma_\ell}\mathcal{T}_p^{\underline{r}+1}\right]\leq \lambda_{q+1}^{r+j}\delta_{q+1}^{1/2}\left[\mathcal{T}_p+\mathcal{T}_p^{\underline{r}}\right]\\
        &\lesssim \lambda_{q+1}^{r+j}\mathcal{T}_p\delta_{q+1}^{1/2}.
    \end{split}
\end{equation*}
We can proceed similarly for $T_2$, using in addition that $\tau^a\lambda_q\delta_q^{1/2}\delta_{q+1}^{1/2}=\mathcal{T}_p\delta_q^{1/2}$, see \eqref{Tp}, and $(\ell\lambda_q)^{m_0}\leq \lambda_q/\lambda_{q+1}$, we deduce:
\begin{equation*}
    \begin{split}
        T_2&\lesssim \lambda_{q+1}^{r+j}\mathcal{T}_p\delta_{q}^{1/2}
    \end{split}
\end{equation*}
and note that this is, in fact, a worse bound compared to the one for $T_1$, we conclude that:
$$||\partial_t^j[\bar w^p-w^{p,p}]||_r\lesssim \lambda_{q+1}^{r+j}\mathcal{T}_p\delta_{q}^{1/2} \ \text{ for } \ 0\leq r \leq N-j$$
where the implicit constant depends on $r$ and all the parameters but not on $a$. 

\noindent An analogous argument, together with the additional $\tau^a/\tau^c$ smallness associated with the magnetic field and the fact that $(\tau^a/\tau^c)\delta_{q+1}^{1/2}\geq \mathcal{T}_p\delta_q^{1/2}$, see \eqref{Tp}, shows: 
$$||\partial_t^j\bar b^p||_r\lesssim \lambda_{q+1}^{r+j}(\tau^a/\tau^c)\delta_{q+1}^{1/2} \ \text{ for } \ 0\leq r \leq N-j.$$
Note that we are not removing the leading term in the Lie-Taylor expansion here.

\noindent \textbf{Alfv\'en transport bounds.} We first recall the observations made in the preliminary section of the proof of Lemma \ref{estimatesfieldscircp}, see \eqref{thetatoxi}, and also set: 
$$S_1^k=\sum_I\mathcal{L}_{\xi^p_I}^k(\partial_t+\mathcal{L}_{\tilde v_{\ell,I}})\Theta^p_I \ \text{ and } \ S_2^k=\sum_I\mathcal{L}_{\xi^p_I}^k\mathcal{L}_{\tilde v_\ell-\tilde v_{\ell,I}}\Theta_I^p$$ 
so that we can write:
\begin{equation}\label{velocityperturbationalfven}
    \bar w^p-w^{p,p}=\curl\left[\sum_{k=1}^{k_0^p}\frac{(-1)^k}{(k+1)!}S_1^k+\sum_{k=0}^{k_0^p}\frac{(-1)^k}{(k+1)!}S_2^k\right].
\end{equation}
Since at most one index $I$ is non-zero at each space-time point, see Lemma \ref{disjoint}, the estimates in \ref{estimatesperturbationnash}, \ref{remaindersmoll}, \ref{stabilitynash}, give:
\begin{equation*}
    \begin{split}
        ||\mathcal{\tilde A}^\pm_{\ell,I}\curl S_1^k||_r&\lesssim (C)^k\lambda_{q+1}^{r}\lambda_q^{[k-\underline{r}]^+(b-1)\gamma_\ell}\mathcal{T}_p^k1/\tau^a\delta_{q+1}^{1/2},\\
        ||(\tilde z^\pm_q-\tilde z ^\pm_{\ell,I})\cn\curl S_1^k||_r&\lesssim(C)^k\lambda_{q+1}^{r+1}\lambda_q^{[k-\underline{r}]^+(b-1)\gamma_\ell}\mathcal{T}_p^k(\ell\lambda_q)^{m_0}\delta_q^{1/2}\delta_{q+1}^{1/2},\\
        ||\mathcal{\tilde A}^\pm_{\ell,I}\curl S_2^k||_r&\lesssim (C)^k\lambda_{q+1}^{r+1}\lambda_q^{[k-\underline{r}]^+(b-1)\gamma_\ell}\mathcal{T}_p^{k}(\ell\lambda_q)^{m_0}\delta_q^{1/2}\delta_{q+1}^{1/2},\\
        ||(\tilde z^\pm_q-\tilde z ^\pm_{\ell,I})\cn\curl S_2^k||_r&\lesssim (C)^k\lambda_{q+1}^{r+2}\lambda_q^{[k-\underline{r}]^+(b-1)\gamma_\ell}\mathcal{T}_p^k(\ell\lambda_q)^{2m_0}\delta_q\tau^a\delta_{q+1}^{1/2}
    \end{split}
\end{equation*}
and
\begin{equation*}
    \begin{split}
        ||\mathcal{\tilde A}^\pm_{\ell,I}\curl \theta^p_w||_r&\leq||\partial_t\theta^p_{w}||_{r+1}+||\tilde z^\pm_q||_r||\theta^p_w||_2+||\tilde z^\pm_q||_0||\theta^p_w||_{r+2}\\
        &\lesssim \frac{(C)^{k_0+1}}{(k_0+2)!}\lambda_{q+1}^{r+1}\mathcal{M}_p^{r+2}\mathcal{T}_p^{k_0+1}\lambda_q^{[k_0+1-\underline{r}]^+(b-1)\gamma_\ell}\left[[1+\tau^a/\tau^c]\delta_{q+1}^{1/2}+\lambda_{q+1}(\ell\lambda_q)^{m_0}\tau^a\delta_q^{1/2}\delta_{q+1}^{1/2}\right]\\
        &\lesssim \lambda_{q+1}^{r+1}\mathcal{M}_p^{r+2}\left(\lambda_q^{(b-1)\gamma_\ell}\mathcal{T}_p\right)^{k_0+1}\delta_{q+1}^{1/2}\\
        &\leq \lambda_{q+1}^{r} 1/\tau^a\mathcal{T}_p\delta_{q+1}^{1/2} 
    \end{split}
\end{equation*}
for $0\leq r \leq M-1$, where we used our choice of $k_0$ in \eqref{konash} and $(\ell\lambda_q)^{m_0}\leq \lambda_q/\lambda_{q+1}$. 

\noindent Gathering the bounds above, for $0\leq r\leq M-1$, and going back to the rewriting in \eqref{velocityperturbationalfven}, we deduce:
\begin{equation*}
\begin{split}
    &||\mathcal{\tilde A}^\pm [\bar w^p-w^{p,p}]||_r\\
    &\leq \sum_{k=1}^{k_0}\frac{1}{(k+1)!}||\mathcal{\tilde A}^\pm_{\ell,I}\curl S_1^k||_r+\sum_{k=0}^{k_0}\frac{1}{(k+1)!}||\mathcal{\tilde A}^\pm_{\ell,I} \curl S_2^k||_r\\
    &+\sum_{k=1}^{k_0}\frac{1}{(k+1)!}||(\tilde z^\pm_q-\tilde z ^\pm_{\ell,I})\cn\curl S_1^k||_r+\sum_{k=0}^{k_0}\frac{1}{(k+1)!}||(\tilde z^\pm_q-\tilde z ^\pm_{\ell,I})\cn\curl S_2^k||_r\\
    &+||\mathcal{\tilde A}^\pm\curl\theta^p_{w}||_{r}\\
    &\lesssim \underbrace{\sum_{k=1}^{k_0}\frac{(C)^k}{(k+1)!}\lambda_{q+1}^{r}\lambda_q^{[k-\underline{r}]^+(b-1)\gamma_\ell}\mathcal{T}_p^k1/\tau^a\delta_{q+1}^{1/2}}_{T_1}+\underbrace{\sum_{k=0}^{k_0}\frac{(C)^k}{(k+1)!}\lambda_{q+1}^{r+1}\lambda_q^{[k-\underline{r}]^+(b-1)\gamma_\ell}\mathcal{T}_p^{k}(\ell\lambda_q)^{m_0}\delta_q^{1/2}\delta_{q+1}^{1/2}}_{T_2}\\
    &+\sum_{k=1}^{k_0}\frac{(C)^k}{(k+1)!}\lambda_{q+1}^{r+1}\lambda_q^{[k-\underline{r}]^+(b-1)\gamma_\ell}\mathcal{T}_p^k(\ell\lambda_q)^{m_0}\delta_q^{1/2}\delta_{q+1}^{1/2}+\sum_{k=0}^{k_0}\frac{(C)^k}{(k+1)!}\lambda_{q+1}^{r+2}\lambda_q^{[k-\underline{r}]^+(b-1)\gamma_\ell}\mathcal{T}_p^k(\ell\lambda_q)^{2m_0}\delta_q\tau^a\delta_{q+1}^{1/2}\\
    &+\lambda_{q+1}^r\mathcal{T}_p(1/\tau^a)\delta_{q+1}^{1/2}\\
    &\lesssim T_1+T_2+\lambda_{q+1}^r\mathcal{T}_p(1/\tau^a)\delta_{q+1}^{1/2}.
\end{split}
\end{equation*}

\noindent We can proceed with $T_1, \ T_2$ as in the pure derivatives bounds above by splitting the sum according to the $\gamma_\ell$ loss. We deduce that:
$$T_1\lesssim \lambda_{q+1}^r\mathcal{T}_p(1/\tau^a)\delta_{q+1}^{1/2} \text{ and } T_2\lesssim \lambda_{q+1}^r\mathcal{T}_p(1/\tau^a)\delta_{q}^{1/2}$$
and we conclude:
$$||\mathcal{\tilde A}^\pm[\bar w^p-w^{p,p}]||_r\lesssim \lambda_{q+1}^r\mathcal{T}_p(1/\tau^a)\delta_{q}^{1/2} \ \text{ for } \ 0\leq r \leq M-1$$
where the implicit constant depends on $r$ and all the other parameters but not on $a$.

\noindent A similar argument, together with the additional $\tau^a/\tau^c$ smallness associated with the magnetic field, and the fact that $1/\tau^c\delta_{q+1}^{1/2}\geq \mathcal{T}_p1/\tau^a\delta_q^{1/2}$, see \eqref{Tp}, shows: 
$$||\mathcal{\tilde A}^\pm\bar b^p||_r\lesssim \lambda_{q+1}^r1/\tau^c\delta_{q+1}^{1/2}\ \text{ for } \ 0\leq r \leq M-1$$
and the same remark regarding the implicit constant applies here.
\end{proof}

\subsection{Reynolds Stress: \texorpdfstring{$R^p$}{Rp}}\label{Rps}
According to the decompositions \eqref{Rq+1} and \eqref{Rgp}, where in fact $\pi^p=0$, we have:
\begin{equation}\label{Rp}
    \begin{split}
        \ddiv \left[R^p+R^g\right]&=\partial_tv_{q+1}+\ddiv(v_{q+1}\otimes v_{q+1}-B_{q+1}\otimes B_{q+1})+\nabla p_{q+1}\\
        &=\underbrace{\partial_tw^p+\ddiv\left[(v_q+w^{g})\otimes w^p+w^p\otimes (v_q+w^{g})-(B_q+b^{g})\otimes b^p-b^p\otimes(B_q+b^{g})\right]}_{\ddiv \ R^{p,lin}}\\
        &+\underbrace{\ddiv \left[w^p\otimes w^p-b^p\otimes b^p -\sum_I g_I^2 \tilde{A}_I\right]}_{\ddiv \ R^{p,qua}}+\ddiv\underbrace{\left[\sum_I g_I^2 (\tilde{A}_I-A_I)\right]}_{R^{p,crt}}\\
        &=\ddiv \left[R^{p,lin}+R^{p,qua}+R^{p,crt}\right]
    \end{split}
\end{equation}
where, following \eqref{Rg} we set 
$$R^p=R^{p,lin}+R^{p,qua}+R^{p,crt}$$
and enforced the identities
\begin{equation}\label{enforced}
    \begin{split}
        \ddiv R^{p,lin}&=\partial_tw^p+\ddiv\left[(v_q+w^{g})\otimes w^p+w^p\otimes (v_q+w^{g})-(B_q+b^{g})\otimes b^p-b^p\otimes (B_q+b^{g})\right]\\
        &=\partial_tw^p+\ddiv\left[\tilde v_q\otimes w^p+w^p\otimes \tilde v_q-\tilde B_q\otimes b^p-b^p\otimes \tilde B_q\right],\\
        \ddiv R^{p,qua}&=\ddiv \Biggr[w^p\otimes w^p-b^p\otimes b^p -\sum_I g_I^2 \tilde{A}_I\Biggr]
    \end{split}
\end{equation}
so that $R^{p,lin}$ corresponds to the linear error regarding $(w^p, \ b^p)$ as a perturbation of $(\tilde v_q, \ \tilde B_q)$ and $R^{p,qua}$
corresponds to the quadratic one; their precise definitions will be given after further rewriting. $R^{p,crt}$ was defined already in Lemma \ref{chartstability} and is due to the charts update procedure in the construction. 

\noindent \textbf{Linear Errors.} To shorten the notation, we define:
\begin{equation}\label{potentialsnashw1}
    w_I=\curl[\Theta_{w,I}], \ r_{w}=\curl[\theta_{w}]
\end{equation}
where
\begin{equation}\label{potentialsnashw}
    \begin{cases}
        \Theta_{w,I}=\sum_{k=0}^{k_0^p}\frac{(-1)^k}{(k+1)!}\mathcal{L}_{\xi^p_I}^k(\partial_t+\mathcal{L}_{\tilde v_{\ell,I}})\Theta^p_I+\sum_{k=0}^{k_0^p}\frac{(-1)^k}{(k+1)!}\mathcal{L}_{\xi^p_I}^k\mathcal{L}_{\tilde v_q-\tilde v_{\ell,I}}\Theta_I^p,
        \\
        \theta_{w}=\theta_{w}^p+\mathring\theta_{w}^p
    \end{cases}
\end{equation}
and 
\begin{equation}\label{potentialsnashb1}
    b_I=\curl[\Theta_{b,I}], \ r_{b}=\curl[\theta_{b}]
\end{equation}
where
\begin{equation}\label{potentialsnashb}
    \begin{cases}
        \Theta_{b,I}=\sum_{k=0}^{k_0^p}\frac{(-1)^k}{(k+1)!}\mathcal{L}_{\xi^p_I}^k\mathcal{L}_{\tilde B_{\ell,I}}\Theta^p_I+\sum_{k=0}^{k_0^p}\frac{(-1)^k}{(k+1)!}\mathcal{L}_{\xi^p_I}^k\mathcal{L}_{\tilde B_q-\tilde B_{\ell,I}}\Theta_I^p,
        \\
        \theta_{b}=\theta_{b}^p+\mathring\theta_{b}^p.
    \end{cases}
\end{equation}
This is nothing but a rearrangement of the decompositions in \eqref{decompnashw}, \eqref{decompnashb} and from the calculations there, we deduce:
$$w^p=\sum_Iw_I+r_w=\curl\left[\sum_I \Theta_{w,I}+\theta_w\right], \ b^p=\sum_Ib_I+r_b=\curl\left[\sum_I \Theta_{b,I}+\theta_b\right].$$
We now manipulate \eqref{enforced} as in Subsection \ref{Rgs}, namely
\begin{equation}\label{linen}
    \begin{split}
        \ddiv[R^{p,lin}]&=\partial_tw^p+\ddiv\left[\tilde v_q\otimes w^p+w^p\otimes \tilde v_q-\tilde B_q\otimes b^p-b^p\otimes \tilde B_q\right]\\
        &=\sum_I\left[\partial_tw_{I}+\ddiv\left[\tilde v_q\otimes w_{I}+w_{I}\otimes \tilde v_q-\tilde B_q\otimes b_{I}-b_{I}\otimes \tilde B_q\right]\right]\\
        &+\partial_tr_{w}+\ddiv\left[\tilde v_q\otimes r_{w}+r_{w}\otimes \tilde v_q-\tilde B_q\otimes r_{b}-r_{b}\otimes\tilde B_q\right]\\
        &=\sum_I\left[(\partial_t +\mathcal{L}_{\tilde v_{\ell,I}})w_I-\mathcal{L}_{\tilde B_{\ell,j}}b_I\right]+2\sum_I\left[w_I\cn \tilde v_{\ell,I}-b_I\cn \tilde B_{\ell,I}\right]\\
        &+\sum_I\ddiv\left[(\tilde v_q-\tilde v_{\ell,I})\otimes w_I+w_I\otimes (\tilde v_q-\tilde v_{\ell,I})\right]\\
        &-\sum_I\ddiv\left[(\tilde B_q-\tilde B_{\ell,I})\otimes b_I-b_I\otimes (\tilde B_q-\tilde B_{\ell,I})\right]\\
        &+\partial_tr_{w}+\ddiv\left[\tilde v_q\otimes r_{w}+r_{w}\otimes \tilde v_q-\tilde B_q\otimes r_{b}-r_{b}\otimes \tilde B_q\right]\\
        &=\ddiv\left[R^{p,tr}+R^{p,na}+R^{p,mo}+R^{p,rm}\right]
    \end{split}
\end{equation}
where we defined: 
$$R^{p,lin}=R^{p,tr}+R^{p,na}+R^{p,mo}+R^{p,rm}$$ 
with
\begin{equation*}
    \begin{split}
        R^{p,mo}&=\sum_I\left[(\tilde v_q-\tilde v_{\ell,I})\otimes w_I+w_I\otimes (\tilde v_q-\tilde v_{\ell,I})\right]-\sum_I\left[(\tilde B_q-\tilde B_{\ell,I})\otimes b_I+b_I\otimes (\tilde B_q-\tilde B_{\ell,I})\right],\\
        R^{p,na}&=2\mathcal{R}\sum_I\left[w_I\cn \tilde v_{\ell,I}-b_I\cn \tilde B_{\ell,I}\right]\\
        &=2\mathcal{R}\ddiv\sum_I\left[\Theta_{w,I}\times\nabla\tilde v_{\ell,I}-\Theta_{b,I}\times\nabla \tilde B_{\ell,I}\right]^\top,\\
        R^{p,tr}&=\mathcal{R}\sum_I\left[(\partial_t +\mathcal{L}_{\tilde v_{\ell,I}})w_I-\mathcal{L}_{\tilde B_{\ell,I}}b_I\right]\\
        &=\mathcal{R}\curl\sum_I\left[(\partial_t +\mathcal{L}_{\tilde v_{\ell,I}})\Theta_{w,I}-\mathcal{L}_{\tilde B_{\ell,I}}\Theta_{b,I}\right],\\
        R^{p,rm}&=\mathcal{R}\left[\partial_tr_{w}+\ddiv\left[\tilde v_q\otimes r_{w}+r_{w}\otimes \tilde v_q-\tilde B_q\otimes r_{b}-r_{b}\otimes \tilde B_q\right]\right]\\
        &=\mathcal{R}\curl[\partial_t\vartheta_{w}]+\mathcal{R}\ddiv\left[\tilde v_q\otimes \curl \ \theta_{w}+\curl \ \theta_{w}\otimes \tilde v_q-\tilde B_q\otimes \curl \ \theta_{b}-\curl \ \theta_{b}\otimes \tilde B_q\right]
    \end{split}
\end{equation*}
and we refer the reader to Subsection \ref{Rgs} for more details.

\subsubsection{Estimates on the Transport Error}

\begin{lemma}[Estimates on $R^{p,tr}$]\label{Rptr} Under the choice of parameters in \ref{choiceofparameters}, the following bounds hold.
\begin{equation*}
    \begin{split}
        ||R^{p,tr}||_{r+\alpha}&\lesssim \lambda_{q+1}^r(\tau^c/\tau^a)\frac{\lambda_q\delta_q^{1/2}\delta_{q+1}^{1/2}}{\lambda_{q+1}^{1-2\alpha}} \ \text{for} \ 0\leq r\leq N-k_0^p-m_0-3,\\
        ||\mathcal{\tilde A}^\pm R^{p,tr}||_{r+\alpha}&\lesssim\lambda_{q+1}^r(1/\tau^a)(\tau^c/\tau^a)\frac{\lambda_q\delta_q^{1/2}\delta_{q+1}^{1/2}}{\lambda_{q+1}^{1-2\alpha}} \ \text{ for } 0\leq r\leq N-k_0^p-m_0-4.\\
    \end{split}
\end{equation*}
The implicit constants depend on $r$ and all the other parameters, but not on $a$. In particular, the bounds hold for $0\leq r \leq M$ and $0\leq r \leq M-1$ respectively.
\end{lemma}
\begin{remark} The additional $\lambda_{q+1}^\alpha$ loss in the first estimate is just out of convenience to make the ratio $\tau^c/\tau^a$ appear, and comes from the definition of $\tau^c$ in \eqref{tauc}.
\end{remark}
\begin{proof}[Proof of Lemma \ref{Rptr}] The type of arguments and the splitting are the same as in the proof of Lemma \ref{Rgtr}. We first rewrite $R^{p,tr}$ as:
$$R^{p,tr}=T_1+T_{2,a}+T_{2,b}+T_3+T_4+T_{5,a}+T_{5,b}$$
where
\begin{equation}\label{Rptrequation}
    \begin{split}
        T_1&=\mathcal{R}\curl \sum_{I\in \mathcal{I}}\sum_{k=0}^{k_0^p}\mathcal{L}_{\xi_I^p}^k\left[\left(\partial_t+\mathcal{L}_{\tilde v_{\ell,I}}\right)^2\Theta_I^p-\mathcal{L}_{\tilde B_{\ell,I}}^2\Theta_I^p\right],\\
        T_{2,w}&=\mathcal{R}\curl\sum_{I\in \mathcal{I}}\sum_{k=1}^{k_0^p}\frac{(-1)^k}{(k+1)!}\sum_{i=0}^{k-1}\left[\mathcal{L}_{\xi_I}^i\mathcal{L}_{\partial_t\xi_I^p+[\tilde v_{\ell,I},\xi_I]} \mathcal{L}_{\xi_I^p}^{(k-1)-i}\left(\partial_t+\mathcal{L}_{\tilde v_{\ell,I}}\right)\Theta_I^p \right],\\
        T_{2,b}&=-\mathcal{R}\curl\sum_{I\in \mathcal{I}}\sum_{k=1}^{k_0^p}\frac{(-1)^k}{(k+1)!}\sum_{i=0}^{k-1}\left[\mathcal{L}_{\xi_I^p}^i\mathcal{L}_{[\tilde B_{\ell,I},\xi_I^p]} \mathcal{L}_{\xi_I^p}^{(k-1)-i}\mathcal{L}_{\tilde B_{\ell,I}}\Theta_I^p\right],\\
        T_3&=\mathcal{R}\curl\sum_{I\in \mathcal{I}}\sum_{k=0}^{k_0^p}\frac{(-1)^k}{(k+1)!}\mathcal{L}_{\xi_I^p}^k\left[\mathcal{L}_{\tilde v_q-\tilde v_{\ell,I}}\left(\partial_t+\mathcal{L}_{\tilde v_{\ell,I}}\right)\Theta_I^p-\mathcal{L}_{\tilde B_q-\tilde B_{\ell,I}}\mathcal{L}_{\tilde B_{\ell,I}}\Theta_I^p\right],\\
        T_4&=\mathcal{R}\curl\sum_{I\in \mathcal{I}}\sum_{k=0}^{k_0^p}\frac{(-1)^k}{(k+1)!} \mathcal{L}_{\xi_I^p}^k\left[\mathcal{L}_{\partial_t(\tilde v_q-\tilde v_{\ell,I})+[\tilde v_{\ell,I},\tilde v_q-\tilde v_{\ell,I}]}\Theta_I^p-\mathcal{L}_{[\tilde B_{\ell,I},\tilde B_q-\tilde B_{\ell,I}]}\Theta_I^p\right],\\
        T_{5,w}&=\mathcal{R}\curl\sum_{I\in \mathcal{I}}\sum_{k=1}^{k_0^p}\frac{(-1)^k}{(k+1)!}\sum_{i=0}^{k-1}\left[\mathcal{L}_{\xi_I^p}^i\mathcal{L}_{\partial_t\xi_I^p+[\tilde v_{\ell,I},\xi_I^p]} \mathcal{L}_{\xi_I^p}^{(k-1)-i}\mathcal{L}_{\tilde v_q-\tilde v_{\ell,I}}\Theta_I^p\right],\\
        T_{5,b}&=-\mathcal{R}\curl\sum_{I\in \mathcal{I}}\sum_{k=1}^{k_0^p}\frac{(-1)^k}{(k+1)!}\sum_{i=0}^{k-1}\left[\mathcal{L}_{\xi_I^p}^i\mathcal{L}_{[\tilde B_{\ell,I},\xi_I^p]} \mathcal{L}_{\xi_I^p}^{(k-1)-i}\mathcal{L}_{\tilde B_q-\tilde B_{\ell,I}}\Theta_I^p\right].\\
    \end{split}
\end{equation}
We deal with each $T_i$ separately, the key tool being the Inductive Lemma \ref{inductive}, which we now need in its full generality, compared to Lemma \ref{Rgtr}. We borrow the notation from there. 

\noindent\textbf{Estimates on $T_1$.} Set: 
$$\underline{r}=M-2m_0-k_0^g-8.$$
The top term reads: 
$$S=(\partial_t+\mathcal{L}_{\tilde v_{\ell,I}})^2\Theta_I^p-\mathcal{L}_{\tilde B_{\ell,I}}^2\Theta_I^p,$$
using the definition of $\Theta_I^p$ in \eqref{ansatz} and the Lie-transport properties of the charts in Lemma \ref{chartprop}, we compute:
\begin{equation*}
    \begin{split}
        S&=\frac{1}{\tau^a\lambda_{q+1}}\alpha_I^{''} \tilde a_I \tilde\Psi_I^{1*}(\varphi\nu_I)+\frac{\tau^a}{\lambda_{q+1}}(\partial_t+\tilde v_{\ell,I}\cn)^2( \tilde a_I)\alpha_I \tilde\Psi_I^{1*}(\varphi\nu_I)+\frac{2}{\lambda_{q+1}}(\partial_t+\tilde v_{\ell,I}\cn)( \tilde a_I)\alpha_I'\tilde\Psi_I^{1*}(\varphi\nu_I)\\
        &-\frac{\tau^a}{\lambda_{q+1}}(\tilde B_{\ell,I}\cn)^2( \tilde a_I)\alpha_I \tilde\Psi_I^{1*}(\varphi\nu_I)\\
        &=\frac{1}{\lambda_{q+1}}\left[\frac{1}{\tau^a}\alpha_I^{''} \tilde a_I +(\mathcal{\tilde A}^+_{\ell,I}+\mathcal{\tilde A}^-_{\ell,I})( \tilde a_I)\alpha_I'+\tau^a(\mathcal{\tilde A}^+_{\ell,I}\mathcal{\tilde A}_{\ell,I}^-)( \tilde a_I)\alpha_I\right]\tilde\Psi_I^{1*}(\varphi\nu_I)\\
        &=: a \tilde\Psi_I^{1*}(\varphi\nu_I)=F_0
    \end{split}
\end{equation*}
The estimates contained in Lemmas \ref{chartprop} and \ref{slowcoeffestimates} give:
\begin{equation*}
    \begin{split}
        ||a_1||_{r}&\lesssim \frac{1}{\lambda_{q+1}}\left[\frac{1}{\tau^a}\sup_t|\alpha_I''|\ ||\tilde a_I||_r +\sup_t|\alpha_I'|\ ||(\mathcal{\tilde A}^+_{\ell,I}+\mathcal{\tilde A}^-_{\ell,I})( \tilde a_I)||_r+\tau^a\sup_t|\alpha_I|\ ||(\mathcal{\tilde A}^+_{\ell,I}\mathcal{\tilde A}_{\ell,I}^-)( \tilde a_I)||_r\right]\\
        &\lesssim \frac{1}{\lambda_{q+1}}\lambda_q^{r+[r-\underline{r}]^+(b-1)\gamma_\ell}\delta_{q+1}^{1/2}\left[\frac{1}{\tau^a}+\frac{1}{\tau^c}+\frac{1}{\ell_t}\tau^a/\tau^c\right]\\
        &\lesssim\lambda_q^{r+[r-\underline{r}]^+(b-1)\gamma_\ell}1/\tau^a\frac{\delta_{q+1}^{1/2}}{\lambda_{q+1}}\\
    \end{split}
\end{equation*}
for all $r\geq 0$, where we used $\tau^a=\ell_t$. Moreover,
\begin{equation*}
    \begin{split}
        \mathcal{\tilde A}^\pm_{\ell, I} a_1&= \frac{1}{\lambda_{q+1}}\mathcal{\tilde A}^\pm_{\ell, I}\left[\frac{1}{\tau^a}\alpha_I^{''} \tilde a_I +(\mathcal{A}^+_{\ell,I}+\mathcal{A}^-_{\ell,I})( \tilde a_I)\alpha_I'+\tau^a(\mathcal{\tilde A}^+_{\ell,I}\mathcal{\tilde A}_{\ell,I}^-)( \tilde a_I)\alpha_I\right]\\
        &=\frac{1}{(\tau^a)^2\lambda_{q+1}}\alpha_I^{'''} \tilde a_I +\frac{1}{\tau^a\lambda_{q+1}}\alpha_I^{''}\mathcal{\tilde A}_{\ell,I}^\pm( \tilde a_I) \\
        &+\frac{1}{\lambda_{q+1}\tau^a}(\mathcal{A}^+_{\ell,I}+\mathcal{A}^-_{\ell,I})( \tilde a_I)\alpha_I''+\frac{1}{\lambda_{q+1}}\mathcal{\tilde A}_{\ell,I}^\pm(\mathcal{A}^+_{\ell,I}+\mathcal{A}^-_{\ell,I})( \tilde a_I)\alpha_I'\\
        &+\frac{1}{\lambda_{q+1}}(\mathcal{\tilde A}^+_{\ell,I}\mathcal{\tilde A}_{\ell,I}^-)( \tilde a_I)\alpha_I'+\frac{\tau^a}{\lambda_{q+1}}(\mathcal{\tilde A}_{\ell,I}^\pm\mathcal{\tilde A}^+_{\ell,I}\mathcal{\tilde A}_{\ell,I}^-)( \tilde a_I)\alpha_I
    \end{split}
\end{equation*}
and from the same lemmas we deduce:
\begin{equation*}
    \begin{split}
        ||\mathcal{\tilde A}^\pm_{\ell,I}a_1||_{r}&\lesssim\lambda_q^{r+[r-\underline{r}]^+(b-1)\gamma_\ell}\frac{\delta_{q+1}^{1/2}}{\lambda_{q+1}}\left[\frac{1}{(\tau^a)^2}+\frac{1}{\tau^a\tau^c}+\frac{1}{\ell_t\tau^c}+\frac{1}{\ell_t^2}\tau^a/\tau^c\right]\\
        &\lesssim\lambda_q^{r+[r-\underline{r}]^+(b-1)\gamma_\ell}1/(\tau^a)^2\frac{\delta_{q+1}^{1/2}}{\lambda_{q+1}}\\
    \end{split}
\end{equation*}
for all $r\geq 0$, and we used $ \tau^a=\ell_t$ again. We can now apply the Inductive Lemma \ref{inductive} with:
\begin{equation*}
    \begin{split}
        &z^\pm=\tilde z^\pm_{\ell,I},\\
        &\sigma =\sigma_i=\xi_I,\\
        &\bar \varsigma=\tau^a\delta_{q+1}^{1/2}, \ \bar \varsigma_{\mathcal{A}}=\delta_{q+1}^{1/2},\\
        &A_1=\frac{\delta_{q+1}^{1/2}}{\lambda_{q+1}}\frac{1}{\tau^a}, \ A_{1,\mathcal{A}}=\frac{\delta_{q+1}^{1/2}}{\lambda_{q+1}}\frac{1}{(\tau^a)^2},\\
        &A_3=A_{3,\mathcal{A}}=0,\\
        &L^\varsigma(r)=L^1(r)=L^1_{\mathcal{A}}(r)=\lambda_q^{[r-\underline{r}]^+(b-1)\gamma_\ell}.\\
    \end{split}
\end{equation*}
Note that the needed estimates on $\sigma$ and the charts can be found in \ref{slowcoeffestimates} and \ref{chartprop}. Moreover, we could set any $\bar N$, since every term is mollified; we fix $\bar N=N+k_0^p$ for definiteness. This is not really an issue here, but it will become one when discussing the terms containing non-mollified objects.  The Inductive Lemma gives decompositions:
\begin{equation*}
    \begin{split}
        \mathcal{L}_{\xi_I^p}^kS&=F_k\simeq a_{1,k}\tilde \Psi^{1*}_I(\phi_{1,k}\nu_I),\\
        (\partial_t+\mathcal{L}_{\tilde z^\pm_{\ell,I}})\mathcal{L}_{\xi_I^p}^kS&=(\partial_t+\mathcal{L}_{z^\pm_{\ell,I}})F_k\simeq a_{1,k,\mathcal{A}}\tilde \Psi^{1*}_I(\phi_{1,k}\nu_I)\\
    \end{split}
\end{equation*}
and for $0\leq r+k\leq \bar N$, estimates: 
\begin{equation*}
    \begin{split}
        ||a_{1,k}||_{r}&\leq C'(C)^k \lambda_{q}^{r+k}L(r+k)A_1\prod_{i=1}^k\bar \varsigma_i\\
        &\lesssim (C)^k\lambda_q^{r+[r+k-\underline{r}]^+(b-1)\gamma_\ell}\mathcal{T}^k_p1/\tau^a\frac{\delta_{q+1}^{1/2}}{\lambda_{q+1}},
        \\
        ||\mathcal{\tilde A}^\pm_{\ell,I}a_{1,k}||_r &\leq C'(C)^k \lambda_{q}^{r+k}\left[L_{\mathcal{A}}(k+r)A_{1,\mathcal{A}}\prod_{i=1}^k\bar\varsigma_i +L(k+r)A_1\max_{j\in\{1,\dots,k\}}\bar\varsigma_{j,\mathcal{A}}\prod_{i=1,\dots,k \ i\neq j}\bar\varsigma_i\right]\\
        &\lesssim (C)^k   \lambda_q^{r+[r+k-\underline{r}]^+(b-1)\gamma_\ell}\mathcal{T}^k_p1/(\tau^a)^2\frac{\delta_{q+1}^{1/2}}{\lambda_{q+1}},
    \end{split}
\end{equation*}
note that we keep the same letter for the fast coefficients, for which we have bounds:
$$||\Psi^{1*}_I(\phi_{1,k}\nu_I)||_r\lesssim \lambda_{q+1}^r.$$

\noindent \textit{Conclusion.} From the above, we deduce that:
\begin{equation*}
    \begin{split}
        ||\mathcal{L}_{\xi_I^p}^kS||_{r+\alpha}&\lesssim ||a_{1,k}||_{r+\alpha}||\tilde \Psi^{1*}_I(\phi_{1,k}\nu_I)||_{0}+||a_{1,k}||_{0}||\tilde \Psi^{1*}_I(\phi_{1,k}\nu_I)||_{r+\alpha}\\
        &\lesssim(C)^k\ell^{-\alpha}\lambda_q^{r+[r+k-\underline{r}]^+(b-1)\gamma_\ell}\mathcal{T}^k_p1/\tau^a\frac{\delta_{q+1}^{1/2}}{\lambda_{q+1}}+(C)^k\lambda_{q+1}^{r+\alpha}\lambda_q^{[k-\underline{r}]^+(b-1)\gamma_\ell}\mathcal{T}^k_p1/\tau^a\frac{\delta_{q+1}^{1/2}}{\lambda_{q+1}}\\
        &\lesssim (C)^k\lambda_{q+1}^{r+\alpha}\lambda_q^{[k-\underline{r}]^+(b-1)\gamma_\ell}\mathcal{T}^k_p1/\tau^a\frac{\delta_{q+1}^{1/2}}{\lambda_{q+1}}
        \\
        ||\mathcal{\tilde A}^\pm_{\ell,I}\mathcal{L}_{\xi_I^p}^kS||_{r+\alpha}&\leq ||(\partial_t+\mathcal{L}_{\tilde z^\pm_{\ell,I}})\mathcal{L}_{\xi_I^p}^kS||_{r+\alpha}+||(\mathcal{L}_{\xi_I^p}^kS)\cn \tilde z^\pm_{\ell,I}||_{r+\alpha}\\
        &\lesssim ||a_{1,k,\mathcal{A}}||_{r+\alpha}||\tilde \Psi^{1*}_I(\phi_{1,k}\nu_I)||_{0}+||a_{1,k,\mathcal{A}}||_{0}||\tilde \Psi^{1*}_I(\phi_{1,k}\nu_I)||_{r+\alpha}\\
        &+(C)^k\lambda_{q+1}^{r+\alpha}\lambda_q^{[k-\underline{r}]^+(b-1)\gamma_\ell}\mathcal{T}^k_p1/(\tau^a\tau^c)\frac{\delta_{q+1}^{1/2}}{\lambda_{q+1}}\\
        &\lesssim (C)^k\lambda_q^{r+[r+k-(\underline{r}-1)]^+(b-1)\gamma_\ell}\ell^{-\alpha}\mathcal{T}^k_p1/(\tau^a)^2\frac{\delta_{q+1}^{1/2}}{\lambda_{q+1}}+(C)^k\lambda_{q+1}^{r+\alpha}\lambda_q^{[k-(\underline{r}-1)]^+(b-1)\gamma_\ell}\mathcal{T}^k_p1/(\tau^a)^2\frac{\delta_{q+1}^{1/2}}{\lambda_{q+1}}\\
        &+(C)^k\lambda_{q+1}^{r+\alpha}\lambda_q^{[k-\underline{r}]^+(b-1)\gamma_\ell}\mathcal{T}^k_p1/(\tau^a\tau^c)\frac{\delta_{q+1}^{1/2}}{\lambda_{q+1}}\\
        &\lesssim (C)^k\lambda_{q+1}^{r+\alpha}\lambda_q^{[k-(\underline{r}-1)]^+(b-1)\gamma_\ell}\mathcal{T}^k_p(1/\tau^a)^2\frac{\delta_{q+1}^{1/2}}{\lambda_{q+1}}
    \end{split}
\end{equation*}
The constants depend on $\bar N$ but not on $a$ and are uniform in $I$, $0\leq r+k\leq \bar N$. 

\noindent \textit{Final estimate.} Using Proposition \ref{czstuff} to deal with $\mathcal{R}\curl$ and the fact that there is at most one $\Theta_I^p$ that is non-zero at each space-time point (see Lemma \ref{disjoint}), we conclude that:
\begin{equation*}
    \begin{split}
        ||T_1||_{r+\alpha}&\leq \sum_{k=0}^{k_0^p} \frac{1}{(k+1)!}||\mathcal{R}\curl\sum_I \mathcal{L}_{\xi^p_I}^kS||_{r+\alpha}\lesssim \sum_{k=0}^{k_0^p} \frac{1}{(k+1)!}||\sum_I\mathcal{L}_{\xi^p_I}^kS||_{r+\alpha}\\
        &\lesssim\sup_I\sum_{k=0}^{k_0^p} \frac{1}{(k+1)!}||\mathcal{L}_{\xi^p_I}^kS||_{r+\alpha}\\
        &\lesssim\lambda_{q+1}^{r+\alpha}(1/\tau^a)\frac{\delta_{q+1}^{1/2}}{\lambda_{q+1}}\sum_{k=0}^{k_0^p} \frac{C^{k}}{(k+1)!}\mathcal{T}_p^k\lambda_q^{[k-\underline{r}]^+(b-1)\gamma_\ell}\\
        &\leq\lambda_{q+1}^{r+2\alpha}(\tau^c/\tau^a)\frac{\lambda_q\delta_q^{1/2}\delta_{q+1}^{1/2}}{\lambda_{q+1}}\underbrace{\sum_{k=0}^{k_0^p} \frac{C^k}{(k+1)!}\left(\mathcal{T}_p\lambda_q^{(b-1)\gamma_\ell}\right)^k}_{\lesssim 1}\\
        &\lesssim \lambda_{q+1}^{r}(\tau^c/\tau^a)\frac{\lambda_q\delta_q^{1/2}\delta_{q+1}^{1/2}}{\lambda_{q+1}^{1-2\alpha}}
    \end{split}
\end{equation*}
for $0\leq r \leq   \bar N-k_0^p=N$. The implicit constant depends on $r, \ \bar N$ and all the parameters, but not on $a$. 

\noindent Using again Proposition \ref{czstuff} to deal with $\mathcal{R}\curl, \ [\mathcal{\tilde A}^\pm_{\ell,I},\mathcal{R}\curl]$ and Lemmas \ref{standardmollnash}, \ref{localcorrn} and \ref{stabilitynash} to correct $\mathcal{\tilde A}^\pm$ depending on $I$, we deduce:
\begin{equation}\label{transportT1n}
    \begin{split}
        ||\mathcal{\tilde A}^\pm T_1||_{r+\alpha}&\leq \sum_{k=0}^{k_0^p}\frac{1}{(k+1)!} \left[||\mathcal{R}\curl \sum_I\mathcal{\tilde A}^\pm_{\ell,I}\mathcal{L}_{\xi^{p}_I}^kS||_{r+\alpha}+||[\mathcal{\tilde A}^\pm,\mathcal{R}\curl] \mathcal{L}_{\xi^{p}_I}^kS||_{r+\alpha}\right]\\
        &+\sum_{k=0}^{k_0^p}\frac{1}{(k+1)!}\left[||\mathcal{R}\curl \sum_I(\tilde z_q^\pm-z^\pm_{\ell,I})\cn\mathcal{L}_{\xi^{p}_I}^kS||_{r+\alpha}\right]\\
        &\lesssim \sup_I\sum_{k=0}^{k_0^p}\frac{1}{(k+1)!} \underbrace{\left[||\tilde z^\pm_q||_{r+1+\alpha}||\mathcal{L}_{\xi^{p}_I}^kS||_{\alpha}+||\tilde z^\pm_q||_{1+\alpha}||\mathcal{L}_{\xi^{p}_I}^kS||_{r+1+\alpha}\right]}_{\lambda_{q+1}^{r+\alpha}\ell^{-\alpha}\lambda_q^{[k-\underline{r}]^+(b-1)\gamma_\ell}\mathcal{T}^k_p\frac{\lambda_q\delta_q^{1/2}\delta_{q+1}^{1/2}}{\lambda_{q+1}}\frac{1}{\tau^a}}\\
        &+ \sup_I\sum_{k=0}^{k_0^p}\frac{1}{(k+1)!}\underbrace{\left[||\tilde z_q^\pm-\tilde z^\pm_{\ell,I}||_{r+\alpha}||\mathcal{L}_{\xi^{p}_I}^kS||_{1}+||\tilde z_q^\pm-\tilde z^\pm_{\ell,I}||_{0}|| \mathcal{L}_{\xi^{p}_I}^kS||_{r+1+\alpha}\right]}_{\lambda_{q+1}^{r+1+\alpha}\lambda_q^{[k-\underline{r}]^+(b-1)\gamma_\ell}(\ell\lambda_q)^{m_0}\mathcal{T}^k_p\frac{\delta_q^{1/2}\delta_{q+1}^{1/2}}{\lambda_{q+1}}\frac{1}{\tau^a}}\\
        &+ \sup_I\sum_{k=0}^{k_0^p}\frac{1}{(k+1)!}\underbrace{||\mathcal{\tilde A}^\pm_{\ell,I}\mathcal{L}_{\xi^{p}_I}^kS||_{r+\alpha}}_{\lambda_{q+1}^{r+\alpha}\lambda_q^{[k-(\underline{r}-1)]^+(b-1)\gamma_\ell}\mathcal{T}^k_p\frac{1}{(\tau^a)^2}\frac{\delta_{q+1}^{1/2}}{\lambda_{q+1}}}\\
        &\lesssim\lambda_{q+1}^{r+\alpha}(1/\tau^a)^2\frac{\delta_{q+1}^{1/2}}{\lambda_{q+1}}\underbrace{\sum_{k=0}^{k_0^p} \frac{(C)^{k}}{(k+1)!}\left(\mathcal{T}_p\lambda_q^{(b-1)\gamma_\ell}\right)^k}_{\lesssim 1}\\
        &\lesssim \lambda_{q+1}^r(1/\tau^a)(\tau^c/\tau^a)\frac{\lambda_q\delta_q^{1/2}\delta_{q+1}^{1/2}}{\lambda_{q+1}^{1-2\alpha}}
    \end{split}
\end{equation}
for $0\leq r\leq N-m_0-1$. Here we used $\gamma_a-\gamma_{CZ}\geq 0$, see \eqref{misc1}. The implicit constant depends on $r, \ \bar N$ and all the parameters, but not on $a$. 

\noindent \textbf{Estimates on $T_{4}$.} The top term can be written as:
$$\mathcal{L}_{(\partial_t+\mathcal{L}_{\tilde v_{\ell,I}})(\tilde v_q-\tilde v_{\ell,I})-\mathcal{L}_{\tilde B_{\ell,I}}(\tilde B_q-\tilde B_{\ell,I})}\Theta_I^p$$
We can argue exactly as in the corresponding $T_4$ in the proof of Lemma \ref{Rgtr}, now using the bounds in Proposition \ref{recap}. In the notation of Lemma \ref{Rgtr}, we have:
$$S=(\partial_t+\mathcal{L}_{\tilde v_{\ell,I}})(\tilde v_q-\tilde v_{\ell,I})-\mathcal{L}_{\tilde B_{\ell,I}}(\tilde B_q-\tilde B_{\ell,I})=S_1+S_2+S_3,$$
all the bounds stay the same upon shifting the upper bound for the good estimates $\underline{r}$, with the exception of:
\begin{equation*}
    \begin{split}
        &||\mathcal{\tilde A}_{\ell,I}^\pm(\tilde R_q-\tilde R_\ell)||_{r}\lesssim \lambda_q^{r}(\ell\lambda_q)^{m_0}\ell^{-\alpha}(1/\tau^a)\delta_{q+1}  \ \text{ for } \  0\leq r \leq M-2m_0-k_0^g-7,\\
        &||\mathcal{\tilde A}_{\ell,I}^\pm(\tilde p_q-\tilde p_\ell)||_{r}\lesssim \lambda_q^{r}(\ell\lambda_q)^{m_0}\ell^{-\alpha}(1/\tau^a)\delta_{q+1}  \ \text{ for } \  1\leq r \leq M-2m_0-6,\\
    \end{split}
\end{equation*}
this changes the bound on $\mathcal{\tilde A}_{\ell,I}^\pm S_2$. 

\noindent Corresponding to \eqref{T4g} in the proof of Lemma \ref{Rgtr}, we now set $\bar L= \delta_{q}^{-1/2}, \ \bar{\bar L}=\delta_q^{-1}, \ \underline{r}=M-2m_0-k_0^g-9$ and define the loss functions:
\begin{equation}\label{lossupdate}
    \begin{split}
    L(r)&=1_{r\leq \underline{r}}+1_{r\geq \underline{r}+1}\lambda_q^{[r-\underline{r}](b-1)\gamma_{\ell}}\bar L,\\
    L_{\mathcal{A}}(r)&= 1_{r\leq \underline{r}-1}+1_{r\geq \underline{r}}\lambda_q^{[r-(\underline{r}-1)](b-1)\gamma_{\ell}}\bar{\bar L},\\
    \end{split}
\end{equation} 
note the additional $-1$ shift in $\underline{r}$ compared to the estimates on the slow coefficients and charts in Lemmas \eqref{chartprop}, \eqref{slowcoeffestimates}. With the changes in the estimates mentioned above in mind, we deduce, as in Lemma \ref{Rgtr}, the bounds:
\begin{equation}\label{ST4nash}
    \begin{split}
        ||S||_{r}& \lesssim \lambda_q^{r+1}L(r-1)(\ell\lambda_q)^{m_0}\ell^{-\alpha}\delta_q \ \text{ for } \ 0\leq r\leq N-m_0-1,\\
        ||\mathcal{\tilde A}^\pm_{\ell,I}S||_{r}&\lesssim \lambda_q^{r+1}L_{\mathcal{A}}(r-1)\ell^{-\alpha}(\ell\lambda_q)^{m_0}(1/\tau^a)\delta_{q+1} \ \text{ for } \ 0\leq r\leq N-m_0-2.\\
    \end{split}
\end{equation}

\noindent Intending to apply the Inductive Lemma \ref{inductive}, we compute explicitly the leading term of $T_4$. We set the notation to match that of Lemma \ref{inductive} and provide estimates on the various objects.
\begin{equation*}
    \begin{split}
        \mathcal{L}_S\Theta_I^p&=S\cn \Theta_I^p+ \DD S^\top [\Theta_I^p]\\
        &=S\cn [\frac{\tau^a}{\lambda_{q+1}} \tilde a_I\alpha_I \tilde \Psi^{1*}_I(\varphi_{\lambda,k_I}\nu_I)]+\DD S^\top\left[\frac{\tau^a}{\lambda_{q+1}} \tilde a_I\alpha_I \tilde \Psi^{1*}_I(\varphi_{\lambda,k_I}\nu_I)\right]\\
        &=\underbrace{\tau^a (S\cdot \Psi^{1*}_I k_I) \tilde a_I\alpha_I}_{a_1} \tilde \Psi^{1*}_I(\varphi_{\lambda,k_I}'\nu_I)\\
        &+\Psi^{1*}_I(\varphi_{\lambda,k_I})\underbrace{\frac{\tau^a}{\lambda_{q+1}}  \left[S\cn(\tilde a_I\alpha_I\tilde \Psi^{1*}_I\nu_I)+ \DD S^\top\left[\tilde a_I\alpha_I\tilde \Psi^{1*}_I\nu_I\right]\right]}_{=\frac{\tau^a}{\lambda_{q+1}}\mathcal{L}_S(\tilde a_I\alpha_I\tilde \Psi^{1*}_I\nu_I)=:Y}\\
        &=: a_1\Psi^{1*}_I(\phi_1\nu_I)+\Psi^{1*}_I(\phi_3)Y=F_0
    \end{split}
\end{equation*}
The estimates in \eqref{ST4nash} together with the transport properties and bounds in Lemmas \ref{chartprop}, \ref{slowcoeffestimates} give:
\begin{equation*}
    \begin{split}
        ||a_1||_r&\lesssim \tau^a \supp_t|\alpha_I| \ ||(S\cdot \Psi^{1*}_I k_I) \tilde a_I||_r\lesssim \lambda_q^{r}L(r)(\ell\lambda_q)^{m_0}(\tau^a/\tau^c)\delta_q^{1/2}\delta_{q+1}^{1/2},
        \\
        ||Y||_r&\lesssim \frac{\tau^a}{\lambda_{q+1}}\sup_t|\alpha_I|\ || \mathcal{L}_{S}(\tilde a_I\tilde \Psi^{1*}_I\nu_I)||_r\lesssim \frac{\lambda_q}{\lambda_{q+1}}\lambda_q^{r}L(r)(\ell\lambda_q)^{m_0}(\tau^a/\tau^c)\delta_q^{1/2}\delta_{q+1}^{1/2},
        \\
        ||\mathcal{\tilde A}^\pm_{\ell,I} a_1||_r&\lesssim \sup_t|\alpha_I'| \ ||(S\cdot \Psi^{1*}_I k_I) \tilde a_I||_r+\sup_t|\alpha_I| \ \tau^a||(S\cdot \Psi^{1*}_I k_I)\mathcal{\tilde A}^\pm_{\ell,I}( \tilde a_I)||_r\\
        &+\tau^a \sup_t|\alpha_I| \ ||(\mathcal{\tilde A}^\pm_{\ell,I}S\cdot \Psi^{1*}_I k_I) \tilde a_I\alpha_I||_r+\tau^a \sup_t|\alpha_I| \ ||(S\cdot \mathcal{\tilde A}^\pm_{\ell,I}\Psi^{1*}_I k_I) \tilde a_I||_r\\
        &\lesssim\lambda_q^{r}L_{\mathcal{A}}(r)(\ell\lambda_q)^{m_0}(1/\tau^c)\delta_q^{1/2}\delta_{q+1}^{1/2},
        \\
        ||(\partial_t+\mathcal{L}_{\tilde z^\pm_{\ell,I}}) Y||_r&\lesssim\frac{1}{\lambda_{q+1}}\sup_t|\alpha_I'| \ || \mathcal{L}_{S}(\tilde a_I\tilde \Psi^{1*}_I\nu_I)||_r+\frac{\tau^a}{\lambda_{q+1}}\sup_t|\alpha_I| \ || \mathcal{L}_{S}(\mathcal{\tilde A}^\pm_{\ell,I}( \tilde a_I)\tilde \Psi^{1*}_I\nu_I)||_r\\
        &+\frac{\tau^a}{\lambda_{q+1}}\sup_t|\alpha_I|\ || \mathcal{L}_{(\partial_t+\mathcal{L}_{\tilde z^\pm_{\ell,I}})S}(\tilde a_I\tilde \Psi^{1*}_I\nu_I)||_r\\
        &\lesssim \frac{\lambda_q}{\lambda_{q+1}}\lambda_q^{r}L_{\mathcal{A}}(r)(\ell\lambda_q)^{m_0}(1/\tau^c)\delta_q^{1/2}\delta_{q+1}^{1/2}+\frac{1}{\lambda_{q+1}}\lambda_q^{r+2}L_{\mathcal{A}}(r)\ell^{-\alpha}(\ell\lambda_q)^{m_0}\delta_{q+1}^{3/2}\\
        &\lesssim \frac{\lambda_q}{\lambda_{q+1}}\lambda_q^{r}L_{\mathcal{A}}(r)(\ell\lambda_q)^{m_0}(1/\tau^c)\delta_q^{1/2}\delta_{q+1}^{1/2}.
    \end{split}
\end{equation*}
for $0\leq r\leq N-m_0-2$, where we used \eqref{trick} to commute the lie derivatives and compute:
\begin{equation*}
    \begin{split}
        (\partial_t+\mathcal{L}_{\tilde z^\pm_{\ell,I}})Y&=\frac{\tau^a}{\lambda_{q+1}}(\partial_t+\mathcal{L}_{\tilde z^\pm_{\ell,I}})\left[\mathcal{L}_{S}(\tilde a_I\alpha_I\tilde \Psi^{1*}_I\nu_I)\right]\\
        &=\frac{\tau^a}{\lambda_{q+1}}\alpha_I(\partial_t+\mathcal{L}_{\tilde z^\pm_{\ell,I}})\left[\mathcal{L}_{S}(\tilde a_I\tilde \Psi^{1*}_I\nu_I)\right]+\frac{1}{\lambda_{q+1}}\alpha_I'\mathcal{L}_{S}(\tilde a_I\tilde \Psi^{1*}_I\nu_I)\\
        &=\frac{\tau^a}{\lambda_{q+1}}\alpha_I\mathcal{L}_{(\partial_t+\mathcal{L}_{\tilde z^\pm_{\ell,I}})S}(\tilde a_I\tilde \Psi^{1*}_I\nu_I)+\frac{\tau^a}{\lambda_{q+1}}\alpha_I\mathcal{L}_{S}\Biggr[\underbrace{(\partial_t+\mathcal{L}_{\tilde z^\pm_{\ell,I}})\tilde \Psi^{1*}_I\nu_I}_{=0}+\mathcal{\tilde A}^\pm_{\ell,I}( \tilde a_I)\tilde \Psi^{1*}_I\nu_I\Biggr]\\
        &+\frac{1}{\lambda_{q+1}}\alpha_I'\mathcal{L}_{S}(\tilde a_I\tilde \Psi^{1*}_I\nu_I),
    \end{split}
\end{equation*}
note that the difference in the transport estimate for $S$ compared to the one in Lemma \ref{Rgtr} gives in fact a lower order term in $(\partial_t+\mathcal{L}_{\tilde z^\pm_{\ell,I}})Y$. We now follow the same strategy as in $T_1$ and apply the Inductive Lemma \ref{inductive} with: 
\begin{equation*}
    \begin{split}
        &z^\pm=\tilde z^\pm_{\ell,I},\\ 
        &\sigma_i=\sigma=\xi_{I},\\
        &\bar\varsigma_i=\bar\varsigma=\tau^a\delta_{q+1}^{1/2},\ \ \bar\varsigma_{i,\mathcal{A}}=\bar\varsigma_{\mathcal{A}}=\delta_{q+1}^{1/2},\\
        &A_1=(\ell\lambda_q)^{m_0}(\tau^a/\tau^c)\delta_q^{1/2}\delta_{q+1}^{1/2}, \ \ A_{1,\mathcal{A}}=(\ell\lambda_q)^{m_0}(1/\tau^c)\delta_q^{1/2}\delta_{q+1}^{1/2},\\
        &A_3=\frac{\lambda_q}{\lambda_{q+1}}A_1, \ \ A_{3,\mathcal{A}}=\frac{\lambda_q}{\lambda_{q+1}}A_{1,\mathcal{A}},\\
        &L^1(r)=L(r), \ \ L^1_{\mathcal{A}}(r)=L_{\mathcal{A}}(r), \ \ L^\varsigma(r)=\lambda_q^{[r-\underline{r}]^+(b-1)\gamma_\ell},\\
    \end{split}
\end{equation*} 
this time, compared to $T_1$, we are forced to set:
$$\bar N=N-m_0-2$$
and finally we let $\bar C$ be the largest implicit constant in the above estimates for $0\leq r \leq \bar N$.
The inductive Lemma gives decompositions: 
\begin{equation*} 
    \begin{split}
        \mathcal{L}_{\xi^{p}_I}^k\mathcal{L}_S\Theta_I^{p}&\simeq a_{1,k}\tilde \Psi^{1*}_I(\phi_{1,k}\nu_I)+a_{2,k} \tilde \Psi^{1*}_I(\phi_{2,k}k_I)+\tilde \Psi^{*}_I(\phi_{3,k}) Y_k,
        \\
        (\partial_t+\mathcal{L}_{\tilde z^\pm_{\ell,I}})\mathcal{L}_{\xi^{p}_I}^k\mathcal{L}_S\Theta_I^{p}&\simeq a_{1,k,\mathcal{A}}\tilde \Psi^{1*}_I(\phi_{1,k}\nu_I)+a_{2,k,\mathcal{A}} \tilde \Psi^{1*}_I(\phi_{2,k}k_I)+\tilde \Psi^{*}_I(\phi_{3,k}) Y_{k,\mathcal{A}}\\
    \end{split}
\end{equation*}
where, abusing notation, we use the same letter for the fast coefficients and estimates:
\begin{equation*}
    \begin{split}
    ||a_{1,k}||_r&\leq C'(C)^k\lambda_{q}^{r+k}L^1(r+k)A_1\prod_{i=1}^k\bar \varsigma_i\lesssim (C)^k\lambda_q^{r}L(k+r)\mathcal{T}_p^k(\ell\lambda_q)^{m_0}(\tau^a/\tau^c)\delta_q^{1/2}\delta_{q+1}^{1/2},
    \\
        ||a_{2,k}||_r&\leq C'(C)^k \lambda_{q+1}\lambda_{q}^{r+k-1}L(r+k-1)A_3\prod_{i=1}^k\bar \varsigma_i\lesssim (C)^k\lambda_q^{r}L(r+k-1)\mathcal{T}_p^k(\ell\lambda_q)^{m_0}(\tau^a/\tau^c)\delta_q^{1/2}\delta_{q+1}^{1/2},
        \\
        ||Y_k||_r&\leq C'(C)^k  \lambda_{q}^{r+k}L(r+k)A_3\prod_{i=1}^k\bar \varsigma_i\lesssim(C)^k\frac{\lambda_q}{\lambda_{q+1}}\lambda_q^{r}L(r+k)\mathcal{T}_p^k(\ell\lambda_q)^{m_0}(\tau^a/\tau^c)\delta_q^{1/2}\delta_{q+1}^{1/2}
    \end{split}
\end{equation*}
for $k\geq 1, \ 0\leq r+k\leq N-m_0-2$ together with: 
\begin{equation*}
    \begin{split}
        ||a_{1,k,\mathcal{A}}||_r&\leq C'(C)^k\lambda_{q}^{r+k}\left[L_{\mathcal{A}}^1(k+r)A_{1,\mathcal{A}}\prod_{i=1}^k\bar\varsigma_i +L^1(k+r)A_1\max_{j\in\{1,\dots,k\}}\bar\varsigma_{j,\mathcal{A}}\prod_{i=1,\dots,k \ i\neq j}\bar\varsigma_i\right]\\
        &\lesssim (C)^k\lambda_q^{r}L_{\mathcal{A}}(r+k)\mathcal{T}_p^k(\ell\lambda_q)^{m_0}(1/\tau^c)\delta_q^{1/2}\delta_{q+1}^{1/2},
        \\
        ||a_{2,k,\mathcal{A}}||_r&\leq C'(C)^k \lambda_{q+1}\lambda_{q}^{r+k-1}\left[L_{\mathcal{A}}(r+k-1)A_{3,\mathcal{A}}\prod_{i=1}^k\bar\varsigma_i +L(r+k-1) A_3\max_{j\in\{1,\dots,k\}}\bar\varsigma_{j,\mathcal{A}}\prod_{i=1,\dots,k \ i\neq j}\bar\varsigma_i\right]\\
        &\lesssim(C)^k\lambda_q^{r}L_{\mathcal{A}}(r+k-1)\mathcal{T}_p^k(\ell\lambda_q)^{m_0}(1/\tau^c)\delta_q^{1/2}\delta_{q+1}^{1/2},
        \\
        ||Y_{k,\mathcal{A}}||_r&\leq C'(C)^k\lambda_{q}^{r+k}\left[L_{\mathcal{A}}(k+r)A_{3,\mathcal{A}}\prod_{i=1}^k\bar\varsigma_i +L(k+r)A_3\max_{j\in\{1,\dots,k\}}\bar\varsigma_{j,\mathcal{A}}\prod_{i=1,\dots,k \ i\neq j}\bar\varsigma_i\right]\\
        &\lesssim(C)^k\frac{\lambda_q}{\lambda_{q+1}}\lambda_q^{r}L_{\mathcal{A}}(r+k)\mathcal{T}_p^k(\ell\lambda_q)^{m_0}(1/\tau^c)\delta_q^{1/2}\delta_{q+1}^{1/2}
    \end{split}
\end{equation*}
for $k\geq 1, \ 0\leq r+k\leq N-m_0-3$. Note that here the loss functions are not admissible in the sense of \ref{admissible}; this is why there is a $-1$ shift in \eqref{lossupdate}, making the estimate cleaner. 

\noindent The inductive Lemma finally gives: 
\begin{equation*}
    ||\Psi^{1*}_I(\phi_{1,k}\nu_I)||_r, \ ||\Psi^{1*}_I(\phi_{2,k}k_I)||_r, \ ||\Psi^{1*}_I(\phi_{3,k})||_r \lesssim \lambda_{q+1}^r.
\end{equation*}
For later use, we add a second parameter to the loss functions in \eqref{lossupdate}, and shift them by one to avoid troubles when interpolating the $C^{r+\alpha}$ norms. Namely, for $r,k\geq 0$ we set:
\begin{equation}\label{lossupdate2}
    \begin{split}
    L(r,k)&=L(k+1)+\left(\frac{\lambda_q}{\lambda_{q+1}}\right)^rL(r+k+1),\\
    L_{\mathcal{A}}(r,k)&= L_{\mathcal{A}}(k+1)+\left(\frac{\lambda_q}{\lambda_{q+1}}\right)^rL_{\mathcal{A}}(r+k+1).\\
    \end{split}
\end{equation}
From the bounds above, the fact that the worst estimates are associated with $a_{1,k}$ and interpolation, we conclude:
\begin{equation}\label{T4nashk}
    \begin{split}
        ||\mathcal{L}_{\xi^{p}_I}^k\mathcal{L}_S\Theta_I^{p}||_{r+\alpha}&\lesssim ||a_{1,k}||_{r+\alpha}||\tilde \Psi^{1*}_I(\phi_{1,k}\nu)||_0+||a_{1,k}||_0||\tilde \Psi^{1*}_I(\phi_{1,k}\nu)||_{r+\alpha}\\
        &+||a_{2,k}||_{r+\alpha}||\tilde \Psi^{1*}_I(\phi_{2,k}k)||_0+||a_{2,k}||_0||\tilde \Psi^{1*}_I(\phi_{2,k}k)||_{r+\alpha}\\
        &+||\tilde \Psi^{*}_I(\phi_{3,k})||_{r+\alpha}||Y_k||_0+||\tilde \Psi^{*}_I(\phi_{3,k})||_0||Y_k||_{r+\alpha}\\
        &\lesssim (C)^k \lambda_q^{r+\alpha}L(k+r+1)\mathcal{T}_p^k(\ell\lambda_q)^{m_0}(\tau^a/\tau^c)\delta_q^{1/2}\delta_{q+1}^{1/2}\\
        &+(C)^k\lambda_{q+1}^{r+\alpha}L(k+1)\mathcal{T}_p^k(\ell\lambda_q)^{m_0}(\tau^a/\tau^c)\delta_q^{1/2}\delta_{q+1}^{1/2}\\
        &\leq (C)^k \lambda_{q+1}^{r+\alpha}L(r,k)\mathcal{T}_p^k(\ell\lambda_q)^{m_0}(\tau^a/\tau^c)\delta_q^{1/2}\delta_{q+1}^{1/2}
    \end{split}
\end{equation}
for $0\leq r+k\leq N-m_0-3$. The constants depend on all the parameters, but not on $a$ and are uniform in $r, k, I$. Note that we lose one additional derivative to interpolate the $C^{r+\alpha}$ norm.

\noindent From the bounds above, and the fact that the worst estimates are associated with $\mathcal{\tilde A}_{\ell,I}^\pm a_{1,k}$, we first deduce:
\begin{equation*}
    \begin{split}
        ||(\partial_t+\mathcal{L}_{\tilde z^\pm_{\ell,I}})\mathcal{L}_{\xi^{p}_I}^k\mathcal{L}_S\Theta_I^{p}||_{r+\alpha}&\lesssim ||a_{1,k,\mathcal{A}}||_{r+\alpha}||\tilde \Psi^{1*}_I(\phi_{1,k}\nu_I)||_0+||a_{1,k,\mathcal{A}}||_0||\tilde \Psi^{1*}_I(\phi_{1,k}\nu_I)||_{r+\alpha}\\
        &+||a_{2,k,\mathcal{A}}||_{r+\alpha}||\tilde \Psi^{1*}_I(\phi_{2,k}k_I)||_0+||a_{2,k,\mathcal{A}}||_0||\tilde \Psi^{1*}_I(\phi_{2,k}k_I)||_{r+\alpha}\\
        &+||\tilde \Psi^{*}_I(\phi_{3,k})||_{r+\alpha}||Y_{k,\mathcal{A}}||_0+||\tilde \Psi^{*}_I(\phi_{3,k})||_0||Y_{k,\mathcal{A}}||_{r+\alpha}\\
        &\lesssim (C)^k\lambda_q^{r+\alpha}L_{\mathcal{A}}(r+k+1)\mathcal{T}_p^k(\ell\lambda_q)^{m_0}1/\tau^c\delta_q^{1/2}\delta_{q+1}^{1/2}\\
        &+(C)^k\lambda_{q+1}^{r+\alpha}L_{\mathcal{A}}(k+1)\mathcal{T}_p^k(\ell\lambda_q)^{m_0}1/\tau^c\delta_q^{1/2}\delta_{q+1}^{1/2}\\
        &\leq  (C)^k \lambda_{q+1}^{r+\alpha}L_{\mathcal{A}}(r,k)\mathcal{T}_p^k(\ell\lambda_q)^{m_0}1/\tau^c\delta_q^{1/2}\delta_{q+1}^{1/2},
    \end{split}
\end{equation*}
and conclude from the estimates in Lemma\ref{localcorrn} and \eqref{T4nashk} that:
\begin{equation*}
    \begin{split}
        ||\mathcal{\tilde A}^\pm_{\ell,I}\mathcal{L}_{\xi^{p}_I}^k\mathcal{L}_S\Theta_I^{p}||_{r+\alpha}&\leq ||(\partial_t+\mathcal{L}_{\tilde z^\pm_{\ell,I}})\mathcal{L}_{\xi^{p}_I}^k\mathcal{L}_S\Theta_I^{p}||_{r+\alpha}+||\mathcal{L}_{\xi^{p}_I}^k\mathcal{L}_S\Theta_I^{p}\cn z^\pm_{\ell,I}||_{r+\alpha}\\
        &\lesssim (C)^k\lambda_{q+1}^{r+\alpha}L_{\mathcal{A}}(r,k)\mathcal{T}_p^k(\ell\lambda_q)^{m_0}1/\tau^c\delta_q^{1/2}\delta_{q+1}^{1/2}+(C)^k\lambda_{q+1}^{r+\alpha}L(r,k)\mathcal{T}_p^k(\ell\lambda_q)^{m_0}(\tau^a/\tau^c)\lambda_q\delta_q\delta_{q+1}^{1/2}\\
        &\lesssim (C)^k\lambda_{q+1}^{r+\alpha}L_{\mathcal{A}}(r,k)\mathcal{T}_p^k(\ell\lambda_q)^{m_0}1/\tau^c\delta_q^{1/2}\delta_{q+1}^{1/2}\\
    \end{split}
\end{equation*}
for $0\leq r+k\leq N-m_0-4$, where we used $L\leq L_{\mathcal{A}}$. The same remarks about the constants and interpolation hold.

\noindent \textit{Conclusion.} Proceeding as in $T_1$, we conclude that:
\begin{equation*}
    \begin{split}
        ||T_4||_{r+\alpha}&\leq \sum_{k=0}^{k_0^p}\frac{1}{(k+1)!} ||\mathcal{R}\curl\sum_I \mathcal{L}_{\xi^{p}_I}^k[\mathcal{L}_{S}\Theta_I^{p}||_{r+\alpha}\lesssim  \sum_{k=0}^{k_0^p}\frac{1}{(k+1)!} ||\sum_I\mathcal{L}_{\xi^{p}_I}^k[\mathcal{L}_{S}\Theta_I^{p}||_{r+\alpha}\\
        &\lesssim \sup_I\sum_{k=0}^{k_0^p}\frac{1}{(k+1)!} ||\mathcal{L}_{\xi^{p}_I}^k[\mathcal{L}_{S}\Theta_I^{p}||_{r+\alpha}\\
        &\lesssim\sum_{k=0}^{k_0^p}\frac{(C)^k}{(k+1)!} \lambda_{q+1}^{r+\alpha}L(r,k)\mathcal{T}_p^k(\ell\lambda_q)^{m_0}(\tau^a/\tau^c)\delta_q^{1/2}\delta_{q+1}^{1/2}\\
        &\leq \lambda_{q+1}^{r}(\tau^a/\tau^c)\frac{\lambda_q\delta_q^{1/2}\delta_{q+1}^{1/2}}{\lambda_{q+1}^{1-\alpha}}\Biggr[\underbrace{\sum_{k=0}^{k_0^p}\frac{\left(C\mathcal{T}_p\right)^k}{(k+1)!}L(k+1)}_{\mathcal{S}_1}+\underbrace{\sum_{k=0}^{k_0^p}\frac{\left(C\mathcal{T}_p\right)^k}{(k+1)!}\left(\frac{\lambda_q}{\lambda_{q+1}}\right)^{r}L(r+k+1)}_{\mathcal{S}_2}\Biggr]\\
    \end{split}
\end{equation*}
for $0\leq r\leq N-k_0^p-m_0-3$, where the implicit constant depends on $r$ and all the other parameters, but not on $a$, and we used $(\ell\lambda_q)^{m_0}\leq \lambda_q/\lambda_{q+1}$. 

\noindent Computing as in $T_1$ to correct the transport operator and deal with the commutators, see \eqref{transportT1n}, we get:
\begin{equation*}
    \begin{split}
        &||\mathcal{\tilde A}^\pm T_4||_{r+\alpha}\\
        &\leq \sum_{k=0}^{k_0^p} \frac{1}{(k+1)!}\left[||\mathcal{R}\curl \sum_I\mathcal{\tilde A}^\pm_{\ell,I}\mathcal{L}_{\xi^p_I}^k\mathcal{L}_{S}\Theta_I^{p}||_{r+\alpha}+||[\mathcal{\tilde A}^\pm,\mathcal{R}\curl]\sum_I \mathcal{L}_{\xi^{p}_I}^k\mathcal{L}_{S}\Theta_I^{p}||_{r+\alpha}\right]\\
        &+ \sum_{k=0}^{k_0^p} \frac{1}{(k+1)!}\left[||\mathcal{R}\curl\sum_I (\tilde z_q^\pm-\tilde z^\pm_{\ell,I})\cn\mathcal{L}_{\xi^{p}_I}^k\mathcal{L}_{S}\Theta_I^{p}||_{r+\alpha}\right]\\
        &\lesssim \sup_I\sum_{k=0}^{k_0^p}\frac{1}{(k+1)!} \left[||\tilde z^\pm_q||_{r+1+\alpha}||\mathcal{L}_{\xi^{p}_I}^k\mathcal{L}_{S}\Theta_I^{p}||_{\alpha}+||\tilde z^\pm_q||_{1+\alpha}||\mathcal{L}_{\xi^{p}_I}^k\mathcal{L}_{S}\Theta_I^{p}||_{r+\alpha}\right]\\
        &+ \sup_I\sum_{k=0}^{k_0^p}\frac{1}{(k+1)!}\left[||\tilde z_q^\pm-\tilde z^\pm_{\ell,I}||_{r+\alpha}||\mathcal{L}_{\xi^{p}_I}^k\mathcal{L}_{S}\Theta_I^{p}||_{1}+||\tilde z_q^\pm-\tilde z^\pm_{\ell,I}||_{0}|| \mathcal{L}_{\xi^{p}_I}^k\mathcal{L}_{S}\Theta_I^{p}||_{r+1+\alpha}\right]\\
        &+ \sup_I\sum_{k=0}^{k_0^p}\frac{1}{(k+1)!}|| \mathcal{\tilde A}^\pm_{\ell,I}\mathcal{L}_{\xi^{p}_I}^k\mathcal{L}_{S}\Theta_I^{p}||_{r+\alpha}\\
        &\lesssim \sum_{k=0}^{k_0^p}\frac{1}{(k+1)!}\left[\lambda_{q+1}^{r+\alpha}(C)^kL(r,k)\mathcal{T}_p^k\lambda_q^{\alpha}(\ell\lambda_q)^{m_0}(\tau^a/\tau^c)\lambda_q\delta_q\delta_{q+1}^{1/2}\right]\\
        &+\sum_{k=0}^{k_0^p}\frac{1}{(k+1)!}\left[\lambda_{q+1}^{r+1+\alpha}(C)^kL(r,k)\mathcal{T}_p^k\lambda_q^{\alpha}(\ell\lambda_q)^{2m_0}(\tau^a/\tau^c)\delta_q\delta_{q+1}^{1/2}\right]\\
        &+\sum_{k=0}^{k_0^p}\frac{1}{(k+1)!}\left[\lambda_{q+1}^{r+\alpha}(C)^kL_{\mathcal{A}}(r,k)\mathcal{T}_p^k(\ell\lambda_q)^{m_0}1/\tau^c\delta_q^{1/2}\delta_{q+1}^{1/2}\right]\\
    \end{split}
\end{equation*}
for $0\leq r\leq N-m_0-k_0^p-4$, where the implicit constant depends on $r$  and all the other parameters, but not on $a$. From $(\ell\lambda_q)^{m_0}\leq \frac{\lambda_q}{\lambda_{q+1}}$ and $L\leq L_{\mathcal{A}}$ and their definition in \eqref{lossupdate2} we obtain:
\begin{equation*}
    \begin{split}
        ||\mathcal{\tilde A}^\pm T_4||_{r+\alpha}&\lesssim \sum_{k=0}^{k_0^p}\frac{1}{(k+1)!}\left[\lambda_{q+1}^{r+2\alpha}(C)^kL_{\mathcal{A}}(r,k)\mathcal{T}_p^k(\ell\lambda_q)^{m_0}1/\tau^c\delta_q^{1/2}\delta_{q+1}^{1/2}\right]\\
        &= \sum_{k=0}^{k_0^p}\frac{1}{(k+1)!}(C)^k\lambda_{q+1}^{r+2\alpha}\left[L_{\mathcal{A}}(k+1)+\left(\frac{\lambda_q}{\lambda_{q+1}}\right)^{r}L_{\mathcal{A}}(r+k+1)\right]\mathcal{T}_p^k(\ell\lambda_q)^{m_0}1/\tau^c\delta_q^{1/2}\delta_{q+1}^{1/2}\\
        &\leq \lambda_{q+1}^r(1/\tau^c)\frac{\lambda_q\delta_q^{1/2}\delta_{q+1}^{1/2}}{\lambda_{q+1}^{1-2\alpha}}\Biggr[\underbrace{\sum_{k=0}^{k_0^p}\frac{\left(C\mathcal{T}_p\right)^k}{(k+1)!}L_{\mathcal{A}}(k+1)}_{\mathcal{S}_{1, \mathcal{A}}}+\underbrace{\sum_{k=0}^{k_0^p}\frac{\left(C\mathcal{T}_p\right)^k}{(k+1)!}\left(\frac{\lambda_q}{\lambda_{q+1}}\right)^{r}L_{\mathcal{A}}(r+k+1)}_{\mathcal{S}_{2,\mathcal{A}}}\Biggr].\\
    \end{split}
\end{equation*}

\noindent We now bound the sums $\mathcal{S}_1, \ \mathcal{S}_2, \ \mathcal{S}_{1,\mathcal{A}}, \ \mathcal{S}_{2,\mathcal{A}}$. Since the ideas are the same, and in terms of the parameters, the constraints are more stringent, we will show the explicit calculations only for the $\mathcal{A}$ versions. 

\noindent We begin with $\mathcal{S}_{2,\mathcal{A}}$ which we split according to the number of derivatives $r$ we are considering and the definition of $L_{\mathcal{A}}$ in \eqref{lossupdate}, we have:
\begin{equation*}
    \begin{split}
        \mathcal{S}_{2,\mathcal{A}}&=\sum_{k=0}^{k_0^p}\frac{\left(C\mathcal{T}_p\right)^k}{(k+1)!}\left(\frac{\lambda_q}{\lambda_{q+1}}\right)^{r}L_{\mathcal{A}}(r+k+1)\\
        &= \sum_{0\leq k+r\leq \underline{r}-2, \ 0\leq k\leq k_0^p}\frac{\left(C\mathcal{T}_p\right)^k}{(k+1)!}\left(\frac{\lambda_q}{\lambda_{q+1}}\right)^{r}+\left(\frac{\lambda_q}{\lambda_{q+1}}\right)^{r-[r-(\underline{r}-2)]^+\gamma_\ell}\bar{\bar L}\sum_{\underline{r}-1\leq k+r, \ 0\leq k\leq k_0^p }\frac{\left(C\lambda_q^{(b-1)\gamma_\ell}\mathcal{T}_p\right)^k}{(k+1)!}\\
        &= \sum_{k=0}^{[\underline{r}-1-r]^+-1}\frac{\left(C\mathcal{T}_p\right)^k}{(k+1)!}\left(\frac{\lambda_q}{\lambda_{q+1}}\right)^{r}+\left(\frac{\lambda_q}{\lambda_{q+1}}\right)^{r-[r-(\underline{r}-2)]^+\gamma_\ell}\bar{\bar L}\sum_{k=[\underline{r}-1-r]^+ }^{k_0^p}\frac{\left(C\lambda_q^{(b-1)\gamma_\ell}\mathcal{T}_p\right)^k}{(k+1)!}\\
        &\lesssim 1+\left(\frac{\lambda_q}{\lambda_{q+1}}\right)^{r-[r-(\underline{r}-2)]^+\gamma_\ell}\bar{\bar L}\left(\lambda_q^{(b-1)\gamma_\ell}\mathcal{T}_p\right)^{[\underline{r}-1-r]^+}\\
        &\overset{0\leq r\leq \underline{r}-2}{=} 1+\left(\frac{\lambda_q}{\lambda_{q+1}}\right)^{r}\bar{\bar L}\left(\lambda_q^{(b-1)\gamma_\ell}\mathcal{T}_p\right)^{\underline{r}-1-r}=1+\left(\frac{\lambda_q}{\lambda_{q+1}}\right)^{r\underbrace{[1-(\gamma_a+\beta+\gamma_{CZ})+\gamma_\ell]}_{\geq 0}}\underbrace{\bar{\bar L}\left(\lambda_q^{(b-1)\gamma_\ell}\mathcal{T}_p\right)^{\underline{r}-1}}_{\leq 1}\lesssim 1\\
        &\overset{r\geq  \underline{r}-1}{=}1+\left(\frac{\lambda_q}{\lambda_{q+1}}\right)^{r-[r-(\underline{r}-2)]\gamma_\ell}\bar{\bar L}\leq 1+\left(\frac{\lambda_q}{\lambda_{q+1}}\right)^{r(1-\gamma_\ell)}\bar{\bar L}\leq 1+\left(\lambda_q^{(b-1)\gamma_{\ell}}\mathcal{T}_p\right)^{\underline{r}-1}\bar{\bar L}\lesssim 1\\
    \end{split}
\end{equation*}
where we used, see \eqref{constraintadmissibility} and \eqref{M}, that:
$$1-(\gamma_a+\gamma_\ell+2\beta+2\gamma_{CZ})\geq 0,\ \ \ \ \underline{r}-1\geq k_0^g, \ \ \ \ \left(\lambda_q^{(b-1)\gamma_{\ell}}\mathcal{T}_p\right)^{k_0^g}\bar{\bar L}\leq 1 $$
from which it also follows, given the definition of $\mathcal{T}_p$ in \eqref{Tp}, that:
$$\mathcal{T}_p=\lambda_q\tau^a\delta_{q+1}^{1/2}=\left(\frac{\lambda_q}{\lambda_{q+1}}\right)^{\gamma_a+\beta+\gamma_{CZ}}\geq \frac{\lambda_q}{\lambda_{q+1}}.$$
Similarly, we can bound $\mathcal{S}_{1,\mathcal{A}}$:
\begin{equation*}
    \begin{split}
        \mathcal{S}_{1,\mathcal{A}}&=\sum_{k=0}^{k_0^p}\frac{\left(C\mathcal{T}_p\right)^k}{(k+1)!}L_{\mathcal{A}}(k+1)= \sum_{k=0}^{\underline{r}-2}\frac{\left(C\mathcal{T}_p\right)^k}{(k+1)!}+\bar{\bar L}\sum_{k=\underline{r}-1}^{k_0^p}\frac{\left(C\lambda_q^{(b-1)\gamma_\ell}\mathcal{T}_p\right)^k}{(k+1)!}\\
        &\lesssim 1+ \left(\lambda_q^{(b-1)\gamma_\ell}\mathcal{T}_p\right)^{\underline{r}-1}\bar{\bar L}\leq 1+ \left(\lambda_q^{(b-1)\gamma_\ell}\mathcal{T}_p\right)^{k_0^g}\bar{\bar L}\lesssim 1
    \end{split}
\end{equation*}
and we conclude that: 
$$\mathcal{S}_1+\mathcal{S}_2\lesssim 1, \  \mathcal{S}_{1, \mathcal{A}}+\mathcal{S}_{2,\mathcal{A}}\lesssim 1.$$ 

\noindent Going back to the original estimates, we obtain:
\begin{equation*}
    \begin{split}
        ||T_4||_{r+\alpha}&\lesssim \lambda_{q+1}^{r}(\tau^a/\tau^c)\frac{\lambda_q\delta_q^{1/2}\delta_{q+1}^{1/2}}{\lambda_{q+1}^{1-\alpha}} \ \text{ for } \ 0\leq r\leq N-k_0^p-m_0-3,\\
        ||\mathcal{\tilde A}^\pm T_4||_{r+\alpha}&\lesssim \lambda_{q+1}^{r}(1/\tau^c)\frac{\lambda_q\delta_q^{1/2}\delta_{q+1}^{1/2}}{\lambda_{q+1}^{1-2\alpha}} \ \text{ for } \ 0\leq r\leq N-k_0^p-m_0-4\\
    \end{split}
\end{equation*}
where the implicit constants depend on $r$ and all the parameters, but not on $a$.

\noindent \textbf{Estimates on $T_{2,a},\ T_{2,b}, \ T_3, \ T_{5,a},\ T_{5,b}$.} The proof of the following bounds is just a slight variation of the techniques and ideas contained in the proof of Lemma \ref{estimatesperturbationnash} and the ones above. We refer to that and state only the final bounds. 
\begin{equation*}
    \begin{split}
        &||T_{2,w}||_{r+\alpha}\lesssim \lambda_{q+1}^r\frac{\lambda_q\delta_{q+1}}{\lambda_{q+1}^{1-\alpha}} \ \text{ for } r\geq 0,\\
        &||\mathcal{\tilde A}^\pm T_{2,w}||_{r+\alpha}\lesssim \lambda_{q+1}^r1/\tau^a\frac{\lambda_q\delta_{q+1}}{\lambda_{q+1}^{1-2\alpha}} \ \text{ for } 0\leq r \leq N-m_0-1,\\
        &||T_{2,b}||_{r+\alpha}\lesssim \lambda_{q+1}^r(\tau^a/\tau^c)^2\frac{\lambda_q\delta_{q+1}}{\lambda_{q+1}^{1-\alpha}} \ \text{ for } r\geq 0,\\
        &||\mathcal{\tilde A}^\pm T_{2,b}||_{r+\alpha}\lesssim \lambda_{q+1}^r(1/\tau^a)(\tau^a/\tau^c)^2\frac{\lambda_q\delta_{q+1}}{\lambda_{q+1}^{1-2\alpha}} \ \text{ for } 0\leq r\leq N-m_0-1,\\
        &||T_{3}||_{r+\alpha}, \lesssim \lambda_{q+1}^r\frac{\lambda_q\delta_q^{1/2}\delta_{q+1}^{1/2}}{\lambda_{q+1}^{1-\alpha}} \ \text{ for } \ 0\leq r\leq N-k_0^p-m_0-2,\\
        &||\mathcal{\tilde A}^\pm T_{3}||_{r+\alpha}\lesssim \lambda_{q+1}^r(1/\tau^a)\frac{\lambda_q\delta_q^{1/2}\delta_{q+1}^{1/2}}{\lambda_{q+1}^{1-2\alpha}} \ \text{ for } \ 0\leq r \leq N-k_0^p-m_0-3,\\
        &||T_{5,w}||_{r+\alpha}\lesssim \lambda_{q+1}^r\tau^a/\tau^c\frac{\lambda_q\delta_{q+1}}{\lambda_{q+1}^{1-\alpha}} \ \text{ for } \ 0\leq r\leq N-k_0^p-m_0-2,\\
        &||\mathcal{\tilde A}^\pm T_{5,w}||_{r+\alpha} \lesssim \lambda_{q+1}^r1/\tau^c\frac{\lambda_q\delta_{q+1}}{\lambda_{q+1}^{1-2\alpha}} \ \text{ for } \ 0\leq r\leq N-k_0^p-m_0-3,\\
        &||T_{5,b}||_{r+\alpha}\lesssim \lambda_{q+1}^r(\tau^a/\tau^c)^2\frac{\lambda_q\delta_{q+1}}{\lambda_{q+1}^{1-\alpha}} \ \text{ for }\ 0\leq r\leq N-k_0^p-m_0-2,\\
        &||\mathcal{\tilde A}^\pm T_{5,b}||_{r+\alpha}\lesssim\lambda_{q+1}^r(1/\tau^a)(\tau^a/\tau^c)^2\frac{\lambda_q\delta_{q+1}}{\lambda_{q+1}^{1-2\alpha}} \ \text{ for }\ 0\leq r\leq N-k_0^p-m_0-3.
    \end{split}
\end{equation*}
The implicit constants depend on $r$ and all the other parameters, but not on $a$.

\noindent \textbf{Conclusion.} We gather all the estimates and observe that $T_1$ has the worst, while $T_4$ has the worst range of derivatives we can allow. We conclude that: 
\begin{equation*}
    \begin{split}
        ||R^{p,tr}||_{r+\alpha}&\lesssim \lambda_{q+1}^r(\tau^c/\tau^a)\frac{\lambda_q\delta_q^{1/2}\delta_{q+1}^{1/2}}{\lambda_{q+1}^{1-2\alpha}} \ \text{ for } \ 0\leq r\leq N-k_0^p-m_0-3,\\
        ||\mathcal{\tilde A}^\pm R^{p,tr}||_{r+\alpha}&\lesssim\lambda_{q+1}^r(1/\tau^a)(\tau^c/\tau^a)\frac{\lambda_q\delta_q^{1/2}\delta_{q+1}^{1/2}}{\lambda_{q+1}^{1-2\alpha}} \ \text{ for } \ 0\leq r \leq N-m_0-k_0^p-4\\
    \end{split}
\end{equation*}
where the implicit constants depend on $r$ and all the parameters, but not on $a$. From our choice of parameters in \ref{choiceofparameters}, see \eqref{constraintremainder}, it follows that the bounds hold in particular for $0\leq r\leq M$ and $0\leq r\leq M-1$.
\end{proof}


\subsubsection{Estimates on the Nash, Mollification and Remainder Errors}
\begin{lemma}[Estimates on $R^{p,na}$] \label{Rpna} Under the choice of parameters in \ref{choiceofparameters}, the following bounds hold.
        \begin{equation*}
        \begin{split}
            &||R^{p,na}||_{r+\alpha}\lesssim \lambda_{q+1}^{r}\frac{\lambda_q\delta_q^{1/2}\delta_{q+1}^{1/2}}{\lambda_{q+1}^{1-\alpha}} \ \text{ for } \ 0\leq r\leq N-k_0^p-m_0-2,\\
            &||\mathcal{\tilde A}^\pm R^{p,na}||_r\lesssim  \lambda_{q+1}^{r}(1/\tau^a)\frac{\lambda_q\delta_q^{1/2}\delta_{q+1}^{1/2}}{\lambda_{q+1}^{1-2\alpha}} \ \text{ for } \ 0\leq r\leq N-k_0^p-m_0-3.\\
        \end{split}
    \end{equation*}
    The implicit constants depend on $r$ and all the other parameters, but not on $a$. In particular, the bounds hold for $0\leq r \leq M$ and $0\leq r \leq M-1$, respectively.
\end{lemma}

\begin{proof}[Proof of Lemma \ref{Rpna}]
Using the definition contained in \eqref{linen}, and the decompositions \eqref{potentialsnashw}, \eqref{potentialsnashb}, we can write:
\begin{equation*}
    \begin{split}
        R^{p,na}&=2\mathcal{R}\ddiv\sum_I\left[\Theta_{w,I}\times\nabla\tilde v_{\ell,I}-\Theta_{b,I}\times\nabla \tilde B_{\ell,I}\right]^\top\\
        &=2\mathcal{R}\ddiv\sum_I\sum_{k=0}^{k_0^p}\frac{(-1)^j}{(j+1)!}\Biggr[\underbrace{\left[\left[\mathcal{L}_{\xi_I^p}^k(\partial_t+\mathcal{L}_{\tilde v_{\ell,I}})\Theta_I^p\right]\times \nabla \tilde v_{\ell,I}\right]^\top}_{S_{1,w}^k}+\underbrace{\left[\left[\mathcal{L}_{\xi_I^p}^k\mathcal{L}_{\tilde B_{\ell,I}}\Theta_I^p\right]\times \nabla \tilde B_{\ell,I}\right]^\top}_{S_{1,b}^k}\Biggr]\\
        &+2\mathcal{R}\ddiv\sum_I\sum_{k=0}^{k_0^p}\frac{(-1)^k}{(j+1)!}\Biggr[\underbrace{\left[[\mathcal{L}_{\xi_I^p}^k\mathcal{L}_{\tilde v_q-\tilde v_{\ell,I}}\Theta_I^p]\times \nabla \tilde v_{\ell,I}\right]^\top}_{S_{2,w}^k}+\underbrace{\left[[\mathcal{L}_{\xi_I^p}^k\mathcal{L}_{\tilde B_q-\tilde B_{\ell,I}}\Theta_I^p]\times \nabla \tilde B_{\ell,I}\right]^\top}_{S_{2,b}^k}\Biggr],
    \end{split}
\end{equation*}
note that the $S_i^k$ we just defined in principle depend also on $I$, but we omit this dependence.  

\noindent Let $\underline{r}=M-2m_0-k_0^g-8$. The estimates in Lemmas \ref{estimatesperturbationnash} with $\bar N=N-m_0-1$ and \ref{localcorrn}, together with Proposition \ref{czstuff} to deal with the operator $\mathcal{R}\ddiv$, allow us to bound:
\begin{equation*}
    \begin{split}
        ||\mathcal{R}\ddiv \sum_I S_{1,w}^k||_{r+\alpha}&\lesssim ||\sum_I S_{1,w}^k||_{r+\alpha}\lesssim \sup_I||S_{1,w}^k||_{r+\alpha}\\
        &\lesssim \sup_I\left[||\mathcal{L}_{\xi_I^p}^k(\partial_t+\mathcal{L}_{\tilde v_{\ell,I}})\Theta_I^p||_{r+\alpha}||\tilde v_{\ell,I}||_{1}+||\mathcal{L}_{\xi_I^p}^k(\partial_t+\mathcal{L}_{\tilde v_{\ell,I}})\Theta_I^p||_{0}||\tilde v_{\ell,I}||_{r+1+\alpha}\right]\\
        &\lesssim (C)^k\lambda_{q+1}^{r-1+\alpha}\lambda_q^{[k-\underline{r} ]^+(b-1)\gamma_\ell}\mathcal{T}_p^k\lambda_q\delta_q^{1/2}\delta_{q+1}^{1/2}\\
        &+(C)^k\lambda_q^{r}\lambda_q^{[r-\underline{r}]^+(b-1)\gamma_\ell}\lambda_q^{[k-\underline{r}]^+(b-1)\gamma_\ell}\ell^{-\alpha}\mathcal{T}_p^k\lambda_{q+1}^{-1}\lambda_q\delta_q^{1/2}\delta_{q+1}^{1/2}\\
        &\lesssim (C)^k\lambda_{q+1}^{r}\lambda_q^{[k-\underline{r} ]^+(b-1)\gamma_\ell}\mathcal{T}_p^k\frac{\lambda_q\delta_q^{1/2}\delta_{q+1}^{1/2}}{\lambda_{q+1}^{1-\alpha}}
    \end{split}
\end{equation*}
for $0\leq r+k\leq N-m_0-2$, similarly,
\begin{equation*}
    ||\mathcal{R}\ddiv \sum_I S_{1,b}^k||_{r+\alpha}\lesssim (C)^k\lambda_{q+1}^{r}\lambda_q^{[k-\underline{r} ]^+(b-1)\gamma_\ell}\mathcal{T}_p^k(\tau^a/\tau^c)\frac{\lambda_q\delta_q^{1/2}\delta_{q+1}^{1/2}}{\lambda_{q+1}^{1-\alpha}} \ \text{ for } \ 0\leq r+k\leq N-m_0-2,
\end{equation*}
moreover, from Lemma \ref{estimatesperturbationnash2} and \ref{localcorrn}, \ref{czstuff} again, we deduce:
\begin{equation*}
    \begin{split}
        ||\mathcal{R}\ddiv \sum_I(S_{2,w}^k+S_{2,b}^k)||_{r+\alpha}&\lesssim||\sum_I(S_{2,w}^k+S_{2,b}^k)||_{r+\alpha}\\
        &\lesssim \sup_I\left[||\mathcal{L}_{\xi_I^p}^k\mathcal{L}_{\tilde v_q-\tilde v_{\ell,I}}\Theta_I^p||_{r+\alpha}||\tilde v_{\ell,I}||_{1}+||\mathcal{L}_{\xi_I^p}^k\mathcal{L}_{\tilde v_q-\tilde v_{\ell,I}}\Theta_I^p||_{0}||\tilde v_{\ell,I}||_{r+1+\alpha}\right]\\
        &+ \sup_I\left[ ||\mathcal{L}_{\xi_I^p}^k\mathcal{L}_{\tilde B_q-\tilde B_{\ell,I}}\Theta_I^p||_{r+\alpha}||\tilde B_{\ell,I}||_1+||\mathcal{L}_{\xi_I^p}^k\mathcal{L}_{\tilde B_q-\tilde B_{\ell,I}}\Theta_I^p||_{0}||\tilde B_{\ell,I}||_{r+1+\alpha}\right]\\
        &\lesssim(C)^k\lambda_{q+1}^{r+\alpha}L_p(k)\mathcal{T}_p^k(\ell\lambda_q)^{m_0}\lambda_q\tau^a\delta_q\delta_{q+1}^{1/2}\\
        &\leq (C)^k\lambda_{q+1}^{r}L_p(k)\mathcal{T}_p^k(\tau^a/\tau^c)\frac{\lambda_q\delta_q^{1/2}\delta_{q+1}^{1/2}}{\lambda_{q+1}}
    \end{split}
\end{equation*}
for $0\leq r \leq N-m_0-2$.

\noindent In the bounds above the implicit constants depend on $r$ and all the parameters but not on $a$, and are uniform in $r, \ k$, moreover, we lost one derivative from the $N-m_0-1$ at our disposal to deduce the $C^{r+\alpha}$ norm by interpolation; we also used the fact that at most one $\Theta_I^p$ is non-zero at each space-time point, see Lemma \ref{disjoint} and $(\ell\lambda_q)^{m_0}\leq \lambda_q/\lambda_{q+1}$. We won't mention these facts again.

\noindent In addition to the Lemmas used above, Proposition \ref{czstuff} and Lemma \ref{stabilitynash} allow us to deal with commutators and correct the transport operator, we deduce:
\begin{equation*}
    \begin{split}
        ||\mathcal{\tilde A}^\pm \mathcal{R}\ddiv \sum_IS_{1,w}^k||_{r+\alpha}&\leq ||\mathcal{R}\ddiv\sum_I\mathcal{\tilde A}^\pm_{\ell,I}S_{1,w}^k||_{r+\alpha}+||[\mathcal{\tilde A}^\pm,\mathcal{R}\ddiv]\sum_I S_{1,w}^k||_{r+\alpha}+||\mathcal{R}\ddiv\sum_I(\tilde z_q^\pm-z^\pm_{\ell,I})\cn S_{1,w}^k||_{r+\alpha}\\
        &\lesssim  \underbrace{\sup_I\left[||\tilde z^\pm_q||_{r+1+\alpha}||S_{1,w}^k||_{\alpha}+||\tilde z^\pm_q||_{1+\alpha}||S_{1,w}^k||_{r+\alpha}\right]}_{\lambda_{q+1}^{r-1+\alpha}\ell^{-\alpha}\lambda_q^{[k-\underline{r} ]^+(b-1)\gamma_\ell}\mathcal{T}_p^k\lambda_q^2\delta_q\delta_{q+1}^{1/2}}\\
        &+\underbrace{\sup_I\left[||\tilde z_q^\pm-\tilde z^\pm_{\ell,I}||_{r+\alpha}||S_{1,w}^k||_{1}+||\tilde z_q^\pm-\tilde z^\pm_{\ell,I}||_{0}|| S_{1,w}^k||_{r+1+\alpha}\right]}_{\lambda_{q+1}^{r+\alpha}\ell^{-\alpha}\lambda_q^{[k-\underline{r} ]^+(b-1)\gamma_\ell}\mathcal{T}_p^k(\ell\lambda_q)^{m_0}\lambda_q\delta_q\delta_{q+1}^{1/2}}\\
        &+\sup_I\Biggr[\underbrace{||\left[\mathcal{A}^\pm_{\ell,I}\mathcal{L}_{\xi_I^p}^k(\partial_t+\mathcal{L}_{\tilde v_{\ell,I}})\Theta_I^p\right]\times \nabla \tilde v_{\ell,I}||_{r+\alpha}}_{\lambda_{q+1}^{r-1+\alpha}\lambda_q^{[k-\underline{r}]^+(b-1)\gamma_\ell}\mathcal{T}_p^k\left[\frac{1}{\tau^a}+\frac{\tau^a}{\tau^c \ell_t}\right]\lambda_q\delta_q^{1/2}\delta_{q+1}^{1/2}}+\underbrace{||\left[\mathcal{L}_{\xi_I^p}^k(\partial_t+\mathcal{L}_{\tilde v_{\ell,I}})\Theta_I^p\right]\times \nabla \mathcal{A}^\pm_{\ell,I}\tilde v_{\ell,I}||_{r+\alpha}}_{\lambda_{q+1}^{r-1+\alpha}\lambda_q^{[k-\underline{r}]^+(b-1)\gamma_\ell}\mathcal{T}_p^k\lambda_q^2\delta_q\delta_{q+1}^{1/2}}\Biggr]\\
        &+\sup_I\Bigg[\underbrace{||\left[\mathcal{L}_{\xi_I^p}^k(\partial_t+\mathcal{L}_{\tilde v_{\ell,I}})\Theta_I^p\right]\times (\DD \tilde z^\pm_{\ell,I})^\top (\DD\tilde v_{\ell,I})^\top||_{r+\alpha}}_{\lambda_{q+1}^{r-1+\alpha}\lambda_q^{[k-\underline{r}]^+(b-1)\gamma_\ell}\mathcal{T}_p^k\lambda_q^2\delta_q\delta_{q+1}^{1/2}}\Biggr]\\
        &\lesssim (C)^k \lambda_{q+1}^{r}\lambda_q^{[k-\underline{r}]^+(b-1)\gamma_\ell}\mathcal{T}_p^k(1/\tau^a)\frac{\lambda_q\delta_q^{1/2}\delta_{q+1}^{1/2}}{\lambda_{q+1}^{1-2\alpha}}
    \end{split} 
\end{equation*}
for $0\leq r \leq N-m_0-2$, where we used $(\ell\lambda_q)^{m_0}\leq \lambda_q/\lambda_{q+1}$ and $\tau^a=\ell_t$.  Taking into account the additional $\tau^a/\tau^c$ smallness of the magnetic part, one can show:
$$||\mathcal{\tilde A}^\pm \mathcal{R}\ddiv \sum_I S_{1,b}^k||_{r+\alpha}\lesssim (C)^k \lambda_{q+1}^{r}\lambda_q^{[k-(\underline{r}-1)]^+(b-1)\gamma_\ell}\mathcal{T}_p^k(1/\tau^c)\frac{\lambda_q\delta_q^{1/2}\delta_{q+1}^{1/2}}{\lambda_{q+1}^{1-2\alpha}} \ \text{ for }\ 0\leq r \leq N-m_0-1.$$
Similarly, we compute:
\begin{equation*}
    \begin{split}
        ||\mathcal{\tilde A}^\pm \mathcal{R}\ddiv \sum_I S_{2,w}^k||_{r+\alpha}&\leq ||\mathcal{R}\ddiv\sum_I\mathcal{\tilde A}^\pm_{\ell,I}S_{2,w}^k||_{r+\alpha}+||[\mathcal{\tilde A}^\pm,\mathcal{R}\ddiv] \sum_I S_{2,w}^k||_{r+\alpha}+||\mathcal{R}\ddiv \sum_I(\tilde z_q^\pm-z^\pm_{\ell,I})\cn S_{2,w}^k||_{r+\alpha}\\
        &\lesssim \underbrace{\sup_I\left[||\tilde z^\pm_q||_{r+1+\alpha}||S_{2,w}^k||_{\alpha}+||\tilde z^\pm_q||_{1+\alpha}||S_{2,w}^k||_{r+\alpha}\right]}_{\lambda_{q+1}^{r+\alpha}\ell^{-\alpha}L_p(k)\mathcal{T}_p^k(\ell\lambda_q)^{m_0}\lambda_q^2\tau^a\delta_q^{3/2}\delta_{q+1}^{1/2}}\\
        &+ \underbrace{\sup_I\left[||\tilde z_q^\pm-\tilde z^\pm_{\ell,I}||_{r+\alpha}||S_{2,w}^k||_{1}+||\tilde z_q^\pm-\tilde z^\pm_{\ell,I}||_{0}||S_{2,w}^k||_{r+1+\alpha}\right]}_{\lambda_{q+1}^{r+1+\alpha}L_p(k)\mathcal{T}_p^k(\ell\lambda_q)^{2m_0}\lambda_q\tau^a\delta_q^{3/2}\delta_{q+1}^{1/2}}\\
        &+ \sup_I\Biggr[\underbrace{||\left[\mathcal{A}^\pm_{\ell,I}\mathcal{L}_{\xi_I^p}^k\mathcal{L}_{\tilde v_q-\tilde v_{\ell,I}}\Theta_I^p\right]\times \nabla \tilde v_{\ell,I}||_{r+\alpha}}_{\lambda_{q+1}^{r+\alpha}L_{p,\mathcal{A}}(k)\mathcal{T}_p^k(\ell\lambda_q)^{m_0}\lambda_q\delta_q\delta_{q+1}^{1/2}}+\underbrace{||\left[\mathcal{L}_{\xi_I^p}^k\mathcal{L}_{\tilde v_q-\tilde v_{\ell,I}}\Theta_I^p\right]\times \nabla \mathcal{A}^\pm_{\ell,I}\tilde v_{\ell,I}||_{r+\alpha}}_{\lambda_{q+1}^{r+\alpha}L_p(k)\mathcal{T}_p^k(\ell\lambda_q)^{m_0}\lambda_q^2\tau^a\delta_q^{3/2}\delta_{q+1}^{1/2}}\Biggr]\\
        &+\sup_I\Biggr[\underbrace{||\left[\mathcal{L}_{\xi_I^p}^k\mathcal{L}_{\tilde v_q-\tilde v_{\ell,I}}\Theta_I^p\right]\times (\DD \tilde z^\pm_{\ell,I})^\top (\DD\tilde v_{\ell,I})^\top||_{r+\alpha}}_{\lambda_{q+1}^{r+\alpha}L_p(k)\mathcal{T}_p^k(\ell\lambda_q)^{m_0}\lambda_q^2\tau^a\delta_q^{3/2}\delta_{q+1}^{1/2}}\Biggr]\\
        &\lesssim (C)^k\lambda_{q+1}^{r}L_{p,\mathcal{A}}(k)\mathcal{T}_p^k(1/\tau^c)\frac{\lambda_q\delta_q^{1/2}\delta_{q+1}^{1/2}}{\lambda_{q+1}^{1-\alpha}}\\
    \end{split} 
\end{equation*}
for $0\leq r+k \leq N-m_0-3$, where we used $(\ell\lambda_q)^{m_0}\leq \lambda_q/\lambda_{q+1}$, $\tau^a=\ell_t$ and $L_p\leq L_{p,\mathcal{A}}$. The same estimate holds for $\mathcal{\tilde A}^\pm \mathcal{R}\ddiv \sum_IS_{2,b}^k$.

\noindent Gathering the above estimates, we conclude:
\begin{equation*}
    \begin{split}
       ||R^{p,na}||_{r+\alpha}&\leq \sum_{k=0}^{k_0^p} \frac{1}{(k+1)!}||\mathcal{R}\ddiv \sum_I\left[S^k_{1,w}+S^k_{1,b}+S^k_{2,w}+S^k_{2,b}\right]||_{r+\alpha}\\
        &\lesssim \sum_{k=0}^{k_0^p} \frac{(C)^k}{(k+1)!}\lambda_{q+1}^{r}\lambda_q^{[k-\underline{r} ]^+(b-1)\gamma_\ell}\mathcal{T}_p^k\frac{\lambda_q\delta_q^{1/2}\delta_{q+1}^{1/2}}{\lambda_{q+1}^{1-\alpha}}\\
        &=\lambda_{q+1}^{r}\frac{\lambda_q\delta_q^{1/2}\delta_{q+1}^{1/2}}{\lambda_{q+1}^{1-\alpha}}\underbrace{\sum_{k=0}^{k_0^p} \frac{\left(C\lambda_q^{(b-1)\gamma_\ell}\mathcal{T}_p\right)^k}{k!}}_{\lesssim 1}\\
        &\lesssim \lambda_{q+1}^{r}\frac{\lambda_q\delta_q^{1/2}\delta_{q+1}^{1/2}}{\lambda_{q+1}^{1-\alpha}}\\
    \end{split}
\end{equation*}
for $0\leq r\leq N-m_0-k_0^p-2$. The implicit constant depends on $r$ and all the other parameters, but not on $a$. From our choice of parameters in \ref{choiceofparameters}, see \eqref{constraintremainder}, it follows that the bound holds in particular for $0\leq r\leq M$.

\noindent We proceed similarly for the Alfv\'en transport bound. Recall the definition of the loss functions in \eqref{lossparameters}. We compute:
\begin{equation*}
    \begin{split}
        ||\mathcal{\tilde A}^\pm R^{p,na}||_{r+\alpha}&\leq \sum_{k=0}^{k_0^p} \frac{1}{(k+1)!}\left[||\mathcal{\tilde A}^\pm \mathcal{R}\ddiv \sum_I S^k_{1,w}||_{r+\alpha}+||\mathcal{\tilde A}^\pm \mathcal{R}\ddiv \sum_IS^k_{1,b}||_{r+\alpha}\right]\\
        &+ \sum_{k=0}^{k_0^p} \frac{1}{(k+1)!}\left[||\mathcal{\tilde A}^\pm \mathcal{R}\ddiv \sum_IS^k_{2,w}||_{r+\alpha}+||\mathcal{\tilde A}^\pm \mathcal{R}\ddiv \sum_IS^k_{2,b}||_{r+\alpha}\right]\\
        &\lesssim \underbrace{\sum_{k=0}^{k_0^p} \frac{1}{(k+1)!}\left[(C)^k \lambda_{q+1}^{r}\lambda_q^{[k-\underline{r}]^+(b-1)\gamma_\ell}\mathcal{T}_p^k(1/\tau^a)\frac{\lambda_q\delta_q^{1/2}\delta_{q+1}^{1/2}}{\lambda_{q+1}^{1-2\alpha}}\right]}_{T_1}\\
        &+ \underbrace{\sum_{k=0}^{k_0^p} \frac{1}{(k+1)!}\left[ (C)^k\lambda_{q+1}^{r}L_{p,\mathcal{A}}(k)\mathcal{T}_p^k(1/\tau^c)\frac{\lambda_q\delta_q^{1/2}\delta_{q+1}^{1/2}}{\lambda_{q+1}^{1-\alpha}}\right]}_{T_2},
    \end{split}
\end{equation*}
we bound $T_1$ as above:
\begin{equation*}
    \begin{split}
        T_1&=\lambda_{q+1}^{r}(1/\tau^a)\frac{\lambda_q\delta_q^{1/2}\delta_{q+1}^{1/2}}{\lambda_{q+1}^{1-2\alpha}}\underbrace{\sum_{k=0}^{k_0^p} \frac{\left(C\lambda_q^{(b-1)\gamma_\ell}\mathcal{T}_p\right)^k}{(k+1)!}}_{\lesssim 1}\lesssim \lambda_{q+1}^{r}(1/\tau^a)\frac{\lambda_q\delta_q^{1/2}\delta_{q+1}^{1/2}}{\lambda_{q+1}^{1-2\alpha}} \ \text{ for } \ 0\leq r\leq N-k_0^p-m_0-3,
    \end{split}
\end{equation*}
we can bound $T_2$ using the explicit formula for $L_{p,\mathcal{A}}$ and splitting the sum accordingly, recall that:
$$L_{p,\mathcal{A}}(k)=1_{r\leq \underline{r}-k_0^g-2}+1_{\underline{r}-k_0^g-1\leq r \leq \underline{r}-1}\left[1+\left(\lambda_q^{(b-1)\gamma_\ell}\mathcal{T}_g\right)^{\underline{r}-1-r}\bar L\right]+1_{r\geq \underline{r}}\lambda_q^{[r-(\underline{r}-1)](b-1)\gamma_\ell}\bar L,$$
we now rewrite:
\begin{equation*}
    \begin{split}
        T_2&= \lambda_{q+1}^{r}(1/\tau^c)\frac{\lambda_q\delta_q^{1/2}\delta_{q+1}^{1/2}}{\lambda_{q+1}^{1-\alpha}}\left[\sum_{k=0}^{\underline{r}-1}\frac{\left(C\mathcal{T}_p\right)^k}{(k+1)!}+\bar L\sum_{k=\underline{r}-k_0^g-1}^{\underline{r}-1}\frac{\left(C\mathcal{T}_p\right)^k}{(k+1)!}\left(\lambda_q^{(b-1)\gamma_\ell}\mathcal{T}_g\right)^{\underline{r}-1-k}\right]\\
        &+\lambda_{q+1}^{r}(1/\tau^c)\frac{\lambda_q\delta_q^{1/2}\delta_{q+1}^{1/2}}{\lambda_{q+1}^{1-\alpha}}\left[\lambda_q^{-(\underline{r}-1)(b-1)\gamma_\ell}\bar L\sum_{k=\underline{r}}^{ k_0^p}\frac{\left(C\lambda_q^{(b-1)\gamma_\ell}\mathcal{T}_p\right)^k}{(k+1)!}\right],\\
    \end{split}
\end{equation*}
where we use the convention that the sum over an empty set is zero. Arguing as in the proof of Lemma \ref{estimatesfieldscircp}, see \eqref{circtransportestimate}, we deduce:
$$T_2\lesssim \lambda_{q+1}^{r}(1/\tau^c)\frac{\lambda_q\delta_q^{1/2}\delta_{q+1}^{1/2}}{\lambda_{q+1}^{1-\alpha}} \ \text{ for } \ 0\leq r\leq N-k_0^p-m_0-3.$$

\noindent Gathering the bounds for $T_1$ and $T_2$, we conclude:
$$||\mathcal{\tilde A}^\pm R^{p,na}||_{r+\alpha}\lesssim \lambda_{q+1}^{r}(1/\tau^a)\frac{\lambda_q\delta_q^{1/2}\delta_{q+1}^{1/2}}{\lambda_{q+1}^{1-2\alpha}} \ \text{ for } \ 0\leq r\leq N-k_0^p-m_0-3$$
where the implicit constant depends on $r$ and all the parameters, but not on $a$. From our choice of parameters in \ref{choiceofparameters}, see \eqref{constraintremainder}, it follows that the bound holds in particular for $0\leq r\leq M-1$.
\end{proof}


\begin{lemma}[Estimates on $R^{p,rm}$]\label{Rprm} Under the choice of parameters in \ref{choiceofparameters}, the following bounds hold.
    \begin{equation*}
        \begin{split}
            &||R^{p,rm}||_r\lesssim \lambda_{q+1}^r\frac{\lambda_q\delta_q^{1/2}\delta_{q+1}^{1/2}}{\lambda_{q+1}^{1-\alpha}} \ \text{ for } \ 0\leq r \leq M,\\
            &||\mathcal{\tilde A}^\pm R^{p,rm}||_r\lesssim \lambda_{q+1}^r\lambda_q\delta_q^{1/2}\frac{\lambda_q\delta_q^{1/2}\delta_{q+1}^{1/2}}{\lambda_{q+1}^{1-\alpha}} \ \text{ for } \ 0\leq r \leq M-1.
        \end{split}
    \end{equation*}
    The implicit constants depend on $r$ and all the parameters, but not on $a$. 
\end{lemma}
\begin{proof}[Proof of Lemma \ref{Rprm}] Recall that:
\begin{equation*}
    \begin{split}
        R^{p,rm}&=\mathcal{R}\curl[\partial_t\vartheta_{w}]+\mathcal{R}\ddiv\left[\tilde v_q\otimes \curl \  \theta_{w}+\curl \  \theta_{w}\otimes \tilde v_q-\tilde B_q\otimes \curl \  \theta_{b}-\curl \  \theta_{b}\otimes \tilde B_q\right]\\
    \end{split}
\end{equation*}
where
$$\theta_{w}=\theta_w^p+\mathring \theta_w^p\ \ \text{ and } \ \ \theta_{b}=\theta_b^p+\mathring \theta_b^p$$
and their definitions in \eqref{decompnashw}, \eqref{decompnashb}. The bounds required can then be found in Proposition \ref{recap} and Lemmas \ref{remaindersmoll}, \ref{estimatesperturbationnash2}. 

\noindent The proof is based on the fact that, according to our choice of parameters in \ref{choiceofparameters}, see \eqref{constraintremainder}, we have:
\begin{equation}\label{remainderconstraintreynolds}
    k_0^p=(N-m_0-5)-M, \ \ \ \ \mathcal{M}_p^{M}\left(\lambda_q^{(b-1)\gamma_\ell}\mathcal{T}_p\right)^{k_0^p}\left(\frac{\lambda_{q+1}}{\lambda_q}\right)^3\lambda_{q+1}\leq \lambda_q\frac{\lambda_q}{\lambda_{q+1}}\delta_q.
\end{equation}
We begin with the estimates. Proposition \ref{czstuff} allows us to deal with $\mathcal{R}\curl$ and we deduce:
\begin{equation*}
    \begin{split}
        ||R^{p,rm}||_{r+\alpha}&\lesssim ||\partial_t\theta_w||_{r+\alpha}+||\tilde v_q||_{r+\alpha}||\theta_w||_1+||\tilde v_q||_{0}||\theta_w||_{r+1+\alpha}+||\tilde B_q||_{r+\alpha}||\theta_b||_1+||\tilde B_q||_{0}||\theta_b||_{r+1+\alpha}\\
        &\lesssim \lambda_{q+1}^{r+\alpha}\mathcal{M}_p^{r+2}\mathcal{T}_p^{k_0^p+1}\lambda_q^{[k_0^p+1-\underline{r}]^+(b-1)\gamma_\ell}\left[\delta_{q+1}^{1/2}+\lambda_{q+1}(\ell\lambda_q)^{m_0}\tau^a\delta_q^{1/2}\delta_{q+1}^{1/2}\right]\\
        &\lesssim \lambda_{q+1}^{r+\alpha}\left(\lambda_q^{(b-1)\gamma_\ell}\mathcal{T}_p\right)^{k_0^p+3}\mathcal{M}_p^{r}\left(\frac{\lambda_{q+1}}{\lambda_q}\right)^2\delta_{q+1}^{1/2}\\
        &\lesssim \lambda_{q+1}^{r+\alpha}\left(\lambda_q^{(b-1)\gamma_\ell}\mathcal{T}_p\right)^{k_0^p+3}\mathcal{M}_p^{M}\left(\frac{\lambda_{q+1}}{\lambda_q}\right)^2\delta_{q+1}^{1/2}
    \end{split}
\end{equation*}
for $0\leq r\leq N-(k_0^p+1)-m_0-3, \ \ 0\leq r\leq M$ and 
\begin{equation*}
    \begin{split}
        ||\mathcal{\tilde A}^\pm R^{rm}||_{r+\alpha}&\lesssim||\partial_tR^{rm}||_{r+\alpha}+||\tilde z_q^\pm||_{r+\alpha}||R^{rm}||_1+||\tilde z_q^\pm||_{0}||R^{rm}||_{r+1+\alpha}\\
        &\lesssim ||\partial_t^2\theta_w||_{r+\alpha}+||[\partial_t \tilde v_q\otimes \curl \ \theta _{w}]^{sym}||_{r+\alpha}+||[\tilde v_q\otimes \curl\  \partial_t\theta_{w}]^{sym}||_{r+\alpha}\\
        &+||[\partial_t \tilde B_q\otimes \curl \ \theta_{b}]^{sym}||_{r+\alpha}+||[\tilde B_q\otimes \curl \ \partial_t \theta_{b}]^{sym}||_{r+\alpha}\\
        &+||\tilde z_q^\pm||_{r+\alpha}||R^{rm}||_1+||\tilde z_q^\pm||_{0}||R^{rm}||_{r+1+\alpha}\\
        &\lesssim \lambda_{q+1}^{r+1+\alpha}\mathcal{M}_p^{r+3}\mathcal{T}_p^{k_0^p+1}\lambda_q^{[k_0^p+1-\underline{r}]^+(b-1)\gamma_\ell}\left[\delta_{q+1}^{1/2}+\lambda_{q+1}(\ell\lambda_q)^{m_0}\tau^a\delta_q^{1/2}\delta_{q+1}^{1/2}\right]\\
        &\lesssim \lambda_{q+1}^{r+1+\alpha}\left(\lambda_q^{(b-1)\gamma_\ell}\mathcal{T}_p\right)^{k_0^p+4}\mathcal{M}_p^{r}\left(\frac{\lambda_{q+1}}{\lambda_q}\right)^3\delta_{q+1}^{1/2}\\
        &\lesssim \lambda_{q+1}^{r+1+\alpha}\left(\lambda_q^{(b-1)\gamma_\ell}\mathcal{T}_p\right)^{k_0^p+4}\mathcal{M}_p^{M}\left(\frac{\lambda_{q+1}}{\lambda_q}\right)^3\delta_{q+1}^{1/2}
    \end{split}
\end{equation*}
for $0\leq r\leq N-(k_0^p+1)-m_0-4, \ \ 0\leq r \leq M$. 

\noindent From \eqref{remainderconstraintreynolds} we deduce:
$$N-(k_0^p+1)-m_0-4=M$$
and thus with quite some margin we can bound:
\begin{equation*}
        \begin{split}
            &||R^{p,rm}||_r\lesssim \lambda_{q+1}^r\frac{\lambda_q\delta_q^{1/2}\delta_{q+1}^{1/2}}{\lambda_{q+1}^{1-\alpha}} \ \text{ for } \ 0\leq r \leq M,\\
            &||\mathcal{\tilde A}^\pm R^{p,rm}||_r\lesssim \lambda_{q+1}^r\lambda_q\delta_q^{1/2}\frac{\lambda_q\delta_q^{1/2}\delta_{q+1}^{1/2}}{\lambda_{q+1}^{1-\alpha}} \ \text{ for } \ 0\leq r \leq M-1
        \end{split}
    \end{equation*}
where the implicit constants depend on $r$ and all the other parameters, but not on $a$.
\end{proof}


\begin{lemma}[Estimates on $R^{p,mo}$]\label{Rpmo} Under the choice of parameters in \ref{choiceofparameters}, the following bounds hold.
    \begin{equation*}
        \begin{split}
            &||R^{p,mo}||_r\lesssim \lambda_{q+1}^r\frac{\lambda_q\delta_{q}^{1/2}\delta_{q+1}^{1/2}}{\lambda_{q+1}} \ \text{ for } \ 0 \leq r \leq M,\\
            &||\mathcal{\tilde A}^\pm R^{p,mo}||_r\lesssim \lambda_{q+1}^r(1/\tau^a)\frac{\lambda_q\delta_{q}^{1/2}\delta_{q+1}^{1/2}}{\lambda_{q+1}} \ \text{ for } \ 0\leq r \leq M-1.\\
        \end{split}
    \end{equation*}
The implicit constants depend on $r$ and all the other parameters, but not on $a$.
\end{lemma}

\begin{proof}[Proof of Lemma \ref{Rpmo}] Recall that:
\begin{equation*}
    \begin{split}
        R^{p,mo}&=\sum_I\left[(\tilde v_q-\tilde v_{\ell,I})\otimes w_I+w_I\otimes (\tilde v_q-\tilde v_{\ell,I})\right]-\left[(\tilde B_q-\tilde B_{\ell,I})\otimes b_I+b_I\otimes (\tilde B_q-\tilde B_{\ell,I})\right],\\
    \end{split}
\end{equation*}
we will use the bounds in Lemmas \ref{leadingtermsnash}, \ref{estimatesfieldscircp} and \ref{estimatesfieldsp}, even if technically speaking those are for $\sum_Iw_I, \ \sum_Ib_I$ (see the definitions in \eqref{decompnashw}, \ref{decompnashb}, \ref{potentialsnashw1}, \eqref{potentialsnashb1}), together with the ones in \ref{standardmollnash}, \ref{stabilitynash}, for $\tilde v_q-\tilde v_{\ell,I}, \ B_q-\tilde B_{\ell,I}$. 

\noindent We first deduce:
\begin{equation*}
    \begin{split}
        ||R^{p,mo}||_r&\lesssim \sup_I\left[||\tilde v_q-\tilde v_{\ell,I}||_r||w_I||_0+||\tilde v_q-\tilde v_{\ell,I}||_0||w_I||_r+||\tilde B_q-\tilde B_{\ell,I}||_r||b_I||_0+||\tilde B_q-\tilde B_{\ell,I}||_0||b_I||_r\right]\\
        &\lesssim \lambda_{q+1}^r(\ell\lambda_q)^{m_0}\delta_q^{1/2}\delta_{q+1}^{1/2}\\
        &\lesssim \lambda_{q+1}^r\frac{\lambda_q\delta_q^{1/2}\delta_{q+1}^{1/2}}{\lambda_{q+1}}\\
    \end{split}
\end{equation*}
for $0\leq r\leq M$, where we used $(\ell\lambda_q)^{m_0}\leq \lambda_q/\lambda_{q+1}$ and the fact that at most one $w_I, \ b_I$ is non zero at each space-time point, see Lemma \ref{disjoint}.

\noindent We move to the Alfv\'en transport bounds. We first compute:
\begin{equation*}
    \begin{split}
        \mathcal{\tilde A}^\pm\left[(\tilde v_q-\tilde v_{\ell,I})\otimes w_I\right]&= [\mathcal{\tilde A}^\pm(\tilde v_q-\tilde v_\ell)+(\tilde z^\pm_q-\tilde z^\pm_{\ell,I})\cn (\tilde v_\ell-\tilde v_{\ell,I})+\mathcal{\tilde A}^\pm_{\ell,I}(\tilde v_\ell-\tilde v_{\ell,I})]\otimes w_I\\
        &+(\tilde v_q-\tilde v_{\ell,I})\otimes [(\tilde z^\pm_q-\tilde z^\pm_{\ell,I})\cn w_I+\mathcal{\tilde A}^\pm_{\ell,I}w_I].
    \end{split}
\end{equation*}
Before proceeding, for $\underline{r}=M-m_0-k_0^g-6$ using Proposition \ref{deepmollification} with the estimates in Proposition \ref{recap}, we bound:
\begin{equation*}
    \begin{split}
         ||\mathcal{\tilde A}^\pm(\tilde v_q-\tilde v_\ell)||_r\lesssim\lambda_q^{r+1}(\ell\lambda_q)^{m_0}\delta_q\left[1_{r\leq \underline{r}-m_0-2}+1_{r\geq\underline{r}-m_0-1}\left(\frac{\lambda_{q+1}}{\lambda_{q}}\right)^{r-(\underline{r}-m_0-2)}\frac{\lambda_{q}}{\lambda_{q+1}(\ell\lambda_q)^{m_0}}\right] 
    \end{split}
\end{equation*}
for $0\leq r \leq M-1$, where the second range of derivative comes from the fact that for $\underline{r}-m_0-2\leq r\leq M-1$ instead of using the mollification error bound, we estimate each term separately, trade a good derivative $\lambda_q$, which we have left for a bad $\lambda_{q+1}$ to compensate for the $1/(\ell\lambda_q)^{m_0}$ loss and use the bound \ref{boundlossfunction} for $L_{g,\mathcal{A}}$. The derivative count is as follows: we have $\underline{r}-1$ good derivatives on the transport of $\tilde v_q$ (see \ref{recap} and the definition of $L_{g,\mathcal{A}}$ in \ref{estimatesfieldsringg}), lose $m_0$ in the mollification error estimate from \ref{deepmollification}, and a final one as just explained. 

\noindent From this, we deduce:
\begin{equation*}
    \begin{split}
        ||\mathcal{\tilde A}^\pm(\tilde v_q-\tilde v_\ell)\otimes w_I||_r&\lesssim ||\mathcal{\tilde A}^\pm(\tilde v_q-\tilde v_\ell)||_r||w_I||_0+||\mathcal{\tilde A}^\pm(\tilde v_q-\tilde v_\ell)||_0||w_I||_r\\
        &\lesssim\lambda_q^{r+1}(\ell\lambda_q)^{m_0}\delta_q\delta_{q+1}^{1/2}\left[1_{r\leq \underline{r}-m_0-2}+1_{r\geq\underline{r}-m_0-1}\left(\frac{\lambda_{q+1}}{\lambda_{q}}\right)^{r-(\underline{r}-m_0-2)}\frac{\lambda_{q}}{\lambda_{q+1}(\ell\lambda_q)^{m_0}}\right] \\
        &+\lambda_{q+1}^{r}(\ell\lambda_q)^{m_0}\lambda_q\delta_q\delta_{q+1}^{1/2}\\
        &\lesssim \lambda_{q+1}^{r}\lambda_q\delta_q^{1/2}\frac{\lambda_q\delta_q^{1/2}\delta_{q+1}^{1/2}}{\lambda_{q+1}}
    \end{split}
\end{equation*}
for $0\leq r \leq M-1$.

\noindent The other terms are less troublesome, and from the Lemmas recalled above, we deduce:
\begin{equation*}
    \begin{split}
        ||\mathcal{\tilde A}^\pm\left[(\tilde v_q-\tilde v_{\ell,I})\otimes w_I\right]||_r&\lesssim ||\mathcal{\tilde A}^\pm(\tilde v_q-v_\ell)\otimes w_I||_r+||[(\tilde z^\pm_q-\tilde z^\pm_{\ell,I})\cn +\mathcal{\tilde A}^\pm_{\ell,I}](\tilde v_\ell-\tilde v_{\ell,I})\otimes  w_I||_r\\
        &+||(\tilde v_q-\tilde v_{\ell,I}) \otimes [(\tilde z^\pm_q-\tilde z^\pm_{\ell,I})\cn w_I+\mathcal{\tilde A}^\pm_{\ell,I}w_I]||_r\\
        &\lesssim\lambda_{q+1}^{r}\lambda_q\delta_q^{1/2}\frac{\lambda_q\delta_q^{1/2}\delta_{q+1}^{1/2}}{\lambda_{q+1}}+\lambda_{q+1}^r(\ell\lambda_q)^{m_0}(1/\tau^a)\delta_q^{1/2}\delta_{q+1}^{1/2}+\lambda_{q+1}^r(\ell\lambda_q)^{m_0}\mathcal{T}_p(1/\tau^a)\delta_q\\
        &\lesssim \lambda_{q+1}^r(1/\tau^a)\frac{\lambda_q\delta_{q}^{1/2}\delta_{q+1}^{1/2}}{\lambda_{q+1}}
    \end{split}
\end{equation*}
for $0\leq r \leq M-1$. 

\noindent Similar bounds hold for $\mathcal{\tilde A}^\pm\left[(\tilde B_q-\tilde B_{\ell,I})\otimes b_I\right]$ and from the fact that at most one $w_I, \ b_I$ is non zero at each space-time point, see Lemma \ref{disjoint}, we conclude:
\begin{equation*}
    ||\mathcal{\tilde A}^\pm R^{p,mo}||_r\lesssim \lambda_{q+1}^r(1/\tau^a)\frac{\lambda_q\delta_{q}^{1/2}\delta_{q+1}^{1/2}}{\lambda_{q+1}} \ \text{ for } \ 0\leq r \leq M-1
\end{equation*}
where the implicit constant depends on $r$ and all the parameters, but not on $a$. 
\end{proof}


\subsubsection{Estimates on the Quadratic Error}
\begin{lemma}[Estimates on $R^{p,qua}$]\label{Rpqua} Under the choice of parameters in \ref{choiceofparameters}, the following bounds hold:
\begin{equation*}
    \begin{split}
        ||R^{p,qua}||_r&\lesssim \lambda_{q+1}^r\frac{\lambda_q\delta_{q+1}}{\lambda_{q+1}^{1-\alpha}}+\lambda_{q+1}^r(\tau^a/\tau^c)\delta_{q+1} \ \text{ for } \ 0\leq r \leq M,\\
        ||\mathcal{\tilde A}^\pm R^{p,qua}||_r&\lesssim  \lambda_{q+1}^{r+\alpha}1/\tau^a\delta_{q+1} \ \text{ for } \ 0\leq r \leq M-1.\\
    \end{split}
\end{equation*}
The implicit constants depend on $r$ and all the parameters, but not on $a$.    
\end{lemma}

\begin{remark}[Regularity Bottleneck] Note that:
$$\lambda_{q+1}^r\frac{\lambda_q\delta_{q+1}}{\lambda_{q+1}^{1-\alpha}}\leq\lambda_{q+1}^r(\tau^a/\tau^c)\delta_{q+1}$$
we state the result this way to show the first and more familiar bound as well, the leading one comes from $T_2$ below, see \eqref{quadraticsplittingn}, and is the same as the one for the cut-off error $R^{cut}$ from the Galbrun stage given in Lemma \ref{Rcut}. One might be able to correct $R^{cut}$ by introducing a full Newton iteration, as in the work of Giri-Radu \cite{GR}. 
\end{remark}

\begin{proof}[Proof of Lemma \ref{Rpqua}]
We first recall the decomposition \eqref{Rp} and \eqref{enforced}, then we rewrite:
\begin{equation*}
\begin{split}
    \ddiv \ R^{p,qua}&=\ddiv\left[w^p\otimes w^p-b^p\otimes b^p -\sum_I g_I^2 \tilde{A}_I\right]\\
    &=\ddiv\left[\sum_I\left[w^{p,o}_I\otimes w^{p,o}_I-g_I^2 \tilde{A}_I\right]+\sum_I(w^p-w^{p,o}_I)\otimes w^p+w^p\otimes (w^p-w^{p,o}_I)-b^p\otimes b^p\right]
\end{split}
\end{equation*}
where the definition of $w^{p,o}_I$ was given in \eqref{leadingtermnashexpansion}, see also \eqref{decompnashw} and then set: 
\begin{equation}\label{quadraticsplittingn}
\begin{split}
    R^{p,qua}&=\underbrace{\mathcal{R}\ddiv\sum_I\left[w^{p,o}_I\otimes w^{p,o}_I-g_I^2 \tilde{A}_I\right]}_{T_1}+\underbrace{\sum_I(w^p-w^{p,o}_I)\otimes w^p+w^p\otimes (w^p-w^{p,o}_I)-b^p\otimes b^p}_{T_2}.
\end{split}
\end{equation}

\noindent \textbf{Estimates on $T_1$.} This is the key term in all convex integration schemes. We use the standard space-averaging method and apply $\mathcal{R}\ddiv$ to the remaining zero-average fast-oscillating terms. This achieves the necessary smallness, as those were designed to be a pressure-less stationary solution of the momentum equation in the fast variables. We will need to adjust the proof of the standard Stationary Phase Lemma to the case at hand.

\noindent Fix an index $I\in \mathcal{I}$. We begin with the space-averaging, namely
\begin{equation}\label{ddivRqua}
    \begin{split}
        &\ddiv\left[w^{p,o}_I\otimes w^{p,o}_I-g_I^2 \tilde{A}_I\right]\\
        &=\ddiv\left[\chi^2_I \tilde a^2_Ig^2_I(\tilde \Psi^{*}_I(\varphi'_{\lambda_{q+1},k_I}))^2\tilde\Psi^{2*}_I\zeta_I\otimes \tilde \Psi^{2*}_I\zeta_I-\chi^2_I \tilde a^2_Ig^2_I\tilde \Psi^{2*}_I\zeta_I\otimes \tilde \Psi^{2*}\zeta_I\right]\\
        &=\ddiv\left[\chi^2_I \tilde a^2_Ig^2_I\left((\tilde \Psi^{*}_I(\varphi'_{\lambda_{q+1},k_I}))^2-1\right)\tilde\Psi^{2*}_I\zeta_I\otimes \tilde \Psi^{2*}_I\zeta_I\right]\\
        &=\chi^2_I \tilde a^2_Ig^2_I\left((\tilde \Psi^{*}_I(\varphi'_{\lambda_{q+1},k_I}))^2-1\right)\ddiv(\tilde\Psi^{2*}_I\zeta_I) \tilde \Psi^{2*}_I\zeta_I+g^2_I\tilde\Psi^{2*}_I\zeta_I\cn\left[\left((\tilde\Psi^{*}_I(\varphi'_{\lambda_{q+1},k_I}))^2-1\right)\chi^2_I \tilde a^2_I\tilde \Psi^{2*}_I\zeta_I\right]\\
        &=\underbrace{\left((\tilde \Psi^{*}_I(\varphi'_{\lambda_{q+1},k_I}))^2-1\right)}_{=\phi(\lambda_{q+1}\tilde \Psi_I)}\underbrace{g^2_I\tilde\Psi^{2*}_I\zeta_I\cn(\chi^2_I \tilde a^2_I\tilde \Psi^{2*}_I\zeta_I)}_{=S_I}\\
        &=\phi(\lambda_{q+1}\tilde \Psi_I)S_I.
    \end{split}
\end{equation}
where we used the properties in \eqref{DG0} and \eqref{DG5} to compute:
\begin{equation*}
        \ddiv(\tilde\Psi^{2*}_I\zeta_I)=\dd \tilde\Psi^{2*}_I\zeta_I=\tilde\Psi^{*}_I\dd \zeta=0 
\end{equation*}
and
\begin{equation*}
    \begin{split}
        \tilde\Psi^{2*}_I\zeta_I\cn\left((\tilde \Psi^{*}_I(\varphi'_{\lambda_{q+1},k_I}))^2-1\right)&=\tilde\Psi^{2*}_I\zeta_I\cn (\tilde \Psi^{*}_I(\varphi'_{\lambda_{q+1},k_I}))^2\\
        &=\tilde\Psi^{2*}_I\zeta_I\cn (\tilde \Psi^{*}_I((\varphi'_{\lambda_{q+1},k_I})^2))\\
        &=\tilde\Psi^{2*}_I\zeta_I\cdot\tilde \Psi^{1*}_I(\nabla (\varphi'_{\lambda_{q+1},k_I})^2))\\
        &=2\lambda_{q+1}\tilde \Psi^{*}_I(\varphi_{\lambda_{q+1},k_I}'\varphi_{\lambda_{q+1},k_I}'')\tilde \Psi^{2*}_I\zeta_I\cdot\tilde \Psi^{1*}_Ik_I\\
        &=2\lambda_{q+1}\tilde \Psi^{*}_I(\varphi_{\lambda_{q+1},k_I}'\varphi_{\lambda_{q+1},k_I}'')[\det\DD\tilde\Psi_I]\tilde \Psi^{*}_I\zeta_I\cdot\tilde \Psi^{1*}_Ik_I\\
        &=2\lambda_{q+1}\tilde \Psi^{*}_I(\varphi_{\lambda_{q+1},k_I}'\varphi_{\lambda_{q+1},k_I}'')[\det\DD\tilde\Psi_I]\tilde \Psi^{*}_I(\zeta_I\cdot k_I)=0.
    \end{split}
\end{equation*}

\noindent By definition, $\varphi=\sqrt{2}\sin$ and thus: 
\begin{equation*}
\begin{split}
    \phi(\lambda_{q+1}\tilde \Psi_I)&=\left(\tilde\Psi^{*}_I(\varphi'_{\lambda_{q+1},k_I}))^2-1\right)\\
    &=2\cos^2(\lambda_{q+1}\tilde\Psi_I\cdot k_I)-1=\cos(2\lambda_{q+1}\tilde\Psi_I\cdot k_I)\\
    &=\frac{e^{i2\lambda_{q+1}\tilde\Psi_I\cdot k_I}+e^{-i2\lambda_{q+1}\tilde\Psi_I\cdot k_I}}{2}.
\end{split}
\end{equation*}
We now want to apply the operator $\mathcal{R}$ in \eqref{invdiv} to \eqref{ddivRqua} and gain smallness by a stationary phase argument; see the classical \cite[Lemma 2.2]{DaSz}. Since $\tilde\Psi_I$ is not globally defined on $\mathbb{T}^3$, we need to adapt its proof to this case, note however that it is defined for $I=(\zeta,j,j')$ on $B_{\tau^c}(x_{j})\times B_{\tau^c}(t_j)=Q_J$ and by definition:
\begin{equation}\label{supportRqua}
    \supp_{x,t}S_I\subset \supp_{x,t} \tilde a_I\Subset Q_J
\end{equation}
see \eqref{supportslow} in Lemma \ref{slowcoeffestimates}. Recall from Lemma \ref{chartprop} that 
\begin{equation}\label{alphaquadratic}
    |\DD \tilde \Psi_I(x,t)-\IId|\lesssim \lambda_{q+1}^{-\alpha} \ \text{ for } \ (x,t)\in Q_J
\end{equation}
in particular, we can choose $a$ sufficiently large so that
$|\DD \tilde \Psi_I(x,t)-\IId|< 1/2 \ \text{ for } \ (x,t)\in Q_J$
this will be useful in a moment. By construction, we have:
\begin{equation}\label{zeroavg}
    \int_{\mathbb{T}^3}\phi(\lambda_{q+1}\tilde \Psi_I)S_I=\int_{\mathbb{T}^3}\ddiv \left[w^{p,o}_I\otimes w^{p,o}_I-g_I^2 \tilde{A}_I\right]=0
\end{equation}
and mimicking the proof of \cite[Lemma 2.2]{DaSz}, for some $\bar N$ to be chosen later we rewrite: 
\begin{equation}\label{rhsRqua}
    \begin{split}
     \sum_I\phi(\lambda_{q+1}\tilde \Psi_I)S_I&=\sum_IS_I\frac{e^{i2\lambda_{q+1}\tilde\Psi_I\cdot k_I}+e^{-i2\lambda_{q+1}\tilde\Psi_I\cdot k_I}}{2}\\
     &=\ddiv\left[\sum_I\sum_{n=0}^{\bar N-1}S_{I,n}\frac{\nabla(\tilde\Psi_I\cdot k_I)}{|\nabla(\tilde\Psi_I\cdot k_I)|^2}\frac{e^{i2\lambda_{q+1}\tilde\Psi_I\cdot k_I}+e^{-i2\lambda_{q+1}\tilde\Psi_I\cdot k_I}}{2}\frac{1}{(i2\lambda_{q+1})^{n+1}}+R_{\bar N}\right]\\
     &=\ddiv\left[\sum_I\sum_{n=0}^{\bar N-1}S_{I,n}\frac{\nabla(\tilde\Psi_I\cdot k_I)}{|\nabla(\tilde\Psi_I\cdot k_I)|^2}\phi(\lambda_{q+1}\tilde \Psi_I)\frac{1}{(i2\lambda_{q+1})^{n+1}}+R_{\bar N}\right]\\
    \end{split}
\end{equation}
where $S_{I,n}$ is recursively defined from $S_{I,0}=S_I$ by 
$$S_{I,n+1}=-\ddiv\left[S_{I,n}\frac{\nabla(\tilde\Psi_I\cdot k_I)}{|\nabla(\tilde\Psi_I\cdot k_I)|^2}\right] \ \text{ for }\  0\leq n\leq \bar N$$
and $R_{\bar N}$ solves:
\begin{equation}\label{diveq}
    \begin{cases}
        \ddiv \ R_{\bar N}= \frac{1}{(i2\lambda_{q+1})^{\bar N}}\sum_I S_{I,\bar N}\frac{e^{i2\lambda_{q+1}\tilde\Psi_I\cdot k_I}+e^{-i2\lambda_{q+1}\tilde\Psi_I\cdot k_I}}{2}=\frac{1}{(i2\lambda_{q+1})^{\bar N}}\sum_I S_{I,\bar N}\phi(\lambda_{q+1}\tilde \Psi_I)\\
        \int_{\mathbb{T}^3} R_{\bar N}=0
    \end{cases}
\end{equation}
From \eqref{zeroavg}, we deduce that the right-hand side also has zero average, and a unique solution exists. Moreover, from \eqref{supportRqua} and  \ref{alphaquadratic} we can bound from below:
\begin{equation*}
    \begin{split}
        \nabla(\tilde\Psi_I\cdot k_I)= k_I+(\DD \tilde\Psi_I-\IId)^\top[k_I] \Longrightarrow |\nabla(\tilde\Psi_I\cdot k_I)|> 1-1/2=1/2
    \end{split}
\end{equation*}
so that the expressions above are all well-defined. 

\noindent Now let $\underline{r}=M-2m_0-k_0^g-8$, from the bounds in Lemmas \ref{chartprop}, \ref{slowcoeffestimates}, and the observation above we get:
\begin{equation}\label{Sn}
    \begin{split}
        ||S_{I,0}||_{r+\alpha}
        &\lesssim\lambda_q^{r+1}\ell^{-\alpha}\lambda_q^{[r-(\underline{r}-1)]^+(b-1)\gamma_\ell}\delta_{q+1},\\
        ||\nabla(\tilde\Psi_I\cdot k_I)/|\nabla(\tilde\Psi_I\cdot k_I)|^2||_{r+\alpha}&\lesssim\lambda_q^{r}\ell^{-\alpha}\lambda_q^{[r-\underline{r}]^+(b-1)\gamma_\ell}\delta_{q+1}\\
        ||\phi(\lambda_{q+1}\tilde \Psi_I)||_{r}&\lesssim \lambda_{q+1}^r\\
    \end{split}
\end{equation}
from this and induction, we also deduce:
\begin{equation}\label{Sn1}
    ||S_{I,n}||_{r+\alpha}\lesssim \lambda_q^{r+n+1}\ell^{-\alpha}\lambda_q^{[r+n+1-\underline{r}]^+(b-1)\gamma_\ell}\delta_{q+1},
\end{equation}
finally, from standard estimates for the divergence equation in \eqref{diveq} and the bounds in \eqref{Sn},\eqref{Sn1} we obtain: 
$$||R_{\bar N}||_{r+\alpha}\lesssim ||\frac{1}{(i2\lambda_{q+1})^{\bar N}}\sum_I S_{I,\bar N}\phi(\lambda_{q+1}\Psi_I)||_{r+\alpha}\lesssim \lambda_{q+1}^{r+\alpha}\left[\lambda_{q+1}\left(\frac{\lambda_q}{\lambda_{q+1}}\right)^{\bar N(1-\gamma_\ell)}\right]\frac{\lambda_q\delta_{q+1}}{\lambda_{q+1}}$$
and given $\gamma_\ell,\ b$ we can choose $\bar N$ so large that:
\begin{equation}\label{remainderquadratic}
    \lambda_{q+1}\left(\frac{\lambda_q}{\lambda_{q+1}}\right)^{\bar N(1-\gamma_\ell)}\leq 1\ \Longrightarrow \ ||R_{\bar N}||_{r+\alpha}\lesssim \lambda_{q+1}^{r+\alpha}\frac{\lambda_q\delta_{q+1}}{\lambda_{q+1}}
\end{equation}

\noindent Now recall the definition of $\mathcal{R}$ in \eqref{invdiv}, it follows that for $u$, solving on $\mathbb{T}^3$:
\begin{equation*}
    \begin{cases}
        \Delta u=\sum_I\phi(\lambda_{q+1}\tilde\Psi_I)S_I,\\
        \int_{\mathbb{T}^3} u=0
    \end{cases}
\end{equation*}
for every fixed time and $\mathbb{P}$ the Leray projector we have: 
$$T_1=\mathcal{R}\sum_I\phi(\lambda_{q+1}\tilde \Psi_I)S_I=\frac{1}{4}\left[\DD \mathbb{P}u+(\DD \mathbb{P}u)^\top\right]+\frac{3}{4}\left[\DD u+(\DD u)^\top\right]-\frac{1}{2}(\ddiv \ u)\IId$$ 
where the Poisson problem admits a unique solution since \eqref{zeroavg} holds.

\noindent We wish to show that $T_1$ actually gained smallness compared to $S_I$, to do so, we use the rewriting of the right-hand side of the Poisson equation \eqref{rhsRqua}. Using Proposition \ref{czstuff} to deal with $\nabla\Delta^{-1}\ddiv$, the fact that at each space-time point at most one $S_I$ is non-zero, see Lemma \ref{disjoint}, and the bounds in \eqref{Sn}, \eqref{Sn1}, \eqref{remainderquadratic} we deduce:
\begin{equation*}
    \begin{split}
        ||\nabla u||_{r+\alpha}&=||\nabla\Delta^{-1}\ddiv\left[\sum_I\sum_{n=0}^{\bar N-1}S_{I,n}\frac{\nabla(\tilde\Psi_I\cdot k_I)}{|\nabla(\tilde\Psi_I\cdot k_I)|^2}\phi(\lambda_{q+1}\tilde \Psi_I)\frac{1}{(i2\lambda_{q+1})^{n+1}}+R_{\bar N}\right]||_{r+\alpha}\\
        &\lesssim||\sum_I\sum_{n=0}^{\bar N-1}S_{I,n}\frac{\nabla(\tilde\Psi_I\cdot k_I)}{|\nabla(\tilde\Psi_I\cdot k_I)|^2}\phi(\lambda_{q+1}\tilde \Psi_I)\frac{1}{(i2\lambda_{q+1})^{n+1}}+R_{\bar N}||_{r+\alpha}\\
        &\lesssim \sup_I\sum_{n=0}^{\bar N-1}\frac{1}{(\lambda_{q+1})^{n+1}}||S_{I,n}\frac{\nabla(\tilde\Psi_I\cdot k_I)}{|\nabla(\tilde\Psi_I\cdot k_I)|^2}\phi(\lambda_{q+1}\tilde \Psi_I)||_{r+\alpha}+||R_{\bar N}||_{r+\alpha}\\
        &\lesssim  \lambda_{q+1}^{r}\frac{\lambda_q\delta_{q+1}}{\lambda_{q+1}^{1-\alpha}}\sum_{n=0}^{\bar N-1} \left(\frac{\lambda_q}{\lambda_{q+1}}\right)^{n(1-\gamma_\ell)}+\lambda_{q+1}^{r}\frac{\lambda_q\delta_{q+1}}{\lambda_{q+1}^{1-\alpha}}\\
        &\lesssim \lambda_{q+1}^{r}\frac{\lambda_q\delta_{q+1}}{\lambda_{q+1}^{1-\alpha}}
    \end{split}
\end{equation*}
where the implicit constant depends on $r, \ \bar N$ and all the other parameters but not on $a$ and we conclude that:
$$||T_1||_{r+\alpha}\lesssim \lambda_{q+1}^r\frac{\lambda_q\delta_{q+1}}{\lambda_{q+1}^{1-\alpha}} \ \text{ for } \ r\geq 0.$$

\noindent\textit{Alfv\'en transport.} Thanks to \eqref{ddivRqua} we can write:
$$T_1=\mathcal{R}\ddiv\sum_I\left[\phi(\lambda_{q+1}\tilde \Psi_I)g_I^2\tilde A_I\right]$$
and with the help of Proposition \ref{czstuff} to deal with the CZ operator $\mathcal{R}\ddiv$ and the commutators, Lemmas \ref{standardmollnash}, \ref{stabilitynash} to correct the transport operators, the estimates in Lemmas \ref{localcorrn}, \ref{slowcoeffestimates}, and the transport properties in Lemma \ref{chartprop}, we deduce:
\begin{equation*}
    \begin{split}
        ||\mathcal{\tilde A}^\pm T_1||_{r+\alpha}&\leq||\mathcal{R}\ddiv\sum_I\mathcal{\tilde A}^\pm_{\ell,I}[\phi(\lambda_{q+1}\tilde \Psi_I)g^2_I\tilde A_I]||_{r+\alpha}+||[\mathcal{\tilde A}^\pm,\mathcal{R}\ddiv]\sum_I\phi(\lambda_{q+1}\tilde \Psi_I)g^2_I\tilde A_I||_{r+\alpha}\\
        &+||\mathcal{R}\ddiv\sum_I(\tilde z^\pm_q-\tilde z^\pm_{\ell,I})\cn[\phi(\lambda_{q+1}\tilde \Psi_I)g^2_I\tilde A_I]||_{r+\alpha}\\
        &\lesssim\sup_I\Biggr[\underbrace{2/\tau^a\sup_t|g_Ig_I'|||\phi(\lambda_{q+1}\tilde \Psi_I)\tilde A_I||_{r+\alpha}}_{\lambda_{q+1}^{r+\alpha}(1/\tau^a)\delta_{q+1}}+\underbrace{\sup_tg_I^2||\phi(\lambda_{q+1}\tilde \Psi_I)\mathcal{\tilde A}^\pm_{\ell,I}\tilde A_I]||_{r+\alpha}}_{\lambda_{q+1}^{r+\alpha}(1/\tau^c)\delta_{q+1}}\Biggr]\\
        &+\sup_I\Biggr[\underbrace{\sup_tg_I^2||\tilde z^\pm_q||_{r+1+\alpha}||\phi(\lambda_{q+1}\tilde \Psi_I)\tilde A_I]||_{\alpha}+\sup_tg_I^2||\tilde z^\pm_q||_{1+\alpha}||\phi(\lambda_{q+1}\tilde \Psi_I)\tilde A_I]||_{r+\alpha}}_{\lambda_{q+1}^{r+\alpha}\lambda_q\ell^{-\alpha}\delta_q^{1/2}\delta_{q+1}}\Biggr]\\
        &+\sup_I\Biggr[\underbrace{\sup_tg_I^2||\tilde z^\pm_q-\tilde z^\pm_{\ell,I}||_{r+\alpha}||\phi(\lambda_{q+1}\tilde \Psi_I)\tilde A_I||_{1}+\sup_tg_I^2||\tilde z^\pm_q-\tilde z^\pm_{\ell,I}||_{0}||\phi(\lambda_{q+1}\tilde \Psi_I)\tilde A_I||_{r+1+\alpha}}_{\lambda_{q+1}^{r+1+\alpha}(\ell\lambda_q)^{m_0}\delta_q^{1/2}\delta_{q+1}}\Biggr]\\
        &\lesssim \lambda_{q+1}^{r+\alpha}1/\tau^a\delta_{q+1}
    \end{split}
\end{equation*}
for $0\leq r\leq N-m_0-1$, where we used $\gamma_a\geq \gamma_{CZ}, \ (\ell\lambda_q)^{m_0}\leq \lambda_q/\lambda_{q+1}$ and the fact that at most one $g_I\tilde A_I$ is non-zero at each space-time point, see Lemma \ref{disjoint}.

\noindent \textbf{Estimates on $T_2$.} The bounds follow immediately from Lemmas \ref{leadingtermsnash}, \ref{estimatesfieldscircp}, \ref{estimatesfieldsp} and we omit the details:
\begin{equation*}
    \begin{split}
        ||T_2||_r&\lesssim \lambda_{q+1}^r(\tau^a/\tau^c)\delta_{q+1} \text{ for }\ 0\leq r \leq M,\\
        ||\mathcal{\tilde A}^\pm T_2||_r&\lesssim \lambda_{q+1}^r(1/\tau^c)\delta_{q+1}\ \text{ for } \ 0\leq r \leq M-1.\\
    \end{split}
\end{equation*}

\noindent \textbf{Conclusion.} Gathering the bounds on $T_1, \ T_2$ above we conclude that:
\begin{equation*}
    \begin{split}
        ||R^{p,qua}||_r&\lesssim \lambda_{q+1}^r\frac{\lambda_q\delta_{q+1}}{\lambda_{q+1}^{1-\alpha}}+\lambda_{q+1}^r(\tau^a/\tau^c)\delta_{q+1} \ \text{ for } \ 0\leq r \leq M,\\
        ||\mathcal{\tilde A}^\pm R^{p,qua}||_r&\lesssim  \lambda_{q+1}^{r+\alpha}(1/\tau^a)\delta_{q+1} \ \text{ for } \ 0\leq r \leq M-1\\
    \end{split}
\end{equation*}
where the implicit constants depend on $r$ and all the other parameters, but not on $a$.
\end{proof}


\section{Conclusion}\label{conclusion}
\subsection{Proof of the Iterative Proposition \ref{iterative}.} \label{proofiterative} We need to ensure that the bounds in \eqref{inductiveassumptionsgeneral} are verified with $q$ replaced by $q+1$. We begin with the estimates for the perturbation and $(v_{q+1}, B_{q+1})$, then move to those for the new Reynolds stress $R_{q+1}$.

\noindent \textbf{$||\cdot||_r$ estimates on the perturbation.} Recall from the decompositions \eqref{decompnashw} and \eqref{leadingtermnashexpansion} that we can write: 
$$w^p=\bar w^p+\mathring w^p=w^{p,p}+(\bar w^p-w^{p,p})+\mathring w^p=\sum_I\left[w^{p,o}_I+w^{p,c}_I\right]+(\bar w^p-w^{p,p})+\mathring w^p.$$ 
We now want to compute the lowest implicit constant $\bar C_0$ we can allow. From Lemma \ref{leadingtermsnash} we get a constant $C_N$ independent of $C_0$ and from the same Lemma together with \ref{estimatesfieldscircp}, \ref{estimatesfieldsp} we get a possibly different constant $C$ which might depend on $C_0$, such that
\begin{equation}\label{convergence}
    \begin{split}
        ||\partial_t^jw^p||_r&\leq \sup_I\left[||\partial_t^jw^{p,o}_I||_r+||\partial_t^jw^{p,c}_I||_r\right]+||\partial_t^j(\bar w^p-w^{p,p})||_r+||\mathring w^p||_r\\
        &\leq C_N\lambda_{q+1}^{r+j}\delta_{q+1}^{1/2}+C\lambda_{q+1}^{r+j}[(\tau^a/\tau^c)\delta_{q+1}^{1/2}+\mathcal{T}_p\delta_q^{1/2}]\\
        &\leq C_N\lambda_{q+1}^{r+j}\delta_{q+1}^{1/2}+2C\lambda_{q+1}^{r+j}(\tau^a/\tau^c)\delta_{q+1}^{1/2}\\
    \end{split}
\end{equation}
for $j=0,1,2$ and $0\leq r \leq N-j$. Given any $b, \gamma_a$ we can pick $a$ sufficiently large that:
$$2C(\tau^a/\tau^c)=2C\left(\frac{\lambda_q}{\lambda_{q+1}}\right)^{\gamma_a}\leq C_N$$
It follows that for any $\bar C_0$ satisfying: 
\begin{equation}\label{barC01}
    \bar C_0\geq 4C_N,
\end{equation}
we have:
$$||\partial_t^j w^p||_r\leq \lambda_{q+1}^{r+j}\frac{\bar C_0}{2}\delta_{q+1}^{1/2} \ \text{ for } \ j=0,1,2 \ \text{ and } \ 0\leq r \leq N-j$$
we will fix the actual value of $\bar C_0$ afterwards. Using, in addition, the bounds in Lemmas \ref{estimatesfieldsringg}, \ref{estimatesfieldsbarg}, we obtain for some possibly different $C$:
\begin{equation}\label{Crperturbation}
    \begin{split}
        ||\partial_t^j(v_{q+1}-v_q)||_r&\leq ||\partial_t^jw^p||_r+||\partial_t^jw^g||_r\\
        &\leq \lambda_{q+1}^{r+j}\frac{\bar C_0}{2}\delta_{q+1}^{1/2}+\lambda_{q+1}^r\lambda_q^jC\left[\lambda_q\ell^{-\alpha}\tau^a\delta_{q+1}+\mathcal{T}_g(\ell\lambda_q)^{m_0}\delta_q^{1/2}\right]\\
        &\leq \bar C_0 \lambda_{q+1}^{r+j}\delta_{q+1}^{1/2}
    \end{split}
\end{equation}
for $0\leq r\leq N-j$. Here we choose $a$ sufficiently large so that the last inequality holds, given any $b, \ \gamma_a, \ \beta$, this is always possible as:
$$C\frac{\lambda_q^{j+1}\ell^{-\alpha}\tau^a\delta_{q+1}}{\lambda_{q+1}^j\delta_{q+1}^{1/2}}\leq C(\tau^a/\tau^c)(\delta_{q+1}/\delta_q)^{1/2}=C\left(\frac{\lambda_q}{\lambda_{q+1}}\right)^{\beta+\gamma_a}$$
and
$$C\frac{\lambda_q^j\mathcal{T}_g(\ell\lambda_q)^{m_0}\delta_q^{1/2}}{\lambda_{q+1}^j\delta_{q+1}^{1/2}}\leq C\left(\frac{\lambda_q}{\lambda_{q+1}}\right)^{1+\beta+\gamma_a+\gamma_{CZ}}$$
we will use this type of argument several times in what follows, without repeating the full details.

\noindent The additional smallness we have in $b^g, \ b^p$ allows us to use the same constant from above, upon choosing $a$ sufficiently large:
\begin{equation*}
    \begin{split}
        ||\partial_t^j(B_{q+1}-B_q)||_r&\leq ||\partial_t^jb^p||_r+||\partial_t^jb^g||_r\\
        &\leq C\lambda_{q+1}^{r+j}(\tau^a/\tau^c)\delta_{q+1}^{1/2}+C\lambda_{q+1}^r\lambda_q^j\left[\lambda_q\ell^{-\alpha}\tau^a\delta_{q+1}+\mathcal{T}_g(\ell\lambda_q)^{m_0}\delta_q^{1/2}\right]\\
        &\leq \bar C_0 \lambda_{q+1}^{r+j}\delta_{q+1}^{1/2}.
    \end{split}
\end{equation*}

\noindent \textbf{Full $||\cdot||_0$ estimates.} We now move to the non-vanishing condition for the magnetic field, from Proposition \ref{recap} and Lemmas \ref{estimatesfieldscircp}, \ref{estimatesfieldsp}, we deduce:
\begin{equation*}
    \begin{split}
        |B_{q+1}|&=|\tilde B_q+b^p|\\
        &\geq |\tilde B_q|-|b^p|\\
        &\geq c_0(1+2\delta_{q+1}^{1/2})-C(\tau^a/\tau^c)\delta_{q+1}^{1/2}\\
        &\geq c_0(1+\delta_{q+1}^{1/2})
    \end{split}
\end{equation*}
where we used the definition of $\mathcal{T}_p$ in \eqref{Tp} and where $C$ comes from the bounds in Lemma \ref{estimatesfieldscircp}, \ref{estimatesfieldsp}, and the last inequality can be ensured by choosing $a$ sufficiently large. 

\noindent For $\bar C_0$ as in \eqref{barC01} and the same Lemmas, we deduce:
\begin{equation*}
    \begin{split}
        ||v_{q+1}||_0&\leq ||\tilde v_{q}||_0+||w^p||_0\\
        &\leq C_0(1-2\delta_{q+1}^{1/2})+\frac{\bar C_0}{2}\delta_{q+1}^{1/2}\\
        &\leq C_0(1-\delta_{q+1}^{1/2})
    \end{split}
\end{equation*}
which can be ensured for any $C_0\geq \bar C_0/2$. One can argue similarly for $||B_{q+1}||_0$.

\noindent \textbf{Full $||\cdot||_r$ estimate.} From the Iterative Assumptions \eqref{inductiveassumptionsgeneral} and \eqref{Crperturbation} it follows that:
\begin{equation*}
    \begin{split}
        ||\partial_t^jv_{q+1}||_r&\leq ||\partial_t^jv_{q}||_r+||\partial_t^j(v_{q+1}-v_q)||_r\\
        &\leq C_0\lambda_q^{r+j}\delta_q^{1/2}+\lambda_{q+1}^{r+j}\frac{\bar C_0}{2}\delta_{q+1}^{1/2}\\
        &\leq C_0\lambda_{q+1}^{r+j}\delta_{q+1}^{1/2}
    \end{split}
\end{equation*}
where to ensure that the last inequality holds, we choose $a$ sufficiently large so that:
$$C_0\lambda_q\delta_q^{1/2}\leq (C_0-\bar C_0)\lambda_{q+1}\delta_{q+1}^{1/2},$$
which we can always find if we fix the constraint $C_0>\bar C_0$. One can argue similarly for $B_{q+1}$.

\noindent \textbf{Alfv\'en Transport Bounds.} From Proposition \ref{recap} and Lemmas \ref{estimatesfieldscircp}, \ref{estimatesfieldsp}, and the bounds above, we deduce:
\begin{equation*}
    \begin{split}
        ||\mathcal{A}^\pm_{q+1}v_{q+1}||_r&\leq ||\mathcal{\tilde A}^\pm w^p||_r+||(w^p\pm b^p)\cn w^p||_r+||\mathcal{\tilde A}^\pm \tilde v_q||_r+||(w^p\pm b^p)\cn \tilde v_q||_r\\
        &\leq C\lambda_{q+1}^r1/\tau^a\delta_{q+1}^{1/2}+C_N^2\lambda_{q+1}^{r+1}\delta_{q+1}+C\lambda_{q+1}^r\lambda_q\delta_q\\
        &+C C_N\lambda_{q+1}^r\lambda_q\delta_q^{1/2}\delta_{q+1}^{1/2}\\
        &\leq 2C_N^2 \lambda_{q+1}^{r+1}\delta_{q+1}^{1/2}.
    \end{split}
\end{equation*}
Abusing notation, the constants $C$ are possibly different and may depend on $C_0$ and the other parameters, while $C_N$ is the same as in \eqref{convergence}. The last inequality follows upon choosing $a$ sufficiently large, given all the other parameters.

\noindent We can now set: 
$$\bar C_0=4C_N^2$$
which still ensures the inequality in \eqref{barC01}, as without loss of generality $C_N\geq 1$, and conclude that:
$$||\mathcal{A}^\pm_{q+1}v_{q+1}||_r\leq \bar C_0\lambda_{q+1}^{r+1}\delta_{q+1}^{1/2}\leq  C_0\lambda_{q+1}^{r+1}\delta_{q+1}^{1/2},$$
this suffices as one can argue similarly for $\mathcal{A}^\pm_{q+1}B_{q+1}$.

\noindent \textbf{Support estimates - Vector Fields.} Eventually adding the correct index $p,g$ we see that the Lagrangian perturbations defined in \eqref{lagrangiang}, \eqref{lagrangianp} satisfy:
\begin{equation*}
    X_s(x,t)-\IId=\int_0^s\xi(X_s(x,t),t)\dd s.
\end{equation*}
In particular, $X^p=X_1^p, \ X^g=X^g_1$ are the identity map $\IId:\mathbb{T}^3\to \mathbb{T}^3$ on the complement of $\supp _t \xi^p,\ \supp _t \xi^g$ respectively. From this and the definitions \eqref{split1g}, \eqref{split1p}, namely
\begin{equation*}
    \begin{split}
        \tilde v_q&=\partial_t X^g\circ (X^g)^{-1}+(X^g)_*v_q,\\
        v_{q+1}&=\partial_t X^p\circ (X^p)^{-1}+(X^p)_*\tilde v_q,\\
    \end{split}
\end{equation*}
we deduce that:
\begin{equation}\label{supportfinal1}
    \supp_t \ (v_{q+1}-v_q)\subset \supp_t \ \xi^p\cup \supp_t \ \xi^g.
\end{equation}
Now, from the support estimates for $\tilde a_I$ in Lemma \ref{slowcoeff} and the Ansatz in \eqref{ansatz} we have:
\begin{equation}\label{supportfinal2}
    \supp_t \ \xi^p\subset \bigcup_I \supp_t \ \tilde a_I \subset B_{C\ell_t}(\supp_t R_\ell),
\end{equation}
moreover, from the construction in Subsection \ref{gnldf} we have: 
$$\xi^g=\sum_j\xi_j^g=\sum_j\curl \ \tilde \eta_j \Theta_j^g$$ 
and 
$$\Theta_j^g\equiv0\iff\sum_{I:I_t=j}f_{p(I)}\mathbb{T} A_I\equiv0 \ \text{ on }\  B_{\tau^c}(t_j)\iff \supp_t R_q \cap B_{\tau^c}(t_j)=\emptyset$$ 
where we used the support properties of $A_I$ from Lemma \ref{slowcoeffgalb}. We deduce that:
\begin{equation}\label{supportfinal3}
    \supp_t \ \xi^g \subset B_{2\tau^c}(\supp_t R_q).
\end{equation}
We now plug \eqref{supportfinal2}, \eqref{supportfinal3} in \eqref{supportfinal1}
and we conclude by means of the Iterative Assumptions \eqref{inductiveassumptionsgeneral} that:
\begin{equation}\label{finalsupport}
    \begin{split}
        \supp_t \ (v_{q+1}-v_q)&\subset B_{C\ell_t}(\supp_t R_\ell)\cup B_{2\tau^c}(\supp_t R_\ell)\\
        &\subset B_{2\tau^c}(\supp_t R_\ell)\\
        &\subset ((1+\delta_q^{1/2})1/2-2\tau^c,(1-\delta_q^{1/2})5/2+2\tau^c)\\
        &\subset ((1+\delta_{q+1}^{1/2})1/2,(1-\delta_{q+1}^{1/2})5/2)
    \end{split}
\end{equation}
where the inclusions follow by taking $a$ sufficiently large, given the definitions of $\tau^c$ and $\ell_t=\tau^a$ in \eqref{tauc}, \eqref{taua} respectively. In particular,
$$\supp_t \ (v_{q+1}-v_q)\subset (1/2,5/2)$$
as claimed. A similar argument for $B_{q+1}$ leads to: 
$$\supp_t \ (B_{q+1}-B_q)\subset ((1+\delta_{q+1}^{1/2})1/2,(1-\delta_{q+1}^{1/2})5/2)\subset (1/2,5/2).$$

\noindent \textbf{Pressure.} Note that $p_{q+1}=\tilde p_q$, the estimates are provided in Proposition \ref{recap}, and one can argue as above to deal with the implicit constant.

\noindent \textbf{Reynolds Stress.} From Proposition \ref{recap} and Lemmas \ref{Rcrt}, \ref{Rptr}, \ref{Rpna}, \ref{Rprm}, \ref{Rpmo}, \ref{Rpqua}, we deduce:
\begin{equation*}
    \begin{split}
        ||R_{q+1}||_r&\leq ||R^{p,tr}||_r+||R^{p,na}||_r+||R^{p,qua}||_r+||R^g||_r+||R^{p,crt}||_r+||R^{p,mo}||_r+||R^{p,rm}||_r\\
        &\lesssim \lambda_{q+1}^r(\tau^c/\tau^a)\frac{\lambda_q\delta_q^{1/2}\delta_{q+1}^{1/2}}{\lambda_{q+1}^{1-2\alpha}}+\lambda_{q+1}^{r}\frac{\lambda_q\delta_q^{1/2}\delta_{q+1}^{1/2}}{\lambda_{q+1}^{1-\alpha}}+\lambda_{q+1}^r\frac{\lambda_q\delta_{q+1}}{\lambda_{q+1}^{1-\alpha}}+\lambda_{q+1}^r(\tau^a/\tau^c)\ell^{-\alpha}\delta_{q+1}\\
        &\lesssim \lambda_{q+1}^{r+2\alpha}\left[(\tau^c/\tau^a)\frac{\lambda_q\delta_q^{1/2}\delta_{q+1}^{1/2}}{\lambda_{q+1}}+(\tau^a/\tau^c)\delta_{q+1}\right]\\
        &\lesssim \lambda_{q+1}^{r+2\alpha}\sqrt{\frac{\lambda_q\delta_q^{1/2}\delta_{q+1}^{3/2}}{\lambda_{q+1}}}\\
        &\leq \lambda_{q+1}^{r-2\alpha}\delta_{q+2}
    \end{split}
\end{equation*}
for $0\leq r \leq M$, where the last two inequalities are guaranteed by the optimisation procedure used to find the value of $\gamma_a$, see \eqref{constraintR1} and our choice of parameters in \ref{choiceofparameters}, see \eqref{constraintR2}. 

\noindent From the same Proposition and Lemmas, and the bounds above, we also obtain:
\begin{equation*}
    \begin{split}
        ||\mathcal{A}^\pm_{q+1}R_{q+1}||_r&\leq ||\mathcal{\tilde A}^\pm R^{p,tr}||_r+||\mathcal{\tilde A}^\pm R^{p,na}||_r+||\mathcal{\tilde A}^\pm R^{p,qua}||_r+||\mathcal{\tilde A}^\pm R^g||_r+||\mathcal{\tilde A}^\pm R^{p,crt}||_r+||\mathcal{\tilde A}^\pm R^{p,mo}||_r\\
        &+||\mathcal{\tilde A}^\pm R^{p,rm}||_r +||(w^p\pm b^p)\cn R_{q+1}||_r\\
        &\lesssim \lambda_{q+1}^r(1/\tau^a)(\tau^c/\tau^a)\frac{\lambda_q\delta_q^{1/2}\delta_{q+1}^{1/2}}{\lambda_{q+1}^{1-2\alpha}}+\lambda_{q+1}^{r}(1/\tau^a)\frac{\lambda_q\delta_q^{1/2}\delta_{q+1}^{1/2}}{\lambda_{q+1}^{1-2\alpha}}+\lambda_{q+1}^{r+\alpha}1/\tau^a\delta_{q+1}\\
        &+\lambda_{q+1}^r1/\tau^c\ell^{-\alpha}\delta_{q+1}+\lambda_{q+1}^{r+1-2\alpha}\delta_{q+1}^{1/2}\delta_{q+2}\\
        &\lesssim \lambda_{q+1}^{r+2\alpha}(1/\tau^a)\left[(\tau^c/\tau^a)\frac{\lambda_q\delta_q^{1/2}\delta_{q+1}^{1/2}}{\lambda_{q+1}}+(\tau^a/\tau^c)\delta_{q+1}\right]+\lambda_{q+1}^{r+1-2\alpha}\delta_{q+1}^{1/2}\delta_{q+2}\\
        &\lesssim \lambda_{q+1}^{r+1-2\alpha}\delta_{q+1}^{1/2}\delta_{q+2}
    \end{split}
\end{equation*}
for $0\leq r \leq M-1$, where we used that according to our choice of parameters in \ref{choiceofparameters}, see \eqref{constraintadmissibility} and \eqref{taua}, we have: 
$$1-(\beta+\gamma_a+\gamma_{CZ})\geq 0\Longrightarrow 1/\tau^a\leq \lambda_{q+1}\delta_{q+1}^{1/2}$$
from which using \eqref{constraintR1} followed by \eqref{constraintR2} we deduce: 
$$\lambda_{q+1}^{r+2\alpha}(1/\tau^a)\left[(\tau^c/\tau^a)\frac{\lambda_q\delta_q^{1/2}\delta_{q+1}^{1/2}}{\lambda_{q+1}}+(\tau^a/\tau^c)\delta_{q+1}\right]\leq 2\lambda_{q+1}^{r+1-2\alpha}\delta_{q+1}^{1/2}\delta_{q+2}$$
and 
$$\lambda_{q+1}^{r+\alpha}1/\tau^a\delta_{q+1}=\lambda_{q+1}^{r+1+2\alpha}\delta_{q+1}^{1/2}\tau^c/\tau^a\frac{\lambda_q\delta_{q}^{1/2}\delta_{q+1}^{1/2}}{\lambda_{q+1}}\leq2 \lambda_{q+1}^{r+1-2\alpha}\delta_{q+1}^{1/2}\delta_{q+2}.$$
Now let $C$ be the largest of the implicit constants in the bounds for $R_{q+1}, \ \mathcal{A}^\pm_{q+1}R_{q+1}$ above for $0\leq r\leq M$ and $0\leq r\leq M-1$, respectively.  Since $C$ depends on all the parameters but not on $a$, for fixed $b, \ \gamma_{CZ}$, upon choosing $a$ sufficiently large we can ensure: 
$$C\lambda_{q+1}^{-\alpha}=C\left(\frac{\lambda_{q}}{\lambda_{q+1}}\right)^{\gamma_{CZ}}\leq 1$$
and conclude:
\begin{equation*}
    \begin{split}
        ||R_{q+1}||_r&\leq\lambda_{q+1}^{r-\alpha}\delta_{q+2} \ \text{ for } 0\leq r \leq M,\\
        ||\mathcal{A}^\pm_{q+1}R_{q+1}||_r&\leq \lambda_{q+1}^{r+1-\alpha}\delta_{q+1}^{1/2}\delta_{q+2} \ \text{ for } 0\leq r \leq M-1.
    \end{split}
\end{equation*}

\noindent \textbf{Support estimates - Reynolds Stress.} Inspection of the proofs of the Lemmas in Subsections \ref{Rgs}, \ref{Rps} shows that: 
$$\supp_t \ R_{q+1}\subset \bigcup_I\supp_t\ \Theta_I^p\cup \bigcup_j\supp_t\ \Theta_j^g\cup \supp_t R_q.$$
Now, arguing as in the support estimate for the vector fields, see \eqref{finalsupport}, we deduce:
$$\supp_t \ R_{q+1}\subset ((1+\delta_{q+1}^{1/2})1/2,(1-\delta_{q+1}^{1/2})5/2)$$
as claimed.


\subsection{Proof of the Main Theorem \ref{main}.} \label{proofmain} We first initialise the scheme by producing explicit subsolutions exhibiting the desired non-conservation properties and finally prove convergence of the iterative scheme.

\noindent Let $\eta,\ \chi:\mathbb{R}\to [0,1]$ be smooth non-negative functions such that: 
\begin{equation*}
        \eta(t)=\begin{cases}
        1 \ \text{ for } \  t \leq 2/3,\\
        0 \ \text{ for } \  t\geq 1 \\
    \end{cases} \ \ \text{ and } \ \
    \chi(t)=\begin{cases}
        1 \ \text{ for } \  |t| \leq 2/3,\\
        0 \ \text{ for } \  |t|\geq 1.\\
    \end{cases}
\end{equation*}
We now fix $1\leq t_0< t_1< 2-1/4$, with $t_1-t_0>1/4$, and define for $\lambda^\parallel_0>8$ arbitrarily large  but independent of the parameters:
\begin{equation}\label{f(t)}
    f(t)=\chi(\lambda^\parallel_0(t-t_0))-\chi(\lambda^\parallel_0(t-t_1)).
\end{equation}
To sample the types of subsolutions this scheme can handle, we consider the following three examples of $(v_0, \ B_0)$. In what follows, $X=X_t(x)$ denotes the Lagrangian flow of $v_0$ at the identity at $t=0$. 

\noindent \textbf{First Example.} The following example leads to a solution that does not conserve energy and which possesses gradients parallel to the magnetic field.
\begin{equation*}
    \begin{cases}
        \bar v(x)=\cos (\lambda_0\delta_0^{1/2}x\cdot e_2)e_1,\\
        \bar{\bar v}(x)=\cos(\lambda_0^\parallel x\cdot e_3)\cos ( \lambda_0\delta_0^{1/2}x\cdot e_2)e_1,\\
        \bar B=e_3
    \end{cases} \leadsto 
    \begin{cases}
        v_0=\eta \bar v+f\bar{\bar v},\\
        B_0=X_*\bar B.
    \end{cases}
\end{equation*}

\noindent \textbf{Second Example.} The following example leads to a solution that does not conserve energy and cross helicity.
\begin{equation*}
    \begin{cases}
        \bar v(x)=\cos (\lambda_0\delta_0^{1/2}x\cdot e_1)e_2,\\
         \bar B(x)=e_3+\cos (\lambda_0\delta_0^{1/2}x\cdot e_1)e_2
    \end{cases}\leadsto
    \begin{cases}
        v_0=\eta \bar v,\\
        B_0=X_*\bar B.
    \end{cases}
\end{equation*}

\noindent\textbf{Third Example.} The following example leads to a solution that does not conserve energy and has non-trivial magnetic helicity.
\begin{equation*}
    \begin{cases}
        \bar v(x)=\cos (\lambda_0\delta_0^{1/2}x\cdot e_3)e_1+\sin (\lambda_0\delta_0^{1/2}x\cdot e_3)e_2,\\
         \bar B(x)=\sin (\lambda_0\delta_0^{1/2}x\cdot e_3)e_1-\cos (\lambda_0\delta_0^{1/2}x\cdot e_3)e_2\\
    \end{cases}\leadsto
    \begin{cases}
        v_0=\eta \bar v,\\
        B_0=X_*\bar B.
    \end{cases}
\end{equation*}

\noindent Note that one should, in principle, take the integer part in the arguments of the functions above to make the objects well defined on $\mathbb{T}^3$; we omit this detail. We remark that the above setups cannot arise as reformulations of ideal MHD as a $2\frac{1}{2}$- Euler flow, either because of parallel gradients or the geometry of the vector fields themselves. We now check that these choices of initial data give rise to subsolutions satisfying the Iterative Assumptions \eqref{inductiveassumptionsgeneral} and the claimed non-conservation of energy/cross helicity. Combinations of the above are also possible. We give explicit computations for the first and second examples only.

\noindent \textit{Details of the first example.} Note that: 
$$e_3\cn \bar v=0 \ \text{ and } \ (e_3\cn \bar{\bar v})(x)=-\lambda^\parallel\sin(\lambda_0^\parallel x_3)\cos (\lambda_0\delta_0^{1/2}x_1)e_2.$$
We will use $\bar v$ to show that the energy is not conserved, while  $\bar{\bar v}$ is there to show that one can accommodate for the presence of parallel gradients, which can be as large as the time we apply the vector field for, in the sense that:
$$f'\sim \lambda_0,$$
see the definition of $f(t)$ in \eqref{f(t)}. 

\noindent We can compute $X_t$ explicitly, indeed
$$X_t=\IId+\eta^{[1]}\bar v+f^{[1]}\bar{\bar v} \ \text{ and } \ X_t^{-1}=\IId-\eta^{[1]}\bar v-f^{[1]}\bar{\bar v}$$
where 
$$\eta^{[1]}(t)=\int_0^t\eta(s)\dd s \ \text{ and } \ f^{[1]}(t)=\int_0^tf(s)\dd s.$$
In particular, we have:
\begin{equation*}
    \begin{split}
        &f^{[1]}(t)=0 \ \text{ for } t\leq 1 \ \& \ t\geq 2 \ \text{ and } \ |f^{[1]}(t)|\leq 1/\lambda_0^\parallel,\\
        &\eta^{[1]}(t)=-t \ \text{ for } t\leq 0 \ \text{ and } \eta^{[1]}(t)=\text{const} \ \text{ for } t\geq 1.
    \end{split}
\end{equation*}
By construction, $B_0$ solves: 
$$\partial_t B_0+\curl[B_0\times v_0]=0, \ \ddiv \ B_0=0.$$
Moreover, we can explicitly compute:
\begin{equation*}
    \begin{split}
        \DD X_t(x)&=\IId-\lambda_0\delta_0^{1/2}\sin (\lambda_0\delta_0^{1/2}x_2)\eta^{[1]}(t)e_1\otimes e_2\\
        &-\lambda_0\delta_0^{1/2}\cos(\lambda_0^\parallel x_3)\sin ( \lambda_0\delta_0^{1/2}x_2)f^{[1]}(t)e_1\otimes e_2-\lambda_0^\parallel\sin(\lambda_0^\parallel x_3)\cos ( \lambda_0\delta_0^{1/2}x_2)f^{[1]}(t)e_1\otimes e_3,\\
        B_0(x,t)&=\DD X_t\circ X_t^{-1}(x)[e_3]=e_3-\sin(\lambda_0^\parallel x_3)\cos ( \lambda_0\delta_0^{1/2}x_2)\lambda_0^\parallel f^{[1]}(t)e_1.
    \end{split}
\end{equation*}
By inspection, we deduce that as soon as $a$ is chosen sufficiently large to have: 
$$\lambda_0^\parallel,\ \sup_t |\eta'|\leq \lambda_0\delta_0, \ \sup_t |\eta''|\leq \lambda_0^2, \ \delta_0^{1/2}<1/2,$$ 
we can ensure that:
\begin{equation}\label{r0}
    \begin{split}
        ||B_0||_0, \ ||v_0||_0&\leq C_0(1-\delta_q^{1/2}),\\
        |B_0|&\geq c_0(1+\delta_0^{1/2}),\\
        ||\partial_t^jB_0||_r,||\partial_t^j v_0||_r&\leq C_0\lambda_0^{r+j}\delta_0^{1/2} \ \text{ for } \ j=0,1,2, \ \& \ \ 0\leq r \leq N-j\ \ \& \ \  (r,j)\neq (0,0),\\
        ||\mathcal{A}^\pm_0 B_0||_r, \ ||\mathcal{A}^\pm_0 v_0||_r&\leq C_0\lambda_0^{r+1}\delta_0 \text{ for } 0\leq r \leq M-1,\\
    \end{split}
\end{equation}
for some $C_0>\bar{C}_0$, a constant sufficiently large, and $c_0<1/2$, another constant, sufficiently small, both independent of the parameters.  

\noindent Plugging $(v_0, B_0)$ as above in the ideal MHD system \eqref{MHD}, we see that the errors can be reabsorbed in a Reynolds stress $R_0$ as in \eqref{relaxedMHD} of the form: 
\begin{equation*}
    \begin{split}
        R_0&=f'\frac{1}{\lambda_0 \delta_0^{1/2}}\cos(\lambda_0^\parallel x_3)\sin(\lambda_0\delta_0^{1/2}x_2)(e_1\otimes e_2+e_2\otimes e_1)+\eta'\frac{1}{\lambda_0\delta_0^{1/2}}\sin(\lambda_0\delta_q^{1/2}x_2)(e_1\otimes e_2+e_2\otimes e_1)\\
        &+\frac{\lambda^\parallel_0 \lambda_0^\parallel f^{[1]}(t)}{\lambda_0 \delta_0^{1/2}}\cos(\lambda_0^\parallel x_3)\sin(\lambda_0\delta_0^{1/2}x_2)(e_1\otimes e_2+e_2\otimes e_1)
    \end{split}
\end{equation*}
where the first row comes from the time derivative of the cut-offs in $v_0$ while the last one comes from $-B_0\cn B_0$. We remark that:
$$|f'|\lesssim \lambda_0^\parallel, \ |f''|\lesssim (\lambda_0^\parallel)^2 \ \text{ and } \ \lambda_0^\parallel |f^{[1]}(t)|\lesssim 1,$$ 
moreover, our choice of $b, \ \gamma_{CZ}$ from Subsection \ref{choiceofparameters}, see \eqref{misc1}, guarantees that 
$$\beta(2b+1)<1-2\gamma_{CZ}$$
and we deduce:
\begin{equation*}
\begin{split}
    ||R_0||_{r}&\leq C (\lambda^\parallel_0+\sup_t |\eta'|)\lambda_0^{r-2\alpha}\delta_1 \ \text{ for } \ 0\leq r \leq M,\\
    ||\mathcal{A}^{\pm}R_0||_{r}&\leq C ((\lambda^\parallel_0)^2+\sup_t |\eta''|)\lambda_0^{r-2\alpha}\delta_1 \ \text{ for } \ 0\leq r \leq M-1,
\end{split}
\end{equation*}
for some constant $C$ that may depend on $M$ but not on $r$. We can now choose $a$ so large that: 
$$C((\lambda^\parallel_0)^2+\sup_t |\eta''|)\lambda_0^{-\alpha}\leq \lambda_0\delta_0^{1/2} \ \text{ and } \ C(\lambda^\parallel_0+\sup_t |\eta''|)\lambda_0^{-\alpha}\leq 1$$
and we conclude:
\begin{equation}\label{r1}
\begin{split}
    ||R_0||_{r}&\leq \lambda_q^{r-\alpha}\delta_1 \ \text{ for } \ 0\leq r \leq M,\\
    ||\mathcal{A}^{\pm}R_0||_{r}&\leq \lambda_q^{r+1-\alpha}\delta_0^{1/2}\delta_1 \ \text{ for } \ 0\leq r \leq M -1.
\end{split}
\end{equation}

\noindent With $p_0=0$, we have that $(v_0, \ B_0,\ p_0, \ R_0)$ solves \eqref{relaxedMHD} and
\begin{equation}\label{r2}
    v_0=\begin{cases}
            \bar v \ \text{ for }\ t\leq 2/3\\
            0 \ \text{ for }\ t\geq 2\\
        \end{cases}, \ B_0=e_3  \ \text{ for } t\leq 1-1/4 \ \& \ t\geq 2, \
        R_0=0 \ \text{ for }\ t\leq 2/3 \ \& \ t\geq 2.
\end{equation}
From this, we deduce that $(v_0, \ B_0,\ p_0)$ is an exact solution in $t\leq 2/3 \ \& \ t\geq 2$ and 
$$\mathcal{E}_{(v_0,B_0)}(0)>\mathcal{E}_{(v_0,B_0)}(3).$$ 

\noindent \textit{Details of the second example.} As above we can compute $X_t$ explicitly, indeed
$$X_t=\IId+\eta^{[1]}\bar v \ \ \text{ and }\ \ X_t^{-1}=\IId-\eta^{[1]}\bar v.$$
By construction $B_0$ is a solution of \ref{FH} and
\begin{equation*}
    \begin{split}
        \DD X_t(x)&=\IId-\lambda_0\delta_0^{1/2}\sin (\lambda_0\delta_0^{1/2}x_1)\eta^{[1]}(t)e_2\otimes e_1\\
        B_0(x,t)&=\DD X_t\circ X_t^{-1}(x)[\bar B(X_t^{-1}(x))],\\
        &=\left[\IId-\lambda_0\delta_0^{1/2}\sin (\lambda_0\delta_0^{1/2}X_t^{-1}(x)\cdot e_1)\eta^{[1]}(t)e_2\otimes e_1\right][e_3+\cos (\lambda_0\delta_0^{1/2}X_t^{-1}(x)\cdot e_1)e_2]\\
        &=\left[\IId-\lambda_0\delta_0^{1/2}\sin (\lambda_0\delta_0^{1/2}x_1)\eta^{[1]}(t)e_2\otimes e_1\right][e_3+\cos (\lambda_0\delta_0^{1/2}x_1)e_2]\\
        &=e_3+\cos (\lambda_0\delta_0^{1/2}x_1)e_2,
    \end{split}
\end{equation*}
moreover, $$(\partial_t+v_0\cn)B_0=B_0\cn B_0=B_0\cn v_0=0$$
and
\begin{equation*}
    \begin{split}
        (\partial_t+v_0\cn)v_0(x,t)&=\eta'(t)\cos (\lambda_0\delta_0^{1/2}x_1)e_2\\
        &=\ddiv\left[\eta'(t)\frac{1}{\lambda_0\delta_0^{1/2}}\sin(\lambda_0\delta_0^{1/2}x_1)(e_1\otimes e_2+e_2\otimes e_1)\right]\\
        &=\ddiv \ R_0(x,t),
    \end{split}
\end{equation*}
where we defined: 
$$R_0(x,t)=\eta'(t)\frac{1}{\lambda_0\delta_0^{1/2}}\sin(\lambda_0\delta_0^{1/2}x_1)(e_1\otimes e_2+e_2\otimes e_1).$$
From the above, we deduce that, with $p_0=0$, the quadruple $(v_0, \ B_0,\ p_0, \ R_0)$ solves \eqref{FH} and
\begin{equation}\label{r22}
    v_0=\begin{cases}
            \bar v \ \text{ for }\ t\leq 2/3\\
            0 \ \text{ for }\ t\geq 2\\
        \end{cases}, \ B_0=\bar B  , \
        R_0=0 \ \text{ for }\ t\leq 2/3 \ \& \ t\geq 2,
\end{equation}
that is, $(v_0, \ B_0,\ p_0)$ is an exact solution in $t\leq 2/3 \ \& \ t\geq 2$. Moreover, a simple computation shows that the total energy and cross helicity are not conserved: 
$$\mathcal{E}_{(v_0,B_0)}(0)>\mathcal{E}_{(v_0,B_0)}(3) \ \text{ and } \ \mathcal{H}^\times_{(v_0,B_0)}(0)>\mathcal{H}^\times_{(v_0,B_0)}(3).$$ 
We can now argue as above and conclude that, also in this example, \eqref{r0} and \eqref{r1} are satisfied. 

\noindent \textbf{Convergence of the scheme.} With \eqref{r0}, \eqref{r1}, \eqref{r2}, \eqref{r22} at hand and given $0\leq \gamma<1/5$ we can apply the Iterative Proposition \ref{iterative} with $\beta$ satisfying $0\leq \gamma <\beta <1/5$ and obtain a sequence of sub-solutions $(v_q, \ B_q,\ p_q, \ R_q)$. From \eqref{iterativeeq}, and \eqref{iterativeRq}, we deduce:
\begin{equation*}
    \begin{split}
        ||v_{q+1}-v_q||_\gamma&\lesssim ||v_{q+1}-v_q||_0^{1-\gamma}||v_{q+1}-v_q||_1^\gamma\lesssim \delta_{q+1}^{1/2}\lambda_{q+1}^\gamma\lesssim \lambda_{q+1}^{\gamma-\beta}.
    \end{split}
\end{equation*}
Similarly, 
$$||B_{q+1}-B_q||_\gamma\lesssim\lambda_{q+1}^{\gamma-\beta} \ \ \text{ and } \ \ 
        ||R_q||_\gamma\lesssim\lambda_{q+1}^{\gamma-2\beta}.$$
The additional use of the Poisson problem determining $p_q$, namely
$$\Delta p_q=\ddiv\ddiv \left[R_q-v_q\otimes v_q+B_q\otimes B_q\right],$$
and standard Schauder estimates allow us to conclude that $(v_q, \ B_q,\ p_q, \ R_q)$ is a Cauchy sequence in $C_tC_x^\gamma$, in particular: 
$$(v_q, \ B_q,\ p_q, \ R_q)\overset{C_tC_x^\gamma}{\longrightarrow} (v, \ B, \ p, \ 0)$$ 
and passing to the limit in the distributional formulation of \ref{relaxedMHD} we see that $(v, B, p)$ solves \ref{MHD} weakly. Moreover, from \eqref{iterativeeq} it follows that
$$\sup_t [v-v_0], \ \sup_t [B-B_0]\subset (1/2,5/2)$$
and we deduce:
$$\mathcal{E}_{(v,B)}(0)=\mathcal{E}_{(v_0,B_0)}(0)>\mathcal{E}_{(v_0,B_0)}(3)=\mathcal{E}_{(v,B)}(3),$$
which shows that the energy is not conserved. Additionally, in the second example, we have:
$$\mathcal{H}^\times_{(v,B)}(0)=\mathcal{H}^\times_{(v_0,B_0)}(0)>\mathcal{H}^\times_{(v_0,B_0)}(3)=\mathcal{H}^\times_{(v,B)}(3),$$
which shows that the cross-helicity is not conserved.


\section{Lie-Taylor Series Expansion} \label{lt}
Here we develop the theory of perturbations led by a Lagrangian Displacement Field (LDF) $\xi$. This object appears in a different form and with a different goal compared to our work in Newcomb \cite{Newcomb} and Vladimirov, Moffat, Ilin \cite{VMI}. In these works, the authors study the stability of constrained energy minimisers and use this object to perform energy variations. Being interested in linear stability theory, the knowledge of the leading term of the perturbation suffices in that case, but not in ours.

\subsection{Perturbation and Lie-Taylor Series}
We study smooth perturbations $(v,B)$ of a smooth solution $(\bar v, \bar B)$ of the Faraday-Ohm system on $\mathbb{\mathbb{T}}^3\times \mathbb{R}$, namely
\begin{equation}\label{fo}
    \begin{cases}
        \partial_t \bar B+\curl[\bar B \times \bar v]=0,\\
        \bar B|_{t=t_0}=B_0,\\
        \ddiv \ \bar{B}=\ddiv \ \bar{v}=0.
    \end{cases}
\end{equation}
We recall that for sufficiently regular $\bar v, B_0$, this is equivalent to:
$$\bar B(x,t)=(\bar X_t)_*B_0(x), \ \ddiv \ B_0=0$$
see \eqref{duhamel}  with the correct identifications and \eqref{DG0}, where $\bar X$ is the Lagrangian flow map of $\bar v$, that is, a solution of:
\begin{equation*}
    \begin{cases}
        \partial_t \bar X_t(x)=\bar{v}(\bar X_t(x),t),\\
        \bar X|_{t=0}=\IId.
    \end{cases}
\end{equation*}
We require the perturbation $(v, B)$ to also be a solution of \eqref{fo}. A natural way to construct it is then to perturb the flow of $\bar v$ via some volume-preserving time-dependent diffeomorphism $X(\cdot,t):\mathbb{T}^3\to \mathbb{T}^3$ close to the identity and such that 
\begin{equation}\label{icX}
    X(\cdot,0)=\IId
\end{equation}
namely
$$\bar X \leadsto X\circ \bar{X},$$
we then set $v$ to be the divergence-free vector field having $X\circ \bar X$ as its Lagrangian flow map:
$$v=\partial_t(X\circ \bar X)\circ(X\circ \bar X)^{-1},$$
where the composition happens in space, that is, as maps $\mathbb{T}^3\to \mathbb{T}^3$. The associated perturbed magnetic field $B$ having to solve \eqref{fo}, will then be:
\begin{equation}\label{newb}
    B=(X\circ \bar{X})_*B_0=X_*\bar{X}_*B_0=X_* \bar{B}.
\end{equation}
Note that \eqref{icX} is there only to ensure the initial condition for the magnetic field is preserved:
$$B|_{t=0}=B_0$$ 
and one might want to remove it in general, moreover, $v$ has a pushforward structure as well, indeed:
\begin{equation}\label{newv}
    \begin{split}
        v&=\partial_t(X\circ \bar X)\circ(X\circ \bar X)^{-1}\\
            &=\partial_tX(\bar X\circ (X\circ \bar X)^{-1}) +\DD X(\bar X\circ(X\circ \bar X)^{-1})\partial_t\bar X\circ(X\circ \bar X)^{-1}\\
            &=\partial_tX\circ (X)^{-1}+\DD X(X^{-1})\bar{v}(X^{-1})\\
            &=\partial_tX\circ X^{-1}+X_*\bar{v}.
    \end{split}
\end{equation}
We rewrote the problem of constructing a perturbation of \eqref{fo} to that of constructing a time-dependent volume-preserving diffeomorphism $X(\cdot,t)$ close to the identity. Passing to the Lie Algebra of such a group, that is, divergence-free time-dependent vector fields, we can construct a continuous family indexed by $s\in[-1,1]$ of these as: 
$$X_s(\cdot,t)=\exp_{\IId}^{s\xi(\cdot,t)}:\mathbb{T}^3\to \mathbb{T}^3,$$
for $\xi(x,t):\mathbb{T}^3\times\mathbb{R}\to \mathbb{R}^3$ some divergence-free vector field, satisfying $\xi(\cdot,0)=0$ so that \eqref{icX} is retained. Explicitly, we solve:
\begin{equation*}
    \begin{cases}
        \partial_s X_s(x,t)=\xi(X_s(x,t),t),\\
        X|_{s=0}(x,t)=x
    \end{cases}
\end{equation*}
and think of $X$ as $X=X_1$ and $X^{-1}=X_{-1}$. We call the vector field $\xi$ the Lagrangian Displacement Field (LDF) after \cite{VMI}, where a similar object is named the first-order displacement field. This gives us a continuous family of divergence field vector fields $(v_s, B_s)$ starting at $(\bar{v},\bar{B})$ solving \eqref{fo}. The identities \eqref{newv} and \eqref{newb} applied with $X=X_s$ give us:
\begin{equation*}
    \begin{cases}
        v_s=\partial_tX_s\circ (X_s)^{-1}+(X_s)_*\bar{v},\\
        B_s=(X_s)_*\bar{B}.
    \end{cases}
\end{equation*}
We now study the derivative of $v_s$ with respect to the parameter $s$ we just introduced. We compute:
\begin{equation}\label{vs}
    \begin{split}
        \partial_s v_s&=\partial_s[\partial_tX_s\circ (X_s)^{-1}]+\partial_s[(X_s)_*\bar{v}]\\
        &=(\partial_s\partial_t X_s)\circ (X_s)^{-1}+\DD (\partial_tX_s)((X_s)^{-1})[\partial_s(X_s)^{-1}]-\mathcal{L}_\xi [(X_s)_*\bar{v}]\\
        &=(\partial_t\partial_s X_s)\circ (X_s)^{-1}-\DD (\partial_tX_s)((X_s)^{-1})\DD((X_s)^{-1})[\xi]-\mathcal{L}_\xi [(X_s)_*\bar{v}]\\
        &=\partial_t(\xi(X_s))\circ (X_s)^{-1}-\DD (\partial_tX_s\circ (X_s)^{-1})[\xi]-\mathcal{L}_\xi [(X_s)_*\bar{v}]\\
        &=\partial_t\xi+\DD\xi[\partial_tX_s\circ (X_s)^{-1}]-\DD (\partial_tX_s\circ (X_s)^{-1})[\xi]-\mathcal{L}_\xi [(X_s)_*\bar{v}]\\
        &=\partial_t \xi-\mathcal{L}_\xi(\partial_tX_s\circ (X_s)^{-1})-\mathcal{L}_\xi [(X_s)_*\bar{v}]\\
        &=\partial_t \xi-\mathcal{L}_\xi v_s\\
    \end{split}
\end{equation}
where we used the equations solved by $X_s, \ X^{-1}_s, \ (X_s)_*\bar{v}$, see \eqref{lagrangianvstransport} with the correct identifications and \eqref{duhamel}. Since $\xi$ is independent of $s$ we deduce that for any $k\geq 0$ that:
\begin{equation}\label{k+1}
    \begin{split}
        \partial_s^{k+1}v_s&=\partial_s^{k}(\partial_t\xi-\mathcal{L}_\xi v_s)\\
        &=-\partial_s^{k-1}\mathcal{L}_\xi \partial_sv_s\\
        &=(-1)^1\partial_s^{k-1}\mathcal{L}_\xi (\partial_t\xi-\mathcal{L}_\xi v_s)\\
        &=\dots\\
        &=(-1)^{k}\mathcal{L}_\xi^{k} (\partial_t\xi-\mathcal{L}_\xi v_s)
    \end{split}
\end{equation}
Since for $s=0$ we have $v_0=\bar{v}$, we can write a truncated Lie-Taylor series expansion in the parameter $s$ at $s=0$ of the form:
\begin{equation}\label{LDF1}
\begin{split}
    v_s&=\bar{v}+\sum_{k=1}^{k_0}\frac{s^k}{k!}\partial_s^kv_s|_{s=0}+r^{k_0,s}_w\\
    &=\bar{v}+s\left(\partial_t \xi-\mathcal{L}_\xi \bar{v}\right)-\frac{s^2}{2}\mathcal{L}_\xi\left(\partial_t \xi-\mathcal{L}_\xi \bar{v}\right)+\dots+r^{k_0,s}_w\\
    &=\bar v+\sum_{k=0}^{k_0}\frac{(-1)^ks^{k+1}}{(k+1)!}\mathcal{L}_{\xi}^k(\partial_t+\mathcal{L}_{\bar{v}})\xi\\
\end{split}
\end{equation}
where the remainder $r^{k_0,s}_w$ can be written as:
\begin{equation}\label{remainderimplicit}
    r^{k_0,s}_w=\frac{1}{k_0!}\int_0^s\partial_{s'}^{k_0+1}v_{s'}(s-s')^{k_0}\dd s'.
\end{equation}
We will study this in more detail in Subsection \ref{remainder}. Now recall that the magnetic field associated with $v_s$ is $B_s$ given by the pushforward formula:
$$B_s=(X_s\circ \bar X)_*B_0=(X_s)_*[(\bar X)_*B_0]=(X_s)_*\bar{B},$$
from \eqref{duhamel} or arguing as above we find expressions for $\partial_s^kB_s$ for any $k\geq 1$, namely
$$\partial_s^kB_s=(-1)^k\mathcal{L}_\xi^kB_s,$$
we can then write a truncated Lie-Taylor expansion for the parameter $s$ around zero of the form:
\begin{equation}\label{LDF2}
    \begin{split}
        B_s&=\bar B+\sum_{k=1}^{k_0}\frac{s^k}{k!}\partial_s^kB_s|_{s=0}+r^{k_0,s}_b\\
        &=\bar{B}-s\mathcal{L}_\xi \bar{B}+\frac{s^2}{2}\mathcal{L}_\xi^2\bar{B}-\dots+r^{k_0,s}_b\\
        &=\bar{B}+\sum_{k=0}^{k_0}\frac{(-1)^ks^{k+1}}{(k+1)!}\mathcal{L}_{\xi}^k\mathcal{L}_{\bar{B}}\xi\\
    \end{split}
\end{equation}
where, as before, the remainder $r^{k_0,s}_b$ has the explicit expression:
$$r^{k_0,s}_b=\frac{1}{k_0!}\int_0^s\partial_{s'}^{k_0+1}B_{s'}(s-s')^{k_0}\dd s'.$$

\subsection{Remainders}\label{remainder} 
In this section, we derive explicit expressions for the remainder terms in the Lie-Taylor series above. This will allow estimates of $v, B$ in terms of $\bar v, \bar B$ and $X$,  where the ones on $X$ follow from the ones on $\xi$ and Gr{\"o}nwall Lemma type of arguments.

\noindent\textit{Velocity Field:} We want to find a suitable representation for $r^{k_0,s}_w$. The problem here is that Taylor's Theorem gives \eqref{remainderimplicit} and we don't have an explicit expression or estimates for $v_s$.  We define $w$ to be: 
$$\partial_t\xi+\mathcal{L}_{\bar{v}}\xi=w$$ 
and compute  
\begin{equation*}
    \begin{split}
        \partial_s(v_s-\bar{v})&=\partial_s v_s\\
        &=\partial_t\xi-\mathcal{L}_{\xi}v_s\\
        &=\partial_t\xi-\mathcal{L}_{\xi}\bar{v}-\mathcal{L}_{\xi}(v_s-\bar{v})\\
        &=\underbrace{\partial_t\xi+\mathcal{L}_{\bar{v}}\xi}_{=w}-\mathcal{L}_{\xi}(v_s-\bar{v})\\
        &=w-\mathcal{L}_{\xi}(v_s-\bar{v})
    \end{split}
\end{equation*}
where we used \eqref{vs}. This shows that the difference $F=v_s-\bar{v}$ solves the forced Lie-transport equation:
\begin{equation*}
    \begin{cases}
        \partial_sF+\mathcal{L}_{\xi}F=w,\\
        F|_{s=0}=0.
    \end{cases}
\end{equation*}
Note that since $\xi$ is independent of $s$ we have $\partial_s\mathcal{L}_{\xi}=\mathcal{L}_{\xi}\partial_s$ and thus applying $\mathcal{L}_{\xi}^k$ to the equation above and commuting we get:
\begin{equation*}
    \begin{cases}
        \partial_s\mathcal{L}_{\xi}^kF+\mathcal{L}_{\xi}\mathcal{L}_{\xi}^kF=\mathcal{L}_{\xi}^kw\\
        \mathcal{L}_{\xi}^kF|_{s=0}=0
    \end{cases}
\end{equation*}
for any $k\geq 0$. In particular, denoting $\Phi_s=X_s^{-1}$  the inverse flow of $\xi$ we have an explicit representation of $\mathcal{L}_{\xi}^kF=\mathcal{L}_{\xi}^k(v_s-\bar{v}_q)$ as the solution of this forced Lie transport given by the Duhamel formula, see \eqref{duhamel}. Namely:
\begin{equation*}
    \begin{split}
        \mathcal{L}_{\xi}^k(v_s-\bar{v})&=\int_0^s\Phi_s^*X_{s'}^*(\mathcal{L}_{\xi}^kw)\dd s'\\
        &=\int_0^s(X_{s'}\circ\Phi_s)^*(\mathcal{L}_{\xi}^kw)\dd s'\\
        &=\int_0^s(X_{s'}\circ X_{-s})^*(\mathcal{L}_{\xi}^kw)\dd s'\\
        &=\int_0^sX_{s'-s}^*(\mathcal{L}_{\xi}^kw)\dd s'\\
    \end{split}
\end{equation*}
here we used that since $\xi$ does not depend on $s$ we have $\Phi_s=X_{-s}$ and $X_s\circ X_{s'}=X_{s+s'}$. With this and \eqref{k+1} at hand, we can rewrite the velocity remainder as follows:
\begin{equation}\label{LDF3}
    \begin{split}
        r^{k_0,s}_w&=\frac{1}{k_0!}\int_0^{s}\partial_{s'}^{k_0+1} v_{s'}(s-s')^{k_0}\dd s'\\
        &=\frac{(-1)^{k_0}}{k_0!}\int_0^{s}\mathcal{L}_{\xi}^{k_0}(\partial_t\xi -\mathcal{L}_{\xi}v_{s'})(s-s')^{k_0}\dd s'\\
        &=\frac{(-1)^{k_0}}{k_0!}\int_0^{s}\mathcal{L}_{\xi}^{k_0}(\partial_t\xi -\mathcal{L}_{\xi}\bar{v}-\mathcal{L}_{\xi}(v_{s'}-\bar{v}))(s-s')^{k_0}\dd s'\\
        &=\frac{(-1)^{k_0}}{k_0!}\int_0^{s}\mathcal{L}_{\xi}^{k_0}(w-\mathcal{L}_{\xi}(v_{s'}-\bar{v}))(s-s')^{k_0}\dd s'\\
        &=\frac{(-1)^{k_0}s^{k_0+1}}{(k_0+1)!}\mathcal{L}_{\xi}^{k_0}w-\frac{(-1)^{k_0}}{k_0!}\int_0^{s}\mathcal{L}_{\xi}^{k_0+1}(v_{s'}-\bar{v})(s-s')^{k_0}\dd s'\\
        &=\frac{(-1)^{k_0}s^{k_0+1}}{(k_0+1)!}\mathcal{L}_{\xi}^{k_0}w+\frac{(-1)^{k_0+1}}{k_0!}\int_0^{s}\left[\int_0^{s'} X_{l-s'}^*\mathcal{L}_{\xi}^{k_0+1}w\dd l\right](s-s')^{k_0}\dd s'\\
        &=\frac{(-1)^{k_0}s^{k_0+1}}{(k_0+1)!}\mathcal{L}_{\xi}^{k_0}w+\frac{(-1)^{k_0+1}}{(k_0+1)!}\int_0^{s}X_{-s'}^*\mathcal{L}_{\xi}^{k_0+1}w(s-s')^{k_0+1}\dd s'\\
        &=\frac{(-1)^{k_0}s^{k_0+1}}{(k_0+1)!}\mathcal{L}_{\xi}^{k_0}w+\frac{(-1)^{k_0+1}}{(k_0+1)!}\int_0^{s}(X_{s'})_*\mathcal{L}_{\xi}^{k_0+1}w(s-s')^{k_0+1}\dd s'\\
    \end{split}
\end{equation}
where we used Fubini's theorem to compute with $f(u)=X_{u}^*\mathcal{L}_{\xi}^{k_0+1}w$ the multiple integral
\begin{equation*}
    \begin{split}
        \int_0^{s}\left[\int_0^{s'}f(l-s')\dd l\right](s-s')^{k_0}\dd s'&=\int_0^{s}\left[\int_{-s'}^{0}f(u)\dd u\right](s-s')^{k_0}\dd s'\\
        &=\int_{-s}^{0}\left[\int_{-u}^{s}(s-s')^{k_0}\dd s'\right]f(u)\dd u\\
        &=-\frac{1}{k_0+1}\int_{-s}^0(s-s')^{k_0+1}|_{s'=-u}^{s'=s}f(l)\dd l\\
        &=\frac{1}{k_0+1}\int_{-s}^0(s+u)^{k_0+1}f(u)\dd u\\
        &=\frac{1}{k_0+1}\int_0^sf(-s')(s-s')^{k_0+1}\dd s'
    \end{split}
\end{equation*}
to switch the integral, one sees that the integration happens on $0\leq s'\leq s$ and $-s'\leq u\leq 0$, that is, on the triangle $0\leq -u\leq s'\leq s$, and thus for fixed $0\leq -u\leq s$ we have $s'$ in the range $-u\leq s'\leq s$.

\noindent\textit{Magnetic Field}: Computing as above with 
$$b=-\mathcal{L}_\xi\bar{B}=\mathcal{L}_{\bar{B}}\xi$$ replacing the role of $w$, one can show that:
\begin{equation}\label{LDF4}
    r^{k_0,s}_b=\frac{(-1)^{k_0}s^{k_0+1}}{(k_0+1)!}\mathcal{L}_{\xi}^{k_0}b+\frac{(-1)^{k_0+1}}{(k_0+1)!}\int_0^{s}(X_{s'})_*\mathcal{L}_{\xi}^{k_0+1}b(s-s')^{k_0+1}\dd s'.
\end{equation}

\subsection{Lagrangian Perturbation Lemma}
We have constructed a continuous family $(v_s, B_s)$ of solutions of \eqref{fo}  sharing the same initial condition and starting at $(\bar v,\bar B)$. It has been constructed by means of an auxiliary divergence-free vector field $\xi$. We now set $s=1$, write $(v_1, B_1)=(v,B),\ X_1=X$ and collect the work above in the following lemma. 
\begin{lemma}[Lagrangian Perturbation Lemma] \label{lpl} Let $\xi$ be a time-dependent divergence-free vector field on $\mathbb{T}^3\times \mathbb{R}$ such that $\xi|_{t=0}=0$, let $(\bar{v},\bar{B})$ be a smooth solution of \eqref{fo}. Let $X=\exp_{\IId}^\xi$ we call $(w,b)$ given by:
\begin{equation*}
    \begin{cases}
        w=\partial_tX\circ X^{-1}+(X_*-\IId_*)\bar v\\
        b=(X_*-\IId_*)\bar{B}
    \end{cases}
\end{equation*}
a perturbation of $(\bar{v},\bar{B})$ along the LDF $\xi$. We have that $(\bar{v}+w,\bar{B}+b)$ is still a solution of \eqref{fo} with initial condition $\bar{B}|_{t=0}$ and can be expanded in a truncated at level $k_0$ Lie-Taylor series of the following form:
\begin{equation*}
    \begin{split}
        w&=(\partial_t+\mathcal{L}_{\bar{v}})\xi+\sum_{k=1}^{k_0}\frac{(-1)^k}{(k+1)!}\mathcal{L}_{\xi}^k(\partial_t+\mathcal{L}_{\bar{v}})\xi+r^{k_0}_w,\\
        b&=\mathcal{L}_{\bar{B}}\xi+\sum_{k=1}^{k_0}\frac{(-1)^k}{(k+1)!}\mathcal{L}_{\xi}^k\mathcal{L}_{\bar B}\xi+r_b^{k_0},\\
        r^{k_0}_{w}&=\frac{(-1)^{k_0+1}}{(k_0+1)!}\int_0^{1}(X_s)_*\left[\mathcal{L}_{\xi}^{k_0+1}(\partial_t+\mathcal{L}_{\bar{v}})\xi\right](1-s)^{k_0+1}\dd s,\\
        r^{k_0}_{b}&=\frac{(-1)^{k_0+1}}{(k_0+1)!}\int_0^{1}(X_s)_*\left[\mathcal{L}_{\xi}^{k_0+1}\mathcal{L}_{\bar B}\xi\right](1-s)^{k_0+1}\dd s.
    \end{split}
\end{equation*}
Moreover, assuming $\xi=\curl \ \Theta$ this can be rewritten as:
\begin{equation*}
    \begin{split}
        w&=\curl\left[(\partial_t+\mathcal{L}_{\bar{v}})\Theta+\sum_{k=1}^{k_0}\frac{(-1)^k}{(k+1)!}\mathcal{L}_{\xi}^k(\partial_t+\mathcal{L}_{\bar{v}})\Theta+\theta^{k_0}_w\right],\\
        b&=\curl\left[\mathcal{L}_{\bar{B}}\Theta+\sum_{k=1}^{k_0}\frac{(-1)^k}{(k+1)!}\mathcal{L}_{\xi}^k\mathcal{L}_{\bar B}\Theta+\theta_b^{k_0}\right],\\
        \theta^{k_0}_{w}&=\frac{(-1)^{k_0+1}}{(k_0+1)!}\int_0^{1}(X_s)_*\left[\mathcal{L}_{\xi}^{k_0+1}(\partial_t+\mathcal{L}_{\bar{v}})\Theta\right](1-s)^{k_0+1}\dd s,\\
        \theta^{k_0}_{b}&=\frac{(-1)^{k_0+1}}{(k_0+1)!}\int_0^{1}(X_s)_*\left[\mathcal{L}_{\xi}^{k_0+1}\mathcal{L}_{\bar B}\Theta\right](1-s)^{k_0+1}\dd s
    \end{split}
\end{equation*}
where the Lie derivatives happen in the 1-form sense.
\end{lemma}

The first assertion follows immediately from \eqref{LDF1}, \eqref{LDF2}, \eqref{LDF3}, \eqref{LDF4}. The statement about the expansion in terms of the potential $\Theta$ follows from the properties of the pushforward and Lie derivative recalled in Appendix \ref{diffgeom}, namely \eqref{DG0} and \eqref{DG1}, given that $\bar v, \ \bar B, \ \xi$ are all divergence-free, with the Lie derivations occurring in the 1-form sense.


\section{Forced Galbrun's Equation} \label{gn}
This subsection is inspired by the work of Giri-Radu in \cite{GR}; the equation at hand, however, is of second order, and the 3D setting requires additional effort to write a formulation for the potential.
\subsection{Rewriting} \label{rew}
In the following, we consider the linearization of the MHD momentum equation set on $\mathbb{T}^3\times(t_0-\tau^c, t_0+\tau^c)$ at $(v,B)$ with forcing $f\mathbb{P}\ddiv F$. Namely, we look for time-dependent vector fields $(w,b)$ and pressure $\pi$ solving the problem:
\begin{equation}\label{linearized}
\begin{cases}
    \partial_t w+\left[v\cn w+w\cn v\right]-\left[ B\cn b +b\cn B\right]+\nabla \pi=f\mathbb{P}\ddiv F,\\
    \ddiv \ w = \ddiv \ b = 0.\\
\end{cases}
\end{equation}
Here $f=f(t/\tau^a)$ is a fixed one-periodic zero-average function rescaled using a parameter $\tau^a$ e.g. the ones constructed in Lemma \ref{timeprofiles}, in particular, the non-rescaled $f$ admits a bounded primitive $f^{[1]}$ defined explicitly as: 
$$f^{[1]}(t)=\int_0^t f(s)\dd s,$$
we think of them as functions $\mathbb{R}\to \mathbb{R}$ and always rescale them with $\tau^a$, moreover, $F$ is a given matrix field, e.g. $A_I$ in \eqref{slowcoeffgalb}. We also assume that on $\mathbb{T}^3\times(t_0-\tau^c,t_0+\tau^c)$ the quadruple $(v, \ B, \ p, \ R)$ solves: 
\begin{equation}\label{backgroundgalbrun}
    \begin{cases}
        \partial_t v+v\cn v+B\cn B+\nabla p=\ddiv \ R,\\
        \partial_t B+\curl [B\times v]=0,\\
        \ddiv \ v =\ddiv \ B =0.\\
    \end{cases}
\end{equation}
To simplify notation, we set:
\begin{equation*}
        D_t=\partial_t+v\cn,\ z^\pm=v\pm B, \ \mathcal{A}^\pm=\partial_t+z^\pm \cn.
\end{equation*}
We now gather the precise assumptions on the H{\"o}lder norms of $v, \ B, \ p, \ R, \ F$ and their transport properties.
\begin{assumptions}\label{standingass}
Fix $\underline{r}\geq 2$ an integer. We assume the following estimates on the vector fields:
\begin{subequations}
    \begin{align}
        &||v||_{r},||B||_{r}\lesssim \lambda_q^r\lambda_q^{[r-\underline{r}]^+(b-1)\gamma_\ell}\delta_q^{1/2}\ \text{ for } \ r\geq 1,\\
        &||\mathcal{A}^\pm v||_{r},||\mathcal{A}^\pm B||_{r}\lesssim \lambda_q^{r+1}\lambda_q^{[r-\underline{r}]^+(b-1)\gamma_\ell}\delta_q \ \text{ for }\ r\geq 0,\\
        &||p||_{r}\lesssim \lambda_q^{r}\lambda_q^{[r-\underline{r}]^+(b-1)\gamma_\ell}\delta_q \ \text{ for }\ r\geq 1,\\
        &||\mathcal{A}^\pm p||_{r}\lesssim \lambda_q^{r+1}\lambda_q^{[r-\underline{r}]^+(b-1)\gamma_\ell}\delta_q^{3/2} \ \text{ for } \ r \geq 1
    \end{align}
\end{subequations}
and the following bounds on the forcings:
\begin{subequations}
    \begin{align}
        &||R||_{r},||F||_{r}\lesssim \lambda_q^r\lambda_q^{[r-\underline{r}]^+(b-1)\gamma_\ell}\delta_{q+1} \ \text{ for }\ r\geq 0,\\
        &||\mathcal{A}^\pm R||_{r}\lesssim \lambda_q^{r+1}\lambda_q^{[r-\underline{r}]^+(b-1)\gamma_\ell}\delta_q^{1/2}\delta_{q+1} \ \text{ for }\ r\geq 0,\\
        &||\mathcal{A}^\pm F||_{r}\lesssim\lambda_q^{r}\lambda_q^{[r-\underline{r}]^+(b-1)\gamma_\ell}1/\tau^c\delta_{q+1} \ \text{ for }\ r\geq 0,\\
        &\sup_t|f|, \ \sup_t|f^{[1]}|\lesssim 1
    \end{align}
\end{subequations}
where all the implicit constants depend only on $r$ (in the applications, also on the parameters, but not on $a$). Note that the fractional H{\"o}lder norm bounds can be obtained from these by interpolation, and to simplify the notation, we set a common lower bound for the good derivative range $\underline{r}$, but it is possible to work in a more general setting, see also Remark \ref{lossapriori}.

\noindent We will restrict the time domain of the solution to $|t-t_0|< \tau^c$ that is $t\in (t_0-\tau^c,\ t_0+\tau^c)$ with $\tau^c$ given by \eqref{tauc}. For later applications, we assume that $F\equiv 0$ on: 
$$I^{cut}=(t_0-\tau^c, \ t_0-2/3\tau^c]\cup[t_0+2/3\tau^c, \ t_0+\tau^c).$$
\end{assumptions}
Going back to the linearization, we think of $(w,b)$ as the leading terms of a Lagrangian perturbation along some vector field $\xi$ with 1-form potential $\Theta$ as in Lemma \ref{lpl}, this gives: 
\begin{equation}\label{firstapprox}
        w=(\partial_t+\mathcal{L}_{v})\xi, \
        b=\mathcal{L}_{B}\xi, \
        \xi=\curl[\Theta].
\end{equation}
We are interested in rewriting \ref{linearized} above in terms of the 1-form potential $\Theta$, we first compute: 
\begin{equation}\label{rewriting1}
    \begin{split}
        f\mathbb{P}\ddiv F&=\partial_t w+\left[v\cn w+w\cn v\right]-\left[ B\cn b +b\cn B\right]+\nabla \pi\\
        &=\partial_t\left[(\partial_t+\mathcal{L}_{v})\xi\right]+v\cn\left[(\partial_t+\mathcal{L}_{v})\xi\right]+\left[(\partial_t+\mathcal{L}_{v})\xi\right]\cn v +\nabla \pi\\
        &-\left[ B\cn \mathcal{L}_{B}\xi  + \mathcal{L}_{B}\xi\cn B\right]\\
        &=(\partial_t+\mathcal{L}_{v})^2\xi-(\mathcal{L}_{B})^2\xi+2(\partial_t+\mathcal{L}_{v})\xi\cn v-2\mathcal{L}_{B}\xi\cn B +\nabla \pi\\
        &=\curl[\underbrace{(\partial_t+\mathcal{L}^1_{v})^2\Theta-(\mathcal{L}^1_{B})^2\Theta}_{T_1}]+\underbrace{2\ddiv\left[(\partial_t+\mathcal{L}^1_{v})\Theta\times \nabla v-\mathcal{L}^1_{B}\Theta\times \nabla B\right]^\top+\nabla \pi}_{T_2}\\
    \end{split}
\end{equation}
where we used the identities \eqref{DG1} and \eqref{DG2}. Taking the divergence of \eqref{rewriting1} we see that:
\begin{equation*}
        \Delta \pi = -2\ddiv\ddiv\left[(\partial_t+\mathcal{L}^1_{v})\Theta\times \nabla v-\mathcal{L}^1_{B}\Theta\times \nabla B\right]^\top
\end{equation*}
and writing the Leray projector as $\mathbb{P}=\IId-\nabla\Delta^{-1}\ddiv$ we end up with:
\begin{equation}\label{rewriting2}
    \begin{split}
        T_2&=2\mathbb{P}\ddiv\left[(\partial_t+\mathcal{L}^1_{v})\Theta\times \nabla v-\mathcal{L}^1_{B}\Theta\times \nabla B\right]\\
        &=\mathbb{P}\ddiv\left[\mathcal{A}^+\Theta \times \nabla z^-+\mathcal{A}^-\Theta \times \nabla z^++2\DD v^\top \Theta \times \nabla v-2\DD B^\top \Theta \times \nabla B\right],
    \end{split}
\end{equation}
we now rewrite $T_1$ in \eqref{rewriting1} as:
\begin{equation}\label{rewriting3}
    \begin{split}
        T_1&=(\partial_t+\mathcal{L}^1_{v})(\DD_t\Theta+\DD v^\top[\Theta])-\mathcal{L}^1_{B}(B\cn \Theta+\DD B^\top[\Theta])\\
        &=\DD_t^2\Theta+\DD_t(\DD v^\top[\Theta])+\DD v^\top[\DD_t\Theta+\DD v^\top[\Theta]]\\
        &-(B\cn)^2\Theta -B\cn(\DD B^\top[\Theta])-\DD B^\top[B\cn\Theta+\DD B^\top[\Theta]]\\
        &=\DD_t^2\Theta+\DD(\DD_tv)^\top[\Theta]-(\DD v\DD v)^\top[\Theta]+\DD v^\top\DD_t\Theta+\DD v^\top[\DD_t\Theta+\DD v^\top[\Theta]]\\
        &-(B\cn)^2\Theta -\DD(B\cn B)^\top[\Theta]+(\DD B \DD B)^\top[\Theta]-\DD B^\top[B\cn\Theta]-\DD B^\top[B\cn\Theta+\DD B^\top[\Theta]]\\
        &=\DD_t^2\Theta-(B\cn)^2\Theta+\DD(\DD_tv-B\cn B)^\top[\Theta]+2\DD v^\top[\DD_t\Theta]-2\DD B^\top[B\cn\Theta]\\
        &=\mathcal{A}^+\mathcal{A}^-\Theta+(\DD[\DD_tv-B\cn B])^\top[\Theta]+(\DD z^-)^\top[\mathcal{A}^+\Theta]+(\DD z^+)^\top[\mathcal{A}^-\Theta]
    \end{split}
\end{equation}
where we used that:
$$\DD_t B\cn=(\partial_tB+[v,B])\cn+B\cn D_t $$ 
and since by assumption $\partial_tB+[v,B]=0$, see \eqref{backgroundgalbrun}, we can write: 
\begin{equation}\label{pmcommutation}
    \begin{split}
        \mathcal{A}^-\mathcal{A}^+&=\left(D_t-B\cn\right)(D_t+B\cn)\\
        &=D_t^2-(B\cn)^2+\DD_t B\cn-B\cn D_t\\
        &=D_t^2-(B\cn)^2 +(\partial_tB+[v,B])\cn\\
        &=D_t^2-(B\cn)^2.
    \end{split}
\end{equation}
Note, in passing, that this also shows the Alfv\'en transport operators commutation: 
$$\mathcal{A}^+\mathcal{A}^-=\mathcal{A}^-\mathcal{A}^+.$$
Substituting \eqref{rewriting2} in \eqref{rewriting1} and applying $\curl^{-1}=-\Delta^{-1}\curl$, we finally reach 
\begin{equation*}
    \begin{split}
        &\curl^{-1}\curl[T_1]+\curl^{-1}\mathbb{P}\ddiv\left[\mathcal{A}^+\Theta \times \nabla z^-+\mathcal{A}^-\Theta \times \nabla z^++2\DD v^\top \Theta \times \nabla v+2\DD B^\top \Theta \times \nabla B\right]^\top\\
        &=f\curl^{-1}\mathbb{P}\ddiv F,
    \end{split}
\end{equation*}
where we fix the gauge of $\curl^{-1}$ to be Coulomb, i.e. requiring the output to be divergence-free. We now set $\mathbb{T}$ to be the CZ-type operator appearing in the rewriting, namely
$$\mathbb{T}=\curl^{-1}\mathbb{P}\ddiv$$
and substitution of the expression of $T_1$ found in \eqref{rewriting3} then leads to: 
\begin{equation*}
    \begin{split}
        &\mathcal{A}^+\mathcal{A}^-\Theta+(\DD[\DD_tv-B\cn B])^\top[\Theta]+(\DD z^-)^\top[\mathcal{A}^+\Theta]+(\DD z^+)^\top[\mathcal{A}^-\Theta]+\\
        &+\mathbb{T}\left[\mathcal{A}^+\Theta \times \nabla z^-+\mathcal{A}^-\Theta \times \nabla z^++2\DD v^\top \Theta \times \nabla v-2\DD B^\top \Theta \times \nabla B\right]^\top+\nabla \tilde \pi\\
        &=f\mathbb{T}\ddiv F,
    \end{split}
\end{equation*}
for some $\tilde \pi$ coming from the Gauge condition, namely
$$\curl^{-1}\curl[T_1]=T_1+\nabla \tilde \pi.$$
We forget about this gradient, as we are not interested in keeping any divergence-free constraint on $\Theta$ and solve instead: 
\begin{equation}\label{rew1}
    \begin{split}
        &\mathcal{A}^+\mathcal{A}^-\Theta+(\DD[\DD_tv-B\cn B])^\top[\Theta]+(\DD z^-)^\top[\mathcal{A}^+\Theta]+(\DD z^+)^\top[\mathcal{A}^-\Theta]\\
        &+\mathbb{T}\left[\mathcal{A}^+\Theta \times \nabla z^-+\mathcal{A}^-\Theta \times \nabla z^++2\DD v^\top \Theta \times \nabla v-2\DD B^\top \Theta \times \nabla B\right]^\top\\
        &=f\mathbb{T}\ddiv F,
    \end{split}
\end{equation}
the two possibly different solutions for the potential lead to $(w,b)$ solving the same equation, which is what really matters to us. Note in addition that one of the terms involving $\Theta$ sees the fact that $(v, \ B, \ p, \ R)$ solve \eqref{backgroundgalbrun}, indeed
$$\DD[\DD_tv-B\cn B])^\top[\Theta]=\DD(\ddiv \ R- \nabla p)^\top[\Theta]$$
this will be important later when doing estimates. We now set for any input vectors $a, b, c\in \mathbb{R}^3$:
\begin{equation}\label{Hdefinition}
    \begin{split}
        &H_1(a)=(\DD[\DD_tv-B\cn B])^\top[a], \ \tilde{H}_1(a)=2\left[\DD v^\top a \times \nabla v+\DD B^\top a \times \nabla B\right]^\top,\\
        &H_2(b,c)=(\DD z^-)^\top[b]+(\DD z^+)^\top[c], \ \tilde{H}_2(b,c)=\left[b \times \nabla z^-+c \times \nabla z^+\right]^\top.\\
    \end{split}
\end{equation}
note that apart from the parameters $a,b,c$, all the objects depend also on $(x,t)$ and we think of them as matrices acting on vectors. To simplify the notation further, we group the local and non-local parts into two single operators, one of zeroth order with respect to $\Theta$ and the other at first order in the Alfv\'en transport of $\Theta$. Namely, we set:
\begin{equation*}
        \mathbb{H}_1=H_1+\mathbb{T}\tilde{H}_1, \ 
        \mathbb{H}_2=H_2+\mathbb{T}\tilde{H}_2.
\end{equation*}
For later use, we remark that solving \eqref{rew1} is the same as solving:
\begin{equation}\label{secondlietransportgalbrun}
    \begin{split}
        (\partial_t+\mathcal{L}^1_{v})^2\Theta-(\mathcal{L}^1_{B})^2\Theta+\mathbb{T}\tilde H_1(\Theta)+\mathbb{T}\tilde H_2(\mathcal{A}^+\Theta,\mathcal{A}^-\Theta)=f\mathbb{T}F.
    \end{split}
\end{equation}
We formalise the above manipulations of \eqref{rew1} in the following Lemma.
\begin{lemma}[Galbrun Equation] The forced linearization \eqref{linearized} of the MHD momentum equation solved on $\mathbb{T}^3\times(t_0-\tau^c,t_0+\tau^c)$ under the Standing Assumptions \ref{standingass} and \eqref{backgroundgalbrun}, using the first-order approximation \eqref{firstapprox} of a perturbation along a Lagrangian displacement field $\xi$ with potential $\Theta$, see Lemma \ref{lpl}, corresponds to the second-order initial value problem:
    \begin{equation}\label{rewritinggr}
    \begin{cases}
        \mathcal{A}^-\mathcal{A}^+\Theta+\mathbb{H}_1[\Theta]+\mathbb{H}_2[\mathcal{A}^+\Theta,\mathcal{A}^-\Theta]=f\mathbb{T} F,\\
        \Theta|_{t=t_0}=G,\\
        \mathcal{A}^+\Theta|_{t=t_0}=H.\\
    \end{cases}
\end{equation}
\end{lemma}
We now prove some bounds on the CZ operators $\mathbb{H}_1$ and $\mathbb{H}_2$, which will be needed later to prove a priori estimates on the solution. To do this, it is useful to define: 
\begin{equation*}
    \begin{split}
    h_1^r&=\lambda_q^r\lambda_q^{[r-(\underline{r}-2)]^+(b-1)\gamma_\ell}(\lambda_q^{2}\ell^{-\alpha}\delta_q)\\
    h_{1,t}^r&=\lambda_q^r\lambda_q^{[r-(\underline{r}-3)]^+(b-1)\gamma_\ell}(\lambda_q^{2}\ell^{-\alpha}\delta_q)\\
    \end{split} \ \text{ and } \
    \begin{split}
        h_2^r&=\lambda_q^r\lambda_q^{[r-(\underline{r}-1)]^+(b-1)\gamma_\ell}(\lambda_q\ell^{-\alpha}\delta_q^{1/2})\\
    h_{2,t}^r&=\lambda_q^r\lambda_q^{[r-(\underline{r}-2)]^+(b-1)\gamma_\ell}(\lambda_q\ell^{-\alpha}\delta_q^{1/2})
    \end{split}
\end{equation*}
these correspond to the bounds for the coefficients in the operators $\mathbb{H}_1, \ \mathbb{H}_2$, and we remark the loss of two good derivatives in $h_1$.
\begin{lemma}[Operator Norms]\label{operatornorms} 
Under the Standing Assumptions \ref{standingass}, we have the following estimates:
\begin{subequations}
    \begin{align}
        ||\mathbb{H}_1\Theta||_{r+\alpha}&\lesssim h_1^{r}||\Theta||_0+h_1||\Theta||_{r+\alpha},\\
        ||\mathbb{H}_2(\Theta_1,\Theta_2)||_{r+\alpha}&\lesssim h_2^{r}(||\Theta_1||_0+||\Theta_2||_0)+h_2^0(||\Theta_1||_{r+\alpha}+||\Theta_2||_{r+\alpha}),\\
        ||[\mathcal{A}^\pm,\mathbb{H}_1]\Theta||_{r+\alpha}&\lesssim \lambda_q\ell^{-\alpha}\delta_q^{1/2}\left[h_{1}^{r}||\Theta||_\alpha+h_{1}^0||\Theta||_{r+\alpha}\right],\\
        ||[\mathcal{A}^\pm,\mathbb{H}_2](\Theta_1,\Theta_2)||_{r+\alpha}&\lesssim \lambda_q\ell^{-\alpha}\delta_q^{1/2}\left[h_{2}^{r}(||\Theta_1||_\alpha+||\Theta_2||_\alpha)+h_{2}^0(||\Theta_1||_{r+\alpha}+||\Theta_2||_{r+\alpha})\right],\\
        ||[\partial_t,\mathbb{H}_1]\Theta||_{r+\alpha}&\lesssim \lambda_q\left[h_{1,t}^{r}||\Theta||_0+h_{1,t}^0||\Theta||_{r+\alpha}\right],\\
        ||[\partial_t,\mathbb{H}_2](\Theta_1,\Theta_2)||_{r+\alpha}&\lesssim \lambda_q\left[h_{2,t}^{r}(||\Theta_1||_0+||\Theta_2||_0)+h_{2,t}^0(||\Theta_1||_{r+\alpha}+||\Theta_2||_{r+\alpha})\right].
    \end{align}
\end{subequations}
    The implicit constants depend on $r$ and $\alpha$.
\end{lemma}
\begin{remark}
    The commutator with the Alfv\'en transport has a slightly different bound because we use the estimate from Proposition \ref{czstuff}, which loses an additional $\alpha$.
\end{remark}
\begin{proof}[Proof of Lemma \ref{operatornorms}] The proof is just a lengthy sequence of computations. We will use Proposition \ref{czstuff} several times to deal with the operator $\mathbb{T}$ and commutators. We show explicit bounds only for $\mathbb{H}_1$; the ones for $\mathbb{H}_2$ follow similarly. We have:
\begin{equation*}
    \begin{split}
        ||\mathbb{H}_1\Theta||_{r+\alpha} &\lesssim||\mathbb{T}\left[\DD v^\top \Theta \times \nabla v+\DD B^\top \Theta \times \nabla B\right]^\top||_{r+\alpha}+||\DD[\DD_tv-B\cn B][\Theta]||_{r+\alpha}\\
        &\lesssim||\left[\DD v^\top \Theta \times \nabla v+\DD B^\top \Theta \times \nabla B\right]^\top||_{r+\alpha}+||\DD[\ddiv R-\nabla p][\Theta]||_{r+\alpha}\\
        &\lesssim ||\Theta||_{0}(||v||_{r+1+\alpha}||v||_1+||B||_{r+1+\alpha}||B||_1)+||\Theta||_{r+\alpha}(||v||_1^2+||B||_1^2)\\
        &+||\Theta||_{r+\alpha}(||R||_{2}+||p||_2)+||\Theta||_{0}(||R||_{r+2+\alpha}+||p||_{r+2+\alpha})\\
        &\lesssim h_1^0 ||\Theta||_{r+\alpha} +h_1^{r}||\Theta||_0.\\
    \end{split}
\end{equation*}
The commutator estimates are done similarly, using the assumptions on the transport derivatives:
\begin{equation*}
    \begin{split}
        &||[\mathcal{A}^\pm,\mathbb{H}_1]\Theta||_{r+\alpha} \\
        &\lesssim||\mathcal{A}^\pm\mathbb{T}\left[\DD v^\top \Theta \times \nabla v+\DD B^\top \Theta \times \nabla B\right]^\top-\mathbb{T}\left[\DD v^\top \mathcal{A}^\pm\Theta \times \nabla v+\DD B^\top \mathcal{A}^\pm\Theta \times \nabla B\right]^\top||_{r+\alpha}\\
        &+||\mathcal{A}^\pm\DD[\DD_tv-B\cn B][\Theta]-\DD[\DD_tv-B\cn B][\mathcal{A}^\pm\Theta]||_{r+\alpha}\\
        &\lesssim ||\mathbb{T}\DD (\mathcal{A}^\pm v)^\top \Theta \times \nabla v||_{r+\alpha}+ ||\mathbb{T}\DD v^\top \Theta \times \nabla (\mathcal{A}^\pm v)||_{r+\alpha}\\
        &+||\mathbb{T}\DD(z^\pm)^\top\DD v^\top \Theta \times \nabla v||_{r+\alpha}+ ||\mathbb{T}\DD v^\top \Theta \times\DD (z^\pm)^\top \DD v^\top||_{r+\alpha}\\
        &+||\mathbb{T}\DD (\mathcal{A}^\pm B)^\top \Theta \times \nabla B||_{r+\alpha}+ ||\mathbb{T}\DD B^\top \Theta \times \nabla (\mathcal{A}^\pm B)||_{r+\alpha}\\
        &+||\mathbb{T}\DD(z^\pm)^\top\DD B^\top \Theta \times \nabla B||_{r+\alpha}+ ||\mathbb{T}\DD B^\top \Theta \times\DD (z^\pm)^\top \DD B^\top||_{r+\alpha}\\
        &+||[\mathcal{A}^\pm,\mathbb{T}]\left[\DD v^\top \Theta \times \nabla v+\DD B^\top \Theta \times \nabla B\right]^\top||_{r+\alpha}\\
        &+||\DD[\ddiv [\mathcal{A}^\pm R][\Theta]||_{r+\alpha}+||[\mathcal{A}^\pm,\DD\ddiv] (R)[\Theta]||_{r+\alpha}+||\DD[\nabla\mathcal{A}^\pm p][\Theta]||_{r+\alpha}+||[\mathcal{A}^\pm,\DD \nabla](p)[\Theta]||_{r+\alpha}\\
        &\lesssim \lambda_q\ell^{-\alpha}\delta_q^{1/2}\left[h_{1}^{r}||\Theta||_\alpha+h_{1}^0||\Theta||_{r+\alpha}\right]
    \end{split}
\end{equation*}
where we used that according to \eqref{backgroundgalbrun}, we have $\DD_tv-B\cn B=\ddiv \ R - \nabla p$, together with the following computations:
\begin{equation*}
    \begin{split}
        &||[\mathcal{A}^\pm,\mathbb{T}]\left[\DD v^\top \Theta \times \nabla v+\DD B^\top \Theta \times \nabla B\right]^\top||_{r+\alpha}\\
        &\lesssim ||z^\pm||_{1+\alpha} ||\DD v^\top \Theta \times \nabla v+\DD B^\top \Theta \times \nabla B||_{r+\alpha}+||z^\pm||_{r+1+\alpha} ||\DD v^\top \Theta \times \nabla v+\DD B^\top \Theta \times \nabla B||_{\alpha}\\
        &\lesssim\lambda_q^{r+3}\lambda_q^{[r-(\underline{r}-1)]^+(b-1)\gamma_\ell}\ell^{-2\alpha}\delta_q^{3/2}||\Theta||_\alpha+\lambda_q^{3}\ell^{-2\alpha}\delta_q^{3/2}||\Theta||_{r+\alpha}\\
        &=\lambda_q\ell^{-\alpha}\delta_q^{1/2}\left[h_1^{r}||\Theta||_\alpha+h_1^0||\Theta||_{r+\alpha}\right].
    \end{split}
\end{equation*}
and
\begin{equation*}
    \begin{split}
        \mathcal{A}^\pm\DD \ddiv R&=\DD \mathcal{A}^\pm\ddiv R-(\DD \ddiv R)(\DD z^\pm)=\DD \left[\ddiv \mathcal{A}^\pm R-\nabla v:\nabla R\right]-(\DD \ddiv R)(\DD z^\pm)\\
    \end{split}
\end{equation*}
where $(\nabla v:\nabla R)^i=\partial_j v^k\partial_kR^{ij}$, from which we deduce:
$$||[\mathcal{A}^\pm,\DD\ddiv] (R)[\Theta]||_{r+\alpha}\lesssim \lambda_q\ell^{-\alpha}\delta_q^{1/2}\left[h_{1}^{r}||\Theta||_\alpha+h_{1}^0||\Theta||_{r+\alpha}\right].$$
Similar bounds hold for the terms involving $p$.

\noindent We are left with the commutator with the pure time-derivative, which is non-zero as the vector fields are time-dependent. We first compute:
\begin{equation*}
    \begin{split}
        ||\partial_t[\DD_tv-B\cn B]||_{r+\alpha}&=||\partial_t[\ddiv R-\nabla p]||_{r+\alpha}\\
        &=||\mathbb{P} \ddiv \partial_t R + \nabla \Delta^{-1}\ddiv[\partial_t[v\cn v-B\cn B]]||_{r+\alpha}\\
        &=||\mathbb{P} \ddiv \partial_t R + 2\nabla \Delta^{-1}\ddiv[[\partial_tv\cn v-\partial_tB\cn B]]||_{r+\alpha}\\
        &\lesssim ||\partial_t R||_{r+1+\alpha}+||[\partial_tv\cn v-\partial_tB\cn B]||_{r+\alpha}\\
        &\lesssim \lambda_q^{r+2}\lambda_q^{[r-(\underline{r}-2)]^+(b-1)\gamma_\ell}\ell^{-\alpha} \delta_{q+1}+\lambda_q^{r+2}\lambda_q^{[r-(\underline{r}-2)]^+(b-1)\gamma_\ell}\ell^{-\alpha}\delta_{q}\\
        &\lesssim \lambda_q^{r+2}\lambda_q^{[r-(\underline{r}-2)]^+(b-1)\gamma_\ell}\ell^{-\alpha}\delta_{q}\\
    \end{split}
\end{equation*}
where we used that: 
$$\ddiv \ a= \ddiv \ b =0 \Longrightarrow\ddiv [a\cn b]=\ddiv [b\cn a]$$
and deduced the pure time derivatives of $v, \ R$ by writing:
$$\partial_t v=\mathcal{A}^\pm v-z^\pm\cn v,$$
similarly, for $B,\ R$. We use this bound in the following estimate:
\begin{equation*}
    \begin{split}
        &||[\partial_t,\mathbb{H}_1]\Theta||_{r+\alpha} \\
        &\lesssim||\partial_t\mathbb{T}\left[\DD v^\top \Theta \times \nabla v+\DD B^\top \Theta \times \nabla B\right]^\top-\mathbb{T}\left[\DD v^\top \partial_t\Theta \times \nabla v+\DD B^\top \partial_t\Theta \times \nabla B\right]^\top||_{r+\alpha}\\
        &+||\partial_t\DD[\DD_tv-B\cn B][\Theta]-\DD[\DD_tv-B\cn B][\partial_t\Theta]||_{r+\alpha}\\
        &\lesssim ||\partial_tv||_{r+1+\alpha}||\Theta||_0||v||_1+||\partial_tv||_{1}||\Theta||_{r+\alpha}||v||_1 +||\partial_tv||_{1}||\Theta||_0||v||_{r+1+\alpha}\\
        &+ ||\partial_tB||_{r+1+\alpha}||\Theta||_0||B||_1+||\partial_tB||_{1}||\Theta||_{r+\alpha}||B||_1 +||\partial_tB||_{1}||\Theta||_0||B||_{r+1+\alpha}\\
        &+||\DD \partial_t[\DD_tv-B\cn B]||_{r+\alpha}||\Theta||_0+||\DD \partial_t[\DD_tv-B\cn B]||_{0}||\Theta||_{r+\alpha}\\
        &\lesssim \lambda_q\left[h_{1,t}^{r}||\Theta||_0+h_{1,t}^0||\Theta||_{r+\alpha}\right].\\
    \end{split}
\end{equation*}

\noindent The implicit constants in the above bounds depend on $r$ and $\alpha$.
\end{proof}

\subsection{Existence and Uniqueness} 
\begin{lemma}[Existence and Uniqueness]\label{EU} Let $r\geq 1$ integer, $0<\alpha<1$ and $0<T<\infty$. Let $z^\pm\in C([-T,T],C^{r+1+\alpha}(\mathbb{T}^3))$ vector fields and $F\in C([-T,T],C^{r+\alpha}(\mathbb{T}^3))$ a 1-form, define $\mathcal{A}^\pm=\partial_t+z^\pm\cn$ and assume
$$\mathcal{A}^+\mathcal{A}^-=\mathcal{A}^-\mathcal{A}^+,$$
moreover, let $$\mathbb{H}_1:C^{r'+\alpha}(\mathbb{T}^3,\mathbb{R}^3)\to C^{r'+\alpha}(\mathbb{T}^3,\mathbb{R}^3), \ \ \ \ \mathbb{H}_2:C^{r'+\alpha}(\mathbb{T}^3,\mathbb{R}^3)\times C^{r'+\alpha}(\mathbb{T}^3,\mathbb{R}^3)\to C^{r'+\alpha}(\mathbb{T}^3,\mathbb{R}^3)$$
be bounded linear operators for $0\leq r'\leq r$.  Consider the initial value problem
\begin{equation}\label{rewritinggr1}
    \begin{cases}
        \mathcal{A}^-\mathcal{A}^+\Theta+\mathbb{H}_1[\Theta]+\mathbb{H}_2[\mathcal{A}^+\Theta,\mathcal{A}^-\Theta]=F,\\
        \Theta|_{t=0}=G,\\
        \mathcal{A}^+\Theta|_{t=0}=H\\
    \end{cases}
\end{equation}
for a 1-form $\Theta$. Then for any pair $G,H\in C([-T,T],C^{r+\alpha}(\mathbb{T}^3)$ the problem \eqref{rewritinggr1} admits a unique solution such that:
    $$\Theta \in C([-T,T],C^{r+\alpha}(\mathbb{T}^3)).$$
\end{lemma}
\begin{proof}[Proof of Lemma \ref{EU}] We adopt the notation and definitions in Subsection \ref{apriorie}. We begin with the following claim:
\begin{equation}\label{timecommute}
    X_{s'}^+\circ\Phi^+_s\circ X_s^-\circ \Phi_t^-=X_{s'}^-\circ\Phi^-_{s'-s+t}\circ X_{s'-s+t}^+\circ \Phi_t^+
\end{equation}
whenever all the objects are well defined, equivalently
$$X_{s',s}^+\circ X_{s,t}^-=X_{s',s'-s+t}^-\circ X^+_{s'-s+t,t}.$$
Recall from Subsection \ref{timemollification} that
$$\mathbb{X}^\pm_s(x,t)=(X^\pm_{t+s,t}(x),t+s)$$
and since $\mathcal{A}^\pm$ commute by assumption, we have: 
$$\mathcal{A}^+\mathcal{A}^-=\mathcal{A}^-\mathcal{A}^+ \Longrightarrow\mathbb{X}^+_a\circ \mathbb{X}_{b}^-=\mathbb{X}^-_b\circ \mathbb{X}_{a}^+$$
for any $a,b$ such that the composition is well defined, we now apply it with $a=s'-s, \ b=s-t$ and obtain:
\begin{equation*}
    \begin{split}
        (X_{s',s}^+\circ X_{s,t}^-(x),s')&=(X_{s',\cdot}^+\circ \mathbb{X}_{s-t}^-(x,t),s')\\
        &=\mathbb{X}^+_{s'-s}\circ \mathbb{X}_{s-t}^-(x,t)\\
        &=\mathbb{X}_{s-t}^-\circ \mathbb{X}^+_{s'-s}(x,t)\\
        &=\mathbb{X}_{s-t}^-(X^+_{s'-s+t,t}(x),s'-s+t)\\
        &=(X_{s',s'-s+t}^-\circ X^+_{s'-s+t,t}(x),s')
    \end{split}
\end{equation*}
where $\circ$ can be either a composition in space or space-time. This is the claim in \eqref{timecommute}. It follows that:
\begin{equation}\label{timecommuteint}
    \begin{split}
        \int_{t_0}^t\left[\int_{t_0}^s(F)(X_{s'}^-\circ\Phi^-_s\circ X_s^+\circ \Phi_t^+,s')\dd s'\right]\dd s&=\int_{t_0}^t\left[\int_{t_0}^s(F)(X_{s'}^-\circ\Phi^-_{s'-s+t}\circ X_{s'-s+t}^+\circ \Phi_t^+,s')\dd s'\right]\dd s\\
        &=\int_{t_0}^t\left[\int_{s'}^t(F)(X_{s'}^-\circ\Phi^-_{s'-s+t}\circ X_{s'-s+t}^+\circ \Phi_t^+,s')\dd s\right]\dd s'\\
        &=\int_{t_0}^t\left[\int_{s'}^t(F)(X_{s'}^-\circ\Phi^-_{u}\circ X_{u}^+\circ \Phi_t^+,s')\dd u\right]\dd s'\\
        &=\int_{t_0}^t\left[\int_{t_0}^u(F)(X_{s'}^+\circ\Phi^+_u\circ X_u^-\circ \Phi_t^-,s')\dd s'\right]\dd u\\
        &=\int_{t_0}^t\left[\int_{t_0}^s(F)(X_{s'}^+\circ\Phi^+_s\circ X_s^-\circ \Phi_t^-,s')\dd s'\right]\dd s\\
    \end{split}
\end{equation}
where we used the fact that the integration is over the triangle $t_0\leq s'\leq s\leq t$ to switch the order of integration using Fubini's theorem.

\noindent \textbf{Local Existence and Uniqueness.} As we will prove detailed a priori estimates later, we will only be brief here. We want to find a solution of \eqref{rewritinggr1} on $\mathbb{T}^3$ for some short time interval $[t_0, t_0+\tau]$ (or $[t_0-\tau, t_0]$) first. We will do so by a Banach Fixed Point argument. A solution of \eqref{rewritinggr} by the Duhamel formula, see \eqref{duhamel} for the scalar case applied twice, corresponds to a fixed point of the map:
\begin{equation}\label{fixedpointmap}
    \begin{split}
        \mathbb{F}[\Theta]&=\underbrace{G(\Phi_t^+)+\int_{t_0}^tH(\Phi^-_s\circ X_s^+\circ \Phi_t^+)\dd s+\int_{t_0}^t\left[\int_{t_0}^s(F)(X_{s'}^-\circ\Phi^-_s\circ X_s^+\circ \Phi_t^+,s')\dd s'\right]\dd s}_{I}\\
        &+\underbrace{\int_{t_0}^t\left[\int_{t_0}^s(-\mathbb{H}_1[\Theta]-\mathbb{H}_2[\mathcal{A}^+\Theta,\mathcal{A}^-\Theta])(X_{s'}^-\circ\Phi^-_s\circ X_s^+\circ \Phi_t^+,s')\dd s'\right]\dd s}_{II}.\\
    \end{split}
\end{equation}
We now choose a space to set the argument, namely, the closed ball:
$$\mathbb{B}_{R,\tau,t_0}=\left\{\Theta \in C([t_0,t_0+\tau];C^{r+\alpha}(\mathbb{T}^3):\mathcal{A}^\pm\Theta \in C([t_0,t_0+\tau];C^{r+\alpha}(\mathbb{T}^3), \ ||\Theta||_{\mathbb{B}} \leq R\right\}$$
where the norm is defined as:
$$||\Theta||_{\mathbb{B}}:=\sup_{t\in[t_0,t_0+\tau]}\left[||\Theta(\cdot,t)||_{r+\alpha}+||\mathcal{A}^+\Theta(\cdot,t)||_{r+\alpha}+||\mathcal{A}^-\Theta(\cdot,t)||_{r+\alpha}\right],$$
in particular,  $\mathbb{B}_{R,\tau,t_0}$ is a complete metric space. Let $R\geq 1$ be so large that for any $0< \tau< 1$ for $I$ as in \eqref{fixedpointmap} we have: 
$$||I||_{\mathbb{B}}\leq R/2$$
this depends on the $C^r$ norm of the initial data and the $C^{r+1}$ norm of $z^\pm$ but not on $\tau$.

\noindent Now from \eqref{timecommuteint} we deduce:
$$\mathcal{A}^{\pm}II=\int_{t_0}^t(\mathbb{H}_1[\Theta]+\mathbb{H}_2[\mathcal{A}^+\Theta,\mathcal{A}^-\Theta])(X^\mp_s\circ\Phi_t^\mp,s)\dd s$$
where we used that, eventually commuting the flows using \eqref{timecommuteint}, when $\mathcal{A}^\pm$ hit the integrand, we get zero, leaving us with the time derivatives of the integral itself. 

\noindent From this and \eqref{fixedpointmap} we deduce: 
$$||\mathbb{F}\Theta||_{\mathbb{B}}\leq R/2+\tau C||\mathbb{H}_1||_{op}||\Theta||_{r+\alpha}+\tau C'||\mathbb{H}_2||_{op}(||\mathcal{A}^+\Theta||_{r+\alpha}+||\mathcal{A}^-\Theta||_{r+\alpha})$$
for some $C,C'$ depending on the $C^{r+1}$ norm of $z^\pm$ and we can choose $\tau$ small enough to ensure
$$\tau [C||H_1||_{op}+C'||H_2||_{op}]\leq 1/2$$ 
and we conclude that for $\Theta \in \mathbb{B}_{R,\tau}$ we have:
\begin{equation}\label{mapping}
    ||\mathbb{F}\Theta||_{\mathbb{B}}\leq R \Longrightarrow \mathbb{F}(\mathbb{B}_{R,\tau,t_0})\subset \mathbb{B}_{R,\tau,t_0}.
\end{equation}

\noindent We move to the contraction property. We set $\Delta=\Theta_2-\Theta_1$ and compute:
\begin{equation}\label{contraction}
    \begin{split}
        \mathbb{F}[\Theta_1]-\mathbb{F}[\Theta_2]&=-\int_{t_0}^t\left[\int_{t_0}^s(\mathbb{H}_1[\Theta_1]+\mathbb{H}_2[\mathcal{A}^+\Theta_1,\mathcal{A}^-\Theta_1])(\dots,s')\dd s'\right]\dd s\\
        &+\int_{t_0}^t\left[\int_{t_0}^s(\mathbb{H}_1[\Theta_2]+\mathbb{H}_2[\mathcal{A}^+\Theta_2,\mathcal{A}^-\Theta_2])(\dots,s')\dd s'\right]\dd s\\
        &=\int_{t_0}^t\left[\int_{t_0}^s(\mathbb{H}_1[\Delta]+\mathbb{H}_2[\mathcal{A}^+\Delta,\mathcal{A}^-\Delta])(\dots,s')\dd s'\right]\dd s\\
    \end{split}
\end{equation}
from which we deduce:
\begin{equation*}
    \begin{split}
        \mathcal{A}^{\pm}[\mathbb{F}[\Theta_1]-\mathbb{F}[\Theta_2]]&=\int_{t_0}^t(\mathbb{H}_1[\Delta]+\mathbb{H}_2[\mathcal{A}^+\Delta,\mathcal{A}^-\Delta])(X^\mp_s\circ\Phi_t^\mp,s)\dd s\\
    \end{split}
\end{equation*}
and we conclude that: 
\begin{equation*}
        ||\mathbb{F}[\Theta_1]-\mathbb{F}[\Theta_2]||_{\mathbb{B}}\leq C \tau ||\Theta_1-\Theta_2||_{\mathbb{B}}
\end{equation*}
where $C$ is some constant that depends on the $C^{r+1}$ norm of $z^\pm$ and the operator norms of $\mathbb{H}_1, \ \mathbb{H}_2$ but not on the initial data, eventually making $\tau$ smaller we can ensure:
\begin{equation}\label{contraction1}
    ||\mathbb{F}[\Theta_1]-\mathbb{F}[\Theta_2]||_{\mathbb{B}}\leq \frac{1}{2} ||\Theta_1-\Theta_2||_{\mathbb{B}}
\end{equation}
and the desired contraction property follows for $\mathbb{F}$ in the space $\mathbb{B}_{R,\tau,t_0}$. Note that $\tau$ does not depend on $t_0, \ R$.

\noindent Existence then follows from \eqref{mapping}, \eqref{contraction1} and an application of Proposition \ref{CMT}. 

\noindent \textbf{Uniqueness.} Given the regularity of the solution, uniqueness is a consequence of \eqref{contraction} and Gr{\"o}nwall Lemma. 

\noindent\textbf{Long Time Existence:} We can now cover $[-T,T]$ with finitely many intervals of length $\tau$ of the form $[t_j,t_j+\tau]$. We then construct an initial solution given the data $(G, H)$ on $[0, \tau]$ and extend it forward and backwards finitely many times, using the value of the solution at $t_{j-1}+\tau$ as initial condition. This requires increasing the value of $R=R_j$ in the fixed point argument above finitely many times to ensure the iteratively constructed initial data lies in the space $\mathbb{B}_{\tau, R_j, t_j}$, uniqueness guarantees the existence of a global solution given any finite $T$.
\end{proof}


\subsection{A Priori Estimates}\label{apriorie}
With the existence and uniqueness of the problem \eqref{rewritinggr} at hand, we now want to prove a priori estimates. The key observations are the following:
\begin{itemize}
    \item The commutation \eqref{pmcommutation} of the Alfv\'en transport operators $\mathcal{A}^\pm$ allows for the propagation of $C^{r+\alpha}$ estimates on the solution. Despite its geometric nature, \eqref{rewritinggr} effectively behaves as a system of first-order PDEs.
    \item The particular choice of forcing allows, basically by integration by parts (IBP) in time, to trade a $\tau^a$ gain of smallness on the solution for a transport derivative of the forcing $\mathcal{A}^\pm \ddiv F$. Given the second-order nature of the equation, it is possible to gain a second $\tau^a$; one then needs to ensure that good second-order transport properties of the forcing are available. We do not pursue this here and gain only $\tau^c$ coming from the restriction of the time-domain of the solution $|t-t_0|< \tau^c$. This gives a total gain of $\tau^a\tau^c$ on the solution compared to the $C^\alpha$ norm of the right-hand side. Given this IBP procedure, it is convenient to set: 
    \begin{equation}\label{initialdata}
            G=0, \ H=\tau f^{[1]}\mathbb{T} F|_{t=t_0}
    \end{equation}
    as initial data for \eqref{rewritinggr}.
    \item The restriction on  the time interval on which we solve the equation, also guarantees that:
        $$|t-t_0|||z^\pm||_{1+\alpha}\lesssim \tau^c \lambda_q\ell^{-\alpha}\delta_q^{1/2}=\left(\frac{\ell^{-1}}{\lambda_{q+1}}\right)^\alpha$$
and given any $\alpha,\ b$ we can always choose $a$ sufficiently large to reabsorb the implicit constant, and bound:
$$|t-t_0|||z^\pm||_{1+\alpha}\leq 1.$$
This will be used several times in the Duhamel formula + Gr{\"o}nwall type of arguments we are about to perform, and we will do so without further mention.
\end{itemize}
Given any solution $\Theta$, we split it into two different components,
$\Theta=\Theta^p+\Theta^h$. The rationale behind this is the following:
\begin{itemize}
    \item $\Theta^p$ will be the principal part of the solution, taking care of the forcing term and the initial conditions. Estimates will be done explicitly.
    \item $\Theta^h$ will be the part of the solution taking care of the homogeneous problem with zero initial conditions. Estimates will be implicit and given by subsequent Gr{\"o}nwall arguments.
\end{itemize}
More precisely, they solve:
\begin{equation}\label{splitproblems}
    \begin{cases}
        \mathcal{A}^-\mathcal{A}^+\Theta^h=-\mathbb{H}_1[\Theta]-\mathbb{H}_2[\mathcal{A}^+\Theta,\mathcal{A}^-\Theta],\\
        \Theta^h|_{t=t_0}=0,\\
        \mathcal{A}^+\Theta^h|_{t=t_0}=0
    \end{cases} \ \text{ and } \ 
        \begin{cases}
        \mathcal{A}^-\mathcal{A}^+\Theta^p=f\mathbb{T} F,\\
        \Theta^{p}|_{t=t_0}=0,\\
        \mathcal{A}^+\Theta^{p}|_{t=t_0}=\tau f^{[1]}\mathbb{T}F|_{t=t_0}.
    \end{cases}
\end{equation}
We remark here that for any solution of the above problem, we have:
$$\mathcal{A}^-\Theta|_{t=t_0}=\mathcal{A}^+\Theta|_{t=t_0}-2B\cn\Theta|_{t=t_0}$$
and thus for initial conditions such that $\Theta\equiv0$ at $t=t_0$ one has $B\cn\Theta|_{t=t_0}=v\cn\Theta|_{t=t_0}=0$
and we conclude that:
$$\Theta|_{t=t_0}=0 \Longrightarrow \mathcal{A}^-\Theta|_{t=t_0}=\mathcal{A}^+\Theta|_{t=t_0}=\partial_t\Theta|_{t=t_0}.$$
We now set $X^{\pm}, \ \Phi^\pm$ to be the Lagrangian and inverse flow map of $z^\pm=v\pm B$, respectively. That is, we solve:
\begin{equation*}
    \begin{cases}
        \partial_tX^\pm=z^\pm(X^\pm),\\
        X^\pm|_{t=t_0}=\IId
    \end{cases} \ \text{ and } \
    \begin{cases}
        \partial_t\Phi^\pm+z^\pm\cn \Phi^\pm=0,\\
        \Phi^\pm|_{t=t_0}=\IId.
    \end{cases}
\end{equation*}
It follows from the Standing Assumptions \ref{standingass}, Propositions \ref{standardlagrangianestimate}, \ref{standardtransportestimate}  and the remarks above that:
\begin{equation}\label{flowmapsgalb}
\begin{split}
    &||\DD X^\pm-\IId||_{\alpha}, \ ||\DD \Phi^\pm-\IId||_{\alpha} \leq 1\\
    &||\DD X^\pm||_{r+\alpha}, \ ||\DD\Phi^\pm||_{r+\alpha}\lesssim \lambda_q^r\lambda_q^{[r-(\underline{r}-1)]^+(b-1)\gamma_\ell} \ \text{ for } \ r\geq 1.
\end{split}    
\end{equation}

\noindent The goal of this Section is to prove the following Lemma.
\begin{lemma}[A Priori Estimates]\label{fullestimates} Let $r\geq 0$ an integer. Under the Standing Assumptions \ref{standingass} we have:
    \begin{subequations}
        \begin{align}
            &||\Theta||_{r+\alpha}\lesssim \lambda_q^{r}\lambda_q^{[r-(\underline{r}-2)]^+(b-1)\gamma_\ell}\ell^{-\alpha}\tau^a\tau^c\delta_{q+1},\\
            &||\mathcal{A}^\pm \Theta||_{r+\alpha}\lesssim \lambda_q^{r}\lambda_q^{[r-(\underline{r}-2)]^+(b-1)\gamma_\ell}\ell^{-\alpha}\tau^a\delta_{q+1},\\
            &||\partial_t^j\mathcal{A}^+\mathcal{A}^-\Theta||_{r+\alpha},||\partial_t^j(\mathcal{A}^\pm)^2\Theta||_{r+\alpha}\lesssim \lambda_q^{r+j}\lambda_q^{[r+j-(\underline{r}-2)]^+(b-1)\gamma_\ell}\ell^{-\alpha}\delta_{q+1} \ \text{ for } \ j=0,1.
        \end{align}
    \end{subequations}
   Moreover, on 
    $$I^{cut}=[t_0-\tau^c, \ t_0-2/3\tau^c]\cup[t_0+2/3\tau^c, \ t_0+\tau^c]$$ 
    we have the improved second-order transport estimate:
    \begin{equation}
        \begin{split}
            &||\mathcal{A}^+\mathcal{A}^-\Theta|_{I^{cut}}||_{r+\alpha}\lesssim \lambda_q^{r}\lambda_q^{[r-(\underline{r}-2)]^+(b-1)\gamma_\ell}\ell^{-\alpha}(\tau^a/\tau^c)\delta_{q+1}.\\
        \end{split}
    \end{equation}
    The implicit constants depend on $r, \ \alpha$ and the ones in the data from the Standing Assumptions \ref{standingass}.
\end{lemma}
\begin{remark}\label{lossapriori}
    As expected, given the second-order nature of the equation, one loses two good derivatives from the $\underline{r}$ we assume in the Standing Assumptions \ref{standingass}. This loss, however, occurs because of the coefficients in the operator $\mathbb{H}_1$, see Lemma \ref{operatornorms} and not because of the second-order transport. To avoid an additional loss in the special transport $(\mathcal{A}^\pm)^2\Theta$ estimate, we fix the same $\underline{r}$ for the transport estimate in \ref{standingass}, although a more natural choice would be $\underline{r}-1$. In the applications, this will not matter as $R, \ B, \ v$ enjoy much better lower bounds compared to $F$, and $F$ has this property by construction, see Lemma \ref{slowcoeffgalbrun}. 
\end{remark}
\begin{proof}[Proof of Lemma \ref{fullestimates}] We will deal with each part of the solution separately and leave the `special' second-order transport $(\mathcal{A}^\pm)^2$ estimates at the end. Those are not coming from the equation itself but follow as a consequence of the transport bounds we assume on the vector fields and forcings. At the very end, we handle the pure time-derivative bound. In the following, $||\cdot||_{r+\alpha}$ will mean the $C^{r+\alpha}_x$ H{\"o}lder norm for each fixed time, moreover for notational convenience we always assume $t_0\leq t<t_0+\tau^c$; the estimates for $t_0+\tau^c<t\leq t_0$ can be done similarly, as the Duhamel representation formulas we are about to use do not change.


\noindent \textbf{A Priori Estimates - $\Theta^p$.} Representing $\Theta^p$ via the Duhamel formula, see \eqref{duhamel} for the scalar case applied twice, we obtain:
\begin{equation*}
    \begin{split}
        \Theta^p&=\int_{t_0}^tH(\Phi^-_s\circ X_s^+\circ \Phi_t^+)\dd s+\int_{t_0}^t\left[\int_{t_0}^s(f\mathbb{T} F)(X_{s'}^-\circ\Phi^-_s\circ X_s^+\circ \Phi_t^+,s')\dd s'\right]\dd s.\\
    \end{split}
\end{equation*}
Set $\Psi=\Phi^-_s\circ X_s^+\circ \Phi_t^+$. We look first at the term coming from the forcing:
\begin{equation*}
    \begin{split}
        &\int_{t_0}^t\left[\int_{t_0}^s(f\mathbb{T}F)(X_{s'}^-\circ\Psi,s')\dd s'\right]\dd s\\
        &=\int_{t_0}^t\left[\tau^a (f^{[1]}\mathbb{T}F)(X_{s'}^-\circ\Psi,s')|_{s'=t_0}^{s'=s}\right]\dd s-\int_{t_0}^t\left[\int_{t_0}^s\tau^a (f^{[1]}\mathcal{A}^-\mathbb{T}F)(X_{s'}^-\circ\Psi,s')\dd s'\right]\dd s\\
        &=\int_{t_0}^t\left[\tau^a (f^{[1]}\mathbb{T}F)(X_s^+\circ \Phi_t^+,s)-\tau^a (f^{[1]}\mathbb{T}F)(\Psi,t_0)\right]\dd s-\int_{t_0}^t\left[\int_{t_0}^s\tau^a (f^{[1]}\mathcal{A}^-\mathbb{T}F)(X_{s'}^-\circ\Psi,s')\dd s'\right]\dd s\\
        &=\int_{t_0}^t\left[\tau^a (f^{[1]}\mathbb{T}F)(X_s^+\circ \Phi_t^+,s)-H(\Psi)\right]\dd s-\int_{t_0}^t\left[\int_{t_0}^s\tau^a (f^{[1]}\mathcal{A}^-\mathbb{T}F)(X_{s'}^-\circ\Psi,s')\dd s'\right]\dd s\\
    \end{split}
\end{equation*}
and we conclude that:
\begin{equation*}
    \begin{split}
        \Theta^p&=\int_{t_0}^t\left[\tau^a (f^{[1]}\mathbb{T}F)(X_s^+\circ \Phi_t^+,s)\right]\dd s-\int_{t_0}^t\left[\int_{t_0}^s\tau^a (f^{[1]}\mathcal{A}^-\mathbb{T}F)(X_{s'}^-\circ\Phi^-_s\circ X_s^+\circ \Phi_t^+,s')\dd s'\right]\dd s.\\
    \end{split}
\end{equation*}
This cancellation is why we chose \eqref{initialdata} as the initial data. With this representation at hand, the bounds in the Standing Assumptions \ref{standingass} and \eqref{flowmapsgalb}, the composition estimates in Propositions \ref{compestimates}, and Proposition \ref{czstuff} to deal with $\mathbb{T}$, we can compute:
\begin{equation}\label{thetap1}
    \begin{split}
        &||\Theta^p||_{r+\alpha}\\
        &\lesssim\tau^c\tau^a \sup_{|s-t_0|,|t-t_0|\leq \tau^c}||(f^{[1]}\mathbb{T}F)(X_s^+\circ \Phi_t^+,s)||_{r+\alpha}\\
        &+(\tau^c)^2\tau^a \sup_{|s-t_0|,|s'-t_0|,|t-t_0|\leq \tau^c}||(f^{[1]}\mathcal{A}^-\mathbb{T}F)(X_{s'}^-\circ\Psi,s')||_{r+\alpha}\\
        &\overset{r=0}{\leq} \tau^c\tau^a \sup_{|s-t_0|\leq \tau^c}||(f^{[1]}\mathbb{T}F)(\cdot,s)||_{r+\alpha}+(\tau^c)^2\tau^a \sup_{|s'-t_0|\leq \tau^c}||(f^{[1]}\mathcal{A}^-\mathbb{T}F)(\cdot,s')||_{r+\alpha}\\
        &\overset{r\geq1}{\lesssim} \tau^c\tau^a\sup_{|s-t_0|,|t-t_0|\leq \tau^c}\left[||f^{[1]}\mathbb{T}F||_{1+\alpha}||X_s^+\circ \Phi_t^+||_{r+\alpha}+||f^{[1]}\mathbb{T}F||_{r+\alpha}||X_s^+\circ \Phi_t^+||_{1+\alpha}^r\right]||X_s^+\circ \Phi_t^+||_{1}^\alpha\\
        &+(\tau^c)^2\tau^a\sup_{|s-t_0|,|s'-t_0|,|t-t_0|\leq \tau^c}||f^{[1]}\mathcal{A}^-\mathbb{T}F||_{1+\alpha}||X_{s'}^-\circ\Psi||_{r+\alpha}||X_{s'}^-\circ\Psi||_{1+\alpha}^r||X_{s'}^-\circ\Psi||_{1}^\alpha\\
        &+(\tau^c)^2\tau^a\sup_{|s-t_0|,|s'-t_0|,|t-t_0|\leq \tau^c}||f^{[1]}\mathcal{A}^-\mathbb{T}F||_{r+\alpha}||X_{s'}^-\circ\Psi||_{1+\alpha}^r||X_{s'}^-\circ\Psi||_{1}^\alpha\\
        &\overset{r\geq0}{\lesssim} \lambda_q^r\lambda_{q}^{[r-(\underline{r}-1)]^+(b-1)\gamma_\ell}\ell^{-\alpha}\tau^c\tau^a\delta_{q+1}
    \end{split}
\end{equation}
where we used $\sup_t |f^{[1]}|\lesssim 1$ to get rid of $f^{[1]}$ and the computation:
\begin{equation*}
    \begin{split}
        ||f^{[1]}\mathcal{A}^-\mathbb{T}F||_{r+\alpha}&\lesssim ||\mathbb{T}\mathcal{A}^- F||_{r+\alpha}+||[\mathcal{A}^-,\mathbb{T}]F||_{r+\alpha}\\
        &\lesssim ||\mathcal{A}^- F||_{r+\alpha}+||z^\pm||_{r+1+\alpha}||F||_{\alpha}+||z^\pm||_{1+\alpha}||F||_{r+\alpha}\\
    \end{split}
\end{equation*}
together with $\tau^c\ell^{-\alpha}\lambda_q\delta_q^{1/2}\leq 1$. 

\noindent \textit{Transport Estimates.} Applying the operator $\mathcal{A}^+$ to this representation formula and renaming the dummy variable $s'$ into $s$ we see that:
\begin{equation}\label{explicit+}
    \begin{split}
        \mathcal{A}^+\Theta^p&=\tau^a f^{[1]}\mathbb{T}F-\int_{t_0}^t\tau^a (f^{[1]}\mathcal{A}^-\mathbb{T}F)(X_{s}^-\circ\Phi_{t}^-,s)\dd s.\\
    \end{split}
\end{equation}
We claim a similar representation for the $\mathcal{A}^-$ transport. The commutativity $\mathcal{A}^+\mathcal{A}^-=\mathcal{A}^-\mathcal{A}^+$, the Duhamel representation formula, together with the fact remarked above that: 
$$\mathcal{A}^-\Theta^p|_{t=t_0}=\mathcal{A}^+\Theta^p|_{t=t_0}=H$$ 
and repeating the integration by parts, gives indeed:
\begin{equation*}
    \begin{split}
        \mathcal{A}^{-}\Theta^p&=H(\Phi_t^+)+\int_{t_0}^t(f \mathbb{T}F)(X_{s}^+\circ\Phi_{t}^+,s)\dd s\\
        &=H(\Phi_t^+)+\tau^a(f^{[1]} \mathbb{T}F)(X_{s}^+\circ\Phi_{t}^+,s)|_{s=t_0}^{s=t}-\int_{t_0}^t\tau^a f^{[1]}(s)( \mathcal{A}^+\mathbb{T}F)(X_{s}^+\circ\Phi_{t}^+,s)\dd s\\
        &=\tau^a f^{[1]}\mathbb{T}F-\int_{t_0}^t\tau^a (f^{[1]}\mathcal{A}^+\mathbb{T}F)(X_{s}^+\circ\Phi_{t}^+,s)\dd s.\\
    \end{split}
\end{equation*}
With the representations at hand, arguing as in \eqref{thetap1} above, we can bound:
\begin{equation*}
    \begin{split}
        ||\mathcal{A}^{\pm}\Theta^p||_{r+\alpha}&\lesssim\tau^a ||f^{[1]}\mathbb{T}F||_{r+\alpha}+\tau^a\tau^c \sup_{|s-t_0|\leq \tau^c}||(f^{[1]}\mathcal{A}^\mp\mathbb{T}F)(X_{s}^\mp\circ\Phi_{t}^\mp,s)||_{r+\alpha}\\
        &\lesssim \lambda_q^r\lambda_{q}^{[r-(\underline{r}-1)]^+(b-1)\gamma_\ell}\ell^{-\alpha}\tau^a\delta_{q+1}.
    \end{split}
\end{equation*}
The standard second-order transport estimates can be read directly from the equation, given the commutation of the Alfv\'en operators, indeed
\begin{equation*}
        ||\mathcal{A}^+\mathcal{A}^-\Theta^p||_{r+\alpha}=||\mathcal{A}^-\mathcal{A}^+\Theta^p||_{r+\alpha}=||f\mathbb{T} F||_{r+\alpha}\lesssim ||\mathbb{T} F||_{r+\alpha}\lesssim ||F||_{r+\alpha}\lesssim \lambda_q^r\lambda_{q}^{[r-\underline{r}]^+(b-1)\gamma_\ell}\ell^{-\alpha}\delta_{q+1}.
\end{equation*}
Additionally, note that this is zero on $I^{cut}$ as $F$ is according to the Standing Assumptions \ref{standingass}.

\noindent \textit{Special Transport.} We are left with $(\mathcal{A}^\pm)^2\Theta$. We will do the $+$ estimate the other one can be deduced similarly. Using again the Alfv\'en commutation we see that $(\mathcal{A}^+)^2\Theta$ solves:
\begin{equation*}
        \mathcal{A}^-(\mathcal{A}^+)^2\Theta^p=\mathcal{A}^+(f\mathbb{T} F)=\frac{1}{\tau^a}f'\mathbb{T} F+f\mathcal{A}^+\mathbb{T} F\\
\end{equation*}
with initial conditions:
\begin{equation}\label{ICspecial}
    \begin{split} 
    (\mathcal{A}^+)^2\Theta^{p}|_{t=t_0}&=\mathcal{A}^-\mathcal{A}^+\Theta^{p}|_{t=t_0}+2B\cn \mathcal{A}^+\Theta^{p}|_{t=t_0}\\
    &=f\mathbb{T}F|_{t=t_0}+2B|_{t=t_0}\cn [\mathcal{A}^+\Theta^{p}|_{t=t_0}]\\
    &=f\mathbb{T}F|_{t=t_0}+(\mathcal{A}^+-\mathcal{A}^-)|_{t=t_0}[\tau^a f^{[1]}\mathbb{T}F]|_{t=t_0}\\
    &=f\mathbb{T}F|_{t=t_0}+\tau^a f^{[1]}(\mathcal{A}^+-\mathcal{A}^-)\mathbb{T}F|_{t=t_0}.\\
    \end{split}
\end{equation}
Representing the solution by the Duhamel formula and integrating by parts the `bad term', we get:
\begin{equation*}
    \begin{split}
        (\mathcal{A}^+)^2\Theta^p=f\mathbb{T}F+\tau^a f^{[1]}(\mathcal{A}^+-\mathcal{A}^-)\mathbb{T}F|_{t=t_0}(\Phi^-_t)+\int_{t_0}^t[f(\mathcal{A}^+-\mathcal{A}^-)\mathbb{T} F](X^-_s\circ\Phi^-_t,s)\dd s
    \end{split}
\end{equation*}
now estimating as above, we conclude:
\begin{equation*}
    \begin{split}
        ||(\mathcal{A}^+)^2\Theta^p||_{r+\alpha}&\lesssim\sup_{|s-t_0|,|t-t_0|\leq \tau^c}\left[||F||_{r+\alpha}+\tau^a ||(\mathcal{A}^+-\mathcal{A}^-)\mathbb{T}F|_{t=t_0}(\Phi^\mp_t)||_{r+\alpha}+\tau^c||[(\mathcal{A}^+-\mathcal{A}^-)\mathbb{T} F](X^-_s\circ\Phi^-_t)||_{r+\alpha}\right]\\
        &\lesssim \lambda_q^r\lambda_{q}^{[r-(\underline{r}-1)]^+(b-1)\gamma_\ell}\ell^{-\alpha}\delta_{q+1}.
    \end{split}
\end{equation*}

\noindent Gathering all the results we get:
    \begin{subequations}\label{galbrunprincipal}
        \begin{align}
            &||\Theta^p||_{r+\alpha}\lesssim \lambda_q^r\lambda_{q}^{[r-(\underline{r}-1)]^+(b-1)\gamma_\ell}\ell^{-\alpha}\tau^a\tau^c\delta_{q+1},\\
            &||\mathcal{A}^\pm \Theta^p||_{r+\alpha}\lesssim \lambda_q^{r}\lambda_{q}^{[r-(\underline{r}-1)]^+(b-1)\gamma_\ell}\ell^{-\alpha}\tau^a\delta_{q+1},\\
            &||\mathcal{A}^+\mathcal{A}^-\Theta^p||_{r+\alpha},||(\mathcal{A}^\pm)^2\Theta^p||_{r+\alpha}\lesssim \lambda_q^{r}\lambda_{q}^{[r-(\underline{r}-1)]^+(b-1)\gamma_\ell}\ell^{-\alpha}\delta_{q+1},\\
            &\mathcal{A}^+\mathcal{A}^-\Theta^p|_{I^{cut}}\equiv0
        \end{align}
    \end{subequations}
for $r\geq 0$.

\noindent \textbf{A Priori Estimates - $\Theta^h$.} We recall that $\Theta^h$ solves the problem \eqref{splitproblems} and argue similarly as above. The Duhamel formula, the commutation $\mathcal{A}^+\mathcal{A}^-=\mathcal{A}^-\mathcal{A}^+$ and the zero initial conditions $\mathcal{A}^+\Theta^h|_{t=t_0}=\mathcal{A}^-\Theta^h|_{t=t_0}=0$ allow us to write the following representation formula:
\begin{equation*}
    \begin{split}
        \mathcal{A}^{\pm}\Theta^h&=-\int_{t_0}^t\left(\mathbb{H}_1[\Theta]+\mathbb{H}_2[\mathcal{A}^+\Theta,\mathcal{A}^-\Theta]\right)(X^\mp_s\circ\Phi^\mp_t,s)\dd s\\
        &=-\int_{t_0}^t\left(\mathbb{H}_1[\Theta^p+\Theta^h]+\mathbb{H}_2[\mathcal{A}^+\Theta^p+\mathcal{A}^+\Theta^h,\mathcal{A}^-\Theta^p+\mathcal{A}^-\Theta^h]\right)(X^\mp_s\circ\Phi^\mp_t,s)\dd s\\
        &=\underbrace{-\int_{t_0}^t\left(\mathbb{H}_1[\Theta^p]\right)(X^\mp_s\circ\Phi^\mp_t,s)\dd s-\int_{t_0}^t\left(\mathbb{H}_2[\mathcal{A}^+\Theta^p,\mathcal{A}^-\Theta^p]\right)(X^\mp_s\circ\Phi^\mp_t,s)\dd s}_{I_1^\pm}\\
        &\underbrace{-\int_{t_0}^t\left(\mathbb{H}_1[\Theta^h]\right)(X^\mp_s\circ\Phi^\mp_t,s)\dd s-\int_{t_0}^t\left(\mathbb{H}_2[\mathcal{A}^+\Theta^h,\mathcal{A}^-\Theta^h]\right)(X^\mp_s\circ\Phi^\mp_t,s)\dd s}_{I_2^\pm}.\\
    \end{split}
\end{equation*}
We now bound $I_1^\pm$ explicitly as they come from the principal part of the solution. Using the properties of $\mathbb{H}_2$ from Lemma \ref{operatornorms}, the estimates in \eqref{galbrunprincipal}, the composition estimates in Proposition \ref{compestimates} and the bounds on the flow maps in \eqref{flowmapsgalb}, we compute:
\begin{equation}\label{I_1}
    \begin{split}
        ||I_1^\pm||_{r+\alpha}&\overset{r=0}{\lesssim}\int_{t_0}^t||\mathbb{H}_1[\Theta^p]||_{\alpha}\dd s+\int_{t_0}^t||\mathbb{H}_2[\mathcal{A}^+\Theta^p,\mathcal{A}^-\Theta^p]||_{\alpha}\dd s\\
        &\lesssim \int_{t_0}^th_1^0||\Theta^p||_{\alpha}\dd s+\int_{t_0}^th_2^0\left(||\mathcal{A}^+\Theta^p||_\alpha+||\mathcal{A}^-\Theta^p]||_{\alpha}\right)\dd s\\
        &\lesssim \ell^{-2\alpha}(\tau^c)^2\tau^a \lambda_q^2\delta_q\delta_{q+1}+\ell^{-2\alpha}\tau^c\tau^a \lambda_q\delta_q^{1/2}\delta_{q+1}\\
        &\lesssim\ell^{-\alpha}\tau^a\delta_{q+1}\\
        &\overset{r\geq 1}{\lesssim} \lambda_q^{r-1}\lambda_q^{[r-\underline{r}]^+(b-1)\gamma_\ell}\int_{t_0}^t||\mathbb{H}_1[\Theta^p]||_{1+\alpha}\dd s+\int_{t_0}^t||\mathbb{H}_1[\Theta^p]||_{r+\alpha}\dd s\\
        &+\lambda_q^{r-1}\lambda_q^{[r-\underline{r}]^+(b-1)\gamma_\ell}\int_{t_0}^t||\mathbb{H}_2[\mathcal{A}^+\Theta^p,\mathcal{A}^-\Theta^p]||_{1+\alpha}\dd s+\int_{t_0}^t||\mathbb{H}_2[\mathcal{A}^+\Theta^p,\mathcal{A}^-\Theta^p]||_{r+\alpha}\dd s\\
        &\lesssim \lambda_q^{r-1}\lambda_q^{[r-\underline{r}]^+(b-1)\gamma_\ell}\int_{t_0}^th_1^{1}||\Theta^p||_{\alpha}+h_1^0||\Theta^p||_{1+\alpha}\dd s+\int_{t_0}^t h_1^{r}||\Theta^p||_{\alpha}+h_1^0||\Theta^p||_{r+\alpha}\dd s\\
        &+\lambda_q^{r-1}\lambda_q^{[r-\underline{r}]^+(b-1)\gamma_\ell}\int_{t_0}^th_2^{1}(||\mathcal{A}^+\Theta^p||_{\alpha}+||\mathcal{A}^-\Theta^p||_{\alpha})+h_2^0(||\mathcal{A}^+\Theta^p||_{1+\alpha}+||\mathcal{A}^-\Theta^p||_{1+\alpha})\dd s\\
        &+\int_{t_0}^t h_2^{r}(||\mathcal{A}^+\Theta^p||_{\alpha}+||\mathcal{A}^-\Theta^p||_{\alpha})+h_2^0(||\mathcal{A}^+\Theta^p||_{r+\alpha}+||\mathcal{A}^-\Theta^p||_{r+\alpha})\dd s\\
        &\lesssim \lambda_q^{r+2}\lambda_q^{[r-(\underline{r}-2)]^+(b-1)\gamma_\ell}\ell^{-2\alpha}(\tau^c)^2\tau^a\delta_q\delta_{q+1}+\lambda_q^{r+1}\lambda_q^{[r-(\underline{r}-1)]^+(b-1)\gamma_\ell}\ell^{-2\alpha}\tau^c\tau^a\delta_q^{1/2}\delta_{q+1}\\
        &\lesssim \lambda_q^{r}\lambda_q^{[r-(\underline{r}-2)]^+(b-1)\gamma_\ell}\ell^{-\alpha}\tau^a \delta_{q+1}.
    \end{split}
\end{equation}
The integrals $I_2^\pm$ will be estimated implicitly via a double Gr{\"o}nwall argument. Computing as above, we first deduce:
\begin{equation}\label{I_2}
    \begin{split}
        ||I^\pm_2||_{r+\alpha}&\overset{r=0}{\lesssim}\int_{t_0}^t||\mathbb{H}_1[\Theta^h]||_{\alpha}\dd s+\int_{t_0}^t||\mathbb{H}_2[\mathcal{A}^+\Theta^h,\mathcal{A}^-\Theta^h]||_{\alpha}\dd s\\
        &\lesssim \int_{t_0}^t \lambda_q^{2}\ell^{-\alpha}\delta_q||\Theta^h||_{\alpha}+\lambda_q\ell^{-\alpha}\delta_q^{1/2}(||\mathcal{A}^+\Theta^h||_{\alpha}+||\mathcal{A}^-\Theta^h||_{\alpha})\dd s\\
        &\overset{r\geq1}{\lesssim} \lambda_q^{r-1}\lambda_q^{[r-\underline{r}]^+(b-1)\gamma_\ell}\int_{t_0}^th_1^{1}||\Theta^h||_{\alpha}+h_1^0||\Theta^h||_{1+\alpha}\dd s+\int_{t_0}^t h_1^{r}||\Theta^h||_{\alpha}+h_1^0||\Theta^h||_{r+\alpha}\dd s\\
        &+\lambda_q^{r-1}\lambda_q^{[r-\underline{r}]^+(b-1)\gamma_\ell}\int_{t_0}^th_2^{1}(||\mathcal{A}^+\Theta^h||_{\alpha}+||\mathcal{A}^-\Theta^h||_{\alpha})+h_2^0(||\mathcal{A}^+\Theta^h||_{1+\alpha}+||\mathcal{A}^-\Theta^h||_{1+\alpha})\dd s\\
        &+\int_{t_0}^t h_2^{r}(||\mathcal{A}^+\Theta^h||_{\alpha}+||\mathcal{A}^-\Theta^h||_{\alpha})+h_2^0(||\mathcal{A}^+\Theta^h||_{r+\alpha}+||\mathcal{A}^-\Theta^h||_{r+\alpha})\dd s\\
        &\lesssim\lambda_q^{r+2}\lambda_q^{[r-(\underline{r}-2)]^+(b-1)\gamma_\ell}\ell^{-\alpha}\delta_q\int_{t_0}^t||\Theta^h||_{\alpha}\dd s\\
        &+\lambda_q^{r+1}\lambda_q^{[r-\underline{r}]^+(b-1)\gamma_\ell}\ell^{-\alpha}\delta_q\int_{t_0}^t||\Theta^h||_{1+\alpha}\dd s+\lambda_q^{2}\ell^{-\alpha}\delta_q\int_{t_0}^t||\Theta^h||_{r+\alpha}\dd s\\
        &+\lambda_q^{r+1}\lambda_q^{[r-(\underline{r}-1)]^+(b-1)\gamma_\ell}\ell^{-\alpha}\delta_q^{1/2}\int_{t_0}^t(||\mathcal{A}^+\Theta^h||_{\alpha}+||\mathcal{A}^-\Theta^h||_{\alpha})\dd s\\
        &+\lambda_q\ell^{-\alpha}\delta_q^{1/2}\int_{t_0}^t(||\mathcal{A}^+\Theta^h||_{r+\alpha}+||\mathcal{A}^-\Theta^h||_{r+\alpha})\dd s\\
        &+\lambda_q^{r}\lambda_q^{[r-\underline{r}]^+(b-1)\gamma_\ell}\ell^{-\alpha}\delta_q^{1/2}\int_{t_0}^t(||\mathcal{A}^+\Theta^h||_{1+\alpha}+||\mathcal{A}^-\Theta^h||_{1+\alpha})\dd s.
    \end{split}
\end{equation}
We now complete the argument for $r=0$, and use induction afterwards. 

\noindent \textit{Base case: $r=0$.} Collecting the bounds in \eqref{I_1}, \eqref{I_2} and summing the $\pm$ estimates we obtain:
\begin{equation*}
    \begin{split}
        \sum_{\sigma \in \{+,-\}}||\mathcal{A}^{\sigma}\Theta^h||_{\alpha}&\lesssim\ell^{-\alpha}\tau^a \delta_{q+1} +\int_{t_0}^t \lambda_q^{2}\ell^{-\alpha}\delta_q||\Theta^h||_{\alpha}+\lambda_q\ell^{-\alpha}\delta_q^{1/2}\sum_{\sigma \in \{+,-\}}||\mathcal{A}^{\sigma}\Theta^h||_{\alpha}\dd s,\\
    \end{split}
\end{equation*}
an application of Gr{\"o}nwall's Lemma together with the standing assumption $|t-t_0|\leq \tau^c$ leads to:
\begin{equation}\label{basecase}
    \begin{split}
        \sum_{\sigma \in \{+,-\}}||\mathcal{A}^{\sigma}\Theta^h||_{\alpha}&\lesssim  \left[\tau^a \ell^{-\alpha}\delta_{q+1}+\int_{t_0}^t\lambda_q^{2}\ell^{-\alpha}\delta_q||\Theta^h||_{\alpha}\dd s\right] e^{\tau^c \lambda_q\ell^{-\alpha}\delta_q^{1/2}}\\
        &\lesssim\tau^a \ell^{-\alpha}\delta_{q+1}+\lambda_q^{2}\ell^{-\alpha}\delta_q\int_{t_0}^t||\Theta^h||_{\alpha}\dd s.
    \end{split}
\end{equation}
We now integrate along the flow, using the fact that $\Theta^h|_{t=t_0}=0$. Namely, we write:
\begin{equation}\label{alongtheflow}
    \begin{split}
        \Theta^h(t)&=\Theta^h(X^+_t\circ \Phi^+_t,t)=\Theta^h(X^+_t\circ \Phi^+_t,t)-\Theta^h(\Phi^+_t,t_0)\\
        &=\int_{t_0}^t\partial_s\left[(\Theta^h)(X^+_s\circ \Phi^+_t)\right]\dd s\\
        &=\int_{t_0}^t(\mathcal{A}^+\Theta^h)(X^+_s\circ \Phi^+_t)\dd s
    \end{split}
\end{equation}
and using \eqref{basecase} together with $|t-t_0|\leq \tau^c$ we deduce:
\begin{equation*}
    \begin{split}
        ||\Theta^h||_\alpha&\lesssim\ell^{-\alpha}\tau^c\tau^a \delta_{q+1}+\lambda_q^{2}\ell^{-\alpha}\delta_q\int_{t_0}^t \int_{t_0}^{s} ||\Theta^h||_{\alpha}\dd s'\dd s\\
        &\lesssim \ell^{-\alpha}\tau^c\tau^a \delta_{q+1}+\tau^c \lambda_q^{2}\ell^{-\alpha}\delta_q\int_{t_0}^t ||\Theta^h||_{\alpha}\dd s'\\
    \end{split}
\end{equation*}
and a final Gr\"onwall argument with subsequent use of \eqref{basecase} gives:
\begin{equation}
        ||\Theta^h||_\alpha\lesssim \ell^{-\alpha}\tau^c\tau^a \delta_{q+1}, \ \ \ \ ||\mathcal{A}^\pm\Theta^h||_\alpha\lesssim \ell^{-\alpha}\tau^a \delta_{q+1},
\end{equation}
this concludes the base case of the induction. 

\noindent \textit{Induction Step: $r\geq 1$.} Now, assume we proved for $1\leq r'\leq r-1$ that:
\begin{equation*}
    \begin{split}
        &||\Theta^h||_{r'+\alpha}\lesssim \lambda_q^{r'}\lambda_q^{[r'-(\underline{r}-2)]^+(b-1)\gamma_\ell}\ell^{-\alpha}\tau^a\tau^c\delta_{q+1},\\
        &||\mathcal{A}^\pm \Theta^h||_{r'+\alpha}\lesssim \lambda_q^{r'}\lambda_q^{[r'-(\underline{r}-2)]^+(b-1)\gamma_\ell}\ell^{-\alpha}\tau^a\delta_{q+1}  \\
    \end{split}
\end{equation*}
where the implicit constants depend on $\alpha$ and $r'$.  From \eqref{I_1}, \eqref{I_2} and arguing as above we deduce:
\begin{equation*}
    \begin{split}
        \sum_{\sigma \in \{+,-\}}||\mathcal{A}^{\sigma}\Theta^h||_{r+\alpha}&\leq \sum_{\sigma \in \{+,-\}, \ i=1,2}||I_i^\sigma||_{r+\alpha}\\
        &\lesssim \lambda_q^{r}\lambda_q^{[r-(\underline{r}-2)]^+(b-1)\gamma_\ell}\ell^{-\alpha}\tau^a \delta_{q+1}+\lambda_q\ell^{-\alpha}\delta_q^{1/2}\sum_{\sigma \in \{+,-\}}\int_{t_0}^t||\mathcal{A}^{\sigma}\Theta^h||_{r+\alpha}\dd s\\
        &+\lambda_q^{2}\ell^{-\alpha}\delta_q\int_{t_0}^t||\Theta^h||_{r+\alpha}\dd s,
    \end{split}
\end{equation*}
we now apply Gr{\"o}nwall lemma
\begin{equation}\label{dergn}
    \begin{split}
        \sum_{\sigma \in \{+,-\}}||\mathcal{A}^{\sigma}\Theta^h||_{r+\alpha}&\lesssim \left[\lambda_q^{r}\lambda_q^{[r-(\underline{r}-2)]^+(b-1)\gamma_\ell}\ell^{-\alpha}\tau^a \delta_{q+1}+\lambda_q^{2}\ell^{-\alpha}\delta_q\int_{t_0}^t||\Theta^h||_{r+\alpha}\dd s\right]e^{\tau^c\lambda_q\ell^{-\alpha}\delta_q^{1/2}}\\
        &\lesssim \lambda_q^{r}\lambda_q^{[r-(\underline{r}-2)]^+(b-1)\gamma_\ell}\ell^{-\alpha}\tau^a \delta_{q+1}+\lambda_q^{2}\ell^{-\alpha}\delta_q\int_{t_0}^t||\Theta^h||_{r+\alpha}\dd s\\
    \end{split}
\end{equation}
and integrating along the flow $\mathcal{A}^+\Theta^h$ as in \eqref{alongtheflow}, together with the bounds in \eqref{galbrunprincipal}, we deduce:
\begin{equation*}
    \begin{split}
        ||\Theta^h||_{r+\alpha}&\lesssim\lambda_q^{r}\lambda_q^{[r-(\underline{r}-2)]^+(b-1)\gamma_\ell}\ell^{-\alpha}\tau^c\tau^a \delta_{q+1}+\int_{t_0}^t||(\mathcal{A}^+\Theta^h)(X^+_s\circ \Phi^+_t)||_{r+\alpha}\dd s\\
        &\lesssim\lambda_q^{r}\lambda_q^{[r-(\underline{r}-2)]^+(b-1)\gamma_\ell}\ell^{-\alpha}\tau^c\tau^a \delta_{q+1}+\int_{t_0}^t||\mathcal{A}^+\Theta^h||_{r+\alpha}\dd s\\
        &\lesssim\lambda_q^{r}\lambda_q^{[r-(\underline{r}-2)]^+(b-1)\gamma_\ell}\ell^{-\alpha}\tau^c\tau^a \delta_{q+1}+\lambda_q^{2}\ell^{-\alpha}\delta_q\int_{t_0}^t\int_{t_0}^s||\Theta^h||_{r+\alpha}\dd s'\dd s\\
    &\lesssim\lambda_q^{r}\lambda_q^{[r-(\underline{r}-2)]^+(b-1)\gamma_\ell}\ell^{-\alpha}\tau^c\tau^a \delta_{q+1}+\lambda_q^{2}\ell^{-\alpha}\delta_q\tau^c\int_{t_0}^t||\Theta^h||_{r+\alpha}\dd s'\\
    \end{split}
\end{equation*}
and we conclude as above that:
\begin{equation}\label{thetah}
    ||\Theta^h||_{r+\alpha}\lesssim \lambda_q^{r}\lambda_q^{[r-(\underline{r}-2)]^+(b-1)\gamma_\ell}\ell^{-\alpha}\tau^c\tau^a \delta_{q+1}
\end{equation}
and 
\begin{equation}\label{firsttransportapriorih}
    ||\mathcal{A}^\pm\Theta^h||_{r+\alpha}\lesssim \lambda_q^{r}\lambda_q^{[r-(\underline{r}-2)]^+(b-1)\gamma_\ell}\ell^{-\alpha}\tau^a \delta_{q+1}.
\end{equation}

\noindent \textbf{A Priori Estimates - $\Theta$.} Since $\Theta=\Theta^p+\Theta^h$, collecting the bounds in \eqref{galbrunprincipal} and \eqref{thetah}, \eqref{firsttransportapriorih},  results in:
\begin{equation}\label{firsttransportapriori}
        \begin{split}
            &||\Theta||_{r+\alpha}\lesssim \lambda_q^r\lambda_{q}^{[r-(\underline{r}-2)]^+(b-1)\gamma_\ell}\ell^{-\alpha}\tau^a\tau^c\delta_{q+1},\\
            &||\mathcal{A}^\pm \Theta||_{r+\alpha}\lesssim \lambda_q^{r}\lambda_{q}^{[r-(\underline{r}-2)]^+(b-1)\gamma_\ell}\ell^{-\alpha}\tau^a\delta_{q+1}\\
        \end{split}
    \end{equation}
for $r\geq 0$.


\noindent \textbf{Special Transport $(\mathcal{A}^\pm)^2\Theta$.} We write the details for $(\mathcal{A}^+)^2\Theta$ only. Applying $\mathcal{A}^+$ to \eqref{rewritinggr} and commuting we obtain:
\begin{equation}\label{specialtransportrewriting}
    \begin{split}
        \frac{1}{\tau^a}f'\mathbb{T}F+f \mathcal{A}^+\mathbb{T}F&=\mathcal{A}^+f\mathbb{T}F\\
        &=\mathcal{A}^+\mathcal{A}^-\mathcal{A}^+\Theta+\mathcal{A}^+\mathbb{H}_1[\Theta]+\mathcal{A}^+\mathbb{H}_2[\mathcal{A}^+\Theta,\mathcal{A}^-\Theta]\\
        &=\mathcal{A}^-(\mathcal{A}^+)^2\Theta+\mathbb{H}_1[\mathcal{A}^+\Theta]+\mathbb{H}_2[(\mathcal{A}^+)^2\Theta,\mathcal{A}^+\mathcal{A}^-\Theta]\\
        &+[\mathcal{A}^+,\mathbb{H}_1]\Theta+[\mathcal{A}^+\mathbb{H}_2][\mathcal{A}^+\Theta,\mathcal{A}^-\Theta].
    \end{split}
\end{equation}
Let $$H=(\mathcal{A}^+)^2\Theta|_{t_0}=f \mathbb{T}F|_{t=t_0}+\tau^a(\mathcal{A}^+-\mathcal{A}^-)\mathbb{T}F|_{t=t_0}$$ 
be the initial condition for $(\mathcal{A}^\pm)^2\Theta$, where the second equality follows as in \eqref{ICspecial}. Integrating along the flow of $\mathcal{A}^-$ we obtain the following representation formula:
\begin{equation*}
    \begin{split}
        (\mathcal{A}^+)^2\Theta&=H(\Phi^-)-\int_{t_0}^t \left[\mathbb{H}_1[\mathcal{A}^+\Theta]+\mathbb{H}_2[(\mathcal{A}^+)^2\Theta,\mathcal{A}^+\mathcal{A}^-\Theta]\right](X_s^-\circ \Phi^-_t)\dd s\\
        &-\int_{t_0}^t\left[[\mathcal{A}^+,\mathbb{H}_1]\Theta+[\mathcal{A}^+\mathbb{H}_2][\mathcal{A}^+\Theta,\mathcal{A}^-\Theta]\right](X_s^-\circ \Phi^-_t)\dd s+\int_{t_0}^t \left[\frac{1}{\tau^a}f'\mathbb{T}F+f \mathcal{A}^+\mathbb{T}F\right](X_s^-\circ \Phi^-_t)\dd s\\
        &=\underbrace{\tau^a(\mathcal{A}^+-\mathcal{A}^-)\mathbb{T}F|_{t=t_0}(\Phi^-)+f\mathbb{T}F+\int_{t_0}^t \left[f (\mathcal{A}^+-\mathcal{A}^-)\mathbb{T}F-\mathbb{H}_1[\mathcal{A}^+\Theta]\right](X_s^-\circ \Phi^-_t)\dd s}_{I_1}\\
        &-\int_{t_0}^t \left[\mathbb{H}_2[(\mathcal{A}^+)^2\Theta,\mathcal{A}^+\mathcal{A}^-\Theta]\right](X_s^-\circ \Phi^-_t)\dd s-\underbrace{\int_{t_0}^t\left[[\mathcal{A}^+,\mathbb{H}_1]\Theta+[\mathcal{A}^+,\mathbb{H}_2][\mathcal{A}^+\Theta,\mathcal{A}^-\Theta]\right](X_s^-\circ \Phi^-_t)\dd s}_{I_2}.\\
    \end{split}
\end{equation*}
We focus on $I_1$ and $I_2$ separately.

\noindent \textit{Estimates on $I_2$.} The estimates on the commutator $[\mathcal{A}^\pm,\mathbb{H}_2], \ [\mathcal{A}^+,\mathbb{H}_1]$ involving the Alfv\'en transport derivatives of the coefficients of the operators $\mathbb{H}_1,\ \mathbb{H}_2$ were given in Lemma \ref{operatornorms}. We have:
\begin{equation*}
    \begin{split}
        ||I_2||_{r+\alpha}&\lesssim \int_{t_0}^t||\left[[\mathcal{A}^+,\mathbb{H}_1]\Theta+[\mathcal{A}^+,\mathbb{H}_2][\mathcal{A}^+\Theta,\mathcal{A}^-\Theta]\right](X_s^-\circ \Phi^-_t)||_{r+\alpha}\dd s\\
        &\overset{r=0}{\lesssim}\int_{t_0}^t||[\mathcal{A}^+,\mathbb{H}_1]\Theta+[\mathcal{A}^+,\mathbb{H}_2][\mathcal{A}^+\Theta,\mathcal{A}^-\Theta]||_{\alpha}\dd s\\
        &\lesssim\tau^c\lambda_q^2\ell^{-2\alpha}\delta_q(||\mathcal{A}^+\Theta||_\alpha+||\mathcal{A}^-\Theta||_\alpha)+\tau^c\lambda_q^3\ell^{-2\alpha}\delta_q^{3/2}||\Theta||_\alpha\\
        &\overset{r\geq 1}{\lesssim} \int_{t_0}^t||\left[[\mathcal{A}^+,\mathbb{H}_1]\Theta+[\mathcal{A}^+,\mathbb{H}_2][\mathcal{A}^+\Theta,\mathcal{A}^-\Theta]\right]||_{1+\alpha}||X_s^-\circ \Phi^-_t||_{r+\alpha}||(X_s^-\circ \Phi^-_t)||_{1}^\alpha\dd s\\
        &+ \int_{t_0}^t||\left[[\mathcal{A}^+,\mathbb{H}_1]\Theta+[\mathcal{A}^+,\mathbb{H}_2][\mathcal{A}^+\Theta,\mathcal{A}^-\Theta]\right]||_{r+\alpha}||(X_s^-\circ \Phi^-_t)||_{1+\alpha}^r||(X_s^-\circ \Phi^-_t)||_{1}^\alpha\dd s\\
        &\lesssim \lambda_q^{r-1}\lambda_q^{[r-\underline{r}]^+(b-1)\gamma_\ell} \tau^c\lambda_q\ell^{-\alpha}\delta_q^{1/2}\left(h_1^{1}||\Theta||_\alpha+h_1^0||\Theta||_{1+\alpha}+h_2^{1}(||\mathcal{A}^+\Theta||_\alpha+||\mathcal{A}^-\Theta||_\alpha)\right)\\
        &+\lambda_q^{r-1}\lambda_q^{[r-\underline{r}]^+(b-1)\gamma_\ell} \tau^c\lambda_q\ell^{-\alpha}\delta_q^{1/2}h_2^0\left(||\mathcal{A}^+\Theta||_{1+\alpha}+||\mathcal{A}^-\Theta||_{1+\alpha}\right)\\
        &+\tau^c\lambda_q\ell^{-\alpha}\delta_q^{1/2}\left(h_1^{r}||\Theta||_\alpha+h_1^0||\Theta||_{r+\alpha}\right)\\
        &+\tau^c\lambda_q\ell^{-\alpha}\delta_q^{1/2}\left[h_2^{r}(||\mathcal{A}^+\Theta||_\alpha+||\mathcal{A}^-\Theta||_\alpha)+h_2^0(||\mathcal{A}^+\Theta||_{r+\alpha}+||\mathcal{A}^-\Theta||_{r+\alpha})\right]\\
        &\overset{r\geq0}{\lesssim} \lambda_q^{r+1}\lambda_q^{[r-(\underline{r}-2)]^+(b-1)\gamma_\ell}\ell^{-2\alpha}\tau^a\delta_q^{1/2}\delta_{q+1}
    \end{split}
\end{equation*}

\noindent \textit{Estimates on $I_1$.} We first use the commutator estimates in Proposition \ref{czstuff} and the Alfv\'en transport estimates in the Standing Assumpions \ref{standingass} to bound: 
\begin{equation*}
    \begin{split}
        ||\mathcal{A}^\pm\mathbb{T}F||_{r+\alpha}&\leq||[\mathcal{A}^\pm, \mathbb{T}]F||_{r+\alpha}+||\mathbb{T}\mathcal{A}^\pm F||_{r+\alpha}\\
        &\lesssim||[\mathcal{A}^\pm, \mathbb{T}]F||_{r+\alpha}+||\mathcal{A}^\pm F||_{r+\alpha}\\
        &\lesssim ||z^\pm||_{r+1+\alpha}||F||_{\alpha}+||z^\pm||_{1+\alpha}||F||_{r+\alpha}+||\mathcal{A}^\pm F||_{r+\alpha}\\
        &\lesssim\lambda_q^{r}\lambda_q^{[r-(\underline{r}-1)]^+(b-1)\gamma_\ell}1/\tau^c\ell^{-\alpha}\delta_{q+1}
    \end{split}
\end{equation*}
and from this, arguing as above, we conclude:
\begin{equation*}
    \begin{split}
        ||I_1||_{r+\alpha}&\lesssim \tau^a||\mathcal{A}^+\mathbb{T}F|_{t=t_0}(\Phi^-)||_{r+\alpha}+||\mathbb{T}F||_{r+\alpha}+\int_{t_0}^t ||\left[f \mathcal{A}^-\mathbb{T}F+f \mathcal{A}^+\mathbb{T}F-\mathbb{H}_1[\mathcal{A}^+\Theta]\right](X_s^-\circ \Phi^-_t)||_{r+\alpha}\dd s\\
        &\lesssim \lambda_q^{r}\lambda_q^{[r-(\underline{r}-2)]^+(b-1)\gamma_\ell}\ell^{-\alpha}\delta_{q+1}\\
    \end{split}
\end{equation*}

\noindent \textit{Conclusion.} Gathering the estimates on $I_1, \ I_2$ and proceeding as before we obtain:
\begin{equation*}
    \begin{split}
        ||(\mathcal{A}^+)^2\Theta||_{r+\alpha}&\lesssim \lambda_q^{r}\lambda_q^{[r-(\underline{r}-2)]^+(b-1)\gamma_\ell}\ell^{-\alpha}\delta_{q+1}+\int_{t_0}^t ||\left[\mathbb{H}_2[(\mathcal{A}^+)^2\Theta,\mathcal{A}^+\mathcal{A}^-\Theta]\right](X_s^-\circ \Phi^-_t)||_{r+\alpha}\dd s\\
        &\overset{r=0}{\lesssim}
        \ell^{-\alpha}\delta_{q+1}+\lambda_q\ell^{-\alpha}\delta_q^{1/2}\int_{t_0}^t ||(\mathcal{A}^+)^2\Theta||_\alpha\dd s\\
        &\overset{r\geq 1}{\lesssim}
        \lambda_q^{r}\lambda_q^{[r-(\underline{r}-2)]^+(b-1)\gamma_\ell}\ell^{-\alpha}\delta_{q+1}\\
        &+\int_{t_0}^t\lambda_q^{r+1}\lambda_q^{[r-(\underline{r}-1)]^+(b-1)\gamma_\ell}\ell^{-\alpha}\delta_q^{1/2}||(\mathcal{A}^+)^2\Theta||_\alpha+\lambda_q^{r}\lambda_q^{[r-\underline{r}]^+(b-1)\gamma_\ell}\ell^{-\alpha}\delta_q^{1/2}||(\mathcal{A}^+)^2\Theta||_{1+\alpha}\dd s\\
        &+\int_{t_0}^t\lambda_q\ell^{-\alpha}\delta_q^{1/2}||(\mathcal{A}^+)^2\Theta||_{r+\alpha}\dd s\\
    \end{split}
\end{equation*}
and arguing inductively by means of Gr{\"o}nwall Lemma, as in the estimates for $\Theta^h$, we conclude that:
$$||(\mathcal{A}^+)^2\Theta||_{r+\alpha}\lesssim \lambda_q^{r}\lambda_q^{[r-(\underline{r}-2)]^+(b-1)\gamma_\ell}\ell^{-\alpha}\delta_{q+1}.$$

\noindent \textbf{Improved Estimates on $I^{cut}$.}
We can read the  second-order transport estimate for $\Theta^h$ directly from  \eqref{splitproblems} and \eqref{firsttransportapriori}. Indeed: 
\begin{equation*}
    \begin{split}
        ||\mathcal{A}^+\mathcal{A}^-\Theta^h||_{r+\alpha}&\lesssim ||\mathbb{H}_1[\Theta]||_{r+\alpha}+||\mathbb{H}_2[\mathcal{A}^+\Theta,\mathcal{A}^-\Theta]||_{r+\alpha}\\
        &\lesssim h_1^r||\Theta||_{\alpha}+h_1||\Theta||_{r+\alpha}+h_2^r(||\mathcal{A}^+\Theta||_\alpha+||\mathcal{A}^-\Theta]||_{\alpha})+h_2(||\mathcal{A}^+\Theta||_{r+\alpha}+||\mathcal{A}^-\Theta]||_{r+\alpha})\\
        &\lesssim \lambda_q^{r}\lambda_q^{[r-(\underline{r}-2)]^+(b-1)\gamma_\ell}\ell^{-\alpha}(\tau^a/\tau^c)\delta_{q+1}
    \end{split}
\end{equation*}
and given the improved estimate for $\mathcal{A}^+\mathcal{A}^-\Theta^p|_{I^{cut}}=0$ see \eqref{galbrunprincipal}, we conclude:
\begin{equation*}
    ||\mathcal{A}^+\mathcal{A}^-\Theta|_{I^{cut}}||_{r+\alpha}\lesssim \lambda_q^{r}\lambda_q^{[r-(\underline{r}-2)]^+(b-1)\gamma_\ell}\ell^{-\alpha}(\tau^a/\tau^c)\delta_{q+1}.
\end{equation*}

\noindent \textbf{Pure time derivative.} From the full equation \eqref{rewritinggr} we obtain:
\begin{equation*}
    \begin{split}
        \partial_t[\mathcal{A}^-\mathcal{A}^+\Theta]&=\partial_t\left[f\mathbb{T} F-\mathbb{H}_1[\Theta]-\mathbb{H}_2[\mathcal{A}^+\Theta,\mathcal{A}^-\Theta]\right]\\
        &=(\partial_tf)\mathbb{T} F+ f\mathbb{T} \partial_t F-\mathbb{H}_1[\partial_t\Theta]-\mathbb{H}_2[\partial_t\mathcal{A}^+\Theta,\partial_t\mathcal{A}^-\Theta]\\
        &-[\partial_t,\mathbb{H}_1][\Theta]-[\partial_t,\mathbb{H}_2][\mathcal{A}^+\Theta,\mathcal{A}^-\Theta],
    \end{split}
\end{equation*}
now, note that to get the bounds on the pure time derivatives, we can write:
\begin{equation*}
        \partial_t \Theta=\mathcal{A}^+\Theta -z^+\cn \Theta, \ \partial_t \mathcal{A}^\pm\Theta=\mathcal{A}^+\mathcal{A}^\pm \Theta -z^+\cn \mathcal{A}^\pm\Theta
\end{equation*}
and estimate each term separately, given the bounds in \eqref{firsttransportapriori}. We first deduce:
\begin{equation*}
    \begin{split}
        &||\partial_t \Theta||_{r+\alpha}\lesssim \lambda_q^{r+1}\lambda_q^{[r+1-(\underline{r}-2)]^+[b-1]\gamma_\ell}\ell^{-\alpha}\tau^a\tau^c\delta_{q+1},\\
        &||\partial_t \mathcal{A}^\pm\Theta||_{r+\alpha}\lesssim \lambda_q^{r+1}\lambda_q^{[r+1-(\underline{r}-2)]^+[b-1]\gamma_\ell}\ell^{-\alpha}\tau^a\delta_{q+1}\\
    \end{split}
\end{equation*}
and then conclude: 
$$||\partial_t[\mathcal{A}^-\mathcal{A}^+\Theta]||_{r+\alpha}\lesssim \lambda_q^{r+1}\lambda_q^{[r+1-(\underline{r}-2)]^+[b-1]\gamma_\ell}\ell^{-\alpha}\delta_{q+1}$$
where we used the bounds in Lemma \ref{operatornorms} to estimate the commutator terms. 

\noindent Finally, from \eqref{specialtransportrewriting} we obtain:
\begin{equation*}
    \begin{split}
        \partial_t \left[(\mathcal{A}^+)^2\Theta\right]&=\mathcal{A}^-(\mathcal{A}^+)^2\Theta-z^-\cn (\mathcal{A}^+)^2\\
        &= \frac{1}{\tau^a}f'\mathbb{T}F+f \mathcal{A}^+\mathbb{T}F -\mathbb{H}_1[\mathcal{A}^+\Theta]-\mathbb{H}_2[(\mathcal{A}^+)^2\Theta,\mathcal{A}^+\mathcal{A}^-\Theta]\\
        &-[\mathcal{A}^+,\mathbb{H}_1]\Theta-[\mathcal{A}^+\mathbb{H}_2][\mathcal{A}^+\Theta,\mathcal{A}^-\Theta]-z^-\cn (\mathcal{A}^+)^2
    \end{split}
\end{equation*}
and conclude that also in this case we have:
$$||\partial_t\left[(\mathcal{A}^+)^2\Theta\right]||_{r+\alpha}\lesssim \lambda_q^{r+1}\lambda_q^{[r+1-(\underline{r}-2)]^+[b-1]\gamma_\ell}\ell^{-\alpha}\delta_{q+1}.$$
\end{proof}


namely a map

\section{Chart Construction}\label{chart}
In this section, we want to construct a chart $\Psi$ from $Q_{\tau^c}(x_0)=B_{\tau^c}(x_0)\times B_{\tau^c}(t_0)\subset \mathbb{T}^3\times \mathbb{R}$, that given $(v,B)$ two divergence-free time-dependent vector fields solving on $\mathbb{T}^3\times B_{\tau^c}(x_0)$ the Faraday-Ohm system \eqref{FH}, `straightens' simultaneously $B,\ \partial_t+v\cn$. We will call such a chart $\Psi$ \textit{adapted} to $(v, B)$. 

\noindent We now clarify what we mean by `straightens' and all the properties we will require on $\Psi$. 
\begin{enumerate}
    \item \label{chart1}\textit{Chart property:} we think of $\Psi$ as a map: 
    $$\Psi:Q_{\tau^c}(x_0)\to \mathbb{R}^3$$
    such that the restriction to each time slice $$\Psi(\cdot,t):B_{\tau^c}(x_0)\to \mathbb{R}^3$$
is a diffeomorphism onto its image, close to the identity and in particular $\DD \Psi$ is invertible. 
\item \label{chart2}\textit{Magnetic field straightening}: the chart sends the magnetic field to a constant vector parallel to $e_3$, namely 
$$\Psi(\cdot,t)^{*}e_3=c B(\cdot,t)$$ 
where the pullback is meant as a map $$B_{\tau^c}(x_0)\to \mathbb{R}^3 \ \ x\mapsto \Psi(x,t)$$ for each fixed $t$ and we write $e_1, \ e_2, \ e_3$ for the coordinate vectors in $\mathbb{R}^3$ and $c_0/C_0\leq c\leq C_0/c_0$ is a constant.
\item \label{chart3}\textit{Transport straightening:} the chart sends the transport along $v$ operator to a constant vector parallel to $e_4$, namely
$$\Psi^{*}e_4=\partial_t+v\cn$$ 
where the pullback is meant as a map $$Q_{\tau^c}(x_0)\to \mathbb{R}^4 \ \ (x,t)\mapsto (\Psi(x,t),t)$$ and we write $e_4$ for the coordinate vector corresponding to $\mathbb{R}$.
\end{enumerate}

\noindent To make the construction quantitative we assume the following bounds: 
\begin{subequations}\label{nonvanishing}
    \begin{align}
        &||v||_{r},||B||_{r}\lesssim \lambda_q^rL(r)\delta_q^{1/2} \text{ for } \ r\geq 1, \\
        &||\partial_t v||_{r}\lesssim \lambda_q^{r+1}L(r+1)\delta_q^{1/2} \text{ for } \ r\geq 0, \\
         &c_0\leq |B(x,t)|\leq C_0, \ |v(x,t)|\leq C_0  \ \ \forall (x,t)\in \mathbb{T}^3\times (t_0-\tau^c, \ t_0+\tau^c)
    \end{align}
\end{subequations}
where $L:\mathbb{N}_{\geq 0}\to \mathbb{R}_{\geq 1}$ is a loss function in the sense of Definition \ref{admissible} and $c_0, \ C_0$ are constants. The pure derivative bound is needed for the applications we have in mind, since we want to estimate second-order pure time derivatives of the chart; it wouldn't be needed otherwise.

\noindent The assumed differential condition \eqref{FH} corresponds to the commutation of $\partial_t+v\cn, \ B\cn$ as space-time vector fields. The following Lemma states that this is a sufficient condition for achieving the above properties of the chart, and it is just a specific instance of the Frobenius Theorem in the context of magnetohydrodynamics. Since we need estimates on $\Psi$, we need to make the standard construction quantitative. We direct the reader to Appendix \ref{diffgeom} for the notation and conventions used. 
\begin{lemma}[Chart Construction]\label{chartconstr} Given divergence-free vector fields $(v,B)$ solving \eqref{FH} on $\mathbb{T}^3\times (t_0-\tau^c, t_0+\tau^c)$ and satisfying the estimates in \eqref{nonvanishing}, there exist a map
$$\Psi:Q_{\tau^c}(x_0)\to \mathbb{R}^3$$
for which the properties \ref{chart1}, \ref{chart2}, \ref{chart3} hold with quantitative bounds:
\begin{equation*}
    \begin{split}
        &||\DD \Psi-\IId||_0, \ ||(\DD\Psi)^{-1}-\IId||_0, \ ||\det[\DD\Psi]-1||_0\lesssim \lambda_{q+1}^{-\alpha},\\
        &||\partial_t^j\DD \Psi||_r, \ ||\partial_t^j(\DD\Psi)^{-1}||_{r}\lesssim \lambda_q^{r+j}L(r+j+1) \ \text{ for } j=0,1,2 \ \text{ and } \ r \geq 0
    \end{split}
\end{equation*}
where the implicit constants depend on $r$, the data in \eqref{nonvanishing} and the norm is taken on $Q_{\tau^c}(x_0)$.

\noindent Moreover, let $e$ be a fixed unit vector, which we think of as either a constant 1-form or 2-form. We have:
    \begin{equation*}
\begin{cases}
    (\partial_t+\mathcal{L}_{v}^1)\Psi^{1*}e=0,\\
    \mathcal{L}_B^1\Psi^{1*}e=0\\
\end{cases} \ \text{ and } \
\begin{cases}
    (\partial_t+\mathcal{L}_{v})\Psi^{2*}e=0,\\
    \mathcal{L}_B\Psi^{2*}e=0.\\
\end{cases}
\end{equation*}
In addition, for any scalar function $\varphi:\mathbb{R}\to \mathbb{R}$, $k\in e_3^\perp$ define 
$\varphi_k(x)=\varphi(x\cdot k)$, which we regard as a function $\mathbb{R}^3\to \mathbb{R}$. We have:
\begin{equation*}
    \begin{cases}
    (\partial_t+v\cn)\Psi^*\varphi_k=0,\\
    B\cn \Psi^{*}\varphi_k=0.\\
\end{cases}  
\end{equation*}
In particular, if $(k,\nu,\zeta)$ is an orthonormal basis with $k\in e_3^\perp$, then 
$$\xi=\curl[\Psi^{1*}(\varphi_k\nu)]=\varphi'\Psi^{2*}\zeta$$
is a well-defined vector field on $Q_{\tau^c}(x_0)$ satisfying:
\begin{equation*}
    (\partial_t+\mathcal{L}_{v})\xi=0, \
    \mathcal{L}_{B}\xi=0, \
    \ddiv \ \xi=0.
\end{equation*}
\end{lemma}
\begin{remark}
    The Lie derivative as a two-form and as a vector field coincide, as $v, \ B$ are divergence-free.
\end{remark}
\begin{remark}
    The size of the chart domain is optimal in the sense that:
    $$||v||_1 \tau^c, \ ||B||_1\tau^c\sim 1$$
\end{remark}

\begin{proof}[Proof of Lemma \ref{chartconstr}] The construction of $\Psi$ is done in two steps. We first freeze $B$ at $t_0$ and normalize it at $x_0$, we write: $$\bar{B}(x)=\frac{B(t_0,x)}{|B(t_0,x_0)|},$$ 
this is why $c=1/|B(t_0,x_0)|$ appears in the property \ref{chart2} above. We construct a chart $\Upsilon:\Omega \to \mathbb{R}^3$ straightening $\bar B$, in the sense that:
$$\Upsilon^*e_3=\bar B \ \text{ everywhere in } \ \Omega,$$
here $\Omega$ is a special $2C_0\tau^c$ neighborhood of $B_{\tau^c}(x_0)$ which will be constructed later, we then solve on $\mathbb{T}^3 \times(t_0-\tau^c,t_0+\tau^c)$ for the inverse and Lagrangian flow of $v$, that is
\begin{equation*}
    \begin{cases}
        \partial_t \Phi+v\cn \Phi=0,\\
        \Phi|_{t=t_0}=\IId
    \end{cases}
    \text{ and }
    \begin{cases}
        \partial_t X=v(X),\\
        X|_{t=t_0}=\IId
    \end{cases}
\end{equation*}
and set 
$$\Psi(x,t)=\Upsilon\circ \Phi_t(x)\in \mathbb{R}^3.$$
Note that the composition $\Upsilon\circ \Phi$ is well defined on $B_{\tau^c}(x_0)$. Indeed, given \ref{nonvanishing} points in $B_{\tau^c}(x_0)$ can move at most $2C_0\tau^c$ from their initial position and won't thus leave $\Omega$ on which $\Upsilon$ is defined.

\noindent This takes care of the straightening of the magnetic field. Indeed, solving \eqref{FH} in the smooth category according to \eqref{duhamel} is equivalent to $cB=X_*\bar B$, and we can compute:
\begin{equation*}
    \begin{split}
        \DD \Psi [cB]&=\DD \Upsilon(\Phi)\DD \Phi[\DD X(X^{-1})\bar B(X^{-1})]\\
        &=\DD \Upsilon(\Phi)\DD \Phi\DD (X^{-1})^{-1}\bar B(\Phi)\\
        &=\DD \Upsilon(\Phi)\DD \Phi(\DD \Phi)^{-1}\bar B(\Phi)\\
        &=\DD \Upsilon(\Phi)\bar B(\Phi)\\
        &=e_3.
    \end{split}
\end{equation*}
Here we used that $\Phi(\cdot,t)=X^{-1}(\cdot,t)$, see \eqref{lagrangianvstransport}, and $\bar B=\Upsilon^*e_3=(\DD\Upsilon)^{-1}e_3$.  This is equivalent to
$\Psi(\cdot,t)^{*}e_3=c B$.

\noindent To see the straightening of $\partial_t +v\cn$ as a space-time vector field, we first note that:
$$(\partial_t+v\cn)\Psi=\DD \Upsilon(\Phi)[(\partial_t+\cn v)\Phi]=0$$
and then, thinking of $\Psi$ as a mapping $(x,t)\mapsto (\Psi(x,t),t)$, compute: 
$$\DD_{x,t}\Psi[\partial_t+v\cn]=\begin{bmatrix}
\DD_x\Psi & \partial_t\Psi\\
0 & 1 
\end{bmatrix}\begin{bmatrix}
v \\
1 
\end{bmatrix}=\begin{bmatrix}
(\partial_t+v\cn)\Psi \\
1 
\end{bmatrix}=\begin{bmatrix}
0 \\
1 
\end{bmatrix}=e_4$$
and we conclude that $\Psi^{*}e_4=(\partial_t+v\cn)$ under the proper identifications. We are then left with the construction of $\Upsilon$ and $\Omega$. 

\noindent \textbf{Construction of $\Omega$.} We want $\Omega$ to be a flow-box  for $\bar{B}$ centred at $x_0$ and containing a $2C_0\tau^c$ neighborhood of $B_{\tau^c}(x_0)$. We call the plane passing through $x_0$ and orthogonal to $\bar B(x_0)$, $\pi_0$. Eventually applying a fixed rotation, we may assume that $\bar{B}(x_0)=B(x_0,t_0)/|B(x_0,t_0)|=e_3$. We then let $Y$ solve:
\begin{equation*}
    \begin{cases}
        \partial_sY=\bar B\left(Y\right),\\
        Y|_{s=0}=\IId
    \end{cases}
\end{equation*}
fix $r_0,\ s_0\geq 0$, and set 
$$\Omega=\{Y_s(B_{r_0\tau^c}(x_0)\cap\pi_0):|s|\leq s_0\tau^c\}.$$
We claim that for $r_0, \ s_0$ sufficiently large we have:
\begin{equation}\label{chain}
        B_{\tau^c}(x_0)\subset B_{2C_0\tau^c}(B_{\tau^c}(x_0))\subset B_{(1+2C_0)\tau^c}(x_0)\subset \Omega
\end{equation}
where the first two inequalities are always satisfied independently of $r_0, \ s_0$. To see this, we expand $Y_s$ as follows:
\begin{equation*}
    \begin{split}
        Y_s(x)&=x+s\bar B(x)+\int_0^s(\bar B\cn \bar B)(Y_{s'})(s-s')\dd s'\\
        &= x+s\bar B(x_0)+\underbrace{s\int_0^1\DD \bar B(sx+(1-s)x_0)[x-x_0]+\int_0^s(\bar B\cn \bar B)(Y_{s'})(s-s')\dd s'}_{r(x,s)}\\
        &=x+se_3+r(x,s).
    \end{split}
\end{equation*}
Note that for $|s|\leq s_0\tau^c$ and $x\in B_{r_0\tau^c}(x_0)\cap\pi_0\subset B_{r_0\tau^c}(x_0)$, given our assumptions \eqref{nonvanishing} and the definition of $\tau^c$ in \eqref{tauc}, we have:
$$|r(x,s)|\leq  \tau^c\frac{C(s_0^2+s_0r_0)}{\lambda_{q+1}^\alpha}$$
where the constant $C$ is independent of $r_0,s_0$. Now fix $r_0=s_0=(1+2C_0)+2$ for any $\alpha, \ b$ we can pick $a$ so large that $$\frac{C(s_0^2+s_0r_0)}{\lambda_{q+1}^\alpha}< 1\Longrightarrow |r(x,s)|<  \tau^c$$
we deduce
$$B_{(1+2C_0)\tau^c+1}(x_0)\subset \Omega$$
and conclude that \eqref{chain} holds for this choice of $\Omega$.

\noindent \textbf{Construction of $\Upsilon$.} We now set:
\begin{equation}\label{U}
    U=\{x_0+(x_1,x_2,s)\in \mathbb{R}^3:\sqrt{|x_1|^2+|x_2|^2}\leq r_0\tau^c, |s|\leq s_0\tau^c\}
\end{equation}
and we put new coordinates on $\Omega$ by means of the map:
$$\Gamma:\mathbb{R}^3\supset U\ni(x_1,x_2,s)\longrightarrow Y_s(x_0+(x_1,x_2,0))\in\Omega \subset \mathbb{T}^3.$$
The key point here is that its differential can be made uniformly close to the identity and has the right geometric properties. First, write:
\begin{equation*}
    \begin{split}
        \DD \Gamma&=\left[\partial_1\Gamma,\partial_2\Gamma,\partial_s \Gamma\right]=\left[\partial_1Y,\partial_2Y,\partial_s Y\right]=\left[\partial_1Y,\partial_2Y,\bar{B}(Y)\right]=\left[\partial_1Y,\partial_2Y,\bar{B}(\Gamma)\right]
    \end{split}
\end{equation*}
and assume for a moment we showed $\DD \Gamma$ is invertible everywhere in $U$ so that by the Inverse Function Theorem \ref{ift}, $\Gamma$ also is, then for $\Upsilon=\Gamma^{-1}$ we have:
\begin{equation*}
    \begin{split}
        \Upsilon^{*}e_3=(\DD \Upsilon)^{-1}[e_3]=\DD (\Upsilon^{-1})(\Upsilon)[e_3]=\DD \Gamma(\Gamma^{-1})[e_3]=\left[\partial_1Y(\Gamma^{-1}), \ \partial_2Y(\Gamma^{-1}), \ \bar{B}\right][e_3]=\bar B
    \end{split}
\end{equation*}
and the straightening property \ref{chart2} is proven. 

\noindent In what follows, we show invertibility and higher-order derivative estimates for $\Gamma, \Upsilon$, we define
$$\DD \Gamma-\IId=[r_1,r_2,r_3]$$
and bound each remainder $r_i$. We begin with $r_3$:  
\begin{equation*}
    \begin{split}
        \bar{B}(Y_s(x_0+(x_1,x_2,0))&=\bar{B}(x_0+(x_1,x_2,0))+\int_0^s \partial_{s'}(\bar B(Y_{s'}))\dd s'\\
        &=\bar{B}(x_0+(x_1,x_2,0))+\int_0^s (\bar{B}\cn \bar{B})(Y_{s'})\dd s'\\
        &=\bar{B}(x_0)+\int_0^1\DD \bar{B}(s'(x_1,x_2,0)+x_0)[(x_1,x_2,0)^\top]\dd s'+\int_0^s (\bar{B}\cn \bar{B})(Y_{s'})\dd s'\\
        &=e_3+\int_0^1\DD \bar{B}(s'(x_1,x_2,0)+x_0)[(x_1,x_2,0)^\top]\dd s'+\int_0^s (\bar{B}\cn \bar{B})(Y_{s'})\dd s'\\
        &=e_3+r_3
    \end{split}
\end{equation*}
and it follows that for $(x_1,x_2,s)\in U$, we have:
$$|r_3(x_1,x_2,s)|\lesssim \tau_c\lambda_qL(1)\delta_q^{1/2}=\lambda_{q+1}^{-\alpha}$$
where we used $L(1)=1$. The estimates on $r_1,r_2$ follow immediately from those for the Lagrangian flow $Y$ in terms of $\bar B$, see Proposition \ref{standardlagrangianestimate} and we deduce:
\begin{equation}\label{Dgammazero}
    ||\DD \Gamma-\IId||_0\lesssim \lambda_{q+1}^{-\alpha}
\end{equation}
where the sup norm is though on $U$. The higher-order derivatives also follow from Proposition \ref{standardlagrangianestimate} and the assumption \eqref{nonvanishing} upon noticing that for any $r$ the derivatives in the coordinate $s$ correspond to:
$$\partial_s^{r}[Y_s(x_0+(x_1,x_2,0))]=\bar B\cn^{r-1}\bar B(Y_s(x_0+(x_1,x_2,0))),$$
and we conclude:
\begin{equation}\label{gamma}
||\DD \Gamma||_{r} \lesssim \lambda_q^rL(r+1)  \ \text{ for } \ r\geq 1
\end{equation}
where the derivatives are taken with respect to the coordinates $(x_1,x_2,s)$ on $U$, we won't mention this point again in the following. The implicit constant depends only on the data \eqref{nonvanishing} and $r$.  

\noindent In particular, given any $\alpha, b$ we can choose $a$ sufficiently large to ensure that: 
$$|\DD \Gamma|_{op}\geq 1-|\DD \Gamma-\IId|_{op}\geq 1-C\lambda_{q+1}^{-\alpha}>1/2 \ \text{ everywhere in } \ U$$
and an application of the Inverse Function Theorem \ref{ift} then gives the existence of a smooth inverse $\Upsilon=\Gamma^{-1}: \Omega \to U$ and the validity of the formula: 
$$\DD (\Gamma^{-1})=(\DD \Gamma)^{-1}(\Gamma^{-1}).$$
We now establish estimates for this inverse map. Recall that for any matrix $A$ with $|A|_{op}<1$ we have
$$|(\IId-A)^{-1}-\IId|_{op}\lesssim |A|_{op}$$
and thus from \eqref{Dgammazero} we deduce:
$$||\DD \Upsilon-\IId||_0=||(\DD \Gamma)^{-1}(\Gamma^{-1})-\IId||_0\leq||(\DD \Gamma)^{-1}-\IId||_0\lesssim\lambda_{q+1}^{-\alpha}$$
where the $\sup$ norm is though in $\Omega=\Gamma(U)$.

\noindent We now move to higher-order derivatives $r\geq 1$. We first claim that:
\begin{equation}\label{gamma-1}
    ||(\DD\Gamma)^{-1}||_{r}\lesssim \lambda_q^r L(r+1)
\end{equation}
where the implicit constant depends on $r$ and the data \eqref{nonvanishing}. Set $A=\DD\Gamma$, for any index $i=1,2,3$ and multi index $\sigma$, with $|\sigma|=r-1$ we have:
\begin{equation*}
    \partial_\sigma\partial_iA^{-1}=-\partial_\sigma (A^{-1}\partial_iAA^{-1}),\\
\end{equation*}
we can then estimate:
\begin{equation*}
    \begin{split}
        ||\partial_\sigma\partial_i A^{-1}||_0&= ||\partial_\sigma (A^{-1}\partial_iAA^{-1})||_0\leq ||A^{-1}\partial_iAA^{-1}||_{r-1}\\
        &\lesssim ||A^{-1}||_{r-1}||\partial_iA||_0||A^{-1}||_{0}+||A^{-1}||_{0}^2||\partial_iA||_{r-1}\\
    \end{split}
\end{equation*}
and the claim follows by induction together with $L(2)=1$ and the fact that $L(r)$ is non-decreasing.

\noindent This, together with composition estimates in Proposition  \ref{compestimates}, gives us: 
\begin{equation*}
    \begin{split}
        [\DD\Upsilon]_r&=[(\DD\Gamma)^{-1}(\Upsilon)]_r\lesssim [(\DD\Gamma)^{-1}]_1||\DD \Upsilon||_{r-1}+||(\DD\Gamma)^{-1}||_{r}||\DD \Upsilon||^r_0\\
    \end{split}
\end{equation*}
and by induction together with the estimates on $\DD \Gamma$ in \eqref{gamma} we conclude:
\begin{equation*}
    \begin{split}
        ||\DD\Upsilon-\IId||_0\lesssim \lambda_{q+1}^{-\alpha}, \ \ \ \  ||\DD\Upsilon||_r\lesssim \lambda_q^rL(r+1)
    \end{split}
\end{equation*}
where the implicit constants depend only on $r$ and the data in \eqref{nonvanishing}, and
the $\sup$ norm is thought on $\Omega$. 

\noindent \textbf{Conclusion.} We now restrict the space domain of $\Phi$ to $B_{\tau^c}(x_0)$ so that the composition $\Psi=\Upsilon\circ \Phi$ is well defined for $|t-t_0|< \tau^c$ and proceed with proving derivative bounds. We first rewrite:
\begin{equation*}
    \begin{split}
        \DD \Psi-\IId&= \DD \Upsilon(\Phi)\DD\Phi-\IId\\
        &=\DD \Upsilon(\Phi)[\DD\Phi-\IId]+\DD \Upsilon(\Phi)-\IId\\
        &=[\DD \Upsilon(\Phi)-\IId][\DD\Phi-\IId]+[\DD\Phi-\IId]+[\DD \Upsilon(\Phi)-\IId]\\\\
        (\DD \Psi)^{-1}-\IId&= (\DD\Phi)^{-1}(\DD \Upsilon)^{-1}(\Phi)-\IId\\
        &= \DD X(\Phi)\DD \Gamma(\Upsilon\circ\Phi)-\IId\\
        &=[\DD X(\Phi)-\IId][\DD \Gamma(\Upsilon\circ\Phi)-\IId]+[\DD X(\Phi)-\IId]+[\DD \Gamma(\Upsilon\circ\Phi)-\IId]\\\\
    \end{split}
\end{equation*}
here we used that $\Phi^{-1}=X$, see \eqref{lagrangianvstransport}, and conclude from the estimates on $\Phi, X$ coming from Propositions \ref{standardtransportestimate}, \ref{standardlagrangianestimate}, and the definition of $\tau^c$ in \eqref{tauc}, that:
\begin{equation*}
    \begin{split}
        ||\DD \Psi-\IId||_0
        &\leq||\DD \Upsilon-\IId||_0||\DD\Phi-\IId]||_0+||\DD\Phi-\IId||_0+||\DD \Upsilon-\IId||_0,\\
        &\lesssim \lambda_{q+1}^{-\alpha}\\
        ||(\DD \Psi)^{-1}-\IId||_0&\leq ||\DD X-\IId||_0||\DD \Gamma-\IId||_0+||\DD X-\IId||_0+||\DD \Gamma-\IId||_0\\
        &\lesssim \lambda_{q+1}^{-\alpha}\\
    \end{split}
\end{equation*}
and from
$|\det[\IId+A]-1|\lesssim |A|_{op}$
we also deduce:
$$||\det [\DD \Psi] -1||_0, \ ||\det [(\DD \Psi)^{-1}] -1||_0 \lesssim \lambda_{q+1}^{-\alpha}.$$
For $r\geq 1$ by means of the composition estimates in Proposition \ref{compestimates} we conclude:
\begin{equation*}
    \begin{split}
        ||\DD \Psi||_r&\lesssim ||\DD \Upsilon(\Phi)||_r||\DD\Phi||_0+ ||\DD \Upsilon(\Phi)||_0||\DD\Phi||_r\\
        &\lesssim \left(||\DD \Upsilon||_1||\DD \Phi||_{r-1}+||\DD \Upsilon||_r||\DD \Phi||_{0}^r\right)||\DD\Phi||_0+||\DD \Upsilon||_0||\DD\Phi||_r\\
        &\lesssim \lambda_q^rL(r+1)
    \end{split}
\end{equation*}
where we used the fact that $L(2)=1$ and the implicit constant depends on the data and $r$. The bounds for $(\DD \Psi)^{-1}$ follow similarly.

\noindent \textbf{Lie-Transport properties.} It follows from \eqref{duhamel} that
$$(\partial_t+\mathcal{L}_{v}^1)\Psi^{1*}e=(\partial_t+\mathcal{L}_{v}^1)\Phi^{1*}(\Upsilon^{1*}e)=0,$$
similarly,
$$(\partial_t+\mathcal{L}_{v})\Psi^{2*}e=(\partial_t+\mathcal{L}_{v})\Phi^{*}(\Upsilon^{2*}e)=0,$$
the same computations hold for the scalar case and we deduce $(\partial_t+v\cn)\Psi^*\varphi_k=0$, here the direction $k$ plays no role. 

\noindent From the properties in \eqref{DG3} we deduce:
\begin{equation*}
        \mathcal{L}_{B}^1\Psi^{1*}e=1/c\mathcal{L}_{\Psi^{*}e_3}\Psi^{1*}e=1/c\Psi^{1*}(\mathcal{L}_{e_3}e)=0,
\end{equation*}
similarly, from \eqref{DG0} and \eqref{DG4} we obtain:
\begin{equation*}
    \begin{split}
        \mathcal{L}_{B}\Psi^{2*}e&=[B,\Psi^{2*}e]=\curl[B\times\Psi^{2*}e]\\
        &=1/c\curl[\Psi^{*}e_3\times\Psi^{2*}e]\\
        &=1/c\curl[\Psi^{1*}(e_3\times e)]\\
        &=\Psi^{2*}\curl[e_3\times e]=0\\
    \end{split}
\end{equation*}
as both $e_3$ and $e$ are constant. For the scalar case, we use that $k\cdot e_3=0$ and the property in \eqref{DG5} to compute:
\begin{equation*}
    \begin{split}
        B\cn \Psi^{*}\varphi_k=B\cdot \dd\Psi^{*}\varphi_k=B\cdot \Psi^{1*}\dd\varphi_k=1/c\Psi^*e_3\cdot \Psi^{1*}\nabla\varphi_k=1/c\Psi^{*}(e_3\cn \varphi_k)=\varphi'(k\cdot \Psi)e_3\cdot k=0
    \end{split}
\end{equation*}
as $d=\nabla$ in our setting, see \eqref{identifications}. Combining the above with $\ddiv\curl=0$ and \eqref{DG1} to commute the Lie operators with $\curl$, we deduce the properties for $\xi$.

\noindent \textbf{Pure Time derivatives.} From the Lie transport properties for constant 1-forms above, we deduce: 
\begin{equation*}
    \begin{split}
        \partial_t (\DD \Psi)^\top=&-v\cn\DD \Psi^\top-(\DD  v)^\top\DD \Psi^\top,\\
        \partial_t^2 \DD \Psi^\top=&-(\partial_t v)\cn\DD \Psi^\top- v\cn\partial_t\DD \Psi_I^\top\\
        &-(\DD \partial_tv)^\top\DD \Psi^\top-(\DD  v)^\top\partial_t\DD \Psi_I^\top.\\
    \end{split}
\end{equation*}
Similar identities hold for $\DD \Psi^{-1}$ where $-\DD v^\top$ is substituted by $-\DD v$. From the assumed pure-time derivative bound on $v$ in \eqref{nonvanishing}, we conclude:
$$||\partial_t^j\DD \Psi||_r, \ ||\partial_t^j\DD\Psi^{-1}||_{r}\lesssim \lambda_q^{r+j}L(r+j+1) \ \text{ for } j=0,1,2 \ \text{ and } \ r \geq 0.$$
\end{proof}


\begin{remark}[An Interesting Cancellation]\label{interestin}
    We claim that:
    $$B\cn(\det[\DD \Psi])=0$$
    from which it quickly follows that:
    \begin{equation}\label{transportdeterminant}
        \mathcal{A}^\pm\det[\DD \Psi]=0.
    \end{equation}
    Indeed, we have:   
    \begin{equation*}
        \begin{split}
            (\partial_t+v\cn)(\det[\DD \Psi])&=(\partial_t+v\cn)(\det[\DD \Upsilon(\Phi)\DD \Phi])\\
            &=(\partial_t+v\cn)(\det[\DD \Upsilon](\Phi)\det[\DD \Phi])\\
            &=(\partial_t+v\cn)(\det[\DD \Upsilon](\Phi))\\
            &=\DD (\det[\DD \Upsilon])(\Phi)[(\partial_t+v\cn)\Phi]=0
        \end{split}
    \end{equation*}
    where we used that $\ddiv \ v=0 \leadsto \det[\DD \Phi]=1$.
    
\noindent To prove the claim, we first compute: 
    \begin{equation*}
        \begin{split}
            c B\cn\det[\DD \Psi]&=cB\cn(\det[\DD \Upsilon](\Phi))\\
            &=\Phi^*\bar B\cn (\Phi^*(\det[\DD \Upsilon]))\\
            &=\Phi^*\bar B\cdot \Phi^{1*}\nabla(\det[\DD \Upsilon])\\
            &=(\bar{B}\cdot\nabla\det[\DD \Upsilon])(\Phi)
        \end{split}
    \end{equation*}
    where we used that $cB=X_*\bar B=\Phi^*\bar B$ and the identity in \eqref{DG5}. In particular, to prove the claim, it suffices to show: 
    $$\bar{B}\cdot\nabla\det[\DD \Upsilon]=0.$$
Now set $f=\det[\DD \Upsilon]$, since the $\bar{B}$ derivative on the manifold corresponds to the $\partial_3=\partial_s$ derivative in the chart, we have:
\begin{equation}\label{manifoldderivative}
    \begin{split}
        \bar B\cdot \nabla f&=(\partial_s(f\circ \Gamma))(\Upsilon)\\
        &=[\partial_s(\det[\DD \Upsilon(\Gamma)]](\Upsilon)\\
        &=[\partial_s(\det[(\DD \Gamma)^{-1}])](\Upsilon)\\
        &=\left[\partial_s\frac{1}{\det[\DD \Gamma]}\right](\Upsilon)\\
        &=\left[-\frac{1}{2\det[\DD \Gamma]}\text{trace}\left[(\DD\Gamma)^{-1}\partial_s\DD \Gamma\right]\right](\Upsilon)\\
    \end{split}
\end{equation}
and with this at hand, we need only compute: 
\begin{equation*}
\partial_s\DD\Gamma=\DD\partial_s\Gamma=\DD\left[\bar B(\Gamma)\right]=\DD \bar B(\Gamma)\DD\Gamma
\end{equation*}
we can plug this back in \eqref{manifoldderivative} and by means of the identity $\text{trace}[AB]=\text{trace}[BA]$ for any two matrices $A, B$, deduce:
\begin{equation*}
    \begin{split}
        \bar{B}\cdot \nabla (\det[\DD \Upsilon])&=\left[-\frac{1}{2\det[\DD\Gamma]}\text{trace}\left[(\DD\Gamma)^{-1}\partial_s\DD\Gamma\right]\right](\Upsilon)\\
        &=\left[-\frac{1}{2\det[\DD\Gamma]}\text{trace}\left[(\DD\Gamma)^{-1}\DD \bar B(\Gamma)\DD\Gamma\right]\right](\Upsilon)\\
        &=-\frac{\ddiv \bar{B}}{2\det[\DD \Gamma](\Upsilon)}\\
        &=0
    \end{split}
\end{equation*}
and the claim follows.
\end{remark}


\subsection{Stability Estimates}
We now perturb our original vector fields: $$(v,B)\leadsto (\tilde v,\tilde B)=(v+w,B+b),$$
we still assume that $(\tilde v,\tilde B)$ are divergence-free, solve \eqref{FH} and satisfy \eqref{nonvanishing}. In particular, we can construct a new chart $\tilde \Psi: Q_{\tau^c}(x_0)\to \mathbb{R}^3$ associated with this new pair by means of Lemma \ref{chartconstr}. We want to estimate the difference between the two constructions. We assume that:
\begin{equation}\label{stabilitydata}
    ||w||_r,||b||_r\lesssim \lambda_q^r\tilde L(r)\lambda_q\ell^{-\alpha}\tau^a\delta_{q+1}
\end{equation}
and we also require: 
$$L\leq \tilde L, \ \tilde L(0)=\tilde L(1)=\tilde L(2)=1, \ \tilde L(r)\tilde L(r')\leq \tilde L(r+r') \ \text{ for } \ r, r' \geq 0.$$

\noindent We will show the following stability result for the chart construction.
\begin{lemma}[Chart Stability]\label{chartstability} Let $r\geq0$ be an integer. Under the above assumptions, we have:
\begin{equation*}
    \begin{split}
        &||\tilde \Psi-\Psi||_r\lesssim \lambda_q^r\tilde L(r+1)\frac{1}{\lambda_q}\mathcal{T}_g,\\
        &||(\DD\tilde \Psi)^{-1}-(\DD\Psi)^{-1}||_r\lesssim \lambda_q^r\tilde L(r+2)\mathcal{T}_g\\
    \end{split}
\end{equation*}
where the implicit constants depend on $r$ and the data \eqref{nonvanishing}, \eqref{stabilitydata}.
\end{lemma}
\begin{remark} The estimate \eqref{stabilitydata} on the perturbation is chosen to match that of $(w^g,b^g)$ from the Galbrun stage, but in principle, the result can be adapted to any given one.
\end{remark}
\begin{proof}[Proof of Lemma \ref{chartstability}] We split the proof in 2 steps. Since $\Psi=\Upsilon\circ\Phi$ and $\Upsilon$ is constructed as the inverse of the auxiliary map $\Gamma$, in the first part, we prove stability estimates for $\Gamma$ and $\Phi$. In the second part, we deal with the stability of $\Upsilon$ and the composition. 

\noindent \textbf{First step.} Let $\bar B, \ \bar{\tilde B}, \ Y,\ \tilde Y, \ \Gamma, \ \tilde \Gamma, \ \Upsilon, \ \tilde \Upsilon$ be the intermediate objects in the construction behind the proof of Lemma \ref{chartconstr}. We begin with the following observations.  By direct computation, we see that the difference $\tilde \Phi-\Phi$ solves:
\begin{equation*}
    \begin{cases}
        \partial_t(\tilde\Phi-\Phi)+v\cn(\tilde\Phi-\Phi)=-w\cn \tilde \Phi,\\
        (\tilde\Phi-\Phi)|_{t=t_0}=0
    \end{cases} 
\end{equation*}
given any $b, \ \alpha$, we can always choose $a$ sufficiently large to have $||v||_{1}2\tau^c\leq 1$, then apply Propositions \ref{standardtransportestimate} to deduce for $|t-t_0|\leq \tau^c$:
\begin{equation}\label{stabilityphi}
    \begin{split}
        ||\tilde\Phi_t-\Phi_t||_r&\lesssim \int_{t_0}^t||(w\cn \tilde \Phi)(\cdot,s')||_r+(t-s')[v]_{r}||(w\cn \tilde \Phi)(\cdot,l)||_1\dd s'\\
        &\lesssim  \lambda_q^r\tilde L(r+1)\frac{1}{\lambda_q}\mathcal{T}_g.
    \end{split}
\end{equation}
Similarly, the difference: $$F(x_1,x_2,s)=(\tilde \Gamma-\Gamma)(x_1,x_2,s)=(\tilde Y_s-Y_s)(x_0+(x_1,x_2,0)),$$ 
solves the problem:
\begin{equation}\label{Ystability}
    \begin{cases}
        \partial_sF+A F=G,\\
        F|_{s=0}=0
    \end{cases}
\end{equation}
where 
\begin{equation*}
    \begin{split}
        A(x_1,x_2,s)&=-\int_0^1(\DD \bar{\tilde B})(s'\tilde Y_s(x_0+(x_1,x_2,0))+(1-s')Y_s(x_0+(x_1,x_2,0)))\dd s'\\
                    &=-\int_0^1(\DD \bar{\tilde B})(s'\tilde\Gamma(x_1,x_2,s)+(1-s')\Gamma(x_1,x_2,s))\dd s',\\
        G(x_1,x_2,s)&=(\bar{ \tilde B}-\bar B)(Y_s(x_0+(x_1,x_2,0)))\\
        &=(\bar{ \tilde B}-\bar B)(\Gamma(x_1,x_2,s))
    \end{split}
\end{equation*}
and the solution can be explicitly written as
\begin{equation}\label{stabilityF}
    F(\cdot,s)=e^{-\int_0^sA(\cdot,s')\dd s'}\int_0^se^{\int_0^{s'}A(s'',\cdot)\dd s''}G(\cdot,s')\dd s'.
\end{equation}
We now provide bounds on $A$ and $G$, which also guarantee that the representation \eqref{stabilityF} is well-defined. We first rewrite:
\begin{equation*}
    \begin{split}
        G&=(\bar{ \tilde B}-\bar B)(\Gamma)=\frac{(B+b)(\Gamma)}{|B+b|(x_0,t_0)}-\frac{B(\Gamma)}{|B|(x_0,t_0)}\\
        &=\frac{1}{|B+b|(x_0,t_0)}b(\Gamma)+\left(\frac{|B|-|B+b|}{|B||B+b|}\right)(x_0,t_0)B(\Gamma)
    \end{split}
\end{equation*}
note that given \eqref{nonvanishing} and \eqref{stabilitydata} we have:
$$\frac{1}{|B+b|(x_0,t_0)}\lesssim \frac{1}{c_0-|b|}\lesssim1 \ \text{ and } \ \left(\frac{\left||B|-|B+b|\right|}{|B||B+b|}\right)(x_0,t_0)\lesssim \frac{|b|}{c_0(c_0-|b|)}\lesssim\lambda_q\ell^{-\alpha}\tau^a\delta_{q+1}$$
where the implicit constants depend on the data, and to bound the numerators, we choose $a$ sufficiently large so that $c_0/2>\lambda_q\ell^{-\alpha}\tau^a\delta_{q+1}$. By means of the composition estimates in Proposition \ref{compestimates}, the bounds for $\Gamma$ in \eqref{gamma}, and the fact that $|s|\lesssim \tau^c$ we deduce:
\begin{equation*}
        ||G||_r\lesssim \lambda_q^r\tilde L(r)\lambda_q\ell^{-\alpha}\tau^a\delta_{q+1} \ \text{ and } \ ||A||_r\lesssim \lambda_q^r\tilde L(r+1)\lambda_q\delta_q^{1/2},\\
\end{equation*}
in particular, given $b$ we can always choose $a$ sufficiently large so that $||\int_0^sA\dd s'||_0< 1$. From the representation formula \eqref{stabilityF}, and the bounds for $A, \ G$ we conclude:
\begin{equation}\label{stabilitygamma}
    ||\tilde \Gamma-\Gamma||_r=||F||_r\lesssim \lambda_q^r\tilde L(r+1)\lambda_q\ell^{-\alpha}\tau^c\tau^a\delta_{q+1}=\lambda_q^r\tilde L(r+1)\frac{1}{\lambda_q}\mathcal{T}_g.
\end{equation}

\noindent \textbf{Second step.} Fix two diffeomorphisms $\varphi,\ \psi$ with same domain and image. We use the following standard trick
\begin{equation*}
        \varphi^{-1}-\psi^{-1}=\varphi^{-1}\circ \psi \circ \psi^{-1}-\varphi^{-1}\circ \varphi\circ\psi^{-1}
\end{equation*}
from this deduce
\begin{equation*}
    \begin{split}
        ||\varphi^{-1}-\psi^{-1}||_0&\leq [\varphi^{-1}]_1||\psi \circ \psi^{-1}-\varphi\circ\psi^{-1}||_0\\
        &\leq [\varphi^{-1}]_1||\psi -\varphi||_0\\
    \end{split}
\end{equation*}
and applying this with $\varphi= \tilde \Gamma, \psi=\Gamma$ from the construction of $\tilde \Psi, \Psi$ with the related estimates given by Lemma \ref{chartconstr}, and \eqref{stabilitygamma} above we conclude:
\begin{equation}\label{stability0}
    ||\tilde \Upsilon -\Upsilon||_0\lesssim ||\tilde \Gamma -\Gamma||_0\lesssim \frac{1}{\lambda_q}\mathcal{T}_g
\end{equation}
where the implicit constants depend on the data \eqref{nonvanishing} and \eqref{stabilitydata}. We now deal with higher-order derivatives. We need some preliminary bounds. We will repeatedly use the composition estimates in Proposition \ref{compestimates}. We first compute:
\begin{equation*}
    \begin{split}
        (\DD\varphi)^{-1}(\varphi^{-1})-(\DD\varphi)^{-1}(\psi^{-1})=\int_0^1 \DD [(\DD\varphi)^{-1}](s\varphi^{-1}+(1-s)\psi^{-1})[\varphi^{-1}-\psi^{-1}]\dd s,
    \end{split}
\end{equation*}
from which we deduce:
\begin{equation}\label{stability1}
    \begin{split}
       &||(\DD\varphi)^{-1}(\varphi^{-1})-(\DD\varphi)^{-1}(\psi^{-1})||_r\\
       &\lesssim \int_0^1 ||\DD [(\DD\varphi)^{-1}](s\varphi^{-1}+(1-s)\psi^{-1})||_0||\varphi^{-1}-\psi^{-1}||_r\dd s\\
       &+\int_0^1||\DD [(\DD\varphi)^{-1}](s\varphi^{-1}+(1-s)\psi^{-1})||_r||\varphi^{-1}-\psi^{-1}||_0\dd s\\
       &\overset{r=0}{\lesssim}||(\DD\varphi)^{-1}||_1||\varphi^{-1}-\psi^{-1}||_0\\
       &\overset{r\geq 1}{\lesssim} ||(\DD\varphi)^{-1}||_1||\varphi^{-1}-\psi^{-1}||_r\\
       &+\left[||(\DD\varphi)^{-1}||_2(||\varphi^{-1}||_r+||\psi^{-1}||_r)+||(\DD\varphi)^{-1}||_{r+1}(||\varphi^{-1}||_1+||\psi^{-1}||_1)^r\right]||\varphi^{-1}-\psi^{-1}||_0.
    \end{split}
\end{equation}
Secondly, we compute:
\begin{equation}\label{stability2}
    \begin{split}
        ||(\DD\varphi)^{-1}(\DD\varphi-\DD\psi)(\DD\psi)^{-1}||_r&\lesssim\left(||(\DD\varphi)^{-1}||_r||(\DD\psi)^{-1}||_0+||(\DD\varphi)^{-1}||_0||(\DD\psi)^{-1}||_r\right)[\varphi-\psi]_1\\
        &+||(\DD\varphi)^{-1}||_0||(\DD\psi)^{-1}||_0||\varphi-\psi||_{r+1}.
    \end{split}
\end{equation}
We now go back to the bound we are interested in. We will use that for any two matrices $A, B$ one has $A^{-1}-B^{-1}=-A^{-1}(A-B)B^{-1}$ to compute:
\begin{equation}\label{stability3}
    \begin{split}
        [\DD(\varphi^{-1}-\psi^{-1})]_r&=[(\DD\varphi)^{-1}(\varphi^{-1})-(\DD\psi)^{-1}(\psi^{-1})]_r\\
        &=[(\DD\varphi)^{-1}(\varphi^{-1})-(\DD\varphi)^{-1}(\psi^{-1})]_r+[((\DD\varphi)^{-1}-(\DD\psi)^{-1})(\psi^{-1})]_r\\
        &\overset{r=0}{\leq}||(\DD\varphi)^{-1}(\varphi^{-1})-(\DD\varphi)^{-1}(\psi^{-1})||_0+||(\DD\varphi)^{-1}||_0||\DD\varphi-\DD\psi||_0||(\DD\psi)^{-1}||_0\\
        &\overset{r\geq 1}{\lesssim} [(\DD\varphi)^{-1}(\varphi^{-1})-(\DD\varphi)^{-1}(\psi^{-1})]_r\\
        &+||(\DD\varphi)^{-1}-(\DD\psi)^{-1}||_1||\psi^{-1}||_r+||(\DD\varphi)^{-1}-(\DD\psi)^{-1}||_r[\psi^{-1}]_1^r\\
        &\lesssim [(\DD\varphi)^{-1}(\varphi^{-1})-(\DD\varphi)^{-1}(\psi^{-1})]_r\\
        &+||(\DD\varphi)^{-1}(\DD\varphi-\DD\psi)(\DD\psi)^{-1}||_1||\psi^{-1}||_r+||(\DD\varphi)^{-1}(\DD\varphi-\DD\psi)(\DD\psi)^{-1}||_r[\psi^{-1}]_1^r.
    \end{split}
\end{equation}
We combine \eqref{stability1}, \eqref{stability2}, \eqref{stability3} for $r=0$, with $\varphi= \tilde \Gamma, \psi=\Gamma$ and the bounds in \eqref{gamma-1}, \eqref{stability0} to deduce:
$$||\tilde \Upsilon-\Upsilon||_1\lesssim \tilde L(2)\mathcal{T}_g=\mathcal{T}_g$$
and for $r\geq 1$, we proceed by induction and obtain:
\begin{equation}\label{stabilityupsilon}
        ||\tilde \Upsilon-\Upsilon||_r\lesssim \lambda_q^r\tilde L(r+1)\frac{1}{\lambda_q}\mathcal{T}_g
\end{equation}
where the implicit constants depend on $r$ and the data in \eqref{nonvanishing}, \eqref{stabilitydata}.

\noindent \textbf{Conclusion.} We write
\begin{equation*}
    \begin{split}
        \tilde \Psi-\Psi&=\tilde \Upsilon\circ\tilde\Phi-\tilde \Upsilon\circ\Phi+\tilde \Upsilon\circ \Phi-\Upsilon\circ\Phi\\
        &=\int_0^1\DD\tilde\Upsilon(s\tilde\Phi+(1-s)\Phi)[\tilde \Phi-\Phi]\dd s+(\tilde \Upsilon-\Upsilon)(\Phi).\\
    \end{split}
\end{equation*}
and use the composition estimates in Proposition \ref{compestimates} together with \eqref{stabilityphi}, \eqref{stabilityupsilon} to bound:
\begin{equation*}
    \begin{split}
        ||\tilde \Psi-\Psi||_r&\lesssim ||\tilde \Phi-\Phi||_0\int_0^1||\DD\tilde\Upsilon(s\tilde\Phi+(1-s)\Phi)||_r\dd s+||\tilde \Phi-\Phi||_r\int_0^1||\DD\tilde\Upsilon(s\tilde\Phi+(1-s)\Phi)||_0\dd s\\
        &+||(\tilde \Upsilon-\Upsilon)(\Phi)||_r\\
        &\overset{r=0}{\lesssim}\frac{1}{\lambda_q}\mathcal{T}_g\\
        &\overset{r\geq 1}{\lesssim} \frac{1}{\lambda_q}\mathcal{T}_g\left[||\DD\tilde\Upsilon||_1(||\DD \tilde\Phi||_{r-1}+||\DD\Phi||_{r-1})+||\DD\tilde\Upsilon||_r(||\DD\tilde\Phi||_0+||\DD\Phi||_0)^r\right]\\
        &+||\tilde \Phi-\Phi||_r+||\tilde \Upsilon-\Upsilon||_1||\DD \Phi||_{r-1}+||\tilde \Upsilon-\Upsilon||_r||\DD\Phi||_0^r\\
        &\lesssim \lambda_q^r\tilde L(r+1)\frac{1}{\lambda_q}\mathcal{T}_g
    \end{split}
\end{equation*}
and conclude:
\begin{equation*}
    \begin{split}
        ||(\DD\tilde \Psi)^{-1}-(\DD\Psi)^{-1}||_r&= ||(\DD\tilde \Psi)^{-1}(\DD\tilde \Psi-\DD\Psi)(\DD\Psi)^{-1}||_r\\
        &\lesssim ||(\DD\tilde \Psi)^{-1}||_r||\DD\tilde \Psi-\DD\Psi||_0||(\DD\Psi)^{-1}||_0+||(\DD\tilde \Psi)^{-1}||_0||\DD\tilde \Psi-\DD\Psi||_r||(\DD\Psi)^{-1}||_0\\
        &+||(\DD\tilde \Psi)^{-1}||_0||\DD\tilde \Psi-\DD\Psi||_0||(\DD\Psi)^{-1}||_r\\
        &\lesssim \lambda_q^r\tilde L(r+2)\mathcal{T}_g.\\
    \end{split}
\end{equation*}
The implicit constants in the above bounds depend on $r$ and the data in \eqref{nonvanishing} and \eqref{stabilitydata}.
\end{proof}


\section{Inductive Lemma}\label{inductives}
In this section, given a time-dependent map $F_0:\mathbb{T}^3\times (t_0-\tau^c,t_0+\tau^c)\to \mathbb{R}^3$, which we think of as a time-dependent 1-form, we wish to study the structure and give estimates for:
$$F_k=\mathcal{L}^1_{\sigma_k}\dots\mathcal{L}_{\sigma_1}^1F_0,$$
for an arbitrary $k$, where $\{\sigma_i\}_{i}$ is a family of divergence-free vector fields defined as:
\begin{equation}\label{family}
        \sigma_i=\frac{1}{\lambda_{q+1}}\curl\left[\varsigma_i\Psi^{1*}(\varphi_i\nu)\right]=\varsigma_i\Psi^{2*}(\varphi_i'\zeta)+\frac{1}{\lambda_{q+1}}\nabla \varsigma_i\times \Psi^{1*}(\varphi_i\nu)=\sigma_i^o+\sigma_i^c
\end{equation}
with $$\varphi_i=\varphi_i(\lambda_{q+1} x\cdot k)$$ 
for $\varphi:\mathbb{R}\to\mathbb{R}$ a smooth scalar function, $(k,\nu,\zeta)$ an oriented orthonormal basis e.g. the ones from Lemma \ref{geomlemma}, $\Psi(\cdot,t):\mathbb{T}^3\to \mathbb{R}^3$ is a diffeomorphism onto its image for every $t\in B_{\tau^c}(t_0)$, e.g. the charts constructed in Lemma \ref{chartconstr} and $\varsigma_i:\mathbb{T}^3\times (t_0-\tau^c,t_0+\tau^c)\to \mathbb{R}$ are smooth functions e.g. the slow coefficients in \eqref{slowcoeff}. We also define the following types of geometries for a 1-form $F_0$:
\begin{equation}\label{geometry}
    F_0= \underbrace{a_1\Psi^{1*}(\phi_1\nu)}_{\text{good}}+\underbrace{a_2 \Psi^{1*}(\phi_2k)}_{\text{good}_c}+\underbrace{\Psi^*(\phi_3) Y}_{\text{bad}}
\end{equation}
where as above $\phi_i=\phi_i(\lambda_{q+1} x\cdot k)$ and $Y:\mathbb{T}^3\times (t_0-\tau^c,t_0+\tau^c)\to \mathbb{R}^3$ is another smooth time dependent 1-form. For reasons that will be clear in a moment, we refer to the first two terms as having `good geometry' and the third one as having `bad geometry', in the sense that we are not controlling it. We adopt the symbolic notations $\text{good},\  \text{good}_c, \ \text{bad}$ to indicate this. In the applications, we will always have $a_2=0$, which is why we formulate Lemma \ref{inductive} this way. Note that this is not restrictive, as one could always collect the $\text{good}_c$ term in the one with bad geometry. This term, however, plays a crucial role in the induction argument we are about to make.  

\begin{definition}[Admissible Loss Functions]\label{admissible} We say that a non-decreasing functions $L:\mathbb{N}_{\geq 0}\to \mathbb{R}_{\geq 1}$ is an \textit{admissible loss function} if
$$L(0)=L(1)=L(2)=1 \ \text{ and } \ \left(\frac{\lambda_q}{\lambda_{q+1}}\right)^{r'}L(r+r')\leq L(r) \ \text{ for } \ r,r'\geq 0.$$
\end{definition}

We now present two examples of admissible loss functions that will appear in this work.

\noindent \textit{Example 1.} The basic example is the following:
$$L(r)=\lambda_q^{[r-\underline{r}]^+(b-1)\gamma_\ell}$$
for some $\underline{r}\geq2$.

\noindent \textit{Example 2.} The following will appear as the loss associated with the Alfv\'en transport bound after the Galbrun Stage: 
\begin{equation}\label{admissibilityofloss}
    L_{p,\mathcal{A}}(r)=1_{r\leq \underline{r}-k_0^g-2}+1_{\underline{r}-k_0^g-1\leq r\leq \underline{r}-1} \left(\lambda_q^{(b-1)\gamma_\ell}\mathcal{T}_g\right)^{\underline{r}-1-r}\bar L+1_{ r\geq \underline{r}}\lambda_q^{[r-(\underline{r}-1)](b-1)\gamma_\ell}\bar L,
\end{equation}
one should think of it as showing the existence of three ranges of derivatives:
\begin{itemize}
    \item In the first one, we have the desired estimate.
    \item In the second one, we partially lose the good estimate, but we still have some smallness because of the Lie-Taylor expansion used in the Galbrun Stage.
    \item In the third, we have a full loss, and we even start to see the mollification parameter appear.
\end{itemize}
We refer to \eqref{lossparameters} and the Claim afterwards for a proof that this is an admissible loss function.


We now specify the exact assumptions on the initial data $F_0$ and the sequence $\{\sigma_i\}_i$.
\begin{assumptions}\label{inductivelemmaassumptions} Fix $\bar N\geq 2$ integer, let $A_i,\ A_{i,\mathcal{A}}$ with $i=1, 3$ and $\bar \varsigma_i, \ \bar\varsigma_{i,\mathcal{A}}$ with $i=1,\dots,k$ be positive real numbers and $L, \ L^1, \ L^\varsigma, \ L_{\mathcal{A}}, \ L_{\mathcal{A}}^1:\mathbb{N}_{\geq 0}\to \mathbb{R}_{\geq 1}$ be non-decreasing functions with: 
\begin{subequations}
    \begin{align}
        &L^\varsigma(r')L^1(r)\leq L^1(r+r'), \ L^\varsigma(r')L(r)\leq L(r+r')\\
        &L^\varsigma(r')L^1_{\mathcal{A}}(r)\leq L^1_{\mathcal{A}}(r+r'), \ L^\varsigma(r')L_{\mathcal{A}}(r)\leq L_{\mathcal{A}}(r+r')\\
        &L(j)=L^1(j)=L_{\mathcal{A}}(j)=L^1_{\mathcal{A}}(j)=1 \ \text{ for } j=0,1,2
    \end{align}
\end{subequations}
and assume in addition that $L^\varsigma$ is admissible as in Definition \ref{admissibilityofloss}. Let $F_0$ be as in \eqref{geometry} with $a_2=0$ and the family $\{\sigma_i\}_i$ be as in \eqref{family}. We require the following estimates:
\begin{itemize}
    \item \textbf{Fast Coefficients:}
    \begin{equation}
        ||\varphi_i||_r,||\phi_i||_r\leq 1
    \end{equation}
    for $r \geq 0$ and without $\lambda_{q+1}$ rescaling.
    \item \textbf{Slow Coefficients:}
    \begin{subequations}
        \begin{align}
            &||\partial_t^ja_1||_r\leq \bar C \lambda_q^{r+j}L^1(r+j) A_1 \ \text{ for } \ j=0,1,2 \ \text{ and } \ 0\leq r\leq \bar N-j,\\
            &||\partial_t^j Y||_r\leq \bar C \lambda_q^{r+j}L(r+j) A_3 \ \text{ for } \ j=0,1,2 \ \text{ and } \ 0\leq r\leq \bar N-j,\\
            &||\partial_t^j\varsigma_i||_r\leq \bar C \lambda_{q}^{r+j}L^\varsigma(r+j) \  \bar\varsigma_i \ \text{ for } \ j=0,1,2 \ \text{ and } \ 0\leq r\leq \bar N-j,\\
            &||\mathcal{A}^\pm a_1||_r\leq \bar C  \lambda_q^r L_{\mathcal{A}}^1(r)A_{1,\mathcal{A}} \ \text{ for } \ 0\leq r\leq  \bar N -1,\\
            &||(\partial_t+\mathcal{L}_{z^\pm}) Y||_r \leq \bar C \lambda_q^r L_{\mathcal{A}}(r)A_{3,\mathcal{A}} \ \text{ for } \ 0\leq r\leq  \bar N -1,\\
            &||\mathcal{A}^\pm \varsigma_i||_r\leq \bar C \lambda_q^rL^\varsigma(r)\ \bar\varsigma_{i,\mathcal{A}} \ \text{ for } \ 0\leq r\leq \bar N-1.
        \end{align}
    \end{subequations}
    for some constant $\bar C$.
    \item  \textbf{Charts:} 
    \begin{subequations}
        \begin{align}
            &||\DD\Psi-\IId||_0,||(\DD\Psi)^{-1}-\IId||_0\leq 1/2,\\
            &||\partial_t^j\DD\Psi||_r,||\partial_t^j(\DD\Psi)^{-1}||_r\leq \bar C \lambda_{q}^{r+j}L^\varsigma(r+j)\ \text{ for } \ j=0,1,2 \ \text{ and } \ 0\leq r\leq \bar N-j.
        \end{align}
    \end{subequations}
    for $\bar C$ as above.
    \item \textbf{Lie-Transport Properties:} 
    \begin{equation}
        (\partial_t+\mathcal{L}_{z^\pm})\Psi^*(\phi_i)= (\partial_t+\mathcal{L}_{z^\pm})\Psi^*(\varphi_i)=(\partial_t+\mathcal{L}_{z^\pm})\Psi^{1*}\nu=0.
    \end{equation}
    for any constant vector $\nu$ and $\varphi_i,\ \phi_i$ as above.
\end{itemize}

\noindent \textbf{Notation:}  In the following, we write $F\simeq G$, for $F$ and $G$ 1-forms to mean that $F$ can be written as a finite sum of terms of the same geometry as $G$  (in the sense of \eqref{geometry}) and obeying the same estimates.
\end{assumptions}

With all the preparations done, we can state the following key result.
\begin{lemma}[Inductive Lemma: 1-Forms]\label{inductive}Under the Standing Assumptions \ref{inductivelemmaassumptions} and the notation above, we have:
$$F_k\simeq a_{1,k}\Psi^{1*}(\phi_{1,k}\nu)+a_{2,k} \Psi^{1*}(\phi_{2,k}k)+\Psi^{1*}(\phi_{3,k}) Y_k$$
where the number of terms hidden in $\simeq$ grows at most exponentially in $k$.

\noindent Moreover, there exists constants $C, \ C', \ C''$, depending only on $\bar C, \ \bar N$, such that for $j=0,1,2 \text{ and } 0\leq r+k+j\leq \bar N$, we have:
    \begin{itemize}
    \item \textbf{Fast Coefficients:} $$||\phi_{i,k}||_r\leq C'' \ \text{ for } \ i=1,2,3.$$ 
    \item \textbf{Slow Coefficients:}  
        \begin{subequations}
            \begin{align}
                &||\partial_t^ja_{1,k}||_r\leq C'(C)^k \lambda_{q}^{r+k+j}L^1(r+k+j)A_1\prod_{i=1}^k\bar \varsigma_i, \\
                &||\partial_t^j a_{2,k}||_r\leq C'(C)^k \frac{\lambda_{q+1}}{\lambda_q}\lambda_{q}^{r+k+j}L(r+k+j-1)A_3\prod_{i=1}^k\bar \varsigma_i, \\
                &||\partial_t^jY_k||_r\leq C'(C)^k \lambda_{q}^{r+k+j}L(r+k+j) A_3\prod_{i=1}^k\bar \varsigma_i. 
            \end{align}
        \end{subequations}
\end{itemize}
In addition, for possibly different $\phi_{i,k}$, we have: 
$$(\partial_t+\mathcal{L}_{z^\pm})F_k\simeq a_{1,k,\mathcal{A}}\Psi^*(\phi_{1,k}\nu)+a_{2,k,\mathcal{A}} \Psi^{1*}(\phi_{2,k}k)+\Psi^*(\phi_{3,k}) Y_{k,\mathcal{A}}$$
where the number of terms hidden in $\simeq$ grows at most exponentially in $k$.

\noindent Moreover, there exists possibly different constants $C,\ C', \ C''$ such that for $0\leq r+k\leq \bar N-1$, we have:
\begin{itemize}
    \item \textbf{Fast Coefficients:} $$(\partial_t+\mathcal{L}_{z^\pm})\Psi^*(\phi_{i,k})=0 \ \text{ and } \ ||\phi_{i,k}||_r\leq C'' \ \text{ for } \ i=1,2,3.$$ 
    \item \textbf{Slow Coefficients:} 
    \begin{subequations}
        \begin{align}
            ||a_{1,k,\mathcal{A}}||_r
            &\le C'(C)^k \lambda_q^{r+k}
            \Bigl[
                L_{\mathcal{A}}^1(k+r)A_{1,\mathcal{A}}\prod_{i=1}^k \bar\varsigma_i
                + L^1(k+r)A_1
                \max_{j\in\{1,\dots,k\}}
                \Bigl[
                    \bar\varsigma_{j,\mathcal{A}}
                    \prod_{i=1,\dots,k,\; i\neq j}\bar\varsigma_i
                \Bigr]
            \Bigr], \\
            ||a_{2,k,\mathcal{A}}||_r
            &\le C'(C)^k \frac{\lambda_{q+1}}{\lambda_q}\lambda_q^{r+k}
            \begin{aligned}[t]
            \Bigl[
                &L_{\mathcal{A}}(r+k-1)A_{3,\mathcal{A}}\prod_{i=1}^k \bar\varsigma_i \\
                &\quad + L(r+k-1)A_3
                \max_{j\in\{1,\dots,k\}}
                \Bigl[
                    \bar\varsigma_{j,\mathcal{A}}
                    \prod_{i=1,\dots,k,\; i\neq j}\bar\varsigma_i
                \Bigr]
            \Bigr],
            \end{aligned}
                \\
                ||Y_{k,\mathcal{A}}||_r
                &\le C'(C)^k \lambda_q^{r+k}
                \Bigl[
                    L_{\mathcal{A}}(k+r)A_{3,\mathcal{A}}\prod_{i=1}^k \bar\varsigma_i
                    + L(k+r)A_3
                    \max_{j\in\{1,\dots,k\}}
                    \Bigl[
                        \bar\varsigma_{j,\mathcal{A}}
                        \prod_{i=1,\dots,k,\; i\neq j}\bar\varsigma_i
                    \Bigr]
                \Bigr].
        \end{align}
    \end{subequations}
\end{itemize}
\end{lemma}
\begin{proof}[Proof of Lemma \ref{inductive}] We prove the Lemma by induction on $k$, the case $k=0$ being contained in our Standing Assumptions \ref{inductivelemmaassumptions}, with the understanding that $a_{2,0}=a_{2,0,\mathcal{A}}=0$ so the bounds trivially hold. Note that for $k=0$ we can just pick $C'=\bar C, \ C=C''=1$.

\noindent To simplify the notation, we set:
\begin{equation*}
    \begin{split}
        &A_{1,r,k,j}=\lambda_{q}^{r+k+j}L(r+k+j)A_1\prod_{i=1}^k\bar \varsigma_i, \\
        &A_{2,r,k,j}= \lambda_{q+1}\lambda_{q}^{r+k+j-1}L(r+k+j-1) A_3\prod_{i=1}^k\bar \varsigma_i, \\
        &A_{3,r,k,j}= \lambda_{q}^{r+k+j}L(r+k+j)A_3\prod_{i=1}^k\bar \varsigma_i \\
    \end{split}
\end{equation*}
and we will omit the entry $j$ if $j=0$. We now assume the claimed estimates hold for $k$, where the values of the constants $C,\ C'$ will be defined along the induction step.  We also assume the decomposition:
$$F_k=a_{1,k}\Psi^{1*}(\phi_{1,k}\nu)+a_{2,k} \Psi^{1*}(\phi_{2,k}k)+\Psi^{1*}(\phi_{3,k}) Y_k$$
where for definiteness, we wrote $=$ instead of $\simeq$, it will be clear from the following computations that the number of terms grows at most exponentially in $k$ as claimed and we then split:
$$F_{k+1}=\mathcal{L}_{\sigma_{k+1}}F_k=\sigma_{k+1}\cn F_k+\DD \sigma_{k+1}^\top[F_k]=\underbrace{\curl [F_k]\times \sigma_{k+1}}_{T_1}+\underbrace{\nabla (F_k\cdot \sigma_{k+1})}_{T_2}$$
\textbf{Rewriting of $T_1$.} We begin with a sequence of preliminary computations:
\begin{equation*}
    \begin{split}
        \curl [F_k]&=\curl\left[a_{1,k}\Psi^{1*}(\phi_{1,k}\nu)+a_{2,k} \Psi^{1*}(\phi_{2,k}k)+\Psi^{1*}(\phi_{3,k}) Y_k\right]\\
        &=\lambda_{q+1} a_{1,k}\Psi^{2*}(\phi_{1,k}'\zeta)+\nabla a_{1,k}\times \Psi^{1*}(\phi_{1,k}\nu)+\nabla a_{2,k} \times \Psi^{1*}(\phi_{2,k}k)\\
        &+\lambda_{q+1} \Psi^{1*}(\phi_{3,k}'k)\times Y^k +\Psi^*(\phi_{3,k})\curl \ Y_k\\
    \end{split}
\end{equation*}
where we used $$\curl[\Psi^{1*}(\phi_{2,k}k)]=\Psi^{2*}\curl[\phi_{2,k}k]=\lambda_{q+1}\Psi^{2*}(\phi'_{2,k}k\times k)=0$$
see \eqref{DG0}. Recall that:
$$\sigma_{k+1}=\varsigma_{k+1}\Psi^{2*}(\varphi_{k+1}'\zeta)+\frac{1}{\lambda_{q+1}}\nabla \varsigma_{k+1}\times \Psi^{1*}(\varphi_{k+1}\nu)=\sigma_{k+1}^o+\sigma_{k+1}^c,$$
we now compute the cross product with the principal part of $\sigma_{k+1}$, namely
\begin{equation}\label{sigmao}
    \begin{split}
        \curl [F_k]\times \sigma_{k+1}^o&=\left[\lambda_{q+1} a_{1,k}\Psi^{2*}(\phi_{1,k}'\zeta)+\nabla a_{1,k}\times \Psi^{1*}(\phi_{1,k}\nu)+\nabla a_{2,k} \times \Psi^{1*}(\phi_{2,k}k)\right]\times \varsigma_{k+1}\Psi^{2*}(\varphi_{k+1}'\zeta)\\
        &+\left[\lambda_{q+1} \Psi^{1*}(\phi_{3,k}'k)\times Y^k +\Psi^*(\phi_{3,k})\curl \ Y_k\right]\times \varsigma_{k+1}\Psi^{2*}(\varphi_{k+1}'\zeta)\\
        &=[\varsigma_{k+1}\Psi^{2*}(\varphi_{k+1}'\zeta)\cn a_{1,k}]\Psi^{1*}(\phi_{1,k}\nu)+[\varsigma_{k+1}\Psi^{2*}(\varphi_{k+1}'\zeta)\cn a_{2,k}]\Psi^{1*}(\phi_{2,k}k)\\
        &-\lambda_{q+1} [Y_k\cdot (\varsigma_{k+1}\Psi^{2*}(\varphi_{k+1}'\zeta))]\Psi^{1*}(\phi_{3,k}'k)+\Psi^*(\phi_{3,k})\curl \ Y_k\times\varsigma_{k+1}\Psi^{2*}(\varphi_{k+1}'\zeta)
    \end{split}
\end{equation}
where we used 
\begin{equation}\label{cross1}
    (a\times b)\times c=b(a\cdot c)-a(b\cdot c)
\end{equation}
and that according to \eqref{DG5} we have: 
\begin{equation}\label{cancellation}
    \begin{split}
         \Psi^{1*}(\phi_{1,k}\nu)\cdot \Psi^{2*}(\varphi_{k+1}'\zeta)&=\det[\DD \Psi]\Psi^{1*}(\phi_{1,k}\nu)\cdot \Psi^{*}(\varphi_{k+1}'\zeta)=\det[\DD \Psi]\Psi^{*}(\phi_{1,k}\varphi_{k+1}'\nu\cdot \zeta)=0,\\
         \Psi^{1*}(\phi_{1,k}'k)\cdot \Psi^{2*}(\varphi_{k+1}'\zeta)&=\dots=0.
    \end{split}
\end{equation}
We now compute the cross product with the corrector part of $\sigma_{k+1}$, namely
\begin{equation}\label{sigmap}
    \begin{split}
        \curl [F_k]\times \sigma^c_{k+1}&=\left[\lambda_{q+1} a_{1,k}\Psi^{2*}(\phi_{1,k}'\zeta)+\nabla a_{1,k}\times \Psi^{1*}(\phi_{1,k}\nu)\right]\times\left[\frac{1}{\lambda_{q+1}}\nabla \varsigma_{k+1}\times \Psi^{1*}(\varphi_{k+1}\nu)\right]\\
        &+\nabla a_{2,k} \times \Psi^{1*}(\phi_{2,k}k)\times \left[\frac{1}{\lambda_{q+1}}\nabla \varsigma_{k+1}\times \Psi^{1*}(\varphi_{k+1}\nu)\right]\\
        &+\left[\lambda_{q+1} \Psi^{1*}(\phi_{3,k}'k)\times Y^k +\Psi^*(\phi_{3,k})\curl \ Y_k\right]\times \left[\frac{1}{\lambda_{q+1}}\nabla \varsigma_{k+1}\times \Psi^{1*}(\varphi_{k+1}\nu)\right]\\
        &=-(a_{1,k}\Psi^{2*}\zeta\cn \varsigma_{k+1})\Psi^*(\phi_{1,k}'\varphi_{k+1}\nu)\\
        &+\frac{1}{\lambda_{q+1}}\left[\Psi^{1*}\nu\cdot\left(\nabla a_{1,k}\times \nabla \varsigma_{k+1} \right)\right]\Psi^{1*}(\phi_{1,k}\varphi_{k+1} \nu)\\
        &+(\sigma_{k+1}^c\cn a_{2,k}) \Psi^{1*}(\phi_{2,k}k)-(\Psi^{1*}(\phi_{2,k}k)\cdot \sigma_{k+1}^c) \nabla a_{2,k}\\
        &+\Psi^*(\phi_{3,k})\curl \ Y_k\times \sigma_{k+1}^c+\left[\Psi^{1*}(\phi_{3,k}'k)\times Y^k\right]\times\left[\nabla \varsigma_{k+1}\times \Psi^{1*}(\varphi_{k+1}\nu)\right]
    \end{split}
\end{equation}
where we used the cancellations \ref{cancellation} above and the cross-product rule
$$(a\times b)\times (c\times b)=b\cdot (a\times c)b$$
in addition to \eqref{cross1}.

\noindent Summing \eqref{sigmao} and \eqref{sigmap} we rewrite $T_1$ as: 
\begin{equation*}
    \begin{split}
    T_1&=\curl \ F_k\times (\sigma_{k+1}^o+\sigma_{k+1}^c) \\
    &=\underbrace{(\sigma_{k+1}^o\cn a_{1,k})\Psi^{1*}(\phi_{1,k}\nu)}_{\text{good}}-\underbrace{(a_{1,k}\Psi^{2*}\zeta\cn \varsigma_{k+1})\Psi^*(\phi_{1,k}'\varphi_{k+1}\nu)}_{ \text{good}}+\underbrace{\frac{1}{\lambda_{q+1}}\left[\Psi^{1*}\nu\cdot\left(\nabla a_{1,k}\times \nabla \varsigma_{k+1} \right)\right]\Psi^{1*}(\phi_{1,k}\varphi_{k+1} \nu)}_{\text{good}}\\
    &+\underbrace{(\sigma_{k+1}\cn a_{2,k}) \Psi^{1*}(\phi_{2,k}k)}_{\text{good}_c}-\underbrace{\lambda_{q+1} (Y_k\cdot \sigma_{k+1}^o)\Psi^{1*}(\phi_{3,k}'k)}_{\text{good}_c}-\underbrace{(\Psi^{1*}(\phi_{2,k}k)\cdot \sigma_{k+1}^c) \nabla a_{2,k}}_{\text{bad}}\\
    &+\underbrace{\Psi^*(\phi_{3,k})\curl \ Y_k\times \sigma_{k+1}}_{\text{bad}}+\underbrace{\left[\Psi^{1*}(\phi_{3,k}'k)\times Y_k\right]\times \left[\nabla\varsigma_{k+1}\times \Psi^{1*}(\varphi_{k+1} \nu)\right]}_{\text{bad}}.
\end{split}
\end{equation*}

\textbf{Rewriting of $T_2$.} Similarly, we compute:
\begin{equation*}
    \begin{split}
        T_2&=  \nabla\left[F_k\cdot \left(\varsigma_{k+1}\Psi^{2*}(\varphi_{k+1}'\zeta)+\frac{1}{\lambda_{q+1}}\nabla \varsigma_{k+1}\times \Psi^{1*}(\varphi_{k+1}\nu)\right)\right]\\
        &=\nabla\underbrace{\left[a_{1,k}\Psi^{1*}(\phi_{1,k}\nu)\cdot \left(\varsigma_{k+1}\Psi^{2*}(\varphi_{k+1}'\zeta)+\frac{1}{\lambda_{q+1}}\nabla \varsigma_{k+1}\times \Psi^{1*}(\varphi_{k+1}\nu)\right)\right]}_{=0}\\
        &+\nabla\left[a_{2,k} \Psi^{1*}(\phi_{2,k}k)\cdot \left(\varsigma_{k+1}\Psi^{2*}(\varphi_{k+1}'\zeta)+\frac{1}{\lambda_{q+1}}\nabla \varsigma_{k+1}\times \Psi^{1*}(\varphi_{k+1}\nu)\right)\right]\\
        &+\nabla\left[\Psi^{1*}(\phi_{3,k}) Y_k\cdot \left(\varsigma_{k+1}\Psi^{2*}(\varphi_{k+1}'\zeta)+\frac{1}{\lambda_{q+1}}\nabla \varsigma_{k+1}\times \Psi^{1*}(\varphi_{k+1}\nu)\right)\right]\\
        &=-\frac{1}{\lambda_{q+1}}\nabla\left[\Psi^*(\varphi_{k+1}\phi_{2,k})\left(a_{2,k} \Psi^{2*}\zeta\cn \varsigma_{k+1}\right)\right]\\
        &+\nabla\left[ \Psi^{*}(\phi_{3,k}\varphi_{k+1}')\varsigma_{k+1}Y_k\cdot \Psi^{2*}\zeta+ \frac{1}{\lambda_{q+1}}\Psi^{*}(\phi_{3,k}\varphi_{k+1})Y_k\cdot(\nabla \varsigma_{k+1}\times \Psi^{1*}\nu)\right]\\
        &=-\underbrace{\left[a_{2,k} \Psi^{2*}\zeta\cn \varsigma_{k+1}\right]\Psi^{1*}((\varphi_{k+1}\phi_{2,k})'k)}_{\text{good}_c}-\underbrace{\frac{1}{\lambda_{q+1}}\Psi^*(\varphi_{k+1}\phi_{2,k})\nabla\left[a_{2,k} \Psi^{2*}\zeta\cn \varsigma_{k+1}\right]}_{\text{bad}}\\
        &+\underbrace{\lambda_{q+1}\varsigma_{k+1}Y_k\cdot \Psi^{2*}\zeta\Psi^{1*}((\phi_{3,k}\varphi_{k+1}')'k)+ Y_k\cdot(\nabla \varsigma_{k+1}\times \Psi^{1*}\nu)\Psi^{1*}((\phi_{3,k}\varphi_{k+1})'k)}_{\text{good}_c}\\
        &+\underbrace{\Psi^{*}(\phi_{3,k}\varphi_{k+1}')\nabla\left[\varsigma_{k+1}Y_k\cdot \Psi^{2*}\zeta\right]+ \frac{1}{\lambda_{q+1}}\Psi^{*}(\phi_{3,k}\varphi_{k+1})\nabla\left[Y_k\cdot(\nabla \varsigma_{k+1}\times \Psi^{1*}\nu)\right]}_{\text{bad}}\\
    \end{split}
\end{equation*}
where we used the identity in \eqref{DG4} in addition to the cancellations \eqref{cancellation} above.

\noindent We now collect the terms in $T_1, \ T_2$ with the same geometries and prove bounds using the inductive assumptions on $F_k$ and the Standing Assumptions \ref{inductivelemmaassumptions} on $\sigma_{k+1}$.

\noindent \textbf{$\text{good}_c$-geometry.} The new term with $\text{good}_c$ geometry is given by:
\begin{equation*}
    \begin{split}
        \text{good}_c&=(\sigma_{k+1}\cn a_{2,k}) \Psi^{1*}(\phi_{2,k}k)-\lambda_{q+1} (Y_k\cdot \sigma_{k+1}^o)\Psi^{1*}(\phi_{3,k}'k)\\
        &-\left[a_{2,k} \Psi^{2*}\zeta\cn \varsigma_{k+1}\right]\Psi^{1*}((\varphi_{k+1}\phi_{2,k})'k)\\
        &+\lambda_{q+1}\varsigma_{k+1}Y_k\cdot \Psi^{2*}\zeta\Psi^{1*}((\phi_{3,k}\varphi_{k+1}')'k)+Y_k\cdot(\nabla \varsigma_{k+1}\times \Psi^{1*}\nu)\Psi^{1*}((\phi_{3,k}\varphi_{k+1})'k)\\
        &=(\sigma_{k+1}\cn a_{2,k}) \Psi^{1*}(\phi_{2,k}k)\\
        &-\left[a_{2,k} \Psi^{2*}\zeta\cn \varsigma_{k+1}\right]\Psi^{1*}((\varphi_{k+1}\phi_{2,k})'k)\\
        &+\lambda_{q+1}\varsigma_{k+1}Y_k\cdot \Psi^{2*}\zeta\Psi^{1*}((\phi_{3,k}\varphi_{k+1}'')k)+Y_k\cdot(\nabla \varsigma_{k+1}\times \Psi^{1*}\nu)\Psi^{1*}((\phi_{3,k}\varphi_{k+1})'k)\\
        &=[\varsigma_{k+1}\Psi^{2*}\zeta\cn a_{2,k}] \Psi^{1*}(\phi_{2,k}\varphi_{k+1}'k)+\frac{1}{\lambda_{q+1}}[(\nabla \varsigma_{k+1}\times \Psi^{1*}\nu)\cn a_{2,k}]\Psi^{1*}(\phi_{2,k}\varphi_{k+1}k)\\
        &-\left[a_{2,k} \Psi^{2*}\zeta\cn \varsigma_{k+1}\right]\Psi^{1*}((\varphi_{k+1}\phi_{2,k})'k)\\
        &+\lambda_{q+1}\varsigma_{k+1}Y_k\cdot \Psi^{2*}\zeta\Psi^{1*}((\phi_{3,k}\varphi_{k+1}'')k)+Y_k\cdot(\nabla \varsigma_{k+1}\times \Psi^{1*}\nu)\Psi^{1*}((\phi_{3,k}\varphi_{k+1})'k)\\
        &\simeq a_{2,k+1}\Psi^{1*}(\phi_{2,k+1}k)
    \end{split}
\end{equation*}
and we have can bound:
\begin{equation*}
    \begin{split}
        ||a_{2,k+1}||_r&\lesssim ||\varsigma_{k+1}\Psi^{2*}\zeta\cn a_{2,k}||_r+\frac{1}{\lambda_{q+1}}||(\nabla\varsigma_{k+1}\times\Psi^{1*}\nu)\cn a_{2,k}||_r+||a_{2,k} \Psi^{2*}\zeta\cn \varsigma_{k+1}||_r\\
        &+\lambda_{q+1} ||\varsigma_{k+1} Y_k\cdot \Psi^{2*}\zeta||_r+||Y_k\cdot(\nabla \varsigma_{k+1}\times \Psi^{1*}\nu)||_r.
    \end{split}
\end{equation*}
We remark that terms with bad geometry can become $\text{good}_c$ and might gain a $\lambda_{q+1}$, and when this happens, the slow coefficients are not hit by derivatives. This, together with the fact that for $k=0$ we have $a_2=0$, explains why the loss function is shifted for the $\text{good}_c$ geometry. Moreover, this loss is incurred only on top of a better estimate; see the inductive assumptions on $A_{3,k,r}$. We now provide bounds on the most relevant terms and, afterwards, state the full bound.
\begin{equation*}
    \begin{split}
        ||\varsigma_{k+1}\Psi^{2*}\zeta\cn a_{2,k}||_r&\lesssim ||\varsigma_{k+1}||_r||\Psi^{2*}\zeta||_0||a_{2,k}||_1+||\varsigma_{k+1}||_0||\Psi^{2*}\zeta||_r||a_{2,k}||_1\\
        &+||\varsigma_{k+1}||_0||\Psi^{2*}\zeta||_0||a_{2,k}||_{r+1}\\
        &\lesssim \lambda_q^rL^\varsigma(r)\bar\varsigma_{k+1}A_{2,k,1}+\bar\varsigma_{k+1}A_{2,k,r+1}\\
        &\lesssim A_{2,k+1,r} \ ,
        \\
        \lambda_{q+1} ||\varsigma_{k+1} Y_k\cdot \Psi^{2*}\zeta||_r&\lesssim \lambda_{q+1} ||\varsigma_{k+1}||_r||\Psi^{2*}\zeta||_0||Y_k||_0+\lambda_{q+1} ||\varsigma_{k+1}||_0||\Psi^{2*}\zeta||_r||Y_k||_0\\
        &+\lambda_{q+1}||\varsigma_{k+1}||_0||\Psi^{2*}\zeta||_{0}||Y_k||_r\\
        &\lesssim\lambda_{q+1}\lambda_q^rL^\varsigma(r)\bar\varsigma_{k+1}A_{3,k,0}+\lambda_{q+1}\bar\varsigma_{k+1}A_{3,k,r}\\
        &\lesssim A_{2,k+1,r}\ ,
        \\
        \frac{1}{\lambda_{q+1}}||(\nabla\varsigma_{k+1}\times\Psi^{1*}\nu)\cn a_{2,k}||_r&\lesssim \frac{1}{\lambda_{q+1}}||\varsigma_{k+1}||_{r+1}||\Psi^{1*}\nu||_0||a_{2,k}||_1+\frac{1}{\lambda_{q+1}}||\varsigma_{k+1}||_1||\Psi^{1*}\nu||_r||a_{2,k}||_1\\
        &+\frac{1}{\lambda_{q+1}}||\varsigma_{k+1}||_1||\Psi^{1*}\nu||_0||a_{2,k}||_{r+1}\\
        &\lesssim \frac{\lambda_q}{\lambda_{q+1}}\lambda_q^rL^\varsigma(r+1)\bar \varsigma_{k+1}A_{2,k,1}+\frac{\lambda_q}{\lambda_{q+1}}\lambda_q^rL^\varsigma(r)L^\varsigma(1)\bar \varsigma_{k+1}A_{2,k,1}\\
        &+\frac{\lambda_q}{\lambda_{q+1}}L^\varsigma(1)\bar \varsigma_{k+1}A_{2,k,r+1}\\
        &\lesssim A_{2,k+1,r}
    \end{split}
\end{equation*}
where we used
$$L(r)L^\varsigma(r')\leq L(r+r')\ \text{ and } \ \frac{\lambda_q}{\lambda_{q+1}}L^\varsigma(r+1)\leq L^\varsigma(r)$$ 
and we conclude that:
$$||a_{2,k+1}||_r\lesssim A_{2,k+1,r} \ \text{ for } \ 0\leq r \leq \bar N-(k+1).$$
Commuting space and time derivatives in the formulas above, one can also show:
$$||\partial_t^j a_{2,k+1}||_r\lesssim A_{2,k+1,r,j} \ \text{ for } \ j=0,1,2 \ \text{ and } \ 0\leq r \leq \bar N-(k+1)-j.$$

\noindent \textbf{good-geometry.} The new term with good geometry is given by:
\begin{equation*}
    \begin{split}
        \text{good}&=(\sigma_{k+1}^o\cn a_{1,k})\Psi^{1*}(\phi_{1,k}\nu)-(a_{1,k}\Psi^{2*}\zeta\cn \varsigma_{k+1})\Psi^*(\phi_{1,k}'\varphi_{k+1}\nu)\\
        &+\frac{1}{\lambda_{q+1}}\left[\Psi^{1*}\nu\cdot\left(\nabla a_{1,k}\times \nabla \varsigma_{k+1} \right)\right]\Psi^{1*}(\phi_{1,k}\varphi_{k+1} \nu)\\
        &=(\varsigma_{k+1}\Psi^{2*}\zeta\cn a_{1,k})\Psi^{1*}(\varphi_{k+1}'\phi_{1,k}\nu)-(a_{1,k}\Psi^{2*}\zeta\cn \varsigma_{k+1})\Psi^*(\phi_{1,k}'\varphi_{k+1}\nu)\\
        &+\frac{1}{\lambda_{q+1}}\left[\Psi^{1*}\nu\cdot\left(\nabla a_{1,k}\times \nabla \varsigma_{k+1} \right)\right]\Psi^{1*}(\phi_{1,k}\varphi_{k+1} \nu)\\
        &\simeq a_{1,k+1}\Psi^{1*}(\phi_{1,k+1}\nu).
    \end{split}
\end{equation*}
Note that no indexes 2 and 3  appear. This means that if a term had a bad or $\text{good}_c$ geometry stays that way. In fact, good terms stay good. This is why we can differentiate the loss function and keep different bounds on the terms with good geometry. By construction, we have: 
\begin{equation*}
        ||a_{1,k+1}||_r\lesssim ||\varsigma_{k+1}\Psi^{2*}\zeta\cn a_{1,k}||_r+||a_{1,k}\Psi^{2*}\zeta\cn \varsigma_{k+1}||_r+\frac{1}{\lambda_{q+1}}||\Psi^{1*}\nu\cdot\left(\nabla a_{1,k}\times \nabla \varsigma_{k+1} \right)||_r
\end{equation*}
and we can estimate as above using the properties of the loss functions:
\begin{equation*}
    \begin{split}
        ||a_{1,k}\Psi^{2*}\zeta\cn \varsigma_{k+1}||_r&\lesssim ||a_{1,k}||_r||\Psi^{2*}\zeta||_0||\varsigma_{k+1}||_1+||a_{1,k}||_0||\Psi^{2*}\zeta||_r||\varsigma_{k+1}||_1\\
        &+||a_{1,k}||_0||\Psi^{2*}\zeta||_0||\varsigma_{k+1}||_{r+1}\\
        &\lesssim\lambda_qL^\varsigma(1)\bar \varsigma_{k+1}A_{1,k,r}+\lambda_q^{r+1}L^\varsigma(r) L^\varsigma(1)\bar \varsigma_{k+1}A_{1,k,0}\\
        &+\lambda_q^{r+1}L^\varsigma(r+1)\bar \varsigma_{k+1}A_{1,k,0}\\
        &\lesssim A_{1,k+1,r} \ ,
        \\
        ||\varsigma_{k+1}\Psi^{2*}\zeta\cn a_{1,k}||_r&\lesssim ||\varsigma_{k+1}||_r||\Psi^{2*}\zeta||_0||a_{1,k}||_1+||\varsigma_{k+1}||_0||\Psi^{2*}\zeta||_r||a_{1,k}||_1\\
        &+||\varsigma_{k+1}||_0||\Psi^{2*}\zeta||_0||a_{1,k}||_{r+1}\\
        &\lesssim\lambda_q^{r}L^\varsigma(r)\bar \varsigma_{k+1}A_{1,k,1}+\bar \varsigma_{k+1}A_{1,k,r+1}\\
        &\lesssim A_{1,k+1,r} \ ,
        \\
        \frac{1}{\lambda_{q+1}}||\Psi^{1*}\nu\cdot\left(\nabla a_{1,k}\times \nabla \varsigma_{k+1} \right)||_r&\lesssim \frac{1}{\lambda_{q+1}}||\Psi^{1*}\nu||_r||a_{1,k}||_1||\varsigma_{k+1}||_1+\frac{1}{\lambda_{q+1}}||\Psi^{1*}\nu||_0||a_{1,k}||_{r+1}||\varsigma_{k+1}||_1\\
        &+\frac{1}{\lambda_{q+1}}||\Psi^{1*}\nu||_0||a_{1,k}||_1||
        \varsigma_{k+1}||_{r+1}\\
        &\lesssim \frac{\lambda_q}{\lambda_{q+1}}\lambda_q^rL^\varsigma(r)L^\varsigma(1)\bar \varsigma_{k+1}A_{1,k,1}+\frac{\lambda_q}{\lambda_{q+1}}L^\varsigma(1)\bar \varsigma_{k+1}A_{1,k,r+1}\\
        &+\frac{\lambda_q}{\lambda_{q+1}}\lambda_q^rL^\varsigma(r+1)\bar \varsigma_{k+1}A_{1,k,1}\\
        &\lesssim A_{1,k+1,r}
    \end{split}
\end{equation*}
and we conclude that:
$$||a_{1,k+1}||_r\lesssim A_{1,k+1,r} \ \text{ for } \ 0\leq r \leq \bar N-(k+1).$$
Commuting space and time derivatives in the formulas above, one can also show:
$$||\partial_t^j a_{1,k+1}||_r\lesssim A_{1,k+1,r,j} \ \text{ for } \ j=0,1,2 \ \text{ and } \ 0\leq r \leq \bar N-(k+1)-j$$

\noindent \textbf{bad-geometry:} the new term with bad geometry is given by:
\begin{equation*}
    \begin{split}
        \text{bad}&=-(\Psi^{1*}(\phi_{2,k}k)\cdot \sigma_{k+1}^c) \nabla a_{2,k}+\left[\Psi^{1*}(\phi_{3,k}'k)\times Y_k\right]\times \left[\nabla\varsigma_{k+1}\times \Psi^{1*}(\varphi_{k+1} \nu)\right]+\Psi^*(\phi_{3,k})\curl \ Y_k\times \sigma_{k+1}\\
        &-\frac{1}{\lambda_{q+1}}\Psi^*(\varphi_{k+1}\phi_{2,k})\nabla\left[a_{2,k} \Psi^{2*}\zeta\cn \varsigma_{k+1}\right]\\
        &+\Psi^{*}(\phi_{3,k}\varphi_{k+1}')\nabla\left[\varsigma_{k+1}Y_k\cdot \Psi^{2*}\zeta\right]+ \frac{1}{\lambda_{q+1}}\Psi^{*}(\phi_{3,k}\varphi_{k+1})\nabla\left[Y_k\cdot(\nabla \varsigma_{k+1}\times \Psi^{1*}\nu)\right]\\
        &=\frac{1}{\lambda_{q+1}}\Psi^{*}(\phi_{2,k}\varphi_{k+1})[(\Psi^{2*}\zeta\cn\varsigma_{k+1}) \nabla a_{2,k}]+\Psi^{1*}(\phi_{3,k}'\varphi_{k+1})\left[\Psi^{1*}k\times Y_k\right]\times \left[\nabla\varsigma_{k+1}\times \Psi^{1*} \nu\right]\\
        &+\Psi^*(\phi_{3,k}\varphi_{k+1}')[\curl \ Y_k\times \varsigma_{k+1}\Psi^{2*}\zeta]+\frac{1}{\lambda_{q+1}}\Psi^*(\phi_{3,k}\varphi_{k+1})\left[\curl \ Y_k\times\left(\nabla \varsigma_{k+1}\times \Psi^{1*}\nu\right)\right]\\
        &-\frac{1}{\lambda_{q+1}}\Psi^*(\varphi_{k+1}\phi_{2,k})\nabla\left[a_{2,k} \Psi^{2*}\zeta\cn \varsigma_{k+1}\right]\\
        &+\Psi^{*}(\phi_{3,k}\varphi_{k+1}')\nabla\left[\varsigma_{k+1}Y_k\cdot \Psi^{2*}\zeta\right]+ \frac{1}{\lambda_{q+1}}\Psi^{*}(\phi_{3,k}\varphi_{k+1})\nabla\left[Y_k\cdot(\nabla \varsigma_{k+1}\times \Psi^{1*}\nu)\right]\\
        &=\Psi^{1*}(\phi_{3,k}'\varphi_{k+1})\left[\Psi^{1*}k\times Y_k\right]\times \left[\nabla\varsigma_{k+1}\times \Psi^{1*} \nu\right]\\
        &+\Psi^*(\phi_{3,k}\varphi_{k+1}')[\curl \ Y_k\times \varsigma_{k+1}\Psi^{2*}\zeta]+\frac{1}{\lambda_{q+1}}\Psi^*(\phi_{3,k}\varphi_{k+1})\left[\curl \ Y_k\times\left(\nabla \varsigma_{k+1}\times \Psi^{1*}\nu\right)\right]\\
        &-\frac{1}{\lambda_{q+1}}\Psi^*(\varphi_{k+1}\phi_{2,k})a_{2,k}\nabla\left[ \Psi^{2*}\zeta\cn \varsigma_{k+1}\right]\\
        &+\Psi^{*}(\phi_{3,k}\varphi_{k+1}')\nabla\left[\varsigma_{k+1}Y_k\cdot \Psi^{2*}\zeta\right]+ \frac{1}{\lambda_{q+1}}\Psi^{*}(\phi_{3,k}\varphi_{k+1})\nabla\left[Y_k\cdot(\nabla \varsigma_{k+1}\times \Psi^{1*}\nu)\right]\\
        &\simeq \Psi^*(\phi_{3,k+1})Y_{k+1}
    \end{split}
\end{equation*}
where we used the identity in \eqref{DG4}. By construction, we have: 
\begin{equation*}
    \begin{split}
        ||Y_{k+1}||_r&\lesssim||(\Psi^{1*}k\times Y_k)\times (\nabla\varsigma_{k+1}\times \Psi^{1*}\nu)||_r+||\curl \ Y_k\times \varsigma_{k+1}\Psi^{2*}\zeta||_r+\frac{1}{\lambda_{q+1}}||\curl \ Y_k\times( \nabla \varsigma_{k+1}\times\Psi^{1*}\nu)||_r\\
        &+\frac{1}{\lambda_{q+1}}||a_{2,k}\nabla\left[ \Psi^{2*}\zeta\cn \varsigma_{k+1}\right]||_r+||\nabla\left[\varsigma_{k+1}Y_k\cdot \Psi^{2*}\zeta\right]||_r+\frac{1}{\lambda_{q+1}}||\nabla[Y_k\cdot(\nabla \varsigma_{k+1}\times \Psi^{1*}\nu)]||_r.
    \end{split}
\end{equation*}
The key observation here is that the additional $1/\lambda_{q+1}$ smallness compensates for the presence of a $\lambda_{q+1}$ in the terms switching from $\text{good}_g\mapsto \text{bad}$ geometry, namely the one with index $2$. This allows us to close the induction argument. We give an example computation of this fact for the term: 
$$T=\frac{1}{\lambda_{q+1}}a_{2,k}\nabla\left[ \Psi^{2*}\zeta\cn \varsigma_{k+1}\right],$$ 
the remaining ones can be addressed using the same ideas
\begin{equation*}
    \begin{split}
        ||T||_r&\lesssim\frac{1}{\lambda_{q+1}}||a_{2,k}||_r(||\Psi^{2*}\zeta||_1|| \varsigma_{k+1}||_1+||\Psi^{2*}\zeta||_0|| \varsigma_{k+1}||_2)\\
        &+\frac{1}{\lambda_{q+1}}||a_{2,k}||_0(||\Psi^{2*}\zeta||_{r+1}|| \varsigma_{k+1}||_1+||\Psi^{2*}\zeta||_0|| \varsigma_{k+1}||_{r+2})\\
        &\lesssim \frac{1}{\lambda_{q+1}}A_{2,k,r}\lambda_q^2\bar \varsigma_{k+1}(L^\varsigma(1)^2+L^\varsigma(0)L^\varsigma(2))+\frac{1}{\lambda_{q+1}}A_{2,k,0}(\lambda_q^{r+2}L^\varsigma(r+1)L^\varsigma(1)\bar \varsigma_{k+1}+\lambda_q^{r+2}L^\varsigma(r+2)\bar \varsigma_{k+1})\\
        &\lesssim\frac{1}{\lambda_{q+1}}\left[\lambda_{q+1}\lambda_{q}^{r+k-1}L(r+k-1) A_3\prod_{i=1}^k\bar \varsigma_i\right]\lambda_q^2L^\varsigma(2)\bar \varsigma_{k+1}\\
        &+\frac{1}{\lambda_{q+1}}\left[\lambda_{q+1}\lambda_{q}^{k-1}L(k-1) A_3\prod_{i=1}^k\bar \varsigma_i\right]\lambda_q^{r+2}L^\varsigma(r+2)\bar \varsigma_{k+1}\\
        &\lesssim \lambda_{q}^{r+k+1}L(r+k+1)A_3\prod_{i=1}^{k+1}\bar \varsigma_i\\
        &=A_{3,k+1,r}
    \end{split}
\end{equation*}
where we used the properties of the loss functions as above, to handle the products. We conclude that:
$$||Y_{k+1}||_r\lesssim A_{3,k+1,r} \ \text{ for } \ 0\leq r \leq \bar N-(k+1).$$
Commuting space and time derivatives in the formulas above, one can also show:
$$||\partial_t^jY_{k+1}||_r\lesssim A_{3,k+1,r,j} \ \text{ for } \ j=0,1,2 \ \text{ and } \ 0\leq r \leq \bar N-(k+1)-j.$$

\noindent \textbf{Conclusion.} Gathering everything we have:
$$F_{k+1}\simeq a_{1,k+1}\Psi^*(\phi_{1,k+1}\nu)+a_{2,k+1} \Psi^{1*}(\phi_{2,k+1}k)+\Psi^*(\phi_{3,k+1}) Y_{k+1}$$
with 
\begin{itemize}
    \item Fast Coefficients (without $\lambda_{q+1}$ rescaling): $||\phi_{i,k+1}||_r\lesssim 1 $ for $r\geq 0,i=1,2,3$ .
    \item Slow Coefficients:
        \begin{equation*}
            \begin{split}
                &||\partial_t^j a_{i,k+1}||_r\lesssim A_{i,k+1,r,j} \ \text{ for } \ i=1,2\\
                &||\partial_t^jY_{k+1}||_r\lesssim A_{3,k+1,r,j}
            \end{split}
        \end{equation*}
    for $j=0,1,2$ and $0\leq r\leq \bar N-(k+1)-j$.
\end{itemize}
Inspection of the proof above shows that we can set
\begin{equation}\label{constantsinductivelemma}
    C=n_1n_2\bar C, \ C'=\bar C, \ C''=n_1
\end{equation}
where $n_1$ is the value of the implicit constant in the interpolation inequality in Lemma \ref{prodestimates} for $r=\bar N$, $n_2$ the largest number of terms in the definitions of the new slow coefficients.\\

\noindent \textbf{Alfv\'en transport.} The transport estimates follow from the first part of the lemma proved above. Indeed, the identity \eqref{DG1} gives:
\begin{equation*}
    \left(\partial_t+\mathcal{L}_{z^\pm}\right)\sigma_i=\frac{1}{\lambda_{q+1}}\curl\left[\mathcal{A}^\pm(\varsigma_i)\Psi^{1*}(\varphi_i\nu)\right],
\end{equation*}
which shows that the Lie derivative of $\sigma_i$ still has the same structure as $\sigma_i$ and with the help of the identity \eqref{trick} we obtain:
\begin{equation}\label{trick2}
    \begin{split}
        \left(\partial_t+\mathcal{L}_{z^\pm}\right)\mathcal{L}_{\sigma_k}\dots\mathcal{L}_{\sigma_1}F_0
        &=\mathcal{L}_{\sigma_k}\dots\mathcal{L}_{\sigma_1}\left(\partial_t+\mathcal{L}_{z^\pm}\right)F_0+\sum_{i=1}^{k}\mathcal{L}_{\sigma_{k}}\dots\mathcal{L}_{\sigma_{i+1}}\mathcal{L}_{\partial_t\sigma_{i}+[z^\pm,\sigma_{i}]}\mathcal{L}_{\sigma_{i-1}}\dots\mathcal{L}_{\sigma_1}F_0,
    \end{split}
\end{equation}
in particular, the Lie-transport of $\mathcal{L}_{\sigma_k}\dots\mathcal{L}_{\sigma_1}F_0$ is a sum of terms satisfying the key assumptions of the first part of the Lemma with modified sequences of $\sigma_i$ and $F_0$. The assumed bounds on $\mathcal{A}^\pm\varsigma_i,\ \mathcal{A}^\pm a_i$ and $(\partial_t+\mathcal{L}_{z^\pm})Y$ and the first part of the Lemma guarantee that:
$$\mathcal{L}_{\sigma_k}\dots\mathcal{L}_{\sigma_1}\left(\partial_t+\mathcal{L}_{z^\pm}\right)F_0\simeq a_{1,k}\Psi^*(\phi_{1,k+1}\nu)+a_{2,k} \Psi^{1*}(\phi_{2,k+1}k)+\Psi^*(\phi_{3,k+1}) Y_{k+1}$$
with
\begin{equation*}
    \begin{split}
        &||a_{1,k}||_r\lesssim \lambda_{q}^{r+k}L^1_{\mathcal{A}}(r+k)A_{1,\mathcal{A}}\prod_{i=1}^k\bar \varsigma_i \ \text{ for } \ 0\leq r\leq \bar N-k-1,\\
        &||a_{2,k}||_r\lesssim \lambda_{q+1}\lambda_{q}^{r+k-1}L_{\mathcal{A}}(k-1+r)A_{3,\mathcal{A}}\prod_{i=1}^k\bar \varsigma_i \ \text{ for } \ 0\leq r\leq \bar N-k-1,\\
        &||Y_k||_r\lesssim \lambda_{q}^{r+k}L_{\mathcal{A}}(r+k)A_{3,\mathcal{A}}\prod_{i=1}^k\bar \varsigma_i \ \text{ for } \ 0\leq r\leq \bar N-k-1.\\
    \end{split}
\end{equation*}
Similarly, applying the first part of the Lemma, we get: 
$$\mathcal{L}_{\sigma_{k}}\dots\mathcal{L}_{\sigma_{j+1}}\mathcal{L}_{\partial_t\sigma_{j}+[z^\pm,\sigma_{j}]}\mathcal{L}_{\sigma_{j-1}}\dots\mathcal{L}_{\sigma_1}F_0\simeq a_{1,k}\Psi^*(\phi_{1,k}\nu)+a_{2,k} \Psi^{1*}(\phi_{2,k}k)+\Psi^*(\phi_{3,k}) Y_{k}$$
with
\begin{equation*}
    \begin{split}
        &||a_{1,k}||_r\lesssim \lambda_{q}^{r+k}L^1(r+k)A_1\bar\varsigma_{j,\mathcal{A}}\prod_{i=1,\dots,k \ i\neq j}\bar\varsigma_i \ \text{ for } \ 0\leq r\leq \bar N-k-1,\\
        &||a_{2,k}||_r\lesssim \lambda_{q+1}\lambda_{q}^{r+k-1}L(r+k-1) A_3\bar\varsigma_{j,\mathcal{A}}\prod_{i=1,\dots,k \ i\neq j}\bar\varsigma_i \ \text{ for } \ 0\leq r\leq \bar N-k-1,\\
        &||Y_k||_r\lesssim \lambda_{q}^{r+k}L(r+k) A_3\bar\varsigma_{j,\mathcal{A}}\prod_{i=1,\dots,k \ i\neq j}\bar\varsigma_i \ \text{ for } \ 0\leq r\leq \bar N-k-1.\\
    \end{split}
\end{equation*} 

\noindent Gathering the above, we conclude that:
$$\left(\partial_t+\mathcal{L}_{z^\pm}\right)\mathcal{L}_{\sigma_k}\dots\mathcal{L}_{\sigma_1}F_0\simeq a_{1,k,\mathcal{A}}\Psi^*(\phi_{1,k}\nu)+a_{2,k,\mathcal{A}} \Psi^{1*}(\phi_{2,k}k)+\Psi^*(\phi_{3,k}) Y_{k,\mathcal{A}}$$
with
\begin{equation*}
        \begin{split}
            &||a_{1,k,\mathcal{A}}||_r\lesssim \lambda_{q}^{r+k}\left[L^1_{\mathcal{A}}(k+r)A_{1,\mathcal{A}}\prod_{i=1}^k\bar\varsigma_i +L^1(k+r)A_1\max_{j\in\{1,\dots,k\}}\left[\bar\varsigma_{j,\mathcal{A}}\prod_{i=1,\dots,k \ i\neq j}\bar\varsigma_i\right]\right],\\
            &||a_{2,k,\mathcal{A}}||_r\lesssim \lambda_{q+1}\lambda_{q}^{r+k-1}\left[L_{\mathcal{A}}(r+k-1)A_{3,\mathcal{A}}\prod_{i=1}^k\bar\varsigma_i +L(r+k-1)A_3\max_{j\in\{1,\dots,k\}}\left[\bar\varsigma_{j,\mathcal{A}}\prod_{i=1,\dots,k \ i\neq j}\bar\varsigma_i\right]\right],\\
            &||Y_{k,\mathcal{A}}||_r\lesssim \lambda_{q}^{r+k}\left[L_{\mathcal{A}}(k+r)A_{3,\mathcal{A}}\prod_{i=1}^k\bar\varsigma_i +L(k+r)A_3\max_{j\in\{1,\dots,k\}}\left[\bar\varsigma_{j,\mathcal{A}}\prod_{i=1,\dots,k \ i\neq j}\bar\varsigma_i\right]\right]
        \end{split}
    \end{equation*}
for  $0\leq r\leq \bar N-k-1, \ k\geq 1$.

\noindent Eventually, making $n_2$ larger, we can use the same definitions for the constants $C,\ C', \ C''$, given in \eqref{constantsinductivelemma}. 
\end{proof}


\begin{remark}[Inductive Lemma without the Fast Coefficients]\label{slowinductive} We want to have a version of Lemma \ref{inductive} to use in Gabrun's stage. In that setting, there are no fast coefficients, and we do not need to require any geometry on the vector fields $\sigma_i$ and the 1-form $F_0$. We can directly assume:
\begin{equation*}
        \begin{split}
            &||\partial_t^j F_0||_{r+\alpha}\lesssim \lambda_q^rL(r+j) A \ \text{ for } j=0,1,2\ \text{ and } \ 0\leq r\leq \bar N-j,\\
            &||\partial_t^j\sigma_i||_{r+\alpha}\lesssim \lambda_{q}^rL^\varsigma(r+j) \  \bar\varsigma_i \ \text{ for } j=0,1,2\ \text{ and } \ \text{ for } 0\leq r\leq \bar N-j,\\
            &||(\partial_t+\mathcal{L}_{z^\pm}) F_0||_{r+\alpha}\lesssim  \lambda_q^r L_{\mathcal{A}}(r)A_{\mathcal{A}} \ \text{ for } \ 0\leq r\leq  \bar N -1,\\
            &||(\partial_t+\mathcal{L}_{z^\pm}) \sigma_i||_{r+\alpha}\lesssim \lambda_q^rL^\varsigma(r)\ \bar\varsigma_{i,\mathcal{A}} \ \text{ for } \ 0\leq r\leq \bar N-j\\
        \end{split}
    \end{equation*}
with the same properties for $L, \ L^\varsigma, \ L_{\mathcal{A}}$ as in the Standing Assumptions \ref{inductivelemmaassumptions}. A simple induction argument then leads to the estimate:
$$||\partial_t^j\mathcal{L}_{\sigma_k}\dots \mathcal{L}_{\sigma_1}F_0||_{r+\alpha}\leq C'(C)^k \lambda_{q}^{r+k+j}L(r+k+j)A\prod_{i=1}^k\bar \varsigma_i \ \text{ for } \ j=0,1,2 \  \text{ and }\ 0\leq r+k+j\leq \bar N$$
for some constants $C, \ C'$ which depend on $\bar C, \ \bar N$ but not on $r, \ k, \ j$.
The rewriting \eqref{trick} and another induction argument then gives:
$$||(\partial_t+\mathcal{L}_{z^\pm})\mathcal{L}_{\sigma_k}\dots \mathcal{L}_{\sigma_1}F_0||_{r+\alpha}\leq C'(C)^k \lambda_{q}^{r+k}\left[L_{\mathcal{A}}(k+r)A_{\mathcal{A}}\prod_{i=1}^k\bar\varsigma_i +L(k+r)A\max_{j\in\{1,\dots,k\}}\left[\bar\varsigma_{j,\mathcal{A}}\prod_{i=1,\dots,k \ i\neq j}\bar\varsigma_i\right]\right]$$
for $0\leq r +k \leq \bar N -1$ and some possibly different constants $C, \ C'$ which depend on $\bar C, \ \bar N$ but not on $r,\ k$.
\end{remark}

\appendix

\section{Tools from Differential Geometry}\label{diffgeom} 

\noindent \textbf{Vector calculus.} Let $$U \subset \mathbb{T}^3 \text{ or } U \subset \mathbb{R}^3$$be open, and let$$a = a^i,\quad b = b^i,\quad F = F^i,\quad v = v^i : U \to \mathbb{R}^3$$ be smooth vector fields, and 
$$A = A^{ij} : U \to \mathbb{R}^{3\times 3}$$
a smooth matrix field. We use the following conventions for differential operators:
\begin{subequations}
\begin{align}
(\DD a)^{ij} &= \partial_j a^i, \\
(b \cn a)^i &= (\DD a[b])^i = b^j \partial_j a^i, \\
\ddiv a &= \partial_i a^i, \\
(\ddiv A)^i &= \partial_j A^{ij},
\end{align}
\end{subequations}
together with the identities
\begin{subequations}
\begin{align}
\ddiv(a \otimes b) &= b \cn a + (\ddiv b)\, a, \\
\curl(a \times b) &= b \cn a - a \cn b + (\ddiv b)\, a - (\ddiv a)\, b.\end{align}
\end{subequations}
Additionally, we have:
\begin{equation}\label{DG2}
    \ddiv [(F \times \nabla v)^\top]=\curl[F]\cn v
\end{equation}
where 
\begin{equation*}
    (F\times\nabla v)^{ij}= \varepsilon^{ikl}F^k\partial_lv^j
\end{equation*}
and $\varepsilon^{ikl}$ is the standard Levi-Civita tensor. 

\noindent This can be proved as follows: 
\begin{equation}
    \begin{split}
        \left(\ddiv\left[(F\times\nabla v )^\top\right]\right)^i&=\partial_j(\varepsilon^{jkl}F^k\partial_lv^i)\\
        &=\varepsilon^{jkl}\partial_jF^k\partial_lv^i+\varepsilon^{jkl}F^k\partial_j\partial_lv^i\\
        &=\left(\curl [F]\cn v\right)^i
    \end{split}
\end{equation}
where we used the definition of curl
$$\varepsilon^{jkl}\partial_jF^k\partial_lv^i=\varepsilon^{ljk}\partial_jF^k\partial_lv^i=(\curl F)^l\partial_lv^i$$
and
$$\varepsilon^{jkl}F^k\partial_j\partial_lv^i=0$$
which follows from the anti-symmetry of the Levi-Civita tensor, the symmetry of the second-order derivatives
$$\varepsilon^{jkl}=\varepsilon^{ljk}=-\varepsilon^{lkj} \ \text{ and } \ \partial_j\partial_lv^i=\partial_l\partial_jv^i$$
and summing over the dummy indices.\\

\noindent \textbf{Differential Geometry.} For any vector valued map $F:U\to\mathbb{R}^3$, scalar function $f:U\to\mathbb{R}$, vector field $v:U\to \mathbb{R}^3$ and diffeomorphism $\Psi:U'\to U$, where $U,U'$ are open subsets of $\mathbb{T}^3$ or $\mathbb{R}^3$, we set the following notation conventions:
\begin{itemize}
    \item The pullback along $\Psi$ of $F$ as a vector field is given by:
        $$\Psi^*F=(\DD \Psi)^{-1}F(\Psi).$$
    \item The pullback along $\Psi$ of $f$ is given by:
        $$\Psi^{*}f=f(\Psi).$$
    \item The pullback along $\Psi$ of $F$ as as a 1-form is given by:
        $$\Psi^{1*}F=(\DD \Psi)^{T}F(\Psi).$$
    \item The pullback along $\Psi$ of $F$ as a 2-form is given by:
        $$\Psi^{2*}F=\det[\DD\Psi](\DD \Psi)^{-1}F(\Psi).$$
    \item  The pushforward along $\Psi$ of $F$ as a vector field is given by:
        $$\Psi_*F=(\DD \Psi)(\Psi^{-1})F(\Psi^{-1})=(\Psi^{-1})^*F.$$
    \item  The `Lie-derivative' of $F$ as a vector field in the direction $v$, that is, the commutator of $F$ with a vector field $v$, is given by:
        $$\mathcal{L}_vF=[v,F]=v\cn F-F\cn v.$$
    \item  The Lie-derivative of $F$ as a 1-form in the direction $v$ is given by:
        $$\mathcal{L}^{1}_vF=\curl[F]\times v+\nabla(F\cdot v)=v\cn F+\DD v^\top F.$$
    \item  The Lie-derivative of $F$ as a 2-form in the direction $v$ is given by:
        $$\mathcal{L}^{2f}_vF=\curl[F\times v]+(\ddiv \ F)v=[v,F]+(\ddiv \ v)F.$$
\end{itemize}

\noindent Even when the map $\Psi$ is time-dependent, we always consider the pullbacks and Lie derivatives in space for each fixed time. 

\noindent We use the above notation because, in this work, we will use the identifications coming from the well-known commutative diagram:
\begin{equation}\label{identifications}
    \begin{tikzcd}[row sep=large, column sep=huge] 
0 \arrow[r] 
& C^\infty(\mathbb{T}^3) 
\arrow[r,"d"] 
\arrow[d,"\mathrm{id}"'] 
& \Omega^1(\mathbb{T}^3) 
\arrow[r,"d"]  
& \Omega^2(\mathbb{T}^3)
\arrow[r,"d"]  
& \Omega^3(\mathbb{T}^3) 
\arrow[d,"\mathrm{id}"']
\arrow[r] 
& 0
\\
0 \arrow[r] 
& C^\infty(\mathbb{T}^3) 
\arrow[u]
\arrow[r,"\nabla"] 
& \mathfrak{X}(\mathbb{T}^3)
\arrow[u,"\flat"']
\arrow[r,"\nabla \times"] 
& \mathfrak{X}(\mathbb{T}^3) 
\arrow[r,"\nabla \cdot"] 
\arrow[u,"\iota\dd \text{vol}=*\flat"']
& C^\infty(\mathbb{T}^3) 
\arrow[r]
\arrow[u]
& 0
\end{tikzcd}
\end{equation}
together with the Hodge isomorphism
$$*:\Omega^1(\mathbb{T}^3)\to \Omega^2(\mathbb{T}^3)$$
where the metric considered is the trivial one. In particular, a map $F:\mathbb{T}^3\to \mathbb{R}^3$ (possibly defined on an open subset $U\subset \mathbb{T}^3$ only) can be thought of as either a 1-form $F\in \Omega^1(\mathbb{T}^3)$, a 2-form $F\in \Omega^2(\mathbb{T}^3)$ or a vector field $F\in \mathfrak{X}(\mathbb{T}^3)$. We want, however, to take advantage of the following standard identities, see \cite[Chapters 12-14] {lee2013smooth} and \cite[Chapter 2]{MS}, for differential forms $a, \ b$ and vector fields $v, \ w$:
\begin{subequations}\label{DG-1}
    \begin{align}
        &\Psi^*\dd =\dd \Psi^*,\\
        &\Psi^*a\wedge \Psi^*b=\Psi^*(a\wedge b),\\
        &\Psi^*\mathcal{L}_{v}a=\mathcal{L}_{\Psi^*v}\Psi^*a,\\
        &\Psi^*\iota_{v}a=\iota_{\Psi^*v}\Psi^*a,\\
        &\mathcal{L}_va=\iota_v\dd a+\dd \iota_v a,\\
        &\mathcal{L}_v\mathcal{L}_w-\mathcal{L}_w\mathcal{L}_v=\mathcal{L}_{[v,w]}.
    \end{align}
\end{subequations}
Despite the notation, these depend on the degree of the forms considered; thus, the identifications in \eqref{identifications} need to be done carefully. This explains why we distinguish in what sense the map $F$ is being considered and add indices in the notation above accordingly. We now verify that the definitions and notational conventions are consistent with the identifications.

\noindent \textit{Vector Fields vs 2-Forms.} From the correspondence between 2-forms and vector fields in \eqref{identifications}, namely
$$X\in \mathfrak{X}(\mathbb{T}^3) \iff \iota_X\dd\text{vol}=*X^\flat\in \Omega^2(\mathbb{T}^3),$$
the key identities in \eqref{DG-1} and $\Psi^*\dd\text{vol}=\det[\DD\Psi]\dd\text{vol}$ we compute: 
\begin{equation*}
    \begin{split}
        \Psi^*\iota_X\dd\text{vol}=\iota_{\Psi^*X}(\Psi^*\dd\text{vol})=\iota_{\Psi^*X}(\det[\DD \Psi]\dd\text{vol})=\det[\DD \Psi]\iota_{\Psi^*X}(\dd\text{vol})=\iota_{\det[\DD \Psi]\Psi^*X}(\dd\text{vol})
    \end{split}
\end{equation*}
and deduce that 
$$\Psi^{2*}F=\det[\DD \Psi](\DD \Psi)^{-1}F(\Psi).$$
Moreover, given the geometric characterisation of divergence:
$$\dd\iota_X\dd\text{vol}=\mathcal{L}_X\dd\text{vol}=(\ddiv \ X)\dd \text{vol},$$
the identity
\begin{equation}\label{doublecontraction}
    \iota_v\iota_w\dd\text{vol}=(w\times v)^\flat,
\end{equation}
which can be easily verified and the commutativity of the diagram \eqref{identifications}, we deduce that:
\begin{equation*}
    \begin{split}
        \mathcal{L}_v(\iota_X\dd\text{vol})&=\dd(\iota_v\iota_X\dd\text{vol})+\iota_v\dd(\iota_X\dd\text{vol})\\
        &=\dd((X\times v)^\flat)+(\ddiv \ X)\iota_v\dd \text{vol}\\
        &=\iota_{\curl [X\times v]}\dd \text{vol}+\iota_{(\ddiv \ X)v}\dd \text{vol}\\
        &=\iota_{\curl [X\times v]+(\ddiv \ X)v}\dd \text{vol}\\
        &=\iota_{[v,X]+(\ddiv \ v)X}\dd \text{vol},
    \end{split}
\end{equation*}
and we conclude:
$$\mathcal{L}^{2f}_vF=[v,F]+(\ddiv \ v)F.$$

\noindent \textit{Vector Fields vs 1-Forms.} Recall that:
$$X^\flat v=X\cdot v \ \leadsto \ \iota_v X^\flat=X^\flat v=X\cdot v.$$
Now, from the commutative diagram \eqref{identifications}, Cartan's Formula in \eqref{DG-1} and \eqref{doublecontraction} we get:
\begin{equation*}
    \begin{split}
        \mathcal{L}_v X^\flat&=\dd \iota_v X^\flat+ \iota_v\dd X^\flat\\
        &=\dd(X\cdot v)+\iota_v\iota_{\curl X}\dd\text{vol}\\
        &=[\nabla(X\cdot v)]^\flat+(\curl X\times v)^\flat\\
        &=[\nabla(X\cdot v)+\curl X\times v]^\flat.
    \end{split}
\end{equation*}

\noindent Concerning the pullback, we have:
$$(\Psi^*X^\flat) (v)=X^\flat(\Psi) (\DD\Psi[v])=X(\Psi)\cdot \DD\Psi[v]=(\DD\Psi^\top X(\Psi))\cdot v=(\DD\Psi^\top X(\Psi))^\flat v$$
and we conclude that:
$$\Psi^{1*}F=\DD\Psi^\top F(\Psi) \ \ \text{ and } \ \ \mathcal{L}_v^1F=\nabla(F\cdot v)+\curl [F]\times v.$$

\noindent \textit{Products.} From the fact that:
$$Ma \times Mb=\det [M] M^{-\top}[a\times b],$$
we deduce that:
\begin{equation}\label{DG4}
    \begin{split}
        \Psi^{2*}F\times \Psi^{*}G&=\det[\DD \Psi]\left[(\DD \Psi)^{-1}F(\Psi)\times(\DD \Psi)^{-1} G(\Psi)\right]=(\DD \Psi)^\top (F\times G)(\Psi)=\Psi^{1*}(F\times G),\\
        \Psi^{1*}F\times \Psi^{1*}G&=\left[(\DD \Psi)^{T}F(\Psi)\times(\DD \Psi)^{T} G(\Psi)\right]=\det[\DD \Psi](\DD \Psi)^{-1}[F(\Psi)\times G(\Psi)]=\Psi^{2*}(F\times G).
    \end{split}
\end{equation} 
We remark that these identities are nothing but an instance of
$$\Psi^*a\wedge \Psi^*b=\Psi^*(a\wedge b)$$
from \eqref{DG-1} above.

\noindent A related property is:
\begin{equation}\label{DG5}
        \Psi^{1*}X\cdot \Psi^{*}Y=[(\DD \Psi)^\top X(\Psi)]\cdot[(\DD \Psi)^{-1} Y(\Psi)]=X(\Psi)^\top\DD \Psi(\DD \Psi)^{-1} Y(\Psi)=\Psi^*(X\cdot Y),
\end{equation} 
this is just the standard action of diffeomorphisms on the pairing between 1-forms and vector fields.

\noindent \textit{Some useful identities.} From the above we deduce the following identities: 
\begin{subequations} \label{DG0}
    \begin{align}
        &\curl \Psi^{1*}=\dd\Psi^*=\Psi^{*}\dd=\Psi^{2*}\curl,\\
        &\ddiv \ \Psi^{2*}=\dd\Psi^{*}=\Psi^{*}\dd=\Psi^{*}\ddiv,\\
        &\det[\DD \Psi]=1 \Longrightarrow \Psi^*=\Psi^{2*},\\
        & \ddiv \ v =0 \Longrightarrow \mathcal{L}^{2f}_v=\mathcal{L}_v \ \text{ and } \ \curl [\mathcal{L}_v^{1}\Theta]=\mathcal{L}_v \curl [\Theta]
    \end{align}
\end{subequations}
where the last identity follows from the definitions and the following computation:
\begin{equation}
    \begin{split}
        \curl [\mathcal{L}_v^{1}\Theta]&=\dd\left(\dd\iota_v \Theta+\iota_v\dd \Theta\right)\\
        &=\dd\left(\iota_v\dd \Theta\right)\\
        &=\curl\left[\curl \Theta\times v\right]\\
        &=\left[v,\curl \Theta\right]\\
        &=\mathcal{L}_v \left[\curl \Theta\right].
    \end{split}
\end{equation}
In particular, if $\xi=\curl\ \Theta$ and $v, \ B$ are divergence-free vector fields, we have:
\begin{equation} \label{DG1}
    \begin{cases}
        (\partial_t+\mathcal{L}_v)\xi
        =\curl \left[(\partial_t+\mathcal{L}^1_v)\Theta\right],\\
        \mathcal{L}_B\xi=\curl[\mathcal{L}^1_B\Theta].
    \end{cases}
\end{equation}

\noindent \textbf{Lie-Transport Equations.}  We recall the following standard identity; see \cite[Chapter 2]{MS} for a proof. Let $M$ be either $\mathbb{T}^3$ or $\mathbb{R}^3$, and let $\xi=\xi(x,s)$ be a smooth time-dependent vector field on $M$. Let $X_s$ be the associated flow map, namely the solution of
\begin{equation*}
\begin{cases}
    \partial_s X_s(x)=\xi(X_s(x),s),\\
    X_0(x)=x
    \end{cases}
\end{equation*}
then, for any smooth time-dependent tensor field
$$F(\cdot,s)\in \Gamma\bigl(M,(TM)^{\otimes p}\otimes (T^*M)^{\otimes q}\bigr)$$
one has
\begin{equation}\label{lietransport}    
\partial_s (X_s^*F)=X_s^*\bigl(\partial_sF+\mathcal{L}_\xi F\bigr)
\end{equation}
where the Lie derivative is understood depending on $(p,q)$.

\noindent \textit{Duhamel Formula for Lie-Transport.} From the identity \eqref{lietransport} follows that for any time-independent smooth differential form $F_0$, we have
\begin{equation}\label{duhamel}
    \begin{cases}
    \partial_s F+\mathcal{L}_\xi F=G\\
    F|_{s=0}=F_0
\end{cases} \iff F(\cdot,s)=\Phi_s^*F_0+\int_0^s \Phi_s^*X_{s'}^*[G(\cdot,s')]\dd s'
\end{equation}
where $\Phi_s=X_s^{-1}$ is the inverse flow map of $\xi$, and we implicitly assume that both sides are well defined.

\noindent Indeed assume that $F$ solves the initial value problem. From \eqref{lietransport} we obtain:
\begin{equation*}
        \partial_s (X_s^*F)=X_s^*[\partial_sF+\mathcal{L}_\xi F]=X_s^*G.
\end{equation*}
Integrating this expression from 0 to $s$, we deduce:
$$X_s^*F-F_0=\int_0^sX_{s'}^*G(\cdot ,s')\dd s',$$
now applying $\Phi_s^*=(X_s^{-1})^*$ we conclude that:
$$F(\cdot,s)=\Phi_s^*F_0+\int_0^s\Phi_s^*X_{s'}^*G(\cdot,s')\dd s'.$$

\noindent The converse implication follows by proceeding backward, applying $X^*_s$ to the Duhamel formula, taking $\partial_s$ derivative and using \eqref{lietransport} we see that:
$$X^*_s(\partial_s F+\mathcal{L}_\xi F)=X^*_sG.$$
But $X_s$ is a diffeomorphism, and the solution property follows. The fact that it has the right initial condition is immediate from the formula.

\noindent \textit{Lagrangian Flow vs Transport.} From the above, applied to the scalar case, we obtain: 
\begin{equation}\label{lagrangianvstransport}
    \begin{cases}
        \partial_sX_s=\xi(X_s,s),\\
        X_0=\IId
    \end{cases} \ \ \ \ 
    \overset{\Phi_s=X_s^{-1}}{\iff} \ \ \ \ 
    \begin{cases}
        \partial_s\Phi_s+\xi(\cdot,s)\cn\Phi_s=0,\\
        \Phi_0=\IId.
    \end{cases}
\end{equation}


\noindent \textbf{Pullback vs Lie Derivation.} Let $T$ be any smooth tensor field and $v$ any smooth vector field. We will show that:
$$\Upsilon^*(\mathcal{L}_vT)=\mathcal{L}_{\Upsilon^*v}(\Upsilon^*T).$$
We first recall that for any time-independent diffeomorphism $\Upsilon$, one has:
$$X^{v}_s\circ \Upsilon=\Upsilon\circ X^{\Upsilon^*v}_s,$$
where we denote by $X^v_s$ the Lagrangian flows of $v$ at time $s$ starting at the identity at time $0$. In our case, $v$ will depend on time, but the differentiation in the Lie derivative happens at each fixed time, i.e. we freeze the time $t$ in $v$ and compute $X^v_s=X^{v|_t}_s$ thinking of $v$ as $v|_t$, a time-independent object. With this in mind, the definition of Lie derivative and the fact that $\Psi^*\Upsilon^*=(\Upsilon\circ \Psi)^*$, we compute:
\begin{equation}\label{DG3}
    \begin{split}
        \Upsilon^*(\mathcal{L}_vT)&=\Upsilon^*(\mathcal{L}_{v|_t}T)\\
        &=\Upsilon^*\left(\frac{\dd}{\dd s}\Big|_{s=0}(X^{v|_t}_s)^*T\right)\\
        &=\frac{\dd}{\dd s}\Big|_{s=0}(X^{v|_t}_s\circ \Upsilon)^*T\\
        &=\frac{\dd}{\dd s}\Big|_{s=0}( \Upsilon\circ X^{\Upsilon^*v|_t}_s)^*T\\
        &=\frac{\dd}{\dd s}\Big|_{s=0}(X^{\Upsilon^*{v|_t}}_s)^*(\Upsilon^*T)\\
        &=\mathcal{L}_{\Upsilon^*{v|_t}}(\Upsilon^*T)\\
        &=\mathcal{L}_{\Upsilon^*v}(\Upsilon^*T)
    \end{split}
\end{equation}
and the commutation is proven.

\noindent \textbf{Transport Trick.} From the key commutation relation:
$$\mathcal{L}_F\mathcal{L}_G-\mathcal{L}_G\mathcal{L}_F=\mathcal{L}_{[F,G]}$$
in \eqref{DG-1}, we deduce a useful identity that can be used to deal with the Alfv\'en Lie-transport of high-order terms in the Lie-Taylor expansions, namely
\begin{equation}\label{trick}
    \begin{split}
        \left(\partial_t+\mathcal{L}_{v}\right)\mathcal{L}_\xi^kF&=\mathcal{L}_\xi\left(\partial_t+\mathcal{L}_{v}\right)\mathcal{L}_\xi^{k-1}\left(\partial_t+\mathcal{L}_{v}\right)\Theta+\mathcal{L}_{\partial_t\xi+[v,\xi]}\mathcal{L}_\xi^{k-1}F\\
        &=\dots\\
        &=\mathcal{L}_\xi^k\left(\partial_t+\mathcal{L}_{v}\right)F+\sum_{i=0}^{k-1}\mathcal{L}_{\xi}^i\mathcal{L}_{\partial_t\xi+[v,\xi]}\mathcal{L}_\xi^{k-1-i}F,\\
    \end{split}
\end{equation}
Similarly,  
$$\mathcal{L}_{B}\mathcal{L}_{\xi}^kF=\mathcal{L}_\xi^k\mathcal{L}_{B}F+\sum_{i=0}^{k-1}\mathcal{L}_{\xi}^i\mathcal{L}_{[B,\xi]}\mathcal{L}_\xi^{k-1-i}F.$$

\noindent \textbf{Minimal covers.} For $\rho>0$ a finite open covering $\mathcal{C}$ of a compact manifold $\mathcal{M}$ of dimension $n$ is a minimal cover of diameter $\rho$ with granularity $c_0$ if:
\begin{itemize}
    \item The diameter of each $U_j\in \mathcal{C}$ is less than $\rho$.
    \item $\mathcal{C}$ can be subdivided into $n+1$ subfamilies $\mathcal{F}_i$, each consisting of pairwise disjoint sets.
    \item $U,U'\in \mathcal{F}_i \Longrightarrow \text{dist}(U,U')>c_0$.
\end{itemize}
The following Proposition is a standard fact. We refer to \cite[Section 5.5]{CDS}. 

\begin{prop}[Minimal Coverings] \label{minimal}Let $\mathcal{M}$ be a compact manifold and $\rho>0$. Then $\mathcal{M}$ admits a minimal cover $\mathcal{C}$ of diameter $\rho$. Moreover, there exists a partition of unity $\{\theta_j\}$ subordinate to $\mathcal{C}$, satisfying:
$$\sum_j\theta_j=1, \ \theta_j\in C_c^{\infty}(U_j), \ \theta_j\geq 0, \ ||\theta_j||_r\lesssim \rho^{-r}$$
where the implicit constant depends on $r$ and the specific choice of the covering.
\end{prop}
Note that $c_0$ always exists by compactness, but for our purposes, an explicit control is needed; we will rely on the following explicit construction on $\mathbb{T}^3$, where we set $\rho=4/3\tau^c$ and $c_0=1/3\tau^c$. 

\noindent We think of $\mathbb{T}^3$ as a periodic unit cube in $\mathbb{R}^3$. We now fix the cube: 
$$Q_0=(0,2/3)^3$$
and let $\mathcal{F}_0$ be the family of cubes given by 1-periodically translating $Q_0$. We construct three other periodic families by shifting this initial one in the coordinate directions, namely
$$\mathcal{F}_l=\mathcal{F}_0+\frac{1}{2}e_l \ \ \text{ for } \ \ l=1,2,3.$$
The union $\mathcal{C}=\cup_i \mathcal{F}_i$ gives a periodic tiling of $\mathbb{R}^3$ such that:
\begin{enumerate}
    \item $U\in \mathcal{C}\Longrightarrow \text{diam}(U)=2/3\sqrt{3}$.
    \item $U,U'\in \mathcal{F}_l \Longrightarrow \text{dist}(U,U')\geq 1/3$.
    \item $\mathcal{C}$ is minimal in the sense that each point in $\mathbb{R}^3$ is contained in at most four elements of $\mathcal{C}$.
\end{enumerate}
For each $U_{j'}\in \mathcal{C}$ standard techniques allow to construct a family $\{\theta_{j'}\}$ of cut-off function with $\theta_{j'}\in C_c^{\infty}(U_{j'},[0,1])$ and such that: 
\begin{enumerate}
    \item $\sum_{j'}\theta_{j'}^2\equiv1$, in particular for each point $\mathbb{T}^3$ at most four $\theta_{j'}$ are non-zero.
    \item $$U_{j'_1},U_{j'_2}\in \mathcal{F}_l \Longrightarrow \text{dist}(\supp(\theta_{j_1'}), \ \supp(\theta_{j'_2}))>\text{dist}(U_{j_1'}, \ U_{j'_2})\geq 1/3$$ and $$\text{diam} \ \supp \ \theta_{j'}< \text{diam}(U_{j'})\leq 2/3\sqrt{3}\leq 4/3.$$
    \item $||\theta_{j'}||_r\lesssim 1$ the implicit constant depends on $r$ and the specific choice of the partition of unity.
\end{enumerate}
Appropriately rescaling the construction by $\tau^c$, so that it is still periodic, we can achieve:
\begin{equation*}\text{diam} \ \supp \ \theta_{j'}<4/3\tau^c \ \ \text{ and } \ \ U_{j'_1},U_{j'_2}\in \mathcal{F}_l \Longrightarrow \text{dist}(\supp(\theta_{j_1'}), \ \supp(\theta_{j'_2}))> 1/3\tau^c
\end{equation*}
while the other properties now read:
$$\sum_{j'}\theta_{j'}^2\equiv1 \ \ \text{ and } \ \ ||\theta_{j'}||_r\lesssim (\tau^c)^{-r}.$$

\section{Tools from Analysis}
\noindent \textbf{H{\"o}lder Functions.} Let $r\geq 0$ be an integer and $\alpha\in(0,1)$, for a possibly time dependent function $f:\mathbb{T}^3\times \mathbb{R}\to \mathbb{R}$, we denote: 
\begin{equation*}
    \begin{split}
        &||f||_0=\sup_x|f(x)| \ \text{ or } \ \sup_{x,t}|f(x,t)|,\\
        &[f]_r=\sup_{|\theta|=r}||\partial_\theta f|| \ \text{ or }\ \sup_{t}\sup_{|\theta|=r}||\partial_\theta f(\cdot,t)||,\\
        &[f]_{r+\alpha}=\sup_{|\theta|=r}\sup_{x \neq y}\frac{|\partial_\theta f(x)-\partial_\theta f(y)|}{|x-y|^\alpha} \ \text{ or }\ \sup_{t}\sup_{|\theta|=r}\sup_{x \neq y}\frac{|\partial_\theta f(x,t)-\partial_\theta f(y,t)|}{|x-y|^\alpha},\\
        &||f||_r=\sum _{j=0}^r [f]_j\ \text{ and } \ ||f||_{r+\alpha}=||f||_r+[f]_{r+\alpha}
    \end{split}
\end{equation*}
where $\theta$ denotes a multi-index, and with a slight abuse of notation, we keep the same notation for vector- and matrix-valued functions. 

\noindent The following product and interpolation estimates are classical:
\begin{equation}\label{prodestimates}
        ||fg||_{r}\lesssim ||f||_{r}||g||_0+||f||_0||g||_r
\end{equation}
where $r$ is a non-negative real number, and the implicit constants depend on $r$, and
\begin{equation}
     ||f||_{r}\lesssim ||f||_{r'}^\lambda||f||_{r''}^{1-\lambda} \ \text{ for } \ r=\lambda r'+(1-\lambda)r''\\
\end{equation}
where $r, \ r',  \ r'', \ \lambda$ are non-negative real numbers and the implicit constants depend on $r,r',r'',\lambda$.

\noindent \textbf{Miscellaneous.} We now recall some standard results in Analysis. We give the formulations from \cite[Appendix C.34]{lee2013smooth}.
\begin{prop}[Inverse Function Theorem]\label{ift} Suppose $U$ and $V$ are open subsets
of $\mathbb{R}^n$, and $F:U \to V$ is a smooth function. If $\DD F(x_0)$ is invertible at some point
$x_0\in U$, then there exist connected neighbourhoods $U_0 \subseteq U$ of $x_0$ and $V_0\subseteq V$ of $F(x_0)$ such that $F|_{U_0} : U_0 \to V_0$ is a diffeomorphism.
\end{prop}
Let $(X, d )$ be a metric space. A map $G: X \to X$ is said to be a contraction if there is a constant $c<1$ such that $d (G(x), G(y))\leq c \ d(x, y)$ for all $x,y\in X$. Clearly, every contraction is continuous. A fixed point of a map $G: X \to X$ is a point $x \in X$
such that $G(x)=x$.
\begin{prop}[Contraction Mapping Theorem]\label{CMT}Let $X$ be a nonempty complete metric space.
Every contraction $G: X \to X$ has a unique fixed point.
\end{prop}

We now recall the following composition estimate; the proof is a simple adaptation of the corresponding result in \cite{DS1} to the fractional case.
\begin{prop}[Composition Estimates]\label{compestimates} Let $f:\Omega\to\mathbb{R}$ and $\Psi:\mathbb{R}^n\to \Omega$ be two functions with $\Omega\subset \mathbb{R}^N$. Then for any $0\leq \alpha <1$ and $r\geq 0$ integer, we have:
\begin{subequations}
    \begin{align}
        [f\circ g]_\alpha&\leq [f]_\alpha [g]_1^\alpha,\\
        [f\circ \Psi]_{r+\alpha}&\lesssim \left([f]_{1+\alpha}[\Psi]_{r}+[f]_{r+\alpha}[\Psi]_1^{r}\right)[\Psi]_1^{\alpha}+[f]_{1}[\Psi]_{r+\alpha}+[f]_r[\Psi]_{1+\alpha}^r \ \text{ for } \ r\in \mathbb{N}\setminus\{0\}
    \end{align}
\end{subequations}
where the implicit constants depend on $r, \ N, \ n, \ \alpha$.
\end{prop}
\begin{proof}[Proof of \ref{compestimates}.] Following \cite{DS1}, we can write the Fa' di Bruno's formula symbolically as:
\begin{equation}\label{fadibruno}
    \DD^r (f\circ \Psi)=\sum_{l=1}^r(\DD ^lf)\circ \Psi\sum_\sigma C_{l,\sigma}(\DD \Psi)^{\sigma_1}\dots (\DD^r \Psi)^{\sigma_r}
\end{equation}
for some constant $C_{l,\sigma}$ , where the inner sum is on all the multi indexes $\sigma$ with:
$$\sum_{j=1}^rj\sigma_j=r \ \ \text{ and } \ \ \sum_{j=1}^r\sigma_j=l,$$
from standard interpolation inequalities, we have for $j\geq 1$ and $\alpha\in[0, 1)$:
\begin{equation*}
    \begin{split}
        [\Psi]_{j+\alpha}&\lesssim [\Psi]_{1+\alpha}^{1-\frac{(j+\alpha)-(1+\alpha)}{(r+\alpha)-(1+\alpha)}}[\Psi]_{r+\alpha}^{\frac{(j+\alpha)-(1+\alpha)}{(r+\alpha)-(1+\alpha)}}=[\Psi]_{1+\alpha}^{1-\frac{j-1}{r-1}}[\Psi]_{r+\alpha}^{\frac{j-1}{r-1}},\\
        [f]_{l+\alpha}&\lesssim [f]_{1+\alpha}^{1-\frac{l-1}{r-1}}[f]_{r+\alpha}^{\frac{l-1}{r-1}}\\
    \end{split}
\end{equation*}
and
$$ [(\DD ^lf)\circ \Psi]_\alpha\leq[f]_{l+\alpha}[\Psi]_1^\alpha$$
with this at hand, the product estimates in \eqref{prodestimates} and \eqref{fadibruno}, we can bound:
\begin{equation*}
    \begin{split}
        [f\circ \Psi]_{r+\alpha}&\lesssim \sum_{l=1}^r\sum_\sigma\left[[f]_{l}[\Psi]_{1+\alpha}^{\sigma_1}\dots [\Psi]_{r+\alpha}^{\sigma_r}+[(\DD ^lf)\circ \Psi]_\alpha[\Psi]_1^{\sigma_1}\dots [\Psi]_r^{\sigma_r}\right]\\
        &\leq \sum_{l=1}^r\sum_\sigma\left[[f]_l\prod_j[\Psi]_{j+\alpha}^{\sigma_j}+[f]_{l+\alpha}[\Psi]_1^\alpha\prod_j[\Psi]_{j}^{\sigma_j}\right]\\
        &\lesssim \sum_{l=1}^r\sum_\sigma\left[[f]_{1}^{\frac{r-l}{r-1}}[f]_{r}^{\frac{l-1}{r-1}}\prod_j[\Psi]_{1+\alpha}^{\sigma_j\frac{r-j}{r-1}}[\Psi]_{r+\alpha}^{\sigma_j\frac{j-1}{r-1}}+[\Psi]_1^\alpha[f]_{1+\alpha}^{\frac{r-l}{r-1}}[f]_{r+\alpha}^{\frac{l-1}{r-1}}\prod_j[\Psi]_{1}^{\sigma_j\frac{r-j}{r-1}}[\Psi]_{r}^{\sigma_j\frac{j-1}{r-1}}\right]\\
        &=\sum_{l=1}^r\sum_\sigma\left[[f]_{1}^{\frac{r-l}{r-1}}[f]_{r}^{\frac{l-1}{r-1}}[\Psi]_{1+\alpha}^{\sum_j\sigma_j\frac{r-j}{r-1}}[\Psi]_{r+\alpha}^{\sum_j\sigma_j\frac{j-1}{r-1}}+[\Psi]_1^\alpha[f]_{1+\alpha}^{\frac{r-l}{r-1}}[f]_{r+\alpha}^{\frac{l-1}{r-1}}[\Psi]_{1}^{\sum_j\sigma_j\frac{r-j}{r-1}}[\Psi]_{r}^{\sum_j\sigma_j\frac{j-1}{r-1}}\right]\\
        &\lesssim\sum_{l=1}^r\left[[f]_{1}^{\frac{r-l}{r-1}}[f]_{r}^{\frac{l-1}{r-1}}[\Psi]_{1+\alpha}^{r\frac{l-1}{r-1}}[\Psi]_{r+\alpha}^{\frac{r-l}{r-1}}+[\Psi]_1^\alpha[f]_{1+\alpha}^{\frac{r-l}{r-1}}[f]_{r+\alpha}^{\frac{l-1}{r-1}}[\Psi]_{1}^{r\frac{l-1}{r-1}}[\Psi]_{r}^{\frac{r-l}{r-1}}\right]\\
        &=\sum_{l=1}^r\left[\left([f]_{r}[\Psi]_{1+\alpha}^r\right)^{\frac{l-1}{r-1}}\left([f]_{1}[\Psi]_{r+\alpha}\right)^{\frac{r-l}{r-1}}+[\Psi]_1^\alpha\left([f]_{r+\alpha}[\Psi]_{1}^r\right)^{\frac{l-1}{r-1}}\left([f]_{1+\alpha}[\Psi]_{r}\right)^{\frac{r-l}{r-1}}\right]\\
        &=\sum_{l=1}^r\left[\left([f]_{r}[\Psi]_{1+\alpha}^r\right)^{\frac{l-1}{r-1}}\left([f]_{1}[\Psi]_{r+\alpha}\right)^{1-\frac{l-1}{r-1}}+[\Psi]_1^\alpha\left([f]_{r+\alpha}[\Psi]_{1}^r\right)^{\frac{l-1}{r-1}}\left([f]_{1+\alpha}[\Psi]_{r}\right)^{1-\frac{l-1}{r-1}}\right]\\
        &\leq \sum_{l=1}^r\left(\frac{l-1}{r-1}\right)[f]_{r}[\Psi]_{1+\alpha}^r+\left(1-\frac{l-1}{r-1}\right)[f]_{1}[\Psi]_{r+\alpha}\\
        &+\sum_{l=1}^r\left(\frac{l-1}{r-1}\right)[\Psi]_1^\alpha[f]_{r+\alpha}[\Psi]_{1}^r+\left(1-\frac{l-1}{r-1}\right)[\Psi]_1^\alpha[f]_{1+\alpha}[\Psi]_{r}\\
        &=\frac{r}{2}\left[[f]_{1+\alpha}[\Psi]_{r}[\Psi]_1^{\alpha}+[f]_{r+\alpha}[\Psi]_1^{r+\alpha}+[f]_{1}[\Psi]_{r+\alpha}+[f]_r[\Psi]_{1+\alpha}^r\right]\\
    \end{split}
\end{equation*}
where we used the weighted AMGM inequality.
\end{proof}

\noindent \textbf{Transport equation and Lagrangian Flow.} Fix $t_0, \tau \in \mathbb{R}$ with $\tau>0$, let $g:\mathbb{T}^3\times \mathbb{R}\to \mathbb{R}$ be a smooth time dependent function and $u:\mathbb{T}^3\times \mathbb{R}\to \mathbb{R}$ a smooth time dependent vector field. We consider the forced transport equation on $\mathbb{T}^3\times(t_0-\tau, \ t_0+\tau)$ for a function $f:\mathbb{T}^3\times (t_0-\tau, \ t_0+\tau)\to \mathbb{R}$, initial data $f_0:\mathbb{T}^3\to \mathbb{R}$ and forcing $g:\mathbb{T}^3\times (t_0-\tau, \ t_0+\tau)\to \mathbb{R}$, namely
\begin{equation}\label{standardtransport}
    \begin{cases}
        \partial_tf+u\cn f=g,\\
        f|_{t=t_0}=f_0
    \end{cases}
\end{equation}
and the equation for the Lagrangian flow $X=X_t(x):\mathbb{T}^3\times (t_0-\tau, \ t_0+\tau)\to \mathbb{T}^3$ of $u$, that is:
\begin{equation}\label{lagrangian}
    \begin{cases}
        \partial_tX_t(x)=u(X_t(x),t),\\
        X|_{t=t_0}(x)=x.
    \end{cases}
\end{equation}

\noindent The following result can be found in \cite{BuckmasterDeLellisSzekelyhidiVicol2019}.
\begin{prop}[Transport Estimate]\label{standardtransportestimate} Assume $\tau||u||_1\leq 1$. Given $u, g$, any solution $f$ of \eqref{standardtransport} satisfies:
$$||f(\cdot,t)||_\alpha\leq 2\left[||f_0||_\alpha+\int_{t_0}^t||(g(\cdot,s)||_\alpha\dd s\right],$$
for $\alpha\in[0,1]$. Moreover, for any $r\geq 1$ integer  
$$||f(\cdot,t)||_{r+\alpha}\lesssim||f_0||_{r+\alpha}+|t-t_0|[u]_r[f_0]_1+\int_{t_0}^t[(g(\cdot,s)]_{r+\alpha}+|t-s|[u]_{r+\alpha}[g(\cdot,s)]_1\dd s,$$
where the implicit constant depends on $r$ and $\alpha$. Consequently, the inverse flow $\Phi$ of $u$ at the identity at time $t_0$ satisfies:
\begin{equation*}
    \begin{split}
        ||\DD \Phi-\IId||_\alpha&\lesssim |t-t_0|[u]_{1+\alpha},\\
        [\Phi(t)]_{r+\alpha}&\lesssim |t-t_0|[u]_{r+\alpha} \ \text{ for } \ r\geq 2.
    \end{split}
\end{equation*}
\end{prop}

\noindent Similarly, we have the following classical result.
\begin{prop}[Lagrangian Flow Estimates]\label{standardlagrangianestimate} Let $\alpha\in [0,1)$ and $r\geq 0$ integer. Assume $\tau||u||_1\leq 1$, then any solution $X$ to \eqref{lagrangian} satisfies:
    \begin{equation*}
    \begin{split}
        ||\DD X-\IId||_\alpha&\lesssim |t-t_0|[u]_{1+\alpha},\\
        [X(t)]_{r+\alpha}&\lesssim |t-t_0|[u]_{r+\alpha} \ \text{ for } \ r\geq 2 \ \text{integer},
    \end{split}
\end{equation*}
where the implicit constant depends on $r$ and $\alpha$.
\end{prop}
For a proof, see \cite[Chapter 3]{bertozzimajda}. The higher-order derivative estimates follow as in the proof of Lemma \ref{lagrangiang}.


\noindent \textbf{Mollification Estimates.}  Let $v, \ w, \ f:\mathbb{T}^3\times \mathbb{R}\to \mathbb{R}^3$ be time-dependent vector fields with $v, \ w, \ f, \ g\in C^N(\mathbb{T}^3,\mathbb{R}^3)$, define a transport operator $A$ by 
$$A=\partial_t+w\cn$$ 
and assume in addition that $Av, \ Af, \ Ag\in C^{N-1}(\mathbb{T}^3,\mathbb{R}^3)$.

\noindent Let $\rho\in C_c^\infty(\mathbb{R}^3)$ be a convolution kernel with vanishing moments up to order $\bar m_0\geq 0$, namely, we require:
\begin{equation}\label{momentumconditions}
    \begin{cases}
        \int_{\mathbb{R}^3} y^\sigma\rho(y) \ \dd y=0 \ \text{ for any multi-index }\ \sigma \ \text{ with } 1\leq |\sigma| \leq \bar m_0 \ &\text{ if }\  \bar m_0\geq1,\\
        \text{ no condition } \ &\text{ if } \ \bar m_0=0.
    \end{cases}
\end{equation}

\noindent For any given $\ell>0$ we rescale the kernel as 
$$\rho_\ell(y)=1/\ell^3\rho(y/\ell)$$
and define the mollification $(f)_\ell=(f*\rho_\ell)$ of $f$ at scale $\ell$ to be
\begin{equation}\label{mollnotation}
     (f)_{\ell}(x)=\int_{\mathbb{R}^3} f_{per}(y)\rho_{\ell}(x-y)\dd y
\end{equation}
where $f_{per}$ denotes the periodic extension to $\mathbb{R}^3$ of $f$. 

\noindent The following is a collection of classical results; see, for example, the deep smoothing operators in \cite[Section 2.3.4]{gromov1986partial} and the Constantin, Weinan, Titi (CET) commutator estimates \cite{CET}, where we adopt the formulation from \cite{GR} coming from \cite{CDS}. Note, however, the bounds related to the transport operator $A$. 

\begin{prop}[Deep Mollification]\label{deepmollification}Let $0\leq m_0\leq \bar m_0$ and $r,r'\geq 0$ integers. We have:
\begin{subequations}
    \begin{align}
        ||(f)_\ell||_{r+r'}&\lesssim \ell^{r'}||f||_{r} \ \text{ for }\ 0\leq r\leq N,\label{deep0}\\
        ||f-(f)_\ell||_r&\lesssim \ell^{m_0+1}||f||_{r+m_0+1} \ \text{ for }\ 0\leq r+m_0+1\leq N\label{deep1}
    \end{align}
\end{subequations}
and the following first- and second-order commutator estimates with transport operators:
\begin{equation}\label{deep2}
    \begin{split}
        ||[v\cn,*\rho_\ell]f||_r&\lesssim  \ell^{m_0+1}(||v||_{r+m_0+1}||f||_1+||v||_{1}||f||_{r+m_0+1})\\
    \end{split}
\end{equation}
for $ 0\leq r+m_0+1\leq N$  and 
\begin{equation}
    \begin{split}
        ||[A,[v\cn,*\rho_\ell]]f||_r&\lesssim \ell^{m_0+1}(||Av||_{r+m_0+1}||f||_1+||Av||_{1}||f||_{r+m_0+1})\\
        &+\ell^{m_0+1}(||w||_{r+m_0+1}||v||_{1}||f||_1+||w||_1||v||_{r+m_0+1}||f||_1+||w||_{1}||v||_{1}||f||_{r+m_0+1})\\
        &+\ell^{m_0+2}(||w||_{r+m_0+1}||v||_{1}||f||_2+||w||_1||v||_{r+m_0+1}||f||_2+||w||_{1}||v||_{1}||f||_{r+m_0+2})
    \end{split}
\end{equation}
for $0\leq r+m_0+2\leq N$. 

\noindent Moreover, the following bound on the CET commutator holds:
\begin{equation}\label{deepCWT}
    ||(f)_\ell (g)_\ell-(fg)_\ell||_{r}\lesssim \ell^{2+r'-r}(||f||_{r'+1}||g||_1+||f||_1||g||_{r'+1}) \ \text{ for }\ 0\leq r'+1\leq N
\end{equation}
and assuming in addition $N\geq 2$ and $w\in C^{N'}$ for $N'\geq N$, we have the following associated transport estimate:
\begin{equation}
    \begin{split}
        &||A[(fg)_\ell-(f)_\ell(g)_\ell]||_{r+r'}\\
        &\lesssim \ell^{m_0+2-r'}(||Af||_{r+m_0+1}||g||_1+||Af||_{1}||g||_{r+m_0+1})\\
        &+\ell^{m_0+2-r'}(||f||_{1}||Ag||_{r+m_0+1}+||f||_{r+m_0+1}||Ag||_{1})\\
        &+\ell^{m_0+2-r'}(||w||_{r+m_0+1}||f||_{1}||g||_1+||w||_1||f||_{r+m_0+1}||g||_1+||w||_{1}||f||_{1}||g||_{r+m_0+1})\\
        &+\ell^{m_0+2}(||w||_{r+r'}||f||_{m_0+2}||g||_1+||w||_{r+r'}||f||_1||g||_{m_0+2})
    \end{split}
\end{equation}
for $0\leq r+r'\leq N'$ and $0\leq r+m_0+1\leq N$ and $m_0+2\leq N$.

\noindent The implicit constants depend on the specific choice of kernel $\rho$ and $r, \ r',\ m_0$.
\end{prop}

\begin{remark} Note that requiring no moments on the kernel that is $\bar m_0=m_0=0$, one in particular recovers the classical estimates. 
\end{remark}

\begin{proof}[Proof of Proposition \ref{deepmollification}] Abusing notation, we do not differentiate functions from their periodic extensions. The first two estimates \eqref{deep0} \eqref{deep1} are standard; see \cite[Lemma 7.2]{burczak_szekelyhidi_wu_2023} for a proof. The third estimate is also well known, see \cite[Proposition A.3]{GR}. 

\noindent For any $v:\mathbb{T}^3\times \mathbb{R}\to \mathbb{R}^3$ and $y\in \mathbb{R}^3$ we write:
$$\delta_yv(x)=v(x-y)-v(x).$$

\noindent \textbf{First-order commutator.} Without loss of generality, we may assume that $f$ is a scalar function. This  reads:
\begin{equation}\label{firstordercomm}
    \begin{split}
        \left([v\cn,*\rho_\ell]\right)f(x)&=(v\cn f_\ell)(x)-(v\cn f)_\ell(x)\\
        &=v(x)\cn\int f(y)\rho_\ell (x-y)\dd y-\int (v\cn f)(y)\rho_\ell(x-y)\dd y\\
        &=-v(x)\cn\int f(x-y)\rho_\ell (y)\dd y+\int (v\cn f)(x-y)\rho_\ell(y)\dd y\\
        &=\int \rho_\ell(y)(v(x-y)-v(x))\cn f(x-y)\dd y\\
        &=\int \rho_\ell(y)(\delta_yv)(x)\cn f(x-y)\dd y
    \end{split}
\end{equation}
We will not provide estimates for this at this point, as the second-order commutator will contain a term of similar structure; we refer to $T_0, \ T_1$ below, see \eqref{deepcommutator1}.

\noindent \textbf{Second-order commutator.} Using the expression in \eqref{firstordercomm} we compute: 
\begin{equation*}
    \begin{split}
        (A[v\cn,*\rho_\ell]f)(x)&=\partial_t\int \rho_\ell(y)\delta_yv(x)\cn f(x-y)\dd y+w(x)\cn \int \rho_\ell(y)\delta_yv(x)\cn f(x-y)\dd y\\
        &=\int \rho_\ell(y)(\partial_t\delta_yv)(x)\cn f(x-y)\dd y+\int \rho_\ell(y)(\delta_yv)(x)\cn \partial_tf(x-y)\dd y\\
        &+\int \rho_\ell(y)w^k(x)(\partial_k\delta_yv)^i(x)\partial_i f(x-y)\dd y+\int \rho_\ell(y)w^k(x)(\delta_yv)^i(x)\partial_{i}\partial_kf(x-y)\dd y\\
        &=\underbrace{\int \rho_\ell(y)(\delta_yAv)(x)\cn f(x-y)\dd y}_{([A v\cn,*\rho_\ell]f)(x)}+\underbrace{\int \rho_\ell(y)(\delta_yv)(x)\cn (Af)(x-y)\dd y}_{([v\cn,*\rho_\ell]A f)(x)}\\
        &-\int \rho_\ell(y)[(\delta_yv)(x)\cn w(x-y)+(\delta_yw)(x)\cn v(x-y)]\cn f(x-y)\dd y\\
        &-\int \rho_\ell(y)\DD ^2f(x-y)[(\delta_yv)(x),(\delta_yw)(x)]\dd y\\
    \end{split}
\end{equation*}
where, we used the following two calculations:
\begin{equation*}
    \begin{split}
        &\int \rho_\ell(y)w^k(x)(\partial_k\delta_yv)^i(x)\partial_i f(x-y)\dd y\\
        &=\int \rho_\ell(y)\left[w^k(x-y)\partial_kv^i(x-y)-w^k(x)\partial_kv^i(x)-(\delta_yw)^k(x)\partial_k v^i(x-y)\right]\partial_i f(x-y)\dd y\\
        &=\int \rho_\ell(y)\left[(\delta_y(w\cn v))^i(x)-(\delta_yw)(x)\cn v^i(x-y)\right]\partial_i f(x-y)\dd y\\
        &=\int \rho_\ell(y)\left[\delta_y(w\cn v)(x)-(\delta_yw)(x)\cn v(x-y)\right]\cn f(x-y)\dd y,\\
    \end{split}
\end{equation*}
and
\begin{equation*}
    \begin{split}
        &\int \rho_\ell(y)w^k(x)(\delta_yv)^i(x)\partial_{i}\partial_kf(x-y)\dd y\\
        &=\int \rho_\ell(y)[w^k(x-y)(\delta_yv)^i(x)-(\delta_yw)^k(x)(\delta_yv)^i(x)]\partial_{i}\partial_kf(x-y)\dd y\\
        &=\int \rho_\ell(y)[(\delta_yv)^i(x)\partial_{i}(w^k(x-y)\partial_kf(x-y))-(\delta_yv)^i(x)\partial_{i}w^k(x-y)\partial_kf(x-y)]\dd y\\
        &-\int \rho_\ell(y)\DD ^2f(x-y)[(\delta_yv)(x),(\delta_yw)(x)]\dd y\\
        &=\int \rho_\ell(y)(\delta_yv)(x)\cn  ((w\cn f)(x-y))\dd y-\int \rho_\ell(y)[(\delta_yv)(x)\cn w(x-y)]\cn f(x-y)\dd y\\
        &-\int \rho_\ell(y)\DD ^2f(x-y)[(\delta_yv)(x),(\delta_yw)(x)]\dd y.\\
    \end{split}
\end{equation*}
We deduce that:
\begin{equation*}
    \begin{split}
        ([A,[v\cn,*\rho_\ell]]f)(x)&=(A[v\cn,*\rho_\ell]f)(x)-([v\cn,*\rho_\ell]Af)(x)\\
        &=\underbrace{([Av\cn,*\rho_\ell]f)(x)}_{T_0}\\
        &-\underbrace{\int \rho_\ell(y)[(\delta_yv)(x)\cn w(x-y)+(\delta_yw)(x)\cn v(x-y)]\cn f(x-y)\dd y}_{T_1}\\
        &-\underbrace{\int \DD ^2f(x-y)[(\delta_yv)(x),(\delta_yw)(x)]\dd y}_{T_2}.\\
    \end{split}
\end{equation*}

\noindent \textit{Estimates on $T_1$.} We expand in a Taylor series up to order $m_0$ and write:
\begin{equation*}
    \begin{split}
         h_x(y)&=[(v(x-y)-v(x))\cn w(x-y)+(w(x-y)-w(x))\cn v(x-y)]\cn f(x-y)\\
        &=\sum_{0<|\sigma|\leq m_0} \frac{1}{\sigma!}\partial_\sigma h_x(0)(y)^\sigma+r^{m_0}_{x}(y)
    \end{split}
\end{equation*}
where
\begin{equation*}
    r^{m_0}_{x}(y)=\sum_{|\sigma|= m_0+1}\frac{1}{\sigma!}y^\sigma\int_0^1 \partial_\sigma h_x(sy)(1-s)^{m_0}\dd s,
\end{equation*}
note the missing constant term in the expansion, with the understanding that for $m_0=0$ the sum is just empty. We can bound:
\begin{equation*}
    \begin{split}
        ||\partial_\sigma h_x(y)||_{C^{r}_x}&\leq ||h_x(y)||_{C^{m_0+1+r}_x}\\
        &\lesssim||w||_{r+m_0+1}||v||_{1}||f||_1+||w||_1||v||_{r+m_0+1}||f||_1+||w||_{1}||v||_{1}||f||_{r+m_0+1}\\
    \end{split}
\end{equation*}
and from this we deduce:
$$||r^{m_0}_{x}(y)||_{C^{r}_x}\lesssim |y|^{m_0+1}(||w||_{r+m_0+1}||v||_{1}||f||_1+||w||_1||v||_{r+m_0+1}||f||_1+||w||_{1}||v||_{1}||f||_{r+m_0+1}),$$
given the zero momentum assumption on $\rho$, see \eqref{momentumconditions}, and the missing constant term in the series, we conclude that:
\begin{equation}
    \begin{split}
        ||T_1||_r&\leq ||\underbrace{\int\rho_\ell(y)\sum_{0<|\sigma|\leq m_0} \partial_\sigma h_x(0)y^\sigma\dd y}_{=0}||_r+\int|\rho_\ell(y)|||r^{m_0}_{x}(y)||_{C^r_x}\dd y\\
        &\leq (||w||_{r+m_0+1}||v||_{1}||f||_1+||w||_1||v||_{r+m_0+1}||f||_1+||w||_{1}||v||_{1}||f||_{r+m_0+1})\int|y|^{m_0+1}|\rho_\ell(y)|\dd y\\
        &\lesssim \ell^{m_0+1}(||w||_{r+m_0+1}||v||_{1}||f||_1+||w||_1||v||_{r+m_0+1}||f||_1+||w||_{1}||v||_{1}||f||_{r+m_0+1})\\
    \end{split}
\end{equation}
where the implicit constant depends on the specific choice of kernel $\rho,\ r$ and $m_0$.

\noindent \textit{Estimates on $T_0$.} The same argument we used for $T_1$ can be used here and we deduce:
\begin{equation}\label{deepcommutator1}
    ||T_0||_r\lesssim \ell^{m_0+1}(||Av||_{r+m_0+1}||f||_1+||Av||_{1}||f||_{r+m_0+1})
\end{equation}
where the implicit constant depends on the specific choice of kernel $\rho,\ r$ and $m_0$.

\noindent \textit{Estimates on $T_2$.} Here we have $O(|y|^2)$ decay despite the pure $C^2$ norm of $f$ appearing. Proceeding as above, we first write:
\begin{equation*}
    \begin{split}
         h_x(y)&=\DD ^2f(x-y)[v(x-y)-v(x),w(x-y)-w(x)]\\
        &=\sum_{0<|\sigma|\leq m_0} \frac{1}{\sigma!}\partial_\sigma h_x(0)(y)^\sigma+r^{m_0}_{x}(y),
    \end{split}
\end{equation*}
then given the better decay at $|y|\to0$ deduce:
$$||r^{m_0}_{x}(y)||_{C^{r}_x}\lesssim |y|^{m_0+2}(||w||_{r+m_0+1}||v||_{1}||f||_2+||w||_1||v||_{r+m_0+1}||f||_2+||w||_{1}||v||_{1}||f||_{r+m_0+2})$$
and conclude from \eqref{momentumconditions} that:
\begin{equation}\label{T2above}
    \begin{split}
        ||T_2||_r&\lesssim \ell^{m_0+2}(||w||_{r+m_0+1}||v||_{1}||f||_2+||w||_1||v||_{r+m_0+1}||f||_2+||w||_{1}||v||_{1}||f||_{r+m_0+2})\\
    \end{split}
\end{equation}
where the implicit constant depends on the specific choice of kernel $\rho,\ r$ and $m_0$.

\noindent \textit{Conclusion.} Collecting the bounds above we deduce that for any $0\leq m_0\leq \bar m_0$ we have:
\begin{equation*}
    \begin{split}
        ||[A,[v\cn,*\rho_\ell]]f||_r&\lesssim \ell^{m_0+1}(||Av||_{r+m_0+1}||f||_1+||Av||_{1}||f||_{r+m_0+1})\\
        &+\ell^{m_0+1}(||w||_{r+m_0+1}||v||_{1}||f||_1+||w||_1||v||_{r+m_0+1}||f||_1+||w||_{1}||v||_{1}||f||_{r+m_0+1})\\
        &\lesssim \ell^{m_0+2}(||w||_{r+m_0+1}||v||_{1}||f||_2+||w||_1||v||_{r+m_0+1}||f||_2+||w||_{1}||v||_{1}||f||_{r+m_0+2})
    \end{split}
\end{equation*}
for $0\leq r +m_0+2\leq N$, where the implicit constant depends on the specific choice of kernel $\rho,\ r$ and $m_0$.

\noindent \textbf{Transport of the CWT commutator.} We first recall the classical rewriting:
\begin{equation}\label{CWTrewriting}
    \begin{split}
        ((fg)_\ell-(f)_\ell(g)_\ell)(x)=\underbrace{\int\rho_\ell(y)(\delta_yf)(x)(\delta_yg)(x)\dd y}_{T_2}-\underbrace{[(f)_\ell-f)](x)[(g)_\ell-g](x)}_{T_1},
    \end{split}
\end{equation}
we then have
$$A[(f)_\ell-f]=(Af)_\ell-Af+[w\cn, *\rho_\ell]f$$
and conclude from \eqref{deep1} and \eqref{deep2} that
\begin{equation*}
    \begin{split}
        ||AT_1||_r&\lesssim \ell^{m_0+2}(||Af||_{r+m_0+1}||g||_1+||Af||_{1}||g||_{r+m_0+1})\\
        &+\ell^{m_0+2}(||f||_{1}||Ag||_{r+m_0+1}+||f||_{r+m_0+1}||Ag||_{1})\\
        &+\ell^{m_0+2}(||w||_{r+m_0+1}||f||_{1}||g||_1+||w||_1||f||_{r+m_0+1}||g||_1+||w||_{1}||f||_{1}||g||_{r+m_0+1}),\\
    \end{split}
\end{equation*}
we then compute:
\begin{equation*}
    \begin{split}
        AT_2&=\int\rho_\ell(y)[(\delta_yAf)(x)(\delta_yg)(x)+(\delta_yf)(x)(\delta_yAg)(x)]\dd y\\
        &-\int\rho_\ell(y)[(\delta_yw)(x)\cn f(x-y)(\delta_yg)(x)+(\delta_yf)(x)(\delta_yw)(x)\cn g(x-y)]\dd y,\\
    \end{split}
\end{equation*}
we can now argue as in the proof for the commutator bounds above, to deduce:
\begin{equation}
    \begin{split}
        ||AT_2||_r&\lesssim \ell^{m_0+2}(||Af||_{r+m_0+1}||g||_1+||Af||_{1}||g||_{r+m_0+1})\\
        &+\ell^{m_0+2}(||f||_{1}||Ag||_{r+m_0+1}+||f||_{r+m_0+1}||Ag||_{1})\\
        &+\ell^{m_0+2}(||w||_{r+m_0+1}||f||_{1}||g||_1+||w||_1||f||_{r+m_0+1}||g||_1+||w||_{1}||f||_{1}||g||_{r+m_0+1}).\\
    \end{split}
\end{equation}
and gathering the two bounds, we conclude: 
\begin{equation*}
    \begin{split}
        ||A[(fg)_\ell-(f)_\ell(g)_\ell)]||_r&\leq ||AT_1||_r+||AT_2||_r\\
        &\lesssim \ell^{m_0+2}(||Af||_{r+m_0+1}||g||_1+||Af||_{1}||g||_{r+m_0+1})\\
        &+\ell^{m_0+2}(||f||_{1}||Ag||_{r+m_0+1}+||f||_{r+m_0+1}||Ag||_{1})\\
        &+\ell^{m_0+2}(||w||_{r+m_0+1}||f||_{1}||g||_1+||w||_1||f||_{r+m_0+1}||g||_1+||w||_{1}||f||_{1}||g||_{r+m_0+1})\\
    \end{split}
\end{equation*}
for $0\leq r +m_0+1\leq \bar N$. 

\noindent For higher derivatives, we can follow \cite{CDS} and use the fact that
\begin{equation*}
f(x)*\partial^\theta\rho_\ell=
    \begin{cases}
        f(x) & \ \text{\ if \ } |\theta|=0,\\
        0 & \ \text{\ else \ }
    \end{cases}
\end{equation*}
to write 
\begin{equation}\label{CWTrewriting1}
    \begin{split}
        \partial^\theta[(fg)_\ell-(f)_\ell(g)_\ell]&=[\delta f\delta g]*\partial^\theta\rho_\ell -\sum_{\theta'\leq \theta}\binom{\theta}{\theta'}(\delta f*\partial^{\theta'}\rho_\ell)(\delta g*\partial^{\theta-\theta'}\rho_\ell)
    \end{split}
\end{equation}
where $\theta$ is any multi index with $|\theta|=r'$. We can now commute
\begin{equation*}
    \begin{split}
        \partial^\theta A[(fg)_\ell-(f)_\ell(g)_\ell]= A\partial^\theta[(fg)_\ell-(f)_\ell(g)_\ell]+[\partial^\theta, w\cn][(fg)_\ell-(f)_\ell(g)_\ell]
    \end{split}
\end{equation*}
The second term using \eqref{deepCWT} can be estimated, for $r'\geq 1$, as
\begin{equation*}
    \begin{split}
        ||[\partial^\theta, w\cn][(fg)_\ell-(f)_\ell(g)_\ell]||_{r}&\lesssim ||w||_{r+r'}||(fg)_\ell-(f)_\ell(g)_\ell||_1+||w||_1||(fg)_\ell-(f)_\ell(g)_\ell||_{r+r'}\\
        &\lesssim \ell^{2+m_0-r'}(||w||_1||f||_{r+m_0+1}||g||_1+||w||_1||f||_1||g||_{r+m_0+1})\\
        &+\ell^{2+m_0}(||w||_{r+r'}||f||_{m_0+2}||g||_1+||w||_{r+r'}||f||_1||g||_{m_0+2}).
    \end{split}
\end{equation*}
The estimate on the first term follows, arguing as above, where we use \eqref{CWTrewriting1} instead of \eqref{CWTrewriting}, we deduce
\begin{equation*}
    \begin{split}
        ||A\partial^\theta[(fg)_\ell-(f)_\ell(g)_\ell]||_r&\lesssim \ell^{m_0+2-r'}(||Af||_{r+m_0+1}||g||_1+||Af||_{1}||g||_{r+m_0+1})\\
        &+\ell^{m_0+2-r'}(||f||_{1}||Ag||_{r+m_0+1}+||f||_{r+m_0+1}||Ag||_{1})\\
        &+\ell^{m_0+2-r'}(||w||_{r+m_0+1}||f||_{1}||g||_1+||w||_1||f||_{r+m_0+1}||g||_1+||w||_{1}||f||_{1}||g||_{r+m_0+1})\\
    \end{split}
\end{equation*}
and since the multi-index was arbitrary, summing the estimates, we conclude:
\begin{equation*}
    \begin{split}
        ||A[(fg)_\ell-(f)_\ell(g)_\ell]||_{r+r'}&\lesssim \ell^{m_0+2-r'}(||Af||_{r+m_0+1}||g||_1+||Af||_{1}||g||_{r+m_0+1})\\
        &+\ell^{m_0+2-r'}(||f||_{1}||Ag||_{r+m_0+1}+||f||_{r+m_0+1}||Ag||_{1})\\
        &+\ell^{m_0+2-r'}(||w||_{r+m_0+1}||f||_{1}||g||_1+||w||_1||f||_{r+m_0+1}||g||_1+||w||_{1}||f||_{1}||g||_{r+m_0+1})\\
        &+\ell^{2+m_0}(||w||_{r+r'}||f||_{m_0+2}||g||_1+||w||_{r+r'}||f||_1||g||_{m_0+2}).
    \end{split}
\end{equation*}
for $0\leq r+r'\leq N'$ and $0\leq r+m_0+1\leq N$ and $m_0+2\leq N$.
\end{proof}


\noindent \textbf{Singular Integral Operators.} Let $K:\mathbb{R}^3\to \mathbb{R}$ be a kernel which is smooth away from the origin, is homogeneous of degree $-2$ and has zero mean on circles centred at the origin. Let $K_{per}$ be the periodization of $K$, namely
$$K_{per}(y)=K(y)+\sum_{n\in \mathbb{Z}^3\setminus{0}}\left[K(y+n)-K(n)\right]$$
and define an operator $T_K$ acting on zero-mean periodic functions $f:\mathbb{T}^3\to \mathbb{R}$ as:
\begin{equation}\label{TK}
    T_Kf(y)=p.v.\int_{\mathbb{T}^3} K_{per}(x-y)f(y)\dd y.
\end{equation}
The following result is standard; see, for example, \cite{BuckmasterDeLellisSzekelyhidiVicol2019} for a proof.
\begin{prop}[Properties of Calderon-Zygm\'{u}nd operators]\label{czstuff} Let $\alpha\in(0,1)$ and $r\geq 0$ an integer. Let $T_K$ be a periodic Calderon-Zygm\'{u}nd operator with kernel $K$ as in \eqref{TK}. This is bounded in the space of functions $f$ which are $C^{r+\alpha}(\mathbb{T}^3)$ and have zero mean, namely
$$||T_Kf||_{r+\alpha}\lesssim ||f||_{r+\alpha}$$
where the implicit constant depends on $r, \ \alpha, \ K$. Moreover, let $b\in C^{r+1+\alpha}(\mathbb{T}^3)$ be a vector field, we have:
$$||[T_K,b\cn]g||_{r+\alpha}\lesssim||b||_{r+1+\alpha}||g||_{\alpha}+||b||_{1+\alpha}||g||_{r+\alpha}$$
for $g\in C^{r+\alpha}(\mathbb{T}^3)$, where the implicit constant depends on $r, \ \alpha, \ K$.
\end{prop}

\begin{remark} In this work we will often use that $\nabla\Delta^{-1}\ddiv, \ \mathcal{R}\ddiv, \ \mathcal{R}\curl$ and $\partial_i\mathcal{R}$ for $i=1,2,3$ are all CZ operators, where $\mathcal{R}$ is defined in \eqref{invdiv} below.
\end{remark}


\noindent \textbf{Convex Integration Toolbox.} In this work, we use the following anti-divergence operator: 
\begin{equation}\label{invdiv}
    \begin{split}
        &\mathcal{R}v=\frac{1}{4}\left[\DD \mathbb{P}u+(\DD \mathbb{P}u)^\top\right]+\frac{3}{4}\left[\DD u+(\DD u)^\top\right]-\frac{1}{2}(\ddiv \ u)\IId\\
    \end{split}
\end{equation}
where  $u$ is the unique solution of
\begin{equation*}
    \begin{cases}
        \Delta u=v-\int_{\mathbb{T}^3}v,\\
        \int_{\mathbb{T}^3}u=0.
    \end{cases}
\end{equation*}
This operator maps mean-zero vector fields $v$ on the torus to symmetric trace-free 2-tensors $\mathcal{R}v$. This is a standard tool in convex integration, which originally appeared in \cite{DS}. We refer to \cite[Lemma 2.2]{DaSz} for the classical statement describing its properties and estimates; since our constructions are local in space, we will need to adapt the arguments to this setting and do so directly in the proof of Lemma \ref{Rpqua}, which provides bounds for the oscillation term in the Nash Stage.

\bibliographystyle{plain}
\bibliography{refs}

\end{document}